\documentclass[11pt,reqno,a4paper]{amsart}
\usepackage[T1]{fontenc}
\usepackage[utf8]{inputenc}
\usepackage{lmodern,microtype}
\usepackage[margin=27mm]{geometry}
\usepackage{mathtools,amssymb,mathrsfs}
\usepackage{tikz}
\usepackage{booktabs,longtable,array,enumitem}
\usepackage[hidelinks,bookmarksnumbered=true]{hyperref}
\usepackage{bookmark}

\numberwithin{equation}{section}
\theoremstyle{plain}
\newtheorem{theorem}{Theorem}[section]

\newtheorem{proposition}{Proposition}[section]
\newtheorem{lemma}[proposition]{Lemma}
\newtheorem{corollary}[proposition]{Corollary}

\theoremstyle{definition}

\theoremstyle{remark}
\newtheorem{remark}[proposition]{Remark}

\newcommand{\R}{\mathbb R}
\newcommand{\C}{\mathbb C}
\newcommand{\N}{\mathbb N}
\newcommand{\Z}{\mathbb Z}

\newcommand{\Hh}{\mathcal H}

\DeclareMathOperator{\Cov}{Cov}
\DeclareMathOperator{\Var}{Var}
\DeclareMathOperator{\ind}{ind}

\DeclareMathOperator{\diag}{diag}
\DeclareMathOperator{\sgn}{sgn}
\DeclareMathOperator{\Id}{Id}

\DeclareMathOperator{\Ree}{Re}
\DeclareMathOperator{\Imm}{Im}
\setlist{nosep}
\title[Degree of singular \(SU(3)\) Toda systems]
{Degree counting for singular \(SU(3)\) Toda systems}

\author[Y. Lee]{Youngae Lee}
\address[Youngae Lee]{Department of Mathematical Sciences,
Ulsan National Institute of Science and Technology (UNIST),
Ulsan, Republic of Korea}
\email{youngaelee@unist.ac.kr}

\author[H. Li]{Haoyu Li}
\address[Haoyu Li]{Department of Mathematics, Faculty of Science and Technology, University of Macau, Taipa, Macau, China}
\email{hyli1994@hotmail.com}

\author[J. Wei]{Juncheng Wei}
\address[Juncheng Wei]{Department of Mathematics,
Chinese University of Hong Kong, Shatin, NT, Hong Kong}
\email{wei@math.cuhk.edu.hk}

\author[W. Yang]{Wen Yang}
\address[Wen Yang]{Department of Mathematics, Faculty of Science and Technology,
University of Macau, Taipa, Macau, China}
\email{wenyang@um.edu.mo}

\author[L. Zhang]{Lei Zhang}
\address[Lei Zhang]{Department of Mathematics, University of Florida,
Gainesville, FL 32611, USA}
\email{leizhang@ufl.edu}

\date{}

\begin{document}

\begin{abstract}
We extend our Leray--Schauder degree formula for the regular $SU(3)$
Toda system on a closed surface~\cite{LLWYZ2026} to singular sources
of arbitrary strengths greater than $-1$. Outside the exceptional
sets determined by the local masses, the degree at $(\rho_1,\rho_2)$
is the sum of the coefficients of $x^ay^b$ with $4\pi a<\rho_1$ and
$4\pi b<\rho_2$ in an explicit series with real exponents. On a torus,
the regular degree is $1$, and a single source with positive noninteger
strengths can give degree $-2$.

At mixed singular blow-up points, we choose two coefficients in an
exact local Toda profile to match both exterior fields to first order.
When all outer profile exponents are less than one, we perturb the
first weight away from the marked points to obtain degree $-1$ for
the reduced scale equation on the bubbling side.
\end{abstract}

\maketitle

\tableofcontents

\section{Introduction}
\label{sec:introduction}
Let $(M,g)$ be a smooth closed connected oriented Riemannian surface
of area one and let $S\subset M$ be a finite set of marked points. For $p\in M$, denote by $\delta_p$ the unit Dirac measure supported at $p$.
In this paper, we compute the Leray--Schauder degree of the following singular $SU(3)$ Toda
system for $v=(v_1,v_2)$:
\begin{equation}\label{eq:toda-dirac}
 \left\{
 \begin{aligned}
 &-\Delta_g v_i
 = \sum_{j=1}^2 K_{ij}\rho_j
    \left(
      \frac{h_j^*e^{v_j}}{\displaystyle\int_M h_j^*e^{v_j}\,dV_g}-1
    \right)
    -4\pi\sum_{p\in S}\alpha_{p,i}(\delta_p-1)
    \text{ in }M,\\
 &\int_M v_i\,dV_g = 0,
 \,\,i=1,2.
 \end{aligned}
 \right.
\end{equation}
Here $\Delta_g$ is the Laplace--Beltrami operator with
$-\Delta_g\geq0$ and $dV_g$ is the Riemannian area measure.
The coefficient matrix is the Cartan matrix of type $A_2$ such that
\[
 K=(K_{ij})_{1\leq i,j\leq2}
   =\begin{pmatrix}2&-1\\-1&2\end{pmatrix}.
\]
For $i=1,2$, the weight $h_i^*$ is a smooth positive function on $M$,
$\rho_i>0$ is a parameter and $\alpha_{p,i}>-1$ is the strength of the source at $p\in S$
in the $i$-th component. We write
\[
 \rho=(\rho_1,\rho_2),\qquad
 \alpha=(\alpha_{p,i})_{p\in S,\,i=1,2},
\]
and denote by $\chi=\chi(M)$ the Euler characteristic of $M$.

We keep the set $S$ of marked points fixed as the source strengths vary.
For fixed source strengths, we may omit from $S$ any point $p$ for which
$\alpha_{p,1}=\alpha_{p,2}=0$. For each $y\in M$, let $G(\cdot,y)$ be the Green function of
$-\Delta_g$ defined as
\[
 -\Delta_g G(\cdot,y)=\delta_y-1,\qquad
 \int_M G(x,y)\,dV_g(x)=0.
\]
For $q>1$, let $L^q(M)$ denote the Lebesgue space with respect to
$dV_g$. For $f\in L^q(M)$, define its Green potential by
\[
 (G*f)(x):=\int_M G(x,y)f(y)\,dV_g(y),\qquad x\in M.
\]
With the above normalization,
\[
 -\Delta_g(G*f)=f-\int_M f\,dV_g,\qquad
 \int_M(G*f)\,dV_g=0.
\]

Fix $p\in M$. In a local isothermal coordinate $z$ centered at $p$,
with $z(p)=0$ and $g=e^\varphi|dz|^2$, we write
\[
 G(x,p)=-\frac{1}{2\pi}\log|z-p|+R(x,p),
\]
where $R(\cdot,p)$ is the regular part of $G(\cdot,p)$ in this
coordinate. Let $d_g$ denote the geodesic distance on $M$. The Robin
function is given by
\[
 r_g(p):=\lim_{x\to p}
 \left(G(x,p)+\frac{1}{2\pi}\log d_g(x,p)\right).
\]
Since $d_g(x,p)/|z(x)|\to e^{\varphi(p)/2}$ as $x\to p$, we obtain
$r_g(p)=R(p,p)+\frac{\varphi(p)}{4\pi}$.
Define the weights of the transformed system by
\begin{equation}\label{eq:physical-weights}
h_i:=h_i^*\exp\left(-4\pi\sum_{p\in S}
\alpha_{p,i}G(\cdot,p)\right),
 \qquad i=1,2,
\end{equation}
and introduce the transformed unknowns
\[
u_i:=v_i+4\pi\sum_{p\in S}\alpha_{p,i}G(\cdot,p),
 \qquad i=1,2.
\]
Since $h_i e^{u_i}=h_i^*e^{v_i}$ almost everywhere and
$G(\cdot,p)$ has zero mean, \eqref{eq:toda-dirac} is equivalent to
\begin{equation}\label{eq:toda}
 \left\{
 \begin{aligned}
 &-\Delta_g u_i
=\sum_{j=1}^2 K_{ij}\rho_j
   \left(\frac{h_j e^{u_j}}{\int_M h_j e^{u_j}\,dV_g}-1\right)\text{ in }M,\\
 &\int_M u_i\,dV_g=0,\,\,i=1,2.
 \end{aligned}
 \right.
\end{equation}

Fix $p\in S$ and choose a conformal coordinate $z$ centered at $p$.
For each $i=1,2$, there is a smooth positive function
$\widetilde h_i$ near $p$ such that
\[
 h_i(z)=|z|^{2\alpha_{p,i}}\widetilde h_i(z).
\]
Let $S_i$ be the set of marked points with nonzero source strength in the
$i$th component, and let $\ell$ be the number of marked points:
\[
 S_i:=\{p\in S:\alpha_{p,i}\ne0\}\quad\mbox{for }i=1,2,\qquad \ell:=|S|.
\]
For each $p\in S$, write
\[
 A_p:=1+\alpha_{p,1},\qquad B_p:=1+\alpha_{p,2}.
\]
For fixed source strengths $\alpha$, choose $q>1$ so close to $1$
that $q\alpha_{p,i}>-1$ for all $p\in S$ and $i=1,2$.
The local expression for $h_i$ then gives $h_i\in L^q(M)$ for
$i=1,2$.

Let $C^0(M)$ denote the space of continuous real-valued functions
on $M$, and set
\[
 X_0:=\left\{u=(u_1,u_2)\in C^0(M)^2:
           \int_M u_i\,dV_g=0,\quad i=1,2\right\},
 \qquad
 \|u\|_{X_0}:=\max_{i=1,2}\sup_{x\in M}|u_i(x)|.
\]
For $u\in X_0$, define the two densities by
\[
 f_i(u):=\frac{\rho_i h_i e^{u_i}}
                {\int_M h_i e^{u_i}\,dV_g},
 \qquad i=1,2.
\]
Then $f_i(u)\in L^q(M)$ and $\int_M f_i(u)\,dV_g=\rho_i$.
Define $\mathcal T_\rho:X_0\to X_0$ by
\[
 (\mathcal T_\rho u)_i
 :=\sum_{j=1}^2 K_{ij}\bigl(G*f_j(u)\bigr),
 \qquad i=1,2.
\]
The maps $f_i:X_0\to L^q(M)$ are continuous and bounded on bounded
sets. By the elliptic estimate for the Green operator and the compact
embedding of $W^{2,q}(M)$ into $C^0(M)$, $\mathcal T_\rho$ is compact.

Let $I$ be the identity on $X_0$ and set
$\mathcal F_\rho=I-\mathcal T_\rho$. Its zeros are precisely the
solutions of \eqref{eq:toda}. We write $\deg_{\mathrm{LS}}$ for the
Leray--Schauder degree and $\deg_{\mathrm B}$ for the Brouwer degree.

Fix the marked points, the source strengths and the smooth positive factors.
If the solution set is bounded in $X_0$, choose $R>0$ so that the
open ball $B_R^{X_0}(0)$ contains every zero of $\mathcal F_\rho$.
We define the degree of the singular Toda system by
\begin{equation*}
 d_\alpha(\rho_1,\rho_2)
 :=\deg_{\mathrm{LS}}\bigl(\mathcal F_\rho,B_R^{X_0}(0),0\bigr).
 \tag{D1}\label{def:singular-toda-degree}
\end{equation*}
By excision this value is independent of the choice of such a ball.

The degree counting program for the scalar mean field equation on a
closed surface goes back to Y.~Y.~Li \cite{Li1999}. For a smooth
positive function $h$ on $M$ and a parameter $\lambda>0$, we consider
the equation
\[
 -\Delta_g w
 =\lambda\left(
     \frac{he^w}{\int_M he^w\,dV_g}-1
   \right)
 \quad\text{in }M,\qquad
 \int_M w\,dV_g=0.
\]
Let $X:=\{w\in C^0(M):\int_M w\,dV_g=0\}$ with the norm $\|w\|_X:=\sup_{x\in M}|w(x)|$.  For $\lambda>0$ and a smooth positive function $H>0$ defined on $M$, set
\[
 \mathcal L_{\lambda,H}(w)
 :=w-\lambda G*\left(\frac{He^w}{\int_M He^w\,dV_g}-1\right),
 \qquad w\in X.
\]
For the smooth positive function $h$ above and each integer $m\geq0$,
let $d_m^{\mathrm{sc}}$ denote the common Leray--Schauder degree
of $\mathcal L_{\lambda,h}$ on a large ball in $X$ containing all its
zeros, for $8\pi m<\lambda<8\pi(m+1)$. Since $d_0^{\mathrm{sc}}=1$, Y. Y. Li
proposed to determine all the degrees by computing the jump at each
$8\pi m$ with $m\geq1$. Chen and Lin \cite{CL2002} obtained sharp
estimates for solutions blowing up with $m$ bubbles. In particular,
they derived an expansion of $\lambda-8\pi m$ and the locations of the bubbles in terms of the
heights. Using these estimates, they
computed the contributions of the bubbling solutions to the degree
jump and obtained in  \cite{CL2003}
\[
 d_m^{\mathrm{sc}}
 =\frac{1}{m!}\prod_{j=1}^{m}(j-\chi),
 \qquad m\geq1.
\]
Here, $\chi$ is the Euler characteristic number of $M$. For the scalar equation with singular sources, Chen and Lin
established sharp blow-up estimates in \cite{CL2010} and derived
the degree formula for mean field equations with sources of strength greater than $-1$ in
\cite{CL2015}.

Sharp estimates for blowup solutions also apply to the study of local uniqueness.
Bartolucci, Jevnikar, Lee and Yang \cite{BJLY2019} considered two
sequences of solutions with the same parameters and the same regular
blow-up points. Under their nondegeneracy and certain nonvanishing assumptions,
they used sharp estimates and local Pohozaev identities to show that
the two sequences coincide.
Bartolucci, Yang and Zhang \cite{BYZ2026U} obtained refined expansions
at non-quantized singular sources and proved local uniqueness for
the bubbling solutions that may contain both regular and singular blow-up
points. Their approach involves a higher order expansion based on \cite{CL2002,CL2010,Zhang2009}. The same assumptions also lead us to the nondegeneracy of bubbling solutions  \cite{BYZ2026N}.

Lee, Lin, Wei and Yang \cite{LLWY2018} initiated the study of degree jumps for the $SU(3)$ Toda system by introducing a shadow system. In partial blow-up, the shadow system describes the blow-up location and the distribution of the residual mass.
They proved that, as $\rho_1$ increases through $4\pi$ with
$\rho_2\notin4\pi\mathbb Z_{>0}$ fixed, the jump in the Toda degree
equals the negative of the degree of the corresponding shadow system. For simple singularities, Lee, Lin, Yang and Zhang \cite{LLYZ2020}
computed the degree in the first two strips,
\[
 \rho_1\in(0,4\pi)\cup(4\pi,8\pi),\qquad
 \rho_2\notin4\pi\mathbb Z_{>0}.
\]
For the non-symmetric systems of types $B_2$ and $G_2$, Lee, Lin
and Yan \cite{LLY2022} computed degree jumps involving several
blow-up points in the parameter ranges considered there.

As in the scalar case, the calculation for the $SU(3)$ Toda system
requires precise information on bubbling solutions. For fully bubbling
solutions, Lin, Wei and Zhao \cite{LWZhao2012} obtained sharp
convergence rates and derived first and second order conditions
at the blow-up points. For a positive integer $N$, Wu and Zhang
\cite{WZ2020} studied the case in which the parameters approach
$(4\pi,4\pi N)$ or $(4\pi N,4\pi)$. They determined the leading
terms of the parameter differences and obtained estimates for the
locations and relative heights of the bubbles. When both parameters
approach critical values, Wei, Wu and Zhang \cite{WWZ2023} proved
a spherical Harnack inequality near a fully bubbling point under a curvature assumption. For singular Toda systems, the possible
local mass pairs were studied in \cite{LWZ2015,LWWZ2018}.

In our previous work \cite{LLWYZ2026}, we computed the degree of the
regular $SU(3)$ Toda system on every noncritical parameter rectangle.
For integers $m,k\geq0$, let $d^{\mathrm{reg}}_{m,k}$ denote the
degree of the regular system on
\[
 (4\pi m,4\pi(m+1))\times(4\pi k,4\pi(k+1))
\]
for some non-negative integers $m$ and $k$. By computing the shadow degrees and the local contributions at
ordinary and mixed blow-up points, we obtained the following identity
in the formal variables $x$ and $y$:
\[
 \sum_{m,k\geq0}d^{\mathrm{reg}}_{m,k}x^m y^k
 =(1-x)^{\chi-1}(1-y)^{\chi-1}(1-xy)^\chi.
\]
Here, $\chi$ denotes the Euler characteristic number of the surface $M$. At each marked point $p\in S$, we assume that the source in the $i$th component has strength $\alpha_{p,i}>-1$ for $i=1,2$. For the parameter values specified in Theorem~\ref{thm:main}, we compute the degree of the singular Toda system by summing the coefficients of all terms $x^a y^b$ with $4\pi a<\rho_1$ and $4\pi b<\rho_2$ in a formal series whose exponents need not be integers.

\subsection{Main results}

We first specify the parameter sets on which the degree will be computed. At $p\in S$, set $(A,B)=(A_p,B_p)$ and express the local masses in units of $4\pi$. On the ellipse $s^2-st+t^2=As+Bt$ in $[0,\infty)^2$, start with the six pairs
\[
 (0,0),\ (A,0),\ (0,B),\ (A,A+B),\ (A+B,B),\ (A+B,A+B).
\]
For each pair $(u,v)$ already obtained, adjoin every pair $(s,t)$ on the ellipse satisfying
\[
 s\ge u,\qquad t\ge v,\qquad
 \{s-u,t-v\}\cap\mathbb Z_{\ge0}\ne\varnothing.
\]
Denote
\[
\begin{aligned}
 S_p^{(0)}
 &:=\{(0,0),(A,0),(0,B),(A,A+B),(A+B,B),(A+B,A+B)\},\\
 S_p^{(k+1)}
 &:=\bigcup_{(u,v)\in S_p^{(k)}}
 \left\{(s,t)\in[u,\infty)\times[v,\infty):
 \begin{array}{l}
 s^2-st+t^2=As+Bt,\\
 \{s-u,t-v\}\cap\mathbb Z_{\ge0}\ne\varnothing
 \end{array}
 \right\},
 \qquad k\ge0,\\
 \Sigma_p(\alpha)
 &:=\bigcup_{k=0}^{\infty}
     \left\{(2s,2t):(s,t)\in S_p^{(k)}\right\}.
\end{aligned}
\]
We write $\Sigma_p(\alpha)$ to emphasize its dependence on the two
source strengths $\alpha_{p,1}$ and $\alpha_{p,2}$ at $p$. The sets $S_p^{(k)}$ are nested by construction. For fixed $(A,B)$,
their union is finite by Step~1 of the proof of Lemma~\ref{fnd:coarse-local}.
Hence, there exists $N\ge0$ such that $S_p^{(k)}=S_p^{(N)}$ for all $k\ge N$. The possible masses at limiting source strengths are also contained in this finite set. Define
\begin{equation}\label{eq:coarse-mass-sets}
 \begin{split}
 E_{i,p}(\alpha)&=\{\sigma_i/2:
                    (\sigma_1,\sigma_2)\in\Sigma_p(\alpha)\},\\
 \mathcal C_i(\alpha)&=(0,\infty)\cap
  4\pi\left\{n+\sum_{p\in S}e_p:
          n\in\mathbb Z_{\geq0},\ e_p\in E_{i,p}(\alpha)\right\}.
 \end{split}
\end{equation}
If one component concentrates completely, its total mass belongs to \(\mathcal C_i(\alpha)\). This set is locally finite since it is a finite union of translates of \(4\pi\mathbb Z_{\geq0}\), restricted to positive values. The construction gives necessary conditions for blow-up. We shall also use the smaller sets
\begin{equation}\label{eq:sharp-mass-sets}
\begin{split}
\Gamma_1(\alpha)&=(0,\infty)\cap
  4\pi\left\{n+\sum_{p\in S}e_p:
       n\in\mathbb Z_{\geq0},\ e_p\in\{0,A_p,A_p+B_p\}\right\},\\
\Gamma_2(\alpha)&=(0,\infty)\cap
  4\pi\left\{n+\sum_{p\in S}e_p:
       n\in\mathbb Z_{\geq0},\ e_p\in\{0,B_p,A_p+B_p\}\right\}.
\end{split}
\end{equation}

For arbitrary source strengths $\alpha_{p,i}>-1$,
Proposition~\ref{fnd:signed-compactness}\textup{(i)} gives an a priori
$L^\infty$--bound for all solutions of \eqref{eq:toda},
uniformly on compact subsets of the parameter region
$\{(\alpha,\rho): \rho_i\notin\mathcal C_i(\alpha),\ i=1,2\}$. For fixed source strengths, the same bound holds with $\Gamma_i$ in place of
$\mathcal C_i$ for $i=1,2$, provided that either all source strengths are integers
or $1,(A_p)_{p\in S},(B_p)_{p\in S}$ are linearly independent over
$\mathbb Q$ (Lemma~\ref{fnd:sharp-compactness}).
Thus the degree in \eqref{def:singular-toda-degree} is well defined
under either parameter exclusion.

For \(A,B>0\), we define
\begin{align}
 F(x,y)&=(1-x)(1-y)(1-xy),\label{eq:F}\\
 P_{A,B}(x,y)
 &=1-x^A-y^B+x^Ay^{A+B}+x^{A+B}y^B
       -x^{A+B}y^{A+B}.\label{eq:P}
\end{align}
Especially, $F(x,y)=P_{1,1}(x,y)$. Moreover, we write
\begin{equation}\label{eq:real-jump-series}
 \mathscr J_\alpha(x,y)
 =F(x,y)^{\chi-\ell}\prod_{p\in S}P_{A_p,B_p}(x,y)
 =\sum_{a,b\geq0}c_{a,b}(\alpha)x^ay^b
\end{equation}
and
\begin{equation}\label{eq:real-cumulative}
 \mathscr D_\alpha(r,s)
   =\sum_{\substack{a<r\\b<s}}c_{a,b}(\alpha).
\end{equation}
With these notations, we state our main theorem.
\begin{theorem}[Degree formula]
\label{thm:main}
Let $\alpha_{p,i}>-1$ for all $p\in S$ and $i=1,2$, let
$\rho_1,\rho_2>0$, and let the weights in \eqref{eq:toda} be given
by \eqref{eq:physical-weights}.
Suppose that at least one of the following conditions holds.
\begin{itemize}
\item[\textup{(i)}]
$\rho_1\notin\mathcal C_1(\alpha)$ and
$\rho_2\notin\mathcal C_2(\alpha)$.
\item[\textup{(ii)}]
$\rho_1\notin\Gamma_1(\alpha)$ and $\rho_2\notin\Gamma_2(\alpha)$.
In addition, either all source strengths $\alpha_{p,i}$ are integers or the family
$1,(A_p)_{p\in S},(B_p)_{p\in S}$ is linearly independent over $\mathbb Q$.
\end{itemize}
Then the Leray--Schauder degree of \eqref{eq:toda} is given by
\begin{equation}\label{eq:real-degree}
 d_\alpha(\rho_1,\rho_2)
 =\mathscr D_\alpha\left(\frac{\rho_1}{4\pi},
\frac{\rho_2}{4\pi}\right).
\end{equation}
\end{theorem}

\begin{remark}
Now we have a comparison of the two parameter sets $\mathcal C_i(\alpha)$ and $\Gamma_i(\alpha)$.
\label{rem:mass-sets-strict}
In both \eqref{eq:coarse-mass-sets} and \eqref{eq:sharp-mass-sets},
the sum over $S$ contains one term $e_p$ for each marked point $p$,
and $e_p$ may be zero.
For $\Gamma_i(\alpha)$, this value is the $i$-th coordinate of one
of the six initial pairs, in units of $4\pi$.
In the definition of $\mathcal C_i(\alpha)$, one may also choose $e_p$
to be the $i$--th coordinate of a pair obtained at a later step of the recursion. Thus
\[
  \Gamma_i(\alpha)\subseteq\mathcal C_i(\alpha),
  \qquad i=1,2.
\]
The inclusion can be strict even when all source strengths are integers.

For example, take $S=\{p\}$ and $\alpha_{p,1}=\alpha_{p,2}=1$, so that
$A_p=B_p=2$. The second coordinates of the six initial pairs form
the set $\{0,2,4\}$. Hence
\[
  \Gamma_2(\alpha)
  =(0,\infty)\cap
    4\pi\{n+e:n\in\mathbb Z_{\geq0},\ e\in\{0,2,4\}\}
  =4\pi\mathbb Z_{>0}.
\]
Now start the recursion at $(0,0)$ and increase the first coordinate
by one. On the line $s=1$, the Pohozaev identity (see, for instance, \eqref{eq:first-strip-pohozaev}) gives
\[
  1-t+t^2=2+2t,
  \qquad
  t=\tau:=\frac{3+\sqrt{13}}{2},
\]
where $\tau$ is the positive root. Thus the first step adds $(1,\tau)$.
Taking its second coordinate in \eqref{eq:coarse-mass-sets} and
choosing $n=0$, we obtain
\[
  4\pi\tau\in\mathcal C_2(\alpha)\setminus\Gamma_2(\alpha),
  \qquad 3<\tau<4.
\]
By symmetry, the same value belongs to
$\mathcal C_1(\alpha)\setminus\Gamma_1(\alpha)$.
For these integer source strengths, local masses are integer multiples of $4\pi$
by \cite[Theorem~1.3]{LWWZ2018}. Thus $(1,\tau)$ is a candidate
produced by the recursion and is excluded by mass quantization.

The left panel shows the six initial pairs (filled) and the additional
pair $(1,\tau)$ (open). The right panel compares the five values
$1,2,3,\tau,4$ on the normalized parameter axis $r=\rho_2/(4\pi)$.

\begin{center}
\begin{tikzpicture}[font=\small,line cap=round]
  \pgfmathsetmacro{\massTau}{(3+sqrt(13))/2}

  \begin{scope}[x=0.83cm,y=0.83cm]
    \draw[->] (0,0) -- (4.7,0) node[right] {$s$};
    \draw[->] (0,0) -- (0,4.7) node[above] {$t$};

    \draw[gray,dashed]
      (1,0) -- (1,\massTau) -- (0,\massTau);
    \draw[semithick,domain=-120:120,samples=121,variable=\a]
      plot ({2+2*cos(\a)+2*sin(\a)/sqrt(3)},
            {2+2*cos(\a)-2*sin(\a)/sqrt(3)});

    \foreach \s/\t in {0/0,2/0,0/2,2/4,4/2,4/4}
      \fill (\s,\t) circle (1.7pt);
    \draw[fill=white,thick] (1,\massTau) circle (2.4pt);

    \node[below left] at (0,0) {$0$};
    \foreach \s in {1,2,4} \node[below] at (\s,0) {$\s$};
    \foreach \t in {2,4}   \node[left]  at (0,\t) {$\t$};
    \node[left] at (0,\massTau) {$\tau$};
  \end{scope}

  \begin{scope}[shift={(5.7cm,0.9cm)},x=1.02cm,y=1cm]
    \node[left] at (0.62,1.45) {$\Gamma_2(\alpha)/(4\pi)$};
    \node[left] at (0.62,0.75) {$\mathcal C_2(\alpha)/(4\pi)$};

    \draw[->] (0.70,0) -- (4.45,0) node[right] {$r$};
    \foreach \r in {1,2,3,4}{
      \draw[gray!55,dotted] (\r,0) -- (\r,1.45);
      \foreach \y in {0.75,1.45}
        \fill (\r,\y) circle (1.7pt);
      \node[below] at (\r,0) {$\r$};
    }
    \draw[gray,dashed] (\massTau,0) -- (\massTau,0.75);
    \draw[fill=white,thick] (\massTau,0.75) circle (2.4pt); 
    \node[below] at (\massTau,0) {$\tau$};
  \end{scope}
\end{tikzpicture}
\end{center}
\end{remark}

 Assume that $1,(A_p)_{p\in S},(B_p)_{p\in S}$ are linearly
independent over $\mathbb Q$ and that
$\rho_2/(4\pi)\in\mathbb Q\setminus\mathbb Z$. At a critical
value $\rho_1^0$, choose the smooth factors for the local degree
calculation. Let $Z$ range over the resulting nondegenerate shadow
solutions, with centers of each type unordered. For each $Z$, let
$U(Z)$ count the ordinary regular blow-up points, let $P(Z)$ be
the set of pure singular blow-up points, and let $d_S(Z)$ be the
local degree of the shadow system at $Z$, as defined in
\eqref{def:local-shadow-degree}. For sufficiently small $\varepsilon>0$,
\[
 d_\alpha(\rho_1^0+\varepsilon,\rho_2)
 -d_\alpha(\rho_1^0-\varepsilon,\rho_2)
 =\sum_Z(-1)^{U(Z)+|P(Z)|}d_S(Z).
\]

The smooth factors may be chosen for each critical value. The degree at either noncritical endpoint is independent of this choice.

At fixed source strengths, an a priori bound outside $\Gamma_i$ also holds if, at each $p$, both source strengths are integers or $1,A_p,B_p$ are rationally independent. The alternative may depend on $p$.

More generally, once \eqref{eq:real-degree} holds outside the sets $\mathcal C_i(\alpha)$, it extends to every positive parameter pair having a neighborhood on which all mean-zero solutions are uniformly bounded for the fixed weights and source strengths.

Indeed, approximate the pair from below in each coordinate by pairs outside the two locally finite sets.

Homotopy invariance gives the same degree on a sufficiently small rectangle. Since the support of the series is locally finite and the inequalities in \eqref{eq:real-cumulative} are strict, its partial sum is also the same for all sufficiently close pairs in this sequence.

Fix integer source strengths, so that $\Gamma_i(\alpha)=4\pi\mathbb Z_{>0}$
for $i=1,2$. For $k,l\in\mathbb Z_{\geq0}$, set
\[
 \mathcal C_{k,l}
 =(4\pi k,4\pi(k+1))\times(4\pi l,4\pi(l+1)).
\]
Lemma~\ref{fnd:sharp-compactness} gives a common bound along every
line segment in $\mathcal C_{k,l}$. By homotopy invariance,
$d_\alpha$ is constant there. We denote its common value by
\[
 d_{k,l}:=d_\alpha(\rho_1,\rho_2),
 \qquad (\rho_1,\rho_2)\in\mathcal C_{k,l}.
\]
For nonnegative integers \(a,b\), define
\begin{equation}\label{eq:Q}
 Q_{a,b}(x,y)=\sum_{i=0}^{a}\sum_{j=0}^{b}y^{i+j}
\sum_{r=i}^{a+j}x^r.
\end{equation}
For integer source strengths $P_{a+1,b+1}=FQ_{a,b}$. Set
\begin{align}\label{eq:degree-candidate}
 D_{\mathrm{form}}(x,y)
 &=\frac{F(x,y)^{\chi-\ell}}{(1-x)(1-y)}
       \prod_{p\in S}P_{A_p,B_p}(x,y)\notag\\
 &=(1-x)^{\chi-1}(1-y)^{\chi-1}(1-xy)^\chi
       \prod_{p\in S}Q_{\alpha_{p,1},\alpha_{p,2}}(x,y).
\end{align}
\begin{corollary}[Integer source strengths]
\label{thm:integer-degree}
Suppose $\alpha_{p,i}\in\mathbb Z_{\ge0}$ for all $p,i$.
\begin{itemize}
\item[\textup{(i)}]
The generating function is
\begin{equation}\label{eq:analytic-identification}
 \sum_{k,l\geq0}d_{k,l}x^ky^l=D_{\mathrm{form}}(x,y).
\end{equation}
For $S=\varnothing$ it equals $(1-x)^{\chi-1}(1-y)^{\chi-1}(1-xy)^\chi$.
\item[\textup{(ii)}]
If $\chi\le0$, then $d_{k,l}>0$ for all $k,l\ge0$. Thus \eqref{eq:toda} has a solution for every $\rho_1,\rho_2\in(0,\infty)\setminus4\pi\mathbb N$.
\end{itemize}
\end{corollary}
We next vary the source strengths. We fix the marked set $S$ and the positive smooth functions $h_i^*$. We allow both source strengths at a point of $S$ to vanish.
\begin{theorem}[Continuation to negative source strengths]
\label{thm:signed-continuation}
\leavevmode
\begin{itemize}
\item[\textup{(i)}]
Suppose
\[
 \alpha^0_{p,i}\ge0,\qquad \rho^0_i\in(0,\infty)\setminus\mathcal C_i(\alpha^0),\qquad i=1,2.
\]
There is a neighborhood $\mathcal U$ of $(\alpha^0,\rho^0)$ in $(-1,\infty)^{2\ell}\times(0,\infty)^2$. For every $(\alpha,\rho)\in\mathcal U$ the degree is defined and
\[
 d_\alpha(\rho_1,\rho_2)
 =d_{\alpha^0}(\rho^0_1,\rho^0_2)
 =\mathscr D_\alpha\left(\frac{\rho_1}{4\pi},
\frac{\rho_2}{4\pi}\right).
\]
\item[\textup{(ii)}]
Let $(\alpha(t),\rho(t))$, $0\le t\le1$, be a continuous path in $(-1,\infty)^{2\ell}\times(0,\infty)^2$. Suppose $\alpha_{p,i}(0)\ge0$ and
\[
 \rho_i(t)\notin\mathcal C_i(\alpha(t)),
 \qquad i=1,2,\quad 0\leq t\leq1.
\]
Then, for every $t\in[0,1]$,
\[
 d_{\alpha(t)}(\rho_1(t),\rho_2(t))
 =d_{\alpha(0)}(\rho_1(0),\rho_2(0))
 =\mathscr D_{\alpha(t)}\left(\frac{\rho_1(t)}{4\pi},
                             \frac{\rho_2(t)}{4\pi}\right).
\]
\end{itemize}
\end{theorem}
\begin{corollary}[Small perturbations of zero source strengths]
\label{cor:small-signed-orders}
Let $K$ be a compact subset of
\[
 \bigl((0,\infty)\setminus4\pi\mathbb N\bigr)^2.
\]
There exists $\varepsilon_K\in(0,1)$ such that, if $|\alpha_{p,i}|<\varepsilon_K$ for all $p\in S$ and $i=1,2$, the following hold.
\begin{itemize}
\item[\textup{(i)}]
For every $\rho\in K$ the degree is defined and
\[
 \begin{split}
 d_\alpha(\rho_1,\rho_2)
 &=\mathscr D_\alpha\left(\frac{\rho_1}{4\pi},
                          \frac{\rho_2}{4\pi}\right)\\
 &=[x^ky^l](1-x)^{\chi-1}(1-y)^{\chi-1}(1-xy)^\chi,
 \end{split}
\]
where $k=\lfloor\rho_1/(4\pi)\rfloor$ and $l=\lfloor\rho_2/(4\pi)\rfloor$.
\item[\textup{(ii)}]
If $\chi\le0$, then $d_\alpha(\rho_1,\rho_2)>0$ and \eqref{eq:toda} has a solution for every $\rho\in K$.
\end{itemize}
\end{corollary}
Compactness gives $\varepsilon_K$ without an explicit bound. It may depend on $K$ and the fixed data. For arbitrary source strengths greater than $-1$, put $a_i=\min(\{1\}\cup\{1+\alpha_{p,i}:p\in S\})$. For $(i,j)=(1,2)$ or $(2,1)$, the degree is given by \eqref{eq:real-degree} when
\[
 0<\rho_i<4\pi a_i,\qquad \rho_j\notin\mathcal C_j(\alpha).
\]
The proof deforms $\rho_i$ to zero and applies the scalar degree formula. This gives the initial degree for all source strengths greater than $-1$. For integer source strengths on a torus, $d_{k,l}$ is constant for large $k,l$. For noninteger source strengths, the degree can be negative in one rectangle and zero in another.

Fix $\rho_2\notin\Gamma_2(\alpha)$ satisfying the scalar exclusions for the configurations under consideration. We compute the change in degree as $\rho_1$ crosses a critical value.

At fixed generic source strengths greater than $-1$, each sequence blowing up at this value determines a solution of the shadow system, a scalar equation coupled to the balance equations for the regular blow-up points.

We choose positive smooth factors for which these solutions are nondegenerate.

Near each one we construct bubbling solutions and show that every blowing-up sequence belongs to one of the resulting families. The reduced equations determine their local degrees.

Adding these contributions to the initial scalar degree gives \eqref{eq:real-degree} at generic source strengths. Compactness under variation of the source strengths then gives the formula at arbitrary source strengths greater than $-1$.
\section{Proof of the degree formula}
\label{sec:assembly}
\subsection{Bubbling configurations and notation}
Throughout the local analysis, rational independence means
\begin{equation}\label{fnd:independence}
       1,(A_p)_{p\in S},(B_p)_{p\in S}
                \quad\hbox{are linearly independent over }\mathbb Q.
\end{equation}
\subsubsection*{Compactness and local masses}
Local masses in this subsection are in units of $4\pi$.
\begin{lemma}[Local mass congruence and the centered estimate]
\label{fnd:local-unitary}
On a fixed disk centered at $0$, suppose $F_i^k=e^{U_i^k}$ have uniformly bounded mass and satisfy
\[
 \Delta U_i^k+(KF^k)_i=b_i^k+4\pi\alpha_i^k\delta_0,
 \qquad \alpha_i^k\longrightarrow\alpha_i>-1,
\]
Assume the $b_i^k$ are smooth and uniformly bounded. For each $k$ the regular parts are continuous at $0$, with uniformly bounded oscillation on the boundary. Set
\[
 a=(A,B)=(1+\alpha_1,1+\alpha_2),\qquad
 4\pi x_i=\text{local mass of component }i\text{ at }0.
\]
Let $\lambda_v\in\mathbb R^3$ have sum zero and successive differences $v_1,v_2$. Write $x=(x_1,x_2)$ and $\iota(x)=(x_1,x_2-x_1,-x_2)$.
\begin{itemize}
\item[\textup{(i)}]
The local masses satisfy
\[
 x_1^2-x_1x_2+x_2^2=Ax_1+Bx_2,\qquad
 \lambda_a-\iota(x)=P\lambda_a+n
 \quad(P\text{ a permutation},\ n\in\mathbb Z^3,\ \sum_jn_j=0).
\]
\item[\textup{(ii)}]
If $1,A,B$ are rationally independent, the possible pairs are
\[
 (0,0),\ (A,0),\ (0,B),\ (A,A+B),\ (A+B,B),\ (A+B,A+B).
\]
For these fixed source strengths, if $0$ is the only possible blow-up point, there is a smaller disk $D$ about $0$ such that
\[
 \sup_{k,\,z\in D\setminus\{0\}} |z|^2 F_i^k(z)<\infty.
\]
\item[\textup{(iii)}]
At a regular blow-up point the nonzero pairs are
\[
 (1,0),\ (0,1),\ (1,2),\ (2,1),\ (2,2).
\]
At integer limiting source strengths, $x\in\mathbb Z^2$.
\end{itemize}
\end{lemma}
\begin{lemma}[The local masses at successive scales]
\label{fnd:coarse-local}
Under the assumptions of Lemma~\ref{fnd:local-unitary}, define $\mathscr S(A,B)$ on the nonnegative part of the Pohozaev ellipse. Start with
\[
 (0,0),\ (A,0),\ (0,B),\ (A,A+B),\ (A+B,B),\ (A+B,A+B).
\]
For each pair $(u,v)$ already obtained, include every pair $(s,t)$ on the ellipse satisfying
\[
 s\ge u,\quad t\ge v,\qquad
 s-u\in\mathbb Z_{\ge0}\ \text{or}\ t-v\in\mathbb Z_{\ge0}.
\]
Every limiting local mass pair at $0$ belongs to $\mathscr S(A,B)$.
\end{lemma}
\begin{proposition}[Uniform compactness as the source strengths vary]
\label{fnd:signed-compactness}
Fix $S$, including points where both source strengths vanish. Let $\mathscr K\subset C^3(M)^2$ be a compact family of pairs of positive smooth functions. Assume
\[
 \alpha^0_{p,i}>-1,\qquad \rho_i^0>0,\qquad \rho_i^0\notin\mathcal C_i(\alpha^0)\quad(i=1,2),
\]
where \eqref{eq:coarse-mass-sets} uses the full local mass lists with $\mu_i=1+\alpha^0_{p,i}$. There exist a neighborhood
\[
 U\subset (-1,\infty)^{2|S|}\times(0,\infty)^2
 \quad\text{of }(\alpha^0,\rho^0)
\]
and $C>0$ such that \textup{(i)} and \textup{(ii)} hold.
\begin{itemize}
\item[\textup{(i)}]
For $(\alpha,\rho)\in U$ and $(h_1^*,h_2^*)\in\mathscr K$, every mean-zero solution of \eqref{eq:toda} satisfies
\[
             \|u_1\|_{C^0(M)}+\|u_2\|_{C^0(M)}\leq C.
\]
\item[\textup{(ii)}]
On $(C^0(M)\cap\{\int_Mu=0\})^2$ the corresponding maps form a continuous compact family on bounded sets. Every continuous homotopy with these data is admissible on a common sufficiently large ball.
\item[\textup{(iii)}]
Along each continuous path on a compact interval the solutions are uniformly bounded and the homotopy is admissible on one ball, provided the positive smooth factors vary continuously in $C^3$ and
\[
 \alpha_{p,i}>-1,\qquad \rho_i>0,\qquad \rho_i\notin\mathcal C_i(\alpha)\quad(i=1,2).
\]
\end{itemize}
\end{proposition}
\begin{lemma}[Compactness at fixed source strengths]
\label{fnd:sharp-compactness}
Fix the points and source strengths. At each $p\in S$, assume that both source strengths are nonnegative integers or that $1,A_p,B_p$ are linearly independent over $\mathbb Q$. Consider a compact family of data with
\[
 \rho_i>0,\quad \rho_i\notin\Gamma_i(\alpha)\quad(i=1,2),
 \qquad h_i^*>0\text{ smooth, with the }C^3\text{ topology}.
\]
The mean-zero solutions have a common $C^0$ bound. The corresponding fixed-point maps form a continuous compact family on bounded sets.
\end{lemma}
Consider a sequence of bubbling solutions with $\rho_1$ approaching a critical value and $\rho_2\notin\Gamma_2(\alpha)$ fixed. Pass to a subsequence along which both measures have limits. By the argument of Lemma~\ref{fnd:sharp-compactness}, one component concentrates completely. For the second this would give
\[
 \rho_2\in\Gamma_2(\alpha).
\]
Hence the first component concentrates completely, the residual mass of the second is positive, and its logarithmic normalization constant has a finite limit. At a point with shifts $(A,B)$ and local masses $(x,y)$, Jensen's inequality on circles for the second component gives
\[
 2y-x\le B,
\]
as in the proof of Proposition~\ref{fnd:signed-compactness}. Among the six pairs in Lemma~\ref{fnd:local-unitary}, this leaves $(0,0),(A,0),(A+B,B)$. At a regular point the Pohozaev identity and integrality give the same list with
\[
 A=B=1.
\]
Thus we have $U$ ordinary and $V$ mixed regular points, together with disjoint sets $P,Q\subset S$ of pure and mixed singular points. Their local masses are
\[
\begin{array}{ccl}
\toprule
 \text{location}&\text{state}&(\text{first mass},\text{second mass})\\
\midrule
 M\setminus S&\text{ordinary}&(1,0)\\
 M\setminus S&\text{mixed}&(2,1)\\
 p\in S&\text{unused}&(0,0)\\
 p\in P&\text{pure singular}&(A_p,0)\\
 p\in Q&\text{mixed singular}&(A_p+B_p,B_p)\\
\bottomrule
\end{array}
\]
Write
\begin{align}
 \tau&=U+2V+\sum_{p\in P}A_p+\sum_{p\in Q}(A_p+B_p),\\
 \sigma_{\rm at}&=V+\sum_{p\in Q}B_p,\qquad
 \gamma=\rho_2-4\pi\sigma_{\rm at}>0,\label{eq:residual}\\
 n&=U+V+|P|+|Q|,\qquad v=V+|Q|.
\end{align}
The critical value is \(\rho_1^0=4\pi\tau\), and \(\mu=\rho_1-\rho_1^0\). The weight parameters for the residual scalar equation are
\[
\begin{array}{c|ccccc}
\text{state}&\text{unused }p&\text{pure }p&\text{mixed }p&
             \text{ordinary regular}&\text{mixed regular}\\
\hline
\text{weight parameter}&B_p-1&A_p+B_p-1&A_p-1&1&0 .
\end{array}
\]
Let $\nu_i$ be the sum of the limiting point masses in component $i$, now in the original units, and let $\eta_2$ be the residual measure of the second component. Away from the blow-up points and sources, the limit of $u_2$ is $-G*\nu_1+2G*\nu_2+w$, where $w=2G*\eta_2$. Set
\[
 H_{\mathbf z}=h_2\exp(-G*\nu_1+2G*\nu_2).
\]
Local convergence gives
\[
 \eta_2=e^{c_2}H_{\mathbf z}e^w\,dV_g,\qquad
 c_2=\lim\log\frac{\rho_2}{\int_M h_2 e^{u_2}},
\]
away from these points. Since $\eta_2$ gives zero mass to each point, integration yields $e^{c_2}=\gamma/\int_M H_{\mathbf z}e^w$. Thus the mean-zero potential $w$ solves
\[
 -\Delta_gw=2\gamma\left(\frac{H_{\mathbf z}e^w}
               {\int_MH_{\mathbf z}e^w}-1\right).
\]
Each entry in the last row is greater than $-1$. On disks of small $\eta_2$-mass, apply the Jensen and H\"older estimates in Proposition~\ref{fnd:signed-compactness} to $2G*\eta_2$. The density belongs to $L^p(M)$ for some $p>1$, so
\[
 w\in C^0(M).
\]
The equations for the moving centers follow from the translation Pohozaev identities.

For each configuration write $\alpha_p,\beta_p$ for the effective weight parameters in its two shadow factors $h_1,h_2$. For positive source strengths in the original system
\[
 \beta_p>0.
\]
The three possibilities for $(a_p,b_p)=(1+\alpha_p,1+\beta_p)$ are
\begin{equation}\label{fnd:effective}
 (A_p,B_p),\qquad(-A_p,A_p+B_p),\qquad(-(A_p+B_p),A_p).
\end{equation}
Indeed a source with Toda local mass pair $(s,t)$ in $4\pi$ units adds $4\pi(2s-t)G$ and $4\pi(-s+2t)G$ to the two logarithmic weights. The singular part of $G$ then changes the shifted parameters to $(A_p-2s+t,B_p+s-2t)$.

Substitution of $(0,0)$, $(A_p,0)$ and $(A_p+B_p,B_p)$ gives the
three pairs above. Denote the ordinary free centers by $p_1,\ldots,p_U$ and the mixed free centers by $q_1,\ldots,q_V$. Write $\mathbf z=(p_1,\ldots,p_U,q_1,\ldots,q_V)$ and put
\[
 f=\frac{\gamma h_2e^{-4\pi\sum_iG(\cdot,p_i)}e^w}
              {\int_Mh_2e^{-4\pi\sum_iG(\cdot,p_i)}e^w},
 \qquad -\Delta_gw=2(f-\gamma),\quad\int_Mw=0.
\]
The balance equations are
\begin{align*}
 F_i^A&=F_{i,0}^A-\tfrac12\nabla w(p_i),\\
 F_{i,0}^A&=\nabla(\log h_1+4\pi r_g)(p_i)
       +8\pi\sum_{l\ne i}\nabla G(p_i,p_l)
       +12\pi\sum_j\nabla G(p_i,q_j),\\
 F_j^B&=\nabla(2\log h_1+\log h_2+12\pi r_g)(q_j)
       +12\pi\sum_i\nabla G(q_j,p_i)
       +24\pi\sum_{l\ne j}\nabla G(q_j,q_l).
\end{align*}
We assume
\begin{equation}\label{fnd:scalar-exclusion}
 \gamma\notin4\pi\bigcup_{E\subset S}
                   \left(\mathbb Z+\sum_{p\in E}\beta_p\right).
\end{equation}
The finite configuration family is closed under pair replacement
\[
 (U,V,\gamma)\mapsto(U-2,V+1,\gamma-4\pi)
\]
whenever applicable with positive residual parameter. This operation preserves $U+2V$ and all source states. The increase of the total concentrated mass of the second component by $4\pi$ gives the displayed change in $\gamma$.

The exclusion in \eqref{fnd:scalar-exclusion} persists after pair
replacement because the union of cosets is invariant under this
integer shift.

Choose disjoint oriented coordinate neighborhoods of the free centers
in $M\setminus S$. Use these coordinates for the centers and for the
components of their balance equations. Write $f(w,\mathbf z)$ for
the density above. For $w\in X$, define
\[
 \widehat F_i^A(w,\mathbf z)
 :=F_{i,0}^A(\mathbf z)-
   \int_M\nabla_{p_i}G(p_i,x)f(w,\mathbf z)(x)\,dV_g(x).
\]
At a zero of $\mathcal L_{2\gamma,H_{\mathbf z}}$, the integral
equals $\nabla w(p_i)/2$, so $\widehat F_i^A=F_i^A$. In this
product chart the shadow map, with its scalar component first, is
\[
 \Phi(w,\mathbf z)=
 \left(\mathcal L_{2\gamma,H_{\mathbf z}}(w),
       (\widehat F_i^A(w,\mathbf z))_{i=1}^{U},
       (F_j^B(\mathbf z))_{j=1}^{V}\right).
\]
For an isolated zero $Z$ of $\Phi$, choose a bounded open
neighborhood $\mathcal O_Z$ whose closure lies in the chart and
contains no other zero. Define
\begin{equation*}
 d_S(Z):=\deg_{\mathrm{LS}}(\Phi,\mathcal O_Z,0).
 \tag{D2}\label{def:local-shadow-degree}
\end{equation*}
We also write $d_S$ when $Z$ is fixed. Permutations of centers of
the same type preserve orientation, since each center contributes
two real coordinates.

Let \(s_a\) and \(t_a\) denote the outer and inner profile exponents at a blow-up point, respectively, with
\[
 (s_a,t_a)=(1,0),\ (2,1),\ (A_p,0),\
                       (A_p+B_p,B_p)
\]
in the four nonempty cases. \(R_a\) and \(r_a\) denote the outer and inner scales, respectively. We write \(B_a\in\C\) for the shape parameter of a regular mixed bubble and \(B_p=1+\alpha_{p,2}\) for the shifted parameter at a source \(p\in S\).

The common concentration parameter \(T\) increases toward blow-up. For configurations containing only regular blow-up points, we use
\[
 t=e^{-T/2},
\]
so ordinary and inner mixed widths are comparable to \(t\) and outer mixed widths to \(t^{1/2}\).

In a conformal coordinate $z=x+iy$, let $\partial_z=(\partial_x-i\partial_y)/2$. Identifying $(v_1,v_2)$ with $v_1+iv_2$ gives
\[
 \nabla_{\rm eucl}f=f_x+if_y.
\]
In the equations for the centers and the shape parameters, we take the complex conjugate before passing to real coordinates. The quadratic harmonic term $\Re(Tz^2)$ then contributes $(\Re T,-\Im T)$ to the equation for the shape parameter.

At a regular mixed point $q$ take $z(q)=0$ and $\varphi(q)=0$. The regular part of the first limiting field is
\[
 \psi_1=\log h_1+\varphi+2G*\nu_1-G*\nu_2-w/2+6\log|z|.
\]
Here $h_1$ is the original physical weight. Put
\[
 \mathcal T(q)=(\partial_z^2\psi_1+(\partial_z\psi_1)^2)(q),\qquad
 \kappa(q)=\Delta_g\log h_1(q)+2\rho_1^0-\rho_2-2K_g(q).
\]
We shall impose the following strict inequality at every mixed center:
\begin{equation}\label{fnd:gap}
                  \kappa>4e^{-\varphi}|\mathcal T|.
\end{equation}
\begin{theorem}[Choice of weights satisfying the trace inequalities]
\label{fnd:common-weights}
Fix $A_p,B_p>1$ satisfying \eqref{fnd:independence}. Let $\mathscr F$ be a finite family of configurations closed under pair replacement, each satisfying \eqref{fnd:scalar-exclusion}. There are common positive smooth factors $h_1^*,h_2^*$ such that
\begin{itemize}
\item[\textup{(i)}]
each shadow system in $\mathscr F$ has a finite, compact set of nondegenerate solutions;
\item[\textup{(ii)}]
the linearized scalar operator at fixed centers is invertible at every solution;
\item[\textup{(iii)}]
the inequality \eqref{fnd:gap} holds at every mixed center.
\end{itemize}
All bounds and strict inequalities are uniform over the finite family.
\end{theorem}
For an ordinary configuration put
\[
 \psi_i^A=\log h_1+\varphi+8\pi R(\cdot,p_i)
       +8\pi\sum_{\ell\ne i}G(\cdot,p_\ell)-w/2,\quad
 d_i=\tfrac12 e^{\psi_i^A(p_i)/2},\quad \kappa_i=\Delta_g\psi_i^A(p_i)
\]
in coordinates with $\varphi(p_i)=0$.
\begin{corollary}[Additional conditions on the weights]
\label{fnd:extra-traces}
The common weights may also be chosen to satisfy the following conditions simultaneously.
\begin{itemize}
\item[\textup{(i)}]
At the prescribed points one may impose any finite set of inequalities of the form
\[
 \Delta_g\log h_1^*(p)-2K_g(p)-4\pi\sum_{q\in S}\alpha_{q,1}
+2\rho_1^{\rm wall}-\rho_2>0.
\]
\item[\textup{(ii)}]
For each solution with $U\ge1$ ordinary bubbles and no other bubble, the coefficient of the ordinary mass satisfies
\[
 L_0=2\pi\sum_i d_i^2\kappa_i\ne0.
\]
\end{itemize}
\end{corollary}
For a fixed configuration $(U,V,P,Q)$, write
$\boldsymbol\beta=(\beta_p)_{p\in S}$ for the effective scalar weight parameters.
Fix the free centers. Whenever the scalar solution set is bounded,
choose a ball $B_R^X(0)$ containing it and set
\[
 d_{\rm sc}(\gamma;\boldsymbol\beta,U)
 :=\deg_{\mathrm{LS}}\bigl(\mathcal L_{2\gamma,H_{\mathbf z}},
                           B_R^X(0),0\bigr).
\]
The scalar weight has an additional quadratic zero at each of
the $U$ ordinary centers. In the cases considered below, compactness
along paths of distinct centers in $M\setminus S$ and homotopy
invariance give the same degree for every choice of their positions.

For the chosen weights, if all shadow solutions of this configuration
are isolated and form a compact set, denote their total degree by
\[
 d^{\rm sh}(U,V,P,Q):=\sum_Z d_S(Z),
\]
where the sum counts centers of each type up to permutation.

\begin{proposition}[Degree of the shadow system]
\label{fnd:shadow-degree}
For $A_p,B_p>1$, take the common weights of
Theorem~\ref{fnd:common-weights}. Then
\[
 d^{\rm sh}(U,V,P,Q)
 =\binom{\chi(M)-|S|}{U+V}\binom{U+V}{U}
    d_{\rm sc}(\gamma;\boldsymbol\beta,U).
\]
\end{proposition}
\begin{proposition}[Shadow systems with integrable singular weights]
\label{univ:shadow}
Assume rational independence and
\[
 A_p,B_p>0,\qquad \beta_p>-1 .
\]
Use the effective pairs in \eqref{fnd:effective}. Suppose each configuration in a finite family closed under pair replacement satisfies \eqref{fnd:scalar-exclusion}. Then the following hold.
\begin{itemize}
\item[\textup{(i)}]
The centers stay away from the prescribed sources. One can choose common weights for which the shadow solutions are nondegenerate and the strict trace inequalities hold at every regular mixed point.
\item[\textup{(ii)}]
The homotopy to the decoupled system and the degree formula for the shadow system remain valid.
\item[\textup{(iii)}]
For fixed $p>1$ sufficiently close to one, the scalar solutions belong to $W^{2,p}(M)$ and are smooth away from the prescribed sources.
\end{itemize}
\end{proposition}
\begin{lemma}[The shadow data of a concentrating sequence]
\label{reg:shadow-limit}
Let the prescribed source strengths exceed $-1$. Assume that the following hold.
\begin{itemize}
\item[\textup{(i)}]
The first component concentrates completely. The residual mass of the second component is $\gamma>0$.
\item[\textup{(ii)}]
Each nonzero local mass pair is one of the four pairs in the configuration table.
\end{itemize}
Then the mean-zero scalar remainder and the separated limiting centers satisfy the scalar and center equations above.
\end{lemma}
\begin{proposition}[Continuation in the source strengths]
\label{fnd:continuation}
Fix nonnegative source strengths, positive smooth factors and positive parameters $\rho_i\notin\mathcal C_i(\alpha)$.
\begin{itemize}
\item[\textup{(i)}]
With the factors fixed, the mean-zero solutions are uniformly bounded in a neighborhood of $(\alpha,\rho)$ with $\alpha_{p,i}>-1$. Within this neighborhood one may preserve the degree and choose
\[
 \widehat\rho_i/(4\pi)\in\mathbb Q\setminus\mathbb Z,\qquad
 \widehat\alpha_{p,i}>0\quad\text{satisfying \eqref{fnd:independence}}.
\]
At these data the scalar exclusions in \eqref{fnd:scalar-exclusion} hold. The degree is constant under deformations of the positive smooth factors in a compact family.
\item[\textup{(ii)}]
If the original source strengths satisfy \eqref{fnd:independence} or are all integers, it suffices to assume $\rho_i\notin\Gamma_i(\alpha)$. First choose the parameters in the same component of
\[
 \bigl((0,\infty)\setminus\Gamma_1(\alpha)\bigr)
 \times
 \bigl((0,\infty)\setminus\Gamma_2(\alpha)\bigr).
\]
\item[\textup{(iii)}]
In both cases the data can be chosen to preserve the finite sum of coefficients of $F^{\chi-|S|}\prod_{p\in S}P_{A_p,B_p}$ over exponents $(a,b)$ with
\[
 a<\rho_1/(4\pi),\qquad b<\rho_2/(4\pi).
\]
\end{itemize}
\end{proposition}
\subsection{Reduction to generic source strengths}
For integer source strengths, the linearized singular Liouville equation has additional
kernel functions \cite{DEM2012}. We compute the local degrees at fixed
source strengths $\alpha_{p,i}$ satisfying \eqref{fnd:independence} and use compactness to return
to the prescribed data. The inverse estimates may depend on these fixed source strengths.

\smallskip
\noindent\textit{Step 1. Choice of source strengths and parameters.}
Fix the smooth factors and suppose first that $\alpha_{p,i}\ge0$.
\begin{itemize}
\item If $\rho_i\notin\mathcal C_i(\alpha)$ for $i=1,2$,
Proposition~\ref{fnd:signed-compactness} gives a common $C^0$ bound
in a neighborhood of $(\alpha,\rho)$.
\item If the source strengths are all integers or satisfy \eqref{fnd:independence},
it suffices that $\rho_i\notin\Gamma_i(\alpha)$ for $i=1,2$.
Keep $\alpha$ fixed and first move $\rho$ in the same parameter component
to $\widetilde\rho$ with $\widetilde\rho_i\notin\mathcal C_i(\alpha)$
for $i=1,2$, as in Proposition~\ref{fnd:continuation}\textup{(ii)}.
\end{itemize}
In this neighborhood choose $\widehat\alpha$ satisfying
\eqref{fnd:independence} and $\widehat\rho$ such that
\[
 \frac{\widehat\rho_i}{4\pi}\in\mathbb Q\setminus\mathbb Z,\qquad
 \widehat\alpha_{p,i}>\alpha_{p,i}\ge0.
\]
By Proposition~\ref{fnd:continuation}\textup{(i)--(iii)},
\[
\begin{aligned}
 d_\alpha(\rho_1,\rho_2)
 &=d_{\widehat\alpha}(\widehat\rho_1,\widehat\rho_2),\\
 \mathscr D_\alpha\left(\frac{\rho_1}{4\pi},\frac{\rho_2}{4\pi}\right)
 &=\mathscr D_{\widehat\alpha}\left(\frac{\widehat\rho_1}{4\pi},
                                    \frac{\widehat\rho_2}{4\pi}\right).
\end{aligned}
\]
Each relevant exponent stays on the same side of its cut-off, and
coefficients are added when exponents coincide. At the chosen data,
rational independence gives \eqref{fnd:scalar-exclusion} for the finite
family of configurations below the cut-offs with $\gamma>0$, including
all pair replacements with positive remainder.

\smallskip
\noindent\textit{Step 2. Choice of the smooth factors.}
At the fixed positive generic source strengths, Theorem~\ref{fnd:common-weights}
gives one pair of positive smooth factors for which all these shadow
solutions are nondegenerate and \eqref{fnd:gap} holds at every regular
mixed point. Corollary~\ref{fnd:extra-traces} gives the trace inequalities
at the sources and $L_0\ne0$ for purely ordinary configurations.
Only the auxiliary factors need satisfy these conditions. At noncritical
parameters Lemma~\ref{fnd:sharp-compactness} gives a common bound along
the homotopy back to the factors prescribed in Theorem~\ref{thm:main}.

\smallskip
\noindent\textit{Source strengths greater than $-1$.}
\begin{itemize}
\item Near a nonnegative base point satisfying the coarse exclusions,
or along a path from it with $\alpha_{p,i}(t)>-1$ and
$\rho_i(t)\notin\mathcal C_i(\alpha(t))$,
Proposition~\ref{fnd:signed-compactness} gives an admissible Toda
homotopy even at scalar resonance. No uniform scalar inverse is needed
there. The sum \eqref{eq:signed-finite-sum} is constant by the same
cut-off argument (Theorem~\ref{thm:signed-continuation}).
\item For arbitrary $\alpha_{p,i}>-1$ with
$\rho_i\notin\mathcal C_i(\alpha)$ for $i=1,2$, the local degree
calculation uses nearby generic source strengths greater than $-1$
(Section~\ref{sec:wall-completion}, Step~4 of the proof of
Theorem~\ref{thm:main}). Proposition~\ref{univ:shadow} gives the choice
of factors for integrable singular weights.
Lemma~\ref{universal:wallwise-weights}\textup{(ii)} allows the factors
in Proposition~\ref{univ:wall-jump} to be chosen separately at each
critical value.
\end{itemize}

\smallskip
\noindent\textit{Blow-up estimates used in the reduction.}
\begin{itemize}
\item The local mass classification used in
Proposition~\ref{fnd:signed-compactness} follows from
Lemma~\ref{fnd:coarse-local}. For local masses and compactness, see
also \cite{LWZ2015,LWWZ2018}. For entire Toda solutions with one source,
see \cite{LWY2012}.
\item For pointwise estimates for singular Liouville equations and
for collapsing integer singularities, see \cite{BT2017} and
\cite{Tarantello2023}, respectively. For a cluster with one fixed
noninteger source, Proposition~\ref{fnd:scalar-cosets} gives
\eqref{fnd:two-cosets} by a unitary holonomy argument.
\item For regular and non-simple blow-up, see \cite{WZ2020,WWZ2023}.
The local selection and profile estimates used in the comparison are
given in Lemmas~\ref{reg:elementary-selection} and~\ref{hyb:singular-rate}.
\end{itemize}

Write \(\dot H^1(M)=\{v\in H^1(M):\int_M v\,dV_g=0\}\) and set \(\Hh=(\dot H^1(M))^2\).
We use the inner product and the functional
\[
 \langle u,v\rangle_{\Hh}
 =\sum_{i,j=1}^2(K^{-1})_{ij}\int_M\nabla u_i\cdot\nabla v_j,\qquad
 J_\rho(u)=\frac12\|u\|_{\Hh}^2-
 \sum_{i=1}^2\rho_i\log\int_Mh_ie^{u_i}.
\]

We compute the degree at small $\rho_1$ and add the local degrees as $\rho_1$ increases.
\begin{proposition}[Local degree contribution]
\label{def:complete-transfer}
Choose the generic source strengths and weights as in Theorem~\ref{fnd:common-weights} and Corollary~\ref{fnd:extra-traces}. Fix $\rho_1^0$ critical and $\rho_2\notin\Gamma_2(\alpha)$. Consider all configurations with
\[
 4\pi\tau=\rho_1^0,\qquad \gamma=\rho_2-4\pi\sigma_{\rm at}>0.
\]
Assume \eqref{fnd:scalar-exclusion} for each $\gamma$. Pair replacement with positive remainder preserves the family. Let $Z$ solve the shadow system for $(U,V,P,Q)$, with unordered regular centers, and assume $Z$ is nondegenerate.
\begin{itemize}
\item[\textup{(i)}]
For each sufficiently small $\mu\ne0$ there is a bounded open set $\mathcal U_{Z,\mu}$, possibly empty, with zero-free boundary.
\item[\textup{(ii)}]
Any sequence of solutions $u_j$ with $0\ne\mu_j\to0$ and $\|u_j\|_{\Hh}\to\infty$ has a subsequence with limiting shadow solution $Z_*$ such that $u_j\in\mathcal U_{Z_*,\mu_j}$ for all sufficiently large $j$.
\item[\textup{(iii)}]
The sets for distinct shadow solutions are disjoint. The solutions at $\rho_1=\rho_1^0$ form a bounded set.
\item[\textup{(iv)}]
For each sufficiently small $\mu\ne0$, define the local Toda degree by
\begin{equation*}
 d_{Z,\mu}:=
 \deg_{\mathrm{LS}}\bigl(\mathcal F_{(\rho_1^0+\mu,\rho_2)},
\mathcal U_{Z,\mu},0\bigr). \tag{D3}\label{def:local-toda-degree}
\end{equation*}
We set $d_{Z,\mu}=0$ if $\mathcal U_{Z,\mu}$ is empty. With the
orientation of $I-\mathcal T_\rho$, the contribution of $Z$ to the
jump, for sufficiently small $\epsilon>0$, is
\begin{equation}\label{eq:local-transfer}
 \Delta d_Z:=d_{Z,\epsilon}-d_{Z,-\epsilon}
            =(-1)^{U+|P|}d_S(Z).
\end{equation}
\end{itemize}
The constants are uniform over the finite family of configurations and shadow solutions. By Lemma~\ref{lem:local-degree-spaces}, the local degrees in $\Hh$ agree with those in $C^0$. For zeros in these neighborhoods,
\[
 \mu_j\to0,\quad u_j\in\mathcal U_{Z,\mu_j}
 \quad\Longrightarrow\quad \|u_j\|_{\Hh}\to\infty .
\]
\end{proposition}
\begin{proposition}[The global degree jump]
\label{prop:global-localization}
Under the assumptions of Proposition~\ref{def:complete-transfer}, take $\epsilon>0$ sufficiently small with $\rho_1^0\pm\epsilon\notin\Gamma_1(\alpha)$. Then
\[
 d_\alpha(\rho_1^0+\epsilon,\rho_2)-d_\alpha(\rho_1^0-\epsilon,\rho_2)
 =\sum_{(U,V,P,Q)}\ \sum_Z (-1)^{U+|P|}d_S(Z).
\]
The sums run over the configurations and shadow solutions at $\rho_1^0$, with local contributions given by \eqref{eq:local-transfer}.
\end{proposition}
\begin{proof}
By Lemma~\ref{fnd:sharp-compactness}, the degree is defined at both endpoints since $\Gamma_1(\alpha)$ is locally finite and $\rho_2\notin\Gamma_2(\alpha)$. Part~(iii) of the preceding proposition bounds the zero set at
\[
 \mu=0.
\]
This set is compact since the normalized map is compact on bounded energy sets. For fixed weights, the exponential estimates give joint compactness there for nearby parameters. Choose an open ball $\mathcal B$ containing the compact zero set. If its boundary contained zeros for a sequence
\[
 \mu_j\to0,
\]
joint compactness would give a boundary zero at $\mu=0$.

For small $\epsilon$, all zeros in the finitely many local domains lie outside $\overline{\mathcal B}$. Set
\[
 \mathcal V_{Z,\pm\epsilon}=\mathcal U_{Z,\pm\epsilon}\setminus\overline{\mathcal B},
 \qquad \partial\mathcal V_{Z,\pm\epsilon}
 \subset \partial\mathcal U_{Z,\pm\epsilon}\cup\partial\mathcal B .
\]
Both boundaries on the right are zero-free, and excision preserves the local degrees. The resulting sets are mutually disjoint and contain every zero outside \(\mathcal B\) at \(\mu=\pm\epsilon\). Otherwise, along a sequence
\[
 \epsilon_j\to0,
\]
a bounded subsequence of uncovered zeros would converge to a zero at the critical parameter outside \(\mathcal B\). Along a subsequence with norm tending to infinity, item (ii) would place the zeros in one of the neighborhoods and hence outside \(\overline{\mathcal B}\) in the corresponding \(\mathcal V\).

Both alternatives are impossible. Let \(\mathcal O_\pm\) be bounded open sets containing \(\mathcal B\), the \(\mathcal V_{Z,\pm\epsilon}\), and the entire endpoint zero sets, and abbreviate \(\mathcal F_{\rho_1^0\pm\epsilon,\rho_2}\) by \(\mathcal F_{\pm\epsilon}\). Additivity gives
\[
 \deg_{\mathrm{LS}}(\mathcal F_{\pm\epsilon},\mathcal O_\pm,0)
 =\deg_{\mathrm{LS}}(\mathcal F_{\pm\epsilon},\mathcal B,0)
+\sum_Z\deg_{\mathrm{LS}}(\mathcal F_{\pm\epsilon},\mathcal V_{Z,\pm\epsilon},0).
\]
Homotopy invariance gives the same degree on $\mathcal B$ at the two endpoints. Subtract the identities and apply excision to replace $\mathcal V_{Z,\pm\epsilon}$ by $\mathcal U_{Z,\pm\epsilon}$. The formula follows from \eqref{eq:local-transfer}.
\end{proof}
\begin{lemma}[The first subcritical strip]
\label{lem:base-strip}
Suppose
\[
 \rho_2\in(0,\infty)\setminus\mathcal C_2(\alpha),\qquad
 0<\rho_1<4\pi\min(\{1\}\cup\{A_p:p\in S\}).
\]
Then $d_\alpha(\rho_1,\rho_2)$ equals the degree of the scalar equation for the second component with parameter $2\rho_2$ and weight $h_2$.
\end{lemma}
\begin{proof}
Along $(t\rho_1,\rho_2)$, $0\le t\le1$, the solutions have a common bound. This follows from local compactness and the classification of local masses. Indeed, a local mass pair $(x,y)$ in units of $4\pi$ satisfies
\[
 0\le x<A.
\]
The exclusion on $\rho_2$ gives positive residual mass in the second component, hence $2y-x\le B$. The Pohozaev identity gives
\[
 x(x-A)+y(y-x-B)=0,\qquad y-x-B\le-\tfrac12(x+B).
\]
Both summands are nonpositive, so $x=y=0$. The same argument holds as $t\to0$, with $x=0$. At $t=0$ use
\[
 (a,b)=(u_1+u_2/2,u_2).
\]
In these coordinates the map $I-\mathcal T_\rho$ is the identity in $a$ and the scalar map for
\[
 -\Delta_g b=2\rho_2(h_2e^b/\int_Mh_2e^b-1)
\]
in $b$. Compactness and this change of coordinates give the equality of degrees independently of the numerical scalar degree and of the recurrence across critical values.
\end{proof}
\begin{proposition}[The degree of the scalar equation]
\label{fnd:scalar-degree-internal}
Let $h$ have a positive smooth factor and fixed shifted parameters $a_p>0$ at the points of a finite set $S$. Suppose
\[
 \gamma>0,\qquad \gamma\notin
 4\pi\left\{n+\sum_{p\in P}a_p:n\in\mathbb Z_{\ge0},\ P\subset S\right\}.
\]
The degree of $-\Delta_gw=2\gamma(h e^w/\int h e^w-1)$ is
\[
 \sum_{b<\gamma/(4\pi)}[y^b]\,(1-y)^{\chi(M)-|S|}\prod_{p\in S}(1-y^{a_p}).
\]
\end{proposition}
\begin{theorem}[Degree for generic positive source strengths]
\label{thm:degree-assembly}
Suppose $\alpha_{p,i}>0$ and \eqref{fnd:independence} holds. For $\rho_i\in(0,\infty)\setminus\Gamma_i(\alpha)$, $i=1,2$, we have
\[
 d_\alpha(\rho_1,\rho_2)
 =\sum_{\substack{a<\rho_1/(4\pi)\\ b<\rho_2/(4\pi)}}
             [x^ay^b]\mathscr J_\alpha(x,y),
 \qquad
 \mathscr J_\alpha(x,y)=F(x,y)^{\chi-\ell}
                      \prod_{p\in S}P_{A_p,B_p}(x,y).
\]
\end{theorem}
\begin{proof}
\noindent\textit{Step 1. Choice of parameters and weights.} Set $r=\rho_1/(4\pi)$ and $s=\rho_2/(4\pi)$. Lemma~\ref{fnd:sharp-compactness} applies at the fixed generic source strengths and implies that the degree is constant on each component of
\[
 \bigl((0,\infty)\setminus\Gamma_1(\alpha)\bigr)
       \times\bigl((0,\infty)\setminus\Gamma_2(\alpha)\bigr).
\]
The nonzero coordinates of the support of $\mathscr J_\alpha$ belong to $\Gamma_i(\alpha)/(4\pi)$ in component $i$. The sum in the statement is therefore constant on the same components. Local finiteness of $\mathcal C_i(\alpha)$ allows the choice
\[
 r,s\in\mathbb Q\setminus\mathbb Z,\qquad \rho_i\notin\mathcal C_i(\alpha)\quad (i=1,2).
\]
If \eqref{fnd:scalar-exclusion} failed for a retained configuration, we would have
\[
 s-m_2=n+\sum_{p\in E}\beta_p,\qquad
 m_2=V+\sum_{p\in Q}B_p,\qquad n\in\mathbb Z,\quad E\subseteq S.
\]
Substitution of the three possible values of $\beta_p$ expresses $s$ as an integer plus an integral linear combination of the $A_p,B_p$. By \eqref{fnd:independence} and $s\in\mathbb Q$, all coefficients of the shifts vanish. This would give
\[
 s\in\mathbb Z.
\]
Thus \eqref{fnd:scalar-exclusion} holds for all retained configurations. Only finitely many have $\tau<r$ and $m_2<s$. Replacing two ordinary points by one mixed free point gives
\[
 (U,V,m_2,\tau)\longmapsto(U-2,V+1,m_2+1,\tau).
\]
A positive remainder keeps the new configuration in the family. By Theorem~\ref{fnd:common-weights} and Corollary~\ref{fnd:extra-traces}, choose common smooth weights with invertible scalar linearizations at fixed centers, nondegenerate shadow solutions, and the required trace and radial conditions.

Lemma~\ref{fnd:sharp-compactness} preserves the Toda degree along a positive smooth-weight homotopy to these weights, at the fixed generic source strengths and noncritical endpoint parameters.

\noindent\textit{Step 2. Local masses and degree jumps.} Increase $\rho_1$ from a positive subcritical value below the target. If a sequence blows up along this segment, one component concentrates completely.

The local mass classification and $\rho_2\notin\Gamma_2(\alpha)$ exclude complete concentration of the second component. The positive residual mass gives
\[
 2y-x\le B
\]
for the local masses $(x,y)$ in units of $4\pi$. The three excluded pairs satisfy
\[
\begin{array}{c|ccc}
(x,y)&(0,B)&(A,A+B)&(A+B,A+B)\\
2y-x&2B&A+2B&A+B
\end{array}
\]
with $2y-x>B$ in each column. At a regular point take
\[
 A=B=1.
\]
Only the four listed types remain, and their first masses determine the critical value. Lemma~\ref{reg:shadow-limit} supplies the limiting shadow data.

At a critical value with no such configuration this argument rules out an unbounded sequence and gives zero jump by compactness. For the remaining critical values, the local contributions are given by Proposition~\ref{def:complete-transfer}.

Apply Proposition~\ref{prop:global-localization} and add the jumps below the target to the scalar degree in the first strip.

\noindent\textit{Step 3. Sum over the configurations.} Write $c=\chi-\ell$. For a configuration $(U,V,P,Q)$, let
\[
 \tau=U+2V+\sum_{p\in P}A_p+\sum_{p\in Q}(A_p+B_p),
 \qquad m_2=V+\sum_{p\in Q}B_p.
\]
Proposition~\ref{fnd:shadow-degree}, applied at the selected common weights and at scalar parameter $\gamma=4\pi(s-m_2)>0$, gives the joint degree as the fixed-center scalar degree multiplied by the decoupled center degree
\[
 \binom{c}{U+V}\binom{U+V}{U}.
\]
The scalar weight has shifted parameters \(1+\beta_p\) at the distinct prescribed points and the shifted parameter \(2\) at each of the \(U\) ordinary centers, disjoint from the prescribed points.

In the scalar critical set, the shifts two contribute only integers in units of \(4\pi\). The exclusion \eqref{fnd:scalar-exclusion} therefore excludes that set at
\[
 \gamma=4\pi(s-m_2).
\]
Proposition~\ref{fnd:scalar-degree-internal} gives the scalar jump series at this fixed, positive, noncritical parameter
\[
 S_{U,P,Q}(y)
   =(1-y)^{c-U}(1-y^2)^U
        \prod_{p\in S}(1-y^{1+\beta_p}),
\]
where the shifted weight parameters are
\[
 1+\beta_p=
 \begin{cases}
 B_p,&p\notin P\cup Q,\\
 A_p+B_p,&p\in P,\\
 A_p,&p\in Q.
 \end{cases}
\]
Its coefficients with exponent strictly less than $s-m_2$ sum to the scalar degree. Thus the contribution at the first critical mass $\tau$ is represented by
\[
 (-1)^{U+|P|}
 \binom{c}{U+V}\binom{U+V}{U}
       x^\tau y^{m_2}S_{U,P,Q}(y),
\]
in the sense that the coefficients with second exponent less than $s$ sum to the degree jump there. Configurations with $m_2\geq s$ contribute nothing to this cumulative sum. The configuration
\[
 U=V=0,\qquad P=Q=\varnothing
\]
has $\tau=m_2=0$. We assign its scalar jump series to the formal term $x^0$. By Lemma~\ref{lem:base-strip}, its cumulative sum is precisely the initial scalar degree. It represents the initial value in the summation over first critical masses. Since
\[
 (1-y)^{c-U}(1-y^2)^U=(1-y)^c(1+y)^U,
\]
an ordinary regular point contributes $-x(1+y)$, while a mixed regular point contributes $x^2y$. Summing over the regular points gives
\begin{align*}
 \sum_{U,V\geq0}\binom{c}{U+V}\binom{U+V}{U}
       [-x(1+y)]^U(x^2y)^V
 &=\{1-x(1+y)+x^2y\}^{c}\\
 &=[(1-x)(1-xy)]^c.
\end{align*}
At a prescribed source $p$, the three possibilities give
\[
 \begin{array}{c|c}
 \text{no blow-up}&1-y^{B_p}\\
 \text{pure singular blow-up}&-x^{A_p}(1-y^{A_p+B_p})\\
 \text{mixed singular blow-up}&x^{A_p+B_p}y^{B_p}(1-y^{A_p}).
 \end{array}
\]
Their sum is $P_{A_p,B_p}$. Multiplying the contributions and including the initial scalar degree at $x^0$ gives
\[
 (1-y)^c[(1-x)(1-xy)]^c\prod_{p\in S}P_{A_p,B_p}
 =\mathscr J_\alpha(x,y).
\]
All operations are valid on every bounded rectangle of exponents: there are only finitely many source choices and finitely many regular configurations with first mass below its upper bound, and every scalar series has locally finite support.

Taking the cumulative sum over $a<r$ and $b<s$ gives the degree asserted in the theorem.
\end{proof}
\begin{proof}[The case of nonnegative source strengths in Theorem~\ref{thm:main}]
\noindent\textit{Step 1. Approximation for (i).} Fix nonnegative source strengths $\alpha$ and positive parameters satisfying
\[
 \rho_i\notin\mathcal C_i(\alpha).
\]
Proposition~\ref{fnd:continuation} gives a joint neighborhood with a common bound; at the fixed parameters take a short segment in the space of source strengths in this neighborhood.

Positive endpoints $\alpha^{(\varepsilon)}$ arbitrarily close to $\alpha$ and satisfying \eqref{fnd:independence} can be chosen by avoiding the countably many affine hyperplanes. The compact homotopy on the segment preserves the degree. For sufficiently small perturbations, $\rho_i\notin\Gamma_i(\alpha^{(\varepsilon)})$ as well. Indeed,
\[
 \Gamma_i(\alpha)\subseteq\mathcal C_i(\alpha),
\]
and only finitely many affine expressions defining $\Gamma_i$ can lie near the fixed parameter $\rho_i$ for nearby source strengths. Theorem~\ref{thm:degree-assembly} therefore computes this degree by the rectangular cumulative sum of $\mathscr J_{\alpha^{(\varepsilon)}}$.

\noindent\textit{Step 2. Stability of the partial sum.} We verify that this sum agrees with that of $\mathscr J_\alpha$ for all sufficiently small perturbations. Write the finite source product, before combining equal exponents, as
\[
 \prod_{p\in S}P_{A_p,B_p}(x,y)
   =\sum_{\nu}\epsilon_\nu
                  x^{a_\nu(\alpha)}y^{b_\nu(\alpha)},
 \qquad \epsilon_\nu\in\{-1,1\}.
\]
Each exponent is a nonnegative affine function of the source strengths. The remaining factor $F^{\chi-\ell}$ has fixed integer exponents and coefficients. With
\[
 r=\rho_1/(4\pi),\qquad s=\rho_2/(4\pi),
\]
only finitely many terms of this factor can contribute below the cut-offs $r,s$, uniformly for the chosen perturbations. Each relevant coordinate is consequently of the form
\[
 m+a_\nu(\alpha)\quad\text{or}\quad n+b_\nu(\alpha),
 \qquad m,n\in\mathbb Z_{\geq0}.
\]
When positive, these coordinates belong to $\Gamma_1(\alpha)/(4\pi)$ and $\Gamma_2(\alpha)/(4\pi)$, respectively, and hence to the corresponding sets $\mathcal C_i(\alpha)/(4\pi)$. Neither coordinate equals its cut-off. Continuity of the finitely many exponents therefore preserves every inequality defining the cumulative sum.

The coefficients are unchanged. Terms whose exponents coincide at the limiting source strengths are added before equal exponents are combined. The sums agree, and homotopy invariance under the perturbation of the source strengths gives the formula.

\noindent\textit{Step 3. The exceptional sets in (ii).} If the source strengths satisfy \eqref{fnd:independence}, their nonnegativity implies that they are all positive, and the weaker exclusions follow directly from Theorem~\ref{thm:degree-assembly}.

If all source strengths are integers, Lemma~\ref{fnd:sharp-compactness} gives compact homotopies along parameter segments in each rectangle $\mathcal C_{k,l}$, since these rectangles avoid
\[
 \Gamma_i=4\pi\mathbb N.
\]
Choose a representative there outside the locally finite sets $\mathcal C_i(\alpha)$ and apply the part just proved. Since $\mathscr J_\alpha$ has only integer exponents, its cumulative sum is also constant on that rectangle.

Formula~\eqref{eq:real-degree} therefore holds throughout the rectangle, proving the assertion under the weaker exclusions in the integer case.
\end{proof}
\begin{proposition}
\label{prop:quotient}
For \(a,b\in\Z_{\geq0}\),
\[
 P_{a+1,b+1}=FQ_{a,b},\qquad
 Q_{a,b}(x,y)=Q_{b,a}(y,x),\qquad
 Q_{a,b}(1,1)=\frac{(a+1)(b+1)(a+b+2)}2 .
\]
Every coefficient of \(Q_{a,b}\) is a nonnegative integer, and its constant coefficient is one.
\end{proposition}
\begin{proof}
Sum the finite geometric series in \eqref{eq:Q} after multiplication by \(1-x\):
\[
 (1-x)Q_{a,b}
 =\left(\sum_{i=0}^a(xy)^i\right)
       \left(\sum_{j=0}^by^j\right)
 -x^{a+1}\left(\sum_{i=0}^ay^i\right)
       \left(\sum_{j=0}^b(xy)^j\right).
\]
Multiplication by \((1-y)(1-xy)\) gives
\[
 (1-(xy)^{a+1})(1-y^{b+1})
 -x^{a+1}(1-y^{a+1})(1-(xy)^{b+1}),
\]
which expands to \(P_{a+1,b+1}\). Both \(F\) and \(P\) have the component-exchange symmetry, so cancellation in the integral domain \(\Z[x,y]\) gives the symmetry of \(Q\).

The finite sum in the definition has nonnegative integer coefficients and constant term one. At \((1,1)\), the summand for \(i,j\) has \(a+j-i+1\) terms. Summing their number gives
\[
 (a+1)(b+1)\left(a+1+\frac b2-\frac a2\right)
 =\frac{(a+1)(b+1)(a+b+2)}2.
\]
\end{proof}
\begin{proof}[Proof of Corollary~\ref{thm:integer-degree}]
For integer source strengths, Theorem~\ref{thm:main} applies throughout each rectangle $\mathcal C_{k,l}$. Since $\mathscr J_\alpha$ has only nonnegative integer exponents, it gives
\[
 d_{k,l}=\sum_{a=0}^k\sum_{b=0}^l[x^ay^b]\mathscr J_\alpha(x,y).
\]
Consequently
\[
 \sum_{k,l\geq0}d_{k,l}x^ky^l
 =\frac{\mathscr J_\alpha(x,y)}{(1-x)(1-y)}
 =D_{\mathrm{form}}(x,y).
\]
The second expression in \eqref{eq:degree-candidate} follows from
\[
 P_{a+1,b+1}=FQ_{a,b},
\]
proved in Proposition~\ref{prop:quotient}. This proves \eqref{eq:analytic-identification} throughout every integer parameter chamber. When $S=\varnothing$, the source product is one, giving the stated regular generating function.

For $\chi\leq0$, the factors $(1-x)^{\chi-1}$ and $(1-y)^{\chi-1}$ have strictly positive coefficients. The factor $(1-xy)^\chi$ and each $Q_{a,b}$ have nonnegative coefficients and constant coefficient one. Their product consequently has a strictly positive coefficient at every $x^ky^l$. Thus
\[
 d_{k,l}>0,
\]
and the nonzero-degree property gives the existence assertion.
\end{proof}
\begin{proof}[Proof of Theorem~\ref{thm:signed-continuation}]
\noindent\textit{Step 1. Dependence of the partial sum on the source strengths.} We first express the dependence of the cumulative sum on the source strengths and parameters. For every
\[
 \alpha_{p,i}>-1,
\]
the numbers $A_p,B_p$ are positive, so the source product can be written, before combining equal exponents, as
\[
 \prod_{p\in S}P_{A_p,B_p}(x,y)
 =\sum_\nu\epsilon_\nu
       x^{a_\nu(\alpha)}y^{b_\nu(\alpha)},
 \qquad \epsilon_\nu\in\{-1,1\}.
\]
There are finitely many terms, and the functions $a_\nu,b_\nu$ are continuous affine functions of the source strengths with nonnegative values. Write
\[
 F(x,y)^{\chi-\ell}
 =\sum_{m,n\geq0}f_{m,n}x^my^n.
\]
For $r=\rho_1/(4\pi)$ and $s=\rho_2/(4\pi)$, we therefore have
\begin{equation}\label{eq:signed-finite-sum}
 \mathscr D_\alpha(r,s)
 =\sum_{\nu,m,n}\epsilon_\nu f_{m,n}
   \mathbf 1_{\{m+a_\nu(\alpha)<r\}}
   \mathbf 1_{\{n+b_\nu(\alpha)<s\}}.
\end{equation}
On a bounded set of positive parameters only finitely many indices $m,n$ can contribute, uniformly in the source strengths, since the source exponents are nonnegative. Each positive number $4\pi(m+a_\nu(\alpha))$ belongs to $\Gamma_1(\alpha)$, and each positive number $4\pi(n+b_\nu(\alpha))$ belongs to $\Gamma_2(\alpha)$.

The inclusion $\Gamma_i(\alpha)\subset\mathcal C_i(\alpha)$ remains valid for all source strengths greater than $-1$, by the six initial pairs in the construction of $\Sigma_p(\alpha)$. Consequently, at any point satisfying the coarse exclusions, none of the finitely many relevant exponents equals its cut-off.

All the indicators in \eqref{eq:signed-finite-sum} are unchanged in a sufficiently small joint neighborhood of that point. In particular, the cumulative sum is locally constant there. The argument allows coincident exponents, because the terms are summed before equal exponents are combined.

\noindent\textit{Step 2. Continuation near nonnegative source strengths.} For part~\textup{(i)}, Proposition~\ref{fnd:signed-compactness} gives a joint neighborhood of $(\alpha^0,\rho^0)$ on which all solutions have a common bound and the associated compact maps form a continuous homotopy.

Shrink this neighborhood to a convex one on which the indicators in \eqref{eq:signed-finite-sum} are fixed.

The segment from $(\alpha^0,\rho^0)$ to any point of this neighborhood is admissible on a common ball. Thus its degree is equal to $d_{\alpha^0}(\rho^0_1,\rho^0_2)$.

Theorem~\ref{thm:main} identifies this value with the cumulative sum at the base point, and the preceding argument identifies that sum with the cumulative sum at every point of the neighborhood.

This proves~\textup{(i)}. No further exclusion is needed at the nearby points, since compactness holds throughout the chosen neighborhood.

\noindent\textit{Step 3. Continuation along the path.} For part~\textup{(ii)}, apply Proposition~\ref{fnd:signed-compactness} at every point of the path.

A finite cover of its compact image gives a uniform bound for all solutions along the path. The resulting continuous compact homotopy is admissible on a common ball, and hence the degree is constant.

The cumulative sum is also locally constant along the path by \eqref{eq:signed-finite-sum} and the stated exclusions. It is therefore constant on $[0,1]$. At $t=0$, Theorem~\ref{thm:main} identifies the degree with this sum, which proves the assertion for every $t$.
\end{proof}
\begin{proof}[Proof of Corollary~\ref{cor:small-signed-orders}]
Let $\mathbf0$ denote the array of zero source strengths on the fixed marked set $S$. At these source strengths the regular local mass list gives
\[
 \mathcal C_i(\mathbf0)=4\pi\mathbb N.
\]
For each $\rho^0\in K$, part~\textup{(i)} of Theorem~\ref{thm:signed-continuation} provides a product neighborhood of $(\mathbf0,\rho^0)$ on which the degree and the cumulative sum are both equal to their values at $(\mathbf0,\rho^0)$. Choose the neighborhoods of the zero source strengths to be balls in the maximum norm.

Choose finitely many of the parameter neighborhoods to cover $K$, and let $\varepsilon_K\in(0,1)$ be smaller than all the corresponding radii in the $\alpha$ variables. Applying the same local statement at both $(\alpha,\rho)$ and $(\mathbf0,\rho)$ then gives
\[
 d_\alpha(\rho_1,\rho_2)
 =d_{\mathbf0}(\rho_1,\rho_2)
 =\mathscr D_\alpha\left(\frac{\rho_1}{4\pi},
                         \frac{\rho_2}{4\pi}\right)
\]
for every $\rho\in K$ and every stated perturbation of the source strengths. The identity $P_{1,1}=F$ shows that retaining a marked point with two zero source strengths does not change the series. In particular,
\[
 \mathscr J_{\mathbf0}(x,y)
 =F(x,y)^{\chi-\ell}F(x,y)^\ell=F(x,y)^\chi.
\]
Its cumulative coefficients, or equivalently Corollary~\ref{thm:integer-degree} at zero source strengths, give
\[
 d_{\mathbf0}(\rho_1,\rho_2)
 =[x^ky^l]\frac{F(x,y)^\chi}{(1-x)(1-y)}.
\]
This is the asserted expression. For $\chi\leq0$, the factors $(1-x)^{\chi-1}$ and $(1-y)^{\chi-1}$ have strictly positive coefficients, while $(1-xy)^\chi$ has nonnegative coefficients and constant coefficient one.

Thus every indicated coefficient is positive. The existence assertion follows from the nonzero-degree property.
\end{proof}
\section{The degree jump and completion of the proof}
\label{sec:wall-completion}
\begin{proposition}[The local degree for source strengths greater than $-1$]
\label{univ:wall-jump}
Assume \eqref{fnd:independence} and
\[
 A_p,B_p>0,\qquad \frac{\rho_2}{4\pi}\in\mathbb Q\setminus\mathbb Z.
\]
At each critical value of $\rho_1$ use the weights selected in the proof. For each nondegenerate shadow solution $Z$ with configuration $(U,V,P,Q)$, the local degree jump is
\[
                  \Delta d_Z=(-1)^{U+|P|}d_S(Z).
\]
The choice of weights may depend on the critical value.
\end{proposition}
\begin{proof}[Completion of the proof of Theorem~\ref{thm:main}]
\noindent\textit{Step 1. The initial degree and the critical values.} First take $A_p,B_p>0$ satisfying \eqref{fnd:independence} and rational noninteger normalized parameter cut-offs such that
\[
 \rho_i\notin\mathcal C_i(\alpha),
\]
$i=1,2$. These rational choices are dense, since the sets defining the parameter exclusions are locally finite. At fixed $\rho_2$ start with
\[
 0<\rho_1^{\rm start}<\min(\rho_1,4\pi\min(\{1\}\cup\{A_p:p\in S\})).
\]
In Lemma~\ref{lem:base-strip} the deformation to $\rho_1=0$ and the triangular change of variables at that endpoint give
\[
 d_\alpha(\rho_1^{\rm start},\rho_2)=d_{\rm sc}(\rho_2;\boldsymbol\beta,0),
 \qquad 1+\beta_p=B_p.
\]
The scalar equation has coefficient $2\rho_2$. Its critical set in $\rho_2$ is contained in $\mathcal C_2(\alpha)$. Proposition~\ref{fnd:scalar-degree-internal} now computes the right-hand side from $[x^0]\mathscr J_\alpha(x,y)$.

Increase the first parameter to $\rho_1$. Take disjoint intervals about the finitely many values of $\Gamma_1$ crossed, with endpoints outside $\Gamma_1$. Compactness at fixed source strengths applies off these intervals.

At a critical value admitting no configuration with positive remainder, the local mass classification in the proof of Proposition~\ref{univ:wall-jump} excludes blow-up. Thus the degree can be computed on a common ball in a neighborhood of this value and its jump is zero.

At each remaining critical value apply Proposition~\ref{univ:wall-jump} with its chosen weights. At the two noncritical endpoints, Lemma~\ref{fnd:sharp-compactness} applies to the fixed generic source strengths and to a compact positive smooth-weight homotopy, so each endpoint degree is independent of that weight choice.

The mass estimates for solutions and their local description exclude unbounded sequences at the critical value itself. Choose a ball containing its bounded solution set. Compactness of the fixed-point maps on bounded sets gives a zero-free boundary also for nearby parameters.

A sequence of zeros outside this ball as the parameter approaches the critical value must escape every bounded set, after a subsequence, and the local description puts it in one of the isolating neighborhoods.

The matched widths force the projected zeros to escape the ball. For an outer radial bubble with profile exponent $s>0$ the first component after applying $K^{-1}$ differs from
\[
 V_R(z)=-\log(R^{2s}+|z|^{2s})
\]
by a bounded term in $H^1$. Fix $L>1$ and $d>0$. Then
\[
 |\partial_r V_R|\ge c_s/r\quad(LR<r<d),\qquad
 \int_{LR<|z|<d}|\nabla V_R|^2\,dA\ge c_s'\log\frac{d}{LR}.
\]
At a regular mixed point the outer potential is $-\log(R^4+|z^2-R^2B|^2)$. The same bound holds on this annulus if $L$ is sufficiently large, uniformly for $B$ in the compact set under consideration. The corrections on the complement are $o(1)$ in $H^1$.

For solutions, the normalized equation turns a uniform $C^0$ bound into a uniform $H^1$ bound. Use the disjoint neighborhoods of Proposition~\ref{univ:wall-jump}, counting each unordered configuration of each type once.

Denote these neighborhoods by $\mathcal U_Z$ and the ball containing the solutions at the critical value by $\mathcal B$. Then
\[
 \partial(\mathcal U_Z\setminus\overline{\mathcal B})
 \subset\partial\mathcal U_Z\cup\partial\mathcal B.
\]
Neither boundary on the right contains zeros. Excision preserves all local degrees. Additivity gives their sum outside $\mathcal B$, and the degrees on $\mathcal B$ agree by homotopy.

\noindent\textit{Step 2. Summation over the configurations.} For each retained configuration the formula for the degree of the shadow system applies with
\[
 \gamma=\rho_2-4\pi m_2>0
\]
and shifted parameters $B_p,A_p+B_p,A_p$, according to its three prescribed states, together with $U$ distinct shifts two.

The corresponding scalar parameter is noncritical by the exclusions for the shadow system, since the shifts equal to two contribute only integers. Use Proposition~\ref{fnd:scalar-degree-internal} at this parameter. The inequality for the exponents is
\[
 \beta<\rho_2/(4\pi)-m_2
\]
for a scalar term $y^\beta$. After multiplication by $y^{m_2}$, its exponent is
\[
 b=\beta+m_2<\rho_2/(4\pi).
\]
The unordered center factor is $\binom{\chi-|S|}{U+V}\binom{U+V}{U}$, and the sign of the local jump is $(-1)^{U+|P|}$. Thus the finite sum in Section~\ref{sec:assembly} applies with these shifts and the same signs. It yields
\[
 d_\alpha(\rho_1,\rho_2)
 =\sum_{a<\rho_1/(4\pi),\ b<\rho_2/(4\pi)}
             [x^ay^b]\,
             F(x,y)^{\chi-|S|}\prod_{p\in S}P_{A_p,B_p}(x,y).
\]

\noindent\textit{Step 3. All noncritical parameters at generic source strengths.} For fixed generic source strengths and a positive target outside both $\Gamma_i$, take a small rectangle in the corresponding component of the noncritical parameter set.

Choose there a rational noninteger normalized pair outside both locally finite $\mathcal C_i$ and join it to the target by the straight parameter segment, keeping the source strengths and smooth factors fixed. Lemma~\ref{fnd:sharp-compactness} gives an admissible homotopy.

In each source polynomial the first exponents are $0,A_p,A_p+B_p$ and the second are $0,B_p,A_p+B_p$. The factor $F^{\chi-|S|}$ adds only nonnegative integers.

Thus each positive exponent projection lies in $\Gamma_i/(4\pi)$ by \eqref{eq:sharp-mass-sets}, so neither strict test can change on this rectangle. This proves the weaker exclusion result for the globally independent source strengths.

\noindent\textit{Step 4. Continuation to arbitrary source strengths.} For arbitrary $\alpha_{p,i}>-1$ and
\[
 \rho_i\notin\mathcal C_i(\alpha),
\]
use the joint neighborhood of the source strengths and parameters in Proposition~\ref{fnd:signed-compactness}, with the smooth positive factors fixed. Shrink it to a convex neighborhood with positive shifted parameters and bounded cut-offs.

Use \eqref{eq:signed-finite-sum} before combining coincident exponents. The integer shifts of any potentially relevant term lie in one fixed finite set since all source exponents are nonnegative.

At the target each positive exponent projection lies in $\Gamma_i/(4\pi)\subset\mathcal C_i/(4\pi)$; zero also cannot equal a positive cut-off. There is consequently a smaller neighborhood on which both strict inequalities have the same truth value for every relevant term.

In it choose generic source strengths, then rational noninteger normalized parameters for those source strengths outside the locally finite coarse sets.

The straight joint homotopy stays in the neighborhood of the common bound, so the target degree and coefficient sum equal their generic values. This gives \eqref{eq:real-degree}.

For fixed integer source strengths
\[
 \Gamma_i=4\pi\mathbb N.
\]
In an open rectangle in its complement choose an interior representative outside the two coarse sets and use the formula just proved. The straight parameter segment at the original integer source strengths is admissible by Lemma~\ref{fnd:sharp-compactness}.

All exponents are integers, so neither strict cut-off changes. This gives the weaker integer exclusions and the corresponding integer generating-function identity in Section~\ref{sec:assembly}.
\end{proof}
\section{Compactness and degree of the shadow systems}
\label{fnd:section}
We first establish compactness as the source strengths vary in $(-1,\infty)$. We then reduce the degree calculation for nonnegative source strengths to one at fixed positive source strengths and choose common smooth factors for the finitely many shadow systems below the parameter cut-offs. Write
\[
 A_p=1+\alpha_{p,1},\qquad B_p=1+\alpha_{p,2}.
\]
For positive source strengths we assume $A_p,B_p>1$ and \eqref{fnd:independence}.

The constants may depend on these fixed source strengths. In the shadow-system calculations below, local masses are divided by $4\pi$.

We use $U$ ordinary free centers, $V$ mixed free centers, and disjoint sets $\mathcal P,\mathcal Q\subset S$ at which pure and mixed singular bubbles occur, respectively. The concentrated masses of the first and second components, in these units, are
\[
 \tau=U+2V+\sum_{p\in\mathcal P}A_p+
                   \sum_{p\in\mathcal Q}(A_p+B_p),
 \qquad m_2=V+\sum_{p\in\mathcal Q}B_p.
\]
The remaining mass in the second component is $\gamma=\rho_2-4\pi m_2>0$.
\begin{lemma}[Unitary holonomy and integrated curvature]
\label{fnd:holonomy}
Remove finitely many disjoint closed disks from a disk to obtain $\Omega$. Let $\mathcal A$ be a $C^1$ connection on $\overline\Omega$ with values in $\mathfrak{su}(m)$, $m\ge2$. Define parallel transport by
\[
 P'=-\mathcal A(\dot\ell)P.
\]
Join the inner boundaries by cuts to a common base point on the outer boundary. Denote the outer holonomy by $H_0$ and the counterclockwise inner holonomies at this base point by $H_j$. In the order induced by the cuts,
\[
       \left\|H_0-\prod_jH_j\right\|
               \leq\int_\Omega|d\mathcal A+\mathcal A\wedge\mathcal A|.
\]
Use the operator norm for the matrices and for the curvature coefficient. The bound is independent of the distances between the holes and of the size of $\mathcal A$.
\end{lemma}
\begin{proof}
\smallskip\noindent\textit{Step 1. Parallel transport on a rectangle.}\quad On a rectangle $[0,1]^2$, parallel transport first along the left edge and then along horizontal segments gives a unitary gauge with $\mathcal A_s=0$ and
\[
 \mathcal A_t(0,t)=0.
\]
The curvature is then $\mathcal F_{st}=\partial_s\mathcal A_t$, so $\mathcal A_t(1,t)=\int_0^1\mathcal F_{st}(s,t)\,ds$. Three edge transports are the identity. Unitarity of the fourth and its transport equation give
\[
 \|H_{\partial\Omega}-I\|
       \leq\int_0^1\|\mathcal A_t(1,t)\|\,dt
       \leq\int_0^1\int_0^1\|\mathcal F_{st}\|\,ds\,dt.
\]
For a positively oriented triangle with vertices $v_0,v_1,v_2$, pull back the connection by
\[
 x(s,t)=(1-s)v_0+s\bigl((1-t)v_1+tv_2\bigr).
\]
One edge of the rectangle collapses to $v_0$, and the determinant of the derivative on $s>0$ is
\[
 s\det(v_1-v_0,v_2-v_0)>0.
\]
The rectangle estimate therefore bounds the norm of the triangular boundary holonomy minus the identity by the curvature integral on the triangle. Triangulate a simple polygon by diagonals.

At each diagonal, the two oriented transports cancel when taking the product of the two subpolygon boundary holonomies based at an endpoint of the diagonal. Changing the base point acts by unitary conjugation. Iterating this identity and using
\[
 \|UV-I\|\leq\|U-I\|+\|V-I\|
\]
for unitaries gives the bound on the polygon. For a rounded cut domain, first take a smooth parallel curve just inside its boundary and approximate that curve by inscribed polygons.

The paths converge uniformly and their velocities converge in $L^1$ as the polygon mesh and then the inward displacement tend to zero.

The transport coefficients then converge in $L^1$, and the transport equations give convergence of their endpoint matrices.

The integrals on the interiors converge by boundedness of the curvature and convergence of the areas of the symmetric differences. This proves the estimate for the rounded cut domains used below.

\smallskip\noindent\textit{Step 2. The cut disk.}\quad Cut the punctured disk along the chosen arcs. Its boundary word consists of the outer loop and the inverses of the based inner loops.

Apply the preceding estimate. More formally, remove thin strips about the cuts and round the corners, apply the estimate for positive strip width, and then let the width tend to zero for the fixed connection and fixed hole radii.

The two transports along a cut converge to the required conjugations. The cuts have zero area. Multiplication and conjugation by unitary matrices preserve the norm, giving the assertion.
\end{proof}
\begin{proof}[Proof of Lemma~\ref{fnd:local-unitary}]
\noindent\textit{Step 1. Local estimates.}\quad On a smaller disk, Green's representation expresses the regular part of $U_i^k$, up to a constant, as the sum of the potentials of $K_{ij}F_j^k$ and a function with uniformly bounded gradient.

For this representation use the Dirichlet solution for the bounded smooth forcing, the harmonic extension of the boundary data after subtracting its mean, and the disk Dirichlet Green function.

Its regular kernel has bounded first derivatives on the smaller disk for every integration point. Thus the singular kernel may be taken to be $-(2\pi)^{-1}\log|z-y|$, and the remaining gradients have the asserted uniform bound. The same representation in the coordinates
\[
 z=d_ky,
\]
with fields $U_i^k(d_ky)+2\log d_k$, has remainder gradients $O(d_k)$ and the same source logarithms. Constants can again absorb the divergent constant parts of the kernels, by replacing $\log|y-z|$ with $\log(|y-z|/(1+|z|))$ in the rescaled representation.

Consider a disk of radius comparable to one away from the source. On a smaller disk write the field as the sum of a constant, a function of bounded oscillation and $2P_i-P_j$, where
\[
 P_j\geq0
\]
is the local logarithmic potential of the density of mass $m_j$ on the larger disk. The remainder oscillation is bounded by the separated kernel derivative estimate and by the mass bound. In a fixed scaled disk take
\[
 P_j=(2\pi)^{-1}\int \log(C_0/|z-y|)F_j(y)\,dy,
\]
choosing $C_0$ so the logarithm is nonnegative. Its area average on the smaller disk is $O(m_j)$. Let $c_i$ be the additive constant. Jensen's inequality on the smaller disk and the oscillation bound give
\[
 \int_{\rm smaller\ disk} F_i\geq c e^{c_i}
\]
with a uniform $c>0$. For $m_i>0$, apply Jensen to $F_i(y)\,dy/m_i$:
\[
 \int_{\rm smaller\ disk}e^{2qP_i(z)}\,dz
 \leq \sup_y\int_{\rm smaller\ disk}
               (C_0/|z-y|)^{q m_i/\pi}\,dz.
\]
The right-hand side is bounded if $m_i$ is sufficiently small for the fixed $q$. For $m_i=0$ the potential term is zero. The factor $\exp(-P_j)$ is at most one. This gives an $L^q$ density bound and gives $o(1)$ in $L^q$ when both local masses tend to zero, by the bound on $e^{c_i}$.

With $q>2$, H\"older for the near kernel $|z-y|^{-1}$ gives the same boundedness or vanishing for its gradient contribution. These estimates are invariant under simultaneous dilation of disk and densities.

They apply on fixed compact sets off the explicit source and the finite set of concentration points, and on annuli after scaling their radius to one.

In particular, on compact subsets disjoint from this finite set and the source the densities are uniformly bounded and have subsequential local uniform limits.

To identify which limits are positive, fix a base disk away from that set and use the potential representation and the small-mass bounds on a finite connected chain of disks avoiding the set.

Differences of the fields at points on this chain are uniformly bounded and converge along a subsequence by the local gradient and compactness estimates.

For each component the additive constant at the base point then converges to a finite number or to $-\infty$. The limiting density of a component is therefore strictly positive throughout the complement or identically zero there.

\noindent\textit{Step 2. Limits on circles.}\quad On a circle around an isolated blow-up point, first take $k\to\infty$ at fixed radius $r$, then $r\downarrow0$. The mass away from the center tends to zero on the rescaled annuli of radii comparable to $r$. The preceding estimate on smaller annular disks gives
\[
 r^2F_i^k=o(1)
\]
there in the two successive limits and gives $o(1)$ for the local gradient in these scaled coordinates. For fixed $0<\eta<1/2$, write
\[
 M_{i,k}(\eta r)=\int_{D_{\eta r}}F_i^k\,dA.
\]
Replacing $\mathbf 1_{D_{\eta r}}F_i^k\,dA$ by $M_{i,k}(\eta r)\delta_0$ changes the rescaled gradient on $\partial D_r$ by at most
\[
 C\eta M_{i,k}(\eta r).
\]
The mass in the remaining intermediate annuli tends to zero, and mass beyond $r/\eta$ contributes $O(\eta)$ to the relative gradient by the separated kernel bound. Thus the regular part $V_i^k=U_i^k-2\alpha_i^k\log|z|$ satisfies
\[
 r\nabla V_i^k=-2(Kx)_i e_r+o(1),\qquad r^2 F_i^k=o(1)
 \quad\text{on }\partial D_r
\]
in the same successive limits. The boundary statements follow uniformly from finitely many scaled local disks. Identical estimates apply at a regular point after absorbing the smooth source logarithms into its canonical field and setting the explicit source strengths to zero.

To compute the total mass on the limiting rescaled plane, apply the same argument on $|z|=d_kr$, $d_k\to0$. First let $k\to\infty$ with $r$ fixed, then let
\[
 r\to\infty.
\]
The mass of the limiting rescaled measure in every annulus comparable to $r$ tends to zero.

For fixed $\eta$, its disk mass inside $\eta r$ tends to the total mass on the rescaled plane, the remaining intermediate mass tends to zero, and the mass escaping beyond $r/\eta$ contributes $O(\eta)$ by the separated gradient bound, uniformly before taking
\[
 k\to\infty.
\]
The smooth remainders in the annular coordinates have gradients $O(d_kr)$ and the forcing is $O(d_k^2r^2)$. This proves the same two boundary statements, with $x$ denoting the total mass on the rescaled plane and gradients taken for the rescaled regular parts.

These arguments only use that the remaining mass tends to zero on its small annuli or finite-mass tails. They do not require a strict exponent margin for a surviving residual density.

\noindent\textit{Step 3. The radial Pohozaev identity.}\quad For fixed $k$, multiply the equations for $V$ by $K^{-1}(z\cdot\nabla V)$, integrate over $D_r$, and use the divergence identity for the quadratic gradient terms. The result is
\begin{align*}
 \frac r2\int_{\partial D_r}
 [ (\partial_r V)^t K^{-1}\partial_r V
   -(\partial_{\tan} V)^t K^{-1}\partial_{\tan} V]
   &=2\sum_i(1+\alpha_i^k)\int_{D_r} F_i^k
       -r\int_{\partial D_r}\sum_i F_i^k\\
   &\quad+\int_{D_r}(z\cdot\nabla V)^tK^{-1}b^k.
\end{align*}
For the density term we used $F_i^k z\cdot\nabla V_i^k=z\cdot\nabla F_i^k-2\alpha_i^k F_i^k$. To justify integration through zero, first excise a circle of radius $\epsilon$ for fixed $k$. Continuity of the regular part implies
\[
 F_i^k=O(\epsilon^{2\alpha_i^k})
\]
on the circle. The nearby and separated logarithmic kernel gradient bounds with this power give
\[
 \epsilon|\nabla V_i^k|=o(1),
\]
so both inner boundary terms tend to zero. The forcing term on the right is $O(r^2)$ in physical coordinates because
\[
 \int_{D_r}|z-y|^{-1}\,dz\leq Cr
\]
uniformly in $y$ and the remaining gradients are bounded. On a physical circle of radius $d_k r$ the bound is $O(d_k^2 r^2)$. The boundary limits just obtained now give
\[
 x^tKx=2a\cdot x,
\]
both for a shrinking disk about a point mass at a finite point and for the total finite mass of a rescaled plane.

\noindent\textit{Step 4. Holonomy.}\quad Choose the real triplet $\Phi$ with sum zero and successive differences $U_1,U_2$. Write
\[
 d^cv=-v_y\,dx+v_x\,dy.
\]
On the punctured disk use the skew-Hermitian traceless connection whose diagonal is $(i/2)d^c\Phi_j$, whose neighboring upper entries are $\sqrt{F_j}\,dz/2$, and whose neighboring lower entries are their negative conjugates.

The other entries vanish. For a neighboring upper entry the differential and the diagonal commutator cancel because
\[
 \Phi_j-\Phi_{j+1}=U_j.
\]
Non-neighboring products vanish since
\[
 dz\wedge dz=0.
\]
The diagonal curvature is $i/2$ times the smooth forcing for $\Phi_j$, since $\Delta\Phi_j=-F_j+F_{j-1}$ plus that forcing and $dz\wedge d\bar z=-2i\,dx\wedge dy$. The curvature therefore has integral $O(r^2)$ on a physical disk of radius $r$.

Apply Lemma~\ref{fnd:holonomy} to this connection in a unitary frame of dimension three. For each fixed $k$, first shrink the circle about the source. The off-diagonal transports have norm integrals tending to zero since
\[
 1+\alpha_j^k>0.
\]
The preceding regular gradient bound gives vanishing regular diagonal flux. Up to conjugation, the limiting source holonomy is diagonal with entries $\exp(-2\pi i (\lambda_{\alpha^k})_j)$.

On the enclosing circles used above the off-diagonal integrals also tend to zero, by the density estimate on the annulus.

In the diagonal transported frame their removal changes the holonomy by $o(1)$. The exact distributional flux of the remaining diagonal one-forms is $2\pi i(\lambda_{\alpha^k}-\iota(x^k))$ plus the smooth forcing integral, where $4\pi x^k$ are the enclosed masses.

The boundary curvature bound compares the enclosing and source holonomies up to conjugation. Taking limits gives equality of the unordered three eigenvalues of the two diagonal transports.

The same proof on the physical circles of radius $d_kr$ gives that equality for the total mass on the rescaled plane, by taking $k\to\infty$ before
\[
 r\to\infty.
\]
Because $\lambda_a-\lambda_{\alpha}=(1,0,-1)$ is integral, equality of eigenvalues gives the asserted permutation and lattice vector. At a regular point the connection extends through the center for each fixed $k$ and the inner holonomy is the identity, with the same conclusion at
\[
 \alpha_1=\alpha_2=0.
\]
The Pohozaev identity gives $|\lambda_a-\iota(x)|^2=|\lambda_a|^2$ because $|\iota(x)|^2=x^tKx$ and $\lambda_a\cdot\iota(x)=a\cdot x$. Consequently $2(P\lambda_a)\cdot n+|n|^2=0$. The entries of $\lambda_a$ are rational combinations of $A,B$, so independence of $1,A,B$ forces
\[
 |n|^2=0.
\]
The six permutations give the displayed list. If both limiting source strengths are integers, the entries of $\lambda_a$ are pairwise congruent modulo integers and the congruence implies
\[
 \iota(x)\in\mathbb Z^3.
\]
For $A=B=1$ the ellipse has each coordinate less than $3$. Its integer pairs are obtained by substituting $x_1=0,1,2$, giving precisely the regular list and $(0,0)$.

\noindent\textit{Step 5. The estimate near the source.}\quad We prove the bound on $|z|^2F_i^k$ for fixed generic source strengths by contradiction. If it fails on a smaller disk, take $z_k\to0$ with
\[
 d_k=|z_k|,\qquad d_k^2F_i^k(z_k)\to\infty.
\]
Here the assertion is used on an isolating disk with zero as its only possible concentration point. Rescale by $d_k$ and pass to the finite measures on the rescaled plane. The local estimate gives positive mass at a regular blow-up point obtained as a limit of $z_k/d_k$.

It also shows that the limiting measures are sums of finitely many point masses and the residual densities considered above. At any finite blow-up point with local mass pair $x$, a component with nonzero residual density has
\[
 (Kx)_i\leq 1+\alpha_i,
\]
where the source strength is zero at a regular blow-up point. For a finite measure without a point mass at zero, use the identity
\[
 \frac1{2\pi}\int_0^{2\pi}\log|re^{it}-y|\,dt
         =\log\max\{r,|y|\}.
\]
It follows from the mean-value property after factoring out the larger radius, with the equality at $|y|=r$ obtained by dominated convergence: for $1/2\leq t_0<1$ one has
\[
 |1-t_0e^{it}|^2\geq 2\sin^2(t/2)
\]
and
\[
 |1-t_0e^{it}|\leq2.
\]
For the remaining measure inside a disk of radius $\delta<1$ the absolute value of the mean is bounded by its mass times $|\log r|$ for $r<\delta$. The part outside this disk is $O_\delta(1)$ locally, using the subtracted kernel for a rescaled plane.

Taking $r\downarrow0$ and then $\delta\downarrow0$ gives $o(|\log r|)$. Thus the circular mean of the surviving logarithmic density is $2(\alpha_i-(Kx)_i)\log r+o(|\log r|)$. If
\[
 (Kx)_i>1+\alpha_i,
\]
Jensen's inequality bounds the circular mean of the density below by $r^{-2-\varepsilon}$ for some $\varepsilon>0$ and all small $r$, contradicting integrability. This proves the stated inequality.

If both components survived, multiplying these inequalities by $x_1,x_2$ and using the already proved Pohozaev identity would give
\[
 x^tKx\leq a\cdot x=x^tKx/2
\]
at the regular blow-up point of positive mass. Thus at most one component has residual density. If one survives, the regular list and its integrability inequality show that the other component has positive local mass at every regular blow-up point.

If neither survives, choose a component with positive local mass at the regular blow-up point on $|y|=1$.

In either case there is a component with zero residual density whose total mass on the rescaled plane exceeds its local mass at the origin by a strictly positive integer in $4\pi$ units.

The central pair and the pair of total masses on the rescaled plane both lie on the generic list proved above.

In each coordinate any difference of members of that list is either zero or a nonzero rational combination of $A,B$, and cannot be a nonzero integer. This contradiction proves the centered estimate.
\end{proof}
\begin{lemma}[Finite-mass profiles with one source]
\label{fnd:entire-profiles}
Assume that the regular parts are continuous at $0$ and that the gradients have the logarithmic potential representation, with no additional entire harmonic gradient.
\begin{itemize}
\item[\textup{(i)}]
Let $U$ solve the canonical Toda system on $\mathbb R^2$ with zero forcing, one source and shifted parameters $A,B>0$. If $F_i=e^{U_i}$ have finite masses, then
\[
 \int_{\mathbb R^2}F_1=\int_{\mathbb R^2}F_2=4\pi(A+B).
\]
\item[\textup{(ii)}]
For a finite-mass solution of $\Delta U+2e^U=4\pi(q-1)\delta_0$, $q>0$, the mass is $4\pi q$ and $e^U$ is given by
\[
 \frac{4q^2\ell^2 |z|^{2q-2}}
      {(\ell^2+|z^q-D|^2)^2},\qquad \ell>0.
\]
Here $D\in\mathbb C$ for integer $q$ and $D=0$ otherwise, so only radial powers occur in the latter case. For $q=1$ there is no source and the formula allows translations.
\end{itemize}
\end{lemma}
\begin{proof}
\smallskip\noindent\textit{Step 1. Behavior at infinity.}\quad At infinity the annular small-mass argument in the proof of Lemma~\ref{fnd:local-unitary} applies to the given logarithmic potentials. It gives
\[
 r^2F_i=o(1)
\]
and the relative gradient of the potential of $F_i$ is $-(\int F_i)/(2\pi)$ times the radial unit vector, up to $o(1)$ on the circle. At zero, subtract the explicit source logarithm. The local bound
\[
 F_i\leq C|z|^{2\alpha_i}
\]
from the continuous regular part gives
\[
 r^2F_i=o(1)
\]
and vanishing relative gradient of its potential there. These conclusions are uniform on each scaled compact annulus.

On such an annulus, bounded scaled logarithmic gradients and vanishing scaled densities give vanishing scaled first derivatives of the densities. Interior estimates with the small scaled $C^{0,\theta}$ forcing give $C^{2,\theta}$ convergence modulo constants. The identity
\[
 D^2 F_i=F_i(DU_i\otimes DU_i+D^2U_i)
\]
in scaled variables then makes the scaled first density derivatives small in $C^{0,\theta}$. Interior estimates for the differentiated equation give convergence through order three, also for the scalar potential with coefficient two.

\smallskip\noindent\textit{Step 2. The Toda system.}\quad For the Toda system set $a=((2U_1+U_2)/3)_z$ and $b=((U_1+2U_2)/3)_z$. Then $a_{\bar z}=-F_1/4$, $b_{\bar z}=-F_2/4$ off zero. Set
\[
 c=a_z+ab-a^2,\qquad W_2=b_z+c-b^2,\qquad W_3=c_z-bc.
\]
Both $W_2,W_3$ are holomorphic off zero. Indeed $(F_1)_z=(2a-b)F_1$, $(F_2)_z=(2b-a)F_2$ and $c_{\bar z}=-aF_2/4$. Thus $(W_2)_{\bar z}=(-(2b-a)-a+2b)F_2/4=0$ and
\[
 (W_3)_{\bar z}=(-a_z-a(2b-a)+c+ba)F_2/4=0.
\]
By the annular derivative limits each $z^j W_j$, $j=2,3$, has a finite limit at zero and at infinity. Hence each is constant. The operator of these invariants factors as
\[
 \partial_z^3+W_2\partial_z+W_3
       =(\partial_z-b)(\partial_z+b-a)(\partial_z+a).
\]
With leading fields $a=\xi/z$ and $b=\upsilon/z$, its Euler indicial roots are $-\xi,1+\xi-\upsilon,2+\upsilon$. At zero $\xi=(2\alpha_1+\alpha_2)/3$ and
\[
 \upsilon=(\alpha_1+2\alpha_2)/3,
\]
giving consecutive gaps $A,B$. If the masses at infinity in $4\pi$ units are $(x,y)$, the two coefficients there are $\xi-x,\upsilon-y$. Radial integration gives
\[
 \overline U_i(r)=2(\alpha_i-(K(x,y))_i)\log r+o(\log r).
\]
Circular Jensen and integrability imply $2(\alpha_i-(K(x,y))_i)\le-2$. Thus the differences between consecutive roots at infinity satisfy
\[
 A-2x+y\leq0,\qquad B+x-2y\leq0.
\]
The polynomial is the same at the two ends and its roots are distinct. Thus the roots are reversed, giving $x=y=2+\xi+\upsilon=A+B$.

\smallskip\noindent\textit{Step 3. The scalar equation.}\quad For the scalar equation, $L=U_{zz}-U_z^2/2$ is holomorphic since
\[
 U_{z\bar z}=-e^U/2.
\]
The same two annular estimates give $L=(1-q^2)/(2z^2)$, by boundedness of $z^2L$ at both ends and the leading coefficient at zero. Set
\[
 f=e^{-U/2}.
\]
The equations $f_{zz}+(L/2)f=0$ and its conjugate give, on any branch of the logarithm, a constant combination of the products $z^{(1\pm q)/2}\bar z^{(1\pm q)/2}$.

Indeed the first equation has those two holomorphic basis functions: solving for their coefficients in $(f,f_z)$ by the nonzero determinant of the two solutions and their first derivatives shows those coefficients are antiholomorphic.

The conjugate equation does the same for each of them. Reality gives
\[
 f=|z|^{1-q}(a+c z^q+\bar c\,\bar z^q+b|z|^{2q}).
\]
A full turn eliminates $c$ for noninteger $q$. The equation $ff_{z\bar z}-f_zf_{\bar z}=1/4$ yields $q^2(ab-|c|^2)=1/4$. Continuity of the regular part at zero gives $a>0$, so $b>0$. Complete the square with
\[
 D=-\bar c/b,\qquad \ell^2=(ab-|c|^2)/b^2.
\]
For integer $q$ change variables $w=z^q-D$ with multiplicity $q$. For noninteger $q$ use $t=r^{2q}$ in the radial integral. Both integrations give $4\pi q$.

Conversely, for the stated $\ell,D$ set
\[
 b=(2q\ell)^{-1},\qquad c=-b\bar D,
\]
and $a=b(\ell^2+|D|^2)$. They give $q^2(ab-|c|^2)=1/4$, so the same calculation verifies the scalar equation off zero. At zero the regular part is continuous, with its $r\partial_r$ derivative tending to zero by the explicit formula. The source is therefore $4\pi(q-1)\delta_0$.
\end{proof}
\begin{proof}[Proof of Lemma~\ref{fnd:coarse-local}]
\noindent\textit{Step 1. Finiteness.}\quad A positive integer increment increases the bounded sum $s+t$ by at least one. A zero increment either fixes the pair or selects the other root of the quadratic equation in one coordinate. These reflections are $(s,t)\mapsto(A+t-s,t)$ and
\[
 (s,t)\mapsto(s,B+s-t).
\]
Denote the two reflections by $R_1,R_2$. Their orbits have at most six
points, since
\[
 R_1^2=R_2^2=(R_1R_2)^3=\Id.
\]
For either coordinate the positive integer increment has finitely many
possible values, each giving at most two roots for the other coordinate.
Since the number of positive increments is bounded,
\[
 \#\mathscr S(A,B)<\infty.
\]

\noindent\textit{Step 2. Choice of radii.}\quad
Let $\mathscr H$ be the set of nonnegative pairs $x$ satisfying both
identities in Lemma~\ref{fnd:local-unitary}\textup{(i)}.
For each permutation $P$ occurring there,
\[
 n=\lambda_a-\iota(x)-P\lambda_a\in\mathbb Z^3,\qquad
 |n|\le 2|\lambda_a|,
\]
since $|\lambda_a-\iota(x)|=|\lambda_a|$ by the Pohozaev identity.
Each $(P,n)$ determines $x$, so $\mathscr H$ is finite.
For $v=(v_1,v_2)$ write $\sigma(v)=v_1+v_2$ and set
\[
 x_k(r)=\frac1{4\pi}\int_{D_r}(F_1^k,F_2^k)\,dA.
\]
The pair $(0,0)$ is an initial pair. For a nonzero limiting local mass
pair $m$, a diagonal choice of continuity radii gives $R_k\downarrow0$
with $x_k(R_k)\to m$. Write
\[
 \sigma(\mathscr H)\cap[0,\sigma(m)]=\{c_0,\ldots,c_N\},
 \qquad 0=c_0<\cdots<c_N=\sigma(m).
\]
For $1\le j\le N$, set $b_j=(c_{j-1}+c_j)/2$.
Since $\sigma(x_k(r))$ is continuous in $r$ and vanishes at $r=0$,
for all large $k$ we can choose
\[
 0<r_{j,k}<R_k,\qquad \sigma(x_k(r_{j,k}))=b_j.
\]
Pass to a subsequence on which all radius ratios converge in $[0,\infty]$.
Group the radii by
\[
 j\sim h\quad\Longleftrightarrow\quad
 0<\lim_{k\to\infty}r_{j,k}/r_{h,k}<\infty.
\]
Choose one representative $\delta_{\ell,k}$ from each group, ordered
so that $\delta_{\ell,k}/\delta_{\ell+1,k}\to0$ for $1\le\ell<L$,
where $L$ is the number of groups. After a further subsequence,
define the endpoint pairs by
\[
 \frac{\delta_{\ell,k}^2}{4\pi}
 (F_1^k,F_2^k)(\delta_{\ell,k}y)\,dA_y
 \rightharpoonup\mu_\ell\quad\text{locally on }\mathbb R^2,
 \qquad
 m_\ell^-=\mu_\ell(\{0\}),\quad m_\ell^+=\mu_\ell(\mathbb R^2).
\]
The two circular limits in Steps~2--4 of the proof of
Lemma~\ref{fnd:local-unitary} give $m_\ell^\pm\in\mathscr H$.
Every fixed multiple of $\delta_{\ell,k}$ tends to zero, and successive
scales are separated. Monotonicity of the disk masses therefore gives,
coordinatewise,
\[
 (0,0)\le m_\ell^-\le m_\ell^+\le m,\qquad
 m_\ell^+\le m_{\ell+1}^-\quad(1\le\ell<L).
\]
For $j$ in group $\ell$, small and large fixed continuity radii on that
scale bracket $r_{j,k}/\delta_{\ell,k}$. Neither endpoint sum equals
a midpoint. Hence
\[
 \sigma(m_\ell^-)<b_j<\sigma(m_\ell^+).
\]
Any gap left by the intervals $[\sigma(m_\ell^-),\sigma(m_\ell^+)]$
in $[0,\sigma(m)]$ would have endpoints among the $c_j$ and hence
contain some $b_j$. Hence the intervals cover $[0,\sigma(m)]$,
and the coordinatewise inequalities give
\[
 m_1^-=(0,0),\qquad
 m_\ell^+=m_{\ell+1}^-\quad(1\le\ell<L),\qquad
 m_L^+=m.
\]

\noindent\textit{Step 3. Local masses on each scale.}\quad
Fix a group and write $r_k=\delta_{\ell,k}$ and $m^\pm=m_\ell^\pm$.
Let $x(y)=\mu_\ell(\{y\})$ and let $g_i$ be the residual density of
component $i$. Lemma~\ref{fnd:local-unitary}\textup{(iii)} gives
\[
 x(y)\in\mathbb Z_{\ge0}^2\qquad(y\ne0).
\]
Let $\mathcal B_\ell$ be the finite blow-up set. The fixed-annulus
potential estimate in Step~1 of that proof applies on rescaled disks
away from the source and gives, for each $i$,
\[
 g_i\equiv0\quad\text{or}\quad g_i>0
 \qquad\text{on }\mathbb R^2\setminus(\mathcal B_\ell\cup\{0\}).
\]
At a point $y$ with $x=x(y)\ne0$, take $\alpha_i=0$ if $y\ne0$.
If both residual densities were positive, Jensen's inequality on
circles, as in Step~5 of the same proof, would give
\[
 (Kx)_i\le 1+\alpha_i,\qquad i=1,2.
\]
Multiplying by $x_i$, summing in $i$ and applying the Pohozaev identity give
\[
 x^tKx\leq x^tKx/2,
\]
a contradiction. In the coordinate with vanishing residual density the increase from $m^-$ to $m^+$ is therefore a sum of integer masses at regular blow-up points. This is precisely an edge of the generation rule, also when that increase is zero.

If there is no point mass, then
\[
 m^-=0.
\]
On the rescaled plane the limit has zero forcing. There is no additional entire harmonic gradient. The gradient of the disk remainder is $O(r_k)$ in the new variables.

On each fixed compact set the mass outside $D_L$ contributes at most $C/L$ to the gradient. On each compact set there is convergence of the regular parts of the surviving fields. At the origin choose a uniform $\alpha_*>-1$ below all the source strengths $\alpha_i^k$ and a $p_0>1$ such that
\[
 p_0\alpha_*>-1
\]
if $\alpha_*<0$. The factors $|z|^{2\alpha_i^k}$ are uniformly in $L^{p_0}$ on the small disk. Choose $1<p<p_0$ and $q$ with
\[
 1/p=1/p_0+1/q.
\]
Since the limiting measure has no point mass at the origin, shrink the disk to a continuity radius at which the masses are small enough for the Jensen estimate of Lemma~\ref{fnd:local-unitary} with exponent $q$.

The additive constants are bounded above by matching on a fixed compact annulus off zero.

H\"older then bounds the densities in $L^p$ on a smaller disk. The local Poisson estimate and the compact embedding of $W^{2,p}$ into $C^0$ on still smaller disks give convergence of the continuous regular parts of the surviving fields.

The surviving field constants are finite after fixing a reference point off zero, and the limit equations inherit the source logarithms.

If both fields survive, Lemma~\ref{fnd:entire-profiles} gives the initial pair $(A+B,A+B)$. If one survives, its scalar equation gives $(A,0)$ or $(0,B)$. At least one survives since
\[
 m^+>m^-
\]
in the coordinate sum. Hence each edge goes to an initial pair or follows the stipulated generation rule. The last endpoint belongs to $\mathscr S(A,B)$, as required.
\end{proof}
\begin{lemma}[Local mass at an isolated scalar blow-up point]
\label{fnd:scalar-atom}
On a disk centered at $0$, suppose that $v_k$ is continuous at $0$ and
\[
 -\Delta v_k=H_k=|z|^{2\vartheta_k}b_k e^{v_k},\qquad
 \vartheta_k\to\vartheta>-1,\quad \inf_k\inf b_k>0,\quad \sup_k\|b_k\|_{C^1}<\infty.
\]
Assume uniformly bounded total mass, bounded oscillation on a fixed boundary circle, and local upper bounds for $v_k$ away from $0$. If $H_k\,dA\rightharpoonup\mu$ with $\mu(\{0\})>0$, then
\[
 \mu(\{0\})=8\pi(1+\vartheta),\qquad H_k\longrightarrow0\quad\text{locally away from }0.
\]
\end{lemma}
\begin{proof}
\smallskip\noindent\textit{Step 1. The Green potential.}\quad On the fixed circle write
\[
 v_k=c_k+b_k^{\rm harm}+P_k,
\]
where $P_k$ is the zero-boundary Dirichlet potential of $H_k$, $c_k$ is the mean of $v_k$ on that circle, and the harmonic remainder has mean zero there. The boundary hypotheses give
\[
 c_k\le C
\]
and locally uniform bounds on the harmonic remainders and their derivatives. On compact sets off zero, $H_k$ are uniformly bounded.

Splitting the logarithmic kernel and its first derivative into a small neighborhood of the observation point and its complement yields convergence of $P_k$ and its gradient there along a subsequence after taking the weak measure limit.

The near integrals are bounded by $C\epsilon^2(1+|\log\epsilon|)$ and $C\epsilon$, respectively, by the local bound. Consequently the limit measure is the point mass $m\delta_0$ plus a density, and the density is zero if $c_k\to-\infty$.

If $c_k$ has a finite limit, the positive limit density is $|z|^{2\vartheta-m/(2\pi)}$ times a factor bounded below by a positive constant, because the zero-boundary Dirichlet potential of the remaining measure is nonnegative. Its integrability gives
\[
 2\vartheta-m/(2\pi)>-2.
\]
This strict margin gives $L^p$ integrability of that power for some $p>1$. On a sufficiently small disk the remaining measure has arbitrarily small mass.

Jensen's inequality for the logarithmic potential gives a bound for its exponential in any preassigned finite $L^q$ on a smaller disk, by choosing the disk sufficiently small. H\"older's inequality with the strict margin gives the limit density in some $L^{p_1}$,
\[
 p_1>1,
\]
so its potential is continuous and bounded. Thus the limit density has the same power as an upper bound as well.

In both cases, including when the density vanishes, the regular part of the limiting gradient is $o(r^{-1})$ on $|z|=r$. Indeed for the power bound with exponent $\lambda>-2$ the integral of $|z-y|^{-1}|y|^\lambda$ over $|y|<2r$ is $O(r^{\lambda+1})$ by scaling.

The remaining part is bounded by a constant times $1+\int_{2r}^{r_0}s^\lambda\,ds$, whose product with $r$ tends to zero. The bounded harmonic part has no leading gradient.

\smallskip\noindent\textit{Step 2. The local mass.}\quad For fixed $k$, apply the divergence formula for $(z\cdot\nabla v_k)\nabla v_k-z|\nabla v_k|^2/2$ to a punctured disk. The inner circle can be removed first: $H_k=O_k(r^{2\vartheta_k})$ and the same direct Green gradient estimate gives $r|\nabla v_k|=o(1)$. Hence
\[
 \frac r2\int_{\partial D_r}
 [(\partial_r v_k)^2-(\partial_{\tan}v_k)^2]
 +r\int_{\partial D_r}H_k
 =\int_{D_r}(2+2\vartheta_k+z\cdot\nabla\log b_k)H_k.
\]
This follows by writing
\[
 H_k z\cdot\nabla v_k=z\cdot\nabla H_k- (2\vartheta_k+z\cdot\nabla\log b_k)H_k.
\]
First let $k$ tend to infinity on a fixed circle, then let $r\downarrow0$. The right side tends to $(2+2\vartheta)m$ since the logarithmic gradients of $b_k$ are uniformly bounded.

On the left the leading derivative is $-m/(2\pi r)$ and the remaining relative gradient tends to zero. The density term also vanishes, by the strict power bound or by zero residual density. Thus
\[
 m^2/(4\pi)=(2+2\vartheta)m.
\]
The positive mass $8\pi(1+\vartheta)$ is incompatible with
\[
 2\vartheta-m/(2\pi)>-2.
\]
Hence
\[
 c_k\to-\infty,\qquad v_k\to-\infty
\]
locally off zero inside the fixed circle. To propagate this convergence in the punctured disk, take any ball compactly contained there. Its Dirichlet potential $P_{k,B}$ of $H_k$ is uniformly bounded by the local bound on $H_k$.

For a uniform constant $C$ large enough, $C-v_k+P_{k,B}$ is positive harmonic on this ball. Harnack's inequality on smaller balls along a finite overlapping chain from the interior annulus gives
\[
 v_k\to-\infty
\]
locally throughout the punctured disk. The weight is uniformly bounded on such compact sets, so $H_k\to0$ there.
\end{proof}
\begin{lemma}[Compactness at scalar parameters]
\label{fnd:scalar-parameter}
Let the points of $S$ be distinct and let
\[
 h=h^*\prod_{p\in S}\exp(-4\pi(a_p-1)G(\cdot,p)),\qquad a_p>0.
\]
Suppose $\gamma\ge0$ avoids the positive values of
\[
 4\pi\left\{n+\sum_{p\in P}a_p:
                n\in\mathbb Z_{\geq0},\ P\subset S\right\},
\]
\begin{itemize}
\item[\textup{(i)}]
The mean-zero solutions of
\[
 -\Delta_gw=2\gamma\left(h e^w/\int_M h e^w-1\right)
\]
are bounded.
\item[\textup{(ii)}]
The bound is uniform near the limiting data as the shifted parameters $a_p$ and the scalar parameter $\gamma\ge0$ vary, the points vary in a compact set of distinct configurations, and the positive smooth factors vary in a compact $C^3$ family.
\end{itemize}
\end{lemma}
\begin{proof}
\smallskip\noindent\textit{Step 1. Concentration.}\quad Suppose there is an unbounded sequence with data tending to the given limit. Since $w=0$ when $\gamma=0$, it suffices to consider positive parameters. Write
\[
 g_k=\gamma_k h_ke^{w_k}/\int h_ke^{w_k}.
\]
The zero mean and Jensen's inequality give
\[
 \gamma_k\Big/\int h_ke^{w_k}\le C,\qquad \int g_k=\gamma_k\le C.
\]
Choose $p_0>1$ with $\inf_{p,k}p_0(a_{p,k}-1)>-1$ for the marked points. In local coordinates centered at a moving marked point, separation bounds the other factors, so
\[
 h_k^{p_0}\le C|z|^{2p_0(a_{p,k}-1)}.
\]
Choose $p,q>1$ with $pq\leq p_0$ and write $q'=q/(q-1)$. For a disk $D$ and a smaller disk $D'$, $w_k=2G*g_k$ gives $w_k\leq 2G*(\mathbf1_Dg_k)+C$ on $D'$. Fix $0<m_*<2\pi/(pq')$. If $0<m_k=\int_Dg_k<m_*$, Jensen's inequality for the probability $\mathbf1_Dg_k/m_k$ gives
\[
 \int_{D'}e^{pq'w_k}\,dV_g
 \le C\sup_{y\in D}\int_{D'}|z-y|^{-pq'm_k/\pi}\,dV_g(z)\le C.
\]
The bound also holds if $m_k=0$. H\"older's inequality with the uniform $L^{pq}$ weight bound and the bounded normalization then gives a uniform $L^p(D')$ estimate for $g_k$.

Choose disks whose boundaries have zero limiting mass. Weak convergence gives the preceding estimate away from the finitely many points at which the local mass reaches the threshold, including at marked points with no concentration.

The estimates are independent of small translations of the separated centers. The Green representation gives bounded $w_k$ on such smaller disks. Thus if the solution is unbounded there is a positive point mass.

\smallskip\noindent\textit{Step 2. The limiting masses.}\quad Isolate one blow-up point in a fixed coordinate disk. If a marked point converges to it, center the coordinates at the marked point and restrict to a smaller disk of fixed radius.

The disk contains at most one marked point. The translations tend to the identity and hence preserve the limiting local mass. Writing the conformal factor in these coordinates, the equation is
\[
 -\Delta w_k=2e^\varphi g_k-2\gamma_k e^\varphi.
\]
Add the Dirichlet solution of $\Delta\xi_k=-2\gamma_k e^\varphi$ and the constant $\log(2\gamma_k/\int h_ke^{w_k})$ to $w_k$. The resulting equation is the one in Lemma~\ref{fnd:scalar-atom}, with right-hand density $2e^\varphi g_k$, weight parameter in the translated coordinates
\[
 \vartheta_k=a_{p,k}-1
\]
(or zero at an unmarked point), and uniformly positive $C^1$ coefficient made up of the area factor, smooth weight and $e^{-\xi_k}$.

The fixed smaller circle is disjoint from the blow-up set, so $w_k$ is uniformly bounded there. The Dirichlet correction is uniformly bounded, and the added constant cancels from the boundary oscillation.

Its upper bound from Jensen, together with the local bounds for $w_k$ away from the blow-up point, gives the local upper bounds for the transformed field.

For continuity at its origin, keep a member fixed and shrink a disk until the mass of its $L^1$ density is below the threshold of the preceding Jensen estimate.

The local weight has the same strict power margin, so H\"older and the Green estimate give $W^{2,p}$ and hence continuity there for this member.

Only continuity for fixed $k$ is used at the blow-up point. The local mass in $g_k$ is therefore $4\pi a_p$ or $4\pi$. The lemma also gives convergence of the density to zero on a surrounding annulus.

Since $w_k$ is bounded there and the weight is bounded below there, the normalizing factor tends to zero. The estimates off the finite blow-up set then give zero residual mass everywhere. Hence the limiting total mass is a positive value of the displayed set, a contradiction.

\smallskip\noindent\textit{Step 3. Continuity of the homotopy.}\quad On bounded $C^0$ balls the normalized densities depend continuously in a common $L^p$ space, by convergence of the weights in that space and Jensen's lower bound on their integrals.

For moving centers the $L^p$ convergence follows by first removing a uniformly small ball around the limiting center (which has small $L^p$ mass by the strict power margin) and then using smooth convergence off that ball.

The Green $W^{2,p}$ estimate and compact embedding give the common compact map for degree homotopies.
\end{proof}
\begin{proof}[Proof of Proposition~\ref{fnd:signed-compactness}]
\noindent\textit{Step 1. Local bounds for (i).}\quad For negative source strengths the weights may be unbounded, so we keep the dependence on the source strengths explicit. Suppose, to the contrary, that there are solutions $u^k$ with
\[
 \alpha^k\longrightarrow\alpha^0,\qquad
 \rho^k\longrightarrow\rho^0,\qquad
 (h_1^{*,k},h_2^{*,k})\in\mathscr K,
 \qquad \|u^k\|_{C^0(M)^2}\longrightarrow\infty.
\]
After passing to a subsequence, the smooth factors converge in $C^3$. They have uniform positive lower bounds and uniform $C^3$ bounds. Choose $a_*>-1$ below every source strength under consideration. Set
\[
 I_i^k=\int_Mh_i^ke^{u_i^k}\,dV_g,\qquad
 c_i^k=\log\frac{\rho_i^k}{I_i^k},\qquad
 f_i^k=e^{c_i^k}h_i^ke^{u_i^k}.
\]
The mean-zero normalization and Jensen's inequality imply
\[
 I_i^k\geq
 \exp\left(\int_M\log h_i^k\,dV_g\right)
 =\exp\left(\int_M\log h_i^{*,k}\,dV_g\right)\geq c>0.
\]
Here the logarithms are integrable, including at the marked points, and the Green functions have zero mean. Thus $c_i^k\leq C$, whereas $\int_M f_i^k\,dV_g=\rho_i^k$. Passing to a further subsequence, write
\[
 f_i^k\,dV_g\rightharpoonup\nu_i,\qquad
 c_i^k\longrightarrow c_i\in[-\infty,C].
\]
We first obtain local bounds away from the blow-up points. Choose $p,q>1$ with $pq a_*>-1$ if $a_*<0$; when $a_*\geq0$, any fixed $p,q>1$ can be used. The weights are uniformly bounded in $L^{pq}(M)$. The Green function is bounded below, and
\[
                 u_i^k=\sum_{j=1}^2K_{ij}G*f_j^k.
\]
On a smaller concentric disk $D'$ inside a fixed disk $D$, this gives
\[
             u_i^k\leq 2G*(\mathbf1_Df_i^k)+C.
\]
The constant accounts for the part of $f_i^k$ outside $D$ and the lower bound for $G*f_j^k$, $j\ne i$. Put $q'=q/(q-1)$. The logarithmic kernel and Jensen's inequality give
\[
 \int_{D'}\exp\bigl(2pq'G*(\mathbf1_Df_i^k)\bigr)\,dV_g\leq C
\]
whenever $\int_Df_i^k\,dV_g$ is bounded above by a fixed number strictly smaller than $2\pi/(pq')$. H\"older's inequality,
\[
 c_i^k\leq C,
\]
and the uniform $L^{pq}$ bound for the weights then give an $L^p(D')$ bound for $f_i^k$. It follows that there is a finite set $\mathcal B$ outside which both densities have local uniform $L^p$ bounds.

Indeed, choose a fixed sufficiently small mass threshold and take the points at which one of the limiting measures has a point mass at least that large.

At every other point the preceding estimate applies on a small disk.

To bound both signs of $u_i^k$ on a still smaller disk, apply H\"older to the near kernels of both $G*f_j^k$ using their local $L^p$ bounds, and use the bounded separated kernel and the total masses for the rest. This gives the local $C^0$ bound also at marked points outside $\mathcal B$.

On compact subsets of $M\setminus(\mathcal B\cup S)$, the weights are smooth and bounded, and the same argument followed by elliptic regularity gives uniform $C^1$ bounds. In particular, every point of positive local mass belongs to $\mathcal B$.

If $\mathcal B$ were empty, a finite cover and the Green representation would bound $u^k$ on all of $M$, contrary to the choice of the sequence.

Thus at least one nonzero local mass occurs. After a subsequence the fields converge locally away from $\mathcal B\cup S$ by the interior estimates, and in $L^1(M)$ to
\[
                     u_i=\sum_jK_{ij}G*\nu_j.
\]
For the $L^1$ convergence, cut off $G(x,y)$ within distance $\delta$ of the diagonal. The cut-off kernels are continuous and their potentials converge in $L^1$ by weak convergence and bounded total masses. The discarded kernel has $\sup_y\int_{d_g(x,y)<2\delta}|G(x,y)|\,dV_g(x)=O(\delta^2(1+|\log\delta|))$.

\noindent\textit{Step 2. The local mass pairs.}\quad Fix a disk isolating one point of $\mathcal B$. At a marked point, take a conformal coordinate centered there and write
\[
 g=e^\varphi|dz|^2,\qquad
 h_i^k=|z|^{2\alpha_i^k}a_i^k(z),
\]
where $a_i^k$ are uniformly positive and bounded in $C^3$. At an unmarked point put $\alpha_i^k=0$ and use the same notation. Choose smooth functions $\xi_i^k$ on this fixed disk satisfying
\[
              \Delta\xi_i^k=(K\rho^k)_i e^\varphi,
\]
with zero boundary values. They are uniformly bounded in every fixed smooth norm. The functions
\[
 v_i^k=u_i^k+c_i^k+2\alpha_i^k\log|z|-\xi_i^k
\]
satisfy
\[
 \Delta v_i^k+\sum_jK_{ij}\widehat h_j^ke^{v_j^k}
       =4\pi\alpha_i^k\delta_0,
 \qquad \widehat h_i^k=e^\varphi a_i^ke^{\xi_i^k}.
\]
The coefficients $\widehat h_i^k$ are uniformly positive and bounded in $C^3$, and
\[
 \widehat h_i^ke^{v_i^k}=e^\varphi f_i^k.
\]
The masses are bounded. The boundary oscillations are bounded by the estimates on fixed annuli obtained above; the constants $c_i^k$ do not affect them. The origin is the only possible blow-up point in the disk.

Since the disk isolates a positive point mass, the regular parts $v_i^k-2\alpha_i^k\log|z|$ cannot both be bounded above: otherwise $\widehat h_i^k e^{v_i^k}$ would be uniformly bounded in $L^p$ by the integrability margin for the powers
\[
 |z|^{2\alpha_i^k},
\]
and neither limiting measure could have a point mass there. Use the canonical fields
\[
 U_i^k=v_i^k+\log\widehat h_i^k
\]
and densities
\[
 F_i^k=\widehat h_i^k e^{v_i^k}.
\]
Their forcings
\[
 b_i^k=\Delta\log\widehat h_i^k
\]
are uniformly bounded; the regular parts are continuous at zero and bounded in oscillation on the fixed boundary.

Thus Lemma~\ref{fnd:coarse-local} gives the marked local mass in $\Sigma_p(\alpha^0)$ in $2\pi$ units, while Lemma~\ref{fnd:local-unitary} gives integer masses at an unmarked point in $4\pi$ units. If both limiting source strengths vanish at a marked point, that same integer conclusion holds.

\noindent\textit{Step 3. Concentration in one component.}\quad If $c_i=-\infty$, the local bounds away from $\mathcal B$ give $f_i^k\to0$ there, in $L^1$ on compact subsets, so that $\nu_i$ is supported on $\mathcal B$.

Suppose instead that both $c_i$ are finite. Each $\nu_i$ then has the positive residual density $e^{c_i}h_i^0e^{u_i}$ away from $\mathcal B\cup S$. Consider any nonzero local mass pair $(s,t)$ in units of $4\pi$, and put
\[
 A=1+\alpha_1^0,\qquad B=1+\alpha_2^0
\]
at that point, with $A=B=1$ at an unmarked point. The local Pohozaev identity is
\[
                         s^2-st+t^2=As+Bt.
\]
The integrability of both residual densities implies
\[
                         2s-t\leq A,\qquad 2t-s\leq B.
\]
For completeness, these inequalities follow by averaging the logarithm of the limiting density on circles. The contribution of the point mass to the mean of $u_i$ is $-2(K(s,t))_i\log r$.

For each remaining measure, the circular-mean formula in the proof of Lemma~\ref{fnd:local-unitary} applies to the local logarithmic part of the Green kernel.

The measure has no point mass at the center, and the regular and separated parts of the kernel contribute $O(1)$. Thus its mean is $o(|\log r|)$ for each of the two measures. The circular mean of the logarithm of the residual density is
\[
 \bigl(2\alpha_i^0-2(K(s,t))_i\bigr)\log r
                         +o(|\log r|).
\]
If $(K(s,t))_i>1+\alpha_i^0$, Jensen's inequality on circles contradicts the integrability of this density. Multiplying the two inequalities by $s$ and $t$, respectively, and adding now gives
\[
 2(s^2-st+t^2)\leq As+Bt=s^2-st+t^2,
\]
which is impossible for a nonzero pair. Thus $c_i=-\infty$ for at least one component $i$. For that component all the mass is concentrated at finitely many points.

At the marked points the masses divided by $4\pi$ belong to $E_{i,p}(\alpha^0)$, with zero used when no mass is present; the sum of the masses at unmarked points is an integer in these units. Therefore
\[
                    \rho_i^0=\nu_i(M)\in\mathcal C_i(\alpha^0),
\]
contrary to the hypothesis. This proves the asserted uniform $C^0$ bound in a full neighborhood of $(\alpha^0,\rho^0)$. Only the local mass list at the limiting source strengths enters this argument.

\noindent\textit{Step 4. Compactness for (ii) and (iii).}\quad It remains to show that the fixed-point maps give a continuous compact homotopy in this neighborhood. Shrink the neighborhood so that every source strength in it is bounded below by
\[
 a_*>-1,
\]
and choose a fixed $p>1$ with
\[
 pa_*>-1
\]
if $a_*<0$. In local coordinates at each source, the functions $|z|^{2\alpha_{p,i}}$ vary continuously in $L^p$, by dominated convergence with the integrable majorant $1+|z|^{2pa_*}$ for their $p$th powers.

Thus the physical weights are uniformly bounded and jointly continuous in $L^p(M)$ as the source strengths and smooth factors vary. On a ball
\[
 \|u\|_{C^0(M)^2}\leq R,
\]
their normalized denominators satisfy
\[
 0<c_R\leq\int_Mh_ie^{u_i}\,dV_g\leq C_R.
\]
For example, the lower bound follows from $e^{-R}\int_Mh_i\,dV_g$ and the uniform Jensen lower bound used above. Hence the normalized densities depend continuously in $L^p$ on the source strengths, parameters, smooth factors, and $u$, and are uniformly bounded there on each such ball.

The Green operator maps $L^p(M)$ continuously to $W^{2,p}(M)$, whose embedding into $C^0(M)$ is compact because $p>1$. This proves continuity and compactness of the family. The uniform solution bound excludes zeros on the boundary of one common sufficiently large ball.

Finally, for a compact path satisfying the stated exclusions, apply the neighborhood assertion at each point of the path and take a finite subcover. The smooth factors along the path form a compact family in $C^3$.

The preceding $L^p$ argument applies with a single exponent, since the source strengths along the path are bounded away from $-1$. It gives the common compact homotopy and completes the proof.
\end{proof}
\begin{proof}[Proof of Lemma~\ref{fnd:sharp-compactness}]
For an unbounded sequence, pass to limiting parameters and smooth factors in the compact family. The local bounds and normalizations of Proposition~\ref{fnd:signed-compactness} give the hypotheses of Lemma~\ref{fnd:local-unitary} on each disk.

When $1,A_p,B_p$ are rationally independent, each component of the local mass pair, divided by $4\pi$, is one of the three values defining $\Gamma_i$. At a marked point with fixed integer source strengths, both coordinates are nonnegative integers, as at an unmarked point by the regular mass list.

The circular-average and Pohozaev argument of Proposition~\ref{fnd:signed-compactness} gives a fully concentrating component. Its limiting total mass is therefore in $\Gamma_i(\alpha)$: include all contributions at integer and unmarked points in the integer summand.

This contradicts the limiting parameter exclusion. The continuity and compactness of the maps follow from the same $L^p$ argument as in that proposition.
\end{proof}
\begin{proof}[Proof of Proposition~\ref{fnd:continuation}]
\smallskip\noindent\textit{Step 1. The perturbation of the source strengths in (i).}\quad Proposition~\ref{fnd:signed-compactness} gives the asserted uniform bound and compact homotopy in a full neighborhood of the prescribed nonnegative source strengths and positive parameters.

Choose a small product neighborhood with the preceding bound. Local finiteness permits the choice of $\widehat\rho_i/(4\pi)\in\mathbb Q\setminus\mathbb Z$ in this neighborhood, still outside $\mathcal C_i(\alpha)$. We may then increase the source strengths by arbitrarily small amounts to obtain $\widehat{\alpha}$ satisfying \eqref{fnd:independence}.

Indeed, the exceptional choices form a countable union of proper affine hyperplanes, so they cannot fill any open box of positive perturbations. This also treats a zero target source strength.

The parameter segment and the segment joining $\alpha$ to $\widehat\alpha$ can both be kept inside the neighborhood on which the solutions are uniformly bounded. The compact homotopy from Proposition~\ref{fnd:signed-compactness} therefore proves equality of the degrees at the target and at $(\widehat{\alpha},\widehat\rho)$.

We next keep $\widehat{\alpha}$ and $\widehat\rho$ fixed. The generic local list in Lemma~\ref{fnd:local-unitary} gives compactness away from $\Gamma_i(\widehat{\alpha})$ by Lemma~\ref{fnd:sharp-compactness}. The chosen parameters avoid these sets. In fact, an equality
\[
 \widehat\rho_i\in\Gamma_i(\widehat{\alpha})
\]
would express the rational noninteger $\widehat\rho_i/(4\pi)$ as an integer plus an integer linear combination of the shifted parameters $\widehat A_p,\widehat B_p$. If the coefficient of some $\widehat A_p$ or $\widehat B_p$ were nonzero, this would contradict \eqref{fnd:independence}; otherwise it would make $\widehat\rho_i/(4\pi)$ an integer.

Lemma~\ref{fnd:sharp-compactness} gives a common bound and compact homotopy along the compact family of positive smooth factors. Consequently their deformation preserves the degree. Here $\Gamma_i$ is evaluated at the fixed source strengths. No continuity of $\Sigma_p$ with respect to the source strengths is needed.

At a source the three possible values of the residual weight parameter $\beta_p$ are $\widehat B_p-1,\widehat A_p+\widehat B_p-1,\widehat A_p-1$ in the unused, pure, and mixed states. The stronger scalar exclusions used below are
\[
 \gamma\notin4\pi
       \bigcup_{E\subset S}\left(\mathbb Z+
                            \sum_{p\in E}\beta_p\right).
\]
Writing $\gamma=\widehat\rho_2-4\pi m_2$ with $m_2=V+\sum_{p\in\mathcal Q}\widehat B_p$, a forbidden equality for a subset $E$ would give
\[
 \frac{\widehat\rho_2}{4\pi}
   =n+V-|E|+\sum_{p\in\mathcal Q}\widehat B_p
                   +\sum_{p\in E}(\beta_p+1),
 \qquad n\in\mathbb Z.
\]
Each $\beta_p+1$ is an integer linear combination of $\widehat A_p,\widehat B_p$. A nonzero coefficient of some $\widehat A_p$ or $\widehat B_p$ contradicts \eqref{fnd:independence}; if all these coefficients vanish, the equality makes $\widehat\rho_2/(4\pi)$ an integer.

Hence these exclusions hold for every configuration with positive residual mass, in particular for the finite family retained below the cut-offs.

\smallskip\noindent\textit{Step 2. Fixed source strengths in (ii).}\quad Suppose now that the prescribed source strengths satisfy \eqref{fnd:independence}, or that they are all integers, and assume only
\[
 \rho_i\notin\Gamma_i(\alpha).
\]
At these fixed source strengths, Lemma~\ref{fnd:sharp-compactness} gives a common bound and compact homotopy along each compact parameter segment contained in the displayed component.

Since $\mathcal C_i(\alpha)$ is locally finite, choose a rational nonintegral representative in that component which also avoids both coarse sets.

With the source strengths fixed, the parameter homotopy preserves the degree. Keeping this parameter representative fixed, we can now perturb the source strengths as in Step~1.

We first apply Lemma~\ref{fnd:sharp-compactness} at fixed source strengths. Before varying the strengths, we choose a parameter pair outside both coarse sets and apply Proposition~\ref{fnd:signed-compactness}.

\smallskip\noindent\textit{Step 3. The finite sum in (iii).}\quad Finally, consider the expansion of $F^{\chi-|S|}\prod_{p\in S}P_{A_p,B_p}$ below fixed finite cut-offs. The expansion of the integer power $F^{\chi-|S|}$ has nonnegative integer exponents, while the product of the source factors has finitely many terms.

Choose a bounded rectangle containing all nearby cut-offs. All exponent shifts from the source factors are nonnegative for nearby source strengths, so the integer exponent pair of any contributing term of $F^{\chi-|S|}$ must lie in this same rectangle.

Thus a single finite list contains every term that contributes for at least one nearby pair $(\alpha,\rho)$.

Before coincident terms are combined, each exponent depends continuously and affinely on the source strengths. Each positive coordinate of such an exponent belongs to the corresponding $\Gamma_i(\alpha)/(4\pi)$, and the six initial pairs in the definition of $\Sigma_p$ give
\[
 \Gamma_i(\alpha)\subseteq\mathcal C_i(\alpha).
\]
For target data satisfying $\rho_i\notin\mathcal C_i(\alpha)$ for $i=1,2$, the parameter neighborhood and the perturbation of the source strengths can therefore be chosen so that no retained exponent crosses either cut-off.

For the sharper result at fixed source strengths, the initial parameter segment stays in one component of the complement of the $\Gamma_i$, and hence also preserves these strict inequalities.

Once the chosen representative also avoids $\mathcal C_i(\alpha)$ in both components, the same finite-term argument applies to the perturbation of the source strengths. Coincident limiting exponents are simply combined by adding their coefficients. This proves the last assertion.
\end{proof}
\subsection{Local masses in the scalar equation}
We first determine a restriction on the mass of a cluster of collapsing zeros. The argument allows a fixed weight parameter $\beta\ge0$, possibly noninteger, and arbitrary nonnegative integer multiplicities at the remaining zeros.
\begin{proposition}[Collapsing zeros at a fixed source]
\label{fnd:scalar-cosets}
Let $\beta\ge0$, $m_j\in\mathbb Z_{\ge0}$ and $p_{j,k}\to0$ in a fixed disk. Assume $H_k>0$ is $C^1$, $\|\nabla\log H_k\|_\infty\le C$, and
\[
 -\Delta v_k=2F_k,
 \qquad F_k=H_k|z|^{2\beta}\prod_j|z-p_{j,k}|^{2m_j}e^{v_k},
 \qquad \int F_k\leq C.
\]
For each $k$, $v_k$ is regular at the zeros, with no uniform bound required there. If $F_k\,dx\rightharpoonup\mu$ and $n=\mu(\{0\})/(4\pi)$, then
\begin{equation}\label{fnd:two-cosets}
                     n\in\mathbb Z\ \cup\ (\beta+\mathbb Z).
\end{equation}
The zeros may coincide or approach one another at arbitrary relative rates. The conclusion also holds with uniformly bounded smooth forcing.
\end{proposition}
\begin{proof}
\smallskip\noindent\textit{Step 1. The curvature integral.}\quad Put $h_k=\log H_k$, $U_k=\log F_k$, and
\[
 W_k=U_k-h_k=v_k+2\beta\log|z|
                     +2\sum_jm_j\log|z-p_{j,k}|.
\]
Set $\tau_j=-i\sigma_j/2$, where
\[
 \sigma_1=\begin{pmatrix}0&1\\1&0\end{pmatrix},\quad
 \sigma_2=\begin{pmatrix}0&-i\\i&0\end{pmatrix},\quad
 \sigma_3=\begin{pmatrix}1&0\\0&-1\end{pmatrix}.
\]
Then $[\tau_1,\tau_2]=\tau_3$ and its cyclic permutations hold. On the punctured disk define
\[
 \theta_1=\sqrt{F_k}\,dx,\quad
 \theta_2=\sqrt{F_k}\,dy,\quad
 \omega=\tfrac12(-W_{k,y}\,dx+W_{k,x}\,dy),
 \qquad\mathcal A_k=\theta_1\tau_1+\theta_2\tau_2+\omega\tau_3.
\]
Since $\Delta W_k+2F_k=0$ there, direct differentiation gives
\[
 d\mathcal A_k+\mathcal A_k\wedge\mathcal A_k
 =\tfrac12\sqrt{F_k}(-h_{k,y}\tau_1+h_{k,x}\tau_2)\,dx\wedge dy.
\]
Indeed the third component is $d\omega+\theta_1\wedge\theta_2=(\Delta W_k/2+F_k)dx\wedge dy$. The other components are $d\theta_1-\omega\wedge\theta_2$ and $d\theta_2+\omega\wedge\theta_1$. Consequently
\begin{equation}\label{fnd:curvature-small}
 \int_{D_r}|\mathcal F_k|
       \leq C\int_{D_r}\sqrt{F_k}
       \leq Cr\left(\int_DF_k\right)^{1/2}\leq Cr.
\end{equation}
Thus the curvature estimate depends only on the stated logarithmic gradient bound.

\smallskip\noindent\textit{Step 2. Holonomy at the sources.}\quad For fixed $k$, first shrink the circle about a source. At a source of total strength
\[
 \alpha\geq0,\qquad W_k=2\alpha\log|z-a|
\]
plus a regular function. The off-diagonal connection has loop integral tending to zero, and
\[
 \omega=\alpha d\vartheta
\]
up to an $L^1$-small loop term. The limiting source holonomy is therefore conjugate to $e^{-2\pi\alpha\tau_3}$, with trace $2\cos(\pi\alpha)$. At an integer source it equals the central matrix $(-1)^\alpha I$.

If $M=\sum_jm_j$, the trace of the product of all based source holonomies is thus $2(-1)^M\cos(\pi\beta)$, including coincidences at zero. There is only one possibly noncentral factor, so its conjugation does not change the trace.

For each fixed $k$ use disjoint holes about the distinct zero locations inside the enclosing circle. The connection is $C^1$ on the resulting closed annular region with holes, and the based transports are unitary.

Apply Lemma~\ref{fnd:holonomy}, then shrink each hole, keeping that circle fixed. Since the absolute trace of a matrix in dimension two is bounded by twice its operator norm, \eqref{fnd:curvature-small} implies
\begin{equation}\label{fnd:outer-trace}
 |\operatorname{tr}H_k(r)-2(-1)^M\cos(\pi\beta)|\leq Cr.
\end{equation}
Choose $d_l\downarrow0$ such that the two boundary circles have zero $\mu$-mass and $\mu(D_{2d_l}\setminus D_{d_l})\to0$. A diagonal subsequence $k_l$ can be chosen with
\[
 \max_j|p_{j,k_l}|\leq d_l/l,
\]
vanishing annular mass, and both disk masses tending to $4\pi n$. Fubini supplies $r_l\in(d_l,2d_l)$ with
\[
 r_l\int_{\partial D_{r_l}}F_{k_l}\,ds\to0,
 \qquad\int_{D_{r_l}}F_{k_l}\to4\pi n.
\]
The entire cluster lies well inside these circles. Moreover $\int_{\partial D_{r_l}}\sqrt{F_{k_l}}\,ds\to0$ by Cauchy--Schwarz. In a frame transported by $\omega\tau_3$, the transport equation and unitarity show that removing the off-diagonal part changes holonomy by $o(1)$, independently of the size of $\omega$.

\smallskip\noindent\textit{Step 3. The flux on the outer circle.}\quad To compute the remaining flux, use the distributional identity
\[
 \Delta W_k=4\pi\beta\delta_0+
             4\pi\sum_jm_j\delta_{p_{j,k}}-2F_k,
 \qquad
 \int_{\partial D_{r_l}}\omega
      =2\pi(\beta+M)-\int_{D_{r_l}}F_{k_l}
      \longrightarrow2\pi(\beta+M-2n).
\]
Comparison with \eqref{fnd:outer-trace} yields
\[
 \cos(\pi(\beta+M-2n))=\cos(\pi(\beta+M)).
\]
Equality of these cosines gives $n\in\mathbb Z$ or $n\in\beta+M+\mathbb Z$, proving \eqref{fnd:two-cosets}. A bounded smooth forcing $b_k$ in $\Delta W_k+2F_k=b_k$ adds $\tfrac12\int_{D_r}b_k=O(r^2)$ to the flux and the same order to the curvature bound.
\end{proof}
\begin{remark}
The $SO(3)$ holonomy alone does not distinguish the two parities. For two noninteger sources the holonomies at a common base point need not commute, so the argument does not give the same two cosets.
\end{remark}
\subsection{Blow-up configurations of the shadow system}
\begin{lemma}[Compactness for separated scalar weights]
\label{fnd:separated-scalar-compactness}
Fix $\beta_p>-1$ and $\gamma>0$ satisfying \eqref{fnd:scalar-exclusion}. Suppose the ordinary centers vary in a compact set of distinct configurations in $M\setminus S$.
\begin{itemize}
\item[\textup{(i)}]
The scalar solutions are uniformly bounded.
\item[\textup{(ii)}]
This bound remains uniform for positive smooth factors converging in $C^1$ with a common positive lower bound.
\end{itemize}
\end{lemma}
\begin{proof}
Let the ordinary centers vary in the given compact set and let the positive smooth factors converge in $C^1$ with a uniform positive lower bound. The physical weights are uniformly bounded in some $L^p$, $p>1$, and their logarithms have uniformly bounded integrals.

Jensen's inequality bounds the normalizing factors from above, and the local logarithmic-potential argument in Proposition~\ref{fnd:signed-compactness} gives a finite concentration set and local bounds for $w$ off that set.

A disk isolating a blow-up point contains at most one prescribed or moving zero. In a conformal coordinate, translate that zero to the origin when it moves. Write the local equation as
\[
 -\Delta w_k=2F_k-b_k,
\]
with
\[
 F_k=a_k |z|^{2\vartheta}A_k e^{w_k},\qquad a_k=\gamma/\int_M H_{P_k}e^{w_k}.
\]
The area factor is included in $A_k$, and $\vartheta$ is $0$, $1$, or a fixed $\beta_p$. Solve $\Delta\xi_k=-b_k$ on a fixed surrounding disk with zero boundary data. Then $v_k=w_k+\xi_k+\log(2a_k)$ solves
\[
 -\Delta v_k=|z|^{2\vartheta}A_k e^{-\xi_k} e^{v_k}=2F_k.
\]
Its positive coefficient has uniform upper and lower bounds and a uniform first derivative bound. The constant shift does not change boundary oscillations, which are bounded by the local estimates away from the blow-up point. The positive local mass implies
\[
 \sup v_k\to\infty:
\]
an upper bound would bound $F_k$ in a fixed $L^p$, $p>1$, by
\[
 \vartheta>-1.
\]
Away from this only possible blow-up point, $v_k$ is bounded above, since $a_k$ is bounded above. Thus Lemma~\ref{fnd:scalar-atom}, with right-hand side
\[
 H_k=2F_k,
\]
gives local mass $4\pi(1+\vartheta)$ for $f$ and convergence of the density to zero on a surrounding annulus. The latter forces the normalizing factor to tend to zero, since both the weight and $e^w$ are bounded below there.

The same exterior bounds then give zero residual mass everywhere. Summing these local masses contradicts \eqref{fnd:scalar-exclusion}. In particular, after the centers have been separated no scalar concentration can occur in any of the deformations below. The proof applies to any fixed
\[
 \beta_p>-1.
\]
\end{proof}
\begin{lemma}[Exclusion of collisions at prescribed singularities]
\label{fnd:source-exclusion}
Fix a configuration with residual scalar parameter $\gamma>0$. Under \eqref{fnd:independence} and \eqref{fnd:scalar-exclusion}, every sequence of solutions of the corresponding shadow system satisfies the following.
\begin{itemize}
\item[\textup{(i)}]
No free center approaches $S$.
\item[\textup{(ii)}]
A nontrivial collision at a regular point involves precisely two ordinary centers and scalar local mass $4\pi$.
\end{itemize}
\end{lemma}
\begin{proof}
\noindent\textit{Step 1. The scalar measure.}\quad Write
\[
 f_k=a_k H_{P_k}e^{w_k},
\]
with $H_{P_k}$ the physical scalar weight and
\[
 a_k=\gamma/\int_M H_{P_k}e^{w_k}.
\]
Jensen's inequality and the zero mean of all the Green factors bound $a_k$ above.

Choose $p>1$ for which the fixed source powers belong locally to $L^p$. The factors at the moving ordinary zeros are bounded above, so
\[
 \sup_k\|H_{P_k}\|_{L^p}<\infty.
\]
The local Jensen estimate of Proposition~\ref{fnd:signed-compactness}, applied to
\[
 w_k=2G*f_k,
\]
gives a finite concentration set and local uniform bounds for $w_k$ outside it. If $a_k\to0$, these bounds give $f_k\to0$ in $L^1$ on compact subsets of the complement, so the entire scalar measure is a sum of point masses.

Proposition~\ref{fnd:scalar-cosets} applies on a disk isolating each of these blow-up points and the whole cluster of centers converging to it.

The ordinary centers give integer zeros, a fixed source contributes the factor $\lvert z\rvert^{2\beta}$, and all remote fields enter a positive coefficient with bounded logarithmic gradient.

The bounded smooth forcing is permitted by the last assertion of that proposition. The sum of the local masses, divided by $4\pi$, would then be an integer plus the sum of $\beta_p$ over a subset of $S$, contrary to \eqref{fnd:scalar-exclusion}.

Thus every subsequential limit of $a_k$ is positive, and the scalar limit has positive residual mass.

\noindent\textit{Step 2. Local equations at a collision.}\quad Consider all $u$ ordinary and $v$ mixed centers converging to one point. Write $4\pi n$ for the scalar local mass there. At a fixed source set
\[
 a=1+\alpha,\qquad b=1+\beta;
\]
at a regular point take $a=b=1$. Choose a conformal coordinate in which the cluster limit is $0$, and write
\[
 F_k=f_ke^\varphi.
\]
Multiply each balance equation by its positive conformal metric factor, so that gradients in the following moments are Euclidean. On a disk containing no other full cluster the scalar equation and coefficient have the form
\[
 -\Delta w_k=2F_k-\mathfrak b_k,\qquad
 F_k=V_ke^{w_k},\qquad
 V_k=H_k|z|^{2\beta}\prod_{i=1}^u|z-p_{i,k}|^2,
 \qquad \mathfrak b_k=2\gamma e^\varphi,
\]
where the logarithmic gradients of $H_k$ are uniformly bounded; its spatially constant normalization factor does not affect these gradients. Positivity of the residual mass gives, after passage to a subsequence, a positive limiting normalization and
\[
 w_\infty=-4n\log|z|+\psi,\qquad
 F_0=H_0|z|^{\lambda}e^\psi,
 \qquad \lambda=2(\beta+u)-4n,
\]
with $H_0$ positive and with bounded logarithmic gradient. The function $\psi$ is the Green potential of the remaining measure plus a smooth function, and is bounded below. Integrability of $F_0$ therefore implies
\[
 \lambda>-2,
\]
equivalently
\[
 2n<b+u.
\]
The remaining measure has no point mass at $0$.

On a sufficiently small disk its Green potential consequently satisfies $e^\psi\in L^s$ for every fixed finite $s$: Jensen's inequality reduces this assertion to the integrability of a logarithmic kernel whose exponent is the arbitrarily small mass of that disk.

Choose $s$ sufficiently large and use $\lambda>-2$ to obtain $F_0\in L^p$ for some $p>1$. Local elliptic regularity makes $\psi$ bounded. Hence
\[
 F_0=O(|z|^\lambda),
\]
and the Green gradient estimate gives $|z|\,|\nabla\psi(z)|\to0$. These facts justify the shrinking-circle limits in the radial and angular Pohozaev identities. For a fixed enclosing radius $r$, define
\[
 I_{i,k}(r)=\int_{D_r}F_k(x)
                    \frac{x-p_{i,k}}{|x-p_{i,k}|^2}\,dx,
 \qquad
 \nabla w_k(p_{i,k})=\pi^{-1}I_{i,k}(r)+Q_{i,k}(r).
\]
Green representation bounds $Q_{i,k}(r)$ uniformly in $k$ for each fixed $r$: its terms are the smooth part of the local Green kernel and the field outside $D_r$. In particular,
\[
 p_{i,k}\cdot Q_{i,k}(r)\to0.
\]
For fixed $k$ the local Green estimate and $F_k=O_k(|z|^{2\beta})$ give
\[
 \epsilon|\nabla w_k|=o(1),\qquad \epsilon^2F_k=o(1)
\]
on $|z|=\epsilon$. The radial and angular derivatives of $2\beta\log|z|$ are respectively $2\beta$ and zero under $x\cdot\nabla$ and the positive angular derivative, so the corresponding volume terms are integrable and the inner boundary vanishes.

At the moving ordinary zeros the coefficient vanishes quadratically and the volume terms are locally integrable. The radial identity is
\begin{align*}
 &\int_{\partial D_r}\left(
  \partial_\nu w_k\,x\cdot\nabla w_k
  -\tfrac12|\nabla w_k|^2x\cdot\nu+2F_kx\cdot\nu\right)\\
 &\hspace{15mm}=\int_{D_r}\left(
  4F_k+2F_kx\cdot\nabla\log V_k
                   +\mathfrak b_kx\cdot\nabla w_k\right).
\end{align*}
Using
\[
 x\cdot\nabla\log V_k
 =2\beta+2u+
  2\sum_i p_{i,k}\cdot\frac{x-p_{i,k}}{|x-p_{i,k}|^2}
  +O(|x|)
\]
first expresses its moving-zero terms as $4\pi\sum_i p_{i,k}\cdot\nabla w_k(p_{i,k})+o_k(1)$. On each fixed enclosing circle avoiding the limiting sources, centers and blow-up points, the local bounds and elliptic estimates give $C^1$ convergence of $w_k$ and convergence of $F_k$.

Take $k\to\infty$ there, and then $r\downarrow0$. The boundary limit is $16\pi n^2$ and the remaining mass term is $16\pi(b+u)n$. The smooth-coefficient error is $O(r)$. The background term tends to zero because the global Green representation bounds $\nabla w_k$ in every fixed $L^q$,
\[
 1<q<2.
\]

The Pohozaev identity for the rotation field $Jx$, where $J$ is rotation through $+\pi/2$, has no divergence term and no contraction with its skew derivative.

The fixed source logarithm has zero angular derivative, whereas the moving-zero terms are $4\pi\sum_i Jp_{i,k}\cdot\nabla w_k(p_{i,k})+o_k(1)$. Its boundary limit is zero by the preceding control of $\psi$. Thus the radial and angular moments of the cluster are
\begin{equation}\label{fnd:scalar-cluster-moments}
 \sum_i p_{i,k}\cdot\nabla w_k(p_{i,k})
       \longrightarrow4n(n-b-u),\qquad
 \sum_i Jp_{i,k}\cdot\nabla w_k(p_{i,k})\longrightarrow0.
\end{equation}
Since the moments are taken over the entire cluster, these identities hold for arbitrary relative distances within the cluster.

\noindent\textit{Step 3. The possible collisions.}\quad Take the inner product of each center equation with the position vector of that center and sum over the cluster. The source terms contribute $2(a-1)u+2(2a+b-3)v$ and the pair terms contribute $-2u(u-1)-6uv-6v(v-1)$.

For each pair the two singular forces are opposite, so their combined radial moment is independent of the distance between the centers. Their angular moment vanishes. With the original coefficient $-1/2$ the result is
\begin{equation}\label{fnd:virial}
 m^2-mt+t^2=am+bt,
 \qquad m=u+2v,\quad t=v+n.
\end{equation}
The same estimates control the errors when the smooth factors converge in $C^1$ with uniform positive lower bounds. The estimates use bounded logarithmic gradients of these factors.

\smallskip\noindent\emph{Part (i).} By \eqref{fnd:two-cosets}, $t$ is either integral or $b$ plus an integer. Substituting either form and any of the three pairs \eqref{fnd:effective} into \eqref{fnd:virial} gives a rational linear relation among $1,A_p,B_p$. Its coefficients force $m=0$. For example, when
\[
 (a,b)=(A,B),\qquad t=k
\]
the coefficients of $A,B$ are $-m,-k$; when $t=B+k$ they are $-m,k-m$, again forcing $m=0$. For
\[
 (a,b)=(-A,A+B)
\]
the integral case has coefficients $m-k,-k$, and the shifted case has coefficients $k,k-m$. For $(a,b)=(-(A+B),A)$ these are $m-k,m$ and $k,m$, respectively. Thus $u=v=0$ in every case. With no centers the Pohozaev identity gives $n=0$ or $n=b$, and the latter contradicts $2n<b$.

\smallskip\noindent\emph{Part (ii).} At a regular point $a=b=1$ and $n$ is integral. The nonzero nonnegative integer solutions of \eqref{fnd:virial} are $(m,t)=(1,0),(0,1),(1,2),(2,1),(2,2)$. Together with
\[
 m=u+2v,\qquad t=v+n,
\]
and $2n<u+1$, the only class with at least two centers is $(u,v,n)=(2,0,1)$. If there is at most one center, $u\leq1$ and the same strict bound gives $n=0$.
\end{proof}
\subsection{Collision of two ordinary blow-up points}
\label{fnd:pair-section}
We next consider the collision of two ordinary blow-up points. The first lemma identifies the limiting shadow system under $C^1$ convergence of the smooth coefficients.

We then derive a second-order identity under stronger local regularity. Throughout this subsection
\[
 g=e^\varphi|dz|^2,\qquad u_z=(u_x-iu_y)/2.
\]
The rotation $J$ through $+\pi/2$ acts by $-i$ on this complex encoding of a gradient.
\begin{lemma}[The limiting shadow system at a collision of two points]
\label{fnd:pair-closure}
Suppose the positive smooth factors satisfy
\[
 h_{i,k}^*\to h_i^*\text{ in } C^1,\qquad
 \inf_{i,k}\min_M h_{i,k}^*>0,\qquad \sup_{i,k}\|h_{i,k}^*\|_{C^2}<\infty.
\]
Along a sequence of solutions of the shadow system, group together all centers with the same limit. Assume that there are finitely many collision points, all regular, and that each corresponding group satisfies
\[
 (u,v,n)=(2,0,1).
\]
If the residual scalar mass is positive, replacing each pair by a mixed center gives a solution of the shadow system with fewer ordinary centers. This also holds for
\[
 (I+\eta_kJ)F_{i,0}^A+c_k\nabla w_k(p_i)=0,
 \qquad (I+\eta_kJ)F_j^B=0,
 \qquad c_k\to-\tfrac12,\quad\eta_k\to0.
\]
\end{lemma}
\begin{proof}
\smallskip\noindent\textit{Step 1. The limiting density.}\quad If the pair limits are $q_1,\ldots,q_l$, scalar normalization and Green representation give
\[
 f_kdV_g\rightharpoonup4\pi\sum_a\delta_{q_a}+f_0dV_g,
 \qquad w_k\longrightarrow8\pi\sum_aG(\cdot,q_a)+z
\]
off those points. Each pair gives two quadratic zeros, whose limiting product cancels the singular factor $|z-q_a|^{-4}$ from the Green potential of the point mass.

Hence $z$ solves the scalar equation with parameter $\gamma-4\pi l>0$ and only the surviving ordinary zeros. Near each new mixed center the limiting scalar factor is positive and $C^1$.

The remaining measure has no point mass there, so the local logarithmic Jensen estimate gives $e^z$ in some $L^p$, $p>1$. The Poisson estimate first makes $z$ continuous and bounded.

With the now bounded density it gives $C^{1,\delta}$ regularity for $\delta<1$. The scalar density is locally Lipschitz because the factor is $C^1$, and Schauder estimates then give $C^{2,\theta}$ regularity for every
\[
 \theta<1.
\]
If the limiting factors are smooth, the equation also makes $z$ smooth near the new mixed centers.

\smallskip\noindent\textit{Step 2. The local equations.}\quad In a disk centered at one limiting point $q=0$, include the contributions of all exterior centers and prescribed sources in the positive weights $\widetilde h_i$. Their logarithms converge in $C^1$. Write
\[
 G(x,p)=-\frac1{2\pi}\log|x-p|+R(x,p),\quad
 \Gamma_1=\log\widetilde h_1+4\pi r_g,
 \quad R_x=\nabla_xR(0,0).
\]
Locally $w_k\to8\pi G(\cdot,0)+v$. Use Euclidean gradients in this coordinate disk, and put
\[
 W=\nabla v(0).
\]
Before summing, multiply each intrinsic balance equation by $e^{\varphi(p_i)}$, including its scalar term. The factor commutes with the rotation. In the resulting equations the two singular pair forces are equal and opposite, so their sum gives
\begin{equation}\label{fnd:pair-sum}
 S_k:=\nabla w_k(p_+)+\nabla w_k(p_-)
            \longrightarrow4\nabla\Gamma_1(0)+32\pi R_x.
\end{equation}
For the perturbed equations the sum is $c_k S_k$ plus the rotated smooth force, whose limit is $2\nabla\Gamma_1(0)+16\pi R_x$. Division by $c_k\to-1/2$ gives the same limit in \eqref{fnd:pair-sum}.

\smallskip\noindent\textit{Step 3. The center equations.}\quad Put $F_k=f_ke^\varphi$ and $V_k=F_ke^{-w_k}$. The translation Pohozaev identity for the scalar equation gives
\[
 \int_{\partial D_r}\left(\partial_\nu w_k\nabla w_k
       -\tfrac12|\nabla w_k|^2\nu+2F_k\nu\right)
 =\int_{D_r}\left(2F_k\nabla\log V_k
                         +2\gamma e^\varphi\nabla w_k\right).
\]
For $i=+,-$, define
\[
 I_{i,k}(r)=\int_{D_r}F_k(x)\frac{x-p_{i,k}}{|x-p_{i,k}|^2}\,dx,
 \quad Q_{i,k}(r)=\nabla w_k(p_{i,k})-\pi^{-1}I_{i,k}(r).
\]
Taking first $k\to\infty$ and then $r\downarrow0$ gives $Q_{i,k}(r)\to8\pi R_x+W$. Moreover
\[
 \nabla\log V_k=\nabla(\log\widetilde h_{2,k}+\varphi)
       -4\pi\sum_{i=\pm}\nabla_xR(x,p_{i,k})
       +2\sum_{i=\pm}\frac{x-p_{i,k}}{|x-p_{i,k}|^2}.
\]
The boundary limit of the translation identity is $-8\pi(8\pi R_x+W)$, and its volume limit is
\[
 8\pi\nabla(\log\widetilde h_2+\varphi)(0)
       -128\pi^2R_x+4\pi\lim S_k-8\pi W.
\]
The smooth coefficient term uses only $C^1$ convergence. The bounded background forcing vanishes on shrinking disks because Green representation bounds $\nabla w_k$ in every fixed $L^q$, $1<q<2$. Thus
\[
 \lim S_k=-2\nabla(\log\widetilde h_2+\varphi)(0)+16\pi R_x.
\]
Combining this with \eqref{fnd:pair-sum} and $8\pi R_x=4\pi\nabla r_g-\nabla\varphi$ gives
\[
       \nabla(2\log\widetilde h_1+
                          \log\widetilde h_2+12\pi r_g)(0)=0.
\]
This is the new mixed balance equation. At a surviving ordinary center the pair contributes $16\pi G$ and its scalar point mass contributes $-4\pi G$, leaving $12\pi G$.

At an existing mixed center it contributes $24\pi G$. Expanding the absorbed fields therefore gives every equation of the limiting shadow system.
\end{proof}
\begin{lemma}[Location of scalar blow-up near two colliding points]
\label{fnd:midpoint}
Let the sequence of global shadow solutions be as above. On a fixed disk away from the sources and the other full collisions, assume the weights including the exterior contributions converge in $C^3$. Write
\[
 p_\pm=s\pm d,\qquad \delta=|d|,\qquad d/\delta\to e,\quad |e|=1.
\]
Suppose the balance equations are
\[
 (I+\eta_kJ)F_{i,0}^A+c_k\nabla w(p_i)=0,\qquad c_k\to-1/2,\quad\eta_k\to0.
\]
For a scalar concentration of mass $4\pi$, the rescaled measures in $z=s+\delta y$ converge to $4\pi\delta_0$.
\end{lemma}
\begin{proof}
\smallskip\noindent\textit{Step 1. The local mass.}\quad Subtract a smooth solution of the bounded background equation. The resulting local equation is
\[
 -\Delta v_k=2A_k(z)|z-p_{+,k}|^2|z-p_{-,k}|^2e^{v_k}.
\]
Set $\delta_k=|d_k|$, $e_k=d_k/\delta_k$, and
\[
 V_k(y)=v_k(s_k+\delta_k y)+6\log\delta_k+\log A_k(s_k),
 \qquad a_k(y)=\frac{A_k(s_k+\delta_k y)}{A_k(s_k)}.
\]
Then
\begin{equation}\label{fnd:pair-rescaled-equation}
 -\Delta V_k=2\widehat F_k,\qquad
 \widehat F_k=a_k(y)|y-e_k|^2|y+e_k|^2e^{V_k}.
\end{equation}
On every fixed compact set, $a_k\to1$ with its first derivatives. The change of variables preserves mass. Hence, for every compact $K\subset\mathbb R^2$,
\[
 \limsup_{k\to\infty}\int_K\widehat F_k\leq4\pi.
\]

\smallskip\noindent\textit{Step 2. Rescaling by the distance between the centers.}\quad To locate the scalar blow-up points at scale $\delta$, restrict the Green representation to a fixed disk about the pair. Extend the rescaled density by zero and absorb the constants into $C_k$. On compact sets in $y$,
\[
 V_k(y)= C_k-\frac1\pi\int
       \widehat F_k(\xi)\log\frac{|y-\xi|}{1+|\xi|}\,dA_\xi
       +H_k(y),
\]
where $H_k$ has uniformly bounded oscillation and can be normalized at a fixed point. Indeed, the smooth part of the local Green kernel, the exterior mass and the bounded background have uniformly bounded gradients on a smaller physical disk.

Their rescaled oscillation on each fixed compact set is $O(\delta_k)$. The logarithm is bounded above there uniformly in $\xi$. It follows that
\[
 V_k\ge C_k-C
\]
on a fixed ball avoiding the two zeros. The coefficient in \eqref{fnd:pair-rescaled-equation} is bounded below on that ball and the mass is bounded, so
\[
 C_k\le C.
\]
Near a point at which the weak limit of the rescaled measures has no point mass of size greater than a fixed small threshold, the local logarithmic Jensen estimate gives $\widehat F_k$ uniformly in some $L^p$, $p>1$, on a smaller disk.

The logarithmic representation then bounds $V_k-C_k$ there both above and below. Thus there are only finitely many possible blow-up points on the plane. Off them, after a subsequence, either $C_k$ has a finite limit and $V_k$ converges locally, or
\[
 C_k\to-\infty,\qquad V_k\to-\infty
\]
locally uniformly.

On a small disk about each possible blow-up point, we apply Lemma~\ref{fnd:scalar-atom}. To check its assumptions at the positive zero, translate the actual point $e_k$ to the origin:
\[
 U_k(x)=V_k(e_k+x),\qquad
 -\Delta U_k=|x|^2W_k(x)e^{U_k},\qquad
 W_k(x)=2a_k(e_k+x)|x+2e_k|^2.
\]
For fixed sufficiently small $r<1/2$, the other zero lies outside $B_r$, and
\[
 0<c\leq W_k\leq C,\qquad \|\nabla W_k\|_{L^\infty(B_r)}\leq C.
\]
The same assertion holds at the negative zero. Away from the two zeros, the corresponding coefficient is uniformly positive and has bounded first derivatives.

The global Green representation also bounds the boundary oscillation in Lemma~\ref{fnd:scalar-atom}. After the usual local selection, choose a small circle $\Gamma$ and a surrounding annulus disjoint from the finite limiting blow-up set.

On that annulus $\widehat F_k$ is uniformly bounded. For $y,y'\in\Gamma$, subtracting the Green representations gives
\[
 V_k(y)-V_k(y')
 =-\frac1\pi\int\widehat F_k(\xi)
       \log\frac{|y-\xi|}{|y'-\xi|}\,d\xi
       +H_k(y)-H_k(y'),
\]
where the smooth term has uniformly bounded oscillation. The logarithmic singularities are uniformly integrable on the annulus. Off the annulus the kernel is bounded on compact sets and is $O_\Gamma(|\xi|^{-1})$ at infinity. The uniform total mass bound therefore gives
\[
 \operatorname{osc}_{\Gamma}V_k\leq C.
\]
The disk may be chosen to contain only the prescribed limiting blow-up point, so all the hypotheses of Lemma~\ref{fnd:scalar-atom} hold. Its mass is $8\pi(1+\alpha)$ for the right-hand side of the equation. In our normalization this gives $\widehat F_k$-mass $4\pi$ at a regular point and $8\pi$ at a quadratic zero of the scalar weight.

\smallskip\noindent\textit{Step 3. The position of the blow-up point.}\quad There must be a finite blow-up point on the pair scale. If there were none, either the rescaled potential has an entire limit $V$, or it tends locally to $-\infty$.

In the first case $-\Delta V=2|y-e|^2|y+e|^2e^V$ has total mass at most $4\pi$. Circular averaging gives
\[
 \overline V'(r)=-M(r)/(\pi r),
\]
so $\overline V(r)\geq-4\log r+C$. For large $r$ the polynomial weight is at least $cr^4$, and Jensen's inequality contradicts finite area. In the second case, Green representation gives
\[
 \nabla V_k(y)=-\frac1\pi\int\widehat F_k(\xi)
                   \frac{y-\xi}{|y-\xi|^2}\,d\xi+o(1).
\]
The exterior harmonic gradient is $O(\delta_k)$ after scaling. With no finite blow-up point, the local estimate for the density gives $\widehat F_k=o(1)$ in local $L^q$, $q>2$, in the branch
\[
 C_k\to-\infty.
\]
Thus the gradient integral tends to zero at the moving points $\pm e_k$ (the tail outside $D_R$ is $O(1/R)$), and the positive ordinary balance would be $-2e=0$.

Neither zero can blow up because $8\pi$ exceeds the available mass. There is therefore exactly one regular blow-up point $y_0$ of mass $4\pi$. The Green formula and local elliptic estimates near the two zeros give
\[
 \nabla V_k(y)\to-4(y-y_0)/|y-y_0|^2
\]
there. The balance equation at the positive zero becomes $-2e+2(e-y_0)/|e-y_0|^2=0$. Inversion is injective and $|e|=1$, so $y_0=0$. This already accounts for the entire original local mass.
\end{proof}
\begin{proposition}[A second-order identity for two colliding points]
\label{fnd:pair-identity}
Assume the hypotheses on the pair and coefficients in Lemma~\ref{fnd:midpoint}. At the limiting mixed solution use $\kappa,\mathcal T$ in \eqref{fnd:gap} and set
\[
 \zeta=d^2,\qquad \varphi(q)=0.
\]
\begin{itemize}
\item[\textup{(i)}]
For the original pair,
\begin{equation}\label{fnd:physical-pair}
              -2\mathcal T(q)\zeta-\tfrac12\kappa(q)|\zeta|
                          =o(|\zeta|).
\end{equation}
\item[\textup{(ii)}]
For $(I+\eta J)F_{i,0}^A+c\nabla w(p_i)=0$ with $c=-1/2$, $\eta\to0$, the right-hand side is
\begin{equation}\label{fnd:rotated-pair}
        2i\eta+2\eta^2+O(|\eta||\zeta|)+o(|\zeta|).
\end{equation}
\item[\textup{(iii)}]
Both estimates are uniform on compact families of coefficients and solutions of the limiting shadow system. Thus $\kappa>4|\mathcal T|$ excludes pair collisions for the original system and all sufficiently small rotations.
\end{itemize}
\end{proposition}
\begin{proof}
\noindent\textit{Step 1. The equation near the pair.}\quad Absorb exterior fields as in Lemma~\ref{fnd:pair-closure}, and put $\Gamma_1=\log\widetilde h_1+4\pi r_g$, $t=1-i\eta$. The exact balance equation is
\[
 t\{\Gamma_{1,z}(p_i)+8\pi G_z(p_i,p_j)\}+cw_z(p_i)=0.
\]
For $F=fe^\varphi$ and $U=\log F$,
\[
 \Delta U+2e^U=\mathcal K+4\pi(\delta_{p_+}+\delta_{p_-}),
 \quad\mathcal K=\Delta(\log\widetilde h_2+\varphi)
                                  +(2\gamma-8\pi)e^\varphi.
\]
Set $\mathscr T=U_{zz}-U_z^2/2$. Then $\partial_{\bar z}\mathscr T=(\mathcal K_z-\mathcal K U_z)/4$. On nested fixed disks write
\[
 U_z=B_k(z)+\sum_{i=\pm}\frac1{z-p_i}
                    -\frac1{2\pi}\int\frac{F_k(y)}{z-y}\,dA_y,
\]
where the density integral is over a fixed larger disk whose boundary has zero limiting mass.

The remainder $B_k$ contains $(\log\widetilde h_{2,k}+\varphi)_z$, the regular derivatives of the two source Green functions, and the derivatives of the regular local kernel and the exterior Green potential against the density.

For the last two terms, the kernels and their first $z$ derivatives are uniformly smooth on the smaller disk, with separation from the exterior integration set. Weak convergence of the bounded measures is uniform on each such compact family of kernels.

Together with local convergence of the smooth factors and the source positions this proves convergence of $B_k$ in $C^1$ on the smaller disk. Take its fixed boundary away from the limiting collision point and any other blow-up point. Let
\[
 Ca(z)=\pi^{-1}\int a(\xi)/(z-\xi)\,dA_\xi
\]
and define
\[
 E_{\mathcal K}(z,y)=C\!\left[\frac{\mathcal K(\cdot)-\mathcal K(y)}
                                      {\cdot-y}\right](z),
 \quad
 A_{\mathcal K}(z,y)=\mathcal K(y)\frac{\bar z-\bar y}{z-y}
                                      +E_{\mathcal K}(z,y).
\]
The first quotient in the Cauchy integral is uniformly bounded. Its dependence on $y$ is continuous in every finite $L^p$; the Cauchy map for $p>2$ makes $E_{\mathcal K}$ jointly continuous. The angular quotient has modulus one. Since
\[
 \partial_{\bar z}A_{\mathcal K}=\mathcal K(z)/(z-y),
\]
a solution of $\partial_{\bar z} Q_k=\partial_{\bar z}\mathscr T$ away from $p_\pm$ is
\[
 Q_k=C[(\mathcal K_{k,z}-\mathcal K_kB_k)/4]
       -\tfrac14\sum_{i=\pm}A_{\mathcal K_k}(\cdot,p_i)
       +\frac1{8\pi}\int F_k(y)A_{\mathcal K_k}(\cdot,y)\,dA_y.
\]
At $p_i$ remove its single angular term $-\mathcal K_k(p_i)(\bar z-\bar p_i)/(4(z-p_i))$. The remaining continuous value satisfies
\begin{equation}\label{fnd:finite-part}
           Q_k^{\rm fin}(p_\pm)=Q_0(q)+
                         \frac{\mathcal K(q)\bar d}{4d}+o(1).
\end{equation}
Indeed the continuous $E_{\mathcal K}$ contributions cancel with coefficients $-1/4,-1/4,+1/2$. By Lemma~\ref{fnd:midpoint}, for any fixed $0<\epsilon<1/2$ the disk $|y-s|\leq\epsilon|d|$ has mass tending to $4\pi$. There the angular quotient at $z=p_\pm$ is $\bar d/d+O(\epsilon)$, since
\[
 p_\pm-s=\pm d
\]
and the denominator has magnitude at least $(1-\epsilon)|d|$. The quotient has modulus one on the rest of the disk. After removing
\[
 |y-s|\leq\epsilon|d|,
\]
the remaining mass in a fixed small disk about $q$ is bounded in the limit by the residual mass of that disk. Outside it the angular kernel converges uniformly near the evaluation point. Taking the radius and then $\epsilon$ to zero gives the contribution of the point mass $\mathcal K(q)\bar d/(2d)$, up to $o(1)$.

The other source contributes $-\mathcal K(q)\bar d/(4d)$. The residual part is continuous after integration by the same small-disk split.

\noindent\textit{Step 2. Expansion at the two points.}\quad Near $p_i$ write $U=\log|z-p_i|^2+v_i$ and $a_i=-(v_i)_z(p_i)$. Since $F(p_i)=0$, $(v_i)_{z\bar z}(p_i)=\mathcal K(p_i)/4$, and Taylor expansion gives the expansion
\[
 \mathscr T=-\frac3{2(z-p_i)^2}+\frac{a_i}{z-p_i}
    -\frac{\mathcal K(p_i)}4\frac{\bar z-\bar p_i}{z-p_i}
    -\frac{a_i^2}{2}+o(1).
\]
After subtracting these poles and $Q_k$, the remainder $J_k$ is holomorphic on the entire smaller disk. On its fixed boundary the two double poles converge to $-3/(z-q)^2$. The singular parts of the two simple poles, as calculated below, have sum
\[
 \frac{\alpha(c,\eta)}{p_+-p_-}
 \left(\frac{1}{z-p_+}-\frac{1}{z-p_-}\right)
    =\frac{\alpha(c,\eta)}{(z-p_+)(z-p_-)}
    \longrightarrow \frac{3}{(z-q)^2}.
\]
The bounded simple-pole coefficients $-L_i^{c,\eta}$ tend to zero by the limiting balance equation. In $Q_k$ the two source contributions and the contribution of the point mass cancel on the same fixed boundary.

Uniform convergence of the remaining integral follows by splitting near its diagonal, where the kernel and the local densities are uniformly bounded, and away from it, where the kernel is continuous.

Hence the boundary values of $J_k$ converge to $\mathscr T_0-Q_0$, and the maximum principle gives
\[
 J_k(p_\pm)\to\mathscr T_0(q)-Q_0(q).
\]

\noindent\textit{Step 3. The balance equations.}\quad The balance equation gives
\[
 a_i=\frac{\alpha(c,\eta)}{p_i-p_j}-L_i^{c,\eta},
 \quad\alpha(c,\eta)=-1-\frac{2t}{c},
\]
\[
 L_i^{c,\eta}=(\log\widetilde h_2+\varphi)_z(p_i)
       -\frac tc\Gamma_{1,z}(p_i)-4\pi R_z(p_i,p_i)
       -\left(4\pi+\frac{8\pi t}{c}\right)R_z(p_i,p_j).
\]
Put $C=2\log\widetilde h_1+\log\widetilde h_2+12\pi r_g$. Lemma~\ref{fnd:pair-closure} gives $C_z(q)=0$. For $c=-1/2,\eta=0$, the identity $(r_g)_z=2R_z(z,z)+\varphi_z/(4\pi)$ implies the exact expression
\[
 L_i=C_z(p_i)+12\pi\{R_z(p_i,p_j)-R_z(p_i,p_i)\}.
\]
The signed displacement of the second argument from $p_i$ to $p_j$ is $\mp 2d$ for $i=\pm$. Taylor expansion of the above expression therefore gives
\[
 L_+-L_-=2(C_{zz}-24\pi R_{zp})d+
          2(C_{z\bar z}-24\pi R_{z\bar p})\bar d
          +O((|c+1/2|+|\eta|)|d|)+o(|d|),
\]
with derivatives in the leading terms at the limiting point. The $O$ term follows by taking the difference of the two smooth expressions in the formula for $L_i$ above, using the uniform coefficient derivatives. Adding the two finite parts gives
\[
 \frac{J_k(p_+)+J_k(p_-)+Q_k^{\rm fin}(p_+)+Q_k^{\rm fin}(p_-)}2
 =\frac{3+2\alpha-\alpha^2}{8d^2}
       +\frac{(\alpha-1)(L_+-L_-)}{4d}
       -\frac{L_+^2+L_-^2}{4}.
\]
As $L_\pm\to0$, multiplying by $d^2$ and using \eqref{fnd:finite-part} gives
\[
 \mathcal A\zeta+\mathcal B|\zeta|
       =\frac{3+2\alpha(c,\eta)-\alpha(c,\eta)^2}{8}
              +O((|c+1/2|+|\eta|)|\zeta|)+o(|\zeta|),
\]
where $\mathcal A=\mathscr T_0-C_{zz}+24\pi R_{zp}$ and $\mathcal B=\mathcal K/4-C_{z\bar z}+24\pi R_{z\bar p}$. In particular, this identity controls the constant term without an a priori relation between $\eta$ and $d$.

\noindent\textit{Step 4. The limiting coefficients.}\quad For $w_0=8\pi G(\cdot,q)+v$, put
\[
 \psi_1=\log\widetilde h_1+\varphi+12\pi R(\cdot,q)-v/2,
 \quad\psi_2=\log\widetilde h_2+\varphi+v,
 \quad\mathcal T=\psi_{1,zz}+\psi_{1,z}^2,
\]
\[
 \kappa=\Delta_g\log\widetilde h_1+16\pi-\gamma-2K_g.
\]
Here $\gamma$ is the scalar parameter before the pair collides; the equation for $v$ has residual parameter $\gamma-4\pi$. Thus
\[
 16\pi-\gamma=12\pi-(\gamma-4\pi),
\]
which identifies this trace with the canonical forcing at the resulting mixed point. If $\ell$ separated pairs collide at once, work near one resulting point $q_j$ and put
\[
 G_{\mathrm{rem}}=\sum_{a\ne j}G(\cdot,q_a),\qquad \gamma_{\mathrm{fin}}=\gamma-4\pi\ell.
\]
In the formula for one pair the scalar remainder is
\[
 v=w_{\mathrm{fin}}+8\pi G_{\mathrm{rem}}.
\]
Relative to the weights after all replacements, the two logarithmic weights in that formula contain, respectively, $4\pi G_{\mathrm{rem}}$ and $-8\pi G_{\mathrm{rem}}$. These additions cancel in $\psi_1,\psi_2$. Also
\[
 \Delta_gG_{\mathrm{rem}}=\ell-1
\]
near $q_j$, so the trace in the formula for one pair becomes the final first logarithmic weight's Laplacian plus $12\pi-\gamma_{\mathrm{fin}}-2K_g$.

Hence the same local identity uses the final canonical backgrounds even for simultaneous separated replacements. By symmetry of the Green function and $r_g(z)=R(z,z)+\varphi(z)/(4\pi)$, at $q$ one has
\[
 (r_g)_{zz}=2R_{zz}+2R_{zp}+\varphi_{zz}/(4\pi),\qquad
 (r_g)_{z\bar z}=e^\varphi/2+2R_{z\bar p}+\varphi_{z\bar z}/(4\pi).
\]
The second identity uses $R_{z\bar z}=e^\varphi/4$; the two mixed terms agree on the diagonal by the symmetry of $R$. Then
\[
 2\psi_{1,z}+\psi_{2,z}=C_z(q)=0,\qquad 2\psi_{1,zz}+\psi_{2,zz}-C_{zz}=-24\pi R_{zp}.
\]
Since $\mathscr T_0=\psi_{2,zz}-\psi_{2,z}^2/2$, these give $\mathcal A=-2(\psi_{1,zz}+\psi_{1,z}^2)$. Substitution of $\mathcal K$ in $\mathcal B$ gives
\[
 \mathcal B=-2(\log\widetilde h_1)_{z\bar z}
            -2\varphi_{z\bar z}+(\gamma/2-8\pi)e^\varphi
            =-\tfrac12 e^\varphi\kappa,
\]
where $\varphi_{z\bar z}=-K_ge^\varphi/2$. Finally
\[
 \frac{3+2\alpha-\alpha^2}{8}
     =\frac{-c-1/2+i(c+1)\eta+\eta^2/2}{c^2}.
\]
This proves \eqref{fnd:physical-pair}--\eqref{fnd:rotated-pair}. For the rotated equations, the imaginary part first gives
\[
 |\eta|=O(|\zeta|);
\]
its real part then contradicts $\kappa/2>2|\mathcal T|$. Uniformity follows by taking a convergent subsequence if any uniform assertion failed.
\end{proof}
\subsection{Nondegeneracy and separation of shadow solutions}
Let $\lambda$ index a finite family of configurations closed under pair replacement. The effective weights and the residual parameter depend on $\lambda$, whereas the positive smooth factors $h_1^*,h_2^*$ are common to all these systems.
\begin{lemma}[Propagation of vanishing off the marked set]
\label{fnd:scalar-unique}
Let $V$ be smooth on $M\setminus S$ and let $-\Delta_g\phi=V\phi$ there. If $\phi$ vanishes on a nonempty open disk, then
\[
 \phi=0\text{ on }M\setminus S,\qquad
 \phi=0\text{ on }M\ \text{if }\phi\in C^0(M).
\]
\end{lemma}
\begin{proof}
In a conformal disk centered at a boundary point of the region of open vanishing, $\phi$ is smooth and flat at the center. The equation there is $-\Delta\phi=e^\varphi V\phi$ with bounded coefficient. Put $r=e^{-s}$. For a smooth function $u$ compactly supported in a punctured disk and
\[
 a\in\mathbb Z+\tfrac12,
\]
$a>0$, expansion in $e^{ij\theta}$ and transformation in $s$ give
\[
 \|r^{-a}u\|_{L^2(ds\,d\theta)}
 \leq \frac{2}{a}\|r^{2-a}\Delta u\|_{L^2(ds\,d\theta)}.
\]
Indeed the conjugated cylinder operator is $(\partial_s-a)^2+\partial_\theta^2$. Its symbol factors as $(i\xi-a+j)(i\xi-a-j)$ for $j\in\mathbb Z$. One factor has modulus at least $a$, the other at least $1/2$. Apply the estimate to $\phi$ times a cut-off equal to one on $r<R/2$ and supported in $r<R$. At first also cut off
\[
 r<\epsilon;
\]
flatness makes the extra error tend to zero as $\epsilon\to0$ for each fixed $a$. The potential term is absorbed for large $a$, because it is bounded by $CR^2\|r^{-a}u\|$.

The remaining errors are supported in $R/2\leq r<R$ and have norm at most $C(R/2)^{-a}$. Restricting the left side to $r<R/4$ and multiplying by $(R/4)^a$ proves vanishing there as
\[
 a\to\infty.
\]
Thus the set of points near which $\phi$ vanishes is both open and closed in the connected surface $M\setminus S$. Continuity proves the last assertion.
\end{proof}
\begin{lemma}[Nondegeneracy of the shadow system and its scalar equation]
\label{fnd:genericity}
Assume uniform compactness for finitely many shadow systems on an open set of positive smooth weights. An arbitrarily small perturbation of the common smooth factors gives the following.
\begin{itemize}
\item[\textup{(i)}]
Every solution of every shadow system is nondegenerate.
\item[\textup{(ii)}]
At each solution the linearized scalar operator with fixed centers is invertible.
\item[\textup{(iii)}]
Each system has finitely many solutions.
\end{itemize}
\end{lemma}
\begin{proof}
\smallskip\noindent\textit{Step 1. The full shadow system.}\quad Choose
\[
 0<\theta<\min(1,2\min_{\lambda,a}\beta_{\lambda,a}),
\]
omitting the source minimum when $S$ is empty. We take $\theta$ strictly below an available H\"older exponent so that the corresponding embeddings are compact.

On separated center neighborhoods the map defining the shadow system is $C^1$ from $C^{2,\theta}_0(M)$ times the finite center space to $C^{0,\theta}_0(M)$ times the same center space, and its derivative in the unknowns is Fredholm of index zero.

The scalar coefficients are $C^{0,\theta}$ even when the source weights are not smooth at $S$. Uniform $C^0$ compactness on separated center domains gives $W^{2,q}$ bounds for $q>2$, because the scalar densities are bounded. Using the $C^1$ bound and the power exponent $2\min\beta$ when
\[
 S\ne\varnothing,
\]
the densities are uniformly bounded in $C^{0,\theta'}$ for some
\[
 \theta'>\theta.
\]
Schauder estimates and the compact embedding give compactness in $C^{2,\theta}$, and local regularity away from the fixed sources gives the higher regularity used at the centers. The fixed weight parameters introduce no new parameter in this map.

For $\eta_i=\delta\log h_i^*$, the universal derivatives are
\[
 D_{\eta_1}F_0=0,\quad
 D_{\eta_1}F^A_i=\nabla\eta_1(p_i),\quad
 D_{\eta_1}F^B_j=2\nabla\eta_1(q_j),\quad
 D_{\eta_2}F_0=-2f(\eta_2-\langle\eta_2\rangle_f).
\]
At a single solution, varying $\nabla\eta_1$ independently at its distinct centers shows that every functional in the cokernel of the full derivative has zero center components. The remaining functional annihilates the image of the fixed-center scalar linearization
\[
 L_f=-\Delta_g-2f(I-\langle\cdot\rangle_f).
\]
This operator is Fredholm of index zero on the indicated mean-zero H\"older spaces and symmetric under the $L^2$ pairing.

The pairings against a basis of its kernel space are independent cokernel functionals; the index shows that these exhaust the cokernel. Thus the functional is represented by a kernel function $\phi$, and arbitrary variations $\eta_2$ give
\[
 f(\phi-\langle\phi\rangle_f)=0.
\]
Since $f>0$ outside finitely many points, $\phi$ is constant, and its mean-zero representative is zero. The universal map is therefore surjective for each configuration.

Choose finitely many common smooth parameter directions spanning all cokernels at the original solutions by compactness and openness of surjectivity.

Uniform compactness puts every zero for a small parameter ball in these surjectivity neighborhoods. The zero set of the full map over that ball is a $C^1$ manifold of dimension the number of parameters, by surjectivity and the index-zero property.

At a point of this manifold, its projection onto parameter space is a submersion if and only if the derivative in the unknowns is invertible (use the tangent equation and index zero).

Sard's theorem for equal-dimensional $C^1$ maps, applied to each configuration, gives a full-measure set of regular parameters.

Intersecting these finitely many sets gives invertibility of the derivatives of the full system at all solutions. Compactness and discreteness make each solution set finite.

\smallskip\noindent\textit{Step 2. The scalar linearization.}\quad For the scalar linearization set
\[
 L_f\xi=-\Delta_g\xi-2f(\xi-\langle\xi\rangle_f).
\]
On any connected disk $B$ avoiding $S$ and every center of the finite set of solutions, $f>0$. If the kernel space is trivial there is nothing to prove. Otherwise choose a basis $\xi_1,\ldots,\xi_d$ of the mean-zero kernel space and put
\[
 \phi_i=\xi_i-\langle\xi_i\rangle_f.
\]
We claim that for some $v\in C_c^\infty(B)$ with $\int_Mv=0$ the symmetric matrix
\[
 M(v)=\left(\int_Bv\phi_i\phi_j\right)_{i,j}
\]
is invertible. Suppose otherwise, and choose $M_0$ of maximal rank $r<d$ in this linear space of matrices. Fix a nonzero $e\in\ker M_0$. On $\operatorname{ran}M_0\oplus\mathbb Re$ the matrix $M_0$ is $\operatorname{diag}(A,0)$ with $A$ invertible.

For every admissible $v$, the corresponding $(r+1)$-minor of $M_0+tM(v)$ vanishes identically in $t$, since its rank cannot exceed $r$. Its linear coefficient is $\det(A)\,e^TM(v)e$, so
\[
 \int_Bv\left(\sum_i e_i\phi_i\right)^2=0
 \quad\hbox{for every }v\in C_c^\infty(B),\quad\int_Bv=0.
\]
Therefore $\phi^2$ is constant on $B$, where $\phi=\sum_i e_i\phi_i$. By continuity and connectedness, $\phi$ itself is constant there. The equation
\[
 -\Delta_g\phi=2f\phi,\qquad f>0
\]
force that constant to be zero. Lemma~\ref{fnd:scalar-unique} applies with
\[
 V=2f
\]
on $M\setminus S$, since the ordinary zeros are smooth there. Continuity of the kernel function extends its vanishing to $M$. Thus $\sum_i e_i\xi_i$ is constant; its mean-zero normalization makes it zero, contradicting the choice of a basis of kernel functions. This proves the claim.

\smallskip\noindent\textit{Step 3. Perturbation of the weights.}\quad Let $-\Delta_gz=2v$, $\int_Mz=0$, and vary the common smooth factors by
\[
 h_{2,t}^*=h_2^*(1+tv/f)e^{-tz},\qquad
 h_{1,t}^*=h_1^*e^{tz/2},\qquad w_t=w+tz.
\]
The quotient $v/f$ is smooth and compactly supported in $B$, extended by zero.

The chosen solution retains exactly the same centers: its scalar density becomes $f+tv$, its ordinary potential $\log h_{1,\lambda}-w/2$ is unchanged, and its mixed potential changes by $\log(1+tv/f)$, which vanishes near every center. Since
\[
 \int v=0,
\]
the total mass of $f_t=f+tv$ remains $\gamma$, and direct differentiation of the weighted average gives
\[
 \left\langle\xi_i,
    \left.\frac{d}{dt}\right|_{t=0}L_{f_t}\xi_j\right\rangle_{L^2}
       =-2\int_Bv\phi_i\phi_j=-2M(v)_{ij}.
\]
On the orthogonal complement of its kernel space, $L_f$ is invertible. After eliminating the orthogonal complement, the operator induced by $L_{f_t}$ on the kernel space is therefore
\[
 -2tM(v)+O(t^2),
\]
which is invertible for every sufficiently small nonzero $t$. Hence the fixed-center scalar block is invertible.

Nondegeneracy of the full system persists by openness. The implicit function theorem continues its other nondegenerate solutions. Compactness on the weight neighborhood prevents any new zeros outside those continuation neighborhoods.

Repeat over the finite collection, keeping each previously obtained invertibility by taking subsequent perturbations small. The resulting common weights have both required properties.
\end{proof}
\begin{lemma}[Separation of centers in different shadow solutions]
\label{fnd:cross-separation}
Under the assumptions of Lemma~\ref{fnd:genericity}, the common weights may also be chosen as follows.
\begin{itemize}
\item[\textup{(i)}]
No ordinary center of any solution coincides with a mixed center of any solution.
\item[\textup{(ii)}]
Two solutions with $U$ and $U'\ne U$ ordinary centers have no mixed center in common.
\end{itemize}
For the same $U$ and different source configurations, mixed centers may coincide.
\end{lemma}
\begin{proof}
\smallskip\noindent\textit{Step 1. The scalar equation at fixed centers.}\quad After solving the scalar equation, whose linearization is invertible, the derivative of the reduced center equations is a symmetric, invertible matrix $S$. This follows from the reduced energy $E_c-E_s/(8\pi)$, where
\[
 E_s(w,P)=\tfrac12\int|\nabla w|^2
 -2\gamma\log\int h_{2,\lambda}e^{-4\pi\sum_iG(\cdot,p_i)}e^w,
\]
and $E_c$ sums the potentials depending on a single point in $F^A_{i,0},F^B_j$ and the Green interactions with coefficients $8\pi,12\pi,24\pi$ for an unordered pair of two ordinary points, one point of each type, and two mixed points, respectively.

Along scalar solutions, the $w$-derivative of $E_s$ vanishes and its explicit $p_i$-derivative is
\[
 8\pi\int_M f(x)\nabla_{p_i}G(x,p_i)\,dV_g=4\pi\nabla w(p_i).
\]
Thus the reduced center gradient is $(F^A,F^B)$. The scalar equation does not depend on $Q$, hence
\[
 S_{q,p}=12\pi\nabla_q\nabla_pG(q,p),\qquad
 D_{\text{gradients of }\log h_1^*}(P,Q)=-S^{-1}D,
 \quad D=\operatorname{diag}(I_{2U},2I_{2V}).
\]
Independent variations of the gradient at the centers of a given solution therefore vary all its coordinates independently. To examine the mixed Green derivatives, choose conformal coordinates near $q$ and $p$ and write
\[
 H(q,p)=\bigl(\partial_{q_a}\partial_{p_b}G(q,p)\bigr)_{a,b=1}^2.
\]
For fixed $q$, differentiating the Green equation in $q$ shows that $\partial_{q_a}G(q,\cdot)$ is harmonic away from $q$. Thus the entries of $H(q,p)$, and hence its determinant, are real analytic in the conformal coordinate of
\[
 p\ne q.
\]
In a common coordinate near the diagonal, the leading term of this determinant is $-1/(4\pi^2|p-q|^4)$, so it is not identically zero. Passing from coordinate derivatives to intrinsic gradients only multiplies $H$ by invertible metric matrices.

Consequently the zero set of $\det(\nabla_q\nabla_pG(q,p))$ has no interior, even when the metric is merely smooth. Successive small changes of the weights therefore make every mixed--ordinary Green matrix nonsingular for all solutions in the finite collection.

\smallskip\noindent\textit{Step 2. Variation of the center equations.}\quad Next vary the second derivatives of $\log h_1^*$, keeping its values and first derivatives fixed at all current centers.

All solutions and scalar blocks remain fixed. Under the particular common choice of second derivative $tI$ at every distinct center, $S$ becomes $S+tD$. For large $|t|$,
\[
 ((S+tD)^{-1})_{q,q}=(2t)^{-1}I+O(t^{-2}),\qquad
 ((S+tD)^{-1})_{q,p}=-6\pi t^{-2}\nabla_q\nabla_pG(q,p)+O(t^{-3}).
\]
Each corresponding determinant of $S^{-1}$ is therefore a nonzero rational function of the second derivatives. Even when several solutions share a center, taking all second derivatives equal to $tI$ shows that each numerator is a nonzero polynomial.

The product of these finitely many numerators is not identically zero. Thus arbitrarily small variations of the second derivatives give, for $R=S^{-1}$,
\[
                  \det R_{q,q}\ne0,\qquad \det R_{q,p}\ne0
\]
for every mixed diagonal or mixed--ordinary block under consideration. All previous open conditions persist. Suppose now
\[
 p_i^X=q_j^Y=x
\]
for two solutions, possibly of different source configurations. Vary the gradient of the common second weight by the local perturbation
\[
 \eta_2(z)=\vartheta((z-x)/r)\,\xi\cdot(z-x),
 \qquad \|\eta_2\|_{C^{0,\theta}}=O(r^{1-\theta}|\xi|),
\]
whose support avoids every other center of $X,Y$ and $S$. Its scalar forcing and the resulting scalar correction are $o_r(1)|\xi|$. It has no direct force in the ordinary balance equation of $X$, but has direct force $\xi$ in the marked mixed balance equation of $Y$. Thus
\[
 D_\xi(p_i^X-q_j^Y)=R^Y_{q_j,q_j}+o_r(1),
\]
which is invertible for fixed sufficiently small $r$. This removes all ordinary--mixed coincidences by finite-dimensional transversality. Different residual parameters or effective weight parameters only change the finite constants in this estimate.

\smallskip\noindent\textit{Step 3. Coincident mixed centers.}\quad Finally suppose $q_i^X=q_j^Y$ and $U_X>U_Y$. Since
\[
 U_X>U_Y,
\]
there is a point $p\in P_X\setminus P_Y$. The ordinary--mixed separation just obtained ensures $p\notin Q_Y$ as well.

A local perturbation of the gradient of $\log h_1^*$ at $p$, avoiding all other centers of $X$ and all centers of $Y$, leaves the solution $Y$ exactly fixed: its scalar equation does not involve $h_1^*$ and all derivatives at its centers stay fixed. For $X$,
\[
 D_\xi(q_i^X-q_j^Y)=-R^X_{q_i,p},
\]
which is invertible. In each of the two separation steps, use the finite-dimensional span of these variations for all coincident centers.

Each equality of the two continued coordinates then has codimension two near the original weights, so the complement of their finite union contains arbitrarily small parameters.

Compactness and the implicit function theorem account for all continued solutions, and sufficiently small parameters preserve the earlier open properties.
\end{proof}
\subsection{Choice of the weights and computation of the degree}
Under $z=a w+bw^2+cw^3+O(w^4)$, $a\ne0$, the added term in $\psi_1$ is $2\log|dz/dw|-6\log|z/w|$. Its first complex derivative at zero is $2b/a-3b/a=-b/a$. Its second is $(6c/a-4b^2/a^2)-3(2c/a-b^2/a^2)=-b^2/a^2$. Thus
\[
 \psi_{1,w}=a\psi_{1,z}-b/a,\qquad \psi_{1,ww}=a^2\psi_{1,zz}+2b\psi_{1,z}-b^2/a^2
\]
there. Adding the square of the first identity to the second gives $\mathcal T_w=a^2\mathcal T_z$. Since $e^{\varphi_w}=|a|^2e^{\varphi_z}$ at the center, \eqref{fnd:gap} is invariant. With $\varphi=0$ at the center the inequality is $\kappa>4|\mathcal T|$. We abbreviate $\mathcal T$ by $T$ below.
\begin{proof}[Proof of Theorem~\ref{fnd:common-weights}]
\noindent\textit{Step 1. The induction on the number of ordinary centers.}\quad We argue by induction on $U$, treating all configurations with the same number of ordinary centers at each step.

Replacing a pair decreases $U$ by two and preserves the concentrations at the prescribed points. Thus every configuration produced by a pair collision has already been treated.

At the smallest value of $U$ in the family, Lemma~\ref{fnd:source-exclusion} excludes approaches to a prescribed source and leaves only regular ordinary pairs. The same lemma gives a positive residual mass through the scalar exclusion.

For a pair, Lemma~\ref{fnd:pair-closure} applies since the coefficients are fixed smooth factors, or converge in $C^1$ with bounded $C^2$ norms when testing a small $C^3$ neighborhood. It gives a configuration in $\mathscr F$ with fewer ordinary centers and positive residual parameter.

This contradicts the minimality of $U$. Thus all free centers stay separated. Lemma~\ref{fnd:separated-scalar-compactness} excludes scalar blow-up without center collision.

Its uniform version, together with the collision exclusion for weights converging to the selected weights, also gives compactness on a small smooth-weight neighborhood.

For the induction step, assume compactness and the strict trace inequalities for the configurations already treated. An unbounded sequence of solutions for a new configuration has no cluster at a prescribed source, and every regular collision consists of a pair of ordinary centers.

Lemma~\ref{fnd:pair-closure} and Proposition~\ref{fnd:pair-identity} give a solution with fewer ordinary centers and a mixed collision point at which
\[
 |\kappa|=4e^{-\varphi}|T|.
\]
This contradicts the strict inequality. For separated centers, apply Lemma~\ref{fnd:separated-scalar-compactness} at the fixed residual mass with the current exclusions. All solution sets up to this stage are therefore compact.

The argument also proves uniform compactness in a sufficiently small smooth-weight neighborhood: otherwise a sequence of weights tending to the current weights would give the same forbidden collision limit.

Apply Lemmas~\ref{fnd:genericity}--\ref{fnd:cross-separation} to the entire finite collection considered so far. Smallness preserves the positive gaps of earlier stages.

\noindent\textit{Step 2. Choice of the weights.}\quad Let $x_1,\ldots,x_s$ be the distinct mixed-coordinate values of all solutions in the current stage, combining coincidences between source configurations with the same $U$.

Choose disjoint coordinate disks around them, separated from $S$, every ordinary coordinate of every solution considered so far, and every mixed coordinate of an earlier stage. Such disks exist by Lemma~\ref{fnd:cross-separation}. Normalize
\[
 \varphi(x_l)=0.
\]
With a fixed cut-off $\vartheta=1$ on $B_{1/2}$, put
\[
 b_r(z)=\sum_l L_l|z-x_l|^2\vartheta((z-x_l)/r),\qquad
 h_{1,r}^*=h_1^*e^{b_r/2},\qquad
 h_{2,r}^*=h_2^*e^{-b_r}.
\]
The mixed potential $2\log h_1^*+\log h_2^*$ is exactly unchanged for every source configuration. For fixed $L_l$,
\[
 \|b_r\|_{C^1}=O(r),\qquad
 \|b_r\|_{C^{0,\theta}}=O(r^{2-\theta}),\qquad
 \|b_r\|_{C^2}=O(1).
\]
The weights are unchanged on fixed neighborhoods of $S$ and of all earlier mixed coordinates. Although third derivatives of $b_r$ may grow as $r\downarrow0$, the $C^1$ convergence and uniform $C^2$ bound suffice for Lemma~\ref{fnd:pair-closure}.

\noindent\textit{Step 3. Compactness for the perturbed weights.}\quad We prove compactness uniformly as $r\downarrow0$ for all configurations treated up to this step.

Lemma~\ref{fnd:source-exclusion}, whose estimates use uniform $C^1$ bounds for the scalar factors and bounded center fields, excludes concentration at a prescribed source; the weights near each source are unchanged. Any collision at a regular point would involve a pair of ordinary centers. Since
\[
 h_{i,r}^*\to h_i^*
\]
in $C^1$ with uniform $C^2$ bounds, Lemma~\ref{fnd:pair-closure} first gives a solution with strictly fewer ordinary centers for the unperturbed weights. Each new mixed collision point is therefore an earlier mixed coordinate, so a fixed disk about it avoids all supports of $b_r$.

On that disk the local coefficients are unchanged and the external Green fields converge smoothly. We may therefore apply Proposition~\ref{fnd:pair-identity}: its identity in the original coordinates \eqref{fnd:physical-pair} contradicts the strict trace inequality at the earlier solution.

For a separated sequence the effective weight parameters and the scalar parameter stay fixed, and the positive smooth factors converge in $C^1$ by the weight convergence and separation. Thus Lemma~\ref{fnd:separated-scalar-compactness} applies. This gives uniform compactness of the scalar fields and hence of the separated solutions.

\noindent\textit{Step 4. Persistence of the solutions.}\quad Let $Z_1,\ldots,Z_N$ be the solutions of the finite subfamily before the perturbation, and choose disjoint small coordinate neighborhoods with fixed orderings of the centers.

Uniform compactness implies that all perturbed solutions lie in these neighborhoods for small $r$.

On each neighborhood the ordinary coordinates stay away from the bump supports, the mixed potential is unchanged, and the scalar coefficient changes by $O(r^{2-\theta})$ in $C^{0,\theta}$.

Its derivatives in the moving ordinary centers have the same estimate, since $e^{-4\pi G(x,p)}$ is jointly smooth through $x=p$. Consequently the maps defining the shadow systems satisfy
\[
 \|F_r-F_0\|_{C^1(\mathcal U_l)}\le C r^{2-\theta}.
\]
Fix an old root $Z_l$. On a fixed small ball about $Z_l$, put
\[
 A_l=DF_0(Z_l),\qquad \Phi_r(Z)= Z-A_l^{-1}F_r(Z).
\]
By continuity of $DF_0$ and the displayed $C^1$ estimate,
\[
 \operatorname{Lip}(\Phi_r)\le\tfrac12,\qquad
 \|\Phi_r(Z_l)-Z_l\|=\|A_l^{-1}F_r(Z_l)\|=O(r^{2-\theta}).
\]
For small $r$ the ball is invariant. Its unique fixed point $Z_{l,r}$ satisfies $\|Z_{l,r}-Z_l\|\le2\|\Phi_r(Z_l)-Z_l\|$. In particular,
\[
\|w_{l,r}-w_l\|_{C^{2,\theta}}
 +|P_{l,r}-P_l|+|Q_{l,r}-Q_l|=O(r^{2-\theta})=o(r).
\]
Joint nondegeneracy and scalar invertibility persist. Every continued current mixed coordinate initially at $x_l$ lies in $B_{r/2}(x_l)$. There
\[
 (b_r)_{zz}=0,\qquad (b_r)_z=O(r^{2-\theta}),\qquad
 \Delta_gb_r=4L_l+o(1).
\]
For every current solution at this coordinate, regardless of its source configuration,
\[
 \kappa_{\lambda,r}=\kappa_\lambda(x_l)+2L_l+o(1),\qquad
 T_{\lambda,r}=T_\lambda(x_l)+o(1).
\]
Choose the one number $L_l$ to satisfy
\[
 2L_l>\max_{(\lambda,Z):\,x_l\in Q(Z)}
             \{4|T_\lambda(Z,x_l)|-\kappa_\lambda(x_l)\}+2.
\]
The maximum is finite. The same perturbation therefore gives the strict inequality at every coincident mixed coordinate in this step.

Earlier mixed coordinates avoid the supports, so their traces and coefficients change only by $o(1)$ through the continued positions and scalar fields; their strict gaps persist.

For small fixed $r$ the stage is complete. Repeating over the finitely many $U$ values proves the theorem.
\end{proof}
At fixed $U$, compactness follows from exclusion of collisions at the sources, scalar compactness for separated centers, the limiting shadow system for two colliding ordinary points and the second-order identity at their limit.

This argument uses the strict inequalities only at smaller values of $U$.

Once compactness holds, we choose the weights for nondegeneracy and finally impose the trace inequalities at the current value of $U$. The induction uses neither the construction of projected solutions nor the estimates placing all blowing-up sequences in that family.
\begin{proof}[Proof of Corollary~\ref{fnd:extra-traces}]
\smallskip\noindent\textit{Step 1. The prescribed points in (i).}\quad At each prescribed source $p\in S$ use a quadratic bump in $\log h_1^*$ with zero value and gradient at $p$.

Fix its Laplacian and let its support shrink to $p$, away from all the free centers in the finite solution sets. The bumps tend to zero in $C^1$ and are bounded in $C^2$.

The scalar equation and the balance maps near all selected solutions are unchanged. Lemma~\ref{fnd:source-exclusion} rules out source collisions and gives a positive residual mass.

For hypothetical regular pair limits, Lemma~\ref{fnd:pair-closure} identifies the limiting solution of the unperturbed lower-$U$ system after all replacements, with its mixed centers away from the shrinking supports.

On fixed disks about these centers the perturbations eventually vanish, so Proposition~\ref{fnd:pair-identity} contradicts the original strict gaps.

For separated configurations, Lemma~\ref{fnd:separated-scalar-compactness} gives compactness with the same effective weight parameters and scalar parameter, upgraded as in the proof of Lemma~\ref{fnd:genericity}. The implicit function theorem now excludes additional separated solutions. This gives the required inequalities.

\smallskip\noindent\textit{Step 2. The ordinary centers in (ii).}\quad For the ordinary mass coefficients, it suffices to consider
\[
 U\geq1,
\]
since $U=0$ gives the scalar base case. Choose disjoint small disks at every distinct location in the finite union of the free center sets, avoiding $S$, and multiply $h_1^*$ by
\[
           \exp\left(\tfrac14\sum_jc_j\chi_jd_g(\cdot,x_j)^2\right).
\]
The values and gradients stay fixed at all selected centers. At an ordinary center in the configurations under consideration, take $\varphi(p_i)=0$. Its ordinary background is
\[
 \psi_i^A=\log h_1+\varphi+8\pi R(\cdot,p_i)
            +8\pi\sum_{\ell\ne i}G(\cdot,p_\ell)-w/2 .
\]
The matched factor is
\[
 d_i=e^{\psi_i^A(p_i)/2}/2,
\]
obtained by matching the coefficient $4s_i^2$ of the ordinary density tail to $t^2e^{\psi_i^A(p_i)}$ for $s_i=td_i$. The scalar density vanishes at $p_i$, and here $\gamma=\rho_2$. Applying $\Delta_g$ to the smooth background at $p_i$ gives
\[
 \kappa_i=\Delta_g\log h_1(p_i)+8\pi U-\rho_2-2K_g(p_i).
\]
Thus the selected solutions keep their scalar densities and width coefficients $d_i$, while the ordinary traces increase by the indicated $c_j$.

The quantity $L_0$ of each ordinary solution is affine in the $c_j$, and each equation $L_0=0$ is a proper affine hyperplane: the coefficient at any one of its center locations is a positive sum of $2\pi d_i^2$.

Avoid the finitely many hyperplanes by arbitrarily small $c_j$. Nondegeneracy, invertibility of the scalar linearization and mixed gaps persist by openness; the traces at the prescribed points are unchanged.

Compactness and the implicit function theorem again exclude additional solutions.
\end{proof}
\begin{lemma}[Deformation of the shadow system]
\label{fnd:decoupling}
For the chosen common weights the shadow system admits a homotopy to the decoupled system with compact solution set along the homotopy. At the decoupled endpoint,
\[
 (F_0^A,F^B)=0,
\]
and the scalar equation is solved with the centers fixed.
\end{lemma}
\begin{proof}
\smallskip\noindent\textit{Step 1. The homotopy.}\quad Use
\[
 (I+\eta J)F_{i,0}^A+c\nabla w(p_i)=0,\qquad (I+\eta J)F_j^B=0.
\]
First fix $c=-1/2$ and increase $\eta$ from zero. Collisions as $\eta\downarrow0$ would, by the cluster identities and source exclusion, be ordinary pairs.

Lemma~\ref{fnd:pair-closure} gives limiting solutions with fewer ordinary centers, and the uniform strict inequality at the mixed centers and \eqref{fnd:rotated-pair} exclude these collisions. Hence the solution set remains compact for all sufficiently small $\eta$.

\smallskip\noindent\textit{Step 2. Collisions.}\quad Consider all centers converging to one point. Denote the radial and angular moments of $(F_0^A,F^B)$ by $R_k,A_k$, and those of $\nabla w$ over the ordinary centers of this group by $R_k^w,A_k^w$.

The pair terms cancel in the sum in the manner used for \eqref{fnd:virial}; the source and pair terms have zero angular moment and all smooth terms contribute $o(1)$. Thus $R_k$ converges to the source and pair coefficients of that calculation and
\[
 A_k\to0.
\]
Scalar normalization and \eqref{fnd:scalar-cluster-moments} are independent of the rotation: $R_k^w\to4n(n-b-u)$ and $A_k^w\to0$. The center equations imply
\[
 R_k-\eta_k A_k+c_k R_k^w=0,\qquad
 A_k+\eta_k R_k+c_k A_k^w=0.
\]
If $\eta_k\to0$ and $c_k=-1/2$, the first equation gives
\[
 u^2+3uv+3v^2+n^2-un=a(u+2v)+b(v+n),
\]
which is \eqref{fnd:virial}. Its classification therefore applies to the limiting clusters at the first stage of the homotopy. When instead
\[
 \eta_k\to\eta_*\ne0,\qquad c_k\to c_*,
\]
the second equation forces
\[
 R_k\to0.
\]
At a source the limit of $R_k/2$ is $a(u+2v)+bv-u^2-3uv-3v^2$. For the three choices in \eqref{fnd:effective}, the pairs of coefficients of $A_p,B_p$ are $(u+2v,v)$, $(-u-v,v)$, and $(-u-v,-u-2v)$. Each is nonzero for $u+v\geq1$, so \eqref{fnd:independence} excludes a source cluster of centers.

At a regular point the pair moment $-2u(u-1)-6uv-6v(v-1)$ is negative for $u+v\geq2$. These conclusions hold also for $c_*=0$. For $c_*\ne0$ the first equation additionally gives
\[
 n(n-b-u)=0,
\]
hence $n=0$ by $2n<b+u$. In all cases the scalar equation and \eqref{fnd:scalar-exclusion} are unchanged. Once the centers are separated, Lemma~\ref{fnd:separated-scalar-compactness} excludes scalar blow-up, including at $c=0$.

\smallskip\noindent\textit{Step 3. The endpoint.}\quad We may consequently increase $c$ to zero at a fixed nonzero $\eta$, and finally remove $I+\eta J$ from the balance equations.

These matrices preserve zeros and have positive determinants. The scalar equation has remained in the system throughout. This proves compactness along the required homotopy.
\end{proof}
\begin{lemma}[Index of the decoupled balance equations]
\label{fnd:center-index}
Set $c_0=\chi(M)-|S|$. On unordered configurations the index of $(F_0^A,F^B)$ is
\[
              \binom{c_0}{U+V}\binom{U+V}{U}.
\]
The zeros lie in a fixed compact subset of configuration space throughout the homotopy in the proof.
\end{lemma}
\begin{proof}
\smallskip\noindent\textit{Step 1. The forces at the boundary.}\quad At a source the logarithmic force coefficients are
\[
 d_A=2\alpha,\qquad d_B=2(2\alpha+\beta).
\]
Allow them to be complex, where multiplication by $i$ means rotation through $+\pi/2$ on real tangent vectors, and keep the pair interactions fixed.

To realize this deformation, for each source $p$ choose a cut-off logarithm $\ell_p$ supported in a coordinate disk about $p$ and equal to $\log|z_p|$ near $p$.

Subtract the original real source coefficients times $\ell_p$ from each potential depending on a single point; the remaining potentials of one point are smooth across the sources.

For the coefficient $d=a+ib$ at $p$ and of the relevant center type, add $a\nabla_g\ell_p+bJ\nabla_g\ell_p$ to the resulting smooth field at that point.

This defines a smooth field on configuration space, invariant under permutations of equal types, and changes only the indicated source moments. At a real endpoint the added field is again a gradient. A source cluster with
\[
 u+v\geq1
\]
gives
\begin{equation}\label{fnd:center-resonance}
 ud_A+vd_B=2u(u-1)+6uv+6v(v-1)
\end{equation}
by its combined radial and angular moment. For fixed $U,V$ these are finitely many complex affine hyperplanes in the coefficient space. The initial coefficients avoid them: after division by two the real equation is
\[
 a(u+2v)+bv=u^2+3uv+3v^2,
\]
and each effective pair \eqref{fnd:effective} is nonresonant by \eqref{fnd:independence}. Choose $d_A,d_B<0$ real at the endpoint. Each source moment on the left of \eqref{fnd:center-resonance} is then negative. The complex line through the initial and final coefficients is not contained in any forbidden hyperplane.

Its parameter plane minus finitely many points is path connected, so join the two parameters there. Use the image of that path for the deformation.

The path is compact and avoids every hyperplane. Hence its distance from their union is positive, and no center approaches a prescribed source. A regular collision has moment
\[
 -2u(u-1)-6uv-6v(v-1)<0\qquad(u+v\geq2).
\]
The smooth Green and metric terms vanish in the limiting moments at
every collision scale.

\smallskip\noindent\textit{Step 2. The degree of the gradient.}\quad
At the endpoint the center field is the gradient of a potential $E$.
By symmetry of $G$, its pair interaction is
\[
 8\pi\sum_{i<i'}G(p_i,p_{i'})
 +12\pi\sum_{i,j}G(p_i,q_j)
 +24\pi\sum_{j<j'}G(q_j,q_{j'}).
\]
Write $Y$ for the configuration space modulo permutations within each
type. Every singular term of $E$ is bounded below, and $d_A,d_B<0$ give
\[
 G(z,z')\longrightarrow+\infty\quad(z\to z'),\qquad
 d_A\ell_p(z),\ d_B\ell_p(z)\longrightarrow+\infty\quad(z\to p).
\]
Thus $E$ is proper and bounded below. Its critical set is compact by
Step~1. Make a compactly supported, $C^1$-small Morse perturbation near
this set in $Y$. Keeping the notation $E$, choose a regular value $a$
above all critical values and put
\[
 K_a=\{z\in Y:E(z)\le a\},\qquad
 W=-\frac{\nabla E}{|\nabla E|^2}\quad\text{on }Y\setminus K_a.
\]
The flow $\Phi_s$ of $W$ stays in the compact strip
$\{a\le E\le E(z)\}$ until it reaches $K_a$, since
\[
 \frac{d}{ds}E(\Phi_s(z))=-1,\qquad
 E(\Phi_{E(z)-a}(z))=a\quad(E(z)>a).
\]
Fix $K_a$ and flow each exterior point $z$ for time $t(E(z)-a)$,
$0\le t\le1$. This retracts $Y$ onto $K_a$. Since $\nabla E$ points
outward on $\partial K_a$, homotopy invariance and the
Poincar\'e--Hopf theorem give
\[
 \ind_Y(F_0^A,F^B)=\chi(K_a)=\chi(Y).
\]

\smallskip\noindent\textit{Step 3. The Euler characteristic.}\quad
Put $n=U+V$. Let $C_j$ be the space of $j$ ordered distinct points of
$M\setminus S$, with $C_0=\{\ast\}$. For $1\le j\le n$, forgetting
the last point gives
\[
 \pi_j:C_j\longrightarrow C_{j-1},\qquad
 \pi_j^{-1}(z_1,\ldots,z_{j-1})
 \cong M\setminus\bigl(S\cup\{z_1,\ldots,z_{j-1}\}\bigr).
\]
Cut-off translations in disjoint disks about the retained points locally
trivialize $\pi_j$. The same exhaustion with one type of center gives
each base finite CW type. Filtering a finite CW model by cells gives
multiplicativity for any monodromy, hence
\[
 \chi(C_0)=1,\qquad
 \chi(C_j)=\chi(C_{j-1})\,\chi(\pi_j^{-1}(z_1,\ldots,z_{j-1}))
          =(c_0-j+1)\chi(C_{j-1}).
\]
The covering $C_n\to Y$ has degree $U!V!$, and permutations of
two-dimensional blocks preserve orientation. With the empty product
equal to one, we obtain
\[
 \ind_Y(F_0^A,F^B)
 =\frac{1}{U!V!}\prod_{j=1}^{n}(c_0-j+1)
 =\binom{c_0}{U+V}\binom{U+V}{U}.\qedhere
\]
\end{proof}
\begin{proof}[Proof of Proposition~\ref{fnd:shadow-degree}]
\smallskip\noindent\textit{Step 1. The degree on the scalar space.}\quad Take the scalar component $w-2G*(f-\gamma)$ on the mean-zero space of Lemma~\ref{fnd:genericity}, followed by the center field in local coordinates.

For the fixed-center scalar degree, $2G*(f-\gamma)$ maps bounded $C^0$ sets compactly into $C^{1,\delta}$, for any
\[
 0<\delta<1,
\]
by choosing $q>2/(1-\delta)$ in the $W^{2,q}$ estimate. On a bounded $C^{1,\delta}$ set, multiplication by the fixed weight and division by the positive normalizing integral give a uniform $C^{0,\theta_1}$ bound on $f$ for some
\[
 \theta_1>\theta
\]
in the H\"older margin of Lemma~\ref{fnd:genericity}. Schauder estimates then give a bounded $C^{2,\theta_1}$ image, compact in $C^{2,\theta}$. For the successive pairs of spaces, the stronger space contains both compact images and is continuously embedded into the weaker one.

The finite-dimensional approximation in the proof of Lemma~\ref{lem:local-degree-spaces} identifies the degrees in each successive pair of spaces for the same isolated compact zero set.

\smallskip\noindent\textit{Step 2. The scalar sign and the center degree.}\quad At a selected zero the second derivative of the scalar functional at fixed centers has quadratic form
\[
 \int_M|\nabla v|^2-2\int_M f(v-\langle v\rangle_f)^2.
\]
In the mean-zero Dirichlet inner product the preconditioned scalar derivative is self-adjoint by the symmetry of the displayed bilinear form and is a compact perturbation of the identity.

Its local sign is $(-1)^{m_S}$, where $m_S$ is the number of its negative eigenvalues counted with multiplicity, equivalently the scalar Morse index.

To compare the Dirichlet realization with the $C^0$ realization, use the exponent estimate of Lemma~\ref{fnd:exponential} and the bounded fixed-center weight.

The common compact Green image on bounded sets is then in $C^0\cap H^1$, and the argument in the proof of Lemma~\ref{lem:local-degree-spaces} applies to these two spaces. Together these comparisons give the same local degree in the H\"older space.

For the full shadow system, use block triangular homotopies with identity diagonal to eliminate the two off-diagonal blocks of its derivative. Differentiating the implicit scalar solution gives the remaining block $DF_{\rm red}$. Hence the local shadow degree is
\[
 d_S(Z)=(-1)^{m_S}\operatorname{sgn}\det DF_{\rm red}.
\]
The homotopy in Lemma~\ref{fnd:decoupling} preserves the full degree. During Lemma~\ref{fnd:center-index}'s deformation, the center zeros stay in a compact separated set. With the scalar parameter and the effective weight parameters fixed, Lemma~\ref{fnd:separated-scalar-compactness} gives a uniform bound for the scalar solutions, upgraded to $C^{2,\theta}$ as above.

The same lemma makes the fixed-center scalar degree constant along separated paths of centers and positive smooth factors.

Move the coordinates one at a time to distinct temporary points away from all initial and final positions, then to their targets, joining points in the complement of the other finitely many points at each step.

Perturb the decoupled center field in its compact isolating region to have finitely many transverse zeros, keeping the scalar equation fixed. In a small separated coordinate ball about each zero, freeze the scalar weight there along straight coordinate segments.

The center equation stays fixed, and the uniform scalar bound keeps this homotopy admissible on a common scalar ball.

By the product rule the local degree in each neighborhood is the index of the center equations times the common scalar degree. Excision and homotopy invariance give the formula for the original center field.

This product contains the degree of the scalar equation at fixed centers. Proposition~\ref{fnd:scalar-degree-internal} computes it from the local degree jump for one component.
\end{proof}
\section{Bubbling at regular points and the degree jump}
\label{reg:section}
We consider solutions that blow up only at regular points. Let $P=(p_i)$ and $Q=(q_j)$ denote the $U$ ordinary and $V$ mixed blow-up points, respectively, and put
\[
 n=U+V,\qquad m=U+2V,\qquad
 \rho_1=4\pi m+\mu,\qquad \sigma=\rho_2-4\pi V>0.
\]
We assume $m>0$. Throughout this section the weights are smooth and positive, and $\rho_2\notin4\pi\N$ is fixed. Fix a neighborhood of a solution of the shadow system in which the centers stay separated.

Assume that the scalar linearization at fixed centers is invertible throughout this neighborhood. The shadow scalar branch and its density are
\begin{equation}\label{reg:shadow-scalar}
 H_P=h_2e^{-4\pi\sum_iG(\cdot,p_i)},\quad
 n_*^g=\frac{\sigma H_Pe^{w(P)}}{\int_MH_Pe^{w(P)}},\quad
 -\Delta_gw=2(n_*^g-\sigma),\quad\int_Mw=0.
\end{equation}
Its Morse index is denoted by $m_S$. Write $F=(F^A,F^B)$ for the center field after eliminating $w$; thus
\begin{align}
 F_i^A&=\nabla(\log h_1+4\pi r_g-w/2)(p_i)
       +8\pi\sum_{l\ne i}\nabla_xG(p_i,p_l)
       +12\pi\sum_j\nabla_xG(p_i,q_j),\label{reg:forces}\\
 F_j^B&=\nabla(2\log h_1+\log h_2+12\pi r_g)(q_j)
       +12\pi\sum_i\nabla_xG(q_j,p_i)
       +24\pi\sum_{l\ne j}\nabla_xG(q_j,q_l).\nonumber
\end{align}
We order the scalar equation before the center equations. At a nondegenerate shadow solution $Z$, its degree is $d_S(Z)=(-1)^{m_S}\sgn\det DF(Z)$.
\subsection{Local profiles and the identities for the blow-up parameters}
At a regular blow-up point, the possible local mass pairs in units of $4\pi$ are $(1,0),(0,1),(1,2),(2,1),(2,2)$ by Lemma~\ref{fnd:local-unitary}.

By Proposition~\ref{fnd:signed-compactness}, one component concentrates completely. Since $\rho_2\notin4\pi\N$, this is the first component. In the notation
\[
 (s,t)=(M_1,M_2)/(4\pi),
\]
the second normalization is finite. Applying Jensen's inequality to the circular mean as in that proposition gives
\[
 2t-s\le1
\]
by integrability of the residual density. This argument uses integrability alone at the blow-up point. Thus the local masses are $(4\pi,0)$ at an ordinary point and $(8\pi,4\pi)$ at a mixed point. The total first mass is $4\pi m$.

After subtracting the logarithmic Green potentials of these masses, the equation for the residual second component has a smooth local coefficient, with a quadratic zero at an ordinary point and a positive coefficient at a mixed point.

The remaining Green potential has arbitrarily small mass on small disks. The local exponential integrability estimate and then elliptic regularity give a smooth local scalar remainder. In particular, the residual density is smooth and positive at a mixed point.

In a source-free conformal disk write $v_i=u_i+\log h_i+\varphi+\log(\rho_i/\int_Mh_ie^{u_i})$, so that $e^{v_i}dA$ is the $i$th measure in these coordinates. Take $\varphi(0)=0$ and use $-\Delta=e^\varphi(-\Delta_g)$ to obtain the canonical equation
\begin{equation}\label{reg:canonical}
 -\Delta v=K e^v+f,\qquad
 f_i=-\Delta(\log h_i+\varphi)-(K\rho)_ie^\varphi,
 \qquad K=\begin{pmatrix}2&-1\\-1&2\end{pmatrix}.
\end{equation}
Here $f$ is evaluated at $\rho_1$. At $\rho_1^0=4\pi m$, the forcing is
\begin{equation}\label{reg:wall-forcing}
 f^0=f+K(\mu,0)^{\mathsf T}e^\varphi
     =f+(2\mu,-\mu)^{\mathsf T}e^\varphi.
\end{equation}
We use the exterior solution at $\rho_1^0$ to construct the approximation and use $f$ in the equation at $\rho_1$.
\begin{proof}[Proof of Lemma~\ref{reg:shadow-limit}]
\emph{Step 1. The scalar equation.} Write $\nu_i$ for the sum of the point masses in component $i$, and $g$ for the residual density in the second component, with $\int_M g=\gamma$. The mean-zero Green representation gives
\[
 u^0=K(G\nu_1,G\nu_2)^t+(-w/2,w)^t,\qquad w=2Gg .
\]
Substitution in the second density gives
\[
 g=\frac{\gamma H e^w}{\int_MHe^w},\qquad H=h_2\exp(-G\nu_1+2G\nu_2).
\]
At each prescribed point $p$, the local weight $H$ is comparable to $\lvert z\rvert^{2\beta_p}$ with $\beta_p>-1$. At free ordinary and mixed centers $H$ is smooth with a quadratic zero or a positive value, respectively. The residual measure has arbitrarily small mass on sufficiently small disks.

The exponential estimate for its Green potential and H\"older's inequality, with $H$ integrable to a power greater than one, give $g\in L^p$ for some $p>1$. Hence $w$ is continuous on $M$ and smooth near each free center.

\emph{Step 2. The center equations.} On a source-free isolating disk use the canonical fields of \eqref{reg:canonical}. Set
\[
 a=(2v_1+v_2)/3,\qquad b=(v_1+2v_2)/3,
\]
and put $k_1=-(2f_1+f_2)/12$, $k_2=-(f_1+2f_2)/12$. Then
\[
 \bar\partial a'=k_1-e^{v_1}/4,\qquad \bar\partial b'=k_2-e^{v_2}/4,
\]
where $'=\partial_z$. On differentiating the squares and the product in $W_2$, the coefficients of $\bar\partial a'$ and $\bar\partial b'$ are $-2a'+b'$ and $a'-2b'$. Because $v_1'=2a'-b'$ and
\[
 v_2'=2b'-a',
\]
the exponential terms from these coefficients cancel $-(e^{v_1})'/4-(e^{v_2})'/4$ from the second derivatives. This yields
\[
 W_2=a''+b''-(a')^2+a'b'-(b')^2,\qquad
 \bar\partial W_2=k_1'+k_2'+(-2a'+b')k_1+(a'-2b')k_2.
\]
Its right-hand side has integral $O(d)$ on $D_d$, uniformly along the sequence, since $k_i$ and their first derivatives are bounded and
\[
 \sup_\zeta\int_{D_d}|z-\zeta|^{-1}\,dA_z\le Cd
\]
for the Cauchy fields of the bounded-mass measures. For every fixed small radius $d$ the enclosing circle lies in an annulus away from the limiting concentration points and the marked points.

The local potential bounds and interior elliptic estimates give convergence there of the first and second derivatives of the fields. The normalizing constants have zero spatial derivatives. The limiting fields have
\[
 a'=-x/z+A_0+O(|z|),\qquad b'=-y/z+B_0+O(|z|)
\]
for masses $(x,y)$ in $4\pi$ units. The simple-pole coefficient is $(2x-y)A_0+(2y-x)B_0$. Integration over the disk gives
\[
 \oint_{\partial D_d}W_2\,dz=2i\int_{D_d}\bar\partial W_2\,dA.
\]
Pass first to the sequence limit on the fixed circle, preserving the uniform $O(d)$ bound.

Next let $d$ decrease to zero: the double pole has zero circular integral, the bounded remainder contributes $O(d)$, and the simple pole contributes $2\pi i$ times its coefficient. That coefficient therefore vanishes. It gives
\[
 2A_0-B_0=0
\]
at $(1,0)$ and $3A_0=0$ at $(2,1)$.

The ordinary identity is the vanishing of the regular derivative of $v_1$ and retains $-\nabla w/2$. The mixed identity is that of $2v_1+v_2$, in which the term from the residual density $2(-w/2)+w$ vanishes.

Twice the conjugate complex derivatives give the coordinate gradients. To pass to the intrinsic equations multiply by $e^{-\varphi}$ at the respective center. Write
\[
 G(z,p)=-(2\pi)^{-1}\log|z-p|+R(z,p).
\]
The contribution of the regular part of the Green function at the ordinary point is $8\pi\nabla_xR(p,p)+\nabla\varphi(p)$, which is $4\pi\nabla r_g(p)$ by symmetry of $R$ and
\[
 r_g=R(p,p)+\varphi(p)/(4\pi).
\]
In the mixed equation the same contribution is $24\pi\nabla_xR(q,q)+3\nabla\varphi(q)$, hence $12\pi\nabla r_g(q)$. At other free centers the ordinary and mixed mass pairs in $4\pi$ units are $(1,0)$ and $(2,1)$.

The first row of $K$ applied to them gives $2$ and $3$, and $(2,1)K=(3,0)$ applied to them gives $3$ and $6$. Hence the ordinary row has interactions $8\pi$ and $12\pi$, while the mixed row has $12\pi$ and $24\pi$, as in \eqref{reg:forces}.

At a prescribed point of mass pair $(s,t)$, absorb $4\pi(2s-t)G$ and $4\pi(-s+2t)G$ into the two logarithmic weights, respectively. The resulting effective weight parameters are
\[
 (\alpha_{p,1}-2s+t,\ \alpha_{p,2}+s-2t).
\]
Substitution of $(0,0)$, $(A_p,0)$ and $(A_p+B_p,B_p)$ yields \eqref{fnd:effective}. The balances follow before any shape or scalar inverse is used.
\end{proof}
For $R>0$, $\beta>0$ and $B\in\C$, put
\begin{align}
 F_R&=R^4+\beta R^4|z|^2+|z^2-R^2B|^2,\nonumber\\
 Q_R&=R^4\{\beta R^4+4|z|^2+\beta|z^2+R^2B|^2\},\nonumber\\
 n_1&=4Q_R/F_R^2,\qquad n_2=16\beta R^8F_R/Q_R^2,
 \qquad a^2=\beta R^4(1+|B|^2)/4.\label{reg:gram}
\end{align}
For $g=(R^2,\sqrt\beta R^2z,z^2-R^2B)$, one has $F_R=\|g\|^2$, $Q_R=\|g\wedge g'\|^2$, and $|\det(g,g',g'')|^2=4\beta R^8$. Differentiating the norm squares gives
\[
 \Delta\log F_R=4Q_R/F_R^2=n_1,\qquad \Delta\log Q_R=4F_R|\det(g,g',g'')|^2/Q_R^2=n_2.
\]
Consequently $-\Delta\log n=K n$. Both norm squares have quartic growth, so the boundary flux of each logarithm at infinity is $8\pi$. Thus $\int_\C n_1=\int_\C n_2=8\pi$. The local limits, for fixed $\beta>0$ and $B$, are
\[
 R^2 n_1(Ry)\longrightarrow
 m_B(y):=\frac{16|y|^2}{(1+|y^2-B|^2)^2},
 \qquad
 a^2 n_2(ay)\longrightarrow\frac{4}{(1+|y|^2)^2},
\]
locally uniformly on the respective scaled planes. On compact subsets away from $z=0$, $n_2(z)\to d_0(z)$, where
\[
 d_0=16\beta/(4+\beta|z|^2)^2.
\]
Denote the corresponding outer and inner profiles by $m_{R,B}=R^{-2}m_B(z/R)$ and $m_a=4a^2/(a^2+|z|^2)^2$. The local masses are $8\pi,4\pi$, although the entire exact family has masses $8\pi,8\pi$. Set
\[
 \mathcal M(B)=\int_\C\frac{|w|}{(1+|w-B|^2)^2}\,dw,
 \qquad \mathcal H(B)=\mathcal M(B)-B\cdot\nabla\mathcal M(B).
\]
The map $\nabla\mathcal M$ is an orientation-preserving diffeomorphism onto the open disk of radius $\pi$. Indeed,
\[
 D^2\mathcal M(B)[v,v]
 =\int_\C\frac{|v|^2-(v\cdot\widehat{w+B})^2}{|w+B|}
                       \frac{dw}{(1+|w|^2)^2}>0,
\]
and its radial gradient tends to $\pi$ at infinity. Integration by parts also gives
\begin{equation}\label{reg:H-positive}
 \mathcal H(B)=\frac12\int_\C
             \frac{dw}{(1+|w|^2)|w+B|}>0.
\end{equation}
The change of variable $w=z^2$ gives the exact moments
\begin{equation}\label{reg:moments}
 \int m_{R,B}=8\pi,\quad
 \int z^2m_{R,B}=8\pi R^2B,\quad
 \int |z|^2m_{R,B}=8R^2\mathcal M(B).
\end{equation}
The smooth part of the first component at $q_j$ is
\[
 \psi_j=\log h_1+\varphi+12\pi R_g(\cdot,q_j)
       +8\pi\sum_iG(\cdot,p_i)
       +12\pi\sum_{l\ne j}G(\cdot,q_l)-w/2.
\]
Here $R_g$ is the regular part of the Green function in these coordinates. Define
\begin{equation}\label{reg:kappa-T}
 \kappa_j=\Delta_g\log h_1(q_j)+8\pi m-\rho_2-2K_g(q_j),
 \quad T_j=\psi_{j,zz}(q_j)+\psi_{j,z}(q_j)^2,
 \quad\mathbf t_j=(\Ree T_j,-\Imm T_j).
\end{equation}
In normalized coordinates at a mixed center, $\kappa_j=-f_1^0(q_j)$ by \eqref{reg:canonical} and \eqref{reg:wall-forcing}. The trace of the full first background there is instead
\[
 \Delta_g\psi_j(q_j)=\kappa_j+n_*^g(q_j),
\]
because $-\Delta_g w/2=n_*^g-\sigma$. In particular the trace contains the residual density $n_*^g(q_j)$. To check the effect of a coordinate change, let
\[
 z=f(\zeta)=
a\zeta+b\zeta^2+c\zeta^3+O(\zeta^4),
\]
$a\ne0$, be a holomorphic coordinate transition centered at the mixed point. Subtracting the leading pole $-6\log|\zeta|$ from the transformed canonical field gives
\[
 \widetilde\psi(\zeta)=
 \psi(f(\zeta))+\log|f'(\zeta)|^2
                    -6\log|f(\zeta)/\zeta|.
\]
In particular, $\widetilde\psi_\zeta(0)=a\psi_z(0)-b/a$ and $\widetilde\psi_{\zeta\zeta}(0) =a^2\psi_{zz}(0)+2b\psi_z(0)-(b/a)^2$. One may take $b=a^2\psi_z(0)$ to remove the first derivative. Adding the square of the first derivative gives $\widetilde T=a^2T$. Since
\[
 e^{\widetilde\varphi(0)}=|a|^2e^{\varphi(0)},
\]
the inequality $\kappa>4e^{-\varphi(0)}|T|$ is independent of the coordinate.
\begin{lemma}[Selection and scalar limits on a regular disk]
\label{reg:elementary-selection}
On a disk $D$ let $v^k$ satisfy
\[
 \Delta v_i^k+\sum_j K_{ij}e^{v_j^k}=b_i^k,
 \qquad b_i^k\ \text{smooth and uniformly bounded}.
\]
Assume the following.
\begin{itemize}
\item[\textup{(i)}]
The total masses of $e^{v_i^k}$ and the oscillations of $v_i^k$ on $\partial D$ are uniformly bounded.
\item[\textup{(ii)}]
The densities are bounded on closed subsets away from $0$. At least one local mass at $0$ is positive, and the second local mass is at most $4\pi$.
\end{itemize}
There exist finitely many centers $x_{\ell,k}$ and scales $\lambda_{\ell,k}$ with the following properties on every fixed smaller disk $D'\Subset D$.
\begin{itemize}
\item[\textup{(a)}]
For $x\in D'$, we have
\[
 v_i^k(x)+2\log\operatorname{dist}(x,\{x_{\ell,k}\}_\ell)\leq C
\]
\item[\textup{(b)}]
For each $j\ne\ell$, $\lambda_{\ell,k}=o(|x_{\ell,k}-x_{j,k}|)$. On rescaling about $x_{\ell,k}$, one density converges to an ordinary scalar profile of mass $4\pi$ and the other to zero.
\item[\textup{(c)}]
The oscillations of $v_i^k$ on $\partial D_r(x_{\ell,k})\subset D'$ are uniformly bounded whenever $4r\leq |x_{\ell,k}-x_{j,k}|$ for every $j\ne\ell$.
\end{itemize}
\end{lemma}
\begin{proof}
\smallskip\noindent\textit{Step 1. Selection of the centers.}\quad Denote the centers already selected by $\Sigma_k$. In selecting the next center we also take the boundary into account. If the claimed upper bound fails, take a point $x$ where
\[
 \max_i v_i^k(x)+2\log d\to+\infty,
\]
with $d=\operatorname{dist}(x,\Sigma_k\cup\partial D)$. Maximize $\max_i v_i^k(z)+2\log(d/2-|z-x|)$ on $D_{d/2}(x)$, obtaining a point $z_k$ and a width $\lambda_k=e^{-\max_i v_i^k(z_k)/2}$. The maximum property gives
\[
 \lambda_k=o(d/2-|z_k-x|).
\]
On every fixed disk of the $\lambda_k$ rescaling both shifted fields are bounded above and one has value zero at the center.

After subtracting the Dirichlet potential of their bounded right-hand sides, Harnack's inequality for the remaining harmonic functions shows that each component has a finite limit or tends uniformly to $-\infty$ on compact subsets.

The forcing tends to zero after rescaling. Any further point $x$ chosen in this way satisfies
\[
 |x-x_{\mathrm{old},k}|/\lambda_{\mathrm{old},k}\to\infty
\]
for each previous center, by its rescaled upper bound and $d\leq|x-x_{\mathrm{old},k}|$. The triangle inequality gives $|z_k-x_{\mathrm{old},k}|\geq|x-x_{\mathrm{old},k}|/2\geq d/2$. Thus both $\lambda_{\mathrm{old},k}$ and $\lambda_k$ are $o(|z_k-x_{\mathrm{old},k}|)$.

\smallskip\noindent\textit{Step 2. The scalar limits.}\quad Rescale the Green representation. The gradient of its regular harmonic part tends to zero.

On fixed compact sets the part of the singular gradient arising from outside $D_L$ is $O(L^{-1})$. This gives the representation in Lemma~\ref{fnd:entire-profiles}.

There is no point mass at the rescaled origin because the rescaled densities are bounded on compact sets. If both entire limits survive, Lemma~\ref{fnd:entire-profiles} with
\[
 A=B=1
\]
gives mass $8\pi$ in each component. This exceeds the available second mass. For a single entire component the scalar assertion of that lemma with $q=1$ gives the stated ordinary profile of mass $4\pi$.

Its translated center in the blown-up coordinates is zero: the boundary-distance factor used in the maximization has ratio $1+o(1)$ between the maximizer and any point at a bounded blown-up distance, so the limiting surviving field attains its maximum at zero.

Its value there is zero by the rescaling convention. Thus the scalar profile at this scale is $16/(4+|y|^2)^2$, and its width in the convention $m_s=4s^2/(s^2+|z|^2)^2$ is asymptotic to $2\lambda_k$.

This positive mass and the uniform mass bound stop the selection after finitely many steps. The lower bounds on the separations ensure that the different extracted masses are counted in disjoint disks.

\smallskip\noindent\textit{Step 3. The oscillation on circles.}\quad For the circle estimate use the Green representation on the original disk. Take two points of the circle of radius $r$.

On a neighborhood of radius comparable to $r$ about either point, selection bounds the density by $Cr^{-2}$. A change of variables by $r$ bounds the contribution of the logarithmic kernel difference there uniformly.

On the complement of those neighborhoods the logarithmic ratios are bounded, so the mass bound controls their integrals. The original harmonic part has bounded oscillation by the boundary hypothesis. This bounds the oscillation on the circle.
\end{proof}
\begin{lemma}[An identity for the shape parameter]
\label{reg:stress}
Let $0$ be an isolated mixed blow-up point of a sequence of normalized global solutions as above. Assume that the residual density is smooth and positive. Use the conformal coordinate with $\varphi(0)=0$.
\begin{itemize}
\item[\textup{(i)}]
Two separate ordinary bubbles in the first component and one in the second component require $|\kappa|=4|T|$.
\item[\textup{(ii)}]
Suppose $\kappa>4|T|$. After passing to a subsequence, take the selected inner center as origin. The outer and inner scales satisfy
\[
 R_n\to0,\qquad a_n=o(R_n).
\]
The rescaled densities corresponding to the local masses converge in $L^1$ on expanding disks to
\[
 m_B(y),\qquad \frac{4}{(1+|y|^2)^2},
\]
respectively, with bounded shape parameter $B$. Choose the disks by diagonal truncation of the two local masses. On the outer scale, use local convergence away from $0$ and exhaustion of the first local mass.
\item[\textup{(iii)}]
Under the assumption in \textup{(ii)}, $B$ satisfies
\begin{equation}\label{reg:shape-limit}
 \kappa\nabla\mathcal M(B)+4\pi(\Ree T,-\Imm T)=0.
\end{equation}
The solution is unique, and its determinant of the real derivative is positive.
\end{itemize}
\end{lemma}
\begin{proof}
\emph{Step 1. The selected centers.} Fix a disk containing only the blow-up point $0$. Denote this disk by $D$. For the fields $v_i^{(n)}$ of
\eqref{reg:canonical}, put $F_i^{(n)}=e^{v_i^{(n)}}$. The forcing is
uniformly smooth. The Green representation and the density bounds give
\[
 \operatorname{osc}_{\partial D}v_i^{(n)}
 +\int_D F_i^{(n)}\,dA+\sup_E F_i^{(n)}\le C_E
 \qquad(E\Subset D\setminus\{0\}).
\]
Let $(M_1,M_2)=(8\pi,4\pi)$ be the local mass pair at $0$.
Any sequence of disks $B_n$ shrinking to $0$ is eventually contained
in each fixed $D_\varepsilon(0)$. Hence
\[
 \limsup_{n\to\infty}\int_{B_n}F_i^{(n)}\,dA
 \le\lim_{\varepsilon\downarrow0}\limsup_{n\to\infty}
       \int_{D_\varepsilon(0)}F_i^{(n)}\,dA=M_i.
\]
Since $M_2=4\pi$, Lemma~\ref{reg:elementary-selection} gives the
selected set $\Sigma_n$, the inverse-square estimate and the circular
oscillation bound. Let $N_i$ count the cores active in component $i$.
Each carries $4\pi$, so disjointness gives, after a subsequence,
\[
 N_1\le2,\qquad N_2\le1,\qquad \#\Sigma_n=N_1+N_2\le3.
\]

If $\#\Sigma_n\ge2$, take a closest pair $p_n,q_n\in\Sigma_n$ and set
\[
 \delta_n=|p_n-q_n|\to0,\qquad
 V_{i,n}(y)=v_i^{(n)}(p_n+\delta_n y)+2\log\delta_n,\qquad
 Z_n=(\Sigma_n-p_n)/\delta_n.
\]
Distinct points of $Z_n$ are at least one unit apart. After a subsequence
let $Z$ be the set of their finite limits. The closest pair gives
$\#Z\ge2$. On each compact $E\subset\mathbb R^2\setminus Z$,
Lemma~\ref{reg:elementary-selection}\textup{(a)} gives
\[
 \sup_E e^{V_{i,n}}\le C_E.
\]

Fix $y_0\notin Z$. Subtract the Green representation at $y_0$ from
that at $y$, both in rescaled coordinates. For $y\in E$ and $L$
sufficiently large, the logarithmic kernel satisfies
\[
 \sup_{|\eta|>L}
 \left|\log\frac{|y-\eta|}{|y_0-\eta|}\right|\le\frac{C_{E,y_0}}{L}.
\]
The total masses are bounded, and the densities are bounded near $E$
and $y_0$. Boundary oscillation bounds the gradient of the harmonic
part on a smaller original disk. Thus
\[
 \sup_{y\in E}|V_{i,n}(y)-V_{i,n}(y_0)|\le C_{E,y_0}.
\]
The rescaled forcing is $O(\delta_n^2)$. Interior estimates give local
$C^2$ convergence of the fields after subtraction at $y_0$, along a
subsequence. Passing to limits at $y_0$ gives, on $\mathbb R^2\setminus Z$,
\[
 V_{i,n}\longrightarrow V_i\quad\text{in }C^2_{\rm loc},
 \qquad\text{or}\qquad
 V_{i,n}\longrightarrow-\infty\quad\text{locally uniformly}.
\]
Let $\mu_j$ be a local weak limit of $e^{V_{j,n}}\,dA$ along a common
subsequence. For a finite limit $V_i$ the Green representation holds
distributionally across $Z$, since
\[
 \sup_{\eta\in\mathbb R^2}
 \int_{D_\varepsilon(\eta)}\bigl|\log|y-\eta|\bigr|\,dA_y
 =O\bigl(\varepsilon^2(1+|\log\varepsilon|)\bigr).
\]
In particular, \eqref{reg:canonical} gives
\[
 -\Delta V_i=\sum_{j=1}^2K_{ij}\mu_j
 \qquad\text{in }\mathcal D'(\mathbb R^2).
\]

Fix a small rescaled disk isolating $z\in Z$. The masses and smooth
forcing are bounded, and the Green estimate bounds the boundary
oscillations. With no Dirac source,
Lemma~\ref{fnd:local-unitary}\textup{(iii)} gives
\[
 x(z):=\frac1{4\pi}\bigl(\mu_1(\{z\}),\mu_2(\{z\})\bigr)
 \in\{(0,0),(1,0),(0,1),(1,2),(2,1),(2,2)\}.
\]
Each such center has another at distance $O(\delta_n)$, so its core
width is $o(\delta_n)$ by Lemma~\ref{reg:elementary-selection}\textup{(b)}.
Its core contributes $4\pi$. The shrinking-disk bound thus gives
\[
 x(z)\ne(0,0)\quad(z\in Z),\qquad
 \sum_{z\in Z}x(z)\le(2,1)\quad\text{coordinatewise}.
\]
Since $\#Z\ge2$, these relations leave only
\[
 x(z)\in\{(1,0),(0,1)\}\qquad(z\in Z).
\]
Suppose that exactly two selected centers remain bounded on this scale.
At either point $z$, let $i$ be its active component. If $V_i$ is finite
off $Z$, integration on $D_r(z)$ gives for its circular mean about $z$
\[
 r\frac{d}{dr}\overline V_i(r)
 =-\frac1{2\pi}\sum_{j=1}^2K_{ij}\mu_j(D_r(z))
 =-4+o(1)\qquad(r\downarrow0).
\] Jensen's inequality would give nonintegrable density there. Thus if the two point masses belong to different components, neither component has a finite limit.

If both belong to the first component, that component again tends to $-\infty$. A finite limit $V_2$ of the second component in this case would satisfy
\[
 \Delta V_2+2e^{V_2}=4\pi(\delta_{z_1}+\delta_{z_2}),\qquad
 \int_{\mathbb C}e^{V_2}\leq4\pi .
\]
The source signs follow directly from the off-diagonal entry of $K$. On large circles enclosing the two concentration points, flux gives
\[
 r\overline V_2'(r)=4-\pi^{-1}\int_{D_r}e^{V_2}\geq0,
\]
which again contradicts integrability by Jensen's inequality. Two point masses in the second component would exceed its total mass.

Hence in every two-center case both fields tend to $-\infty$ on a fixed rescaled circle enclosing the pair. The rescaled measures have no residual mass away from these points.

On such circles the first derivatives after multiplication by $K^{-1}$ converge to the Cauchy transforms of the point masses.

This follows from the gradient representation, using the bounded density and its convergence to zero near the contour and the $O(L^{-1})$ gradient bound for sources outside radius $L$.

The contour integral of $zW_2$ tends to zero: $\bar\partial (zW_2)$ contains only the scaled smooth forcing and terms at most first order in the fields, with bounded $L^1$ norm on fixed disks.

The coefficient of $z^{-1}$ in the limiting $zW_2$ is $x+y-x^2+xy-y^2$ by expanding the rational Cauchy fields at infinity, where $x,y$ are the two total masses in $4\pi$ units.

Thus the two enclosed masses satisfy the radial Pohozaev identity at a regular point. Neither $(2,0)$ nor $(1,1)$ satisfies it. This excludes two centers after rescaling by the minimum distance, including the possibility that a third center is at a larger distance.

It follows that in the case of multiple selected points there are exactly three, all pairwise distances are comparable to $\delta_n$, and the three widths are $o(\delta_n)$.

The first component has two pure cores, the second has one, and the three core profiles exhaust the mass of both components on the rescaled compact sets.

The second-component center is the midpoint of the first two to leading order. Here is a local way to see this without a rate for the heights.

Write $z_0$ for the limiting second-component position and $z_1,z_2$ for the other positions. On fixed rescaled contours, the derivatives of the inverse Cartan fields have limits
\[
 a'=-\frac1{z-z_1}-\frac1{z-z_2},\qquad
 b'=-\frac1{z-z_0}.
\]
Indeed the Cauchy kernels converge against the isolated masses off the observation contours. Near each contour the densities tend uniformly to zero by the selection bound, the compactness alternative and exhaustion.

The part of the gradient kernel beyond distance $L$ is bounded by $C/L$ using the mass bound, and the derivative of the harmonic term becomes $O(\delta_n)$. Interior elliptic estimates give the same contour convergence of the next derivatives. The rescaled expression
\[
 W_2=a''+b''-(a')^2+a'b'-(b')^2
\]
has contour integral tending to zero on a circle around $z_0$: its $\bar\partial$ derivative contains only the scaled smooth forcings and their first derivatives multiplied by terms of at most first order in the fields, as in the quadratic identity below.

Their integrals are $o(1)$ since the rescaled forcing is $O(\delta_n^2)$ and the first derivatives are bounded in $L^1$ on fixed disks by the Cauchy representation. At $z_0$ the only residue in the displayed limit is $-a'(z_0)$. Hence
\[
 \frac1{z_0-z_1}+\frac1{z_0-z_2}=0,
 \qquad z_0=\frac{z_1+z_2}{2}.
\]
This is the separated formation. It suffices below that the relative directions of its three ordered pairs agree up to sign.

\emph{Step 2. Rescaling about a single center.} Suppose instead that there is exactly one selected center $p_n$. Let $\ell_n$ be a radius much larger than the innermost scalar width, with both components fast on $\partial D_{\ell_n}(p_n)$ and with masses tending to the innermost mass pair.

Such a radius is obtained by a diagonal choice on the first rescaled plane. It can be chosen so that there is a shrinking radius $d_n\gg\ell_n$ on which the masses tend to $(2,1)$ in $4\pi$ units. To determine the next scale, put
\[
 t=\log r,\qquad
 M_i(t)=\frac1{4\pi}\int_{D_{e^t}(p_n)} F_i^{(n)},\qquad
 G_i(t)=\overline{v_i^{(n)}}(e^t)+2t .
\]
The overline denotes circular mean about $p_n$. The equations and the circular oscillation bound imply
\[
 G_i'(t)=2-2(KM(t))_i+O(e^{2t}),\qquad
 0\leq M_i'(t)\leq C e^{G_i(t)} .
\]
The constants are independent of $n$ throughout the shrinking disk. At $t=\log\ell_n$ both $G_i$ tend to $-\infty$. For the component of the inner bubble, $G_i'\to-2$ at
\[
 t=\log\ell_n.
\]
For the other component, $G_j'\to4$. Fix $\eta>0$ small. While $M$ stays within $\eta$ of its initial value and $r\le d_n$, we have
\[
 G_i'\le-c<0<c\le G_j'
\]
with $c$ independent of $n$. Choose $L$ so that $Ce^{-L}/c<\eta/4$. Until one of the $G_i$ first reaches $-L$, the masses added by the increasing and decreasing components are bounded by
\[
 Ce^{-L}/c,\qquad C e^{G_i(\log\ell_n)}/c=o(1),
\]
respectively. Integrating in $t$ gives these bounds. Thus the masses stay within $\eta$ of their initial values until one of the $G_i$ reaches $-L$. As the masses in $D_{d_n}$ tend to the full mixed masses, the increasing component must reach $-L$ at some first radius
\[
 R_n<d_n.
\]
The boundedness of the derivatives and $G_i(\log\ell_n)\to-\infty$ give $R_n/\ell_n\to\infty$.

Rescaling by $R_n$ at $p_n$, use the same bounded field differences on compact subsets off $0$ as after rescaling by the minimum distance.

On the unit circle the decreasing field has mean tending to $-\infty$ and the other has mean $-L$. Only the increasing component has a finite limit, and only the innermost bubble contributes a point mass. In fact the decreasing component adds $o(1)$ up to $R_n$.

For any fixed $s\in(0,1)$, integrating the positive derivative from $sR_n$ to $R_n$ bounds the new mass of the increasing component inside $D_{sR_n}(p_n)$ by $Cs^c$ for some $c>0$, in addition to its $o(1)$ mass inside $D_{\ell_n}(p_n)$.

The rescaled equations and the Green representation used above therefore give
\[
 \Delta V+2e^V=4\pi\delta_0\quad\text{in }\mathbb C,
 \qquad \int_{\mathbb C}e^V<\infty
\]
for the increasing component. The signs and coefficient of the Dirac mass follow from its coupling to the opposite point mass $4\pi$.

The preceding rescaled Green representation for the surviving field has no remaining harmonic gradient, because the regular part of the original kernel is bounded on the isolating disk and mass escaping beyond $D_L$ gives an $O(L^{-1})$ tail on each fixed compact set.

Subtract its prescribed logarithm at zero; the residual measure has no point mass there and the cumulative mass bound
\[
 \int_{D_s}e^V=O(s^c)
\]
in the preceding construction gives
\[
 s\overline V'(s)=2+O(s^c)
\]
by radial flux. Its angular oscillation stays bounded by the circle estimate. Thus $V-2\log|y|$ is bounded, so $e^V=O(|y|^2)$ locally. Subtracting its continuous Dirichlet potential removes the remaining bounded harmonic singularity. Lemma~\ref{fnd:entire-profiles} with $q=2$ gives
\[
 \int e^V=8\pi,
\]
and, after absorbing a fixed dilation into $R_n$, $e^V=m_B$ for a finite $B\in\mathbb C$. In particular the increasing component must be the first: otherwise its mass would exceed $4\pi$. The inner scalar limit is thus the ordinary second bubble about $p_n$, and its width $a_n$ is $o(R_n)$.

The profiles exhaust both local masses. Together with the separated case this proves the qualitative alternatives.

All the rescalings here use the local equations on the isolating disk, including when prescribed singularities are present away from that disk. No comparison of the heights to $\log \delta_n^{-1}$ is used.

\emph{Step 3. The identity for $W_3$.} To exclude the first alternative, put $a=(2v_1+v_2)/3$, $b=(v_1+2v_2)/3$, and
\[
 k_1=(2\kappa_1^c+\kappa_2^c)/12,
\]
$k_2=(\kappa_1^c+2\kappa_2^c)/12$, where $\Delta v+Ke^v=\kappa^c$. For primes denoting $z$ derivatives,
\[
 W_3=a'''-2a'a''+a'b''+(a')^2b'-a'(b')^2,
\]
\[
 \bar\partial W_3=k_1''-2a'k_1'+a'k_2'
  +k_1(-2a''+b''+2a'b'-(b')^2)
  +k_2((a')^2-2a'b').
\]
All exponential terms cancel by $\bar\partial a'=k_1-e^{v_1}/4$ and $\bar\partial b'=k_2-e^{v_2}/4$. For $W_2=a''+b''-(a')^2+a'b'-(b')^2$, the same differentiation gives
\[
 \bar\partial W_2=k_1'+k_2'+(-2a'+b')k_1+(a'-2b')k_2.
\]
On a disk of radius $d$ the integral of the right-hand side is $O(d)$, uniformly along the sequence, by the Cauchy representation and
\[
 \sup_\zeta\int_{|z|<d}|z-\zeta|^{-1}\,dA_z\le Cd.
\]
In the exterior limit write $a'=-2/z+A_0+O(|z|)$ and $b'=-1/z+B_0+O(|z|)$. The constant term of the integral of $W_2\,dz$ on the shrinking circle is $2\pi i(3A_0)$. Integration over the disk gives $A_0=0$, or $2\psi_{1,z}+\psi_{2,z}=0$. With
\[
 A_1=(2\psi_{1,zz}+\psi_{2,zz})/3,\qquad B_1=(\psi_{1,zz}+2\psi_{2,zz})/3,
\]
the corresponding coefficient of $W_3$ is $4A_1-2B_1+2B_0^2=2T$. Terms containing $\bar z$ in the Taylor expansion at this order have zero circular integral.

For each component let $\mu_i^{(n)}$ be the density measure restricted to a shrinking disk whose mass tends to the local mass. The truncation is chosen after the profile limits just obtained. The angular integrals
\[
 S_{ij}=\lim_{n\to\infty}\iint
 \Theta(x-y)\,d\mu_i^{(n)}(x)d\mu_j^{(n)}(y),\qquad
 \Theta(z)=\overline z/z ,
\]
may be taken along a subsequence. Their integrands are bounded in absolute value by one and their values on the diagonal are irrelevant for the original smooth densities.

To integrate the identity for $W_3$, let $\chi_d$ be a radial smooth cut-off equal to one in $D_d$ and supported in $D_{2d}$. We first let $n\to\infty$ with $d$ fixed and then let $d\downarrow0$. On the left,
\[
 \int\chi_d\bar\partial W_3=-\int(\bar\partial\chi_d)W_3.
\]
Only the coefficient of $1/z$ in the pole expansion on the limiting annulus survives the angular integration. Its contribution is $2\pi T$, since
\[
 \bar\partial(1/z)=\pi\delta_0
\]
and that coefficient is $2T$. Remainders that are bounded on the shrinking annulus give $O(d)$. The zeroth- and first-gradient terms on the right vanish, by boundedness of the smooth coefficients and the Cauchy $L^1(D_{2d})$ bound $O(d)$ for the first fields.

For a second derivative term such as $\int\chi_d k_1 a''$, integrate by parts using $\partial_z$. The term with $\partial_z k_1$ is $O(d)$. For the term with $\partial_z\chi_d$, replace $k_1(z)$ by $k_1(0)$ at cost $O(d)$. On the annulus, after the first limit,
\[
 a'=-2/z+O(1),\qquad \int (\partial_z\chi_d)/z=0
\]
by radiality. The bounded remainder again gives $O(d)$. The $b''$ term is handled by its pole $-1/z$ and the same angular integral. For a smooth compactly supported $k$, set
\[
 \mathcal Ck(x)=\frac1\pi\int_{\mathbb C}
                      \frac{k(z)}{x-z}\,dA_z .
\]
In the quadratic terms, the singular part of each first field is $-(4\pi)^{-1}$ times the integral of $(z-x)^{-1}$ against its density measure. Products involving a bounded smooth part give $O(d)$ in the cut-off pairing. For the remaining pair kernels use
\[
 \int_{\mathbb C}\frac{k(z)}{(z-x)(z-y)}\,dA_z
      =-\pi\,\frac{\mathcal Ck(x)-\mathcal Ck(y)}{x-y}.
\]
For $k=\chi_d k_\ell^{(n)}$, the preceding identity involves the divided difference
\[
 \left|\frac{\mathcal Ck(x)-\mathcal Ck(y)}{x-y}\right|\le C,
\]
with $C$ independent of $n,d$. Scaling to the unit cut-off gives a uniform $C^2$ bound on the coefficients. The kernel estimate bounds the Cauchy transform and its first derivatives after subtraction at the singularity of the differentiated kernel.

At fixed $d$, take the smooth coefficient limit and use
\[
 \frac{\mathcal Ck(x)-\mathcal Ck(y)}{x-y}
   =(\mathcal Ck)_z(0)+k(0)\Theta(x-y)+o(1)
\]
uniformly as $x,y\to0$, since
\[
 \bar\partial\mathcal Ck=k.
\]
The term $(\mathcal Ck)_z(0)$ is $O(d)$: the constant part of $k_\ell$ pairs to zero with the differentiated radial cut-off kernel, and subtracting that constant leaves an $O(|z|)$ factor against $|z|^{-2}$.

Only pairs of point masses contribute to these limits. Fix $\varepsilon>0$ and separate the complement of $D_\varepsilon$ from the origin. If one argument of the divided difference lies outside $D_\varepsilon$, the bound is $O(d/\varepsilon)$ when the argument separation is at least $\varepsilon/2$, since
\[
 \|\mathcal Ck\|_\infty=O(d).
\]
Otherwise the entire connecting segment is away from $0$ and the bound is $O(d^2/\varepsilon^2)$ by differentiating the kernel there.

Inside $D_\varepsilon$, the mass remaining after truncation at the concentration scale tends to zero as $\varepsilon$ shrinks, and the divided difference is bounded uniformly. Thus the contribution of the residual density and all remaining density is zero. Combining the numerical factors gives
\[
 32\pi^2 T=-k_1(2S_{12}-S_{22})-k_2(S_{11}-2S_{12}),
\]
where $k_1,k_2$ are evaluated at the limiting center. Here
\[
 2k_1-k_2=\kappa_1^c(0)/4=\kappa/4,
\]
because $\kappa_1^c=-f_1^0$ at $\mu=0$ and
\[
 \varphi(0)=0.
\]
\emph{Step 4. The equation for $B$.} In each scalar core the density converges, after rescaling, to a radial profile in $L^1$ on the truncated expanding disks.

This follows from local convergence and exhaustion of the mass $4\pi$, using a diagonal cut-off. Simultaneous rotation in both variables gives zero for the angular integral of the profile against itself.

In the separated case take those cut-offs to have radii $o(\delta_n)$. On pairings between distinct cores the relative complex directions converge uniformly on the product of their supports, so the angular factors converge to their values at the limiting separated positions.

In the nested case, on the outer scale the first density converges to $m_B$ in $L^1$ by local convergence away from $0$ and exhaustion of its full mass. The concentrating part of the second measure tends to $4\pi\delta_0$ there.

In their cross pairing one can first remove the region in which the first argument has modulus less than $\varepsilon$, since the mass of this region under $m_B$ tends to zero. This justifies use of $\Theta(y)$ for the cross pairing on the outer scale.

The integral pairing the second density with itself is computed on the inner scale with its full radial profile. For the split formation, each radial density has zero pairing with itself, and the pairings of distinct densities give
\[
 S_{11}=S_{12}=32\pi^2\tau,
\]
$S_{22}=0$, $|\tau|=1$. Thus $T=-\kappa\tau/4$. Hence $\kappa>4|T|$ excludes the three separate bubbles. For the nested profile, direct substitution $w=y^2$ gives
\[
 S_{12}=32\pi\overline{\nabla\mathcal M(B)},\qquad S_{22}=0.
\]
The scalar identity for $A_B=-\partial_y\log(1+|y^2-B|^2)$, $A_B'=A_B^2+A_B/y$, gives $S_{11}=S_{12}$ after integration. More explicitly, the Cauchy representation and symmetrization give
\[
 S_{11}=-8\pi\int_{\C}\overline y\,m_B(y)A_B(y)\,dA_y.
\]
Since $\partial_y m_B=m_B(1/y+2A_B)$, integration by parts gives
\[
 0=\int_{\C}\overline y\,\partial_y m_B\,dA_y
   =\int_{\C}\frac{\overline y}{y}m_B\,dA_y
       +2\int_{\C}\overline y\,m_BA_B\,dA_y.
\]
The boundary terms vanish at zero and infinity. Thus
\[
 S_{11}=4\pi\int(\overline y/y)m_B=S_{12},
\]
including the pairing of the density with itself. This proves \eqref{reg:shape-limit}. Strict convexity and the gradient range prove its last assertions. Mass saturation and selection also give
\[
 r^2e^{v_i}\leq C
\]
and uniformly bounded circular oscillation on the annuli between the inner scale and a fixed radius.
\end{proof}
\subsection{Construction of the approximate solution}
Choose $c_1=2\log t$, $c_2=\log(\sigma/\int H_Pe^w)$. The exterior pair is
\[
 u_1^0=8\pi\sum_iG(\cdot,p_i)+12\pi\sum_jG(\cdot,q_j)-w/2,
 \qquad u_2^0=w-4\pi\sum_iG(\cdot,p_i).
\]
In a conformal coordinate centered at $p_i$, normalized by $\varphi(p_i)=0$, define the ordinary first-component background by
\[
 \psi_i^A(x)=\log h_1(x)+\varphi(x)+8\pi R_g(x,p_i)
       +8\pi\sum_{l\ne i}G(x,p_l)
       +12\pi\sum_jG(x,q_j)-\frac12w(x).
\]
At each mixed center choose conformal coordinates depending smoothly on the parameters while the centers remain separated. Starting from a fixed coordinate translated to the center, take the positive real factor $a$ with
\[
 |a|^2e^{\varphi(q_j)}=1,
\]
and then the quadratic coefficient $b=a^2\psi_{j,z}(q_j)$ of the transition in the calculation preceding Lemma~\ref{reg:stress}. In these coordinates $\varphi(0)=0$ and $\psi_{j,z}(0)=0$ for each parameter value. These backgrounds fix the scales by
\begin{equation}\label{reg:scales}
 4s_i^2=t^2e^{\psi_i^A(p_i)},\qquad
 16\epsilon_j^2=t^2e^{\psi_j(q_j)},\qquad
 R_j^2=\epsilon_j,\quad
 a_j^2=\beta_j\epsilon_j^2(1+|B_j|^2)/4,
 \quad\beta_j=n_*^g(q_j).
\end{equation}
Thus $s_i\asymp a_j\asymp t$ and $R_j\asymp\sqrt t$. At ordinary points add $(D_i,-D_i/2)$ to $u^0$, where $D_i=-2\log(1+s_i^2/|z|^2)$. At mixed points add $\log n-L$, with
\[
 L_1=2\log\epsilon_j+\log16-6\log|z|
                   +\log(1+\beta_j|z|^2/4),\quad
 L_2=\log\beta_j-2\log(1+\beta_j|z|^2/4).
\]
Let $\widehat u^{\rm app}$ be the function obtained by making these additions with fixed disjoint cut-offs, and set
\[
 u_i^{\rm app}=\widehat u_i^{\rm app}
                    -\int_M\widehat u_i^{\rm app}.
\]
The same constant $c_1$ is used at every first-component bubble. Define
\[
 N_i=e^{c_i}h_i e^{\widehat u_i^{\rm app}},\qquad
 M_i=\int_MN_i,\qquad V_i^{\rm app}=\frac{\rho_iN_i}{M_i}.
\]
\begin{lemma}[Residual and operator comparison]
\label{reg:residual}
Fix $1<p<2$ and let $|\mu|\leq C_0t$. If $V\geq1$, assume the $B_j$ are bounded.
\begin{itemize}
\item[\textup{(i)}]
We have
\[
 M_1=4\pi m+O(t),\quad M_2=\rho_2+O(t),\quad
 \|\mathcal R^{\rm app}\|_{H^{-1}}\leq C_pt^{2/p-1}
 \quad(1<p<2).
\]
\item[\textup{(ii)}]
The normalized approximate densities differ from the sum of the separated ordinary, outer, inner and residual model densities by $O(t^{2/p-1})$ in $L^p$.
\end{itemize}
\end{lemma}
\begin{proof}
\emph{Step 1. The equation for the approximation.} On a disk about a blow-up point, write $N_i^{\rm loc}=e^\varphi N_i$ and
\[
 \mathsf v_i^{\rm app}=\log N_i^{\rm loc}
\]
for the canonical approximation. In the following local formulas abbreviate $N_i^{\rm loc}$ to $N_i$. At an ordinary point $N_1=m_{s_i}e^{J_1}$ and $N_2=\gamma_i(x)(s_i^2+|x|^2)$, with
\[
 J_1=F_i^A\cdot x+O(|x|^2).
\]
Define the local residual before mass normalization to be $-\Delta\mathsf v^{\rm app}-KN-f$, with all densities expressed in the conformal coordinate. At an ordinary point it is
\[
 -K(m_{s_i}(e^{J_1}-1),s_i^2\gamma_i)^{\mathsf T}
       +K(\mu,0)^{\mathsf T}e^\varphi.
\]
At a mixed point the chosen coordinate gives $J_1=O(|x|^2)$, $J_2=F_j^B\cdot x+O(|x|^2)$, and the same residual is
\begin{equation}\label{reg:cancelled-residual}
 -K\binom{n_1(e^{J_1}-1)}{(n_2-d_0)(e^{J_2}-1)}
       +K(\mu,0)^{\mathsf T}e^\varphi.
\end{equation}
Indeed, the smooth correction at $\rho_1^0=4\pi m$ satisfies the exact identities
\[
 -\Delta J_1=f_1+2\mu e^\varphi-n_*^{\rm loc}+d_0,\qquad
 -\Delta J_2=f_2-\mu e^\varphi+2n_*^{\rm loc}-2d_0.
\]
The terms containing the residual density cancel, leaving the forcing $f$ at $\rho_1=4\pi m+\mu$. Replacing $N_i$ by
\[
 V_i^{\rm app}=\rho_iN_i/M_i
\]
adds exactly $-K((\rho_i/M_i-1)N_i)_i$ to this residual. Equivalently, if $\Pi$ denotes subtraction of the integral on $M$ and $N$ again denotes the density with respect to $dV_g$, we obtain
\begin{equation}\label{reg:normalized-residual}
 \mathcal R^{\rm app}
 =-\Delta_g u^{\rm app}-K\Pi N
       -K\Pi\bigl((\rho_i/M_i-1)N_i\bigr)_i.
\end{equation}
The identity holds for every parameter in the given interval. \emph{Step 2. The density estimates.} Write $A=R^2$, $u=|z|^2$,
\[
 H=A^2+|z^2-AB|^2\asymp A^2+u^2,\qquad x=u+a^2
\]
and $Y=Q_R/R^4$. Then $F_R=H+\beta A^2u$ and $Y-4u=O(H)$. Subtracting $16A^2u/H^2$ from $4A^2Y/F_R^2$ gives $|n_1-m_{R,B}|\le CA^2/H$, using $A^2u^2\le C H^2$. For the second fraction put $h=Y/(4x)$ and $h_0=1+\beta u/4$. Then
\[
 h-h_0=O(Au/x),\qquad h\asymp 1,
\]
and $S=F_R-A^2(1+|B|^2)-u^2=O(Au)$. The difference is
\[
 n_2-m_a-d_0
  =\frac{4a^2}{x^2}(h^{-2}-1)
   +\frac{\beta u^2}{x^2}(h^{-2}-h_0^{-2})
   +\beta h_0^{-2}(u^2/x^2-1)
   +\frac{\beta S}{x^2 h^2}.
\]
Each term has absolute value at most $CA/x$, by $a^2\asymp A^2$, $a^2u\le x^2/4$ and $1-u^2/x^2\le 2a^2/x$. For $|z|\le R$ one has $A/x\le C m_a$ and outside that set $A/x\le C A/(A+u)$. Polar integration and rescaling therefore give
\[
 \|n_1-m_{R,B}\|_p=O(R^{2/p}),\quad
 \|n_2-m_a-d_0\|_p=O(R^{4/p-2}),
\]
\[
 \||x|^2n_1\|_p=O(R^{2/p}),\quad
 \||x|(n_2-d_0)\|_p=O(R^{4/p-2}),\quad
 \||x|m_s\|_p=O(s^{2/p-1}).
\]
The linear terms in the mass integrals vanish by parity. The ordinary first mass error is $O(t^2|\log t|)$, its second excess $O(t^2)$, and each mixed mass error is $O(t)$. For the second mixed mass, integrate $n_2=\Delta\log Q_R$ on a fixed disk and use
\[
 Q_R/R^4=4|z|^2+\beta|z|^4+O(R^2)
\]
with the same error for the first derivatives on its boundary. Since
\[
 \Delta\log(4|z|^2+\beta|z|^4)=4\pi\delta_0+d_0,
\]
the disk mass differs from $4\pi+\int d_0$ by $O(R^2)$. Replacing $n_2$ by $n_2e^{J_2}$ and $d_0$ by $d_0e^{J_2}$ does not change this order: the linear term in the Taylor expansion integrates to zero against the even difference $n_2-d_0$. Its $|z|^2$-weighted absolute integral is bounded by
\[
 C\int_{D_d}|z|^2(m_a+A/(|z|^2+a^2))=O(a^2(1+|\log a|)+A)=O(R^2).
\]
The cut-off errors are $O(t^2)$ and $O(t)$ respectively. If $V=0$, parity removes the linear first mass term and $\int_{D_d}|z|^2m_{s_i}=O(t^2|\log t|)$ for each
\[
 s_i\asymp t.
\]
The first mass error is thus $O(t^2|\log t|)$, while the second mass error is $O(t^2)$ by the displayed $s_i^2\gamma_i$ term. On the fixed annuli the cut-off corrections and their first two derivatives are $O(t^2)$.

Outside the bubble disks the first unnormalized density is $O(t^2)$ and the second is $n_*^g$. On each disk the $L^p$ error is bounded using
\[
 \||z|m_{s_i}\|_p=O(t^{2/p-1}),\qquad |e^{J_1}-1|\leq C|z|.
\]
This bounds both the unnormalized density difference and the residual at $\mu=0$. It also bounds the $H^{-1}$ error, by $H^1\hookrightarrow L^{p'}$. The normalizing factors have errors $O(|\mu|+t^2|\log t|)$, so the same conclusions follow when all bubbles are ordinary. Hence
\[
 \log(\rho_i/M_i)=O(t).
\]
This factor multiplies an $L^p$ norm at most $Ct^{2/p-2}$. The remaining smooth parameter term is $O(|\mu|)$ in $H^{-1}$, which is smaller than the claimed bound on
\[
 |\mu|\leq C_0t.
\]
This proves the normalized estimates. More generally, at fixed approximation parameters the exact dependence is
\[
 \mathcal R^{\rm app}(\mu)-\mathcal R^{\rm app}(0)
       =-\mu K(\Pi(N_1/M_1),0)^{\mathsf T}.
\]
The probability density $P=N_1/M_1$ is bounded by $Ct^{-2}$. For its Green energy, the kernel on distances at least $t$ is bounded by $C(1+L)$, and on smaller distances its integral against $P$ is bounded by
\[
 Ct^{-2}\int_0^t(1+|\log r|)r\,dr=O(1+L).
\]
Thus the residual change has norm $O(|\mu|\sqrt L)$ in $H^{-1}$. The change in the second derivative is $-\mu$ times the bilinear form $\int uv\,dP-\int u\,dP\int v\,dP$ for the first component. With
\[
 p=1+1/L,\qquad \|P\|_p=O(1)
\]
and $\int_M Pv^2\le C\|v\|_{2(L+1)}^2\le CL\|v\|_{H^1}^2$ by Lemma~\ref{fnd:exponential}. Polarization gives $O(|\mu|L)$. These estimates will also be used for arbitrary blowing-up sequences, before the sharper relation between $\mu$ and $t$ is known.
\end{proof}
\begin{lemma}[Comparison of two bilinear forms]
\label{hyb:covariance}
Let $p_0>1$ and let $v_T,\widetilde v_T$ be probability densities satisfying
\[
 \|v_T-\widetilde v_T\|_1\le Ce^{-aT},\qquad
 \|v_T\|_{p_0}+\|\widetilde v_T\|_{p_0}\le Ce^{bT},
 \qquad a,b>0.
\]
For all mean-zero $\phi,\psi\in H^1(M)$,
\[
 |\Cov_{v_T}(\phi,\psi)-\Cov_{\widetilde v_T}(\phi,\psi)|
 \le CTe^{-aT}\|\phi\|_{H^1}\|\psi\|_{H^1}.
\]
\end{lemma}
\begin{proof}
For large $T$, put $p_T=1+1/T<p_0$. In the interpolation identity $1/p_T=\theta_T+(1-\theta_T)/p_0$ we have
\[
 1-\theta_T=\frac{p_0}{(p_0-1)(T+1)}.
\]
Hence interpolation gives
\[
 \|v_T-\widetilde v_T\|_{p_T}
 \le Ce^{-aT}\exp((a+b)T(1-\theta_T))\le C'e^{-aT}.
\]
Since both densities have mass one, each $L^{p_T}$ norm is at most $C\exp(bT(1-\theta_T))$. Since
\[
 p_T'=T+1,
\]
the two-dimensional Sobolev inequality
\[
 \|u\|_q\le C\sqrt q\|u\|_{H^1}
\]
bounds the product integral by $CTe^{-aT}\|\phi\|_{H^1}\|\psi\|_{H^1}$. Each weighted mean is bounded by $C\sqrt T\|u\|_{H^1}$, and its difference by $C\sqrt Te^{-aT}\|u\|_{H^1}$.

To subtract the products of the means, add and subtract the mixed product. Each term then contains one mean and one difference of means. This proves the estimate for these bilinear forms for probability densities that may vanish.
\end{proof}
In the following estimate
\[
 m_{q,R}(z)=4q^2R^{2q}|z|^{2q-2}(R^{2q}+|z|^{2q})^{-2}.
\]
\begin{lemma}[Centered estimates and interaction of the two scales]
\label{hyb:centered-estimates}
Let $q>0$ be fixed and let $\omega_{q,R}$ be the probability measure obtained by normalizing $m_{q,R}$ times a fixed cut-off equal to one near zero. The following bounds are uniform for small $R$.
\begin{itemize}
\item[\textup{(i)}]
For $L\ge2$ and $\overline v=\int v\,d\omega_{q,R}$,
\begin{align*}
 \int |v-\overline v|^2\,d\omega_{q,R}
     &\leq C\int |\nabla v|^2,\\
 \int_{\{|z|/R\notin[L^{-1},L]\}}
       |v-\overline v|^2\,d\omega_{q,R}
     &\leq C L^{-2q}(1+\log L)\int |\nabla v|^2.
\end{align*}
The bounds also hold for the mixed outer measures with $q=2$ and $B$ in a fixed compact set. For fields with bounded Dirichlet energy, local strong convergence after subtracting the means implies weighted $L^2$ convergence.
\item[\textup{(ii)}]
Suppose $r/R=\delta\to0$. If a bounded outer field $X$ satisfies $|X(y)-X(0)|\leq C\min\{|y|^\beta,1\}$, $\beta>0$, then
\[
 \int |X(z/R)-X(0)|^2\,d\omega_{t,r}
 \leq C\delta^{\min\{2\beta,2t\}}(1+|\log\delta|).
\]
\item[\textup{(iii)}]
For the same $r,R$, if a bounded inner field $Y$ satisfies $|Y(y)-Y(\infty)|\leq C\min\{1,|y|^{-2t}\}$, then
\[
 \int |Y(z/r)-Y(\infty)|^2\,d\omega_{s,R}
 \leq C\delta^{\min\{4t,2s\}}(1+|\log\delta|).
\]
\end{itemize}
The interaction estimates are uniform for the outer kernel functions and inner dilation functions, respectively. With inner decay $|y|^{-\lambda}$, replace $4t$ in \textup{(iii)} by $2\lambda$.
\end{lemma}
\begin{proof}
\smallskip\noindent\textit{Step 1. The weighted estimate.}\quad Scale $R$ to one and write
\[
 z=e^\xi e^{i\theta}.
\]
The unnormalized measure becomes $q^2\operatorname{sech}^2(q\xi)\,d\xi\,d\theta$, whereas Dirichlet energy is unchanged and equals $\int(|v_\xi|^2+|v_\theta|^2)\,d\xi\,d\theta$. Let $a(\xi)$ be the circular mean of $v$. The Poincar\'e inequality on circles gives
\[
 \int_0^{2\pi}|v-a(\xi)|^2\,d\theta
 \leq\int_0^{2\pi}|v_\theta|^2\,d\theta,
 \qquad
 |a(\xi)-a(0)|^2
 \leq\frac{|\xi|}{2\pi}\int |v_\xi|^2.
\]
Integration against the exponentially decaying cylinder weight proves the first inequality, first with $a(0)$ in place of $\overline v$ and then with the weighted mean. On $|\xi|>\log L$, the angular part is bounded using $\operatorname{sech}^2(q\xi)\leq C L^{-2q}$; the mean part uses
\[
 \int_{|\xi|>\log L}(1+|\xi|)
       \operatorname{sech}^2(q\xi)\,d\xi
 \leq C L^{-2q}(1+\log L).
\]
Cauchy--Schwarz for the probability measure and the first estimate about $a(0)$ give
\[
 |\overline v-a(0)|^2\leq C\int|\nabla v|^2.
\]
The constant $\overline v-a(0)$ contributes at most $C L^{-2q}\int|\nabla v|^2$ to either tail. The cut-off normalization is bounded above and below and the cut-off can only decrease the unnormalized integrals up to its fixed supremum.

For bounded $B$, $1+|y^2+B|^2$ is uniformly comparable to $1+|y|^4$, so the same estimates hold for $m_B$ by using $a(0)$ before taking its weighted mean.

For bounded-energy centered fields, the first bound gives a local $L^2$ bound on each fixed cylinder annulus, where the limiting density is bounded above and below.

The compact Sobolev embedding gives strong $L^2$ convergence there along a subsequence. The tail bound also passes to the local limit by exhaustion. Adding the two small tails upgrades local strong convergence to strong convergence in the weighted spaces.

\smallskip\noindent\textit{Step 2. Interaction of the two scales.}\quad For the first interaction estimate, scale by $r$ and integrate $\min\{\delta^{2\beta}|y|^{2\beta},1\}$ against $m_{t,1}$. The ranges $|y|\leq1$, $1\leq |y|\leq\delta^{-1}$ and $|y|\geq\delta^{-1}$ give respectively
\[
 C\delta^{2\beta},\qquad
 C\delta^{2\beta}\int_1^{\delta^{-1}} a^{2\beta-2t-1}\,da,\qquad
 C\delta^{2t}.
\]
For inner decay exponent $\lambda$, scale the second estimate by $R$ and split at $\delta$ and $1$. The three contributions are bounded respectively by
\[
 C\delta^{2s},\qquad
 C\delta^{2\lambda}\int_\delta^1 a^{2s-2\lambda-1}\,da,\qquad
 C\delta^{2\lambda}.
\]
This proves the second estimate with $\lambda=2t$, and its stated extension. Only equality of the compared powers yields the logarithm. The same cut-off comparison applies. The displayed profiles give
\[
 \beta=2s
\]
for the radial outer dilation and $\beta=2$ for a mixed outer field on a compact shape set. The inner dilation has $\lambda=2t$. At a free mixed point the inner translation field has $\lambda=1$.
\end{proof}
\subsection{The linearized operator}
Let $\Hh=\{u\in H^1(M)^2:\int_Mu_i=0\}$ with energy $(u,v)_K=\int\nabla u^{\mathsf T}K^{-1}\nabla v$. For the model densities $V_1=\sum_{i=1}^U m_{s_i}+\sum_{j=1}^V m_{R_j,B_j}$, $V_2=\sum_{j=1}^Vm_{a_j}+n_*^g$, use the common weighted means $\bar u_i=\int V_iu_i/\int V_i$ and the form
\begin{equation}\label{reg:form}
 Q_t(u)=\|u\|_K^2-\sum_{i=1}^2\int_MV_i(u_i-\bar u_i)^2.
\end{equation}
Push each plane model density forward as a measure in the conformal coordinate disk, so its density with respect to $dV_g$ includes $e^{-\varphi}$. Cut off the resulting measures on the fixed disks and normalize them to the displayed masses.

Fix disjoint cut-offs $\chi_a$, equal to one on smaller disks. For a scalar profile $Z_\alpha$ in component $k(\alpha)$ define the source density
\[
 h_\alpha=\chi_aV_{k(\alpha)}
 \left(Z_\alpha-\frac{\int\chi_aV_{k(\alpha)}Z_\alpha}
                         {\int\chi_aV_{k(\alpha)}}\right)e_{k(\alpha)}.
\]
Use $Z=y_j/(1+|y|^2)$ for ordinary and mixed inner translations, and
\[
 Z=\partial_{B_j}\log m_{R,B}
\]
for the two outer shapes. The common radial source is the sum of the separately centered profiles $-2Z_0$ in component one at all first-component bubbles, where $Z_0=(1-|y|^2)/(1+|y|^2)$ ordinarily and
\[
 Z_{0,B}=(1-|y^2-B|^2)/(1+|y^2-B|^2)
\]
at a mixed outer bubble. Each summand has zero integral on its own disk. Set
\[
 S_\alpha=KGh_\alpha,\qquad L_\alpha u=(S_\alpha,u)_K=\int h_\alpha\cdot u.
\]
Then $E_t=\bigcap_\alpha\ker L_\alpha$ has codimension $2U+4V+1$. The sources associated with the approximation use $V^{\rm app}$ in this same definition, preserving cut-off centering exactly.

In scaled center, shape and logarithmic-scale coordinates $\theta$, these sources satisfy
\begin{equation}\label{reg:modulation-bounds}
 \|L_\theta\|_{(C^0)^*}+\|D_\theta L_\theta\|_{(C^0)^*}
 +\|L_\theta\|_{(H^1)^*}+\|D_\theta L_\theta\|_{(H^1)^*}\leq C.
\end{equation}
For a model width $b$, the active density and its first scaled derivatives are bounded by $Cb^{-2}(1+|z/b|^2)^{-2}$. At a mixed point translation is scaled by $a$, so an outer derivative has the additional factor $a/R$. This bound also bounds derivatives through an outer density zero.

The profiles $Z_\alpha$ and their scaled first derivatives are uniformly bounded. The model sources and their scaled derivatives have zero integral on each supporting disk.

Using the weighted variance bound for the displayed radial envelope, equivalently the centered estimate in Lemma~\ref{hyb:centered-estimates}, gives uniform Green energies.

In an inner translation source, the smooth term from the residual density has $L^p$ norm $O(a)$ for $1<p<2$: both $Z_{j,a}$ and its scaled derivatives are bounded by $Ca/(a+|z|)$ on the fixed disk.

Its centering contributes $O(a\sqrt L)$ to the Green norm, since the centered bubble part has vanishing full-plane mean and the fixed cut-off tail has a positive power of $a$.

To compare the differentiated approximation sources, take one scaled derivative $D$ in a mixed disk. The coordinates and matching give
\[
 |Dz|\le Ca,\qquad |DA|\le CA,
\]
$|Da|\le Ca$, $|Dx|\le Cx$ and $|DH|\le CH$ in the notation of the residual proof. In its four-term expansion for $n_2-m_a-d_0$, use
\[
 |DS|\le CAx,\qquad |DE|\le CA
\]
and $|D(a^2/x)|\le Ca^2/x$. They follow by differentiating $S=\beta A^2u-2A\Re(z^2\bar B)$ and
\[
 4xE=\beta u^2-\beta ux+2\beta A\Re(z^2\bar B).
\]
Each differentiated term is again bounded by $CA/x$. Differentiating the two numerator and denominator differences for $n_1-m_{R,B}$ gives $CA^2/H$ by the same bounds. Taylor expansion of the smooth terms gives
\[
 DJ_1=O(|z|^2+a|z|),\qquad DJ_2=O(|z|+a)
\]
in a mixed disk, and
\[
 DJ_1=O(|z|+s)
\]
in an ordinary disk. Cutoff derivatives on fixed annuli have the bounds of the undifferentiated errors. Thus the unnormalized differences and their first derivatives have $L^p$ norm $O(t^{2/p-1})$. Integrating gives this bound also for the differentiated mass errors.

After normalization a derivative can multiply the $O(t^{2/p-2})$ density norm. Taking a fixed $1<p<4/3$ still gives a positive power $O(t^{4/p-3})$ for the differentiated errors. Differentiating cut-off centering subtracts a local mean coefficient times a normalized local density.

The coefficient errors are of the same bound by a positive power of the radius before applying Green, even with the derivative of the centering denominator.

Their Green representatives cost at most an additional $\sqrt L$. This proves the displayed dual bounds. For the comparison of the normalized bilinear forms, one can use Lemma~\ref{hyb:covariance} with the same powers.
\begin{lemma}[Kernel space of the linearized equations on the plane]
\label{reg:jacobi-kernel}
Let $m$ be an ordinary scalar profile or the quadratic profile $m_B$, in coordinates without a source. Suppose $u$ is a real solution of
\[
 \Delta u+2mu=0\quad\text{in }\mathbb R^2
\]
which is bounded or has finite Dirichlet energy. Then $u$ is a linear combination of the following derivatives of $\log m$.
\begin{itemize}
\item[\textup{(i)}]
For the ordinary profile, take the scale derivative and the two translation derivatives.
\item[\textup{(ii)}]
For $m_B$, take the scale derivative and the derivatives with respect to $\Ree B$ and $\Imm B$.
\end{itemize}
\end{lemma}
\begin{proof}
At the origin and in the inverted coordinate $w=1/z$ the transformed density is smooth. It is positive in the ordinary case and has a quadratic zero in the quadratic case.

A bounded solution extends across each puncture after subtracting the Dirichlet potential of $2mu$ and removing the bounded harmonic singularity.

A finite-energy solution is locally in $H^1$ across a puncture. Indeed its circular mean grows by at most $C\sqrt{1+|\log r|}$ by Cauchy--Schwarz on logarithmic radii.

The remaining angular oscillation has a weighted area $L^2$ bound by the circular Poincar\'e inequality and the finite energy.

Testing the equation with logarithmic cut-offs that equal zero near the puncture and have gradients tending to zero in $L^2$ removes the point in the weak equation. Local elliptic regularity now gives smoothness both at zero and in the inverted coordinate.

Away from the zeros put $v=\log m$ and $\mathcal W=u_{zz}-v_z u_z$. The equations $v_{z\bar z}=-m/2$ and $u_{z\bar z}=-mu/2$ give $\partial_{\bar z}\mathcal W=0$: the two terms involving $m u_z$ cancel, as do the terms involving
\[
 m_z=v_zm.
\]
Under a holomorphic coordinate change the added $\log$ of the area factor in $v$ cancels the second derivative of the coordinate in $u_{zz}$, so $\mathcal W\,dz^2$ is a quadratic differential. At each quadratic zero $v_z=1/z+O(1)$ in the local coordinate, and hence $\mathcal W$ has at most a simple pole.

At ordinary points it is regular. In plane coordinates this gives $z\mathcal W$ bounded at zero and tending to zero at infinity (the quadratic differential contributes $z^{-4}$ under inversion). Thus
\[
 \mathcal W=0.
\]
Its conjugate and the linearized equation now give a closed first-order system for $(u,u_z,u_{\bar z})$ off the punctures. Values of these three real data at one ordinary point determine the solution by uniqueness along paths, so the dimension is at most three.

Differentiating the three real parameters $(\ell,\Ree D,\Imm D)$ of Lemma~\ref{fnd:entire-profiles} at $q=1,2$ gives three solutions.

They are independent: the leading $|z|^{2q}$ coefficient in $|z|^{q-1}e^{-v/2}$ fixes the variation of $\ell$, and the coefficients of $z^q,\bar z^q$ then fix the two variations of $D$. This proves the assertion.
\end{proof}
\begin{lemma}[Derivatives of the normalized exponential]
\label{journal:normalized-nonlinearity}
Let $v_T\geq0$, $\int_Mv_T=1$, and suppose that, for fixed $p_0>1$, $b\geq0$ and $C_0<\infty$,
\[
                 \|v_T\|_{p_0}\leq C_0e^{bT},\qquad T\geq T_0.
\]
Put $\Psi_T(u)=\log\int_Mv_Te^u$. Fix $r_0<\infty$ and consider $\sqrt T\|u\|_{H^1}\le r_0$.
\begin{itemize}
\item[\textup{(i)}]
Its derivatives satisfy
\begin{equation}\label{journal:cumulant-estimates}
 \bigl|D^j\Psi_T(u)[z_1,\ldots,z_j]\bigr|
 \leq C T^{j/2}\prod_{k=1}^j\|z_k\|_{H^1},\qquad j=1,2,3.
\end{equation}
\item[\textup{(ii)}]
If one argument is bounded, write it as $a$ in the third slot. Then
\begin{equation}\label{journal:bounded-tangent-estimate}
 |D^3\Psi_T(u)[z_1,z_2,a]|
 \leq CT\|a\|_\infty\|z_1\|_{H^1}\|z_2\|_{H^1}.
\end{equation}
\item[\textup{(iii)}]
For a finite sum of these functionals with bounded coefficients, the Taylor remainder $\mathcal N_T$ of its gradient in the Dirichlet inner product satisfies
\begin{equation}\label{journal:normalized-taylor}
 \|\mathcal N_T(u)-\mathcal N_T(w)\|_{H^1}
 \leq CT^{3/2}(\|u\|_{H^1}+\|w\|_{H^1})\|u-w\|_{H^1}
\end{equation}
on the above ball.
\end{itemize}
These estimates also hold for the two components of the Toda functional with the fixed Cartan inner product.
\end{lemma}
\begin{proof}
\smallskip\noindent\textit{Step 1. Exponential moments.}\quad For large $T$ set $p_T=1+T^{-1}$. Interpolation between $L^1$ and $L^{p_0}$ gives
\[
 \|v_T\|_{p_T}\leq C,
\]
since the interpolation parameter is $p_0/((p_0-1)(T+1))$. The Sobolev and exponential inequalities on the surface give, for $k\geq2$ and mean-zero $z$,
\[
 \|z\|_k\leq C\sqrt{k}\|z\|_{H^1},\qquad
 \|e^{|u|}\|_k\leq C\exp(Ck\|u\|_{H^1}^2).
\]
Apply H\"older first with exponents $p_T,T+1$, and then with $j+1$ equal exponents to the exponential and the $j$ functions. For $1\leq j\leq3$ this proves
\[
 \int_Mv_Te^{|u|}\prod_{k=1}^j|z_k|
 \leq CT^{j/2}e^{CT\|u\|_{H^1}^2}
                         \prod_{k=1}^j\|z_k\|_{H^1}.
\]
Also $|\int v_Tu|\leq C\sqrt T\|u\|_{H^1}$, so Jensen's inequality yields
\[
 \int_Mv_Te^u\geq\exp\left(\int_Mv_Tu\right)\geq e^{-Cr_0}.
\]
The same moment estimates therefore hold for $d\nu_u=v_Te^u\,dV_g/\int_Mv_Te^u$. Differentiating this probability measure gives the exact formulas
\[
 \begin{aligned}
 D\Psi_T(u)[z]&=\int z\,d\nu_u,\\
 D^2\Psi_T(u)[z_1,z_2]&=\operatorname{Cov}_{\nu_u}(z_1,z_2),\\
 D^3\Psi_T(u)[z_1,z_2,z_3]
   &=\int\prod_{k=1}^3\left(z_k-\int z_k\,d\nu_u\right)d\nu_u.
 \end{aligned}
\]

\smallskip\noindent\textit{Step 2. Derivatives of the normalized integral.}\quad The moment bounds prove \eqref{journal:cumulant-estimates}. For a bounded argument,
\[
 |a-\int a\,d\nu_u|\leq2\|a\|_\infty;
\]
Cauchy--Schwarz and the second-moment bound prove \eqref{journal:bounded-tangent-estimate}. To justify the Fr\'echet derivatives, at fixed $T$ expand $e^w$ through order three. The integrated remainder is $O_T(\|w\|_{\Hh}^4)$ by
\[
 |e^w-\sum_{k=0}^3w^k/k!|\leq C|w|^4 e^{|w|}
\]
and H\"older's inequality with the exponents $p_T$ and $T+1$, followed by the exponential and finite-moment bounds above. Apply the same Taylor estimate after inserting products of at most three variation directions, using finite higher H\"older exponents. The inequality
\[
 |e^w-1|\leq |w|e^{|w|}
\]
then gives continuity in the corresponding multilinear operator norms. Differentiating the quotient gives
\[
 D\nu_u[z]=\nu_u(z-\int z\,d\nu_u),
\]
with equality tested against integrable products of the variations. The Jensen lower bound on the denominator justifies taking its logarithm. For the gradient remainder, the fundamental theorem of calculus gives
\[
 D\mathcal N_T(u)[v]
   =\int_0^1D^2\mathcal F_T(B_T+\sigma u)[u,v]\,d\sigma .
\]
Pair both sides with a function of Dirichlet norm one. For the difference in \eqref{journal:normalized-taylor} integrate $D\mathcal N_T$ in the direction $u-w$ along the segment from $w$ to $u$. For the bounded test $a$, the integral remainder is
\[
 \langle\mathcal N_T(u),a\rangle_{\Hh}
    =\int_0^1(1-\sigma)D^3J_T(B_T+\sigma u)[u,u,a]\,d\sigma ,
\]
where $J_T$ is the corresponding Toda energy. Use the bounded-argument estimate for each component of the logarithmic part; the linear Dirichlet gradient contributes no remainder and the fixed Cartan matrix changes only the constants.
\end{proof}
\begin{lemma}[An a priori estimate for the projected equation]
\label{journal:projected-contraction}
Let $P_T$ project orthogonally onto a closed subspace $E_T$ of the space $\Hh$. Assume the gradient $\mathcal F_T$ satisfies \eqref{journal:normalized-taylor} about $B_T$, uniformly in $T$, on
\[
 \sqrt T\|u\|_{\Hh}\le r_0,\qquad r_0>0\ \text{fixed}.
\]
Suppose
\[
 A_T=P_TD\mathcal F_T(B_T)|_{E_T},\qquad
 \|A_T^{-1}\|\leq C_A T.
\]
Choose $\mathfrak r_T>0$ and define $r_T$ by
\[
 \|P_T\mathcal F_T(B_T)\|\le\mathfrak r_T,\qquad
 T^{7/2}\mathfrak r_T\to0,\qquad r_T=2C_AT\mathfrak r_T.
\]
\begin{itemize}
\item[\textup{(i)}]
For sufficiently large $T$ there is a unique solution of
\[
 P_T\mathcal F_T(B_T+\eta)=0,\qquad \eta\in E_T,
\]
in the closed ball $\|\eta\|\leq r_T$.
\item[\textup{(ii)}]
Every solution satisfying $\sqrt T\|\eta\|\leq r_0$ and $C_AC_NT^{5/2}\|\eta\|\leq\tfrac12$ obeys
\begin{equation}\label{journal:absorption-threshold}
                    \|\eta\|\leq2C_AT\mathfrak r_T,
\end{equation}
where $C_N$ is a fixed constant in the quadratic remainder bound. If the concentration radius is exponentially small in $T$, a bound on $\|\eta\|$ by a positive power of that radius suffices.
\end{itemize}
\end{lemma}
\begin{proof}
The equation is equivalent to
\[
 \eta=-A_T^{-1}P_T\mathcal F_T(B_T)
       -A_T^{-1}P_T\mathcal N_T(\eta).
\]
On $\|\eta\|\le r_T$, the right-hand side has norm at most $r_T/2+C_AC_NT^{5/2}r_T^2$, and Lipschitz constant at most $2C_AC_NT^{5/2}r_T$. Both requirements for contraction follow from
\[
 T^{7/2}\mathfrak r_T\to0.
\]
Also $\sqrt T r_T\to0$, as required by the nonlinear estimate. For an arbitrary solution in the stated range, the same inverse equation gives
\[
 \|\eta\|\leq C_AT\mathfrak r_T
                   +C_AC_NT^{5/2}\|\eta\|^2.
\]
The stated smallness condition absorbs the last term and proves \eqref{journal:absorption-threshold}. This argument uses the preliminary smallness of the error; the inequality by itself would not exclude a second, larger range of values of $\|\eta\|$.
\end{proof}
\begin{proposition}[Invertibility and Morse index]
\label{reg:inverse}
Let $L=\log(1/t)$ and let $L_t$ be the restriction of the self-adjoint operator representing \eqref{reg:form} to the projected space, with orthogonal projection onto that space.
\begin{itemize}
\item[\textup{(i)}]
We have
\begin{equation}\label{reg:inverse-index}
 \|L_t^{-1}\|\leq CL,\qquad
 \ind L_t=m_S+U+4V-1.
\end{equation}
\item[\textup{(ii)}]
Exactly $r=U+2V-1$ eigenvalues tend to zero. Each lies in $[c/L,C/L]$.
\item[\textup{(iii)}]
For $|\mu|\leq C_0t$, both conclusions hold for $u^{\rm app}$ and for the small family of projected solutions constructed below.
\end{itemize}
\end{proposition}
\begin{proof}
\emph{Step 1. The projected operator.} For components with local masses $M_a$, the exact variance identity is
\[
 \int\Bigl(\sum_am_a\Bigr)(u-\bar u)^2
 =\sum_a\int m_a(u-\bar u_a)^2
       +\sum_aM_a(\bar u_a-\bar u)^2.
\]
The matrix of the first-component means has rank $n-1$. Since
\[
 \sigma>0,
\]
subtracting the residual mean from each of the $V$ inner means gives $V$ independent differences in the second component. The total rank is $r$.

Write the full operator as $I-B_t-V_tCV_t^*$, where $B_t$ contains separately centered variances, $C>0$ is their exact finite mass matrix, and $V_t$ has as columns the Green representatives of these differences of means.

For the auxiliary dilation at each scale, let $m$ denote the model density of that bubble alone. Its source in the active component is
\[
 \chi m\left(Z-\frac{\int\chi mZ}{\int\chi m}\right).
\]
In particular, use the inner bubble density $m_{a_j}$ for an inner dilation. The residual density enters the second bilinear form in \eqref{reg:form}. The constraints defining $E_t$ are those specified above.

Let $\mathcal K_t$ be the span of the Green representatives of the translation and shape sources and of all the individual dilation sources just defined. Set
\[
 W_t=\mathcal K_t^\perp.
\]
We first prove that
\[
 A_t=P_{W_t}(I-B_t)|_{W_t}
\]
has a uniformly bounded inverse. Lemma~\ref{hyb:centered-estimates} proves the following estimates for the separately centered measures.

In the scaled coordinate of a bubble, let $\widetilde m$ denote its cut-off density, normalized to its prescribed mass, and let $\langle v\rangle_{\widetilde m}$ be the corresponding mean. Circular Poincar\'e and Cauchy--Schwarz for the circular mean give
\begin{align}
 \int\widetilde m|v-\langle v\rangle_{\widetilde m}|^2
 &\le C\int|\nabla v|^2,\notag\\
 \int_{\{|y|<H^{-1}\}\cup\{|y|>H\}}
       \widetilde m|v-\langle v\rangle_{\widetilde m}|^2
 &\le CH^{-\alpha}(1+\log H)\int|\nabla v|^2,
 \qquad H\ge2.
 \label{reg:centered-tails}
\end{align}
Here $\alpha=2$ for ordinary and inner profiles and $\alpha=4$ for outer profiles. The integrals are taken within the coordinate disk. To obtain the second estimate, compare each circular mean with the mean on $|y|=1$; their squared difference is bounded by
\[
 C|\log|y||\int|\nabla v|^2.
\]
Integrating against the radial density gives the stated powers, and the angular part follows from circular Poincar\'e. Subtracting the weighted mean gives the same bound. For bounded $B$ we have
\[
 m_B\asymp m_{2,1},
\]
so the estimates still hold. Compactness of the local Sobolev embedding and \eqref{reg:centered-tails} imply, after passage to a subsequence, strong convergence in the weighted $L^2$ spaces for every energy-bounded sequence of centered rescalings.

Choose a basis $\Psi_{\alpha,t}$ of $\mathcal K_t$ with the preceding convention for the dilations. On an inner disk the contribution of the residual density to a translation source gives a smaller error. Indeed, for $Z_{j,a}(z)=az_j/(a^2+|z|^2)$ and
\[
 1<p<2,\qquad \|\chi n_*^g Z_{j,a}\|_p=O(a),
\]
since $n_*^g$ is bounded there. The change in the centering constant is $O(a)$, by oddness of the bubble part on the coordinate disk. Its product with the bubble density is $O(a\sqrt{1+|\log a|})$ in $H^{-1}$, by the logarithmic Green estimate.

Thus the two conventions for a translation give Green representatives differing by $o(1)$ in Dirichlet norm. The matrix $((\Psi_{\alpha,t},\Psi_{\beta,t})_K)_{\alpha,\beta}$ has a positive definite limit. Indeed, the cut-off of
\[
 (Ke_{k(\alpha)}/2)(Z_\alpha-Z_\alpha(\infty))
\]
on its scale, with mean removed, differs from $\Psi_{\alpha,t}$ by $O(t^cL^C)$ in Dirichlet norm: their Poisson sources differ only by the cut-off terms, the model tail and the error from centering the residual density in the translation equation estimated above.

Its own centered variance recovers the source since
\[
 \int mZ_\alpha=0.
\]
Thus each diagonal block converges to the positive definite matrix of inner products of the corresponding kernel functions. The off-diagonal blocks tend to zero. The interaction bounds in Lemma~\ref{hyb:centered-estimates}, together with the fixed cut-off errors, give
\begin{equation}\label{reg:centered-jacobi-error}
 \|(I-B_t)\Psi_{\alpha,t}\|_K\le Ct^cL^C
\end{equation}
for some $c>0$. For a nested pair, subtract the outer field's value at the inner center and the inner field's limit at infinity.

The remaining weighted integrals tend to zero with a positive power of $a/R$. Cutoff errors, remote bubbles and the residual measure give positive powers of the corresponding radii.

These estimates also give convergence of the off-diagonal blocks to zero independently of the inverse estimate.

\emph{Step 2. A uniform inverse on $W_t$.} Suppose, to the contrary, that $t_k\to0$ and $\phi_k\in W_{t_k}$ satisfy $\|\phi_k\|_K=1$ and $A_{t_k}\phi_k\to0$. Then
\[
 (I-B_{t_k})\phi_k
   =h_k+\sum_\alpha\lambda_{\alpha,k}\Psi_{\alpha,t_k},
 \qquad h_k\in W_{t_k},\qquad \|h_k\|_K\to0.
\]
Pairing with $\Psi_{\beta,t_k}$, using self-adjointness and \eqref{reg:centered-jacobi-error}, gives
\[
 \sum_\alpha\lambda_{\alpha,k}
       (\Psi_{\alpha,t_k},\Psi_{\beta,t_k})_K=o(1).
\]
The uniform inverse bound for $((\Psi_{\alpha,t_k},\Psi_{\beta,t_k})_K)_{\alpha,\beta}$ implies
\[
 \lambda_{\alpha,k}=o(1),
\]
and hence $(I-B_{t_k})\phi_k\to0$ in $\Hh$.

Let $\phi_k\rightharpoonup\phi$ in $\Hh$. Every centered bubble source tends to zero when paired with a smooth test function, by Cauchy--Schwarz and the vanishing variance of that test function on the bubble scale. The limiting equation is therefore the scalar shadow linearization, with
\[
 \phi_1=-\phi_2/2.
\]
Its invertibility implies
\[
 \phi=0;
\]
the variance for the residual measure tends to zero. Subtract the active weighted mean at each bubble before rescaling. The rescaled compactly supported tests have uniformly bounded Dirichlet norms, so they may be used in the residual equation.

At a mixed point on the outer scale, the inner centered source has zero integral.

Its action is bounded by the square root of the centered density variance of the fields times that of the test function in the inner measure. The latter tends to zero, since the test varies on the outer scale and
\[
 a/R\to0.
\]
On the inner scale the outer model mass on each test support tends to zero, by the outer density bound and
\[
 a/R\to0.
\]
Cauchy--Schwarz excludes its contribution there. The remote and bounded residual measures also give zero on the shrinking test supports. The limiting active field $z$ has finite Dirichlet energy and satisfies
\[
 -\Delta z=2mz.
\]
The other component plus one half of the active one is harmonic. Finite energy removes isolated punctures, so no additional logarithmic term occurs at either scale.

By Lemma~\ref{reg:jacobi-kernel}, every finite-energy solution of either scalar equation on the plane is a linear combination of its three kernel functions. The conditions $\phi_k\perp\mathcal K_{t_k}$ pass to
\[
 \int mZ_\alpha z=0
\]
by \eqref{reg:centered-tails}. The positive definiteness of $(\int mZ_\alpha Z_\beta)_{\alpha,\beta}$ gives $z=0$. Strong weighted convergence now makes every separately centered bubble variance tend to zero. Testing $(I-B_{t_k})\phi_k=o(1)$ against $\phi_k$ yields
\[
 1=\sum_{\text{bubbles}}\int m_a
       |\phi_{k,i(a)}-\langle\phi_{k,i(a)}\rangle_a|^2
       +\int n_*^g
       |\phi_{k,2}-\langle\phi_{k,2}\rangle_g|^2+o(1)=o(1),
\]
a contradiction. This proves the uniform inverse bound, including the exclusion of approximate kernels supported between the bubble scales. \emph{Step 3. The index on $W_t$.} For an active scalar component $z$ and the other component $v$, the Dirichlet form splits as
\[
 \nabla(z,v)^{\mathsf T}K^{-1}\nabla(z,v)
   =\tfrac12|\nabla z|^2
       +\tfrac23|\nabla(v+z/2)|^2.
\]
The local negative index is consequently that of the centered scalar form $\tfrac12\int|\nabla z|^2-\int mz^2$. For the radial profile with exponent $q$, use $\xi=\log|y|$, angular mode $j\in\mathbb Z$, $x=\tanh(q\xi)$ and $\nu=|j|/q$. The weighted eigenvalue equation is
\[
 -((1-x^2)f')'+\nu^2 f/(1-x^2)=\Lambda f.
\]
Substituting $f=(1-x^2)^{\nu/2}P$ gives the self-adjoint polynomial operator $-(1-x^2)P''+2(\nu+1)xP'+\nu(\nu+1)P$ in $L^2((1-x^2)^\nu dx)$.

The polynomial modes have finite energy and solve the weak equation by cut-off at the poles. Orthogonal polynomials by degree $k$ diagonalize this operator and are complete by polynomial density, with
\[
 \Lambda=(k+\nu)(k+\nu+1).
\]
After removal of the constant, none is below two when $q=1$, and exactly two real modes are below two when $q=2$ ($k=0$, $j=\pm1$).

For $m_B$ the latter index is unchanged along bounded $B$ paths: the weighted embedding is compact, the forms depend continuously on $B$, and their kernel space remains three-dimensional. The macroscopic negative index is $m_S$.

Cut off and transplant negative scalar fields and negative shadow fields, then project off $\mathcal K_t$. Their quadratic forms and matrices of energy inner products converge to those of the direct sum, so
\[
 \ind A_t\ge m_S+2V.
\]
For the reverse inequality, the uniform inverse bound gives a fixed spectral gap at zero. If the index were larger, choose $m_S+2V+1$ energy-orthonormal negative eigenvectors. Their span satisfies
\[
 (A_tu,u)_K\le-\delta\|u\|_K^2
\]
for a fixed
\[
 \delta>0.
\]
Take macroscopic and centered local limits. The potential terms converge by \eqref{reg:centered-tails}. Lower semicontinuity on disjoint rescaled annuli and outside the bubble disks bounds the sum of limiting energies by the original energy as the annuli exhaust the planes.

Hence every nonzero linear combination of the limiting fields has strictly negative limiting form.

The limit map is injective and gives a negative subspace of dimension $m_S+2V+1$ in the direct sum of the scalar shadow form, the centered plane forms and the nonnegative inactive forms. This contradicts the index just computed. We conclude that
\[
 \ind A_t=m_S+2V.
\]

Before projection, the inner products of the Green representatives of the local means form the matrix
\begin{equation}\label{reg:mean-gram}
 L\left\{\bigoplus_A\frac1\pi\right\}
 \oplus
 L\left\{\bigoplus_B\frac1{4\pi}
                  \begin{pmatrix}2&-1\\-1&4\end{pmatrix}\right\}+O(1).
\end{equation}
Indeed the averaged Green singularity uses $\log(1/\max(r_i,r_j))/(2\pi)$ and the Cartan entries $2,2,-1$; a mixed outer radius is $\sqrt t$ and its inner radius is $t$. After taking differences of the means and projecting onto $W_t$, the corresponding vectors (again denoted by $V_t$) satisfy
\[
 G_t=V_t^*V_t=LG_0+O(1)
\]
with
\[
 G_0>0.
\]
Projection changes their matrix of inner products by $O(1)$. The mean Green potential of a probability bubble of width $b$ on its own disk is $(2\pi)^{-1}\{\log(1/b)-\log(1+|z|/b)\}+O(1)$ by rescaling the logarithmic kernel against the bounded scaled profile and its tails.

On a concentric disk of width $a$, subtracting $(2\pi)^{-1}\log(1/\max(a,b))$ leaves a fluctuation bounded by $C(1+|\log|z/a||)$ in weighted $L^2$; on other disks and against the residual density use separation and the uniformly square-integrable Green kernel. Thus
\[
 \|B_W V_t\|=O(1),\qquad V_t^*B_W V_t=O(1),\qquad
 F_t:=V_t^*A_t^{-1}V_t=G_t+O(1).
\]
Here $B_W=P_{W_t}B_t|_{W_t}$ and $A_t=I-B_W$ on $W_t$. The last identity follows from
\[
 A_t^{-1}=I+B_W+B_WA_t^{-1}B_W.
\]
The formula for the inverse under a finite-rank perturbation now shows that $A_t-V_tCV_t^*$ has a uniform inverse on $W_t$ and adds exactly $r$ negative mean directions. Its inverse maps $V_t$ to vectors of norm $O(L^{-1/2})$.

\emph{Step 4. The $r$ small eigenvalues.} Restore the $r$ radial modes allowed by the common constraint. If their outer and inner amplitudes are $\gamma_a,\beta_j$, put $\theta_j=(1+|B_j|^2)^{-1}$. The limiting differences of their means are
\[
 D(\gamma,\beta)=
 ((\gamma_a-\gamma_n)_{a<n},
                      (\beta_j-\theta_j\gamma_j)_{j\leq V}).
\]
Its kernel space is the common matched dilation. The sum of the source pairings with the radial functions tends to $-2\sum_a(M_a/3)\gamma_a$. The equivalent zero constraint
\[
 \sum_a(M_a/3)\gamma_a=0
\]
excludes that kernel space, so $D$ is invertible on the retained radial space. For the limiting plane potential normalized to vanish at infinity,
\[
 2G_{\mathbb R^2}(mZ)=Z+1
\]
: the source has zero integral, and both sides solve the same Poisson equation and vanish at infinity. This gives the inner constant $-(Z_{0,B}(0)+1)/2=-\theta_B$. Splitting nested weighted integrals at $a/R\asymp\sqrt t$ gives
\[
 \|Z_{0,B}((a/R)\cdot)-Z_{0,B}(0)\|_{L^2(m_i)}=O(a/R),
\]
\[
 \|Z_i(\cdot/(a/R))+1\|_{L^2(m_B)}
   =O((a/R)^2\sqrt{1+|\log(a/R)|}).
\]
The cut-off errors and the interactions with the other bubbles are bounded by positive powers of $t$. These bounds make the centered radial residuals
\[
 O(t^cL^C)=o(L^{-k})
\]
for each fixed $k$. For $r>0$, choose a uniformly nondegenerate frame for the retained radial space in $E_t\cap W_t^\perp$, and let $D_t$ be its matrix of differences of means, so that
\[
 D_t\to D.
\]
Here the columns of $V_t$ are the mean Green vectors projected onto $W_t$. Relative to $W_t$ and this radial frame, the compressed form has matrix
\[
 \begin{pmatrix}
 A_t-V_tCV_t^*&-V_tCD_t\\
 -D_t^{\mathsf T}CV_t^*&-D_t^{\mathsf T}CD_t
 \end{pmatrix}+O(t^cL^C).
\]
Writing
\[
 M_t=A_t-V_tCV_t^*,
\]
the same inverse formula gives the matrix obtained by eliminating the complementary variables in the displayed principal matrix as
\begin{align*}
 &-D_t^{\mathsf T}CD_t
   -D_t^{\mathsf T}CV_t^*M_t^{-1}V_tCD_t\\
 &\hspace{2em}=D_t^{\mathsf T}(F_t-C^{-1})^{-1}D_t.
\end{align*}
The preceding residual bounds and the uniform inverse bound for $M_t$ therefore give the matrix on the remaining variables of the full form:
\[
 D_t^{\mathsf T}(F_t-C^{-1})^{-1}D_t+O(t^cL^C)
       =L^{-1}D^{\mathsf T}G_0^{-1}D+o(L^{-1})>0.
\]
When $r=0$ there is no radial block to restore. For $r>0$, let $Q_t^{\rm rad}$ be the positive definite matrix of energy inner products of the radial basis functions and write the orthogonal decomposition
\[
 E_t=W_t\oplus(E_t\cap W_t^\perp)
\]
in these coordinates as
\[
 L_t=\begin{pmatrix}\widetilde M_t&K_t^{\rm rad}\\
                       (K_t^{\rm rad})^*&J_t\end{pmatrix}.
\]
Here $\widetilde M_t=M_t+O(t^cL^C)$ has a fixed inverse bound, and $(\widetilde M_t)^{-1}K_t^{\rm rad}=O(L^{-1/2})$ by the preceding estimate on $(M_t)^{-1}V_t$. For $|\lambda|$ in a fixed sufficiently small interval, use the reduced matrix
\[
 S_t(\lambda)=J_t-\lambda Q_t^{\rm rad}
 - (K_t^{\rm rad})^*(\widetilde M_t-\lambda I)^{-1}K_t^{\rm rad}.
\]
Its derivative is
\[
 -Q_t^{\rm rad}-(K_t^{\rm rad})^*(\widetilde M_t-\lambda I)^{-2}
K_t^{\rm rad}=-Q_t^{\rm rad}+O(L^{-1})
\]
in operator norm, by the resolvent identity on the fixed interval. Its ordered eigenvalues are strictly decreasing there, by the uniform positive lower bound on $Q_t^{\rm rad}$.

Since $S_t(0)$ is positive with eigenvalues comparable to $L^{-1}$, $S_t(c/L)$ is positive and $S_t(C/L)$ negative for suitable fixed
\[
 0<c<C.
\]
Each of its $r$ ordered eigenvalues crosses zero once, and the kernel space dimensions agree with those of $L_t-\lambda I$. This proves exactly the stated $r$ small positive eigenvalues, the fixed gap for the others, and the index $m_S+2V+r=m_S+U+4V-1$.

\emph{Step 5. Perturbation of the operator.} For fixed $p>1$ an $L^p$ density error $e_p$ changes the uncentered second moments by $O(e_p)$ in operator norm, by H\"older and $H^1\hookrightarrow L^{2p'}$.

The original linear moments have norm $O(\sqrt L)$ by the probability estimate in the residual proof, and perturbations of the unnormalized linear moments and masses are $O(e_p)$.

Expanding the normalized mean terms gives an $O(e_pL)$ change of the variance form. In each source the bounded kernel function and the density error give an $O(e_p)$ Green error before recentering.

The centering coefficient changes by $O(e_p)$, and its constant times a normalized local density costs $O(e_p\sqrt L)$ in Green norm.

The uniform inverse bound for the matrix of inner products of the source representatives therefore gives the same bound on the projection difference.

The isometry of the two nearby compressed spaces differs from the identity by $O(e_p\sqrt L)$, and the full form has operator norm $O(L)$. Hence the compressed form differs by $O(e_pL^{3/2})$. Lemma~\ref{reg:residual} gives
\[
 e_{6/5}=O(t^{2/3}),
\]
which is $o(L^{-1})$ even with these logarithms. Along the segment from the trial field to a correction $\eta$ with
\[
 \|\eta\|_{H^1}=O(t^{2/3}L),
\]
Lemma~\ref{journal:normalized-nonlinearity} bounds the third differential of each normalized logarithm by $CL^{3/2}$. Thus the change of the second derivative is bounded by $CL^{3/2}\|\eta\|_{H^1}=O(t^{2/3}L^{5/2})=o(L^{-1})$. The gap and the index therefore persist.
\end{proof}
\subsection{Projected solutions and their reduced degree}
Use the source densities of the approximation $h_\alpha$ defined above, including the common radial source, and put $q_c=\sum_\alpha c_\alpha h_\alpha$. We write the projected equation as
\begin{equation}\label{reg:projected-equation}
 -\Delta_g(u^{\rm app}+\eta)
 =K\Pi\left(\frac{\rho_i h_i e^{u_i^{\rm app}+\eta_i}}
                        {\int_Mh_i e^{u_i^{\rm app}+\eta_i}}\right)_i
     +Kq_c,\qquad\eta\in E_t.
\end{equation}
Each source has zero integral. In local canonical equations $q_c$ denotes its Euclidean-area density, including the conformal factor. The inverse yields a unique small orthogonal correction for $|\mu|\leq C_0t$:
\begin{equation}\label{reg:range}
 \eta\in E_t,\qquad
 \|\eta\|_{H^1}=O(t^{2/3}L),\qquad |c|=O(t^{2/3}L^2).
\end{equation}
For the multiplier estimate test against the cut-off kernel functions of the scalar linearization in their active components, with their surface means removed. Their Dirichlet norms are uniformly bounded by scaling.

On a common radial source it suffices to test a single first-component bubble with its cut-off $-2Z_0$ (or $-2Z_{0,B}$).

The limiting test matrix is invertible: the translation and shape blocks are their positive weighted moment matrices, the radial entry is a positive weighted square, and the other pairings vanish by angular orthogonality, separation, or scale separation.

On the full space the linear operator norm is $O(L)$ by the moment estimate for the probability measures. The residual, the bound on $\eta$, and the normalized quadratic remainder then bound the projected gradient by
\[
 O(t^{2/3}L^2),
\]
giving the stated estimate on $c$. This construction also applies for $V=0$, with the mixed tests absent.

For the Taylor estimate after normalization, replace each approximate probability measure $\nu$ by the measure with density $e^v/\int\nu e^v$ with
\[
 \|v\|_{H^1}\le C t^{2/3}L.
\]
The upper bound for the density gives $\|\nu\|_{1+1/L}=O(1)$. The denominator is bounded below by $\exp(\int\nu v)$, with
\[
 |\int\nu v|\le C\sqrt L\|v\|_{H^1}.
\]
After H\"older with $1+1/L$, use the exponent estimate of Lemma~\ref{fnd:exponential} at order $4(L+1)$ for each of three variations and for $e^{|v|}$. The third moment of the three variations is then bounded by $CL^{3/2}$ times their three Dirichlet norms.

The centered third moment satisfies the same bound. This proves the Taylor estimate used in Lemma~\ref{journal:normalized-nonlinearity} on the current ball. Lemma~\ref{journal:projected-contraction} applies to the compressed inverse. Alternatively, fix
\[
 0<\varepsilon<2/3.
\]
H\"older's inequality and the exponential Sobolev inequality give
\[
 \|N(v)-N(w)\|_{H^{-1}}
 \leq C_\varepsilon t^{-\varepsilon}
            (\|v\|+\|w\|)\|v-w\|.
\]
The common denominator is retained. For one density let $V=\rho_i\nu$, where $\nu$ is its probability weight, and put $a=\int\nu v$ and $b=\int\nu(e^v-1-v)$. Subtracting $V(1+v-a)$ from $Ve^v/(1+a+b)$ gives the exact remainder
\[
 V\frac{(e^v-1-v)-b-(a+b)(v-a)}{1+a+b}.
\]
For a probability density bounded by $Ct^{-2}$, H\"older's inequality with exponents $1+1/L$ and $L+1$, followed by Lemma~\ref{fnd:exponential}, gives
\begin{equation}\label{reg:weighted-square}
 \int\nu v^2e^{c|v|}\leq C_cL\|v\|^2e^{C_cL\|v\|^2}.
\end{equation}
Indeed, bound the $L^{L+1}$ norm of the product by the $L^{4(L+1)}$ norm of $v$ squared and the $L^{2(L+1)}$ norm of $e^{c|v|}$. The probability density has bounded $L^{1+1/L}$ norm. Thus $|b|\leq CL\|v\|^2$ on the current ball,
\[
 |a|+\int\nu|v-a|\leq C\sqrt L\|v\|,\qquad 1+a+b=\int\nu e^v\geq e^a
\]
by Jensen's inequality. The inverse is a contraction since $L t^{-\varepsilon}(t^{2/3}L)\to0$. The density remainder above has $L^1$ norm
\[
 O(L\|v\|^2)=O(t^{4/3}L^3)=o(t)
\]
and hence the same order against uniformly bounded bubble tests, including the normalization term.
\begin{proposition}[The reduced equations when $V\geq1$]
\label{reg:mixed-map}
Suppose $V\geq1$.
\begin{itemize}
\item[\textup{(i)}]
The remaining equations have the following expansions, uniform in $C^0$:
\begin{align}
 c_i^A&=-3s_iF_i^A+o(t),&c_j^B&=-3a_jF_j^B+o(t),\label{reg:rows}\\
 c_j^{\rm sh}&=-\frac{3\epsilon_j}{16\pi}
              (\kappa_j\nabla\mathcal M(B_j)+4\pi\mathbf t_j)+o(t),
 &c_0&=\frac3{8\pi m}
       \left(\mu-\sum_j\epsilon_j\kappa_j\mathcal H(B_j)\right)+o(t).
 \nonumber
\end{align}
\item[\textup{(ii)}]
Assume the shadow solution $Z$ is nondegenerate and $\kappa_j>4|T_j|$ for every mixed point. In the isolating neighborhood constructed below, the degree for $\mu>0$ is $(-1)^U d_S(Z)$. For $\mu<0$ there are no zeros there.
\end{itemize}
\end{proposition}
\begin{proof}
\emph{Step 1. Derivatives of the approximate solution.} Write $u=u^{\rm app}+\eta$. Let $Y$ denote a derivative of $u^{\rm app}$ with respect to a center coordinate, multiplied by its translation width, a shape derivative, or the derivative in
\[
 \lambda=\log t.
\]
All these tangents are uniformly bounded. An ordinary center tangent equals $(4Z_k,-2Z_k)+O(t)$ on its own disk, where
\[
 Z_k(y)=y_k/(1+|y|^2).
\]
A mixed center tangent on $z=ay$ tends to $(-2Z_k,4Z_k)$ on compact scaled disks. On its entire mixed disk each component is bounded by
\[
 C\left(\frac{a}{R}+\frac{a|z|}{a^2+|z|^2}\right).
\]
Indeed, $F_R\asymp R^4+|z|^4$ and
\[
 Q_R\asymp R^4(a^2+|z|^2)
\]
on the fixed disk, and their spatial logarithmic derivatives are bounded by $C/R$ and $C|z|/(a^2+|z|^2)$, respectively. The derivatives of the smooth backgrounds and of the logarithms of the matching coefficients are bounded before multiplication by $a$.

Center tangents outside their own disks are $O(t)$. The shape tangents are even under $z\mapsto -z$ on their own disk before mean subtraction, bounded on $|z|\le R$, and $O(R^2/|z|^2)$ on $|z|\ge R$. The differentiated residuals satisfy
\[
 \|D\mathcal F_\rho(u^{\rm app})Y\|_K=O(t^{2/3}).
\]
To see this, differentiate \eqref{reg:cancelled-residual} and \eqref{reg:normalized-residual} before taking norms. On each ordinary disk the first differentiated residual is bounded by $Cm_s(|z|+s)$, using
\[
 s|\nabla m_s|\le Cm_s,
\]
while the second is $O(s^2)$. On a mixed disk use the bounds
\[
 C n_1(|z|^2+a|z|),\qquad
 C(|z|+a)\left\{\frac{a^2}{(a^2+|z|^2)^2}
                      +\frac{R^2}{R^2+|z|^2}\right\}.
\]
For the first bound use $|J_1|\le C|z|^2$ and
\[
 a|\nabla n_1|\le Cn_1.
\]
For the second, first subtract $d_0$. On $|z|\le R$ the derivatives of the rational expressions for the densities in scaled center, shape or logarithmic coordinates are bounded by a constant times the sum of the inner profile and a bounded term. On $|z|\ge R$, write
\[
 n_2=d_0(1+E_2),
\]
where the differences of the numerator and denominator from their expressions in terms of the smooth residual density are relatively $O(R^2/|z|^2)$. The same bound holds for their shape and logarithmic derivatives and for a spatial derivative multiplied by $a$, since
\[
 a/|z|\le a/R.
\]
The derivative of a smooth coefficient contributes at most its unscaled spatial derivative multiplied by $a$. These bounds also apply on the fixed cut-off annuli. Their $L^{6/5}$ norms are
\[
 O(t^{2/3})
\]
or smaller, and their absolute masses are $O(t)$, as follows by splitting the integrals at $a$ and $R$.

Integrating the unnormalized residual and its derivatives on $M$ gives $O(t)$ for each mass defect and its derivative. The differentiated normalization factors therefore preserve the $O(t^{2/3})$ bound, since each density and its scaled derivatives have $L^{6/5}$ norm at most $Ct^{-1/3}$.

\emph{Step 2. The center and shape equations.} Self-adjointness and \eqref{reg:weighted-square} now give
\[
 \langle\mathcal F_\rho(u),Y\rangle_K
 =\langle\mathcal F_\rho(u^{\rm app}),Y\rangle_K
     +O(t^{2/3}\|\eta\|_K+L\|\eta\|_K^2)
 =\langle\mathcal F_\rho(u^{\rm app}),Y\rangle_K+o(t).
\]
For an ordinary center the trial pairing is $-4\pi s_iF_i^A+o(t)$, by
\[
 \int m_{s_i}Z_k x_\ell=\pi s_i\delta_{k\ell}.
\]
On a mixed disk the first residual contributes $o(t)$ to the mixed center pairing because the outer profile has finite second moment and, on the fixed disk,
\[
 \int |z|^2|n_1-m_{R,B}|=O(R^4(1+|\log R|))
\]
by the fraction estimate above. In outer coordinates the bounded center tangent tends to zero off the origin.

The second residual contributes $-4\pi a_jF_j^B+o(t)$: on the inner scale use the same moment identity and the tangent $4Z_k$. For the passage to the inner limit, the inner tail beyond $Ha$ has
\[
 \int |z|m_a=O(a/H).
\]
The comparison in the residual proof bounds $|n_2-d_0|$ by $C m_a$ on $|z|\le R$ and by
\[
 C R^2/(R^2+|z|^2)
\]
on $|z|\ge R$. The second bound contributes $o(R^2)$ upon multiplication by $|z||Y_2|$ and integration. Indeed the two terms of the tangent bound give $O(aR)$ and $O(aR^2(1+|\log R|))$, respectively. For every center tangent
\[
 \int_M N_iY_i=o(1)
\]
by oddness on its active scale and the vanishing spectator pairings. The normalization correction to a center test is consequently $o(t)$.

At a mixed point the smooth first logarithmic error is exactly $J_1=\psi_j(z)-\psi_j(0)-\log(1+\beta_j|z|^2/4)$, by \eqref{reg:scales} and the subtraction of $L_1$. Thus
\[
 J_1(0)=J_{1,z}(0)=0.
\]
In these coordinates $\Delta=\Delta_g$ at the center, and $\Delta_gR_g=\Delta_gG=1$ for the regular self part and the other centers. Using $\Delta_g(-w/2)=\beta_j-\sigma$ there and $\Delta_g\varphi=-2K_g$ gives
\[
 \Delta\psi_j(0)=\kappa_j+\beta_j.
\]
After subtracting $\log(1+\beta_j|z|^2/4)$ the quadratic term is $\kappa_j|z|^2/4+\Re(T_jz^2)$. The shape trial pairing, to order $o(t)$, is
\[
 -\int \partial_{B_k}m_{R,B}
          \left\{\frac{\kappa}{4}|z|^2+\Re(Tz^2)\right\}
 =-2R^2\left(\kappa\,\partial_{B_k}\mathcal M
                                  +4\pi\mathbf t_k\right),
\]
where $\mathbf t_k$ is the $k$th component of $(\Re T,-\Im T)$. The equality follows by differentiating \eqref{reg:moments}. The difference $n_1-m_{R,B}$ and its shape and logarithmic derivatives are bounded by $CR^4/(R^4+|z|^4)$, as follows by writing
\[
 F_R=R^4+|z^2-R^2B|^2+\beta R^4|z|^2
\]
and using the displayed comparabilities. Its second moment on the fixed disk is $O(R^4(1+|\log R|))$. The product $(n_2-d_0)\partial_{B_k}\log n_2$ is even, so only the quadratic remainder of $e^{J_2}-1$ contributes to its integral.

Up to $R$ its second moment is bounded by the inner second moment and $CR^4$. Beyond $R$ it is bounded by the second moment of $CR^4/|z|^4$. Thus it contributes $O(R^4(1+|\log R|))$. For normalization, the leading local integrals $\int n_iY_i$ are derivatives of the constant masses
\[
 \int_{\mathbb C}n_i=8\pi,
\]
up to omitted tails with derivatives $O(R^2)$ on the fixed disks.

The background terms and the mean subtraction change these integrals by $o(1)$, by the bounds on the parameter derivatives and the Taylor expansions on each disk. Thus the normalization corrections to the shape pairing are $o(t)$.

\emph{Step 3. The scale equation.} For the scale pairing let a dot denote differentiation in $\lambda$ at fixed centers, shapes and $\mu$. Put $\rho^0=(4\pi m,\rho_2)$ and
\[
 \mathscr R=-\Delta_g\widehat u^{\rm app}-K(N-\rho^0).
\]
Differentiating the normalized functional before expanding either mass gives the exact formula
\[
\begin{split}
 \partial_\lambda J_\rho(u^{\rm app})
 ={}&2\mu-2(M_1-4\pi m)
       +\int_M K^{-1}\mathscr R\cdot\dot{\widehat u}^{\rm app}\\
    &+\sum_i\frac{M_i-\rho_i}{M_i}\dot M_i
       +\mu\int_M\dot{\widehat u}^{\rm app}_1.
\end{split}
\]
This follows from $\dot M_i=\int N_i(\dot c_i+\dot{\widehat u}^{\rm app}_i)$ and integration by parts on $M$, since $\dot c_1=2$ and $\dot c_2=0$. On a mixed disk
\[
 \dot{\widehat u}^{\rm app}_1=-2+\partial_\lambda
\log n_1,
\]
so the constant part of the first residual pairing cancels the corresponding part of $-2(M_1-4\pi m)$. The remaining contribution of the first density is
\[
 -\int \partial_\lambda m_{R,B}
       \left\{\frac{\kappa}{4}|z|^2+\Re(Tz^2)\right\}
 =-2R^2\{\kappa\mathcal M(B)+4\pi\mathbf t\cdot B\}+o(t),
\]
because $\partial_\lambda R^2=R^2$. The factor $(n_2-d_0)\partial_\lambda\log n_2$ is even, so the linear Taylor term of the second component integrates to zero. On $d=|z|\le R$, \eqref{reg:gram} gives
\[
 |\partial_\lambda\log n_2|\le C,\qquad |n_2-d_0|\le C(m_a+1).
\]
For $d\ge R$, write $Q_R=R^4q_R$; $F_R-d^4$, $q_R-(4d^2+\beta d^4)$ and the $\lambda$ derivatives of $F_R,q_R$ are all $O(R^2d^2)$. The denominators give
\[
 |n_2-d_0|+|\partial_\lambda\log n_2|\le CR^2/d^2
\]
there. Thus the weighted integrals of the product against the quadratic remainder are bounded by
\[
 C\int_{d<R}d^2(m_a+1)+CR^4\int_R^{d_0^\flat}dd/d=O(R^4(1+|\log R|)),
\]
with $d_0^\flat$ the fixed disk radius. On ordinary disks parity leaves $O(t^2L)$ by their quadratic moment.

The lost tails of the total first masses are $O(t^2)$, and the residual and scale tangent on each fixed cut-off annulus are both $O(t)$. The mass defects and their derivatives are $O(t)$, and
\[
 \int_M|\dot{\widehat u}^{\rm app}|=O(tL)
\]
by the tail bounds $Cs^2/(s^2+|z|^2)$ and $CR^2/|z|^2$ on the ordinary and mixed disks, respectively. Thus the last two terms of the exact variation are $o(t)$. At the projected solution each tested pairing is exactly $\sum_\alpha c_\alpha L_\alpha Y$. The center blocks tend to $(4\pi/3)I_2$ since
\[
 4\int m_{1,1}Z_k^2=4\pi/3.
\]
The shape block at each mixed point tends to $(32\pi/3)I_2$ by the substitution $w=y^2-B$ in the outer profile. The common radial coefficient tends to $16\pi m/3$, since the separately centered radial sources use $-2Z_0$ and the squared weighted moment of $Z_0$ is one third of the mass.

On a mixed outer scale the first scale tangent, up to its constant, is $-2Z_{0,B}+B\cdot\partial_B\log m_{R,B}$. Thus its pairing with the shape sources yields the off-diagonal shape--scale column $(32\pi/3)B$.

Radial--shape pairings in the other block vanish by the angular integral in $w$, and the remaining off-diagonal pairings vanish by parity, separation, or the ratio
\[
 a/R\to0.
\]
Shape differentiation of the inner width produces an even inner dilation, which has zero pairing with the odd inner translation sources.

The leading matrix is invertible and each pairing at $u^{\rm app}$ is $O(t)$. Hence
\[
 |c|=O(t),
\]
and multiplication by an $o(1)$ matrix gives an $o(t)$ error. The shape and center pairings give the corresponding equations in \eqref{reg:rows}. Substituting the equations for the shape parameters into the common scale pairing replaces $\mathcal M$ by
\[
 \mathcal H=\mathcal M-B\cdot\nabla\mathcal M,
\]
giving the last equation. The exact local radial identity can also be stated in canonical coordinates. Let $\xi$ be the difference of the solution and the canonical profile $\log n$, including the normalizing constants, so $-\Delta\xi=Kn(e^\xi-1)+f+Kq_c$. Write
\[
 E=\epsilon\partial_\epsilon\log n,\qquad E_\infty=(2,0).
\]
Differentiating the profile equation gives $-\Delta[K^{-1}(E-E_\infty)]=(n_iE_i)_i$. Green's identity on the fixed disk $D$ then gives
\begin{equation}\label{reg:radial-exact}
\begin{aligned}
 2\int_D n_1(e^{\xi_1}-1)
 &=\int_D K^{-1}(E-E_\infty)\cdot f
   +\int_D(E-E_\infty)\cdot q_c\\
 &\quad+\sum_i\int_D E_i n_i(e^{\xi_i}-1-\xi_i)\\
 &\quad-\int_{\partial D}
   \bigl(\partial_\nu[K^{-1}(E-E_\infty)]\cdot\xi
               -K^{-1}(E-E_\infty)\cdot\partial_\nu\xi\bigr).
\end{aligned}
\end{equation}
Indeed, in the volume pairing the subtraction of the linearized equation for $E$ gives $\sum_i E_i n_i(e^{\xi_i}-1-\xi_i)- 2n_1(e^{\xi_1}-1)$. At $\rho_1=4\pi m+\mu$ the forcing is
\[
 f=f^0-K(\mu,0)^{\mathsf T}e^\varphi
\]
as in \eqref{reg:wall-forcing}, and local centering gives $\int_D E_\infty\cdot q_c=0$ for a disk containing the source supports. \emph{Step 4. The local degree.} Write
\[
 \epsilon_j=td_j(P,Q),\qquad S=\sum_jd_j\kappa_j\mathcal H(B_j)>0.
\]
For $\mu>0$, put $t=\mu\tau$. Divide the center equations by $-3s_i$ and $-3a_j$, the shape equations by $-3\epsilon_j/(16\pi)$, and the common scale equation by
\[
 3\mu/(8\pi m)>0.
\]
The center and shape sign changes preserve orientation because their dimensions are even. The limiting map is
\[
 (P,Q,B,\tau)\longmapsto
 (F(P,Q),(\kappa_j\nabla\mathcal M(B_j)+4\pi\mathbf t_j)_j,
                                                    1-\tau S).
\]
Choose an isolating center neighborhood, a large fixed product of shape balls, and $0<\tau_-<1/S_+\leq1/S_-<\tau_+$. Take the construction constant $C_0>1/\tau_-$, so the whole scale interval satisfies
\[
 |\mu|<C_0t.
\]
The angular field points outward because
\[
 \kappa_j>4|T_j|;
\]
the radial equation has opposite signs at the two endpoints. Homotopy to a constant positive $S$ and to $\mathbf t_j=0$ has no boundary zero. The degree of this finite-dimensional map is $-\sgn\det DF$.

The matrix obtained by pairing the sources with the derivatives in $(P,Q,B,\log t)$, with translations scaled by their positive inner widths, is block triangular with diagonal blocks $(4\pi/3)I_{2n}$, $(32\pi/3)I_{2V}$ and $16\pi m/3$.

Its upper off-diagonal shape--scale block is $(32\pi/3)(B_j)_j$, and its determinant is positive. Differentiating $L_\theta\eta=0$ gives
\[
 L_\theta D_\theta\eta=-(D_\theta L_\theta)\eta=o(1)
\]
in scaled directions by \eqref{reg:modulation-bounds}, preserving transversality. Fix $\mu$. We compare the degrees on the product neighborhood.

The estimates below for blowing-up solutions show that the parametrization is injective and open and that the neighborhood isolates the solutions. Transport $E_{\theta_0}$ to $E_\theta$ along parameter segments by
\[
 \dot U=[\dot P_E,P_E]U,\qquad U(0)=I,
\]
for the orthogonal projections $P_E$ onto these spaces. The commutator has finite rank, and its norm-continuous integral makes $U-I$ compact. Use target coordinates $\mathscr S_\theta(c,v)=\sum_\alpha c_\alpha S_\alpha(\theta)+Uv$ and domain coordinates
\[
 \Phi(\theta,\zeta)=u^{\rm app}(\theta)+\eta_*(\theta)+U\zeta.
\]
At fixed $\mu$ the parameter base has compact closure. The family $U-I$ is norm-continuous with compact values, so its action on any bounded fiber ball is collectively relatively compact. Indeed, a finite operator-norm net reduces this assertion to finitely many compact operators.

The graph has compact image. After identification by $\mathscr S_{\theta_0}$, it follows that the pullback of $I-\mathcal T_\rho$ differs from the identity by a compact map: the finite parameter terms, the graph terms, $(U-I)\zeta$, and $\mathcal T_\rho(\Phi(\theta,\zeta))$ each have relatively compact image.

In the target, $\mathscr S_{\theta_0}^{-1}\mathscr S_\theta-I$ is a norm-continuous compact family by the same argument, since the source columns have finite-dimensional domain.

Its inverse has the same property. Applying the inverse along the parameter segments gives an admissible homotopy of the target that preserves every zero.

For the domain change use the oriented change-of-variables rule for Leray--Schauder degree under a change of variables of the form identity plus a compact map.

At an interior point its differential has complementary column $U$, which is identified with the identity by the same transport. The finite block, in the source coordinates, is $G_\theta^{-1}L_\theta D_\theta\Phi$, where
\[
 G_\theta=((S_\alpha,S_\beta)_K)_{\alpha,\beta}
\]
is positive definite. In fact
\[
 L_\theta\Phi=L_\theta u^{\rm app}(\theta),
\]
so this block is $G_\theta^{-1}\{L_\theta D_\theta u^{\rm app} -(D_\theta L_\theta)(\eta_*+U\zeta)\}$. The determinant of $L_\theta D_\theta u^{\rm app}$ in the scaled coordinates is positive. The bounded scaled source derivatives and the small correction radius preserve that sign throughout the neighborhood.

The differential therefore has positive orientation, as one also obtains by first eliminating the finite off-diagonal blocks and then using the transport homotopy on the complement. In these coordinates the complementary equation is $A_\theta\zeta+N_\theta(\zeta)$ with
\[
 \|A_\theta^{-1}\|\le CL,\qquad \|N_\theta(\zeta)\|\le CL^{3/2}\|\zeta\|^2.
\]
In a neighborhood of the graph with
\[
 CL^{5/2}\|\zeta\|<1,
\]
scaling away the remainder keeps every complementary zero on the graph. Evaluating the finite-dimensional equations at $\zeta=0$ preserves the same finite boundary condition.

Along the same parameter segments the invertible fiber operators deform to the one at the base point. The sign of $\det(L_\theta D_\theta u^{\rm app})$ gives the orientation of the parametrization, since passing to the source coefficients divides the determinant by
\[
 \det G_\theta>0.
\]
Excision and the product rule multiply the finite degree by the sign in \eqref{reg:inverse-index}, giving $(-1)^Ud_S$. For $\mu\leq0$, the radial equation $\mu-tS+o(t)=0$ excludes zeros. Only uniform $C^0$ remainders were used.
\end{proof}
Here $m_{q,1}$ is the density $4q^2|y|^{2q-2}/(1+|y|^{2q})^2$.
\begin{lemma}[A scalar linear estimate]
\label{hyb:scalar-inverse}
Let $q>0$ denote the profile exponent. Put $\ell(y)=\log(2+|y|)$ and fix $0<\zeta<\min(1,2q)$.
\begin{itemize}
\item[\textup{(i)}]
If $q$ is noninteger, every regular $v$ on $D_L$ with $v(0)=0$ satisfies, for $L$ sufficiently large,
\[
 \sup_{D_L}|v|/\ell\le C\left\{
 \sup_{D_L}(1+|y|)^{2+\zeta}|(-\Delta-2m_{q,1})v|
       +\operatorname{osc}_{\partial D_L}v\right\}.
\]
\item[\textup{(ii)}]
For $m=m_{1,1}$ or $m=m_B$ with bounded $B$, use $-\Delta-2m$ and replace $v(0)=0$ by $\Lambda_jv=0$, where
\begin{equation}\label{hyb:inverse-pins}
 \Lambda_jv=\int_{\C}\chi_*mZ_jv\,dA,\qquad j=0,1,2 .
\end{equation}
Here $Z_j$ are the three bounded kernel functions. Choose $\chi_*$ radial on a sufficiently large fixed disk so that
\[
 (\Lambda_iZ_j)_{i,j=0}^2
\]
is uniformly positive definite on the compact set of shape parameters. For large $L$ the supports lie in $D_L$.
\item[\textup{(iii)}]
In both cases the forcing norm may use the weight
\[
 \begin{cases}|y|^{2-2\kappa},& |y|<1,\\
 (1+|y|)^{2+\zeta},& |y|\ge1,\end{cases}
 \qquad 0<\kappa<\min(1,q),
\]
with fixed $\kappa$. For $m_B$ take $q=2$.
\end{itemize}
\end{lemma}
\begin{proof}
\noindent\textit{Step 1. The homogeneous equation.}\quad For a radial model the regular angular mode of frequency $j\ge1$ is proportional to
\[
 r^j\frac{(j+q)+(j-q)r^{2q}}{1+r^{2q}} .
\]
Indeed in $x=\log r$ the equation is
\[
 u''-j^2u+2q^2\operatorname{sech}^2(qx)u=0;
\]
differentiating $e^{jx}(j-q\tanh(qx))$ gives the formula. Near $x=-\infty$ this solution is a nonzero multiple of $e^{jx}$ to leading order, so reduction of order gives an independent solution growing as $e^{-jx}$. The displayed mode is therefore the only regular one.

If $q$ is noninteger, every such mode grows faster than a logarithm as $x$ tends to infinity. The two radial solutions are $\tanh(qx)$ and $qx\tanh(qx)-1$. Only the dilation, proportional to the first, is bounded at the origin, and its value there is nonzero.

For $q=1$ the same calculation gives the radial dilation and the two regular translation modes. For $m_B$ apply Lemma~\ref{reg:jacobi-kernel} to its quadratic density, with the coefficient
\[
 2m_B=32|y|^2/(1+|y^2-B|^2)^2
\]
in the normalized coordinates. The kernel space consists of the dilation and the two shape derivatives.

We reduce the case of logarithmic growth to the classification of bounded kernel functions. If $w$ is a homogeneous solution with
\[
 |w|\le C\ell,
\]
put $g=2mw$. It is integrable with an absolute moment of any fixed order $0<\alpha<\min(1,2q)$, since $g=O(r^{-2q-2}\log r)$ at infinity (take $q=2$ for $m_B$). For the radial model take its usual dilation $D$; for $m_B$ take the kernel function
\[
 D=4-8/(1+|y^2-B|^2),
\]
obtained from $\log(16\lambda^4/(\lambda^4+|y^2-B|^2)^2)$ by differentiating in $\log\lambda$ at one. Both have nonzero limit $D_\infty$, $D-D_\infty=O(r^{-2q})$, and $\nabla D=O(r^{-2q-1})$. Interior estimates on disks of radius comparable to $r$ give
\[
 |\nabla w|=O(r^{-1}\log r).
\]
The regular solutions $w$ and $D$ satisfy their equations on the whole disk in the distributional sense. Since they are bounded near the origin and
\[
 m=O(r^{2q-2})
\]
there (with $q=2$ for $m_B$), their right-hand sides belong locally to some $L^p$ with $p>1$. They lie in $W^{2,p}_{\rm loc}$ and hence in $C^0_{\rm loc}\cap H^1_{\rm loc}$. Green's identity in its weak form therefore gives zero boundary pairing on every centered disk. Using the preceding derivative bounds and
\[
 \int_{\partial D_d}\partial_\nu w=-\int_{D_d}g,
\]
its limit as $d$ tends to infinity is $-D_\infty\int g$. Thus
\[
 \int g=0.
\]
The logarithmic potential of $g$ on the plane is bounded. For large $|y|$ subtract $\log|y|$ in the integral. On $|z|<|y|/2$ the logarithmic ratio is at most $C(|z|/|y|)^\alpha$. Use the absolute moment of order $\alpha$ there and the pointwise decay on the remaining annuli and at infinity.

The difference between $w$ and that potential is harmonic on the plane and of logarithmic growth, hence constant. Thus the solution is bounded.

The condition at zero or the three moment conditions then eliminate the kernel. In the latter case, linear independence of the three functions gives
\[
 (\Lambda_iZ_j)>0.
\]
Continuity on the compact set of $B$'s and uniform integrability of the tails allow one fixed sufficiently large cut-off.

\noindent\textit{Step 2. The estimate on a disk.}\quad For the disk estimate use
\[
 G_L(y,z)-G_L(0,z)=\frac1{2\pi}
   \left(\log\frac{|z|}{|y-z|}+\log|1-y\overline z/L^2|\right).
\]
We show why its absolute value is bounded after integration against the weights in the estimate. For the first logarithm split into $|z|<|y|/2$,
\[
 |y|/2\le |z|\le2|y|,\qquad |z|>2|y|.
\]
On the first region use $C+\log(2+|y|)+|\log|z||$, on the middle region scale by $|y|$ in the integrable logarithmic singularity, and on the last region use $C|y|/|z|$. When $|y|\le2$, bound the logarithmic singularities in a fixed disk directly. For the reflected term it suffices to consider
\[
 |y|,|z|>L/2,
\]
since otherwise it is bounded. There
\[
 |1-y\overline z/L^2|\ge |y-z|/L,
\]
as follows by squaring. Rescaling by $L$ bounds the logarithmic integral by $CL^2$ times the supremum of the weight on the outer annulus. These estimates give a bound $C\ell(y)$ for the integral against $(1+|z|)^{-2-\zeta}$ uniformly in $y,L$.

For the additional origin weight $|z|^{2\kappa-2}$ on $|z|<1$, use H\"older with a fixed exponent $p>1$ close to one and the uniform integrability of the logarithmic kernel in $L^{p'}$ there. Applied only to $|z|>M$ with weight $m(z)\ell(z)$, dominate this weight by $CM^{-c}(1+|z|)^{-2-\zeta}$ with
\[
 0<c<2q-\zeta.
\]
The resulting integral divided by $\ell(y)$ tends uniformly to zero.

A counterexample normalized to unit logarithmic norm now converges locally to a homogeneous solution by $W^{2,p}$ estimates for a fixed $p>1$ close to one, also for the allowed origin weight.

If the shape varies, pass to a limit in its compact set. The local convergence preserves the normalization conditions, so the limit vanishes.

In the disk representation subtract the value at zero. The harmonic difference is bounded by the boundary oscillation. The forcing integral tends to zero in logarithmic norm by the kernel bound for the forcing weight, including its origin part when used.

For the $2mv$ integral use local convergence on $D_M$, then the uniform small tail outside $D_M$. The value at zero tends to zero by local convergence as well. This contradicts unit logarithmic norm and completes the estimate.
\end{proof}
\begin{lemma}[Interpolation at different concentration scales]
\label{journal:independent-radii}
Let $R_a\to0$, $1\leq a\leq J$, and let $e_a$ be signed densities with
\[
 \|e_a\|_1\leq C R_a^{c_a}(1+|\log R_a|)^N,
 \qquad \|e_a\|_\infty\leq C R_a^{-K_a},
 \qquad c_a>0,\quad K_a\geq0.
\]
There exist $p>1$ and $c>0$, independent of the relative radii, such that, with $D=\max_a R_a$,
\[
 \Big\|\sum_a e_a\Big\|_p\leq CD^c,
 \qquad
 \Big\|G\sum_a e_a\Big\|_{C^0\cap H^1}\leq CD^c.
\]
The exponents form a fixed finite set.
\end{lemma}
\begin{proof}
Choose $p>1$ so close to one that $K_a(p-1)<c_a/2$ for every $a$. Interpolation gives
\[
 \|e_a\|_p^p\leq \|e_a\|_\infty^{p-1}\|e_a\|_1
 \leq C R_a^{c_a-K_a(p-1)}(1+|\log R_a|)^N.
\]
After taking the \(p\)-th root, decrease each positive exponent to absorb the fixed logarithm. Since
\[
 R_a\leq D,
\]
take the minimum of the resulting exponents and sum. The Green estimate follows from \(G=(-\Delta_g)^{-1}\Pi:L^p\to W^{2,p}\) and the two-dimensional embeddings \(W^{2,p}\hookrightarrow C^0\cap H^1\) for \(p>1\).

For a direct estimate of the Green potential, suppose \(\delta=\|e\|_1>0\) and put \(Q=\|e\|_\infty\) and
\[
 \sigma=(\delta/(Q+\delta))^{1/2}.
\]
In a disk of radius \(\sigma\) about the observation point the absolute logarithmic integral is at most \(C Q\sigma^2(1+|\log\sigma|)\). On the complement it is at most
\[
 C\delta(1+|\log\sigma|).
\]
If \(\sigma\) exceeds the fixed coordinate radius, the same bound follows directly with that coordinate radius, since then \(Q=O(\delta)\). Thus
\[
 \|Ge\|_\infty\leq C\delta\left(1+\log(1+Q/\delta)\right).
\]
For fixed \(Q\), the right-hand side is increasing in \(\delta\): the derivative of \(\delta\log(1+Q/\delta)\) is
\[
 \log(1+Q/\delta)-Q/(Q+\delta)\geq0.
\]
Hence the stated upper bounds on \(\delta\) and \(Q\) give a positive power of the individual radius, including when the actual \(\delta\) is much smaller than its bound. Finally
\[
 \|\nabla Ge\|_2^2=\int_M\Pi e\,Ge=\int_Me\,Ge \le\|e\|_1\|Ge\|_\infty,
\]
since \(Ge\) has zero mean. The case \(\delta=0\) is immediate.
\end{proof}
\begin{lemma}[Comparison with the scalar solution]
\label{journal:scalar-neighborhood}
Fix $p>1$. Let the centers $p_a$ range in a compact set of separated configurations. Let $w(p)$ be a branch, continuous in $C^0$, of
\[
 -\Delta_g w=2\gamma\left(\frac{H_pe^w}{\int_MH_pe^w}-1\right),
 \qquad \int_Mw=0,\qquad\gamma>0.
\]
Assume the linearization has a uniformly bounded inverse between the mean-zero $L^p$ and $W^{2,p}$ spaces. Suppose $H_p\ge0$ depends continuously on the centers in $L^p$ and satisfies
\[
 \|H_p\|_{L^p}\le C,\qquad H_p(x)\ge c>0\quad(x\in O)
\]
on a fixed nonempty open set $O$. For a sequence of Toda solutions, let $f_i$ be the normalized densities. Obtain $\nu_i$ by cutting off the solution densities or the fitted profiles, and set
\[
 e_1=f_1-\nu_1,\quad
 Z=K(G\nu_1,G\nu_2)^t,\quad w_u=u_2-Z_2.
\]
\begin{itemize}
\item[\textup{(i)}]
Suppose the fitted centers converge to a point $p_*$ of the branch, $\|e_1\|_p=o(1)$, and
\[
 f_2-\nu_2\longrightarrow
 \frac{\gamma H_{p_*}e^{w(p_*)}}{\int_M H_{p_*}e^{w(p_*)}}
 \quad\text{strongly in }L^p.
\]
Then $w_u$ belongs, eventually, to the $W^{2,p}$ neighborhood of $w(p)$ given by the inverse function theorem.
\item[\textup{(ii)}]
Assume \textup{(i)} and suppose in addition that
\[
 f_2-\nu_2=C_uH_pe^{w_u}+e_2,\qquad
 \|e_1\|_p+\|e_2\|_p+
 \left|\int_M\nu_2-(\rho_2-\gamma)\right|\leq\delta_D,
 \qquad\delta_D\to0.
\]
Then $\|w_u-w(p)\|_{W^{2,p}}\leq C\delta_D$.
\end{itemize}
With these definitions the second Toda equation gives
\begin{equation}\label{journal:scalar-remainder-equation}
 -\Delta_g w_u=2\Pi(f_2-\nu_2)-\Pi e_1.
\end{equation}
In particular, the remainder $w_u$ depends on both densities cut off near the blow-up points through $Z_2=-G\nu_1+2G\nu_2$.
\end{lemma}
\subsection{Description of all mixed blowing-up solutions}
We next prove that every blowing-up sequence with the prescribed limiting shadow solution belongs to the family constructed above.

We first derive local estimates without imposing a relation between the mass parameter and the concentration scale. In the following lemma and its proof, the local canonical variables and densities are
\[
 v_i=u_i+\log\frac{\rho_i h_i}{\int_Mh_i e^{u_i}}+\varphi,
 \qquad F_i=e^{v_i}
       =e^\varphi\frac{\rho_i h_i e^{u_i}}{\int_Mh_i e^{u_i}}.
\]
Thus $F_i\,dA$ is the normalized solution measure in the conformal coordinate. The notation $F_i$ in these local estimates is distinct from the center field $F$ in \eqref{reg:forces}. Let $D=\max_jR_j\vee\max_i\sqrt{s_i}$ before any scale comparison. Write
\[
 \mathfrak n_i=\rho_i h_i e^{u_i}/\int_Mh_i e^{u_i}
\]
for the densities of the solution with respect to the surface measure. Choose the fixed cut-offs $\chi_a$ about the first-component blow-up points and inner cut-offs $b_j$, equal to one on $D_{\tau_j}$ and supported in $D_{2\tau_j}$, where $\tau_j=R_j^{3/2}$. Put
\begin{equation}\label{reg:actual-diffuse}
 \begin{gathered}
 \nu_1=\sum_a\chi_a\mathfrak n_1,\qquad
 b=\sum_jb_j,\qquad\nu_2=b\mathfrak n_2,\qquad
 e_1=\mathfrak n_1-\nu_1,\quad d_D=\mathfrak n_2-\nu_2,\\
 Z_2=-G\nu_1+2G\nu_2,\qquad w_u=u_2-Z_2,\qquad
 -\Delta_gw_u=2\Pi d_D-\Pi e_1.
 \end{gathered}
\end{equation}
The cut-offs are applied to the solution densities. The last equation follows directly from the second Toda equation.
\begin{lemma}[Local power estimates]
\label{reg:local-power}
Let $p$ be an isolated regular blow-up point. Assume that the second component has a finite positive residual limit on the punctured matching disk. Fix any sufficiently small $\varepsilon>0$.
\begin{itemize}
\item[\textup{(i)}]
At a mixed point, assume $\kappa>4|T|$ and take $R,a$ as in Lemma~\ref{reg:stress}. Then $a=R^{2+o(1)}$ and
\[
 F_1\leq CR^{-4+\varepsilon}r^{2-\varepsilon},\quad
 F_2\leq Ca^{2-\varepsilon}r^{-4+\varepsilon}
                 \quad(Ca<r<R/C),
\]
\[
 F_1\leq CR^{4-\varepsilon}r^{-6+\varepsilon},\quad
 F_2\leq Cr^{-\varepsilon}\quad(CR<r<r_0).
\]
\item[\textup{(ii)}]
Under \textup{(i)}, put $\tau=R^{3/2}$. Then
\[
 \text{inner mass at }\tau=4\pi+O(R^{1/4}),\qquad
 \text{local first mass}=8\pi+O(R^{1/4}).
\]
Truncate the inner second density at $\tau$ as in \eqref{reg:actual-diffuse}, retaining its residual part. The mean-zero Green potentials of the concentrating measures differ from those of the fitted profiles by $O(R^{1/4}|\log R|^C)$ uniformly.
\item[\textup{(iii)}]
At an ordinary point with local masses $(4\pi,0)$, on the fixed matching disk we have
\[
 F_1\leq Cs^{-2}(1+r/s)^{-4+\varepsilon},\qquad
 F_2\leq C(s+r)^{2-\varepsilon}.
\]
For $|y|\leq s^{-1/2}$,
\[
 \left|\log\frac{s^2F_1(p+sy)}{4(1+|y|^2)^{-2}}\right|
       \leq Cs^{1/4}\log(2+|y|),
\]
and the error in the Green potential is $O(s^{1/4}|\log s|^C)$.
\end{itemize}
\end{lemma}
\begin{proof}
\emph{Step 1. Estimates on the two annuli (i).} Circular integration of the canonical equations about the selected center gives
\[
 r\overline v_1'(r)=-2m_1(r)+m_2(r)+O(r^2),\qquad
 r\overline v_2'(r)=m_1(r)-2m_2(r)+O(r^2),
\]
where
\[
 m_i(r)=(2\pi)^{-1}\int_{D_r}F_i.
\]
We estimate the derivatives uniformly on the two annuli. Given
\[
 \varepsilon>0,
\]
first choose a large fixed $C$ and a small fixed $r_0$. On $[Ca,R/C]$, the mass at each endpoint is within $O(\varepsilon)$ of $(0,2)$ in $2\pi$ units.

At the inner endpoint this follows by exhausting the inner scalar profile, and at the outer endpoint it follows by excluding a small disk from the outer profile and by the already exhausted second mass. Monotonicity in the radius then gives the estimate throughout that interval.

On $[CR,r_0]$ the masses are similarly within $O(\varepsilon)$ of $(4,2)$: exhaust the first outer scalar profile at $CR$, and take $r_0$ small enough that the remaining mass is within the prescribed tolerance.

The error in the derivatives can thus be taken smaller than $\varepsilon$, upon choosing the mass tolerance by a fixed multiple of $\varepsilon$. The derivatives on the short transition annulus $[R/C,CR]$ are uniformly bounded by the total mass bounds.

On a fixed profile circle at scale $a$, the second mean is $-2\log a+O(1)$. On fixed circles at scale $R$, the first mean is $-2\log R+O(1)$. On the fixed circle
\[
 r=r_0,
\]
the second mean is bounded, by its smooth positive exterior limit. Integrating the second-component derivative on the middle interval and comparing its two endpoints using the outer interval and the transitions yields
\[
 \big|2\log a-4\log R\big|
       \leq C_\varepsilon+
             C\varepsilon\big(|\log a|+|\log R|\big).
\]
Divide by the sum of the two absolute logarithms, pass to the sequence limit and then let $\varepsilon$ decrease to zero. This gives
\[
 a=R^{2+o(1)}.
\]
Finally integrate each derivative in the appropriate direction: for the first component use the outer profile circle for both estimates, and for the second use its inner profile circle for the middle interval and the fixed exterior circle for the outer interval.

Circular oscillation turns the mean estimates into the four stated pointwise bounds. The choices of $C$ and $r_0$ depend on the exponent tolerance, not on $a/R$. On
\[
 R/C\le r\le CR,
\]
the bounded radial derivatives of the circular mean and oscillation, compared with
\[
 F_2\le C R^{-\varepsilon}
\]
at $r=CR$, give the same bound throughout this transition annulus. \emph{Step 2. The local masses (ii).} We now estimate the Cauchy integrals on the two circles. Fix $0<\varepsilon<1/100$. By decreasing the scales, we may suppose that
\[
 R^{2+1/100}\leq a\leq R^{2-1/100}.
\]
The bound on the intermediate annulus for $F_1$ extends to the inner disk in the form
\[
 F_1(z)\leq C R^{-4+\varepsilon}(a+|z|)^{2-\varepsilon}
 \qquad (|z|<R/C).
\]
To extend the estimate, rescale by $a$ on a disk of radius $2C$ in the rescaled coordinates, choosing
\[
 2Ca<R/C.
\]
The inner profile convergence bounds both rescaled densities on this fixed disk, so the rescaled first field has bounded Laplacian there. Its upper bound on $|z|=2Ca$ is the one already obtained on the middle annulus.

The Poisson barrier for a bounded Laplacian extends that boundary upper bound to $|z|<Ca$ up to an additive uniform constant in the field. Exponentiating gives the stated upper bound on the inner scale. For a fixed smaller coordinate disk write
\[
 \mathcal C_i(z)=\frac1\pi\int\frac{F_i(\zeta)}{z-\zeta}\,dA_\zeta,
 \qquad I_i(r)=\int_{D_{r/2}}F_i\,dA.
\]
The remaining harmonic term and the terms from the smooth right-hand side have bounded first and second derivatives. On the annuli used below, bounded total masses and the radial flux formula bound the circular mean variation by $O(1)$ between radii $r/2$ and $2r$.

Together with the circular oscillation bound this controls the field oscillation on every local disk of radius comparable to $r$ within the annulus. In the rescaled equation $r^2F_i$ is bounded. The local Poisson gradient bound and the harmonic gradient estimate therefore give
\[
 r|\nabla v_i|\leq C,\qquad r|\nabla F_i|\leq CF_i.
\]
The inner extension just proved gives
\[
 \sup_{|z|=\tau}
 \{\tau|\mathcal C_1(z)|+\tau^2|\partial_z\mathcal C_1(z)|\}
 \leq C\frac{\tau}{R},\qquad \tau=R^{3/2}.
\]
In fact $|\partial_z\mathcal C_1|\leq CR^{-2}$ on this circle. For the part between $2\tau$ and $R$, the derivative is bounded by
\[
 CR^{-4+\varepsilon}\int_{2\tau}^R r^{1-\varepsilon}\,dr
 \leq CR^{-2}.
\]
The inner part and the principal-value part meeting the observation point are at most $CR^{-4+\varepsilon}\tau^{2-\varepsilon}$, and the remaining mass is at distance comparable to $R$ or farther. These terms are bounded by $CR^{-2}$ as well.

For the second density, its first moment inside $D_{\tau/2}$ is at most $Ca$. Indeed, the inner profile bound on $D_{Ca}$ and the bound $Ca^{2-\varepsilon}|z|^{-4+\varepsilon}$ from $Ca$ to $\tau/2$ both give this order after multiplying by $|z|$ and integrating. Splitting its Cauchy integral at $\tau/2$ and $2\tau$ gives
\begin{align*}
 &\sup_{|z|=\tau}\left\{
 \tau\left|\mathcal C_2(z)-\frac{I_2(\tau)}{\pi z}\right|
 +\tau^2\left|\partial_z\mathcal C_2(z)
                         +\frac{I_2(\tau)}{\pi z^2}\right|\right\}\\
 &\hspace{3em}\leq
 C\left\{\frac a\tau+\left(\frac a\tau\right)^{2-\varepsilon}
                         +\tau\right\}.
\end{align*}
On $D_{\tau/2}$ use
\[
 |(z-\zeta)^{-1}-z^{-1}|\leq C|\zeta|/\tau^2
\]
and its $z$ derivative. On the annulus meeting $z$, the differentiated kernel is interpreted in principal value. Subtract $F_i(z)$ on a disk of radius $\tau/4$ centered at $z$; its constant part has zero angular integral, and the remainder is bounded by $C\tau\sup|\nabla F_i|$.

On the remaining annuli the integral of $a^{2-\varepsilon}r^{-4+\varepsilon}$ gives $C(a/\tau)^{2-\varepsilon}$ after multiplication by the appropriate power of $\tau$. The bounded outer contribution gives $C\tau$. This proves the differentiated estimates without differentiating a weak limit of measures.

Let $\mathfrak a=(2v_1+v_2)/3$ and $\mathfrak b=(v_1+2v_2)/3$. Then $\mathfrak a'=-\mathcal C_1/4+H_1'$ and $\mathfrak b'=-\mathcal C_2/4+H_2'$, with the smooth terms just described. Put $b_\tau=I_2(\tau)/(4\pi)$. The preceding estimates give
\[
 \mathfrak a'=O(R^{-1}),\qquad
 \mathfrak a''=O(R^{-2}),\qquad
 \mathfrak b'=-\frac{b_\tau}{z}+O(R^{1/4}/\tau),\qquad
 \mathfrak b''=\frac{b_\tau}{z^2}+O(R^{1/4}/\tau^2).
\]
Here every error before weakening its exponent is $O(R^{1/3})$ after the indicated multiplication by $\tau$ or $\tau^2$. We now estimate the error in the integral of $W_2$. With $k_1,k_2$ as in the proof of Lemma~\ref{reg:stress}, direct differentiation gives
\[
 \bar\partial W_2=k_1'+k_2'
             +k_1(-2\mathfrak a'+\mathfrak b')
             +k_2(\mathfrak a'-2\mathfrak b').
\]
The exponential terms cancel identically. The $k_i$ have bounded $C^1$ norms, and
\[
 \sup_\zeta\int_{D_r}\frac{|z|}{|z-\zeta|}\,dA_z\leq Cr^2.
\]
Fubini's theorem and the bounded local masses consequently give
\[
 \left|\frac1{2i}\int_{\partial D_r}zW_2\,dz\right|\leq Cr^2.
\]
Substitution on $r=\tau$ yields $b_\tau-b_\tau^2=O(R^{1/4})$. Since $b_\tau\to1$, this proves $b_\tau=1+O(R^{1/4})$. Replacing $D_{\tau/2}$ by $D_\tau$ changes the mass by at most $C(a/\tau)^{2-\varepsilon}=O(R^{1/2})$. On the outer circle $\varrho=R^{1/2}$ the same decomposition gives, for $j=0,1$,
\begin{align*}
 \varrho^{j+1}
 \left|\partial_z^j\left(\mathcal C_1(z)
                      -\frac{I_1(\varrho)}{\pi z}\right)\right|
 &\leq C\left(\frac R\varrho
           +R^{4-\varepsilon}\varrho^{-4+\varepsilon}\right),\\
 \varrho^{j+1}
 \left|\partial_z^j\left(\mathcal C_2(z)
                      -\frac{I_2(\varrho)}{\pi z}\right)\right|
 &\leq C\left(\frac a\varrho+\varrho^{2-\varepsilon}+\varrho\right).
\end{align*}
The second mass between $\tau/2$ and $\varrho/2$ is $O((a/\tau)^{2-\varepsilon}+\varrho^{2-\varepsilon})$. Thus, with $\alpha=I_1(\varrho)/(4\pi)$ and $b=I_2(\varrho)/(4\pi)$, one has $b=1+O(R^{1/4})$ and
\[
 \alpha+b-\alpha^2+\alpha b-b^2=O(R^{1/2}).
\]
At $(\alpha,b)=(2,1)$ the derivative of this polynomial with respect to $\alpha$ is $-2$. The qualitative limit $\alpha\to2$ therefore gives
\[
 \alpha=2+O(R^{1/4}).
\]
The first tail outside $D_{\varrho/2}$ has mass $O(R^{4-\varepsilon}\varrho^{-4+\varepsilon})=O(R^{3/2})$. This establishes the stated first mass error on a fixed disk.

\emph{Step 3. Comparison with the scalar profiles (ii).} On expanding disks we use the scalar estimate
\[
 \sup_{|y|<L}\frac{|z(y)|}{\log(2+|y|)}
 \leq C\left\{\sup_{|y|<L}(1+|y|)^{2+\sigma}
           |(-\Delta-2m)z|+\operatorname{osc}_{\partial D_L}z\right\},
 \quad0<\sigma<2,
\]
after fixing three orthogonality conditions defined by compactly supported test functions. Specifically, choose a fixed large cut-off $\chi$ and impose
\[
 \int\chi mZ_jz=0.
\]
The matrix $(\int\chi mZ_iZ_j)_{i,j}$ is positive definite for a sufficiently large fixed cut-off, uniformly for bounded $B$. This is the $q=1$ and bounded-$B$ case of Lemma~\ref{hyb:scalar-inverse}, whose proof is independent of the present blow-up estimates.

The Green identity with the dilation mode bounds any logarithmically bounded homogeneous limit. Orthogonality then eliminates its three coefficients of the kernel functions. The Green kernel estimate in the lemma excludes concentration near the boundary of the expanding disk.

If $\sigma\geq1$, apply it first with a smaller exponent below one; the displayed right-hand norm only increases. The nonlinear remainder has norm at most $C\|z\|_{\log}^2$ when
\[
 \|z\|_{\log}=o(1).
\]
We first prove a smallness estimate on each expanding scalar disk. Let $v$ be the regular logarithmic part, $v_*$ its fitted plane profile, and $m$ the corresponding model density in rescaled coordinates.

The difference of their circular fluxes is twice the enclosed mass difference plus the integrated weighted forcing.

For a fixed large $M_0$ the inner circle $|y|=M_0$ has $v-v_*=o(1)$ by local convergence. The bound for the scalar density gives a tail bound of order $M_0^{-c}$, uniform up to the outer boundary. The model has the same property.

The weighted forcing is integrable with norm $o(1)$, so the absolute difference of derivatives on $M_0<|y|<L$ is at most $CM_0^{-c}+o(1)$. The uniformly bounded angular oscillations then give
\[
 \limsup\|v-v_*\|_{\log}
 \le CM_0^{-c}+C/\log(2+M_0).
\]
Take $M_0\to\infty$ to obtain the initial $o(1)$ error, before using the inverse.

Fit the model while keeping the circle fixed in the original variables. Up to invertible constant changes of coordinates, the derivative of the three orthogonality equations is
\[
 \left(\int\chi m Z_i Z_j\right)_{i,j}.
\]
The three kernel functions are independent, so the orthogonality conditions can be imposed.

The implicit function theorem uses only the $o(1)$ error on the fixed compact support of the test functions. For the inner core this is local scalar convergence, and for the outer core it follows across $0$ by the Poisson comparison below.

The centers are first fitted at the inner width on a fixed disk in the original coordinates and then recentered by $o(a)$.

The contour bounds persist at comparable radii about that center, and this shift is $o(R)$ in the outer scaling. The two outer coefficients and the dilation are fitted in fixed coordinates about this center.

Apply this inverse first to the inner equation on $D_\tau$. The required boundary estimate follows from the Green representation. At $\tau=R^{3/2}$ the Cauchy field of the first density is $O(R^{-1})$; the second field minus its inner radial monopole is bounded by
\[
 Ca/\tau^2+C\tau^{-1}(a/\tau)^{2-\varepsilon}+C.
\]
Its angular variation is thus bounded by
\[
 C(\tau/R+a/\tau+(a/\tau)^{2-\varepsilon}+\tau)=O(R^{1/3}),
\]
using $a=R^{2+o(1)}$ and a fixed sufficiently small $\varepsilon$. Density subtraction controls the differentiated Cauchy kernel as above. We estimate the forcing in the inner variables. Write
\[
 z=ay,\qquad L_a=\tau/a,
\]
and choose $0<\sigma<1/10$. The preceding extension of the first density gives
\begin{align*}
 &\sup_{|y|<L_a}(1+|y|)^{2+\sigma}a^2F_1(ay)\\
 &\qquad\leq
 C R^{-4+\varepsilon}a^{4-\varepsilon}
           (1+\tau/a)^{4+\sigma-\varepsilon}
 \leq C R^{-4+\varepsilon}a^{-\sigma}
                    \tau^{4+\sigma-\varepsilon}
 \leq CR^{4/3}.
\end{align*}
The last inequality follows from $\tau=R^{3/2}$ and $a\geq R^{2+1/100}$; the exponent before weakening it is
\[
 2-\varepsilon/2-(1/2+1/100)\sigma>4/3.
\]
The smooth forcing contributes at most
\[
 Ca^2(1+\tau/a)^{2+\sigma}\leq CR^{4/3}
\]
as well. The inner boundary oscillation is bounded by the contour estimates above. After imposing the three orthogonality conditions on the fixed compact set, the scalar inverse therefore gives an $O(R^{1/3})$ logarithmic error.

The quadratic term is absorbed only after qualitative convergence and the integrated derivatives have given an $o(1)$ logarithmic norm; no mass rate is used to obtain that initial smallness.

Solve
\[
 \Delta\Phi=F_2-m_a-f_1
\]
with zero boundary on $D_{R^{1/2}}$, where $f_1$ is the current forcing in \eqref{reg:canonical}. The inner logarithmic estimate and finite logarithmic moments give an $L^1$ error $O(R^{1/3})$. Outside the inner estimate disk the tails have mass $O((a/\tau)^{2-\varepsilon}+R^{1-\varepsilon/2})$, including the transition annulus.

The densities are bounded by a fixed power of $R^{-1}$. Splitting the Dirichlet logarithmic kernel at a sufficiently small power of $R$ and including the smooth forcing yields
\[
 \|\Phi\|_\infty\leq CR^{1/3}|\log R|^2=O(R^{1/4}).
\]
Put $\delta=a/R=R^{1+o(1)}$ and
\[
 V(y)=v_1(Ry)+2\log R-\log(|y|^2+\delta^2)-\Phi(Ry).
\]
Then $-\Delta V=2(|y|^2+\delta^2)e^{\Phi(Ry)}e^V$ on $D_{R^{-1/2}}$. The annular bounds and the inner rescaling give
\[
 \begin{aligned}
 e^V&\leq C(|y|+\delta)^{-\varepsilon} &&(|y|\leq1),\\
 e^V&\leq C|y|^{-8+\varepsilon}
       &&(1\leq|y|\leq R^{-1/2}).
 \end{aligned}
\]
Indeed, on $|z|\leq Ca$ the inner profile bounds the scaled second density. The first mass is small, so the exponential-potential estimate and interior elliptic estimates bound the oscillation of the first field on this disk.

Comparison with its values on $|z|=Ca$ extends the bound on the intermediate annulus through the inner disk. For $0<\sigma<1$ and sufficiently small $\varepsilon$, it follows that
\[
 \sup_{|y|\leq R^{-1/2}}(1+|y|)^{2+\sigma}
 \left|\bigl\{|y|^2(e^{\Phi(Ry)}-1)
                 +\delta^2e^{\Phi(Ry)}\bigr\}e^V\right|
 \leq C\bigl(\|\Phi\|_\infty+\delta^{2-\varepsilon}\bigr).
\]
These bounds also justify convergence of the regular part at the origin. For
\[
 \delta\leq\eta<1,
\]
its source on $D_\eta$ is bounded by $C\eta^{2-\varepsilon}$, and its zero-boundary Green potential there has supremum at most $C\eta^{4-\varepsilon}$, uniformly in $\delta$. Qualitative convergence away from zero therefore extends across zero by harmonic extension from $\partial D_\eta$, first letting $R\to0$ and then
\[
 \eta\to0.
\]
Together with the integrated tail derivatives, this gives an initial $o(1)$ logarithmic error relative to $\log\{16/(1+|y^2-B|^2)^2\}$.

The boundary oscillation needed for the outer scalar inverse also follows from the outer Cauchy estimates just proved. The coefficients $I_i(\varrho)$ of the radial terms are real, so they disappear from
\[
 \partial_\theta v_1=2\Ree(iz\partial_zv_1).
\]
Integrating around $|z|=\varrho$ gives
\[
 \operatorname{osc}_{\partial D_\varrho}v_1
 \leq C\left\{\frac R\varrho+\frac a\varrho
 +R^{4-\varepsilon}\varrho^{-4+\varepsilon}
 +\varrho^{2-\varepsilon}+\varrho\right\}
 \leq CR^{1/2}.
\]
On this circle $\Phi=0$ and the logarithm in the definition of $V$ is radial. The outer model has boundary oscillation at most
\[
 C(R/\varrho)^2=CR,
\]
uniformly for bounded $B$. Thus the boundary oscillation of its logarithmic error is $O(R^{1/2})$, which is smaller than the $O(R^{1/4})$ error used in the scalar estimate.

In the initial $R$-scaled coordinates a dilation $\lambda$ changes the regular outer profile to $-4\log\lambda+\log(16/(1+|(y/\lambda)^2-B|^2)^2)$; its two $B$ derivatives and dilation derivative give the derivatives used in the three orthogonality conditions.

Apply the scalar inverse after fitting these parameters and rewriting at the fitted scale. It gives the outer error $O(R^{1/4}|\log R|^2)\log(2+|y|)$.

At the quadratic zero subtract the two densities additively: the main difference is $|y|^2(e^{\Phi(Ry)+V}-e^{V_*})$, and the extra term $\delta^2 e^{\Phi(Ry)+V}$ has $L^1$ size $O(\delta^2)$ by the two upper bounds for $e^V$.

The finite logarithmic moments give $L^1$ error $O(R^{1/4}|\log R|^C)$, including the tails. The logarithmic kernel split and the polynomial density bound then give the stated Green errors.

\emph{Step 4. The ordinary blow-up point (iii).} At an ordinary point the only selected core is the first scalar bubble. Indeed each selected center contributes a full regular scalar mass, and the available masses are $(4\pi,0)$.

The selection bound and the centered Cauchy estimate give bounded angular oscillation, and the first profile exhausts the first mass.

For $Cs<r<d_0$, the derivatives of the two circular means lie within $\varepsilon$ of $-4$ and $2$, by the enclosed mass limits. Integrating the first from its scalar core and the second from a fixed exterior circle gives the stated tails.

On $D_{2Cs}$ the two rescaled densities are uniformly bounded by the selection and the first scalar profile. Thus $\Delta_y[v_2(p+sy)+2\log s]$ is uniformly bounded on $D_{2C}$. At its boundary the outer estimate gives
\[
 v_2(p+sy)+2\log s\leq (4-\varepsilon)\log s+C'.
\]
The Poisson barrier extends this upper bound to $D_C$ up to a constant, yielding
\[
 F_2\leq C s^{2-\varepsilon}
\]
on $D_{Cs}$. The boundary oscillation and local interior gradient estimate, after subtracting the boundary mean, bound the gradient of the rescaled field on $D_C$. Away from this disk use the same gradient estimate on fixed-ratio annuli as above. Hence $|\nabla F_2|\leq C(s+r)^{1-\varepsilon}$.

Take $d=2\sqrt{s}$ and subtract the radial term due to the first mass in $D_{d/2}$. This choice leaves room for the disk of radius $\sqrt{s}$ after adjustment of the center.

Its relative Cauchy error is $O(s/d+(s/d)^{2-\varepsilon})$: the first moment is $O(s)$ by the first tail, the adjoining annulus uses the density bound and integrability of $1/|x-z|$, and the distant tail is integrated radially.

The second potential has uniformly bounded first derivative on a smaller fixed disk since $F_2$ is uniformly bounded there. Thus the oscillation on the boundary is $O(\sqrt s)$ after subtracting the radial term.

On this expanding disk in the rescaled variable (a fixed multiple of $s^{-1/2}$), the weighted supremum of the scalar forcing is bounded by
\[
 C\{s^{4-\varepsilon}(1+s^{-1/2})^{4+\sigma-\varepsilon}
        +s^2(1+s^{-1/2})^{2+\sigma}\}=O(s^{1/2}).
\]
Here the powers before weakening are $2-(\sigma+\varepsilon)/2$ and $1-\sigma/2$, respectively, with $0<\sigma<1$ and sufficiently small
\[
 \varepsilon>0.
\]
The scalar profile exhausts the local first mass, and the upper bound on the annulus bounds the rescaled tail by $O(H^{-2+\varepsilon})$ beyond radius $H$. Together with the weighted forcing bound this verifies the flux hypothesis of the argument giving the initial $o(1)$ estimate.

Choose the three ordinary parameters by the orthogonality conditions on the fixed support described above.

The fitted center and width differ from the initial ones by $o(s)$ and a relative $o(1)$, respectively, so the smaller fitted disk of radius $\sqrt s$ stays in the comparison disk and the preceding contour and forcing estimates retain their powers.

The inverse gives $O(s^{1/2})$ in logarithmic norm there. The estimate in the statement follows on weakening the power. Integration against the plane profile and the fixed-disk tail gives the $L^1$ error for the concentrating density.

Splitting the logarithmic Green kernel at a sufficiently small power of $s$ converts this to the claimed uniform Green error, using the polynomial bound on the density on its own scale.
\end{proof}
The profile estimates in Lemma~\ref{reg:local-power}, integrated against their finite logarithmic moments, compare each $\nu_i$ to the similarly truncated fitted profile in $L^1$ by $C\varrho_a^{1/4}(1+|\log \varrho_a|)^N$, where
\[
 \varrho_a=s_i
\]
or $R_j$. The densities have polynomial bounds $C\varrho_a^{-K}$; at mixed inner bubbles this follows from
\[
 a_j=R_j^{2+o(1)}.
\]
Choose once and for all a fixed $1<p_0<6/5$ sufficiently close to one that interpolation at every individual radius retains a positive power. Lemma~\ref{journal:independent-radii} then bounds the sum of their Green errors in $C^0\cap H^1$ by $CD^c$ for some $c>0$.

This also controls a local mass error times its own logarithmic width, before summing; no ratio between different radii is used.

This coarse exponent $p_0$ is distinct from the construction exponent $6/5$ in Lemma~\ref{reg:residual}. The first tails on the fixed cut-off annuli and exterior Green oscillation give
\[
 \|e_1\|_{p_0}\leq CD^c.
\]
The retained second inner tail has $p_0$th-power bound
\[
 C a_j^{p_0(2-\varepsilon)}
       \tau_j^{2-p_0(4-\varepsilon)}
       =R_j^{3-2p_0+o(1)+O(\varepsilon)},
\]
which tends to zero as a positive power of $R_j$. The outer bound $Cr^{-\varepsilon}$ has uniformly integrable $p_0$th power, and the ordinary second density has its quadratic-zero bound.

Qualitative convergence off the bubbles therefore gives strong $L^{p_0}$ convergence of $d_D$ to the selected shadow density. The exact equation \eqref{reg:actual-diffuse} now places $w_u$ in the mean-zero $W^{2,p_0}$ neighborhood of $w(P)$, before any quantitative inversion. For that inversion set
\[
 A_D=(1-b)h_2e^{Z_2},\qquad
 C_u=\frac{\rho_2}{\int_Mh_2e^{u_2}}.
\]
The identity $d_D=C_uA_De^{w_u}$ is exact. The fitted comparison for an ordinary bubble is estimated only after multiplication by its zero:
\[
 H_P=r^2k(x),\qquad
 r^2k(x)(1+s_i^2/r^2)-r^2k(x)=s_i^2k(x).
\]
To use the fixed-disk model potentials, first extend the truncated fitted inner profile from $D_{2\tau_j}$ to the fixed coordinate disk. Its omitted mass is $O(a_j^2/\tau_j^2)$, and its Green error is
\[
 O\!\left(\frac{a_j^2}{\tau_j^2}
       (1+|\log a_j|+|\log\tau_j|)\right)
       =R_j^{1+o(1)}(1+|\log R_j|).
\]
At a fixed cut-off the errors in the first profile also decay as a positive power of the radius. Thus this change of fitted model preserves the $CD^c$ error. For a mixed bubble the resulting model difference potential is
\[
 2Gm_{a_j}-Gm_{R_j,B_j}
 =\log\frac{R_j^4+|z^2-R_j^2B_j|^2}{(a_j^2+|z|^2)^2}
       +O(R_j^2|\log R_j|).
\]
Outside $D_{\tau_j}$ the ratio in the logarithm differs from one by at most $C\{R_j^2/r^2+R_j^4/r^4\}$, since $a_j\ll R_j$ and
\[
 r\geq \tau_j\gg a_j.
\]
Its two $L^{p_0}$ norms on the fixed disk are bounded by $CR_j^{3/p_0-1}$ and $CR_j^{3/p_0-2}$: integrate the two radial powers starting at
\[
 \tau_j=R_j^{3/2}.
\]
This comparison concerns the exponential ratio, even in the region where the ratio is large. The positive-power additive error in the logarithm multiplies the same bounded-in-$L^{p_0}$ ratio. Inside the deleted disk the bounded limiting coefficient has $L^{p_0}$ norm
\[
 O(\tau_j^{2/p_0}).
\]
The transition with
\[
 0\leq b_j\leq1
\]
satisfies the same bounds. Multiplication by the $C^0$-small difference of the two Green potentials therefore gives
\[
 \|A_D-H_P\|_{p_0}
 +\left|\int_M\nu_2-4\pi V\right|\leq CD^c.
\]
The factors $H_P$ and $h_2$ are included in all coefficient estimates, which therefore remain valid at their zeros. Since $w_u$ is already $C^0$ bounded, a uniform positive integral of $A_D$ on a fixed open set bounds $C_u$ above.

The lower bound follows from $\int_M d_D\to\sigma>0$ and the uniform $L^1$ bound on $A_D$ in the displayed coefficient comparison. Thus $e_2=C_u(A_D-H_P)e^{w_u}$ satisfies
\[
 d_D=C_uH_Pe^{w_u}+e_2,\qquad
 \|e_1\|_{p_0}+\|e_2\|_{p_0}
 +\left|\int_M\nu_2-4\pi V\right|\leq CD^c.
\]
Write $\delta_2=\int_M\nu_2-4\pi V$ and $Q=H_Pe^{w_u}/\int_MH_Pe^{w_u}$. Integration of the second density gives $C_u\int_MH_Pe^{w_u}=\sigma-\delta_2-\int_Me_2$. Thus the normalized scalar defect is exactly
\[
 -\Delta_gw_u-2\sigma\Pi Q
   =\Pi\bigl(2e_2-e_1-2(\delta_2+\textstyle\int_M e_2)Q\bigr).
\]
The weight $H_P$ is bounded in $L^{p_0}$, depends continuously on separated $P$, and has a uniform positive lower bound on an open set away from all center disks. The normalized scalar linearization is the Laplacian plus a compact perturbation from $W^{2,p_0}_0$ to $L^{p_0}_0$.

Its kernel space embeds in $C^0$ and bootstraps into the fixed-center energy kernel space, which is trivial.

Continuity on the separated branch bounds its inverse on a smaller center neighborhood. The previously established $W^{2,p_0}$ convergence and the displayed defect therefore give
\[
 \|w_u-w(P)\|_{W^{2,p_0}}\leq CD^c,
\]
as in Lemma~\ref{journal:scalar-neighborhood} and \eqref{journal:scalar-remainder-equation}. For regular points the matching follows directly from $u_2=Z_2+w_u$ and $u_1=2G\nu_1-G\nu_2-w_u/2+3G e_1/2$. Write
\[
 A=R^2,\qquad H=A^2+|z^2-AB|^2,
\]
$x=a^2+|z|^2$ and $d=1+\beta|z|^2/4$ in a mixed disk. Subtract the potentials of the point masses from the regular model potentials. After mean subtraction the difference equals the cut-off of
\[
 -\log(1+s^2/r^2)\quad\text{at an ordinary point},\qquad
 (-\log(H/r^4),-\log(x/r^2))\quad\text{at a mixed point},
\]
up to $O(D^c)$ in $C^0\cap H^1$. Indeed, applying $-\Delta_g$ leaves the model tail and the derivatives of the cut-off on a fixed annulus. Both are bounded by a positive power of the fitted widths. Set
\[
 c_1=\log(\rho_1/\int h_1e^{u_1}).
\]
Integration of $d_D=C_u A_D e^{w_u}$ and the scalar comparison give
\[
 C_u=\sigma/\int H_Pe^{w(P)}+O(D^c).
\]
After fixing the unit tangent, a quadratic change to the normalized coordinates changes the fitted measures by a positive power of $D$: their relative density and Jacobian changes are $O(r)$, which is controlled by the first moments.

In these coordinates, the ordinary first and mixed second center values give $4/s_i^2=e^{c_1+\psi_i^A(p_i)}s_i^{-4}(1+O(D^c))$ and
\[
 4/a_j^2=\beta_j R_j^4(1+|B_j|^2)a_j^{-4}(1+O(D^c)).
\]
At any fixed $y\ne0$ on the mixed outer scale, $x(Ry)=R^2|y|^2(1+O((a/R)^2))$. Comparing $R^2F_1(Ry)$ to $m_B(y)$ in $u_1$ gives $16R_j^4=e^{c_1+\psi_j(q_j)}(1+O(D^c))$. Thus with $t^2=e^{c_1}$ all of \eqref{reg:scales} hold to relative error $O(D^c)$, with
\[
 s_i,a_j\asymp t,\qquad R_j^2\asymp t
\]
and $D\asymp\sqrt t$. The estimates for the local first masses and the first density outside the disks give
\[
 |\mu|=O(D^c).
\]
The ordinary fitted cut-off correction already has the form of the trial field. In a mixed disk the two corrections differ by
\[
 K(-\log(F_R/H),-\log(Q_R/(4R^4dx)))^{\mathsf T}.
\]
Both logarithms are $O(A)$. The gradient of the first is bounded by $Cr$ for $r\le R$ and $CA^2/r^3$ for $r\ge R$. The gradient of the second is bounded by $CAr/(a^2+r^2)$. These give $O(A\sqrt{1+|\log a|})$ in $H^1$.

To replace the fitted widths by the matched ones, first regard the widths $s_i,a_j,R_j$ in the fitted logarithms as independent.

Each logarithmic width derivative is bounded in $C^0$ and $H^1$: for $s_i$ and $a_j$ this follows by differentiating $-\log(r^2+b^2)$ and integrating $(b^2r/(r^2+b^2)^2)^2r\,dr$; for $R_j$ rescale $z=R_jy$ in $-\log(R_j^4+|z^2-R_j^2B_j|^2)$.

The differentiated outer profile has bounded $y$-gradient on $|y|\leq C$ and $O(|y|^{-3})$ gradient outside, uniformly for bounded $B_j$. Cutoff and mean-subtraction terms obey the same bounds.

The logarithmic width ratios are $O(D^c)$ by the center and outer comparisons. Integrating the width tangents between the two triples and using the fitted density and scalar remainder bounds gives
\begin{equation}\label{reg:capture-seed}
 |\mu|=O(t^{\delta_0}),\qquad
 \|u-u^{\rm app}\|_{C^0\cap H^1}\leq Ct^{\delta_1},\qquad
 \delta_0,\delta_1>0.
\end{equation}
Fix reference widths at the initial fitted parameters $\theta_*$ and set $g_u(\theta)=L_\theta(u-u^{\rm app}(\theta))$ in those scaled coordinates. Its initial value is $O(t^{\delta_1})$. Write
\[
 J_*=L_{\theta_*}D_\theta u^{\rm app}(\theta_*),
\]
which has a uniform inverse bound. On a fixed small reference ball the leading blocks of $L_\theta D_\theta u^{\rm app}$ are continuous functions of the bounded shapes and the width ratios in these coordinates, with uniform error tending to zero.

Their variation can be made uniformly small by reducing that ball. Since
\[
 Dg_u=(D_\theta L_\theta)(u-u^{\rm app})-L_\theta D_\theta u^{\rm app},
\]
the bounded source derivatives and trial tangents give $\|I+J_*^{-1}Dg_u\|\leq1/2$ there for small $t$. The map
\[
 \theta\mapsto\theta+J_*^{-1}g_u(\theta)
\]
moves the center by $O(t^{\delta_1})$ and maps a ball of twice that radius into itself. Its unique fixed point in the reference ball imposes exact orthogonality within $O(t^{\delta_1})$ in the scaled parameters. The trial tangent bounds in both norms preserve the same error estimate.
\begin{theorem}[Description of mixed blowing-up solutions]
\label{reg:mixed-capture}
Let a blowing-up sequence have a separated nondegenerate shadow solution $Z$ with invertible scalar linearization. Suppose $V\geq1$ and $\kappa_j>4|T_j|$ at every mixed point.
\begin{itemize}
\item[\textup{(i)}]
After adjusting the parameters, the sequence belongs to the family of projected solutions in \eqref{reg:range}.
\item[\textup{(ii)}]
We have
\[
 \frac{\mu}{t}\longrightarrow S_*,\qquad
 S_*=\sum_jd_j(Z)\kappa_j\mathcal H(B_j^*)>0.
\]
\end{itemize}
\end{theorem}
\begin{proof}
\emph{Step 1. The projected correction (i).} The model inverse remains valid on the preliminary parameter range
\[
 |\mu|=O(t^{\delta_0})
\]
: at fixed fitted parameters $\rho_1$ rescales all first-component centered sources by the same nonzero factor, so $E_t$ is unchanged. The operator changes by $O(|\mu|L)$ by \eqref{reg:normalized-residual}, and its product with the inverse bound is
\[
 O(|\mu|L^2)=o(1).
\]
Its residual is $O(D^{4/3}+|\mu|\sqrt L)$. For $E=\|u-u^{\rm app}\|_{H^1}$ write $\xi=u-u^{\rm app}$. If $P$ is an approximation probability density, the normalized density
\[
 P_s=P e^{s\xi_i}/\int P e^{s\xi_i}
\]
has
\[
 P_s/P=1+O(\|\xi\|_\infty)
\]
on the support, uniformly for $0\leq s\leq1$. Write the bilinear form as the double integral of $(v(x)-v(y))(w(x)-w(y))/2$ against the product probability. The probability moment bound implies the norm of this form is at most $CL$ on $H^1$.

Thus the change of that form on replacing $P$ by $P_s$ is $O(\|\xi\|_\infty L)$, and its integral along the segment gives a normalized nonlinear error $O(\|\xi\|_\infty LE)$ in $H^{-1}$. The preliminary power estimate absorbs it after inversion, giving
\[
 E\leq C(D^{4/3}L+|\mu|L^{3/2}).
\]
\emph{Step 2. The scale parameter (ii).} To improve the parameter interval, test with the scale tangent of the mean-zero approximation. Write
\[
 \lambda=\log t,
\]
differentiate at fixed centers, shapes and $\mu$, and denote differentiation by a dot. Put $\rho^0=(4\pi m,\rho_2)$ and
\[
 \mathscr R=-\Delta_g\widehat u^{\rm app}-K(N-\rho^0).
\]
With the notation of \eqref{reg:scales}, $\dot c_1=2$ and $\dot c_2=0$. Differentiating the normalized energy before expanding either mass gives
\[
\begin{split}
 \partial_\lambda J_\rho(u^{\rm app})
 ={}&2\mu-2(M_1-4\pi m)
     +\int_M K^{-1}\mathscr R\cdot\dot{\widehat u}^{\rm app}\\
    &+\sum_i\frac{M_i-\rho_i}{M_i}\dot M_i
     +\mu\int_M\dot{\widehat u}^{\rm app}_1 .
\end{split}
\]
This follows by differentiating
\[
 M_i=\int_M e^{c_i}h_i e^{\widehat u_i^{\rm app}}
\]
and integrating the Dirichlet term by parts on $M$. In particular it introduces no circle boundary term. We give the bounds for this test also on the preliminary parameter interval. On ordinary disks, the first contribution to $\|\mathscr R\|_1$ is bounded by
\[
 C\int |z|m_s=O(t)
\]
and the second is $O(s^2)$. On mixed disks the two contributions are bounded by $C\int |z|^2 n_1$ and $C\int |z||n_2-d_0|$, respectively. The first is $O(R^2)$ and the second is $O(a+R^2)$: for the latter use
\[
 |n_2-d_0|\leq C\left\{
   \frac{a^2}{(a^2+|z|^2)^2}+\frac{R^2}{R^2+|z|^2}\right\}
\]
and integrate over the fixed disk. On the fixed cut-off annuli the unnormalized residual is $O(t)$, and outside the disks the only residual density is the first tail $O(t^2)$.

These bounds hold for $\dot{\mathscr R}$ as well: differentiating the rational expressions for $n_1$ and $n_2$ at fixed $z$ preserves the displayed majorants for $n_1$ and $n_2-d_0$, and
\[
 |\dot m_s|\leq 2m_s
\]
on ordinary disks. Integrating $\mathscr R$ and $\dot{\mathscr R}$, and using the invertibility of $K$, gives
\[
 |M_i-\rho_i^0|+|\dot M_i|=O(t).
\]
The derivative in $\lambda$ is uniformly bounded. On the disks its expressions are logarithmic derivatives of the rational profiles and their tails. On mixed disks outside $|z|=R$ its size is at most $C R^2/|z|^2$, and on ordinary disks it is at most $C s^2/(s^2+|z|^2)$. Hence
\[
 \int_M|\dot{\widehat u}^{\rm app}| =O(tL).
\]
The exact energy identity now gives $\partial_\lambda J_\rho(u^{\rm app})=2\mu+O(t+|\mu|tL)$.

Differentiating the residuals before taking their $L^{6/5}$ norms also gives
\[
 \|\partial_\lambda\mathcal F_\rho(u^{\rm app})\|_K
    \le C(t^{2/3}+|\mu|\sqrt L).
\]
Indeed the differentiated disk errors use the same bounds as the undifferentiated ones, while the derivatives of $M_i-\rho_i^0$ and of the densities obey the bounds just established. For the differentiated parameter term put
\[
 P=N_1/M_1.
\]
Its derivative has zero integral and $|\dot P|\leq CP$, since
\[
 |\dot N_1|\leq CN_1,\qquad |\dot M_1|\leq CM_1.
\]
The log-kernel split for the probability energy bounds $\|G\dot P\|_{H^1}$ by $C\sqrt L$, also for this signed source. Thus the factor $|\mu|$ remains in the estimate throughout the preliminary parameter interval. Since
\[
 \mathcal F_\rho(u)=0,
\]
the self-adjoint linearization and \eqref{reg:weighted-square} therefore give
\[
 |2\mu|\le C\{t+|\mu|tL
   +(t^{2/3}+|\mu|\sqrt L)E+LE^2\}.
\]
Substituting the bound for $E$ and using $(x+y)^2\le2x^2+2y^2$, we get, for large $L$,
\[
 |\mu|\le C\bigl(t+t^{2/3}L^3|\mu|+|\mu|^2L^4\bigr).
\]
The power bound absorbs the terms containing $\mu$, yielding $|\mu|=O(t)$ and
\[
 E=O(t^{2/3}L).
\]
Choose $C_0$ larger than the resulting $|\mu|/t$ bound. For $|\mu|\le C_0t$, let $\bar r=C_r t^{2/3}$ bound the projected trial residual and $C_A L$ bound the inverse. The same projected equation gives
\[
 E\le C_A L\bar r+CL^{5/2}E^2\le C_A L\bar r+E/4,
\]
since $L^{5/2}E=o(1)$. Consequently
\[
 E\le(4/3)C_A L\bar r<2C_A L\bar r,
\]
strictly inside the construction ball. Its unique projected correction is therefore the difference between the solution and the approximation at the adjusted parameters. Subtracting the original equation from \eqref{reg:projected-equation} gives
\[
 K\sum_\alpha c_\alpha h_\alpha=0.
\]
The matrix $K$ is invertible, and testing the source sum against the parameter derivatives shows $c_\alpha=0$ for every $\alpha$.

Conversely, setting all $c_\alpha=0$ in \eqref{reg:projected-equation} gives the original normalized Toda equation with the same parameters and means. Dividing the radial equation by $t$ now gives $\mu/t\to S_*>0$.
\end{proof}
\subsection{Configurations containing only ordinary bubbles}
When $V=0$, the leading radial term has order $t^2\log(1/t)$. Let the shadow solution be nondegenerate and write $s_i=td_i(P)$ with $d_i=e^{\psi_i^A(p_i)/2}/2$. In the normalized local coordinate, put
\[
 \kappa_i=\Delta\psi_i^A(p_i)
 =\Delta_g\log h_1(p_i)+8\pi m-\rho_2-2K_g(p_i).
\]
Set
\[
 \ell(P)=\sum_ie^{\psi_i^A(p_i)}
       (\Delta_g\log h_1(p_i)+8\pi m-\rho_2-2K_g(p_i)),
 \qquad L_0=\frac\pi2\ell(P^0).
\]
Assume
\[
 L_0\ne0.
\]
This can be enforced at a finite set of shadow solutions by arbitrarily small smooth perturbations of the second derivatives of $\log h_1$, with zero values and first derivatives at every free center of the finite family.

The variation of each $\ell(P)$ is a nonzero linear combination with positive coefficients. The strict inequalities at the mixed points and nondegeneracy persist.

For $\lambda=\log t$, $L=\log(1/t)$, $f=|F(P)|$ and $M=|\mu|$, the derivative at the approximate solution is
\begin{equation}\label{reg:ordinary-trial-expansion}
 \partial_\lambda J_\rho(u^{\rm app})
 =2\mu-4\pi\sum_i\kappa_i s_i^2\log(1/s_i)
   +O(t^2+t^2Lf^2+Mt^2L+t^4L^2).
\end{equation}
On $|P-P^0|\le C_0t^{1/2}$ and $|\mu|\le C_1t^2L$, the reduced equations are
\begin{equation}\label{reg:ordinary-rows}
 c_i=-3s_iF_i(P)+o(t^{3/2}),\qquad
 c_0=\frac3{8\pi m}
     \left(\mu-2\pi\sum_i\kappa_i s_i^2\log(1/s_i)\right)
        +o(t^2L).
\end{equation}
On $|P-P^0|\leq C_0t^{1/2}$, $|\mu|\leq C_1t^2L$, the ordinary residual is $O(t^{7/6})$. Indeed $|F(P)|=O(t^{1/2})$,
\[
 \|rm_s\|_{6/5}=O(t^{2/3}),\qquad \|r^2m_s\|_{6/5}=O(t^{5/3}),
\]
and normalization adds $O(t^{5/3}L)$. The projected inverse therefore gives a unique correction on a ball of radius $K t^{7/6}L$ for a fixed sufficiently large $K$. On that ball the factor in the inverted quadratic Lipschitz estimate is
\[
 O(t^{7/6}L^{7/2})=o(1).
\]
In particular $E=O(t^{7/6}L)$. On this parameter domain use the global translation test up to the point where the exact solution equation is used.

For $u^{\rm app}$ the same bounds hold for the nonlinear comparison with
\[
 u=u^{\rm app}+\eta.
\]
For the projected solution the left pairing is $\sum_\alpha c_\alpha L_\alpha Y$. For the translation in direction $j$ at $p_i$, it equals $(4\pi/3)c_{i,j}+O(t|c|)$: the diagonal term is
\[
 4\int m_{s_i} Z_j^2=4\pi/3,
\]
the off-diagonal plane pairings vanish by angular integration, and the smooth completions of the tangents are $O(t)$. The source changes in total variation by
\[
 O(t|F|+t^2 L+|\mu|)=O(t),
\]
using its exact centering. The projected equation and the uniform inverse bound for $((S_\alpha,S_\beta)_K)_{\alpha,\beta}$ give
\[
 |c|=O(t^{7/6} L^2),
\]
so $t|c|=o(t^{3/2})$. The comparison error is also $o(t^{3/2})$ by $E=O(t^{7/6}L)$ and $|F|=O(t^{1/2})$. This proves the translation equation in \eqref{reg:ordinary-rows}.

For the scale equation put $\lambda=\log t$. The trial expansion \eqref{reg:ordinary-trial-expansion}, derived by differentiating the normalized functional on the fixed center neighborhood, gives
\[
 \partial_\lambda J_\rho(u^{\rm app})
     =2\mu-4\pi\sum_i\kappa_i s_i^2\log(1/s_i)+o(t^2 L).
\]
Along the projected solution its change is bounded by
\[
 C\{(t^{2/3}|F|+t^{5/3}L+|\mu|\sqrt L)E+LE^2+|c|E\}
   =o(t^2L).
\]
The first two terms are the self-adjoint linearized pairing with $\partial_\lambda\mathcal F_\rho(u^{\rm app})$ and \eqref{reg:weighted-square}. The last term follows by differentiating the exact orthogonality condition:
\[
 L_\alpha\partial_\lambda\eta =-(\partial_\lambda L_\alpha)\eta.
\]
On each ordinary disk the scale tangent is $(-2(Z_0+1),Z_0+1)$ up to the constant from mean subtraction. The source pairings at the approximation therefore give
\[
 (16\pi m/3)c_0+O((t|F|+t^2L)|c|);
\]
the nonradial pairings of the plane profiles vanish and the constants pair to zero. This error is $o(t^2L)$ and gives the second equation in \eqref{reg:ordinary-rows}.
\begin{theorem}[Ordinary blowing-up solutions and their local degree]
\label{reg:ordinary-capture}
For every sequence concentrating at the above ordinary configuration, the following hold.
\begin{itemize}
\item[\textup{(i)}]
We have
\[
 |F(P)|+|P-P^0|=O(t^{5/3}L^3),\quad
 E=O(t^{5/3}L^2),\quad
 \mu=L_0t^2L+O(t^2).
\]
For all sufficiently small $t$, the sequence belongs to the preceding family of projected solutions on the shrinking neighborhoods of the centers.
\item[\textup{(ii)}]
Concentration occurs only for $\mu L_0>0$. On this side the local degree is $(-1)^m\sgn(\mu)d_S$.
\item[\textup{(iii)}]
For either sign of $L_0$, subtracting the local degree for $\mu<0$ from that for $\mu>0$ gives $(-1)^m d_S$.
\end{itemize}
\end{theorem}
\begin{proof}
\emph{Step 1. The preliminary estimates.} The ordinary part of Lemma~\ref{reg:local-power}, the scalar inverse, and evaluation at each first peak give $s_i^2=e^{c_1+\psi_i^A(p_i)}/4\,(1+O(t^c))$, $|\mu|=O(t^{\delta_0})$, and
\[
 \|u-u^{\rm app}\|_{C^0\cap H^1}=O(t^{\delta_1})
\]
for positive exponents, on a fixed center neighborhood. Use \eqref{reg:actual-diffuse} with
\[
 \nu_2=0,\qquad \gamma=\rho_2,
\]
and the same scalar exponent $p_0$ close to one. Before applying the inverse, impose the orthogonality conditions as for \eqref{reg:capture-seed}.

This choice of parameters preserves the error bounded by a positive power of the scale and places the difference in $E_t$. This does not change the construction estimates at
\[
 p=6/5.
\]
No shrinking center condition is used here. To extend the projected inverse first compare at the same fitted parameters with $\mu=0$. Each first-component source is multiplied by
\[
 \rho_1/(4\pi m)>0,
\]
so its kernel constraint is unchanged. The normalized bilinear form changes by $O(|\mu|L)$, and the inverse norm at zero is $O(L)$. The preliminary power estimate gives invertibility by a convergent geometric series with
\[
 O(|\mu|L^2)=o(1)
\]
before the sharp mass law is used. Absorb the normalized nonlinear error using the preliminary $C^0$ bound as in the mixed case. With $M=|\mu|$ and $f=|F(P)|=o(1)$ this gives
\[
 E\leq C(t^{2/3}Lf+t^{5/3}L^2+ML^{3/2}).
\]
The scale test on the surface gives $M\le C(t^2L+LE^2+M^2)$ on this preliminary interval. Indeed \eqref{reg:ordinary-trial-expansion} uses only the explicit approximation on the fixed center neighborhood; since
\[
 f=O(1),
\]
its right-hand side is $2\mu+O(t^2L+Mt^2L)$. Before taking any $H^{-1}$ norm, the differentiated local residual has absolute mass $O(tf+t^2L)$, from
\[
 |\dot m_s|\le 2m_s,
\]
the fixed smooth $J_i$, and integrals of $m_s|z|$ and $m_s|z|^2$.

The differentiated normalization contributes at most $O(t^2L+M)$ in the same absolute mass bound, with each normalized upper bound for the density bounded by a fixed power of $t^{-1}$. H\"older with exponent $1+1/L$ therefore gives
\[
 \|\partial_\lambda\mathcal F_\rho(u^{\rm app})\|_K \le C\sqrt L(tf+t^2L+M).
\]
Pairing the exact solution equation with the bounded scale tangent and using self-adjointness and \eqref{reg:weighted-square} bounds the comparison error by $C\sqrt L(tf+t^2L+M)E+CLE^2$. Young's inequality and $Mt^2L=o(M)$ give the asserted estimate. Since
\[
 M=O(t^{\delta_0})
\]
for some $\delta_0>0$, substitution absorbs $CM^2L^4$ and yields
\[
 M\leq C(t^2L+t^{4/3}L^3f^2),\qquad
 E\leq C(t^{2/3}Lf+t^{5/3}L^2).
\]
Here is a translation test on the fixed center neighborhood; it does not use the shrinking center domain. For a coordinate direction of $p_i$, let
\[
 Y=s_i\partial_{p_i}u^{\rm app},
\]
with $\rho$ and $t$ fixed. In the corresponding normalized coordinates put $Z_j(y)=y_j/(1+|y|^2)$. On the $i$th disk the two components of $Y$ are $(4Z_j,-2Z_j)+O(t)$, by differentiating $-2\log(s_i^2+|z|^2)$ in the translated coordinates, with
\[
 \partial_{p_i}s_i=O(s_i).
\]
Changes in the smooth coordinate factor, background, and surface mean contribute $O(t)$. On other disks and outside the source disk
\[
 Y=O(t),\qquad \int_M |Y|=O(t),
\]
since $|Z_j(z/s_i)|\le C s_i/|z|$. \emph{Step 2. The center equations.} Pairing the exact normalized residual \eqref{reg:normalized-residual} with $Y$ gives
\[
 \langle\mathcal F_\rho(u^{\rm app}),Y\rangle_K
   =-4\pi s_iF_{i,j}^A
     +O(t^2f+t^3L+tM).
\]
In fact the local first term is $-4\int m_{s_i}Z_j(e^{J_i}-1)$ and the full plane moment is
\[
 \int_{\mathbb R^2}m_{s_i}Z_j x_k=\pi s_i\delta_{jk}.
\]
Its tail outside the fixed disk is $O(t^3)$. The quadratic Taylor polynomial of $e^{J_i}-1$ is even and pairs to zero. The cubic remainder is bounded by
\[
 C\int |Z_j||z|^3m_{s_i}=O(t^3L).
\]
Errors from the $O(t)$ part of the tangent are bounded by $C t\int m_s(f|z|+|z|^2)=O(t^2f+t^3L)$. The second excess $s_i^2\gamma_i$ and the residuals on the fixed cut-off annuli give $O(t^3)$. In the normalization terms
\[
 |\int N_iY_i|=O(t)
\]
: for the first component subtract the radial mass in the odd translation pairing, and for the second use its bounded density and $\int|Y|=O(t)$. The two mass coefficients are $O(t^2L+M)$. The differentiated residual obeys
\[
 \|D\mathcal F_\rho(u^{\rm app})Y\|_K
   \le C\sqrt L\, (tf+t^2L+M).
\]
Indeed, write
\[
 J_k(z,P)=\psi_k^A(z,P)-\psi_k^A(0,P)
\]
in coordinates centered at $p_k$. On the separated scalar branch, a center derivative at fixed local $z$ is $O(|z|)$ by the bounded mixed $z,P$ derivatives of $\psi_k^A$. At a fixed surface point there is also the coordinate derivative on the moving disk, bounded by
\[
 C|\nabla_zJ_k|\le C(f+|z|),
\]
since the linear coefficient of $J_k$ is $F_k(P)$. The directions are multiplied by the ordinary width, comparable to $t$. Differentiating the profile costs at most a bounded factor, and the first absolute mass is bounded by
\[
 C(tf+t^2L)
\]
using the first and second moments of its radial density.

The second excess and cut-off derivatives have smaller mass. The scaled derivatives of $M_i$ are $O(t^2L)$ by integration in the translated disks, preserving the parity cancellation of the local linear term.

The differentiated normalization terms have absolute mass $O(t^2L+M)$ because the scaled derivatives of the approximate densities are bounded by a constant times those densities.

Every upper bound for the density in these estimates, divided by its absolute mass bound, has at most a fixed power of $t^{-1}$ as a pointwise bound. Applying H\"older with exponent $1+1/L$ and
\[
 \|v\|_q\le C\sqrt q\|v\|_{H^1}
\]
yields the asserted $H^{-1}$ estimate.

Since $\mathcal F_\rho(u)=0$, self-adjointness and \eqref{reg:weighted-square} give a comparison error $C\sqrt L(tf+t^2L+M)E+CLE^2$. Use $f=O(1)$, $tL=o(1)$ and Young's inequality on $M\sqrt L E$. Summing the resulting estimates for the center components gives
\[
 tf\le C(t^3L+t^2f+t\sqrt L E+LE^2+tM+M^2).
\]
After substitution and division by $t$, one obtains
\begin{align*}
 f\leq {}&Ct^{5/3}L^3
   +C(tL+t^{2/3}L^{3/2})f\\
 &+C(t+t^{1/3}L^3+t^{4/3}L^3)f^2
   +Ct^{5/3}L^6f^4.
\end{align*}
The terms independent of $f$ that are smaller than $t^{5/3}L^3$ have been included in the first term. All remaining coefficients tend to zero; since
\[
 f=o(1),
\]
the linear, quadratic and quartic terms can be absorbed. This proves the asserted center and energy rates and then $M=O(t^2L)$. This gives
\[
 |P-P^0|=o(\sqrt t),\qquad E=o(t^{7/6}L).
\]
The precise parameter bound placing the solution in the parameter range of the construction will follow from the mass law.

For that calculation only the exact equation and the estimates just proved are needed. We prove the stated remainder in the mass law by differentiating the normalized functional.

Put $\lambda=\log t$ and denote differentiation in $\lambda$, with $P$ and $\rho$ fixed, by a dot. Retain the notation $\widehat u^{\rm app}$, $N_i$ and $M_i$ from the construction, and set
\[
 \rho^0=(4\pi m,\rho_2),\qquad
 \mathscr R=-\Delta_g\widehat u^{\rm app}-K(N-\rho^0).
\]
Since $\dot c_1=2$ and $\dot c_2=0$, differentiation of
\[
 J_\rho(u^{\rm app})
 =\tfrac12\|\widehat u^{\rm app}\|_K^2
   -\sum_i\rho_i\log M_i+\sum_i\rho_i c_i
   +\sum_i\rho_i\int_M\widehat u_i^{\rm app}
\]
and integration by parts on $M$ give the exact identity
\begin{align}
 \partial_\lambda J_\rho(u^{\rm app})
 ={}&2\mu-2(M_1-4\pi m)
       +\int_M K^{-1}\mathscr R\cdot\dot{\widehat u}^{\rm app}
       \notag\\
    &+\sum_i\frac{M_i-\rho_i}{M_i}\dot M_i
       +\mu\int_M\dot{\widehat u}^{\rm app}_1.
 \label{reg:ordinary-trial-variation}
\end{align}
Thus both normalizing integrals are retained before any local expansion. The moment estimates in the construction, and their derivatives, imply
\[
 \begin{split}
 |M_1-4\pi m|+|\dot M_1|&\le Ct^2L,\qquad
 |M_2-\rho_2|+|\dot M_2|\le Ct^2,\\
 \int_M|\dot{\widehat u}^{\rm app}|&\le Ct^2L.
 \end{split}
\]
Indeed the linear term in the expansion of the first density has zero integral, its quadratic moment is $O(s_i^2\log(1/s_i))$, and the excess in the second density is $s_i^2\gamma_i(x)$. Also
\[
 \dot D_i=-4s_i^2/(s_i^2+r^2),\qquad |\dot m_{s_i}|\le2m_{s_i}.
\]
These formulas give the same bounds after differentiation; the fixed cut-off annuli contribute $O(t^2)$. The last line of \eqref{reg:ordinary-trial-variation} is consequently $O(Mt^2L+t^4L^2)$.

\emph{Step 3. The scale equation.} On an ordinary bubble disk let
\[
 Z_i=\frac{s_i^2-r^2}{s_i^2+r^2},\qquad
 J_i(x)=\psi_i^A(x)-\psi_i^A(p_i).
\]
There $N_1^{\rm loc}=m_{s_i}e^{J_i}$, $N_2^{\rm loc}=\gamma_i(x)(s_i^2+r^2)$, and
\[
 \dot{\widehat u}^{\rm app}=(-2(Z_i+1),Z_i+1).
\]
By the local calculation in Lemma~\ref{reg:residual}, the Euclidean density of $\mathscr R$ is $-K(m_{s_i}(e^{J_i}-1),s_i^2\gamma_i)^t$. Its first term in the residual pairing cancels the corresponding constant part of $-2(M_1-4\pi m)$. The remaining contribution from the first density is therefore
\[
 2\int m_{s_i}Z_i(e^{J_i}-1).
\]
The second-density contribution is $O(s_i^2)$. The first mass outside the bubble disks and the omitted model tails are $O(t^2)$, while on each fixed cut-off annulus the residual and the scale derivative are both $O(t^2)$.

Thus these regions contribute $O(t^2)$ to \eqref{reg:ordinary-trial-variation}; no boundary term is discarded.

The smooth backgrounds have uniformly bounded third derivatives. Since $\nabla J_i(0)=F_i(P)$ and $\Delta J_i(0)=\kappa_i$, parity, angular integration and the identities
\[
 \int_{D_d}m_sZ_s r^2=-8\pi s^2\log(1/s)+O_d(s^2),
 \qquad
 \int_{D_d}m_s r^3=O_d(s^2)
\]
give, for any fixed bubble-disk radius $d$,
\[
 2\int m_{s_i}Z_i(e^{J_i}-1)
 =-4\pi\kappa_i s_i^2\log(1/s_i)
       +O(t^2+t^2L|F_i(P)|^2).
\]
Here the quadratic term $\tfrac12(F_i\cdot x)^2$ is included in the last error, and the Taylor remainder has integral $O(t^2)$. This proves \eqref{reg:ordinary-trial-expansion}.

It remains to compare this derivative with that at the exact solution $u=u^{\rm app}+\eta$. Direct differentiation of the residual estimates gives
\[
 \|D\mathcal F_\rho(u^{\rm app})\dot u^{\rm app}\|_K
 =\|\partial_\lambda\mathcal F_\rho(u^{\rm app})\|_K
 \le C(t^{2/3}f+t^{5/3}L+M\sqrt L).
\]
To obtain this bound, differentiate the explicit profiles before taking their $L^{6/5}$ norms: the linear and quadratic Taylor terms still have norms $O(t^{2/3}f)$ and $O(t^{5/3})$, respectively.

Differentiating the normalizing integrals contributes $O(t^{5/3}L)$, and the parameter term is bounded as in \eqref{reg:normalized-residual}. The cut-offs are fixed. Self-adjointness,
\[
 \mathcal F_\rho(u)=0,
\]
the uniform $L^\infty$ bound on $\dot u^{\rm app}$, and \eqref{reg:weighted-square} now imply
\[
 |\partial_\lambda J_\rho(u^{\rm app})|
 \le C\{(t^{2/3}f+t^{5/3}L+M\sqrt L)E+LE^2\}.
\]
No derivative of $\eta$ is required in this estimate. Substituting $f=O(t^{5/3}L^3)$, $E=O(t^{5/3}L^2)$ and $M=O(t^2L)$ shows that this bound is $o(t^2)$. All terms other than $O(t^2)$ in the remainder of \eqref{reg:ordinary-trial-expansion} are also $o(t^2)$. Consequently
\[
 \mu=2\pi\sum_i\kappa_i s_i^2\log(1/s_i)+O(t^2).
\]
Finally $s_i=td_i(P)$, with $d_i$ positive and smooth on the fixed center neighborhood, gives
\[
 2\pi\sum_i\kappa_i s_i^2\log(1/s_i)
 =\frac\pi2\ell(P)t^2L+O(t^2)
 =L_0t^2L+O(t^2),
\]
because $\ell$ is smooth and $L|P-P^0|=o(1)$. The mass expansion gives $|\mu|/(t^2L)\to|L_0|$, so $|\mu|< C_1t^2L$ for every $C_1>4|L_0|$ and all sufficiently small $t$. The center ratio $|P-P^0|/\sqrt t$ tends to zero and
\[
 E/(K t^{7/6}L)=O(t^{1/2}L)
\]
tends to zero, so the modulated solution is strictly inside the shrinking center and correction bounds as well. Uniqueness identifies the correction.

Subtracting the original and projected equations gives zero source sum, hence all multipliers vanish by the invertible source test matrix.

Conversely zero multipliers restore the original normalized equation. \emph{Step 4. The local degree.} On the allowed side let $t_\mu$ solve $|\mu|=|L_0|t_\mu^2\log(1/t_\mu)$, write
\[
 P=P^0+t_\mu^{1/2}\xi,\qquad t=t_\mu\tau.
\]
For $1/2\leq\tau\leq 2$ the ratio $|\mu|/(t^2L)$ tends uniformly to $|L_0|/\tau^2$. Choose the construction constant
\[
 C_1>4|L_0|
\]
and the radius of the $\xi$ ball less than $C_0/\sqrt2$ so that the closed product lies strictly inside both construction bounds. Since $t^2\log(1/t)$ is increasing near zero and its ratio at $(1\pm\epsilon)t_\mu$ tends to $(1\pm\epsilon)^2$, the mass law gives
\[
 t/t_\mu\to1,\qquad \xi\to0.
\]
Divide the two real components of each center equation in \eqref{reg:ordinary-rows} by $-3s_i\sqrt{t_\mu}$ and the scale equation by $3|\mu|/(8\pi m)$. The center divisors have positive two-dimensional determinant, and the scale divisor is positive.

The coordinate change from $\log t$ to $\tau=t/t_\mu$ is increasing. The resulting map converges uniformly to $(DF(P^0)\xi,\sgn(\mu)(1-\tau^2))$ on a product of a small fixed center ball and
\[
 1/2<\tau<2.
\]
Its degree is $-\sgn(\mu)\sgn\det DF$. The determinant obtained by pairing the ordinary sources with the parameter derivatives is positive. The index of the projected operator is $m_S+m-1$.

The homotopy of the projected graph used at mixed points applies here as well: its uniform contraction controls the entire complementary ball over the closed scale interval. The product degree is $(-1)^m\operatorname{sgn}(\mu)d_S$.

For $L_0>0$ it is on the positive side. For $L_0<0$ it is $-(-1)^m d_S$ on the negative side, and is subtracted in the difference of the two local degrees. In both cases the mass expansion excludes concentration for $\mu L_0\le0$.
\end{proof}
\begin{lemma}[Uniqueness of the parameters in scaled coordinates]
\label{journal:modulation-uniqueness}
Let $\dot C^0(M)=\{v\in C^0(M):\int_Mv=0\}$. Take $X$ to be one of the spaces
\[
 (\dot C^0(M))^2,\qquad (\dot H^1(M))^2,\qquad
 (\dot C^0(M)\cap\dot H^1(M))^2,
\]
with the supremum norm, the Cartan energy norm or the sum norm, respectively. On balls of a common positive radius in the scaled parameters $\theta$, let $B_\theta$ be the global approximation and
\[
 L_\theta\in\mathcal L(X,\mathbb R^d)
\]
its orthogonality functionals. Assume $B_\theta,L_\theta$ are $C^1$ and
\begin{equation}\label{journal:modulation-hypotheses}
 \|L_\theta\|+\|D_\theta B_\theta\|+\|D_\theta L_\theta\|
       +\|J_\theta^{-1}\|\leq C_0,
 \qquad J_\theta=L_\theta D_\theta B_\theta.
\end{equation}
All norms are operator norms. The constant $C_0$ and the modulus of continuity of $J_\theta$ are independent of concentration.
\begin{itemize}
\item[\textup{(i)}]
If $\|u-B_{\theta_0}\|_X=o(1)$, then in a fixed sufficiently small ball about $\theta_0$ there is a unique representation
\[
 u=B_\theta+\eta,\qquad L_\theta\eta=0,
 \qquad |\theta-\theta_0|+\|\eta\|_X
       \leq C\|u-B_{\theta_0}\|_X.
\]
\item[\textup{(ii)}]
For the separated profiles assume, uniformly as concentration increases,
\[
 \|u-B_\theta\|_X=o(1).
\]
If the scaled parameter distance between any two representations of the same field is $o(1)$, both lie in one ball as in \textup{(i)}. Thus the parametrization is injective for sufficiently concentrated fields.
\end{itemize}
\end{lemma}
\begin{proof}
Put $\delta=\|u-B_{\theta_0}\|_X$ and $\mathfrak m_u(\theta)=L_\theta(u-B_\theta)$. By the operator bound,
\[
 \|\mathfrak m_u(\theta_0)\|\leq C_0\delta,\qquad
 \|J_{\theta_0}^{-1}\mathfrak m_u(\theta_0)\|
                                      \leq C_0^2\delta.
\]
Moreover
\[
 D\mathfrak m_u(\theta)
 =(D L_\theta)(u-B_\theta)-J_\theta.
\]
For $|\theta-\theta_0|\leq r$, the bound on $D_\theta B_\theta$ gives $\|u-B_\theta\|_X\leq\delta+C_0r$. If $\omega_J$ is the common modulus of continuity, then
\[
 \|I+J_{\theta_0}^{-1}D\mathfrak m_u(\theta)\|
 \leq C_0\{C_0(\delta+C_0r)+\omega_J(r)\}.
\]
Choose a common ambient radius $r_0>0$, and then $\delta_0>0$, so that this bound is at most $1/2$ for $r=r_0$ and $\delta\leq\delta_0$. Decrease $\delta_0$ so that $2C_0^2\delta_0<r_0$. The map
\[
 \mathcal N_u(\theta)=\theta+J_{\theta_0}^{-1}\mathfrak m_u(\theta)
\]
is therefore a contraction on the ambient ball, and maps the closed ball of radius $2C_0^2\delta$ about $\theta_0$ into itself. Its fixed point is a zero of $\mathfrak m_u$, satisfies $|\theta-\theta_0|\leq2C_0^2\delta$, and gives
\[
 \|\eta\|_X\leq\delta+C_0|\theta-\theta_0|\leq C\delta.
\]
The contraction estimate also gives uniqueness in the larger ball. For two representations in the stated class, use the first parameter as $\theta_0$.

Its remainder gives the small $\delta$ in the derivative bound; their $o(1)$ scaled parameter distance places the second zero in the same ambient ball. The local uniqueness just proved then identifies the parameters and remainders.
\end{proof}
\subsection{Isolating neighborhoods and the degree jump at regular blow-up points}
\label{reg:isolation}
For each shadow solution take the finite-dimensional isolating product above.

Use the projected uniqueness ball as the correction ball: $2C_A L\bar r$ in the mixed case and $K t^{7/6}L$ in the ordinary case.

The uniform correction is strictly inside these balls, by the contraction estimate with sufficiently large construction constants. We first compare the parameters of two representations of the same function.

For fixed nonzero $\mu$ the ratios of the two common scale parameters are bounded above and below by the chosen scale interval. If two representations agreed along a sequence with
\[
 \mu\to0,
\]
their trial fields would differ by $o(1)$ in $H^1$. On the ordinary widths and the inner mixed widths the respective first and second gradients converge after scaling to the gradients of the ordinary scalar profiles.

All cut-off tails and smooth exterior terms disappear at these scales; at a mixed point the other profile has width $R$ and contributes $O(a/R)$. A center discrepancy unbounded in its rescaled width would make the other scaled gradient vanish on compact sets.

Otherwise, equality of the two limiting gradients fixes the displacement and width ratio by the unique critical point and its nonzero second derivative. Thus all centers agree to $o(1)$ in their scaled units.

At a mixed point the first inverse-Cartan gradient, $(2\nabla u_1+\nabla u_2)/3$, has no contribution from the inner point mass.

On the outer scale its nonconstant part is the gradient of $-\log(1+|y^2-B|^2)$. Equality of these gradients forces equality of the polynomials after normalizing their leading term and then their quadratic coefficients.

Before making this comparison, the matched first or inner width fixes $t$ up to the bounded shape factor; an outer profile comparison can include a bounded ratio $\lambda$ of outer radii.

In that case the two positive quadratic-denominator polynomials are $1+|y^2-B|^2$ and $\lambda^4+|y^2-\lambda^2\widetilde B|^2$. They are monic at degree four, so comparison of their terms of degree two and zero gives
\[
 \lambda^2\widetilde B=B,\qquad \lambda^4=1.
\]
This also fixes the ratio of $t$ values if there is no ordinary point. Consequently all scaled parameter differences tend to zero. The ordinary rational tangents have bounded Dirichlet norm by rescaling with $s$.

The mixed tangents from $F_R$ and the leading factor $4R^4(a^2+|z|^2)$ of $Q_R$ have the same bound after rescaling with $R$ and $a$, respectively: their rescaled gradients are bounded on unit disks and are $O(|y|^{-2})$ outside. In $Q_R$ the remaining factor is $1+\beta(|z|^4+2R^2\Re(z^2\overline B))/(4(a^2+|z|^2))$.

Its differentiated logarithmic tangents have gradients $O(|z|+R^2/(a+|z|))$, uniformly in the bounded shape set, by the quotient rule. These gradients are bounded in $L^2$, as are the smooth and cut-off terms.

Once two parameters agree to $o(1)$ in scaled coordinates, use fixed reference widths along the short segment joining them. For a common represented field $u$, integrate the derivative of
\[
 \theta\mapsto L_\theta(u-u^{\rm app}(\theta))
\]
along this segment. By the $H^1$ tangent bounds,
\[
 u-u^{\rm app}(\theta)=o(1)
\]
there. The scaled source derivatives are bounded, while $L_\theta D_\theta u^{\rm app}$ converges uniformly on the segment to the invertible matrix obtained above by pairing the sources with the derivatives in the scaled coordinates.

The integrated matrix is invertible, so the parameters coincide. The same derivative is invertible near every interior representation. The implicit function theorem for
\[
 L_\theta(u-u^{\rm app}(\theta))=0
\]
recovers $\theta$ continuously from nearby $u$, and then
\[
 \eta=u-u^{\rm app}(\theta).
\]
Thus the image in the original function space is open, as in Lemma~\ref{journal:modulation-uniqueness}, and it is bounded at fixed nonzero $\mu$. For a converging sequence in this image, compactness of the closed parameter product gives a convergent parameter subsequence.

Subtracting its trial fields gives convergence of the corrections in $H^1$ and preserves the orthogonality. A boundary point is consequently on a parameter face or on the correction sphere.

Every zero in the closure lies on the graph of the unique projected correction and satisfies the finite-dimensional equations. The strict correction bound excludes the sphere, and the reduced boundary signs exclude every parameter face.

For the local degree, the normal derivative in source coordinates has matrix $G_\theta^{-1}(L_\theta D_\theta u^{\rm app} +O(\|\eta\|_{H^1}))$, where
\[
 G_\theta=((S_\alpha,S_\beta)_K)_{\alpha,\beta}>0.
\]
This follows by differentiating
\[
 L_\theta\eta=0.
\]
The determinant of $L_\theta D_\theta u^{\rm app}$ is positive on both parameter domains. Transport the complementary space along the same local parameter coordinates as in the projected construction.

At fixed $\mu$ the compact closed product and the strict correction bound allow excision to a smaller neighborhood of the graph of the unique correction.

Choose it small enough that $CL^{5/2}$ times its radius is less than one. The complementary homotopy then keeps all its zeros on the graph, and fixing the correction in the finite-dimensional equations preserves the zero-free faces.

With this orientation, the local degree in the original space is the product of the finite-dimensional degree and the sign on its complement. The estimates above place each blowing-up sequence in the interior of its corresponding neighborhood.

To distinguish neighborhoods of different shadow solutions, first choose their closed center neighborhoods so that the compact images $(P,Q,w(P))$, with centers of each type unordered, are disjoint for distinct solutions.

This is possible by finiteness and continuity of the scalar branches. A limiting parameter of a fixed mixed center product can be anywhere in that closed neighborhood. If two such neighborhoods met along
\[
 \mu\to0,
\]
their trial fields would differ by $o(1)$ in $L^2$, by the mean-zero Poincar\'e inequality and the correction bounds.

Each trial field converges smoothly on compact sets off the union of the possible blow-up sets. Thus their limits agree there in $L^2$ and agree as distributions on $M$, since both limits are locally integrable.

To recover the data, apply $-\Delta_gK^{-1}$ to each limit. If $\nu_1,\nu_2$ are the sums of its point masses and $w$ is its mean-zero scalar remainder, the result is
\[
  \binom{\nu_1-\nu_1(M)}
        {\nu_2-\nu_2(M)+n^g-\sigma},
  \qquad -\Delta_gw=2(n^g-\sigma).
\]
Here $n^g$ is the integrable shadow density. Its measure has no point masses. Thus the two components of this distribution determine each concentration point and its mass pair.

Subtracting those Green potentials identifies $w$. The limiting triples would then lie in two disjoint compact images.

We use one labeling in each neighborhood, since permutations of centers of the same type give the same configurations. Thus the finitely many neighborhoods are disjoint and their constants can be chosen uniformly.

To sum all contributions at the critical parameter, first take the problem without prescribed singularities. Take all configurations with
\[
 U+2V=m,\qquad \rho_2-4\pi V>0,
\]
and use common weights from Theorem~\ref{fnd:common-weights} and Corollary~\ref{fnd:extra-traces} for this family. The scalar exclusions hold because
\[
 \rho_2\notin4\pi\N,
\]
and the family contains every applicable pair replacement with positive remainder. The second component cannot concentrate completely at this parameter.

Thus every escaping sequence has a limiting shadow of one of these configurations. Both mass laws then exclude unbounded sequences of solutions at
\[
 \mu=0.
\]
Solutions in the corresponding neighborhoods leave every fixed ball as
\[
 \mu\to0.
\]
For the trial field, on a concentric annulus between the square root of a local width and a fixed disk radius, its first component is $-4\log r+O(1)$ at an ordinary point and $-6\log r+O(1)$ at a mixed point, up to its common mean constant.

The cut-off corrections have bounded size at both circles. Their circular means therefore differ by at least $cL-O(1)$.

Cauchy--Schwarz in logarithmic radius gives Dirichlet energy at least $c'L$. The $o(1)$ correction in $H^1$ preserves escape by the triangle inequality for the energy norm.

Choose a fixed open ball containing the entire zero set at the critical parameter in its interior. The normalized Green map is jointly compact in the parameter and the function on bounded $H^1$ sets near the critical value.

Hence a hypothetical sequence of zeros on its boundary would converge to a zero at $\mu=0$ on the boundary, and likewise a bounded sequence of zeros outside the ball would converge to a zero at $\mu=0$ outside it.

Both are impossible. An unbounded $H^1$ sequence of solutions is also unbounded in $C^0$, since a uniform $C^0$ bound gives a uniform $L^p$ density bound and then a uniform energy bound by Green.

For sufficiently close noncritical parameters, the estimates for blowing-up solutions therefore place every zero outside the ball in one of the disjoint neighborhoods. The degrees on the ball agree by homotopy invariance. Excision and additivity give the degree jump at
\[
 \rho_1=4\pi m.
\]
Let $k$ be the integer such that $4\pi k<\rho_2<4\pi(k+1)$. Since $S=\varnothing$, for $m\geq1$ this gives
\begin{equation}\label{reg:wall-jump}
d^{\mathrm{reg}}_{m,k}-d^{\mathrm{reg}}_{m-1,k}=\sum_{\substack{U+2V=m\\4\pi V<\rho_2}}
   (-1)^U d^{\rm sh}(U,V,\varnothing,\varnothing).
\end{equation}
\begin{remark}
\label{reg:hybrid-limitation}
\begin{itemize}
\item[\textup{(i)}]
If all blow-up points are regular and the configuration contains mixed points, then
\[
 L\|\eta\|^2=O(t^{4/3}L^3)=o(t),\qquad R_j^2\asymp t\quad(j\text{ mixed}).
\]
\item[\textup{(ii)}]
If the blow-up set meets $S$, the outer radii need not be comparable. The local identities at regular points, the expansion in the common scale, and uniform $C^0$ estimates for the reduced shape equations give the degree.
\end{itemize}
\end{remark}
\section{Regular bubbling in the presence of fixed zeros}
\label{regzero:section}
We now consider configurations for which blow-up occurs away from all prescribed singular points. The weights may vanish at these points.

We prove that the construction and local degree formula of Section~\ref{reg:section} remain valid in this setting. Fix a finite set $S$ of marked points and source strengths $\alpha_{i,p}\geq0$. Assume
\begin{equation}\label{regzero:weights}
 h_i(z)=|z|^{2\alpha_{i,p}}a_{i,p}(z)\quad\hbox{near }p\in S,
 \qquad a_{i,p}>0\hbox{ smooth},
\end{equation}
and that $h_i$ is smooth and positive on $M\setminus S$. Thus $h_i$ is bounded and nonnegative. The free centers vary in a compact set of configurations, uniformly separated from each other and from $S$. We use $U,V,n,m,\sigma,t,s_i,R_j,a_j$ as in Section~\ref{reg:section}; in particular $m=U+2V>0$ and
\[
 \sigma=\rho_2-4\pi V>0.
\]
Assume that the scalar linearization and the solution of the full shadow system are nondegenerate. At every mixed point assume $\kappa_j>4|T_j|$ uniformly. These assumptions concern a fixed neighborhood of the selected solution of the shadow system.
\subsection{Scalar branches and the diffuse part of the inverse}
\begin{lemma}[The scalar branch with fixed zeros]
\label{regzero:scalar}
For the weights \eqref{regzero:weights}, the branch of \eqref{reg:shadow-scalar} satisfies the following.
\begin{itemize}
\item[\textup{(i)}]
It is smooth in $P$ with values in $W^{2,p}(M)$ for $1<p<\infty$. The normalized inverse and the parameter derivatives are uniformly bounded on a smaller compact neighborhood of the centers.
\item[\textup{(ii)}]
The densities $n_*^g$ and their parameter derivatives are bounded on $M$ and uniformly smooth on the disks about the free blow-up points.
\item[\textup{(iii)}]
The densities have mass $\sigma$ and a common positive lower bound on a fixed open set. The scalar Morse index is constant.
\end{itemize}
\end{lemma}
\begin{proof}
Choose $0<\theta<1$ smaller than every positive $2\alpha_{2,p}$. Then
\[
 h_2\in C^{0,\theta}.
\]
Locally $e^{-4\pi G(x,p)}=|z(x)-z(p)|^2$ times a jointly smooth positive factor. Thus $P\mapsto H_P$ is smooth into $C^{0,\theta}$ and each $L^p$; no division by $H_P$ is involved.

On an open set avoiding $S$ and all center disks, $H_P$ has a uniform positive lower bound. The normalizing integral is consequently bounded away from zero near the fixed scalar solution.

Since $W^{2,p}\hookrightarrow C^0$ for $p>1$ in two dimensions, the normalized scalar equation is a smooth map from the mean-zero $W^{2,p}$ space to the mean-zero $L^p$ space.

Its derivative is Fredholm of index zero, since its difference from the mean-zero Laplacian factors through the compact embedding into $C^0$. A $W^{2,p}$ kernel function belongs to $H^1$ by the two-dimensional embedding.

Conversely a weak energy kernel function belongs to every $W^{2,p}$ because the scalar density is bounded, $H^1$ embeds into each finite $L^p$, and the equation then has $L^p$ right-hand side. Thus energy nondegeneracy makes the derivative invertible on
\[
 W^{2,p}_0\to L^p_0
\]
for each fixed $p>1$. Invertibility and the local inverse bounds persist on a smaller compact center neighborhood. The implicit function theorem gives the parameter derivatives.

On the free disks the coefficients are smooth, including the ordinary zero, so local elliptic estimates give the required higher bounds. Continuity and invertibility of the second derivative of the scalar energy keep its Morse index fixed.
\end{proof}
\begin{lemma}[Invertibility with vanishing residual density]
\label{regzero:inverse}
Assume the hypotheses of Proposition~\ref{reg:inverse}, including uniform invertibility of the centered scalar operator along the branch with fixed centers. Allow its density to satisfy
\[
 0\leq n_*^g\in L^\infty(M),\qquad
 \int_M n_*^g>0,\qquad n_*^g>0\ \text{on a fixed open set}.
\]
Assume also that $n_*^g$ is uniformly smooth on the disks about the free blow-up points. In particular, along the branch of Lemma~\ref{regzero:scalar},
\[
 \|L_t^{-1}\|_{H^{-1}\to H^1}\leq C\log(1/t),\qquad
 \ind L_t=m_S+U+4V-1.
\]
The bounds hold for fixed source strengths on a compact neighborhood of the branch.
\end{lemma}
\begin{proof}
The variance decomposition in \eqref{reg:form} remains valid for a nonnegative residual density of positive mass. For an $H^1$-bounded sequence converging strongly in $L^2$,
\[
 \int_M n_*^g|v_k-v|^2\longrightarrow0.
\]
Thus the contribution of $n_*^g$ to the centered operator is compact. The spectral gap and the scalar negative index follow as before. The weighted mean is a bounded functional. Positivity on the fixed open set controls the constants in the decomposition into local means.

The Green operator maps bounded densities to bounded potentials. The residual density therefore contributes only $O(1)$ to the inner products of the vectors representing the local means. This preserves the logarithmic term in \eqref{reg:mean-gram}. More precisely, for each finite $p$,
\[
 \|G*(n_*^g v)\|_{W^{2,p}}\leq C_p\|v\|_{H^1}.
\]
Choose $p>2$. Pairing this $C^{1,\theta}$ potential with a centered source, subtract its value at the source center. The remaining absolute integral is bounded by the first moment of the absolute source density.

The bounds for the profiles on the plane give a positive power of their width, and the part from the term from the residual density in an inner translation source has the same property by its $O(a)$ local $L^p$ bound.

Thus the centered radial and projection errors decay as a positive power of the width and hence are $o(1/L)$. The diffuse mean Green vectors remain bounded in energy and add only $O(1)$ to these inner products.

All rescaled kernels, source moments and local mean contrasts are unchanged. The proof of Proposition~\ref{reg:inverse} still gives
\[
 V^*A^{-1}V=G+O(1)
\]
in its mean block, the same finite-rank inverse formula, and the positive relative radial quadratic form after elimination of size $1/\log(1/t)$. These estimates give the inverse bound and the index, including the signs of the small radial eigenvalues.
\end{proof}
\subsection{Local profiles and the simultaneous regular construction}
\begin{proposition}[Local estimates and projected solutions]
\label{regzero:construction}
At a free point with isolated local masses $(4\pi,0)$ or $(8\pi,4\pi)$, the conclusions of Section~\ref{reg:section} extend to \eqref{regzero:weights} as follows.
\begin{itemize}
\item[\textup{(i)}]
The selection of bubbles, the identity for the shape parameter and the estimates for the profiles and their Green potentials hold.
\item[\textup{(ii)}]
At a mixed point, $\kappa_j>4|T_j|$ excludes splitting and bounds the shape parameter.
\item[\textup{(iii)}]
On the same parameter domains, the simultaneous projected solution exists, is unique and has the same orientation of the matrix of pairings with the parameter derivatives.
\end{itemize}
\end{proposition}
\begin{proof}
Choose the disks about the blow-up points and their cut-off annuli disjoint from $S$. On a slightly larger disk the canonical forcing in \eqref{reg:canonical} has uniform smooth bounds. Subtract its Dirichlet potential.

The resulting local Toda system has positive coefficients in a compact $C^3$ class, bounded local masses and bounded boundary oscillation. These are the hypotheses of the local selection and grouping results used in Lemmas~\ref{reg:stress} and~\ref{reg:local-power}. In the mixed case splitting would imply
\[
 |\kappa_j|=4|T_j|,
\]
contrary to
\[
 \kappa_j>4|T_j|.
\]
Equation \eqref{reg:shape-limit} bounds $B$ because $\nabla\mathcal M$ maps onto the disk of radius $\pi$. The inner and outer profiles exhaust the two local masses, excluding any further bubble.

These local data and the smooth exterior field give the estimates on the annuli between scales, the bound for the scalar inverse with weight $\log(2+|y|)$, and the estimates for the profiles and Green potentials.

The original fixed sources lie outside the disk; their Green fields have bounded derivatives there. The constants are uniform when the other blow-up points remain separated from this disk.

For a configuration containing only regular bubbles, use the same global approximations and the common choice $c_1=2\log t$ in \eqref{reg:scales}. Outside these disks their equation is the scalar equation of Lemma~\ref{regzero:scalar}, and their densities satisfy
\[
 N_1=O(t^2),\qquad N_2=n_*^g.
\]
The bound for $N_1$ uses only
\[
 h_1\in L^\infty(M),
\]
also at its zeros. On ordinary disks the second density is a smooth positive factor times $s_i^2+|z|^2$. On mixed disks the exact diffuse cancellation is still \eqref{reg:cancelled-residual}. Every cut-off derivative is supported on a smooth annulus.

We thus obtain the same $L^{6/5}$ residual and density comparisons and the same normalized mass errors as Lemma~\ref{reg:residual}; the ordinary shrinking-center construction has its original sharper residual as well.

In these statements the smooth correction is defined at $\rho_1^0$, with forcing $f^0$, whereas the forcing at $\rho_1$ satisfies \eqref{reg:wall-forcing}. We include both $(2\mu,-\mu)e^\varphi$ and the error in the normalization in \eqref{reg:normalized-residual}, also on the larger interval of mass parameters.

Lemma~\ref{regzero:inverse} and the density comparison give an error in operator norm of $O(t^{2/3}\log^{3/2}(1/t))=o(1/\log(1/t))$ in the mixed case. The contraction uses bounded normalized masses, the bound
\[
 \|V_i^{\rm app}\|_\infty\leq Ct^{-2}
\]
and exponential Sobolev estimates, all of which use $h_i$ multiplicatively. Thus \eqref{reg:range} and the ordinary construction follow with the same constants and index. Use the locally centered source densities of the approximation and the projected sign convention of \eqref{reg:projected-equation}.

Estimate \eqref{reg:modulation-bounds} applies unchanged: the sources and their scaled derivatives have bounded total variation and bounded Green energy. The entries of the matrix obtained by applying the $L_\alpha$ to the parameter derivatives depend only on the local terms.

The $L_\alpha$ annihilate constants in the smooth correction, and Lemma~\ref{regzero:scalar} controls its remaining derivatives. Its limiting determinant is positive.

The same orthogonality equations therefore determine the local parameters uniquely. Permuting centers of the same type permutes even-dimensional blocks, hence preserves orientation.
\end{proof}
\subsection{Description of all regular blowing-up solutions}
\begin{theorem}[Local degree in the presence of fixed zeros]
\label{regzero:transfer}
Fix the source strengths and assume the above hypotheses with $m=U+2V>0$. If $V=0$, assume $L_0\ne0$, where $L_0$ is the coefficient for $U\geq1$ ordinary points.
\begin{itemize}
\item[\textup{(i)}]
Every sequence concentrating only at regular points with the given limiting shadow solution belongs to the family of Proposition~\ref{regzero:construction}.
\item[\textup{(ii)}]
The bubbling solutions have bounded open neighborhoods with zero-free boundary. The local degree for $\mu>0$ minus that for $\mu<0$ equals
\[
                         (-1)^U d_S.
\]
\end{itemize}
The proof also gives the estimates for the mass parameter.
\end{theorem}
\begin{proof}
\smallskip\noindent\textit{Step 1. Comparison with the scalar solution.} We use the cut-off definitions in \eqref{reg:actual-diffuse}: $\mathfrak n_i$ are the normalized densities of the solution,
\[
 \nu_1=\sum_a\chi_a\mathfrak n_1
\]
on the fixed disks about the first-component blow-up points, $\nu_2=b\mathfrak n_2$ with inner cut-off radii $\tau_j=R_j^{3/2}$, and
\[
 e_1=\mathfrak n_1-\nu_1,\quad d_D=\mathfrak n_2-\nu_2,\quad
 Z_2=-G\nu_1+2G\nu_2,\quad w_u=u_2-Z_2.
\]
Thus $-\Delta_gw_u=2\Pi d_D-\Pi e_1$ exactly. We apply the scalar inverse to this field $w_u$, comparing $Z_2$ with the corresponding Green potentials of the fitted profiles. Put
\[
 D=\max_jR_j\vee\max_i\sqrt{s_i}
\]
before comparing any radii. At each individual radius $\varrho_a=s_i$ or $R_j$, the local profile and mass estimates give an $L^1$ error between the bubble density and its profile of at most $C\varrho_a^{1/4}(1+|\log \varrho_a|)^N$ and a polynomial density bound $C\varrho_a^{-K}$.

Choose a fixed $1<p_0<6/5$ close enough to one to retain a positive power after interpolation for every bubble. Lemma~\ref{journal:independent-radii} then bounds their Green errors by $CD^c$, including each mass error times its own logarithmic width.

The exponent $p_0$ is used for the comparison of the solution with the approximate family; the construction itself uses
\[
 p=6/5.
\]
At a point on the boundary of a fixed disk about a regular blow-up point, the bound on the tail of the first density and the lower bounds for $h_1$ and $e^{u_1}$ give
\[
 e^{c_1}=O(D^c).
\]
On the remaining fixed complement, including the unused sources, the small-mass Green estimate gives bounded $u_1$; since $h_1$ is bounded, this gives
\[
 \|e_1\|_{p_0}\leq CD^c.
\]
Removing the inner parts of the solution density leaves a tail with $p_0$th-power bound
\[
 Ca_j^{p_0(2-\varepsilon)}
       \tau_j^{2-p_0(4-\varepsilon)}
       =R_j^{3-2p_0+o(1)+O(\varepsilon)},
\]
which has a positive exponent. The bound $Cr^{-\varepsilon}$ farther from the inner bubble is uniformly $L^{p_0}$ integrable; ordinary points have the stronger quadratic-zero bound.

On a fixed disk about an unused source there is no concentrating mass. 
In particular,
\[
 \lim_{r\downarrow0}\limsup_{D\to0}
 \int_{B_r(p)}(\mathfrak n_1+\mathfrak n_2)\,dV_g=0
 \qquad(p\in S).
\]
Fix $q>p_0$. Choose $B=B_r(p)$ disjoint from the bubble disks and small
enough that $\int_B\mathfrak n_i\,dV_g\le\pi/q$ for all sufficiently
small $D$. For every fixed $\alpha_{i,p}\ge0$, \eqref{regzero:weights} gives
\[
 c_i:=\log\frac{\rho_i}{\int_M h_i e^{u_i}\,dV_g},\qquad
 \mathfrak n_i=h_i e^{u_i+c_i}
 =|z|^{2\alpha_{i,p}}a_{i,p}(z)e^{u_i+c_i}
 \le C e^{u_i+c_i}\quad\text{on }B.
\]
Let $G_B$ be the Dirichlet Green operator for $-\Delta_g$ on $B$ and
$H_i$ the harmonic extension of $u_i+c_i$ from $\partial B$. The
boundary bounds for $u_i$ and Jensen's bound $c_i\le C$ give
\[
 \sup_B H_i\le C,\qquad
 u_i+c_i
 =H_i+G_B\!\left(\sum_{j=1}^2 K_{ij}(\mathfrak n_j-\rho_j)\right)
 \le 2G_B\mathfrak n_i+C.
\]
This upper bound allows $u_1+c_1\to-\infty$. On a smaller concentric
disk $B'\Subset B$, the logarithmic-kernel estimate in Step~1 of the
proof of Proposition~\ref{fnd:signed-compactness} gives
\[
 \int_{B'}e^{2qG_B\mathfrak n_i}\,dV_g\le C_q,\qquad
 \|\mathfrak n_i\|_{L^q(B')}\le C_q.
\]
Elliptic regularity now bounds the local Green remainders on smaller
disks. Since $q>p_0$, this estimate and the preceding bubble-tail bounds
give uniform integrability of $d_D^{p_0}$. Qualitative convergence away
from the free centers therefore yields
\[
 \|d_D-n_*^g\|_{L^{p_0}(M)}\longrightarrow0.
\]

Subtract \eqref{reg:shadow-scalar} from \eqref{reg:actual-diffuse}.
Both fields have mean zero, so the elliptic estimate and the preceding
bound for $e_1$ give
\[
 \begin{aligned}
 -\Delta_g(w_u-w(P))
   &=2\Pi(d_D-n_*^g)-\Pi e_1,\\
 \|w_u-w(P)\|_{W^{2,p_0}(M)}
   &\le C\bigl(\|d_D-n_*^g\|_{L^{p_0}(M)}
                +\|e_1\|_{L^{p_0}(M)}\bigr)\longrightarrow0.
 \end{aligned}
\]
We may now use the inverse of the normalized scalar linearization from
Lemma~\ref{regzero:scalar}\textup{(i)} in this neighborhood of $w(P)$.
For the coefficient comparison, $u_2=Z_2+w_u$ gives
\[
 A_D=(1-b)h_2e^{Z_2},\qquad
 C_u=\frac{\rho_2}{\int_M h_2e^{u_2}\,dV_g},\qquad
 d_D=(1-b)\mathfrak n_2=C_uA_De^{w_u}.
\]
The coefficient comparison in the proof of \eqref{reg:actual-diffuse} uses only bounded exterior weights. At an ordinary point the model Green potential replaces $r^2$ by $r^2+s_i^2$: both logarithms have the corresponding mass $4\pi$ and
\[
 \Delta\log(r^2+s_i^2)=m_{s_i}.
\]
Define $H_{P,s}$ by replacing $4\pi\delta_{p_i}$ in the definition of $H_P$ by the fitted ordinary profile with its fixed radial cut-off.

In the regular Green remainder the first moments vanish by radial symmetry, the second moment of $m_{s_i}$ on the fixed disk is $O(s_i^2|\log s_i|)$ and the cut-off tail has mass $O(s_i^2)$.

The bounded exterior harmonic terms therefore change by $O(\sum_i s_i^2|\log s_i|)$ in the logarithm of the fitted coefficient. Multiplying $1+s_i^2/r^2$ by
\[
 H_P=r^2k_i
\]
first then gives
\[
 \|H_{P,s}-H_P\|_\infty\leq C\sum_i s_i^2|\log s_i|.
\]
Extending a fitted inner cut-off to the fixed-disk model costs the Green error $O((a_j/\tau_j)^2(1+|\log a_j|+|\log\tau_j|))$
\[
 =R_j^{1+o(1)}(1+|\log R_j|),
\]
as proved after \eqref{reg:actual-diffuse}. At a mixed point the exponential ratio differs from one by at most $C(R_j^2/r^2+R_j^4/r^4)$ for $r\geq \tau_j$. The two polar integrals, starting at
\[
 \tau_j=R_j^{3/2},
\]
give the powers $R_j^{3/p_0-1}$ and $R_j^{3/p_0-2}$. Both are positive at the chosen $p_0$. The deleted disk and cut-off transition have vanishing norm. The bounded factors $h_2$ and $H_P$ are included before taking norms, so the estimate remains valid at their zeros. Consequently
\[
 \|A_D-H_P\|_{p_0}
 +\left|\int_M\nu_2-4\pi V\right|\leq CD^c.
\]
The positive integral on a fixed open set and the bounded $w_u$ give the upper bound on $C_u$. For the lower bound use the uniform $L^1$ bound on $A_D$ and $\int_Md_D\to\sigma>0$ in the exact identity for $d_D$. Hence $e_2=C_u(A_D-H_P)e^{w_u}$ satisfies
\[
 d_D=C_uH_Pe^{w_u}+e_2,\qquad
 \|e_1\|_{p_0}+\|e_2\|_{p_0}
 +\left|\int_M\nu_2-4\pi V\right|\leq CD^c.
\]
With $\delta_2=\int_M\nu_2-4\pi V$ and $Q=H_Pe^{w_u}/\int_MH_Pe^{w_u}$, the mass identity gives $C_u\int_MH_Pe^{w_u}=\sigma-\delta_2-\int_Me_2$. Thus the scalar defect is exactly
\[
 -\Delta_gw_u-2\sigma\Pi Q
     =\Pi\bigl(2e_2-e_1-2(\delta_2+\textstyle\int_Me_2)Q\bigr).
\]
It is $O(D^c)$ in $L^{p_0}$. The branch of Lemma~\ref{regzero:scalar} has a uniform inverse in the neighborhood obtained by strong convergence above. Lemma~\ref{journal:scalar-neighborhood} and \eqref{journal:scalar-remainder-equation} therefore give $\|w_u-w(P)\|_{W^{2,p_0}}\leq CD^c$.

\smallskip\noindent\textit{Step 2. Determination of the scales.} We apply the matching argument preceding \eqref{reg:capture-seed}. We use the two Green representations and the estimate
\[
 \|w_u-w(P)\|_{W^{2,p_0}}\le CD^c.
\]
All cut-off potentials and center evaluations are on disks disjoint from $S$; after replacing the actual densities by the fitted ones, the same local logarithms and smooth remainder estimates apply.

Evaluation at an ordinary center, at the inner mixed center and at a nonzero point on the outer scale gives, respectively,
\[
\begin{aligned}
4s_i^2&=e^{c_1+\psi_i^A(p_i)}(1+O(D^c)),\\
4a_j^2&=\beta_jR_j^4(1+|B_j|^2)(1+O(D^c)),\\
16R_j^4&=e^{c_1+\psi_j(q_j)}(1+O(D^c)).
\end{aligned}
\]
Thus $t^2=e^{c_1}$ and $D\asymp\sqrt t$, so $O(D^c)=O(t^{c/2})$. The total first mass comparison gives
\[
 |\mu|=O(t^c).
\]
The difference between fitted Green corrections and the trial corrections is estimated by the same logarithm and gradient bounds as there, giving
\[
 \|u-u^{\rm app}\|_{C^0\cap H^1}=O(t^{c'}).
\]
This comparison only uses the scalar $W^{2,p_0}$ bound and the globally bounded weights off the free disks. When $V=0$, the two mixed evaluations are omitted.

\smallskip\noindent\textit{Step 3. Identification of the projected correction.} The source measures have bounded total variation, so this preliminary $C^0$ estimate allows us to choose the parameters by orthogonality. The normalized nonlinear remainder satisfies
\[
 \|\mathcal N(\eta)\|_{H^{-1}}
       \leq C\|\eta\|_\infty L\|\eta\|_{H^1}.
\]
The inverse remains valid on the preliminary parameter range $|\mu|=O(t^c)$: its form changes by $O(|\mu|L)$, whose product with the inverse is $o(1)$. Absorb the displayed nonlinear term. For $V\geq1$, writing $E=\|\eta\|_{H^1}$ and $M=|\mu|$, the projected equation gives
\[
 E\leq C(t^{2/3}L+ML^{3/2}).
\]
Apply the normalized energy test from the proof of Theorem~\ref{reg:mixed-capture} with the backgrounds of the present construction. Its disk estimates use only smoothness on the separated regular disks.

On their fixed complement the first density and its scale derivative are $O(t^2)$ in $L^p$, since the physical weight is bounded there, including at its fixed zeros.

The second term from the residual density there is the actual scalar branch density. Thus the same test, with both normalizing integrals retained, gives
\[
 M\leq C\{t+MtL+(t^{2/3}+M\sqrt L)E+LE^2\}.
\]
Inserting the projected estimate and using the preliminary power bound for $M$ absorbs the linear and quadratic terms in $M$ and gives
\[
 M=O(t),\qquad E=O(t^{2/3}L).
\]
Enlarge $C_0$ to contain the resulting $M/t$ bound, as in the smooth calculation. With
\[
 \bar r=C_r t^{2/3},
\]
the projected inequality gives $E\le C_A L\bar r+E/4$ and hence
\[
 E\le (4/3)C_A L\bar r<2C_A L\bar r.
\]
This identifies the unique projected correction. Subtracting the original equation shows that all multipliers vanish, since the matrix of the $L_\alpha$ applied to the parameter derivatives is invertible.

Conversely zero multipliers give the original normalized equation. For $V=0$, put
\[
 f=|F(P)|=o(1).
\]
The same scale and translation tests on the surface for ordinary bubbles as in the proof of Theorem~\ref{reg:ordinary-capture} give the second and third inequalities below; the first follows from the compressed inverse.

In these tests the local background is smooth on all regular disks, and the $O(t^2)$ exterior first density and the diffuse second density have the same bounds as in the preceding mixed calculation.
\begin{align*}
 E&\leq C(t^{2/3}Lf+t^{5/3}L^2+ML^{3/2}),\\
 M&\leq C(t^2L+LE^2+M^2),\\
 tf&\leq C(t^3L+t^2f+t\sqrt L E+LE^2+tM+M^2).
\end{align*}
Substituting the first estimate into the second gives $M\leq C(t^2L+t^{4/3}L^3 f^2+M^2L^4)$. The preliminary power estimate absorbs the last term. Substitution back gives
\[
 E\leq C(t^{2/3}Lf+t^{5/3}L^2),
\]
using $f=o(1)$. After division by $t$, the third estimate has constant term $O(t^{5/3}L^3)$ and the remaining terms are bounded by
\[
 C(t+t^{2/3}L^{3/2})f
 +C(t^{1/3}+t^{4/3})L^3f^2+Ct^{5/3}L^6f^4.
\]
Since $f=o(1)$ these terms can all be absorbed into the left-hand side. This gives
\[
 f+|P-P^0|=O(t^{5/3}L^3),\quad
 E=O(t^{5/3}L^2),\quad M=O(t^2L).
\]
The center and error bounds satisfy
\[
 |P-P^0|=o(\sqrt t),\qquad E=o(t^{7/6}L).
\]
To locate the parameters, apply the equation for the solution, the variation of the normalized energy and \eqref{reg:ordinary-trial-variation}--\eqref{reg:ordinary-trial-expansion}. This gives $\mu=L_0t^2L+O(t^2)$, as for smooth positive weights. Thus
\[
 |\mu|<C_1t^2L,\qquad |P-P^0|<C_0\sqrt t,
\]
and $E<Kt^{7/6}L$ for small $t$, since the last two ratios tend to zero and $C_1>4|L_0|$. Uniqueness identifies the ordinary projected correction, and the same comparison of the projected and original equations gives zero multipliers in both directions.

\smallskip\noindent\textit{Step 4. Local degree.} We compare the reduced equations on the domains of the projected solutions. The translation and scale tests take place on smooth free disks and fixed smooth cut-off annuli. On the complement of the free disks and cut-offs, the exterior fields satisfy
\[
 -\Delta_g\widehat u^{\rm app}=K(-4\pi m,n_*^g-\rho_2)^t,\qquad N_2=n_*^g.
\]
At $\mu=0$ the residual there is $-K(N_1,0)^t$, also after parameter differentiation. The first density and its scale derivative are $O(t^2)$, the scale tangent of $\widehat u^{\rm app}$ vanishes there, and the scaled translation tangents are $O(t)$.

Thus the local moment identities and all omitted-annulus errors have the same orders as in the smooth calculations.

The estimates for the normalized nonlinear terms use the same upper bounds for the probability densities and the same radii for the corrections.

They give $o(t)$ in the mixed equations, and $o(t^{3/2})$ and $o(t^2L)$ in the equations for the ordinary centers and the scale, respectively.

Testing the projected equation retains $\sum_\alpha c_\alpha L_\alpha Y$. The local matrices of the pairings $L_\alpha Y$ have the same invertible limits.

Consequently the mixed and ordinary reduced maps used below apply to these projected families on their whole stated parameter domains. For $V\geq1$, the solutions in part~\textup{(i)} satisfy
\[
 |\mu|\leq C_0t,
\]
and their radial equation gives $\mu/t\to S_*>0$. Set $t=\mu\tau$ on the positive side. Choose the two endpoints $0<\tau_-<\tau_+$ as in Proposition~\ref{reg:mixed-map}, and take $C_0>1/\tau_-$ so that the entire closed interval lies in
\[
 |\mu|<C_0t.
\]
The fixed isolating center neighborhood and shape ball have the strict boundary signs of the same proposition, by the preceding estimates for the reduced equations.

The shape degree is $+1$, the scale degree is $-1$, and the index of the projected operator is $m_S+U+4V-1$. Their product is $(-1)^Ud_S$, and no blowing-up sequence occurs for $\mu\leq0$. For $V=0$, use
\[
 |P-P^0|\leq C\sqrt t,\qquad |\mu|\leq C_1t^2L.
\]
On the allowed side put $P=P^0+\sqrt{t_\mu}\,\xi$ and $t=t_\mu\tau$. For $1/2\leq\tau\leq2$, the ratio $|\mu|/(t^2L)$ tends uniformly to $|L_0|/\tau^2$. Take $C_1>4|L_0|$ and the $\xi$-ball radius less than $C_0/\sqrt2$ so that the closed product lies in the construction domain
\[
 |P-P^0|<C_0\sqrt t,
\]
\[
 |\mu|<C_1t^2L.
\]
On this product the map reduces to $(DF(P^0)\xi,\sgn(\mu)(1-\tau^2))$. The degree on this side is $(-1)^U\sgn(\mu)d_S$.

For either sign of $L_0$, subtracting the degree for $\mu<0$ from that for $\mu>0$ gives $(-1)^Ud_S$. The preceding bounds put every blowing-up sequence strictly inside these domains.

Compare the approximations on the same disjoint disks in the energy space. The first ordinary and second inner mixed gradients determine the centers and inner widths. At a mixed point the first gradient after applying $K^{-1}$ determines the outer width and shape.

The estimate along the segment of parameters gives injectivity, and the local inverse gives openness. The correction balls remain within the projected uniqueness balls. The strict bounds for the correction exclude zeros on its sphere.

The preceding estimates for the reduced equations exclude zeros on the boundary in parameter space. The determinant of the matrix of these pairings and the inverse on the complementary space have the signs used above.

At fixed nonzero $\mu$, the parameter product is compact after adjoining its faces.

For a convergent sequence in this neighborhood, take a subsequence on which the parameters converge. Subtracting the approximate solutions and passing to the limit preserves orthogonality and the closed bound for the corrections.

Its boundary thus has only the parameter faces and correction sphere already excluded. By excision we may decrease the radius around the graph. The bound
\[
 CL^{5/2}\operatorname{radius}<1
\]
keeps every zero of the complementary homotopy on the graph, as in Subsection~\ref{reg:isolation}. Here $L$ is bounded on each fixed closed scale interval. This gives the stated local degrees in bounded open $H^1$ neighborhoods. The fixed weights are bounded and have integrable logarithms.

The compact-map comparison above therefore applies: by H\"older and the exponential bound the normalized sources on bounded $H^1$ or $C^0$ sets are bounded in some $L^q$, $q>1$, with the denominators bounded below by Jensen's inequality.

Their Green images are relatively compact in the common space $C^0\cap H^1$. Finite-dimensional reduction gives the same local degree in the space of continuous functions. This proves part~\textup{(ii)}.
\end{proof}
\begin{remark}
\label{regzero:scope}
At each marked point, the source strength may vanish in either component. The constants may depend on the fixed source strengths.
\begin{itemize}
\item[\textup{(i)}]
Fixed zeros enter the scalar shadow degree and the smooth terms at regular blow-up points. They contribute no additional local index.
\item[\textup{(ii)}]
The construction at prescribed blow-up points treats blow-up at a prescribed singular point.
\end{itemize}
\end{remark}
\section{Bubbling at regular and singular points}
\label{hyb:section}
In this section we consider configurations containing at least one singular bubble. We fix the generic source strengths and the common weights of Section~\ref{fnd:section}.

Let $U$ and $V$ be the numbers of ordinary and mixed regular points, respectively, and let $P$ and $Q$ be the sets of pure and mixed singular points. Put
\[
 n=U+V+|P|+|Q|,\qquad v=V+|Q|.
\]
The outer and inner profile exponents at a bubble are respectively
\[
 (s_a,t_a)=(1,0),\ (2,1),\ (A_a,0),\ (A_a+B_a,B_a).
\]
At the prescribed blow-up points, all positive profile exponents are nonintegers, with
\[
 A_a,B_a>1.
\]
The centers range over compact separated neighborhoods of the selected solutions of the shadow system. For fixed centers, the scalar equation has a compact smooth branch of solutions. Its linearization is invertible. Write
\begin{equation}\label{hyb:masses}
 M=4\pi\sum_a s_a,\qquad N=4\pi\sum_{t_a>0}t_a,
 \qquad \gamma=\rho _2-N>0,\qquad \rho _1=M+\mu.
\end{equation}
The scalar index is denoted by $m_S$. At mixed regular points choose the compact shape set
\[
 |B_b|\le B_0
\]
large enough to contain the limiting shapes of the selected solutions of the shadow system in its interior, as given by Lemma~\ref{reg:stress}. We impose both integral normalizations. Writing the common first normalization as $e^{-T}$, we obtain
\begin{equation}\label{hyb:matching}
 R_a^{2s_a}=d_a^{s_a}e^{-T},\qquad
 r_a^{2t_a}=b_a R_a^{2s_a}\quad(t_a>0).
\end{equation}
The coefficients $d_a,b_a$ are smooth and bounded above and away from zero. At a mixed regular point the second identity reads
\[
 r_b^2=\frac{\beta_b R_b^4(1+|B_b|^2)}4.
\]
We keep the factor $1+|B_b|^2$ when comparing local means and computing the orientation.
\subsection{Approximate solutions and the linearized operator}
For $q>0$ set
\[
 m_{q,R}(z)=\frac{4q^2R^{2q}|z|^{2q-2}}
 {(R^{2q}+|z|^{2q})^2},\qquad
 Z_{q,R}(z)=\frac{R^{2q}-|z|^{2q}}{R^{2q}+|z|^{2q}}.
\]
Thus $\int_{\C}m_{q,R}=4\pi q$ and $-\Delta Z_{q,R}=2m_{q,R}Z_{q,R}$. The mixed regular outer density is
\[
 m_{R,B}(z)=R^{-2}m_B(z/R),\qquad
 m_B(y)=\frac{16|y|^2}{(1+|y^2-B|^2)^2}.
\]
The corresponding exact two-component profile is given by
\begin{align}
 F&=R^4+\beta R^4|z|^2+|z^2-R^2B|^2,
 &Q_0&=\beta R^8+4R^4|z|^2+\beta R^4|z^2+R^2B|^2,\notag\\
 n_1&=4Q_0/F^2,&n_2&=16\beta R^8F/Q_0^2.
 \label{hyb:gram}
\end{align}
Here $Q_0$ denotes the polynomial and $Q$ the set of mixed singular points. The two tails are
\[
 L_1=4\log R+\log16-6\log|z|+\log(1+\beta|z|^2/4),\qquad
 L_2=\log\beta-2\log(1+\beta|z|^2/4).
\]
Let
\[
 H=h_2\exp\left(4\pi\sum_a(2t_a-s_a)G(\cdot,p_a)\right),
 \quad I=\int_MHe^w,\quad -\Delta_gw=2\gamma(He^w/I-1),
\]
and define the exterior field
\[
 u^0=K\binom{4\pi\sum_a s_aG(\cdot,p_a)}
 {4\pi\sum_a t_aG(\cdot,p_a)}+\binom{-w/2}{w}.
\]
Before imposing the two integral normalizations, put $c_1=-T$ and $c_2=\log(\gamma/I)$. All disk cut-offs below are smooth radial functions equal to one on a smaller concentric disk. On a radial disk add
\[
 (D_{s,R},-D_{s,R}/2),\qquad (-D_{t,r}/2,D_{t,r}),\qquad
 D_{q,R}=-2\log(1+(R/|z|)^{2q}),
\]
with fixed radial cut-offs; omit the second vector when there is no inner bubble. On a mixed regular disk add $\chi(\log n_1-L_1,\log n_2-L_2)$. Matching the tails determines \eqref{hyb:matching}.

We add the first-order corrections at the singular points specified below and subtract the surface means. This defines the approximate solution $B_T$.

We write all densities in local coordinates with the conformal area factor included. If $g=e^\varphi|dz|^2$, $\varphi(0)=0$, define
\[
 \psi_{1,a}=\log h_1+\varphi+u^0_1+(2s_a+2)\log|z|.
\]
At each center the explicit logarithm cancels the singular part of $\log h_1+u_1^0$. Matching on each outer disk gives
\[
 4s_a^2R_a^{2s_a}=e^{-T+\psi_{1,a}(0)}.
\]
For positive source strengths, the residual scalar weight is comparable to $\lvert z\rvert^{2(A+B-1)}$ at a pure prescribed center and to $\lvert z\rvert^{2(A-1)}$ at a mixed one. Both exponents are positive. Local Schauder regularity therefore gives
\[
 \psi_{1,a}\in C^{2,\theta}
\]
across such a center for some $\theta>0$. There put
\[
 \kappa_a=\Delta\psi_{1,a}(0).
\]
At a mixed regular point we also subtract the smooth tail of \eqref{hyb:gram}:
\begin{equation}\label{hyb:canonical-trace}
 J_{1,b}=\psi_{1,b}-\psi_{1,b}(0)-\log(1+\beta_b|z|^2/4),
 \qquad \kappa_b=\Delta\psi_{1,b}(0)-\beta_b.
\end{equation}
We choose a local coordinate in which the smooth part of the first component has zero gradient. Indeed, under
\[
 z=\zeta+a\zeta^2,
\]
the $-6\log|z|+2\log|dz/d\zeta|$ terms change its linear part by $-2\Re(a\zeta)$. The linear part of the smooth second component changes by $4\Re(a\zeta)$, preserving the combination of twice the first and the second.

The pullback of the smooth linear term contributes $2\Re(a\psi_{1,z}(0)\zeta^2)$ and the pole and the factor from the change of variables together contribute $-\Re(a^2\zeta^2)$. Choose
\[
 a=\psi_{1,z}(0).
\]
Then the radial subtraction leaves $T_b=\psi_{1,zz}(0)+\psi_{1,z}(0)^2$ unchanged, and
\[
 J_{1,b}=\kappa_b|z|^2/4+\Re(T_bz^2)+o(|z|^2),\qquad
 J_{2,b}=F_b\cdot x+O(|x|^2).
\]
These Taylor coefficients include the exterior fields from all the other singular masses. At $\rho_1=M$, $\kappa_b$ is the first forcing in the convention
\[
 \Delta U+Ke^U=\kappa
\]
used below. It is $-f_1^0$ in the convention of \eqref{reg:canonical}--\eqref{reg:wall-forcing}. The subtraction of $\beta_b$ in \eqref{hyb:canonical-trace} accounts for the residual density in this identification.

Work on the mean-zero space $\Hh\subset H^1(M)^2$ with scalar product
\[
 \langle u,v\rangle_{\Hh}=\int_M\nabla u^tK^{-1}\nabla v.
\]
For a fixed radial cut-off, center every source on its own disk, replacing $\chi mZ$ by
\[
 \chi m(Z-\int\chi mZ/\int\chi m).
\]
Retain the ordinary translations, the inner translations and outer shapes at each mixed regular point, and one outer radial source at a singular reference point. The translation source is $m_{1,R}(z)y_j/(1+|y|^2)$; the shape source is $\partial_{B_j}m_{R,B}$.

Write these centered component sources as $h_\alpha e_{i(\alpha)}$ with respect to $dV_g$; thus their coordinate measures have the Euclidean densities specified above. Define their Cartan Green representatives and functionals by
\begin{align}
 S_\alpha&=K\bigl(G*(h_\alpha e_{i(\alpha)})\bigr),\qquad
 E_T=\bigcap_\alpha\ker L_\alpha,\notag\\
 L_\alpha u&=\langle S_\alpha,u\rangle_{\Hh}
           =\int_Mh_\alpha u_{i(\alpha)}\,dV_g.
 \label{hyb:source-convention}
\end{align}
Accordingly $\mathcal F_\mu(u)=\sum_\alpha c_\alpha S_\alpha$ adds $+Kq$ to the equation for $-\Delta_gu$, where
\[
 q=\sum_\alpha c_\alpha h_\alpha e_{i(\alpha)}
\]
has zero mean in each component. We use this sign for the multipliers in the local identities below. The radial source at the reference point is supported on its disk and has integral zero there.

Let $\omega_a^{\rm out}$ and $\omega_b^{\rm in}$ be the exactly normalized cut-off outer and inner probability measures. Their densities are Euclidean densities on their coordinate disks. Define
\begin{equation}\label{hyb:model-probabilities}
 \nu_{1,T}=\sum_a\frac{4\pi s_a}{M}\omega_a^{\rm out},\qquad
 \nu_{2,T}=\sum_{t_b>0}\frac{4\pi t_b}{\rho _2}\omega_b^{\rm in}
             +\frac{\gamma}{\rho _2}\frac{He^w}{I}\,dV_g.
\end{equation}
The measure $(He^w/I)\,dV_g$ has mass one on $M$. Consider the self-adjoint model operator
\begin{equation}\label{hyb:model-operator}
 \langle\mathcal T_T\phi,\psi\rangle_{\Hh}
 =\langle\phi,\psi\rangle_{\Hh}
 -M\Cov_{\nu_{1,T}}(\phi_1,\psi_1)
 -\rho _2\Cov_{\nu_{2,T}}(\phi_2,\psi_2).
\end{equation}
\begin{theorem}[Invertibility of the model operator]
\label{hyb:model-inverse}
Write $A=P_{E_T}\mathcal T_T|_{E_T}$. Uniformly on the compact set of separated configurations, the following hold.
\begin{itemize}
\item[\textup{(i)}]
$\|A^{-1}\|\le CT$, and the negative index of $A$ is
\begin{equation}\label{hyb:index}
 \nu=m_S+n-1+v
       +2\sum_a(\lceil s_a\rceil-1)
       +2\sum_{t_b>0}(\lceil t_b\rceil-1).
\end{equation}
\item[\textup{(ii)}]
Exactly $n-1+v$ eigenvalues of $A$ tend to zero. They lie in $[c/T,C/T]$. The rest of the spectrum stays uniformly away from zero.
\end{itemize}
\end{theorem}
\begin{proof}
\noindent\textit{Step 1. Estimates on the individual profiles.}\quad We center each local measure separately. For the calculation in cylindrical coordinates and the estimates at the two scales, we use Lemma~\ref{hyb:centered-estimates}.

For each cut-off radial probability, the Poincar\'e inequality on circles and the Cauchy--Schwarz inequality for the circular mean give
\[
 \Var_{\omega_{q,R}}f\le C\int_M|\nabla f|^2.
\]
In scaled coordinates the contribution outside $[L^{-1},L]$ is bounded by $CL^{-2q}(1+\log L)\int|\nabla f|^2$. Extra fixed logarithmic weights only increase the logarithmic power.

For bounded $B$, $m_B$ is uniformly comparable to $m_{2,1}$, so the same estimates hold. Local compactness of the local Sobolev embedding and these tails give strong weighted compactness of centered fields.

We next estimate $L_\alpha$ and its first derivatives. Write a source as $h_\alpha e_i$, so
\[
 L_\alpha v=\int h_\alpha v_i;
\]
its integral vanishes on its own disk. Its total variation is bounded, and after rescaling it is a bounded kernel function times the corresponding density.

Subtract the corresponding weighted mean of $v_i$ and use Cauchy--Schwarz and the weighted variance estimate above. Together with the total-variation bound, this gives
\begin{equation}\label{hyb:modulation-bounds}
 \|L_\alpha\|_{(C^0)^*}+\|L_\alpha\|_{\Hh^*}
 +\|\partial_{\rm sc}L_\alpha\|_{(C^0)^*}
 +\|\partial_{\rm sc}L_\alpha\|_{\Hh^*}\le C.
\end{equation}
Here $\partial_{\rm sc}$ denotes any first parameter derivative with each regular center scaled by its own translation width. For these derivatives, differentiate the exact local centering as well. The differentiated sources have zero integral on each disk.

After rescaling, their densities are bounded by constant multiples of the radial densities.

Use $m_{q,1}$ for radial scale derivatives, $m_{2,1}$ for the derivatives of the outer mixed profile with respect to shape, and $m_{1,1}$ for spatial derivatives of the regular sources (using the outer scale for an outer source).

The latter bound holds also for the derivative of the outer mixed density: it is $O(|y|)$ at zero and decays faster than $m_{1,1}$ at infinity. Differentiating an outer profile by a mixed center scaled by $r$ introduces the bounded factor $r/R$.

The tails produced by derivatives of the cut-offs are bounded by positive powers of the radii, and derivatives of the matching coefficients are bounded on the compact set of separated configurations.

This proves the derivative bounds in \eqref{hyb:modulation-bounds}. Thus cancellation on each disk gives the uniform $H^{-1}$ estimates.

Let $B_T^{\rm cen}$ be the operator on $\Hh$ whose quadratic form is the sum of the individual bubble variances and the variance of the residual measure, each centered separately.

For a local source associated with a kernel function, let $\Psi$ denote its Green representative. Its active component is $Z-Z(\infty)+o(1)$ and its inactive component is minus one half of that expression. At a mixed regular outer bubble the fields are
\[
 Z_{0,B}=\frac{1-|y^2-B|^2}{1+|y^2-B|^2},\quad
 Z_{1,B}=\frac{4\Re(y^2-B)}{1+|y^2-B|^2},\quad
 Z_{2,B}=\frac{4\Imm(y^2-B)}{1+|y^2-B|^2}.
\]
Their squared norms in $L^2(m_B\,dA)$ are $8\pi/3$, $32\pi/3$ and $32\pi/3$, respectively. The radial function is orthogonal to both shape derivatives in this space. The substitution $w=y^2-B$ proves these identities. Dilations of radial profiles with exponent $q$ have squared norm $4\pi q/3$. For all these fields,
\begin{equation}\label{hyb:jacobi-error}
 \|(I-B_T^{\rm cen})\Psi\|_{\Hh}\le Ce^{-cT}T^C.
\end{equation}
Indeed, for a kernel function \(Z\) in component \(i\) the Green representative is approximated in energy by
\[
 \frac12 Ke_i\chi(Z-Z(\infty))
\]
in the scaled coordinates, with the surface mean subtracted. Applying the Laplacian to this expression gives the cut-off source corresponding to $Z$ and fixed-annulus terms containing the tails of \(Z-Z(\infty)\) and its gradient, each bounded by a positive power of the radius.

The adjustment from the cut-off source corresponding to $Z$ to its exactly centered version is a positive power of the radius times the cut-off mass, since the source corresponding to $Z$ has integral zero on the plane.

Its energy-dual norm has at most the additional factor \(C\sqrt T\) from the logarithmic Green mean.

On its own scale the centered active component gives the source corresponding to $Z$; all cut-off errors are bounded by the same estimates.

At a nested pair let \(\delta=r/R\). In inner coordinates, split at
\[
 |y|=\delta^{-1/2}.
\]
On the inner part the outer field changes from its value at zero by \(O(\delta^{\beta/2})\), with some \(\beta>0\); on the other part the inner probability is \(O(\delta^t)\). In outer coordinates, split at
\[
 |y|=\delta^{1/2}.
\]
On the outer part the inner field changes from its limiting value at infinity by \(O(\delta^{\lambda/2})\), with
\[
 \lambda>0
\]
(the inner translations have \(\lambda=1\)); the remaining outer probability is \(O(\delta^s)\). Taking the weighted \(L^2\) norms and using the centered variance estimate gives positive powers in the corresponding energy-dual norms.

On separated disks the cut-off model field vanishes up to its energy error and an additive constant.

For the variance of the residual measure, split the disk at the square root of its local radius. On the smaller disk the bounded residual density has mass bounded by a positive power of this radius.

On the complement the model field minus its value at infinity satisfies such a bound. These estimates prove \eqref{hyb:jacobi-error}, since all radii and nested ratios are exponentially small.

\noindent\textit{Step 2. The operator on the orthogonal complement.}\quad Let $\mathcal K_T$ be the span of the Green representatives of all sources corresponding to local kernel functions, including each dilation separately, and let
\[
 W_T=\mathcal K_T^\perp.
\]
Then
\[
 A_T=P_{W_T}(I-B_T^{\rm cen})|_{W_T}
\]
has a uniformly bounded inverse. If this fails, choose a sequence with unit energy whose image tends to zero. Testing against the kernel functions and using \eqref{hyb:jacobi-error} eliminates its projection multipliers.

Away from the blow-up points the limit solves the homogeneous linearized scalar equation. Invertibility implies that this limit vanishes. At every bubble the active centered limit solves
\[
 -\Delta z=2mz;
\]
the inactive component plus one half of the active component is harmonic. A source of total mass zero at the inner scale tends to zero in distributions on the outer scale. The puncture is removable by finite energy. On the inner scale the outer centered source vanishes.

For each noninteger profile exponent $q$, separation of variables below shows that the dilation spans the finite-energy scalar kernel space.

For the regular model $q=1$ and the nonradial model $m_B$, the kernel spaces follow from Lemma~\ref{reg:jacobi-kernel}. Orthogonality eliminates them.

Strong weighted compactness then makes every variance tend to zero, contradicting the unit energy when the equation is tested against the sequence itself. The index of $A_T$ is
\begin{equation}\label{hyb:centered-index}
 m_S+2\sum_a(\lceil s_a\rceil-1)
       +2\sum_{t_b>0}(\lceil t_b\rceil-1).
\end{equation}
For the spherical metric associated with the radial profile with exponent $q$, the weighted eigenvalues are $\lambda(\lambda+1)$,
\[
 \lambda=k+|\ell|/q,\qquad k\ge0,
\]
$\ell\in\mathbb Z$. Removing the constant, the eigenvalues below two are precisely $k=0$, $0<|\ell|<q$. These eigenvalues can be obtained directly by separating the angular frequency $\ell$. With
\[
 x=q\log r,\qquad \mu_\ell=|\ell|/q,
\]
the equation $-\Delta\phi=\Lambda m_{q,1}\phi$ becomes
\[
 -f_{xx}+\mu_\ell^2f=\Lambda\operatorname{sech}^2x\,f.
\]
Putting $z=\tanh x$ and $f=(1-z^2)^{\mu_\ell/2}p(z)$ gives
\[
 (1-z^2)p''-2(\mu_\ell+1)zp'
       +\{\Lambda-\mu_\ell(\mu_\ell+1)\}p=0.
\]
Set
\[
 w=(1-z^2)^{\mu_\ell},\qquad b=(1-z^2)w.
\]
The operator $-w^{-1}(b p')'$ is symmetric on polynomials for $w\,dz$ and preserves each degree subspace, with leading multiplier $k(k+2\mu_\ell+1)$ in degree $k$. Orthogonalizing $1,z,z^2,\ldots$ gives an eigenpolynomial at each degree.

The corresponding fields have finite spherical energy and are valid tests in the weak weighted equation; polynomials are dense in $L^2(w\,dz)$, so these tests determine every eigenvalue. Hence
\[
 \Lambda=\mu_\ell(\mu_\ell+1)+k(k+2\mu_\ell+1)
        =(k+\mu_\ell)(k+\mu_\ell+1).
\]
The constant is removed by centering and the dilation has eigenvalue two. This gives the number of negative angular directions in the formula above and the additional frequency-$q$ kernel space pair precisely when $q$ is an integer.

The nonradial $q=2$ index stays two along bounded $B$ paths because the weighted embedding is compact and the kernel space remains three-dimensional.

For the lower bound, transplant finite families of negative eigenfunctions in their respective scaled coordinates, taking the inactive component to be minus one half of the active component.

Constants in the weighted eigenfunctions may be subtracted before cutting off, since the energy and centered variances are unchanged by constants.

Truncate the local energy fields and use logarithmic cut-offs at zero and infinity. Their cut-off energies tend to zero for bounded fields. This gives approximation on fixed rescaled annuli, which we enlarge after taking the scaling limit.

For the negative directions of the limiting problem on $M$, use cut-offs about the finitely many centers whose Dirichlet energies tend to zero.

The transplanted supports are then disjoint for large \(T\); their centered variances converge on the indicated scales and vanish on the other scales.

Projecting orthogonally to the kernel functions changes these finite-dimensional test functions by \(o(1)\), by orthogonality to the limiting kernel.

For the upper bound, use the fixed negative spectral gap and pass a basis of any larger negative spectral subspace to all local and macroscopic limits.

On disjoint scaled annuli and the macroscopic complement the energy is lower semicontinuous, and all centered variances converge strongly. On a local component the identity
\[
 \nabla(u,v)^tK^{-1}\nabla(u,v)
   =\tfrac12|\nabla u|^2+\tfrac23|\nabla(v+u/2)|^2
\]
gives the active threshold two. Any nonzero normalized linear combination of the chosen basis has limiting form bounded above by minus the fixed spectral gap.

In particular its limiting vector cannot be zero, so the limiting negative subspace has the same dimension, contradicting the direct-sum index.

\noindent\textit{Step 3. The differences of the local means.}\quad The exact variance decomposition now gives
\[
 \mathcal T_T=I-B_T^{\rm cen}-V_TC V_T^*,
\]
where the components of $V_T^*\phi$ are the differences of means
\[
 \langle\phi_1\rangle_a-\langle\phi_1\rangle_n\ (a<n),
 \qquad \langle\phi_2\rangle_b-\langle\phi_2\rangle_g\ (b\le v).
\]
The two positive matrices are
\[
 C_1=\diag(M_a)_{a<n}-\frac{(M_a)_{a<n}(M_a)^t_{a<n}}{M},\qquad
 C_2=\diag(N_b)-\frac{(N_b)(N_b)^t}{\gamma+\sum_bN_b},
\]
where $M_a=4\pi s_a$ and $N_b=4\pi t_b$. Project $V_T$ to $W_T$. For a vector $z$ in these coordinates, put $c_a=z_{1,a}$ for $a<n$,
\[
 c_n=-\sum_{a<n}z_{1,a},
\]
and $d_a=z_{2,b}$ at the outer point supporting inner bubble $b$, zero otherwise. Logarithmic Green averaging yields
\begin{equation}\label{hyb:contrast-gram}
 V_T^*V_T=T G_0+O(1),\qquad
 z^tG_0z=\frac1{2\pi}
 \left\{\sum_a\frac{c_a^2-c_ad_a}{s_a}
       +\sum_b\frac{d_{a(b)}^2}{t_b}\right\}>0\quad(z\ne0).
\end{equation}
The nested cross term uses $\log(1/R_{a(b)})$; the inner self term uses $\log(1/r_b)$. Completing squares proves positivity. Nonradial logarithmic moments, remote interactions and the projection off $\mathcal K_T$ change only the $O(1)$ term. Centered logarithmic moments give
\[
 V_T^*B_T^{\rm cen}V_T=O(1),\qquad \|B_T^{\rm cen}V_T\|=O(1).
\]
With $B_W=P_WB_T^{\rm cen}|_W$, the exact identity $A_T^{-1}=I+B_W+B_WA_T^{-1}B_W$ therefore implies
\begin{equation}\label{hyb:resolvent}
 F_T:=V_T^*A_T^{-1}V_T=T G_0+O(1).
\end{equation}
Applying the finite-rank resolvent formula to $H_T=A_T-V_TC V_T^*$ gives
\[
 \|H_T^{-1}\|\le C,\qquad
 \ind H_T=\ind A_T+n-1+v,\qquad
 \|H_T^{-1}V_T\|=O(T^{-1/2}),
\]
because $C^{-1}-F_T$ is negative and its inverse is $O(T^{-1})$.

\noindent\textit{Step 4. The small eigenvalues.}\quad Finally restore the relative radial fields, leaving the reference radial field removed. We use the label $n$ for the chosen singular reference outer point, so the retained radial coordinates have
\[
 x_n=0.
\]
Normalize the active-component mean at each bubble to one. For a regular mixed outer field its inactive inner value is $-(1+|B_b|^2)^{-1}$. Their limiting contrast matrix, restricted to $x_n=0$, is
\[
 D_B(x,y)=\big((x_a-x_n)_{a<n},\ (y_b-\tau_bx_{a(b)})_b\big),\qquad
 \tau_b=1\ \hbox{or}\ (1+|B_b|^2)^{-1}.
\]
It is triangular and invertible. The shape derivatives contribute to the differences of the inner means. Projecting them out changes $V_T^*V_T$ by $O(1)$.

Their Dirichlet inner products with the radial functions vanish in the outer model. The sources defining $E_T$ are among these sources corresponding to local kernel functions, so $W_T\subset E_T$.

When $n-1+v>0$, project the Green representatives of the retained individual radial sources onto $E_T$. The projections change the matrix of their inner products by $o(1)$. For different scales or disjoint disks this follows from \eqref{hyb:jacobi-error}.

At the same scale the angular integral of the radial function against each shape derivative vanishes. The remaining radial fields are also asymptotically orthogonal to the omitted reference radial source.

Thus their Green representatives give a uniformly nondegenerate frame of the orthogonal complement of $W_T$ in $E_T$.

Let $D_T$ be the matrix of mean contrasts in this frame. The projections subtract $o(1)$ linear combinations of the removed sources' Green representatives.

Each such representative has bounded weighted means on every disk, by the bounded kernel functions on the plane and the cut-off estimates. Thus the projections change the contrasts by $o(1)$ as well, and
\[
 D_T\to D_B.
\]
Up to the exponentially small errors in \eqref{hyb:jacobi-error}, its radial block is $-D_T^tCD_T$ and its block into $W_T$ is $-V_TCD_T$. The finite-rank resolvent identity gives
\[
 H_T^{-1}V_TC=A_T^{-1}V_T(C^{-1}-F_T)^{-1}.
\]
Eliminating the complementary variables in the principal blocks gives
\begin{align*}
 &-D_T^t\left(C+CV_T^*H_T^{-1}V_TC\right)D_T\\
 &\hspace{12mm}=-D_T^t(C^{-1}-F_T)^{-1}D_T
              =D_T^t(F_T-C^{-1})^{-1}D_T.
\end{align*}
Including these errors gives
\begin{equation}\label{hyb:positive-schur}
 D_T^t(F_T-C^{-1})^{-1}D_T+O(e^{-cT}T^C)
   =T^{-1}D_B^tG_0^{-1}D_B+o(T^{-1}).
\end{equation}
The matrix $G_T^{\rm rad}$ of inner products of the radial functions is uniformly positive definite. Write the principal off-diagonal block as
\[
 K_T^{\rm rad}=-V_TCD_T.
\]
Choose a fixed $\lambda_0>0$ small enough that the $W_T$ block of $\mathcal T_T-\lambda I$ is invertible with a common bound for
\[
 |\lambda|\le\lambda_0.
\]
The uniform inverse bound for $H_T$ and the exponentially small change of that block give this choice by the resolvent identity. On this entire interval its reduced matrix $S_T(\lambda)$ satisfies
\[
 S_T'(\lambda)=-G_T^{\rm rad}
 -(K_T^{\rm rad})^*(H_T-\lambda I)^{-2}K_T^{\rm rad}
 +O(e^{-cT}T^C)
 =-G_T^{\rm rad}+O(T^{-1}).
\]
Here $\|(H_T-\lambda I)^{-1}K_T^{\rm rad}\|=O(T^{-1/2})$ follows from the preceding bound on $H_T^{-1}V_T$ and
\[
 (H_T-\lambda I)^{-1}=(I-\lambda H_T^{-1})^{-1}H_T^{-1};
\]
the exponentially small perturbations preserve this estimate. Its value at zero is positive of order $T^{-1}$ by \eqref{hyb:positive-schur}. It is still positive at $\lambda=c/T$ for sufficiently small $c>0$, and negative at $\lambda=C/T$ for sufficiently large fixed $C$.

Strict monotonicity therefore gives exactly $n-1+v$ positive roots in $[c/T,C/T]$, counted with multiplicity, with no other roots in the fixed interval.

If $n-1+v=0$, there is no relative radial block. This proves the gap and \eqref{hyb:index}, including the sign of the relative radial block.
\end{proof}
\subsection{Estimates for the approximate solutions}
At each singular outer point we add $(P,-P/2)$, where
\begin{equation}\label{hyb:firstjet}
 P=P_{q,R,g}=-\frac{2q}{q-1}\frac{g\cdot x}{1+(|x|/R)^{2q}},
 \qquad (-\Delta-2m_{q,R})P=2m_{q,R}g\cdot x.
\end{equation}
At a mixed singular point also add $(-Q_1/2,Q_1)$ with $Q_1=P_{t,r,\widehat g_2}$ and $\widehat g_2=\nabla\psi_2(0)+s\nabla\psi_1(0)/(s-1)$. Use fixed radial cut-offs and subtract means. Put
\[
 s_* =\max_a s_a>1,\qquad \mathcal L=e^{-T/s_*}.
\]
Choose $p>1$ sufficiently close to one that
\begin{equation}\label{hyb:p-choice}
 p<2,\quad 2/p-1>1/s_*,\quad 2A_a-2+2/p>2\quad(a\in Q),
\end{equation}
and define
\begin{equation}\label{hyb:epsilon}
 \epsilon(T)=\sum_{a\in P\cup Q}R_a^{2/p}
       +\sum_{b\ {\rm mixed\ regular}}R_b^{4/p-2}
       +\sum_{a\ {\rm ordinary}}R_a^{2/p-1}.
\end{equation}
Then $\epsilon(T)^2T^C=o(\mathcal L)$ for every fixed $C$. We shall use this estimate in the common scale equation. We estimate the shape equations at the respective regular points below. For $d=|x|$, set
\begin{equation}\label{hyb:linear-couplings}
 \mathfrak l_{s,R}(x)=-\frac{s}{s-1}(g_1\cdot x)
                  \frac{(d/R)^{2s}}{1+(d/R)^{2s}}.
\end{equation}
\begin{lemma}[Moments and logarithmic scale derivatives]
\label{hyb:moment-estimates}
Fix $q>0$, $d_0>0$, $p\geq1$ and $b\geq0$. Assume $p(b+2q-2)>-2$.
\begin{itemize}
\item[\textup{(i)}]
We have
\[
 \big\||z|^b m_{q,R}\big\|_{L^p(D_{d_0})}
 \leq C
 \begin{cases}
 R^{b-2+2/p},& b<2q+2-2/p,\\
 R^{2q}(1+|\log R|)^{1/p},&b=2q+2-2/p,\\
 R^{2q},&b>2q+2-2/p.
 \end{cases}
\]
The same bounds hold with $m_{q,R}$ replaced by $\partial_T^j m_{q,R}$, $j=1,2$, when $R^{2q}=d^q e^{-T}$ with $d$ independent of $T$.
\item[\textup{(ii)}]
At a mixed singular point, set
\[
 s=A+t,\quad A,t>1,\quad r^{2t}=b_0R^{2s},\quad
 S=m_{s,R},\quad h=(r/|z|)^{2t}.
\]
Take $b_0$ in a compact subset of $(0,\infty)$ and keep $b_0$ and the Taylor coefficients fixed under differentiation. Then
\[
 \|Sh\|_p\leq C R^{2A-2+2/p},\qquad
 \int_{D_{d_0}}Sh\leq C R^{2A}.
\]
Both estimates hold after two $T$ derivatives.
\item[\textup{(iii)}]
Under the assumptions of \textup{(ii)}, the first-order corrections in \eqref{hyb:firstjet} and \eqref{hyb:linear-couplings} satisfy
\[
 \int_{D_{d_0}}|z|
    \bigl(|\partial_T^j(SQ_1)|+
          |\partial_T^j(I_0\mathfrak l_{s,R})|\bigr)
 \leq C R^{2A+2},\qquad j=0,1,2.
\]
\end{itemize}
All constants are uniform on the fixed compact parameter sets.
\end{lemma}
\begin{proposition}[Residual estimates]
\label{hyb:trial-comparison}
Let $\mathcal F_\mu$ be the normalized variational gradient.
\begin{itemize}
\item[\textup{(i)}]
For some $a>0$, the corrected approximate solution $B_T$ satisfies
\[
 \|D\mathcal F_\mu(B_T)-\mathcal T_T\|
       \le C(e^{-aT}+|\mu|)T^C.
\]
At fixed $\mu$, centers and shape parameters,
\begin{equation}\label{hyb:residual}
 \|\partial_T^j\mathcal F_\mu(B_T)\|_{\Hh}
       \le C(\epsilon(T)+|\mu|)T^C,\qquad j=0,1.
\end{equation}
\item[\textup{(ii)}]
The common unnormalized masses are bounded below by a positive constant. The normalized densities of $B_T$ are $O(e^{CT})$ in $L^\infty$.
\end{itemize}
\end{proposition}
\begin{proof}
\noindent\textit{Step 1. The local residual.}\quad We estimate the residual separately on the regular and singular disks. The moment bounds and their logarithmic-scale derivatives used below are proved in Lemma~\ref{hyb:moment-estimates}. For the densities in \eqref{hyb:gram}, put \(U=R^2\), \(u=|z|^2\),
\[
 A_0=U^2+|z^2-UB|^2,\qquad x=u+r^2.
\]
Then \(A_0\asymp U^2+u^2\) on compact shape sets and
\[
 F=A_0+\beta U^2u,\qquad
 \frac{Q_0}{4R^4x}=1+\frac{\beta u}{4}+O(Uu/x),\qquad
 \beta F=\beta u^2+4r^2+O(Uu).
\]
For the first density the numerator term other than \(16R^4u\) is \(O(U^2 A_0)\), while replacing \(F^2\) by \(A_0^2\) in that leading term gives an error of \(O(U^2/A_0)\) pointwise. For the second, divide the last display by \(x^2\) and by the square of the middle expression.

In the term with numerator \(\beta u^2\) the replacements of \(u^2/x^2\) by one and of the middle expression by \(1+\beta u/4\) cost \(O(U/x)\). In the term with numerator \(4r^2\) the replacement of the middle expression by one gives an error of \(O(U/x)\), using
\[
 r^2u/x^2\le 1/4.
\]
Thus, with
\[
 n_{\rm bg}=\beta(1+\beta|z|^2/4)^{-2},
\]
\[
 |n_1-m_{R,B}|\le CU^2/A_0,\qquad
 |n_2-m_{1,r}-n_{\rm bg}|\le CU/x .
\]
The first bound integrates to \(O(U)\), the second to \(O(U(1+|\log U|))\). On \(u\le U\) the latter bound is controlled by \(Cr^2/x^2\), since \(r^2\asymp U^2\); on \(u\ge U\) it is controlled by \(CU/(U+u)\). This gives
\begin{align*}
 \|n_1-m_{R,B}\|_1&\le CR^2,&
 \|n_2-m_{1,r}-n_{\rm bg}\|_1&\le CR^2(1+|\log R|),\\
 |n_2-n_{\rm bg}|&\le C\left\{
  \frac{r^2}{(r^2+|z|^2)^2}+\frac{R^2}{R^2+|z|^2}\right\}.
\end{align*}
The same bounds hold after one logarithmic-scale derivative. In particular,
\[
 \||z|^2 n_1\|_p\le C_pR^{2/p},\qquad
 \||z|(n_2-n_{\rm bg})\|_p\le C_p(r^{2/p-1}+R^2).
\]
The exact residual before normalization on this disk is
\begin{equation}\label{hyb:gram-residual}
 -K\binom{n_1(e^{J_1}-1)}{(n_2-n_{\rm bg})(e^{J_2}-1)}.
\end{equation}
Indeed \(-\Delta L=K(0,n_{\rm bg})^t\) holds off the center, while the residual density there is \(n_{\rm bg}e^{J_2}\). In the patched field at the center the logarithmic part of \(L_1\) cancels the exterior Green pole, so no point mass is left in the residual.

The difference of the two expressions for the residual density has already been included in \((n_2-n_{\rm bg})(e^{J_2}-1)\). With
\[
 J_1=O(d^2),\qquad J_2=O(d),
\]
the bound is $C_pR^{4/p-2}$. At an ordinary point the residual is $-K(m_{1,R}(e^{J_1}-1),R^2\gamma_0(z))^t$, with $J_1=O(|z|)$, and has norm $C_pR^{2/p-1}$.

At a mixed singular point write $S=m_{s,R}$, $I_0=m_{t,r}$, $h=(r/|z|)^{2t}$, and $\psi_2=\log h_2+\varphi+u_2^0-(2A-2)\log|z|$. Put $J_i=\psi_i-\psi_i(0)$ and match $4t^2r^{2t}=e^{c_2+\psi_2(0)}R^{2s}$. Before angular correction the exact densities are
\[
 N_1^0=Se^{J_1}(1+h),\qquad
 N_2^0=I_0e^{J_2}+g_{\rm loc}(1+h)^{-2},\qquad
 g_{\rm loc}=e^{c_2+\psi_2}|z|^{2A-2}.
\]
After \eqref{hyb:firstjet}, the residual before multiplication by $-K$ is
\begin{align*}
 E_1&=S\{(1+h)e^{J_1+P-Q_1/2}-1-P-g_1\cdot x\},\\
 E_2&=I_0\{e^{J_2-P/2+Q_1}-1-Q_1-\widehat g_2\cdot x\}
       +g_{\rm loc}\{(1+h)^{-2}e^{-P/2+Q_1}-1\}.
\end{align*}
There are two linear coupling terms in these expressions. Use \eqref{hyb:linear-couplings}. Then $J_2-P/2-\widehat g_2\cdot x=\mathfrak l_{s,R}+O(d^2)$. The exponents $J_1+P-Q_1/2$ and $J_2-P/2+Q_1$ are uniformly $O(d)$, so Taylor expansion gives
\begin{align*}
 E_1&=-\tfrac12SQ_1+O(Sd^2)+Sh e^{J_1+P-Q_1/2},\\
 E_2&=I_0\mathfrak l_{s,R}+O(I_0d^2)
       +g_{\rm loc}\{(1+h)^{-2}e^{-P/2+Q_1}-1\}.
\end{align*}
Since $\mathfrak l_{s,R}=O(R)$ when $d\asymp R$, we retain this linear term in the following estimates. For $j=0,1$,
\begin{equation}\label{hyb:linear-coupling-bounds}
 \|\partial_T^j(SQ_1)\|_p+
 \|\partial_T^j(I_0\mathfrak l_{s,R})\|_p
       \le C_pR^{2A-1+2/p}.
\end{equation}
Indeed $r^{2t}\asymp R^{2s}$. On $r\le d\le R$, each product is bounded by $Cd^{2A-1}$, whose $L^p$ norm is $O(R^{2A-1+2/p})$. On $d\le r$ the two bounds for the norms are, respectively,
\[ CR^{-2s}r^{2s-1+2/p},\qquad
 Cr^{-2t}R^{-2s}r^{2t+2s-1+2/p}. \]
They have the same smaller order. On $d\ge R$ use
\[
 |SQ_1|\le CR^{2s}r^{2t}d^{-2s-2t-1},\qquad I_0|\mathfrak l_{s,R}|\le Cr^{2t}d^{-2t-1}
\]
to obtain the stated bound. At fixed $x$, logarithmic-scale differentiation preserves these majorants: $|\partial_T^jS|\le C_jS$, $|\partial_T^jI_0|\le C_jI_0$, $|\partial_T^jQ_1|\le C_j|Q_1|$ and
\[
 |\partial_T^j\mathfrak l_{s,R}|\le C_j|\mathfrak l_{s,R}|
\]
for the first two derivatives used here and in the radial variation. Since $A>1$, \eqref{hyb:linear-coupling-bounds} is smaller than $CR^{2/p}$. Splitting the remaining terms at $r$ and $R$ now gives
\[
 \|E_1\|_p+\|E_2\|_p
 \le C_p(R^{2/p}+r^{2/p}+R^{2A-2+2/p}+e^{-T}T^C),
\]
and the same bound after one $T$ derivative. For example
\[
 Sh\le C|z|^{2A-2}/(1+(|z|/R)^{2s})^2.
\]
The linear coupling terms are odd and have zero integral against every radial test on the bubble disk. Their pairings with a first-order correction tangent, which is $O(d)$, obey the stronger moment bound
\begin{equation}\label{hyb:linear-coupling-moments}
 \int_{D_{d_0}}d\bigl(|\partial_T^j(SQ_1)|+
          |\partial_T^j(I_0\mathfrak l_{s,R})|\bigr)\,dA
      \le CR^{2A+2},\qquad j=0,1.
\end{equation}
Integration over the three regions proves this estimate. In the outer region the second moment uses
\[
 \int_R^{d_0} d^{1-2t}\,dd=O(R^{2-2t})
\]
for $t>1$, so its lower endpoint determines the bound. Thus these terms do not change either the integrated mass order or the leading radial coefficient, including its first logarithmic-scale derivative. The terms from the residual density obey
\[
 \|g_{\rm loc}P\|_p\le C_pR^{2A-1+2/p},\qquad
 \|g_{\rm loc}Q_1\|_p\le C_p(r^{2A-1+2/p}+r^{2t})(1+|\log r|)^{1/p}.
\]
In particular the first is estimated at its own outer radius. Pure singular points have bound $C_p(R^{2/p}+R^{2A})$. Fixed-annulus cut-off errors are $O(R^{2s})$, $O(r^{2t})$, or $O(R^2)$, and outside the disks the first density is $O(e^{-T})$ while the second is the residual density.

\noindent\textit{Step 2. The two normalization factors.}\quad Let $N_i=e^{c_i}h_ie^{\widehat B_{T,i}}$, before surface mean subtraction, and $M_i=\int_MN_i$. Integration of the unnormalized residual and its $T$ derivative gives
\[
 |M_i-\rho_i^0|+|\partial_TM_i|\le C_p\epsilon(T),
\]
where $\rho^0=(M,\rho_2)$. The normalized density is exactly $\rho_iN_i/M_i$. H\"older with exponent $1+1/T$ gives
\[
 \|N_i/M_i-1\|_{H^{-1}}\le C\sqrt T,
\]
and a polynomial bound for its $T$ derivative. Thus the normalization error is the mass error in the denominator multiplied by this $H^{-1}$ norm. This proves \eqref{hyb:residual}.

\noindent\textit{Step 3. Comparison of the two operators.}\quad For the total-variation comparison first compare each unnormalized density to the sum of its model concentrating densities and, in component two, the term from the residual density.

On an ordinary or pure singular disk this uses the first absolute radial moment of \(m_{s,R}\) and the second-component tail computed above. On a regular mixed disk, use the two $L^1$ comparisons for $n_1,n_2$ and the bounds for $d^2n_1$ and $d(n_2-n_{\rm bg})$, since
\[
 J_1=O(d^2),\qquad J_2=O(d).
\]
On a mixed singular disk the component-one difference is bounded by \(C(Sd+Sh)\). The component-two difference is bounded by
\[
 C\bigl(I_0d+g_{\rm loc}\min(1,h)
          +g_{\rm loc}(|P|+|Q_1|)\bigr).
\]
Integrating the first moments and splitting the other terms at \(r,R\) as above gives positive powers of the radii, allowing logarithmic factors.

The fixed-annulus errors and the masses lost on replacing the full scalar profiles by the exactly normalized cut-off profiles are also positive powers. Dividing by \(M_i\) and using the positive lower bounds and the mass errors already proved gives
\[
 \|N_i/M_i-\nu_{i,T}\|_1\le Ce^{-aT},
\]
for some \(a>0\), with densities interpreted with respect to \(dV_g\). Both the model densities and \(N_i/M_i\) have exponential supremum bounds. Lemma~\ref{hyb:covariance} therefore applies; replacing \(M\) by \(M+\mu\) gives an error of \(O(T|\mu|)\), by the bounded \(L^{1+1/T}\) norms of the normalized densities.
\end{proof}
\begin{proposition}[Solutions of the projected equation]
\label{hyb:range}
There exists $\bar\delta>0$ such that the following hold for $0<\delta_0<\bar\delta$, $T$ sufficiently large and $|\mu|\le e^{-\delta_0T}$.
\begin{itemize}
\item[\textup{(i)}]
For fixed $C_0$ sufficiently large, there is a unique $\eta\in E_T$, $\|\eta\|_{\Hh}\le e^{-\delta_0T}T^{C_0}$, satisfying
\[
 P_{E_T}\mathcal F_\mu(B_T+\eta)=0,\qquad
 \mathcal F_\mu(B_T+\eta)=\sum_\alpha c_\alpha S_\alpha.
\]
Moreover,
\[
 \|\eta\|_{\Hh}+|c|\le C(\epsilon(T)+|\mu|)T^C.
\]
\item[\textup{(ii)}]
The family depends continuously on the parameters in the compact set and is differentiable for each finite $T$. Its projected second derivative has index \eqref{hyb:index}.
\end{itemize}
\end{proposition}
\begin{proof}
The difference of the projected operators is $o(T^{-1})$. Hence the norm of the inverse at the approximate solution is bounded by \(CT\).

To estimate the nonlinear remainder, write \(P_i\) for the probabilities at the approximate solution. Interpolating their mass and exponential supremum bounds gives
\[
 \|P_i\|_p\le Ce^{C(1-1/p)T},
\]
with \(p>1\) as close to one as needed. For \(u\) in the exponentially small ball, the denominator \(\int P_i e^{u_i}\) is bounded below, by Jensen and
\[
 |\int P_i u_i|\le C\sqrt T\|u_i\|_{H^1}.
\]
Use H\"older with this fixed \(p>1\) for each appearance of \(P_i\) in the derivatives of the normalized exponential. Use fixed finite Sobolev exponents for the perturbations and the test function, and the two-dimensional exponential Sobolev inequality for \(e^{u_i}\).

The second derivative contains at most three weighted integrals, so choosing \(p\) sufficiently close to one bounds the norm of this bilinear map for the energy norms by \(C_\varepsilon e^{\varepsilon T}\). The integral remainder estimate then gives
\[
 \|\mathcal N(\eta)-\mathcal N(\xi)\|
 \le C_\varepsilon e^{\varepsilon T}
        (\|\eta\|+\|\xi\|)\|\eta-\xi\|.
\]
Take $\bar\delta$ below the positive decay rate of $\epsilon(T)$, and choose
\[
 \varepsilon<\delta_0/2.
\]
The inverse residual is $O((\epsilon(T)+|\mu|)T^C)$, strictly below half the chosen radius for sufficiently large $C_0$. The Lipschitz constant after inversion tends to zero. The fixed-point estimate gives the bound on $\eta$.

Projecting onto the sources $S_\alpha$ gives the bound on $c$, since the matrix $(\langle S_\alpha,S_\beta\rangle_{\Hh})$ has a bounded inverse. The change in the second derivative is $o(T^{-1})$, so its index is unchanged.

The implicit function theorem gives continuous dependence on the parameters and differentiability for each finite $T$.
\end{proof}
\subsection{The balance and shape equations}
We first estimate the change in the local identities when the approximate solution is corrected. If $\nu$ is a probability measure of the approximate solution with $\|\nu\|_\infty\le e^{CT}$, $E=\|\eta\|_{H^1}$ and $\sqrt T E\to0$, then for $q=T$
\begin{equation}\label{hyb:measure-stability}
 Z=\int\nu e^\eta=1+o(1),\qquad
 \|e^\eta/Z-1\|_{L^q(\nu)}
       \le C\sqrt T E e^{CTE^2}=o(1).
\end{equation}
Indeed H\"older gives
\[
 \|e^\eta-1\|_{L^q(\nu)}\le\|\nu\|_2^{1/q}\|e^\eta-1\|_{2q},
\]
where $\|\nu\|_2^{1/q}\le e^{CT/(2q)}$ and
\[
 \|e^\eta-1\|_{2q}\le C\sqrt q E e^{CqE^2}
\]
by exponential Sobolev and H\"older. On a fixed matching annulus the trial densities have uniformly bounded derivatives, and the weighted estimate implies ordinary $L^p$ convergence of the density difference for every fixed $p<q$: for a trial density $v$, one uses
\[
 v^p\leq \|v\|_{L^\infty({\rm annulus})}^{p-1}v
\]
there. Subtracting the trial equation from the projected equation adds only the explicit source tails and the trial residual on that annulus.

The local profiles and their derivatives converge there to the exterior fields and their derivatives, and the exterior fields solve the limiting equation off the centers.

Thus the residual and the source tails tend to zero in the norms of smooth functions on fixed smaller annuli. Interior $W^{2,p}$ estimates, using
\[
 \|\eta\|_{L^p}\leq C_p E
\]
for the zero-order term, give $W^{2,p}$ convergence of $\eta$ on these smaller annuli. With $p>2$ this gives $C^{1,\alpha}$ convergence on a still smaller annulus. Substitution into the density difference in the equation gives $C^{1,\alpha'}$ convergence of its right side for
\[
 0<\alpha'<\alpha,
\]
so the interior Schauder estimate gives $C^{3,\alpha'}$ convergence of $\eta$ there. Define
\[
 \mathcal M(B)=\int_\C\frac{|w+B|}{(1+|w|^2)^2}\,dA_w,
 \qquad \Phi_b(B)=\kappa_b\nabla\mathcal M(B)+4\pi t_b^{\rm sh},
 \quad t_b^{\rm sh}=(\Re T_b,-\Imm T_b).
\]
Here $t_b^{\rm sh}$ denotes the vector of harmonic quadratic coefficients; the inner profile exponent remains $t_b$.
\begin{lemma}[Bounded Cauchy kernels]
\label{hyb:cauchy}
Let $f\in C^2(\C)$ satisfy, for some $\delta>0$,
\begin{equation}\label{hyb:cauchy-hypotheses}
 |D^\alpha f(z)|\le A_\alpha(1+|z|)^{-2-\delta-|\alpha|},\qquad |\alpha|\le2.
\end{equation}
Put $C_f(z)=\pi^{-1}\int f(y)/(z-y)\,dA_y$.
\begin{itemize}
\item[\textup{(i)}]
Pairing two Cauchy potentials against $f$ gives the kernel
\begin{equation}\label{hyb:divided-difference}
 K_f(z,w)=-\pi\frac{C_f(z)-C_f(w)}{z-w}.
\end{equation}
The only possible directional discontinuity at $z=w$ comes from
\[ -\pi f(w)(\overline z-\overline w)/(z-w). \]
For $f=\partial_{B_j}m_B$, this term vanishes at $(0,0)$.
\item[\textup{(ii)}]
Uniformly for bounded constants in \eqref{hyb:cauchy-hypotheses},
\[
 |K_f(z,w)|\le C\min\{(1+|z|)^{-1},(1+|w|)^{-1}\}.
\]
This includes $f_j(y)=16y_j/(1+|y|^2)^3$ and the tests for the shape parameters, uniformly for bounded $B$.
\end{itemize}
\end{lemma}
\begin{proof}
\smallskip\noindent\textit{Step 1. The Cauchy integral.}\quad The decay of $f$ implies absolute convergence of the Cauchy integral. Partial fractions give \eqref{hyb:divided-difference}.

In the principal-value formula for $(C_f)_z$, subtract $f(z)$ on a disk centered at $z$; its constant angular integral is zero. On the unit disk the remaining integral is bounded by $C\|\nabla f\|_\infty$, and on its complement by $C\|f\|_1$. Together with
\[
 \bar\partial C_f=f,
\]
this proves a uniform bound for the real first derivatives of $C_f$; their continuity follows from the same subtraction and the $C^2$ bound. The real Taylor formula now gives the diagonal expression and boundedness of $K_f$.

For $|z|=R\ge2$, split the Cauchy integral at $|y|=R/2$ and $|y-z|=R/2$. The decay gives
\[
 |C_f(z)|\le C/R.
\]
In its differentiated integral, subtract $f(z)$ on $|y-z|<R/2$. The local contribution is at most $CR\sup_{|y-z|<R/2}|\nabla f(y)|=O(R^{-2-\delta})$; the remaining pieces are $O(R^{-2})$. Thus $|DC_f(z)|\le C/R^2$. If $|z|\ge|w|$ and
\[
 |z-w|\ge|z|/2,
\]
the uniform bound for $C_f$ bounds the divided difference by $C/(1+|z|)$. Otherwise the joining segment stays outside $|z|/2$, and the derivative bound gives the stronger $C/(1+|z|)^2$. This proves the displayed decay of $K_f$.

\smallskip\noindent\textit{Step 2. The test functions.}\quad The rational test functions used here satisfy, for $|\alpha|\le2$, $D^\alpha f_j=O((1+|y|)^{-5-|\alpha|})$ and
\[
 D^\alpha\partial_{B_j}m_B=O((1+|y|)^{-8-|\alpha|}),
\]
uniformly on bounded shape sets. They are smooth at finite points, and the latter vanishes at zero because it contains $|y|^2$. A radial cut-off preserves \eqref{hyb:cauchy-hypotheses} uniformly in the rescaled variables.

Exact centering of a shape test can add a bounded multiple of the cut-off $m_B$, whose derivatives decay as $O((1+|y|)^{-6-|\alpha|})$. Thus the centered family still satisfies \eqref{hyb:cauchy-hypotheses} with the common choice
\[
 \delta=3,
\]
and still vanishes at zero. The odd translation tests have zero mean under the radial cut-offs.

We need the same bound for smooth test functions supported on shrinking disks. Let $k_d(z)=\chi(z/d)k(z)$, where $\chi$ is a fixed radial smooth cut-off and $k$ has a bounded $C^2$ norm near zero. Write
\[
 k_d(z)=k(0)\chi(z/d)+d h_d(z/d);
\]
the functions $h_d$ have uniformly bounded $C^2$ norms and fixed compact support. Rescaling the Cauchy integral gives
\[
 K_{k_d}(z,w)=k(0)K_\chi(z/d,w/d)
                  +dK_{h_d}(z/d,w/d).
\]
Hence these kernels too are uniformly bounded, although the unscaled derivatives of $k_d$ need not be. Radiality gives
\[
 (C_\chi)_z(0)=0,
\]
so $(C_{k_d})_z(0)=O(d)$. This is the uniform cut-off estimate used when taking the small-disk limit of the smooth cubic defect.
\end{proof}
\begin{proposition}[The balance and shape equations]
\label{hyb:stress-rows}
The following estimates hold in $C^0$, uniformly on the compact parameter set.
\begin{itemize}
\item[\textup{(i)}]
At the regular points,
\begin{align*}
 c_a^{\rm ctr}/R_a&=-3F_a+o(1)
       &&\text{at ordinary points},\\
 c_b^{\rm ctr}/r_b&=-3F_b+o(1)
       &&\text{at mixed regular points}.
\end{align*}
\item[\textup{(ii)}]
If $c^{\rm ctr}=0$, then, for all shape parameters under consideration,
\begin{equation}\label{hyb:shape-row}
 \frac{c_b^{\rm sh}}{R_b^2}=-\frac3{16\pi}\Phi_b(B_b)+o(1).
\end{equation}
\end{itemize}
\end{proposition}
\begin{proof}
\noindent\textit{Step 1. The two local identities.}\quad On a regular disk disjoint from the prescribed points, put $U_i=u_i+\log h_i+\varphi+\log(\rho_i/\int_Mh_ie^{u_i})$. The projected source is retained; write $q_i$ for its Euclidean density from \eqref{hyb:source-convention}. Set
\[
 a=(2U_1+U_2)/3,\qquad b=(U_1+2U_2)/3.
\]
Then $\Delta U+Ke^U=\kappa-Kq$ with $\kappa_i=\Delta(\log h_i+\varphi)+e^\varphi(K\rho)_i$. With primes denoting $\partial_z$, put
\[
 W_2=a''+b''-(a')^2+a'b'-(b')^2,
\quad
 W_3=a'''-2a'a''+a'b''+(a')^2b'-a'(b')^2.
\]
Differentiation cancels the exponential terms in $\bar\partial W_j$. On a cut-off equal to one on the source support, integration by parts gives the exact forcing contributions
\begin{align}
 E_2&=\tfrac14\int\{(2a'-b')q_1+(2b'-a')q_2\},\notag\\
 E_3&=-\tfrac14\int q_1\{b''+2a'b'-(b')^2\}
       +\tfrac14\int q_2\{a''-(a')^2+2a'b'\}.
 \label{hyb:stress-source}
\end{align}
The source cut-offs give no additional local source term; this uses their support and their exact local centering. To make the smooth terms explicit, put
\[
 k_1=(2\kappa_1+\kappa_2)/12,\qquad k_2=(\kappa_1+2\kappa_2)/12.
\]
Their contributions to $\bar\partial W_2$ and $\bar\partial W_3$ are
\begin{align*}
 (\bar\partial W_2)_{\rm sm}
 &=k_1'+k_2'-(2a'-b')k_1-(2b'-a')k_2,\\
 (\bar\partial W_3)_{\rm sm}
 &=k_1''-2a'k_1'+a'k_2'
   +k_1\{-2a''+b''+2a'b'-(b')^2\}
   +k_2\{(a')^2-2a'b'\}.
\end{align*}
Before integration, the full source-dependent identities are
\begin{align*}
 \bar\partial W_2={}&(\bar\partial W_2)_{\rm sm}
 +\tfrac14\{(2a'-b')q_1+(2b'-a')q_2\}
 -\tfrac14\partial_z(q_1+q_2),\\
 \bar\partial W_3={}&(\bar\partial W_3)_{\rm sm}
 -\tfrac14q_1\{b''+2a'b'-(b')^2\}
 +\tfrac14q_2\{a''-(a')^2+2a'b'\}\\
 &+\partial_z\{-q_1'/4+a'q_1/2-a'q_2/4\}.
\end{align*}
The total $z$ derivatives have only boundary terms, with the multiplier sign fixed by \eqref{hyb:source-convention}. For clarity, one may first fix a small circle of radius $d$ and integrate inside it, then take the concentration limit and send $d$ to zero.

Although the sources have a fixed cut-off, they and their first two spatial derivatives on the circle are $O_d(|c|(R^3+r^3+R^4))$ at a mixed disk, with the $R^3$ term sufficient for ordinary disks. These bounds follow directly from the rational source densities and the fixed smooth cut-offs.

The exterior derivatives of the fields are bounded on each such circle by \eqref{hyb:measure-stability} and the annular elliptic estimates. Thus boundary terms from the total derivatives are $o(1)$ for fixed $d$.

Extending the integrated source pairings from that disk to the full local source supports adds another $o(1)$, since the omitted annulus is a positive distance from the bubble and all its fields have bounded derivatives there.

This justifies using \eqref{hyb:stress-source} up to $o(1)$ in the limiting identities.

\noindent\textit{Step 2. Contributions of the smooth terms.}\quad We next compute the smooth terms in the identity for $W_3$, including the quadratic term involving two inner measures. The local Cauchy representations for the projected equation use the signed measures $f_i+q_i$, since
\[
 -\Delta a=f_1+q_1-4k_1,\qquad -\Delta b=f_2+q_2-4k_2.
\]
Their source total variations are
\[
 O(|c|)=o(1).
\]
For a smooth cut-off $k$, the quadratic gradient product is given by \eqref{hyb:divided-difference}. Its uniform boundedness makes that product stable under a total-variation $o(1)$ replacement.

Thus \eqref{hyb:measure-stability} and the total-variation comparison in Proposition~\ref{hyb:trial-comparison} allow us to use the measures in \eqref{hyb:model-probabilities}, each multiplied by its total mass.

On the fixed disk we include the residual measure. It will be removed only in the subsequent small-disk limit. This replacement determines the diagonal angular integral before the inner measures concentrate at their common point.

Identify a real vector $(v_1,v_2)$ with $v^{\mathbb C}=v_1+iv_2$, and put $g_B=(\partial_{B_1}\mathcal M)+i(\partial_{B_2}\mathcal M)$. The complex derivative $\partial_z$ therefore represents real gradients by their complex conjugates. In the notation
\[
 \Theta(z)=\overline z/z,
\]
let $S_{ij}$ be the limiting double angular integrals of the two concentrating measures. The inner radial model gives $S_{22}=0$ by simultaneous rotation. The outer model and $w=y^2-B$ give the explicit identity
\[
 S_{12}=4\pi\int_{\mathbb C}m_B(y)\frac{\overline y}{y}\,dA_y
 =32\pi\int_{\mathbb C}
 \frac{\overline{w+B}}{|w+B|(1+|w|^2)^2}\,dA_w
 =32\pi\overline{g_B}.
\]
To compute the pairing of the first density with itself, write $A_B=-\partial_y\log(1+|y^2-B|^2)$. It is the Cauchy potential
\[
 -(4\pi)^{-1}\int m_B(x)/(y-x)\,dA_x
\]
: the two expressions have $\bar\partial A_B=-m_B/4$ and decay at infinity, so their holomorphic difference is zero. Symmetrizing in $x,y$ gives
\[
 S_{11}=-8\pi\int_{\mathbb C}\bar y\,m_B(y)A_B(y)\,dA_y.
\]
The identity $\partial_y m_B=m_B(1/y+2A_B)$ and $\partial_y\bar y=0$ imply $2\int\bar y m_B A_B=-\int(\bar y/y)m_B$. Indeed the boundary terms vanish since $m_B=O(|y|^2)$ at zero and $O(|y|^{-6})$ at infinity. Hence
\[
 S_{11}=4\pi\int(\bar y/y)m_B=S_{12}.
\]
The smooth terms with at most one gradient vanish when the testing disk shrinks, by uniform local $L^p$ bounds for $p<2$. For the second derivatives, use a radial cut-off $\chi(z/d)$ equal to one near zero.

On its transition annulus take the concentration limit first. The monopole part paired with $\partial_z\chi(z/d)$ integrates to zero, since $\partial_z\chi(z/d)/z$ has angular frequency $-2$.

The regular parts give $O(d)$ there. Replacing $k_i(z)$ by $k_i(0)$ in this calculation leaves terms bounded by the single-gradient estimate, also $O(d)$.

For the quadratic products, the Cauchy transform of a radial cut-off is $\bar z$ near zero. Indeed it is rotation-equivariant of frequency $-1$, and after subtracting $\bar z$ it is holomorphic in the disk where the cut-off is one. Hence its divided difference there is $-\pi\Theta(z-w)$.

The $O(d)$ change of $k_i$ on the testing disk changes the corresponding scaled Cauchy kernel by $O(d)$, by the bounds of Lemma~\ref{hyb:cauchy}. Terms with residual mass vanish as the testing disk shrinks.

If both variables are in a fixed disk of radius $e$, the residual mass there is $O(e^2)$ and the double kernel is bounded.

If either variable lies outside that disk, the scaled kernel is bounded by $Cd/e$ by its decay. First let $d$ tend to zero with $e$ fixed, and then let $e$ tend to zero. The remaining quadratic products are the angular integrals above: their sum is
\[
 -\frac{k_1(2S_{12}-S_{22})+k_2(S_{11}-2S_{12})}{16\pi}
       =-\frac{\kappa_1}{2}\overline{g_B},
\]
because $2k_1-k_2=\kappa_1/4$. The canonical forcing $\kappa_1$ at the bubble tends to the trace in \eqref{hyb:canonical-trace}. Thus the quadratic term is determined by the canonical forcing of the normalized projected equation.

\noindent\textit{Step 3. The off-diagonal coefficients.}\quad We estimate the terms involving both the translations and the shape parameters before setting the center multipliers equal to zero. Put $\lambda=r/R$ and
\[
 w_r=(1+|z|/r)^{-1}.
\]
From the bounds $n_1\le C(1+|z|^2/R^4)$ on $|z|\le R$ and $n_1\le CR^4|z|^{-6}$ on $R<|z|\le d_0$ in the fixed matching disk $D_{d_0}$, direct integration gives
\begin{equation}\label{hyb:translation-weight}
 \int_{D_{d_0}} w_r^p n_1\le C\lambda^p\qquad(1\le p\le3/2).
\end{equation}
Here are the bounds for this cubic coupling, including the locally projected sources. Rescale by the inner width \(r\) and put
\[
 W(y)=(1+|y|)^{-1},
\]
with
\[
 \lambda=r/R
\]
as above. Let \(\mu_i\) be the pushforwards of the local measures \((f_i+q_i)\,dA\) to this scale, and let \(\mu_i^0\) be the pushforwards of the even measures $n_i\,dA$ from the explicit mixed profile on the same centered disk.

The part of the fields supplied by the smooth forcing and by the exterior is \(O(r)\) after multiplying the first derivatives by \(r\). For the first measures,
\[
 \int W\,d|\mu_1-\mu_1^0|=o(\lambda),\qquad
 \int W\,d|\mu_1^0|=O(\lambda).
\]
For the projected correction, apply \eqref{hyb:measure-stability} with $q=T$ and H\"older with conjugate $T/(T-1)\le3/2$ in \eqref{hyb:translation-weight}. The approximate first density in the core is $n_1 e^{J_1}$ times a normalization $1+o(1)$, and $e^{J_1}-1=O(|z|^2)$ there. Since
\[
 W(z/r)|z|^2\le r|z|,\qquad \int |z|n_1=O(R),
\]
this term contributes $O(rR)$. On fixed cut-off transition annuli the bounded mass contributes $O(r)=o(\lambda)$ since
\[
 W(z/r)=O(r)
\]
there. Thus the disk used for the profile $(n_1,n_2)$ and the disk of integration in the Cauchy transform may have different fixed radii. The first projected source on this mixed disk is a linear combination of the outer shape sources.

Each absolute shape source is bounded on the outer scale by a multiple of $m_B$, including its centering term. Its contribution is therefore
\[
 O(|c|\lambda)=o(\lambda).
\]

For the second measure, the projected correction, normalization and inner projected source contribute $o(1)$ in total variation.

Discrepancies on fixed cut-off transition annuli have bounded variation and cost $O(r)$ by the second-variable weight of Lemma~\ref{hyb:cauchy}. In the core the remaining change is the smooth factor
\[
 e^{J_2}-1=O(|z|).
\]
In pairing this change with an outer first measure using the kernel of Lemma~\ref{hyb:cauchy}, split the second measure at $|z|=R$. On $|z|\le R$ take the weight from the first factor. The bound on $n_2-n_{\rm bg}$ in Proposition~\ref{hyb:trial-comparison}, radial integration and $|n_{\rm bg}|\le C$ give
\[
 \int_{|z|\le R}|z|\,n_2
 \le C R^3+C\int_0^R\frac{d^2 r^2}{(r^2+d^2)^2}\,{\rm d}d
 \le C(R^3+r).
\]
The last integral is bounded after setting $d=ru$ and using
\[
 \int_0^\infty u^2(1+u^2)^{-2}\,{\rm d}u<\infty.
\]
The weaker bound $O(R)$ also follows from bounded total mass. On $|z|>R$ take the weight from the second factor and use $W(z/r)|z|\le r$ and the bounded mass. In both regions the pairings are $o(\lambda)$. The same proof applies with the first measure $\mu_1$ in place of $\mu_1^0$.

In the rescaled cubic coefficient for an inner translation, the differentiated linear term can be integrated by parts against the translation density \(h(y)\). Its kernel is a first derivative of \(C_h\) and is bounded by \(C(1+|y|)^{-2}\) by the proof of Lemma~\ref{hyb:cauchy}.

Its change on replacing \(\mu_1\) by \(\mu_1^0\) is therefore \(o(\lambda)\). The two quadratic terms use the bounded kernels \(K_h(y_1,y_2)\), bounded in absolute value by \(C\min\{W(y_1),W(y_2)\}\).

Expand their differences one factor at a time, taking the weight on the changed first factor, or on the unchanged first factor when the second factor carries an \(o(1)\) total-variation error.

For the change from \(e^{J_2}-1\), use the preceding split. This gives an \(o(\lambda)\) change for the quadratic terms as well; contributions involving the smooth fields are
\[
 O(r)=o(\lambda).
\]
The even model densities have odd first fields and even second derivatives, whereas the translation source is odd. Their cubic pairing is zero. Consequently \(r^2L_{3,j}=o(\lambda)\), or
\[
 rR L_{3,j}=o(1),
\]
as required. For the other off-diagonal coefficient, the outer shape source contributes $O(R^{-1})$ to $E_2$.

After rescaling by $R$, each first derivative is paired with the bounded single Cauchy transform of a shape test density, against measures of bounded mass. This is the \(O(1)\) coefficient of \(Y\) in the first row below.

\noindent\textit{Step 4. The diagonal coefficients.}\quad We compute the diagonal entries in the two equations before solving them. In any scaled disk, write a first field as $P[\mu]+H$, where
\[
 P[\mu](y)=-\frac1{4\pi}\int\frac{d\mu(x)}{y-x}
\]
and $H$ is the scaled contribution from the smooth forcing and the exterior. To see the stated smooth-field bounds, represent the field in the original variables on a disk using its logarithmic potential plus the smooth forcing potential and a harmonic remainder.

The boundary oscillation of the remainder is bounded after subtraction of a constant, by taking the representation on a larger fixed matching disk.

Its gradient and its derivatives are then bounded on the smaller disk. In particular $H=O(r)$ on the inner scale and $O(R)$ on the outer scale.

For a scaled source test $g$, the linear, differentiated-linear, and quadratic tests of the singular parts are
\[
 \int gP[\mu]=\tfrac14\int C_g\,d\mu,\qquad
 \int g\partial_yP[\mu]=-\tfrac14\int C_{\partial_y g}\,d\mu,\qquad
 \int gP[\mu]P[\nu]=\tfrac1{16\pi^2}\iint K_g\,d\mu\,d\nu .
\]
The last formula follows by first integrating in $y$ and using \eqref{hyb:divided-difference}; the differentiated formula follows by integration by parts. The cut-offs converge in these tests. On each bounded rescaled disk the test and its first two derivatives converge to those of the limiting test.

The uniform rational tail bounds of Lemma~\ref{hyb:cauchy} give convergence in $L^1$ and $C^1$ of the tests and in $L^1$ and $C^0$ of their first derivatives. The same bounds also apply along convergent sequences of the shape parameter.

The centering of a shape source changes its coefficient by $o(1)$ since the cut-off mass tends to the constant plane mass of $m_B$, and the cut-off integral of $\partial_{B_j}m_B$ tends to zero.

The principal-value subtraction in Lemma~\ref{hyb:cauchy} then gives convergence of the single and differentiated single kernels on compact sets, and of the double kernels uniformly off the diagonal.

The single and double kernels tend to zero outside large compact sets uniformly, so mass escaping to infinity makes no contribution. For the double kernels, first discard the projected source measures whose total variations tend to zero.

On compact sets the products of the remaining nonnegative measures converge weakly, as one sees by approximation with finite sums of separated continuous tests.

A closed neighborhood of any discontinuities of the limiting kernel can be omitted with arbitrarily small error provided the limiting product gives those discontinuities zero mass.

The only possible obstructions are common point masses, where continuity of the test kernel is checked below.

On the inner scale the first measure tends vaguely to zero and the second to $m_{1,1}(y)\,dA_y$. On the outer scale they tend vaguely to $m_B(y)\,dA_y$ and $4\pi\delta_0$, respectively.

These are the limits of $n_1\,dA$ and $n_2\,dA$. The projected measures have the same limits by \eqref{hyb:measure-stability}, the smooth multipliers $e^{J_i}$, and
\[
 |c|=o(1).
\]
For example the second inner core has mass $4\pi$ and
\[
 r/R\to0,
\]
and the part of the second density outside a large multiple of $r$ and inside a small fixed disk has mass bounded by the decaying inner tail and the bounded residual density times the disk area. A single test of
\[
 P[m_{1,1}\,dA]=-\bar y/(1+|y|^2)
\]
in the inner entry of $E_2$, with the translation source temporarily multiplied by four, gives
\[
 -8\int_{\mathbb C}\frac{\bar y\,y_j}{(1+|y|^2)^4}\,dA_y
 =-\frac{2\pi}{3}\begin{cases}1& j=1,\\-i&j=2.\end{cases}
\]
The ordinary translation has the same coefficient at its own width: there the first, instead of the second, measure is the scalar core and the test uses $(2a'-b')/4$. For the outer cubic entry take
\[
 g=\partial_{B_j}m_B.
\]
The part $b''-(b')^2$ has zero limit exactly in the tested sense. In fact, for $4\pi\delta_0$ the differentiated kernel gives $-\pi\partial_y C_g(0)$, while the quadratic kernel gives
\[
 K_g(0,0)=-\pi\partial_yC_g(0).
\]
The divided difference is continuous at $(0,0)$ because $g(0)=0$. For the cross field the limiting product of measures gives no mass to the diagonal. The two limiting fields are
\[
 A_B=-2y\overline{y^2-B}/(1+|y^2-B|^2)
\]
and $-1/y$, so $2A_B(-1/y)=4\bar w/(1+|w|^2)$ with $w=y^2-B$. Since
\[
 \partial_{B_j}m_B=4m_Bw_j/(1+|w|^2),
\]
the entry in $E_3$ after multiplying by $R^2$ is
\[
 -32\int_{\mathbb C}\frac{\bar w\,w_j}{(1+|w|^2)^4}\,dA_w
 =-\frac{8\pi}{3}\begin{cases}1&j=1,\\-i&j=2.\end{cases}
\]
Here the pushforward of $m_B\,dA_y$ is $8(1+|w|^2)^{-2}\,dA_w$. All limits are uniform on compact shape and separated-center sets, by the stated uniform kernel bounds and compactness of these parameters.

\noindent\textit{Step 5. The two limits.}\quad Write
\[
 \widehat c^{\rm ctr}=c^{\rm ctr}/4
\]
for the coefficient of the fourfold translation sources, and set $X=\widehat c^{\rm ctr}/r$, $Y=c^{\rm sh}/R$. At a mixed regular point the two equations are
\[
 (M_{11}+o(1))X+O(1)Y=\mathcal R+o(1),\qquad
 o(1)X+(M_{22}+o(1))Y=O(R),
\]
where conjugation identifies the complex equations with real equations, $M_{11}=-2\pi/3$ and
\[
 M_{22}=-8\pi/3.
\]
In the second row the fixed matching circle bounds the boundary term; the smooth defect is bounded by integration by parts for second derivatives and the bounded Cauchy kernels for quadratic products. This row is multiplied by $R$.

At an ordinary point use the first row alone with width $R_a$ in place of $r$.

Taking the concentration limit first on a fixed circle and then shrinking it, the integral of the smooth contribution to $\bar\partial W_2$ tends to zero by the $O(d)$ estimate for the terms containing one first derivative. The exterior residue is
\[
 3A_0=\overline{F^{\mathbb C}}/2
\]
at a mixed point. At an ordinary point it is
\[
 2A_0-B_0=\overline{(F^A)^{\mathbb C}}/2,
\]
using $a'=-1/z+A_0+O(|z|)$ and $b'=B_0+O(|z|)$ in $W_2$. Integration over the disk therefore gives
\[
 \mathcal R=\pi\overline{F^{\mathbb C}}/2
\]
for the mixed row, with the same coefficient for the ordinary row. Inverting the two rows first gives $X=O(1)$ and $Y=o(1)$. The first row then gives $X=-3F/4+o(1)$ in real coordinates.

Returning to the original sources multiplies this coefficient by four. When $c^{\rm ctr}=0$, the center equations give
\[
 F=o(1).
\]
For the cubic boundary integral take the concentration limit on a fixed small circle first, and then shrink the circle. In the exterior fields $a'=-2/z+A_0$ and
\[
 b'=-1/z+B_0,
\]
the holomorphic Taylor terms give simple-pole coefficient $4A_0'-2B_0'-2A_0B_0-A_0^2+2B_0^2$ in $W_3$, evaluated at zero. The higher poles cancel identically. The terms involving first-order $\bar z$ terms at order $|z|^{-1}$ have zero circular integral. The center equation gives
\[
 A_0(0)=o(1).
\]
The regular gauge in \eqref{hyb:canonical-trace} also gives
\[
 2A_0(0)-B_0(0)=o(1),
\]
because this combination is the first derivative of the smooth first component, which has zero linear term in that coordinate in the limit.

Thus both constant terms in the derivatives are $o(1)$, and the remaining coefficient from the quadratic Taylor term of the entire exterior field is $2T_b+o(1)$.

Including the $O(d)$ remainder on a circle of radius $d$, the iterated Stokes limit is $2\pi T_b+o(1)$. The smooth defect in that same limit is $-\kappa_b\overline{g_{B_b}}/2+o(1)$. The outer cubic entry just computed applies here with
\[
 c^{\rm ctr}=0.
\]
The residue and the smooth defect therefore give
\[
 2\pi T_b=-\frac{\kappa_b}{2}\overline{g_{B_b}}
       -\frac{8\pi}{3R_b^2}\overline{(c_b^{\rm sh})^{\mathbb C}}
       +o(1)+o(1)|c_b^{\rm sh}|/R_b^2.
\]
To obtain the stated error uniformly, one may first fix a sufficiently small matching radius. The bounded boundary and smooth contributions, and the nonzero limiting outer coefficient $-8\pi/3$, give a uniform bound on $|c_b^{\rm sh}|/R_b^2$ there, since the center sources vanish identically.

With this bound the coefficient errors tend to zero. One can then send the matching radius to zero in the remaining residue and smooth defect. Conjugation gives
\[
 \frac{(c_b^{\rm sh})^{\mathbb C}}{R_b^2}
 =-\frac{3}{16\pi}\left(4\pi\overline{T_b}+\kappa_b g_{B_b}\right)+o(1).
\]
In real coordinates this is \eqref{hyb:shape-row}. The convergence obtained here is uniform in $C^0$. This is the convergence required for the degree homotopy below.
\end{proof}
\subsection{Description of blowing-up solutions}
Take disjoint coordinate disks about the limiting blow-up points. On each disk the masses and boundary oscillations are uniformly bounded. Removing the smooth forcing gives the local system with uniformly positive smooth coefficients.

The canonical version has only uniformly bounded smooth forcing. Since the disk isolates the singular point, the estimate below on disks of small mass gives continuity of the regular part of each solution there.

Since the source strengths satisfy \eqref{fnd:independence}, Lemma~\ref{fnd:local-unitary} gives the centered inverse-square density bound. This bound and the prescribed local masses determine the two scales below.
\begin{lemma}[Profiles at a centered singular point]
\label{hyb:centered-singular-profiles}
Let $A,B>0$ with $A,B,A+B\notin\Z$. On a fixed coordinate disk $D$, let $f_i=e^{U_i}$ satisfy
\[
 \Delta U_1+2f_1-f_2=\kappa_1+4\pi(A-1)\delta_0,\qquad
 \Delta U_2-f_1+2f_2=\kappa_2+4\pi(B-1)\delta_0.
\]
Assume, for $i=1,2$, that $\kappa_i$ is smooth and
\[
 \|\kappa_i\|_{C^1(D)}+\int_D f_i+
 \operatorname{osc}_{\partial D}U_i\le C,\qquad |z|^2 f_i(z)\le C.
\]
For the following local masses at zero, in units of $4\pi$, the profiles are as follows.
\begin{itemize}
\item[\textup{(i)}]
For $(A,0)$ there is one profile $m_{A,R}$ in the first component.
\item[\textup{(ii)}]
For $(A+B,B)$ the second component has the inner profile $m_{B,r}$, and the first has the outer profile $m_{s,R}$, where $s=A+B$ and $r/R\to0$.
\end{itemize}
In both cases we obtain the stated profiles after adjusting the radii by fixed factors. The profiles account for the full local masses.
\end{lemma}
\begin{proof}
\noindent\textit{Step 1. The first scale.}\quad Write
\[
 a_1=A,\qquad a_2=B,
\]
and remove $2(a_i-1)\log|z|$ from $U_i$ to obtain $\widehat U_i$. In a disk about zero its distributional equation has only the two integrable densities and smooth forcing. Subtract their Dirichlet potentials, leaving a harmonic function smooth across zero.

For each fixed member of the sequence the density masses in a sufficiently small disk are arbitrarily small.

The upper bound for $\widehat U_i$ on a smaller disk uses only the positive term $2G_D f_i$ and a bounded harmonic and smooth contribution, since the other density has negative sign. Choose $p_0>1$ with
\[
 p_0(a_i-1)>-1,\qquad 1<p,q
\]
with $pq\le p_0$. The weight $|z|^{2(a_i-1)}$ then has finite $L^{pq}$ norm. Jensen applied to the logarithmic potential of the probability $f_i/\int_D f_i$ gives $e^{2G_D f_i}$ in $L^{pq'}$ by choosing the disk so that
\[
 pq'\int_D f_i/\pi<2.
\]
H\"older gives $f_i\in L^p$ there. The corresponding logarithmic potentials are continuous. Thus $\widehat U_i$ is continuous across zero and
\[
 f_i=O(|z|^{2a_i-2})
\]
for each fixed member, with its own constant. Here the solution is fixed. This argument uses neither uniform oscillation nor the centered density estimate.

The uniform estimates below use the assumed inverse-square bound. The largest of the two regular heights divided by $2a_i$ defines a scale $b$ by
\[
 -\log b=\max_i\frac{\max\widehat U_i}{2a_i},
\]
with maxima in a smaller fixed closed disk about the source. The positive local mass and the integrability of $|z|^{2a_i-2}$ imply
\[
 b\to0.
\]
In the coordinates $z=by$, the regular parts of $U_i(by)+2\log b$ are bounded above by zero. A maximizing point of one of these regular parts stays in a bounded scaled disk: if it is nonzero, the centered bound gives
\[
 |y|^{2a_i}\le C
\]
there. On every fixed scaled disk the two densities are bounded by $|y|^{2a_i-2}$, hence in a fixed $L^p$ with $p>1$ close to one. Subtract the bounded Dirichlet potentials of their regular equations.

The harmonic parts are bounded above, so Harnack comparison for positive harmonic functions shows that each regular field either converges locally, including at zero, or tends to $-\infty$ uniformly on compact sets. The component reaching the normalized maximum survives.

A diagonal subsequence gives these conclusions on the plane. The gradient representation follows by rescaling the fixed-disk Green formula; the regular gradient is $o(1)$ and mass outside a rescaled radius $L$ contributes only $O(1/L)$ on each fixed compact set.

The regular parts of the surviving fields are continuous at zero by the fixed $L^p$ estimate just proved.

By Fatou on fixed rescaled disks their plane masses are bounded above by the original local masses: each such rescaled disk lies in any fixed disk in the original variables for large $k$, which can then shrink.

If both survive, apply Lemma~\ref{fnd:entire-profiles} with shifted parameters $(A,B)$ to get mass $4\pi(A+B)$ in each component, exceeding the available mass of the second component in either case.

If just one survives, its scalar equation has shifted parameter $a_i$ and no additional point mass; the same lemma gives the mass $4\pi a_i$ and the radial profile because $a_i$ is noninteger. For $(A,0)$ the first component must survive and exhausts the mass. We now consider $(A+B,B)$.

\noindent\textit{Step 2. The second scale.}\quad We now determine the second scale. Denote the initially active component by $i$ and the other by $j$. Put
\[
 m_i(d)=\frac1{4\pi}\int_{D_d}f_i,\qquad
 w_i(d)=\overline{U_i}(d)+2\log d .
\]
The angular mean is denoted by an overline. By taking a diagonal sequence of slowly growing fixed multiples of $b$, there is a radius $d_-\to0$ with $d_-/b\to\infty$ such that
\[
 m_i(d_-)=a_i+o(1),
\]
\[
 m_j(d_-)=o(1),
\]
and $w_i(d_-),w_j(d_-)$ both tend to $-\infty$. This uses the fast tail of the entire scalar profile and the vanishing other component on every compact set of the first rescaled disk.

The centered bound gives uniformly bounded angular oscillations at all sufficiently small radii. Indeed the local Green representation of the gradient, split at distance $d/2$ from the observation point on $|z|=d$, has a nearby term bounded by
\[
 C d^{-2}\int_{|x|<d/2}|x|^{-1}\,dA
\]
and a remaining term bounded by $C/d$ times the total mass. The source logarithm and the bounded harmonic and smooth-forcing terms have the same bound after multiplication by $d$.

Consequently $dm_k/d\log d$ is comparable, with fixed positive constants, to $e^{w_k(d)}$. Exact circular integration gives
\[
 \frac{dw_i}{d\log d}=2(a_i-2m_i+m_j)+O(d^2),\qquad
 \frac{dw_j}{d\log d}=2(a_j+m_i-2m_j)+O(d^2).
\]
Choose a small fixed $\beta>0$ less than one fourth of each of $a_i,a_j$ and less than the terminal mass of component $j$. Let $b_+$ be the first radius after $d_-$ with
\[
 m_j(b_+)=\beta.
\]
The positive local mass of component $j$ at zero makes this radius exist and tend to zero. On $[d_-,b_+]$ the first derivative is at most $-c$ and the second at least $c$, for some $c>0$. The former bounds $m_i(b_+)-m_i(d_-)$ by
\[
 C e^{w_i(d_-)}=o(1).
\]
The latter bounds $m_j(b_+)-m_j(d_-)$ by $C e^{w_j(b_+)}$, giving a lower bound for $w_j(b_+)$. Its upper bound follows from the centered estimate. All derivatives are uniformly bounded. Since $w_j(d_-)\to-\infty$, it follows that $b_+/d_-\to\infty$.

In the coordinates $z=b_+y$, the first field tends to $-\infty$ uniformly on every compact annulus. Indeed
\[
 w_i(b_+)\to-\infty,
\]
the derivatives are bounded on each annulus with fixed ratio of radii, and the angular oscillations are uniformly bounded. Its measure tends to $4\pi a_i\delta_0$. For a possible point mass of the second measure at zero, integration of its positive derivative up to $b_+$ gives
\[
 m_j(\eta b_+)\le m_j(d_-)+C\eta^c \qquad (0<\eta<1,\ d_-<\eta b_+).
\]
Thus the limiting measure assigns zero mass to the origin after the second rescaling. The centered density bound and the bounds at $|y|=1$ give elliptic compactness for the second field on every compact annulus.

Its limit satisfies the scalar equation with source strength $a_i+a_j-1$: the source logarithm has coefficient $2(a_j-1)$ before rescaling and the collapsed first mass adds $2a_i$ to the radial derivative of the circular mean.

Equivalently, convergence on annuli and the preceding bound on $m_j(\eta b_+)$ give
\[
 d\overline U_j'(d)=2(a_i+a_j-1)-\frac1\pi\int_{D_d}e^{U_j}
\]
for the limiting field. The angular oscillation bound and the mass bound $O(d^c)$ make $U_j-2(a_i+a_j-1)\log|y|$ bounded at the origin: its circular average converges, since the residual derivative is $O(d^c)$.

The punctured scalar equation therefore extends distributionally with just this prescribed logarithm.

One can see this by subtracting the continuous Dirichlet potential of its locally $L^p$ density and removing the bounded harmonic singularity. Its total mass is bounded by the lower limit of the masses in the sequence.

The rescaled Green representation leaves no entire harmonic gradient, by the same estimate for the translated kernel as for the first scale.

The continuous regular part at zero and the shifted parameter $a_i+a_j=A+B$ therefore allow the scalar assertion of Lemma~\ref{fnd:entire-profiles}, which gives mass $4\pi(a_i+a_j)$ and the radial profile.

If $i=1$, this exceeds the terminal second mass $4\pi B$, a contradiction. Thus $i=2$, and the two limiting scalar masses are exactly $4\pi B$ and $4\pi(A+B)$. Since the two masses exhaust what is available at zero, there is no unaccounted concentration there.
\end{proof}
\begin{lemma}[Quantitative scalar profiles at a singular point]
\label{hyb:singular-rate}
Assume the hypotheses of Lemma~\ref{hyb:centered-singular-profiles} and let $R$ be the fitted outer radius. In the mixed case suppose the second canonical field is uniformly bounded on a fixed circle about the source.
\begin{itemize}
\item[\textup{(i)}]
There are $c>0$ and $N<\infty$ such that
\[
 \|f_1-m_{s,R}\|_{L^1(D_{R^{1/2}})}\le CR^c(1+|\log R|)^N ,
\]
where $s=A$ in the pure case and $s=A+B$ in the mixed case. The first mass outside $D_{R^{1/2}}$, up to a fixed smaller disk, satisfies the same bound.
\item[\textup{(ii)}]
In the mixed case, for each fixed $1<\theta<s/B$, $\|f_2-m_{B,r}\|_{L^1(D_{R^\theta})}$ satisfies a bound of the same form. After cutting off the inner concentration, the remaining second density has an integrable tail.
\item[\textup{(iii)}]
For some fixed $p>1$, the cut-off first and inner second density errors are bounded in $L^p$ by a positive power of $R$. Their mean-zero Green potentials satisfy such a bound in $C^0\cap H^1$.
\end{itemize}
\end{lemma}
\begin{proof}
\noindent\textit{Step 1. Bounds between the two scales.}\quad All masses below are in units of $4\pi$. Set $d=|z|$ and write $w_i(d)=\overline{U_i}(d)+2\log d$. Circular integration gives
\[
  w_1'=2(A-2m_1+m_2)+O(d^2),\qquad
  w_2'=2(B+m_1-2m_2)+O(d^2),
\]
where a prime denotes the derivative with respect to $\log d$. The centered bound and the Green representation give bounded angular oscillations, as in Lemma~\ref{hyb:centered-singular-profiles}.

In the mixed case, first choose a large fixed factor $C_\varepsilon$ and a small fixed outer disk. By exhaustion and monotonicity, the masses on
\[
 C_\varepsilon r<d<R/C_\varepsilon
\]
differ from $(0,B)$ by at most $\varepsilon$, and those on $C_\varepsilon R<d<d_0$ differ from $(s,B)$ by at most $\varepsilon$. Intervals omitted near $r$ and $R$ have bounded length in $\log d$, and the derivatives are bounded there. Since
\[
 w_2(r)=O(1),
\]
the fixed-circle bound for $U_2$ gives
\[
 |2B\log r-2s\log R|
 \le C_\varepsilon+C\varepsilon(|\log r|+|\log R|).
\]
First taking the radius to zero and then the exponent loss to zero gives $r=R^{s/B+o(1)}$. Integrating the derivatives with $w_1(R)=O(1)$ and $w_2(r)=O(1)$ gives, with arbitrarily small fixed exponent loss (renamed $\varepsilon$),
\begin{align*}
 f_1&\le C R^{-2s+\varepsilon}d^{2s-2-\varepsilon},&
 f_2&\le C r^{2B-\varepsilon}d^{-2B-2+\varepsilon}
                    &&(C_\varepsilon r\le d\le R/C_\varepsilon),\\
 f_1&\le C R^{2s-\varepsilon}d^{-2s-2+\varepsilon},&
 f_2&\le C d^{2A-2-\varepsilon}
                    &&(C_\varepsilon R\le d\le d_0).
\end{align*}
The upper bounds at the adjoining annuli of fixed logarithmic length follow from the derivative and oscillation bounds. In the inner disk,
\[
 f_1\le C R^{-2s+\varepsilon}r^{2B-\varepsilon}d^{2A-2},
 \qquad f_2\le C r^{-2B}d^{2B-2}.
\]
For the first bound rescale by $r$ and subtract $2(A-1)\log|y|$. The two rescaled densities on each fixed inner disk are bounded by
\[
 C|y|^{2(A-1)},\qquad C|y|^{2(B-1)},
\]
respectively, by the first rescaling in the proof of Lemma~\ref{hyb:centered-singular-profiles}. Subtract their uniformly bounded Dirichlet potentials in some fixed $W^{2,p_0}$,
\[
 p_0>1.
\]
The harmonic maximum principle then propagates the improved regular-part upper bound from $d=C_\varepsilon r$ to the interior. The second bound comes directly from the same choice of the first scale.

\noindent\textit{Step 2. Oscillation on the comparison circles.}\quad We next estimate the angular oscillation. On $d=\tau=R^\theta$ the first Cauchy field has
\[
 \tau|\mathcal C[f_1]|
 \le C(\tau/R)^{\eta_s},\qquad 0<\eta_s<\min(1,2s-\varepsilon).
\]
Indeed the inner disk $|y|<\tau/2$ and the adjoining annulus have bound $C(\tau/R)^{2s-\varepsilon}$. On the adjoining annulus use
\[
 \int_{|z-y|<C\tau}|z-y|^{-1}\,dA_y\le C\tau.
\]
The remaining middle annuli give $C\tau R^{-2s+\varepsilon}\int_{2\tau}^R d^{2s-2-\varepsilon}\,dd$, and the portion beyond $R$ gives $C\tau/R$. For the second density, subtract its actual mass inside $\tau/2$ divided by $\pi z$. For $|z|=\tau$ and $|y|<\tau/2$,
\[
 \tau\left|(z-y)^{-1}-z^{-1}\right|\le C(|y|/\tau)^{\eta_B},
 \qquad 0<\eta_B<\min(1,2B-\varepsilon).
\]
The core bound on $d\le C_\varepsilon r$ gives $\int |y|^{\eta_B}f_2\,dA=O(r^{\eta_B})$ there, and on the rest of $D_{\tau/2}$ the moment is bounded by
\[
 Cr^{2B-\varepsilon}\int_{C_\varepsilon r}^{\tau/2}
 d^{\eta_B-2B-1+\varepsilon}\,dd=O(r^{\eta_B}).
\]
On $\tau/2<d<2\tau$ the local kernel integrates to $O(\tau)$ and the density is bounded by $Cr^{2B-\varepsilon}\tau^{-2B-2+\varepsilon}$. On $2\tau<d<R/C_\varepsilon$ one uses
\[
 |z-y|^{-1}\le C/d.
\]
In the relative Cauchy bound both contributions are at most $C(r/\tau)^{2B-\varepsilon}$. For $d\ge R/C_\varepsilon$ the contribution is bounded by $C\tau\int_{R/C_\varepsilon}^{d_0} d^{2A-2-\varepsilon}\,dd$, hence by $C\tau^{\eta_A}$ for every
\[
 0<\eta_A<\min(1,2A-\varepsilon).
\]
On $[R/C_\varepsilon,C_\varepsilon R]$ the bound $f_2\le Cd^{2A-2-\varepsilon}$ follows from its outer endpoint by the uniform derivative and oscillation bounds. The powers in these annular estimates are positive because $r=R^{s/B+o(1)}$ and $1<\theta<s/B$ for $\tau=R^\theta$.

At
\[
 \sigma=R^{1/2}
\]
subtract both actual masses in $D_{\sigma/2}$ divided by $\pi z$. For the first density, the interior fractional moment is $O(R^{\eta_s})$, and the adjoining and exterior tail give $O((R/\sigma)^{2s-\varepsilon})$ in the relative Cauchy bound.

For the second density, split at $R/C_\varepsilon$. The inner fractional moment gives $O((r/\sigma)^{\eta_B})$. For its remainder, the relative kernel after subtracting the radial term due to the enclosed mass is bounded by $Cd/\sigma$ on
\[
 R/C_\varepsilon<d<\sigma/2.
\]
The resulting integral is at most $C\sigma^{2A-\varepsilon}$. The same bound holds on
\[
 \sigma/2<d<2\sigma
\]
by integrating the local kernel against the bound $C\sigma^{2A-2-\varepsilon}$. For $d>2\sigma$ the bound is $C\sigma\int_{2\sigma}^{d_0}d^{2A-2-\varepsilon}\,dd$. The three bounds are $O(\sigma^{\eta_A})$ with the above strict choice of $\eta_A$. The terms from the enclosed masses are radial and hence have zero angular oscillation.

The harmonic terms and the smooth forcing contribute $O(\tau)$ or $O(\sigma)$. Thus the boundary oscillations of the rescaled regular parts on these circles are bounded by positive powers of $R$. For each comparison below, write the fitted regular field as $W_q+z$ with
\[
 W_q=\log(4q^2/(1+|y|^{2q})^2),\qquad z(0)=0.
\]
For a fixed large $M$, local convergence makes the difference $z$ small on $D_M$. Outside $D_M$ the difference of the radial derivatives of the circular mean is bounded by the discrepancy of the masses on $D_M$, the two tail masses and the integral of the extra forcing.

The latter quantities have a bound tending to zero by first taking the sequence limit and then sending $M$ to infinity. Integration in the logarithmic radius, followed by division by
\[
 \log(2+|y|),
\]
preserves this bound. The bounded angular oscillation adds at most $C/\log(2+M)$. Thus
\[
 \|z/\log(2+|y|)\|_\infty=o(1).
\]
Use the origin weight $|y|^{2-2\kappa}$ and the exterior weight $(1+|y|)^{2+\zeta}$ of Lemma~\ref{hyb:scalar-inverse}, with $\kappa<q$ and $\zeta<\min(1,2q)$. The equation for $z$ has the stated forcing plus $2m_{q,1}(e^z-1-z)$, whose weighted norm is
\[
 O(\|z/\log(2+|y|)\|_\infty^2)
\]
using the model powers at zero and infinity. The scalar inverse absorbs this term because the logarithmic norm is already small. Thus the weighted forcing and the boundary oscillation give the stated logarithmic error in each case.

\noindent\textit{Step 3. Comparison at the inner and outer scales.}\quad On the inner scale the regular part of $U_2(ry)+2\log r$ solves the scalar equation with weight $\lvert y\rvert^{2(B-1)}$ and forcing $-r^2 f_1(ry)-r^2\kappa_2(ry)$.

Fit $r$ by a factor $1+o(1)$ to normalize its value at zero to $W_B(0)$ in the comparison above. In the weighted forcing norm take
\[
 0<\kappa<\min(1,A,B),\qquad \zeta>0
\]
small. For $L=\tau/r$ the forcing bound is
\[
 C\left\{(r/R)^{2s-\varepsilon}L^{2s-\varepsilon+\zeta}
             +r^2L^{2+\zeta}\right\}.
\]
Near zero this follows from the inner disk bound with power $|y|^{2A-2}$ and
\[
 \kappa<A.
\]
The two terms have exponents $(2s-\varepsilon)(\theta-1)-\zeta(s/B-\theta)+o(1)$ and $2\theta-\zeta(s/B-\theta)+o(1)$ in $R$. Both are positive for small $\zeta$. The tail bound gives uniform mass exhaustion on this expanding disk, and we have just proved the boundary oscillation estimate.

With the value at zero fixed, Lemma~\ref{hyb:scalar-inverse} gives the weighted logarithmic error $O(R^c)$ of the inner scalar profile. Since
\[
 \log(2+L)=O(|\log R|),
\]
the uniform error on this entire disk is $O(R^c(1+|\log R|))$. Thus $\|f_2-m_{B,r}\|_{L^1(D_\tau)}$ has the same bound, by integrating against the finite mass profile. The tail between $\tau$ and $\sigma$ contributes at most $C((r/\tau)^{2B-\varepsilon}+\sigma^{2A-\varepsilon})$, including the omitted model tail. Therefore the solution of
\[
 \Delta\Phi=f_2-m_{B,r}+\kappa_1
\]
with zero boundary data on $D_\sigma$ has a supremum norm bounded by a positive power of $R$. Indeed the density bounds provide some fixed $p_0>1$ for which the $L^{p_0}$ difference is at most $CR^{-K}$. Choose $p_0$ so that each origin power $d^{2a-2}$ in the bounds satisfies
\[
 (2a-2)p_0>-2.
\]
The relation $r=R^{s/B+o(1)}$ makes the remaining scale factors polynomial in $R^{-1}$. Write the preceding $L^1$ bound as $CR^{c_0}$, decreasing $c_0>0$ to absorb any logarithms. If
\[
 1/p=\lambda+(1-\lambda)/p_0,
\]
interpolation gives the exponent $\lambda c_0-(1-\lambda)K$. Fix $p>1$ sufficiently close to one to make it positive. The Dirichlet bound is $C\sigma^{2-2/p}$ times the $L^p$ error by dilation, and the smooth forcing contributes $O(\sigma^2)$.

For the outer scale set $a_0=r/R$ and
\[
 v(y)=U_1(Ry)+2\log R-2(A-1)\log|y|
             -\log(|y|^{2B}+a_0^{2B})-\Phi(Ry).
\]
Then the exact equation is
\[
 -\Delta v=2\bigl(|y|^{2s-2}+a_0^{2B}|y|^{2A-2}\bigr)
                 e^{\Phi(Ry)+v}.
\]
The density bound on a fixed scaled disk is
\[
 C(|y|^{2s-2-\varepsilon}+a_0^{2B-\varepsilon}|y|^{2A-2}).
\]
On $|y|<\delta$ its zero-boundary potential is bounded by $C\delta^{2s-\varepsilon}$ when $a_0\le\delta$, by radial integration. Thus convergence on compact annuli extends uniformly across zero by harmonic comparison, and $v$ is uniformly bounded on $|y|\le1$.

In fitting the outer radius, keep $\Phi$ and its physical disk fixed. Replacing $R$ by $\lambda R$ replaces $v(y)$ by $v(\lambda y)+2s\log\lambda$. The condition $v(0)=W_s(0)$ can therefore be enforced with
\[
 \lambda\to1
\]
and the boundary in the original variable remains at $\sigma$. For $|y|\ge1$, the upper bound on the outer scale gives
\[
 e^v\le C|y|^{-4s+\varepsilon}.
\]
The forcing relative to $2|y|^{2s-2}e^v$ has weighted norm at most $C(\|\Phi\|_\infty+a_0^{2B})$ for the above $\kappa$ and $\zeta<2s-\varepsilon$. The singular factor at zero is allowed since
\[
 \kappa<A,
\]
and the powers at infinity decay. Its $L^1$ bound and the upper bound on the outer scale give convergence of the full mass for $|y|^{2s-2}e^v$.

The oscillation on the circle in the original variable $\sigma$ is the one already estimated, since $\Phi$ vanishes there and the other subtracted terms are radial.

The same normalization and Lemma~\ref{hyb:scalar-inverse} give the logarithmic error for the outer profile. Integration against $m_{s,1}$ gives an $L^1$ error bounded by a positive power of $R$.

The additional term with coefficient $a_0^{2B}$ is bounded in $L^1$ by $Ca_0^{2B}$, using the uniform bound for $v$ near zero and its stated decay at infinity.

The first tail outside $\sigma$ is bounded by $C(R/\sigma)^{2s-\varepsilon}$. This proves the claimed estimate for the first mass. Only qualitative convergence of the masses was used in the linear estimate.

\noindent\textit{Step 4. The pure singular point.}\quad At a pure point, the masses on the annulus from a large fixed multiple of $R$ to a small fixed disk differ from $(A,0)$ by an arbitrarily small constant.

Integrating the derivatives from $R$ for $f_1$ and from a fixed exterior circle for the upper bound on $f_2$ gives the following estimates. On that fixed circle the upper bound for $f_2$ follows from the assumed
\[
 |z|^2 f_2\leq C:
\]
\[
 f_1\le C R^{2A-\varepsilon}d^{-2A-2+\varepsilon},\qquad
 f_2\le C d^{2(A+B)-2-\varepsilon}\quad(d\ge C_\varepsilon R).
\]
For $d\le C_\varepsilon R$, the choice of the first scale and harmonic comparison used above give
\[
 f_1\le C R^{-2A}d^{2A-2},\qquad f_2\le C R^{2A-\varepsilon}d^{2B-2}.
\]
At $\sigma=R^{1/2}$ the first fractional moment and the second estimate for the residual density yield a bound on the boundary oscillation by a positive power of $R$ by the same Cauchy integral estimates with exponents $A$ and $A+B$, respectively.

Fit the first radius at zero. On $|y|<\sigma/R$, the weighted norm of $R^2 f_2(Ry)+R^2|\kappa_1(Ry)|$ is then bounded by
\[
 C\left\{R^{A+B-(\varepsilon+\zeta)/2}
            +R^{1-\zeta/2}\right\}
\]
for sufficiently small $\varepsilon,\zeta>0$ and $0<\kappa<\min(1,A,B)$. The first profile has uniform mass exhaustion by the power tail; applying Lemma~\ref{hyb:scalar-inverse} with the value at zero fixed and integrating proves its mass and density rates.

\noindent\textit{Step 5. The Green potentials.}\quad In both cases take $p_0>1$ sufficiently close to one to integrate the origin powers. The remaining scale factors in the $L^{p_0}$ bounds are polynomial in $R^{-1}$, using
\[
 r=R^{s/B+o(1)}
\]
in the mixed case. Each cut-off density error has an $L^1$ bound $CR^{c_j}$ with $c_j>0$, after decreasing the exponent to absorb logarithms, and an $L^{p_0}$ bound $CR^{-K_j}$. In the interpolation identity
\[
 1/p=\lambda+(1-\lambda)/p_0,
\]
choose $\lambda<1$ so close to one that
\[
 \lambda c_j-(1-\lambda)K_j>0
\]
for every one of these errors. The Green map $L^p\to W^{2,p}$ and the two-dimensional embeddings give the $C^0\cap H^1$ bound. The cut-offs of the concentrating densities are supported or vary only on the indicated disks and tail annuli.

On their transition annuli the displayed tails yield a positive power in $L^1$, and the same interpolation applies there.
\end{proof}
At each prescribed point the preceding local reduction supplies the forcing bounds, the local mass bound and the boundary oscillations in Lemma~\ref{hyb:centered-singular-profiles}.

The regularity of the weight factors gives smooth forcing after the source logarithms are subtracted. By \eqref{fnd:independence}, $A,B,A+B\notin\mathbb Z$, and Lemma~\ref{fnd:local-unitary} gives the centered bound. Thus the lemma applies before any quantitative matching.

On a source-free disk at a mixed regular point the mass pair is $(8\pi,4\pi)$ and the residual density of the second component is smooth and positive.

Lemma~\ref{reg:stress} applies with the Taylor coefficients of the exterior field for the selected shadow solution, including the contributions of the other prescribed sources. The strict trace inequality excludes splitting and gives a bounded outer shape $B$.
\begin{proposition}[Local estimates for blowing-up solutions]
\label{hyb:local-capture}
Assume that the local profiles exhaust the masses and that
\[
 \gamma>0,\qquad \kappa_j>4|T_j|\quad\text{at each mixed regular point}.
\]
Adjust the widths by factors $1+o(1)$ and choose the free centers and shape parameters by the orthogonality conditions on the kernel functions with compact cut-off.
\begin{itemize}
\item[\textup{(i)}]
For some $c>0$ and $N_0<\infty$, at each outer radius $R$ the truncated local mass errors, the total variation errors, and the $C^0\cap H^1$ norms of the differences of the corresponding mean-zero Green potentials are bounded by
\[
 CR^c(1+|\log R|)^{N_0}.
\]
\item[\textup{(ii)}]
The logarithmic norm on the expanding comparison disks satisfies the same bound. At each mixed outer scale we use the regular part obtained by subtracting the inner model logarithm and the Dirichlet correction.
\end{itemize}
The comparison includes the residual density of the second component.
\end{proposition}
\begin{proof}
\smallskip\noindent\textit{Step 1. The prescribed points.}\quad At a prescribed point the isolated local profiles and the centered Harnack bound satisfy the assumptions of Lemmas~\ref{hyb:centered-singular-profiles} and~\ref{hyb:singular-rate}.

At a mixed point, on every fixed circle avoiding other bubbles the second canonical field converges to a finite field: its normalizing constant has a finite limit because the second component has positive residual mass, and its Green field is bounded on such circles.

Thus the boundedness assumption on this circle holds. The estimates apply on the two distinct scales, leaving the residual part in the second density.

\smallskip\noindent\textit{Step 2. The regular points.}\quad At a mixed regular point, the corresponding argument is given in Lemma~\ref{reg:local-power}. We indicate explicitly how its estimates apply here. Set $s=2$, $t=1$, $r=R^{2+o(1)}$ and
\[
 \tau=R^{3/2}.
\]
The inner profile is $m_{1,r}$ and the outer profile is $m_{R,B}$. The inner center and radius are chosen by the three compactly supported orthogonality conditions \eqref{hyb:inverse-pins}.

After subtracting the corresponding Dirichlet potential, the outer radius and the two real components of $B$ are chosen by the three conditions for $m_B$, as in the proof of Lemma~\ref{reg:local-power}.

The inner and outer radii are independent in this local argument; their relation to the common parameter is imposed only in the subsequent matching.

The decay bounds and scalar inverse estimate for $m_B$ are uniform for bounded $B$, and the outer regular part is compared with
\[
 \log\frac{16}{(1+|y^2-B|^2)^2},
\]
as in the proof of Lemma~\ref{reg:local-power}. That proof gives truncated mass errors $O(R^{1/4})$ and logarithmic profile and uniform Green errors $O(R^{1/4}|\log R|^C)$. At the endpoint $t=1$, the second moment of the inner profile satisfies
\[
 \int_{|y|<L}|y|^2m_{1,1}(y)\,dA=O(1+\log L).
\]
Thus its contribution introduces only logarithmic factors in the estimates. The exterior second density is bounded by $Cd^{-\varepsilon}$, which is locally integrable for
\[
 0<\varepsilon<2.
\]
Integration of the profile errors and tails gives the total-variation estimate. Interpolate it against the polynomial $L^{p_0}$ bounds of the truncated concentrating densities for a fixed
\[
 p_0>1.
\]
The two local radii have a polynomial relation, so an $L^1$ exponent $c_0>0$ and an $L^{p_0}$ loss exponent $K$ give the positive exponent $\lambda c_0-(1-\lambda)K$ upon fixing
\[
 1/p=\lambda+(1-\lambda)/p_0
\]
sufficiently close to one. There are finitely many local comparisons, so take the same $p>1$ for all of them. The mean-zero Green inverse
\[
 G:L^p\to W^{2,p}\hookrightarrow C^0\cap H^1
\]
gives the asserted Green estimate. Throughout this argument the remainder of the second measure is retained, as in Lemma~\ref{reg:local-power}. The ordinary case is the single-profile argument in the same lemma, with center and radius fixed by its three kernel conditions.
\end{proof}
\subsection{Matching a singular Liouville background}
The next calculation concerns a single prescribed point. Each result specifies its assumptions on the shifted parameter $A$.
\begin{lemma}[Matching through the first derivatives]
\label{univ:liouville-jet}
Let $A,c>0$ and $g(z)=c|z|^{2A-2}e^{u(z)}$ on a disk about $0$, with $u$ continuous, $u(0)=0$, and
\[
 \Delta\log g+2g=\kappa
 \quad\hbox{away from zero},
 \qquad \kappa\in C^\infty.
\]
\begin{itemize}
\item[\textup{(i)}]
Put $\beta=c/(4A^2)$. The function $u+2\log(1+\beta|z|^{2A})$ is differentiable at $0$. Denote its linear part by $2\Re(\zeta z)$.
\item[\textup{(ii)}]
For fixed $p\in\mathbb C$, set
\[
 v=p+\frac{A}{A+1}\zeta,\qquad
 h(z)=\sqrt\beta\,z^A\frac{1+vz}{1+pz},\qquad
 f=\frac{4|h'|^2}{(1+|h|^2)^2}.
\]
On a sufficiently small disk, $f$ is independent of the branch defining $h$. The function $J=\log(g/f)$ extends to a $C^2$ function at $0$, and
\[
 J(0)=0,\qquad \nabla J(0)=0,\qquad
 J(z)=O(|z|^2),\qquad \Delta J(0)=\kappa(0).
\]
The bounds are uniform on compact bounded families of data with $A,c>0$ fixed.
\end{itemize}
\end{lemma}
\begin{proof}
\smallskip\noindent\textit{Step 1. The linear part at zero (i).} Compare $g$ with $f_0=c|z|^{2A-2}(1+\beta|z|^{2A})^{-2}$. The difference $J_0=\log(g/f_0)$ is continuous, vanishes at zero, and satisfies
\[
 \Delta J_0+V_0J_0=\kappa,\qquad
 V_0=2f_0\int_0^1e^{tJ_0}\,dt=O(|z|^{2A-2}).
\]
The equation for $J_0$ holds across zero in distributions. Indeed its right-hand side $\kappa-V_0J_0$ belongs to $L^p$ for some $p>1$, since $A>0$ and $J_0$ is bounded. Subtract its continuous logarithmic potential on a smaller disk.

The difference is bounded and harmonic on the punctured disk, so its singularity is removable. The same argument applies to the matched difference $J$ below. We give the local bootstrap used here. Suppose
\[
 \Delta v=F,\qquad |F(z)|\le C|z|^{\lambda-2}
\]
with $\lambda>0$ noninteger, and put
\[
 m=\lfloor\lambda\rfloor.
\]
The logarithmic potential has a Taylor polynomial of order $m$. Its coefficients are integrable since each kernel derivative of order $j\leq m$ has absolute bound $C|y|^{-j}$ for $j\geq1$. For $|y|>2|z|$, the Taylor error in the kernel integrates to at most
\[
 C|z|^{m+1}\int_{2|z|}^{d}t^{\lambda-m-2}\,dt=O(|z|^\lambda).
\]
For $|y|\le2|z|$, first subtract $\log|y|$ from $\log|z-y|$ so that the additive $\log|z|$ cancels on rescaling; the resulting logarithms are integrable against $|y|^{\lambda-2}$ in the scaled disk. The remaining Taylor terms of orders $1\leq j\leq m$ cost
\[
 C|z|^j\int_0^{2|z|}t^{\lambda-j-1}\,dt=O(|z|^\lambda).
\]
Add the Taylor polynomial of the harmonic part. At an integer exponent the same estimate holds with any slightly smaller noninteger exponent. The constants are uniform on a smaller disk.

Initially $J_0$ is bounded, so this estimate gives some positive H\"older exponent. If
\[
 J_0=O(|z|^\sigma),
\]
its equation improves the available vanishing order to any noninteger order below $\min\{2,\sigma+2A\}$, up to a possible linear harmonic term. After finitely many repetitions this proves differentiability at zero and $J_0=2\Re(\zeta z)+O(|z|^{1+\epsilon})$ for some $\epsilon>0$.

\smallskip\noindent\textit{Step 2. Matching the first derivative (ii).} The choice of $v-p$ makes $\log(f/f_0)=2\Re(\zeta z)+O(|z|^{1+\epsilon_1})$, where $\epsilon_1>0$. Indeed
\[
 h'=A\sqrt\beta\,z^{A-1}
 \left\{1+\frac{A+1}{A}(v-p)z+O(|z|^2)\right\},
\]
and the difference between the two logarithmic denominator terms is $O(|z|^{2A+1})$. Consequently $J=\log(g/f)$ has zero value and gradient and satisfies
\[
 \Delta J+VJ=\kappa,\qquad
 V=2f\int_0^1e^{tJ}\,dt=O(|z|^{2A-2}).
\]
Apply the preceding bootstrap again, now with no constant or linear term. It gives $J=O(|z|^{2-\varepsilon})$ for any sufficiently small $\varepsilon>0$. Choose
\[
 \varepsilon<2A.
\]
The right side of $\Delta J=\kappa-VJ$ is then $\kappa(0)+O(|z|^{\delta})$ for some
\[
 \delta>0.
\]
We also verify the derivative estimates needed for a $C^2$ extension. On $r/2<|z|<2r$, rescale $J$ by $r^{2-\varepsilon}$. The rescaled coefficient $r^2V$ is bounded by $Cr^{2A}$, and the rescaled smooth forcing is bounded. Interior $W^{2,p}$ estimates with any fixed $p>2$ therefore give
\[
 |\nabla J|\le Cr^{1-\varepsilon}
\]
on a smaller annulus. The explicit Liouville density satisfies $|\nabla f|\le Cr^{2A-3}$. Consequently the function
\[
 F=\kappa-\kappa(0)-2f(e^J-1)
\]
satisfies $|F|\le Cr^\delta$ and $|\nabla F|\le Cr^{\delta-1}$, after decreasing $\delta>0$ if necessary. Subtract $\kappa(0)|z|^2/4$ and use its logarithmic potential. On $|y-z|<|z|/2$, subtract $F(z)$ in the second-derivative kernel and use the displayed gradient bound. On the remaining inner region use
\[
 |F(y)|\le C|y|^\delta;
\]
on $|y|>2|z|$ subtract the Taylor polynomial of the kernel at zero. These estimates give a quadratic Taylor polynomial $Q$ and, for any fixed $0<\delta'<\min\{\delta,1\}$,
\[
 |D^j(J-Q)(z)|\le C|z|^{2+\delta'-j},\qquad j=0,1,2.
\]
For the H\"older bound across zero, these estimates give $F\in C^{0,\delta'}$ on a smaller disk. For $|z-z'|\ge \frac12\max(|z|,|z'|)$ use $|F|\le Cr^\delta$ at each point. Otherwise the segment between the points stays at a comparable radius $r$, and
\[
 |F(z)-F(z')|\le Cr^{\delta-1}|z-z'| \le Cr^{\delta-\delta'}|z-z'|^{\delta'}.
\]
Interior Schauder estimates on the extended equation $\Delta J=\kappa(0)+F$ now give $J\in C^{2,\delta'}$ and
\[
 \Delta J(0)=\kappa(0).
\]
In particular no $|z|^2\log|z|$ term occurs: the nonconstant part of the order-zero forcing tends to zero with a positive power.
\end{proof}
The matched Taylor coefficients and the remainder are $C^1$ in the finite-dimensional background parameters. More precisely, suppose $c$, $\kappa$, and the original branch of scalar shadow solutions depend $C^1$ on such parameters, the last in $C^0$ and in its weighted elliptic space.

These properties follow from the scalar implicit function theorem on a compact nondegenerate branch. Differentiate the equation for $J_0$ before extracting its linear part:
\[
 \Delta\dot J_0+2g\dot J_0
       =\dot\kappa-2\dot f_0(e^{J_0}-1).
\]
Its right side is $O(1)+O(r^{2A-1})$, since
\[
 J_0=O(r),\qquad \dot f_0=O(r^{2A-2}).
\]
The same uniform weighted Poisson estimates show that $\dot J_0$ has a linear expansion, and apply also to the parameter difference quotients. Thus $\zeta$ depends $C^1$ on the parameters. After matching it, differentiate the equation for $J$:
\[
 \Delta\dot J+2g\dot J
       =\dot\kappa-2\dot f(e^J-1).
\]
Now $\dot J$ has zero value and gradient and the additional forcing is $O(r^{2A})$. Repeating the argument for a remainder with zero value and gradient, and using the derivative estimates on annuli, proves uniform $C^{2,\delta'}$ bounds for $\dot J$, with a common sufficiently small
\[
 \delta'>0.
\]
The same estimates hold for the difference quotients in the parameters. Their value and gradient at zero vanish by the matching conditions. Subtraction gives a zeroth-order coefficient $O(r^{2A-2})$.

The right-hand side has the same positive exponent. The uniform $C^{2,\delta'}$ bounds, together with convergence away from zero, give convergence in $C^{2,\delta''}$ for every
\[
 0<\delta''<\delta'.
\]
Thus the asserted parameter differentiability holds in this common smaller H\"older space. For $J_1$ below the sum $J_1+J_2/2$ solves a smooth Poisson equation: the density terms in the equations for $\Delta J_2$ and $\Delta J_1$ cancel. The bounded sum extends across zero.

Its smooth forcing and boundary values on a fixed small circle depend $C^1$ on the background parameters. The interior Poisson estimates, followed by subtraction of $J_2/2$, give the same uniform Taylor remainder and parameter differentiability for $J_1$.
\subsection{An exact local Toda family with prescribed first derivatives}
Fix $A,B>0$, put $s=A+B$, and suppose $s\ne1$ and $B\ne1$. Let $k>0$, $\tau=R^{2s}$, and choose complex numbers $a_0,a_1$. On a small punctured disk define the holomorphic vector
\[
 \nu_R(z)=\bigl(R^s(1+a_0z),\ \sqrt{k}R^s z^A(1+a_1z),\ z^s\bigr),
\]
and set
\[
 F_R=|\nu_R|^2,\qquad Q_R=|\nu_R\wedge\nu_R'|^2,\qquad
 W_R=\det(\nu_R,\nu_R',\nu_R''),
\]
\[
 n_{1,R}=\frac{4Q_R}{F_R^2},\qquad
 n_{2,R}=\frac{4|W_R|^2F_R}{Q_R^2}.
\]
Analytic continuation once around zero multiplies the components of $\nu_R$ by a diagonal unitary matrix. Thus both densities are single-valued. The Pl\"ucker identities give
\[
 \Delta\log F_R=n_{1,R},\qquad \Delta\log Q_R=n_{2,R},
\]
so the densities solve the Toda system exactly on the punctured disk. Their source strengths are $A-1$ and $B-1$. The disk can be chosen independently of $R$, since
\[
 W_R=\sqrt{k}\,\tau z^{A+s-3}W_0(z),\qquad W_0(0)=AsB\ne0,
\]
where $W_0$ is a polynomial independent of $R$. Here are explicit expressions convenient for matching. Put
\begin{align*}
 P_{01}(z)&=A+\{(A-1)a_0+(A+1)a_1\}z+Aa_0a_1z^2,\\
 P_{02}(z)&=s+(s-1)a_0z,\qquad
 P_{12}(z)=B+(B-1)a_1z,\\
 F_0(z)&=|1+a_0z|^2+k|z|^{2A}|1+a_1z|^2,\\
 D(z)&=|P_{02}(z)|^2+k|z|^{2A}|P_{12}(z)|^2,\\
 H_R(z)&=k\tau|z|^{-2B}|P_{01}(z)|^2.
\end{align*}
Then
\[
 F_R=\tau F_0+|z|^{2s},\qquad
 Q_R=\tau|z|^{2s-2}(D+H_R).
\]
Direct expansion of the determinant also gives
\[
 W_0=A P_{02}P_{12}
            +z(P_{12}'P_{02}-P_{12}P_{02}').
\]
The exterior limits of $\log n_{1,R}$ and $n_{2,R}$ are
\[
 L_{1,R}=\log(4\tau)-(2s+2)\log|z|+\log D(z),
 \qquad
 f=\frac{4k|z|^{2A-2}|W_0(z)|^2}{D(z)^2}.
\]
Equivalently,
\[
 f=\frac{4|h'|^2}{(1+|h|^2)^2},\qquad
 h(z)=\sqrt{k}\,z^A\frac{P_{12}(z)}{P_{02}(z)}.
\]
In particular $f\sim c|z|^{2A-2}$, where $c=4A^2B^2k/s^2$. Put $\beta=B^2k/s^2$ and write the linear parts after subtraction of the radial fractional terms as $2\Re(g_1z)$ and $2\Re(g_2z)$. They are
\[
 g_1=\frac{s-1}{s}a_0,\qquad
 g_2=\frac{A+1}{A}
       \left(\frac{B-1}{B}a_1-\frac{s-1}{s}a_0\right).
\]
Thus arbitrary prescribed $g_1,g_2$ are obtained by choosing
\[
 a_0=\frac{s}{s-1}g_1,\qquad
 a_1=\frac{B}{B-1}
       \left(g_1+\frac{A}{A+1}g_2\right).
\]
This construction uses only the nonresonance conditions
\[
 s,B\ne1;
\]
the model inverse may require the stronger noninteger assumptions. Suppose the canonical exterior fields are $L_{1,R}+J_1$ and $\log f+J_2$, with their constant terms and these first derivatives matched. Lemma~\ref{univ:liouville-jet} gives $J_2\in C^2$ and
\[
 J_2=O(|z|^2).
\]
The sum of the actual first exterior field and half the second is smooth after subtraction of its logarithmic term. Since $D=|P_{02}|^2(1+|h|^2)$ and $P_{02}$ is nonzero on the disk, $\Delta\log D=4|h'|^2/(1+|h|^2)^2=f$ there away from zero. Also
\[
 \Delta\log f=-2f,
\]
so the exact sum $L_{1,R}+\tfrac12\log f$ is harmonic after subtracting its logarithmic term and extends by boundedness. Matching its constant and linear terms to those of the smooth actual sum and subtracting $J_2/2$ gives $J_1\in C^2$ and
\[
 J_1=O(|z|^2).
\]
More explicitly, if the canonical first exterior field $E_1$ satisfies
\[
 \Delta E_1=g+\kappa_1
\]
away from zero, where $\kappa_1$ is smooth, then $\Delta L_{1,R}=f$ gives
\[
 \Delta J_1=\kappa_1+f(e^{J_2}-1),\qquad
 \Delta J_1(0)=\kappa_1(0).
\]
Here $f(e^{J_2}-1)=O(r^{2A})$ tends to zero for every $A>0$. Thus the finite trace used in the scale equation is the trace of the matched difference $J_1$, even when the uncorrected exterior field itself is not $C^2$. In a normalized conformal coordinate it is
\[
 \kappa_1(0)=\Delta\log h_{1,\mathrm{reg}}(0)
                   +2M-\rho_2-2K_g(p),
 \qquad h_1=|z|^{2(A-1)}h_{1,\mathrm{reg}}.
\]
In particular a perturbation of the second derivative of $\log h_1^*$ with zero value and first derivative changes this trace without changing the matched first derivatives.
\subsection{Density comparisons and the local residual}
Let $S=m_{s,R}$, $I=m_{B,\ell}$, where $\ell^{2B}=kA^2\tau/s^2$, and let $d>0$ be a sufficiently small fixed disk radius. All constants below may depend on the fixed profile exponents and on a compact set of Taylor coefficients. Write $r=|z|$.
\begin{lemma}[Uniform estimates for the exact local family]
\label{univ:profile-estimates}
On $D_d$ the following estimates hold.
\begin{itemize}
\item[\textup{(i)}]
The local densities satisfy
\[
 n_{1,R}\le C S\{1+r^{2A}+(\ell/r)^{2B}\},
\]
\[
 |n_{2,R}-f|
 \le C\left\{I+r^{2A-2}\min\{1,(\ell/r)^{2B}\}\right\}.
\]
For some $\eta>0$,
\[
 \|n_{1,R}-S\|_1+\|n_{2,R}-I-f\|_1\le CR^\eta.
\]
\item[\textup{(ii)}]
For each fixed $p>1$,
\begin{align*}
 \|r^2n_{1,R}\|_p
   &\le C_p(R^{2/p}+\tau)(1+|\log R|),\\
 \|r^2(n_{2,R}-f)\|_p
   &\le C_p(\ell^{2/p}+\tau)(1+|\log R|).
\end{align*}
\item[\textup{(iii)}]
The estimates in \textup{(i)} and \textup{(ii)} also hold after one differentiation in $\log R$, with possibly different constants and smaller $\eta$. If $ps<1$ and $pB<1$, the two weighted norms in \textup{(ii)} and their derivatives in $\log R$ are $O(\tau)$.
\end{itemize}
\end{lemma}
\begin{proof}
\smallskip\noindent\textit{Step 1. Pointwise estimates in (i).}\quad On a sufficiently small disk, $F_0\asymp1+r^{2A}$, $D\asymp1+r^{2A}$, $H_R\asymp(\ell/r)^{2B}$, and $|W_0|\asymp1$. The formula for $n_1$ gives the pointwise bound and the estimate
\[
 |n_{1,R}-S|\le CS\{r+r^{2A}+(\ell/r)^{2B}\}.
\]
The integrals of the last two terms are $O(R^{2A})$. The integral of the first is $O(R^{\min\{1,2s\}}(1+|\log R|))$. For the second density the exact formula gives
\[
 n_{2,R}
 =\frac{4k\tau r^{-2B-2}|W_0|^2F_0}{(D+H_R)^2}
   +f\frac{D^2}{(D+H_R)^2}.
\]
The first summand is $I\{1+O(r+r^{2A})\}$. Since $H_R/D\asymp(\ell/r)^{2B}$, the difference between the second summand and $f$ is bounded by $Cr^{2A-2}\min\{1,(\ell/r)^{2B}\}$. This proves the pointwise bound for $n_{2,R}-f$ and gives
\[
 \|n_{2,R}-I-f\|_1
 \le C\ell^{\min\{1,2A,2B\}}(1+|\log\ell|).
\]
Both comparisons are therefore bounded by a positive power of $R$.

\smallskip\noindent\textit{Step 2. The integral estimates.}\quad For the $L^p$ estimates we split the radial integrals at $R$ or $\ell$. The scalar estimate is
\[
 \|r^2m_{q,a}\|_{L^p(D_d)}
 \le C_p(a^{2/p}+a^{2q})(1+|\log a|).
\]
The two extra terms in the first density are bounded respectively by
\[
 C_p(R^{2A+2/p}+\tau)(1+|\log R|),\qquad
 C_p(R^{2A+2/p}+\tau^2)(1+|\log R|).
\]
For the second density the additional term has norm at most
\[
 C_p\{\ell^{2A+2/p}+\ell^{2B}\}(1+|\log\ell|).
\]
Use $\ell^{2B}\asymp\tau$ to obtain the claimed bounds. Suppose now that $ps<1$ and $pB<1$. On $r<R$ the term $r^2S(\ell/r)^{2B}$ is bounded by $Cr^{2A}$, since
\[
 \ell^{2B}\asymp R^{2s}=\tau.
\]
On $r>R$ that term is bounded by $Cr^2S$. The other first-component terms have the same bound as $r^2S$ up to a fixed factor. Outside $R$ or $\ell$, the remaining majorants are
\[
 C\tau r^{-2s}
\]
for $r^2S$, $C\tau r^{-2B}$ for $r^2I$, and $C\tau r^{2A-2B}$ for the second diffuse difference. Their $L^p$ norms are $O(\tau)$ under the two strict inequalities.

On the inner disks the integrals give $O(R^{2/p}+R^{2A+2/p}+\ell^{2/p}+\ell^{2A+2/p})$, which is again $O(\tau)$. To differentiate in $\log R$, fix $z,A,B,k,a_0,a_1$. Put
\[
 a=\tau/r^{2s},\qquad b=(\ell/r)^{2B},
\]
so $\partial_{\log R}a=2sa$ and $\partial_{\log R}b=2sb$. Write $D=s^2(1+D_1)$ and $H_R=s^2b(1+H_1)$. The explicit polynomials give
\[
 F_0-1,D_1,H_1=O(r+r^{2A}).
\]
For any such positive multiplier $U=1+O(r+r^{2A})$ independent of $R$,
\[
 \frac{xU}{1+xU}-\frac{x}{1+x}
  =\frac{x(U-1)}{(1+xU)(1+x)} .
\]
This and its logarithmic derivative are bounded by $C|U-1|x/(1+x)$ for $x=a,b$. Use $U=F_0$ and $U=(1+H_1)/(1+D_1)$, respectively, in the two density fractions. The differentiated diffuse fraction satisfies
\[
 |\partial_{\log R}((1+H_R/D)^{-2}-1)| \leq C\min(1,b).
\]
The other factors are independent of $R$ or have bounded logarithmic derivatives. This preserves the difference bounds and the individual majorants in the preceding integrals.
\end{proof}
Replace the exterior fields on the inner disk by adding $(\log n_{1,R}-L_{1,R},\log n_{2,R}-\log f)$, and use a fixed cut-off to match to the exterior fields. Before global normalization, the exact local residual, apart from its Cartan matrix factor, is
\[
 -\binom{n_{1,R}(e^{J_1}-1)}
              {(n_{2,R}-f)(e^{J_2}-1)}.
\]
Indeed
\[
 -\Delta(L_{1,R},\log f)^t=K(0,f)^t
\]
on the punctured disk, whereas the exterior residual density is $fe^{J_2}$.

The identification introduces no point mass in the physical residual: the local canonical fields $\log n_{i,R}+J_i$ have logarithmic coefficients $2(A-1)$ and $2(B-1)$, exactly those subtracted in passing from the canonical fields to the regular physical fields.

The residual identity therefore holds on the full disk in distributions. It proves that the local residual has $L^p$ norm bounded by
\[
 C_p\{R^{2/p}+\ell^{2/p}+\tau\}(1+|\log R|).
\]
The same estimate holds after a logarithmic-scale derivative. Set $U=\tau F_0/r^{2s}$ and $V=H_R/D$. Applying $K^{-1}$ to the two logarithmic corrections gives $(-\log(1+U),-\log(1+V))$ exactly. At fixed matched coefficients with
\[
 \tau=d_0^se^{-T}
\]
their $T$ derivatives are $(U/(1+U),V/(1+V))$. On the fixed cut-off annulus both corrections and these derivatives are $O(\tau)$ with any fixed number of spatial derivatives. Thus the cut-off contributes only $O(\tau)$.

Integrating $\Delta\log F_R$ and $\Delta\log Q_R$ on the fixed disk gives masses $4\pi s+O(\tau)$ and
\[
 4\pi B+\int f+O(\tau)
\]
: the outer coefficients of $\log r$ are $2s$, $2s-2$ with the additional term $\log D$ in the second potential, and the inner coefficients are $0$, $2A-2$, with
\[
 \Delta\log D=f.
\]
The same flux calculation controls the first scale derivatives of the errors. The displayed residual bounds control the remaining mass defects after inserting $e^{J_i}$; each normalization error costs at most $O(\sqrt T)$ times its mass defect in $H^{-1}$.

For finitely many fixed positive shifts choose $p>1$ near one so the squared residual and normalization errors, even times fixed powers of $T$, are $o(\mathcal L)$, where
\[
 \mathcal L=e^{-T/s_*}
\]
if $s_*>1$ and $e^{-T}$ if $s_*\le1$. Indeed the rate is $e^{-T\min(1,1/(ps))}$ up to polynomial factors, with $p<2$ and $s\le s_*$. The other kinds of disks require their corresponding residual bounds in the global construction.
\begin{lemma}[Identification with the projected family]
\label{journal:strict-capture}
Let $u=B_T+\eta$ be a solution. After matching and imposing the orthogonality conditions, assume $L_\alpha\eta=0$ for every retained source and
\[
 |\mu|+\|u-B_T\|_{H^1}\leq C e^{-cT}T^N,
 \qquad c>0.
\]
For some $0<\delta_0<c$, assume the projected solution exists on the parameter domain and is unique in the ball given respectively by
\[
 |\mu|\leq e^{-\delta_0T},\qquad
 \|\eta\|_{H^1}\leq e^{-\delta_0T}T^{C_0}.
\]
For $T$ sufficiently large the two inequalities are strict for $u$, and $u$ equals the projected solution.
\end{lemma}
\begin{proof}
The ratios to the parameter radius and correction radius are bounded, respectively, by \(Ce^{-(c-\delta_0)T}T^N\) and \(Ce^{-(c-\delta_0)T}T^{N-C_0}\). Both tend to zero. The mean-zero Toda solution satisfies the projected equation and the parameters were chosen so that
\[
 L_\alpha\eta=0.
\]
Uniqueness in this same ball identifies it with the projected solution. Whenever the projected residual is written as a sum of source representatives, the independent sources make its vanishing equivalent to vanishing of every multiplier.

The argument is uniform when the constants and exponents are uniform. For finitely many shadow solutions take the minimum positive exponent before choosing \(\delta_0\).
\end{proof}
\subsection{Matching the blow-up parameters}
We return to $A_p,B_p>1$ and the separated configurations under consideration.
\begin{lemma}[Matching the parameters of solutions and of the projected equation]
\label{hyb:capture-matching}
Assume the stated local masses at the critical first parameter and $\gamma>0$. Use the fitted centers for the scalar branch. Let $D$ be the largest fitted outer radius and set
\[
 \delta=CD^c(1+|\log D|)^N,\qquad c>0,
\]
decreasing $c$ if necessary. Assume the following estimates.
\begin{itemize}
\item[\textup{(a)}]
At each bubble the fitted concentrating densities satisfy the stated $L^1$ and Green potential estimates by positive powers of their radii. The ratio of the inner to the outer radius is bounded by a positive power of the outer radius.
\item[\textup{(b)}]
For each inner cut-off radius $\tau$,
\[
 r/\tau\le CR^\eta,\qquad \tau\ge cR^K,\qquad \eta,K>0\ \text{fixed}.
\]
The cut-offs are one on the corresponding fixed rescaled annuli.
\item[\textup{(c)}]
For some fixed $p>1$, the cut-off remainders satisfy
\[
 \|f_1-\nu_1\|_{L^p}+\|f_2-\nu_2-C_uH_pe^{w_u}\|_{L^p}=O(\delta),\quad
 \int_M\nu_2=\rho_2-\gamma+O(\delta),\quad
 \|w_u-w(p)\|_{W^{2,p}}=O(\delta).
\]
Here $w_u=u_2+G\nu_1-2G\nu_2$ and $C_u=\rho_2/\int_Mh_2e^{u_2}$.
\end{itemize}
Then the following conclusions hold.
\begin{itemize}
\item[\textup{(i)}]
All fitted outer and inner radii satisfy the matching relations up to relative error $O(\delta)$, with the same constant at every point,
\[
 c_1=\log\bigl(\rho_1/\textstyle\int_M h_1e^{u_1}\bigr).
\]
\item[\textup{(ii)}]
In this section and for general source strengths, set $T=-c_1$. Then $|\mu|=O(e^{-c_0T})$ for some $c_0>0$. With the exactly matched radii,
\[
 \|u-B_T\|_{C^0\cap H^1}\le C e^{-c_0T}(1+T)^N .
\]
For regular bubbles alone the corresponding estimates use $t^2=e^{c_1}$ and powers of $t$.
\end{itemize}
\end{lemma}
\begin{proof}
\noindent\textit{Step 1. The local Green potentials.}\quad Let $f_i=\rho_i h_i e^{u_i}/\int_M h_i e^{u_i}$ denote the densities with respect to $dV_g$. Let $\nu_i$ be the cut-off concentrating densities used in the scalar comparison,
\[
 Z=K(G\nu_1,G\nu_2)^t,\qquad w_u=u_2-Z_2.
\]
The exact Green identity is
\[
 u=Z+(-w_u/2,w_u)^t+\bigl(3G(f_1-\nu_1)/2,0\bigr)^t .
\]
By the quantitative comparisons, $w_u$ can be replaced by the selected scalar branch $w(p)$, with error $O(\delta)$ in $W^{2,p}$, $p>1$. If $I=\int H_pe^{w(p)}$, integration of $f_2-\nu_2=C_uH_pe^{w_u}+e_2$ gives
\[
 \log C_u=\log(\gamma/I)+O(\delta),\qquad
 C_u=\rho_2/\int_M h_2e^{u_2},
\]
since $\gamma>0$ and $I$ is bounded below. In particular, the actual second normalization may be replaced by the shadow normalization
\[
 c_2=\log(\gamma/I)
\]
with this error. The estimate for the concentrating densities allows us to replace $G\nu_i$ by the potentials of the fitted profiles with an $O(\delta)$ error. On a smaller fixed disk about a fitted center, the potential of a cut-off radial profile with exponent $q$ equals
\[
 4\pi qG(z,p)-\log\{1+(R/|z|)^{2q}\}+O(R^{\eta})
\]
with some $\eta>0$ in $C^0\cap H^1$. For the regular mixed outer profile this formula is
\[
 8\pi G(z,p)
  -\log\{(R^4+|z^2-R^2B|^2)/|z|^4\}+O(R^{\eta}).
\]
To check the constants, the planar logarithmic potentials are $-\log(R^{2q}+|z|^{2q})$ and $-\log(R^4+|z^2-R^2B|^2)$, respectively. Their negative Laplacians are precisely the profiles and their differences from $-2q\log|z|$ and $-4\log|z|$ tend to zero at infinity.

In the surface Green representation, integrate the smooth kernel remainder against the profile after subtracting its value at the center. The error and its first spatial derivative are bounded by a fractional moment
\[
 \int |z|^{\eta}m\,dA=O(R^\eta)
\]
for a sufficiently small
\[
 \eta>0.
\]
Power decay gives the same bound for the omitted planar tail and its logarithmic potential on a smaller disk. Outside that disk use separation from the singularity and smoothness of the kernel. Each concentration radius is treated separately.

If an inner profile had previously been cut off on an intermediate disk, extending its model to the fixed disk produces an error of at most $C(r/\tau)^{2q}\tau^{-(2-2/\widehat p)}$ in $L^{\widehat p}$ by integrating the radial tail. Choose
\[
 1<\widehat p\leq p\quad\Longrightarrow\quad K(2-2/\widehat p)<2q\eta
\]
at every such disk. Then
\[
 r/\tau\leq CR^\eta,\qquad \tau\geq cR^K
\]
give the positive power $2q\eta-K(2-2/\widehat p)$ of the corresponding outer radius. The mean-zero Green map sends this error to $C^0\cap H^1$ with the same bound.

In particular these potential expansions give, after mean subtraction, $u^0$ plus the potential corrections of the independent fitted profiles, at an error of order $\delta$.

\noindent\textit{Step 2. Matching on the rescaled annuli.}\quad We compare on the fixed scaled annulus $\mathcal A=\{1<|y|<2\}$. Write $F_i=e^\varphi f_i$ for the canonical densities and put $P_q(y)=|y|^{2q-2}(1+|y|^{2q})^{-2}$. For a radial outer bubble use
\[
 \psi_1=\log h_1+\varphi+u_1^0+(2s+2)\log|z|.
\]
On $z=Ry$ the outer Green factor in $F_1$ is $(1+|y|^{-2s})^{-2}$, and the inner factor when present differs from one by a positive power. Since $\psi_1$ is uniformly H\"older by the smooth factors and $w(p)\in W^{2,p}$, this gives
\[
 R^2F_1(Ry)=e^{c_1+\psi_1(0)}R^{-2s}P_s(y)\{1+O(\delta)\}.
\]
The fitted density on this annulus is $4s^2P_s$. Its integral differs from that of $R^2F_1(Ry)$ by $O(\delta)$ by the fitted $L^1$ comparison, since the outer cut-off equals one there. The integral of $P_s$ on $\mathcal A$ is positive. Dividing the two integrated expressions gives $4s^2R^{2s}=e^{c_1+\psi_1(0)}\{1+O(\delta)\}$.

At a mixed prescribed source put $\psi_2=\log h_2+\varphi+u_2^0-(2A-2)\log|z|$. On $z=ry$ the outer factor is $R^{2s}(r|y|)^{-2s}\{1+O(\delta)\}$ and the inner factor is $(1+|y|^{-2B})^{-2}$. Hence
\[
 r^2F_2(ry)=e^{c_2+\psi_2(0)}R^{2s}r^{-2B}
                  P_B(y)\{1+O(\delta)\}.
\]
This annulus lies inside the inner cut-off. Comparing its integral with that of $4B^2P_B$ gives $4B^2r^{2B}=e^{c_2+\psi_2(0)}R^{2s}\{1+O(\delta)\}$. In the exact singular family of Lemma~\ref{univ:profile-estimates} the coefficient
\[
 e^{c_2+\psi_2(0)}=4A^2B^2k/s^2
\]
gives the same relation with designated inner radius $\ell^{2B}=kA^2R^{2s}/s^2$. For a mixed regular bubble use $u_1^0=-6\log|z|+O(1)$ and put
\[
 P_B^{\rm out}(y)=|y|^2/(1+|y^2-B|^2)^2.
\]
On $z=Ry$ the nonradial Green factor in the first density is $|y|^8/(1+|y^2-B|^2)^2$. Thus
\[
 R^2F_1(Ry)=e^{c_1+\psi_1(0)}R^{-4} P_B^{\rm out}(y)\{1+O(\delta)\}.
\]
Its integral differs from that of $16P_B^{\rm out}$ by $O(\delta)$. The integral of $P_B^{\rm out}$ is bounded uniformly away from zero on the compact shape set, yielding
\[
 16R^4=e^{c_1+\psi_1(0)}\{1+O(\delta)\}.
\]
On $z=ry$ the limiting second field has no logarithm and the outer correction contributes $R^4(1+|B|^2)(r|y|)^{-4}\{1+O(\delta)\}$. The inner correction is $(1+|y|^{-2})^{-2}$. Writing
\[
 \beta=e^{c_2+\psi_2(0)},
\]
the annular comparison is now between $\beta R^4(1+|B|^2)r^{-2}P_1(y)\{1+O(\delta)\}$ and $4P_1(y)$. It gives
\[
 4r^2=\beta R^4(1+|B|^2)\{1+O(\delta)\}.
\]
To make each annular comparison quantitative, use a fixed annulus in that bubble's own rescaled variables.

If its profile is $P$, the Green and shadow errors give the rescaled density $aP(1+O(\delta))$ there, and the comparison of concentrating densities gives $bP+O_{L^1}(\delta)$, with $b$ the model coefficient.

On the compact shape set $\int P$ is bounded below on this annulus. Its fixed positive normalization and the mass bound first give
\[
 a=O(1).
\]
Subtract the two annular integrals and divide by $\int P$ to get $a-b=O(\delta)$. This also justifies using an $L^1$ estimate on the error of the concentrating densities when matching the common scale.

\noindent\textit{Step 3. The error after matching.}\quad The coefficients are bounded above and below, so
\[
 D\asymp e^{-T/(2s_*)},
\]
where $s_*$ is the largest outer profile exponent. Integrating the first density on the individual disks and using the bound on its exterior mass by a positive power of the radius proves
\[
 |\mu|=O(\delta).
\]
The replacement of each independent fitted logarithmic radius by the exactly matched one has size $O(\delta)$. For the radial correction its derivative is
\[
 \partial_{\log R}\bigl[-2\log(1+(R/|z|)^{2q})\bigr]
 =-\frac{4qR^{2q}}{R^{2q}+|z|^{2q}}.
\]
It is uniformly bounded in $C^0$ and in $H^1$ after cut-off. The Dirichlet integral in the latter estimate is uniformly bounded by setting
\[
 z=Ry.
\]
The same computation for the regular mixed nonradial correction uses the smooth rescaled denominator $1+|y^2-B|^2$ and its power decay, uniformly on the bounded shape set. Thus replacing the radii produces an error of $O(\delta)$ in both norms.

It remains to compare the Poisson sources of the potential corrections and the localized exact profiles. At a mixed regular point, the profile pairs are the two scalar densities and the linear combination of the density $\beta(1+\beta |z|^2/4)^{-2}$ supplied by the two smooth tails.

By the estimates on the local disks in Proposition~\ref{hyb:trial-comparison}, the difference from the exact densities $n_1,n_2$ is bounded in $L^1$ by a positive power of the radius and in $L^{p_0}$ by a power of its reciprocal, for some
\[
 p_0>1.
\]
For the exact family at a mixed prescribed point use the comparison in Lemma~\ref{univ:profile-estimates}, with $f$ supplying the tails. For positive source strengths we use the fitted radial pair at a mixed prescribed point, together with the first-order corrections below.

The singular logarithms agree in all these comparisons. The cut-off potential of the source difference is bounded by a positive power of the radius in both norms. Indeed, at each radius interpolate the $L^1$ and $L^{p_0}$ bounds, then apply
\[
 G:L^p\to W^{2,p}\hookrightarrow C^0\cap H^1
\]
for a fixed $p>1$ close to one. On the fixed matching annulus the exact corrections and their spatial derivatives are small by the explicit formulas. The remaining harmonic correction has small $C^1$ boundary data on the fixed circle by the annular estimates.

Harmonic extension gives the supremum bound by the maximum principle and the energy bound by minimizing against a radial cut-off extension of that trace. Thus cutting off and subtracting the surface means preserves the comparison in both norms.

At a pure point the correction for the matched linear term is a cut-off multiple of $g\cdot z/(1+(|z|/R)^{2q})$ with $q>0$. Its supremum is bounded by
\[
 CR^{\min(1,2q)}.
\]
The gradient is bounded by $C$ for $|z|\le R$ and by $C(R/|z|)^{2q}$ between $R$ and the fixed cut-off radius $d_0$. Thus its squared Dirichlet norm is at most
\[
 C\left(R^2+R^{4q}\int_R^{d_0}d^{1-4q}\,{\rm d}d\right),
\]
which gives $O(R^{\min(1,2q)}(1+|\log R|)^{1/2})$ in $C^0\cap H^1$. First-order corrections for positive source strengths have the same form, with the inner radius in the inner terms. This proves the claimed field estimate.
\end{proof}
\begin{theorem}[Description of all blowing-up solutions]
\label{hyb:capture}
Consider a separated blowing-up sequence at a selected solution of the shadow system, with at least one singular blow-up point.
\begin{itemize}
\item[\textup{(i)}]
For all sufficiently large indices the solutions belong to the family of Proposition~\ref{hyb:range}. The parameters are chosen by $L_\theta(u-B_\theta)=0$, and the shape parameters are bounded.
\item[\textup{(ii)}]
The conclusion holds before the sharp mass expansion, including when the shape equations have multiple zeros.
\end{itemize}
\end{theorem}
\begin{proof}
\noindent\textit{Step 1. The residual density.}\quad Let $\nu_i$ be the sums of fitted, truncated concentrating measures, and put
\[
 Z=K(G\nu_1,G\nu_2)^t,
\]
\[
 w_u=u_2-Z_2,\qquad D=\max_aR_a.
\]
Outside the blow-up disks the $L^p$ norm of the first density is bounded by a positive power of $D$ for some fixed $p>1$.

On compact subsets away from the blow-up points, the estimate for small local mass in Proposition~\ref{fnd:signed-compactness} bounds the mean-zero fields, including on the disks of unused sources.

Evaluate the first density on one fixed profile-disk boundary, away from zeros of $h_1$. The local tail estimate there, the lower bound for $u_1$, and
\[
 f_1=e^{c_1}h_1e^{u_1}
\]
bound $e^{c_1}$ by a positive power of the corresponding local radius. Multiplying this same constant by the bounded exterior $e^{u_1}$ and the fixed $L^p$ weight gives the required exterior bound before any quantitative matching. The exact identity
\[
 u_1-Z_1+\tfrac12w_u=\tfrac32G(f_1-\nu_1)
\]
then has a right-hand side bounded by a positive power of $D$ in $C^0\cap H^1$. The measure errors are also small in some fixed $L^p$, $p>1$. Indeed an individual density error with $L^1$ norm $CR_a^{c_a}|\log R_a|^N$ and supremum $CR_a^{-K_a}$ has $L^p$ norm at most $CR_a^{[c_a-K_a(p-1)]/p}|\log R_a|^N$.

Choose $p$ close to one for the finite list, and only then sum; no smallest radius is substituted. Choose the cut-off of each inner concentration to be one on $D_\tau$ and supported in $D_{2\tau}$. Its retained actual tail on the intervening annulus has $L^p$ norm bounded by
\[
 C r^{2t-\varepsilon}\tau^{\,2/p-2t-2+\varepsilon},\qquad
 \tau=R^\theta .
\]
Indeed the pointwise bound there is
\[
 Cr^{2t-\varepsilon}d^{-2t-2+\varepsilon}
\]
and its radial integral starts at $\tau$. Here $1<\theta<s/t$ and $r=R^{s/t+o(1)}$ by the local profile estimates. At a mixed regular point take $s=2$, $t=1$, and $\theta=3/2$. At $(p,\varepsilon)=(1,0)$ the exponent $(s/t)(2t-\varepsilon)+\theta(2/p-2t-2+\varepsilon)$ equals
\[
 2s-2t\theta>0.
\]
For the finite list choose $1<p<3/2$, close enough to one also for the interpolation above, so these exponents are positive at $\varepsilon=0$. Then choose $\varepsilon>0$ small enough to preserve positivity,
\[
 p(2t+2-\varepsilon)>2,\qquad p\varepsilon<2,
\]
and $p(2A-2-\varepsilon)>-2$. The exterior bounds $Cd^{2A-2-\varepsilon}$ and $Cd^{-\varepsilon}$ then have uniformly integrable $p$-th powers. At an ordinary point use the quadratic-zero bound of Lemma~\ref{reg:local-power}. At a pure prescribed point use the pure outer estimate with
\[
 2B-2>0.
\]
Strong $L^p$ convergence of the retained density follows by first using convergence outside fixed disks and then letting their radii shrink, using these power bounds inside the disks. The fitted error $bf_2-\nu_2$ is already small by the interpolation above. The exact remainder equation
\[
 -\Delta_g w_u=2\Pi(f_2-\nu_2)-\Pi(f_1-\nu_1),\qquad
 \Pi v=v-\int_M v\,dV_g,
\]
then gives \(w_u\to w_*\) in \(W^{2,p}\). This supplies the scalar inverse-function neighborhood. For the quantitative coefficient comparison, let \(b\) be the sum of the inner cut-offs, each supported on the corresponding intermediate disk, and set
\[
 A_D=(1-b)h_2e^{Z_2},\qquad C_u=\frac{\rho_2}{\int_M h_2e^{u_2}}.
\]
The already established independent-radius interpolation gives \(\|bf_2-\nu_2\|_p\le CD^c|\log D|^N\), with \(p>1\) close enough to one. In particular
\[
 f_2-\nu_2=C_u A_De^{w_u}+bf_2-\nu_2.
\]
In comparing \(A_D\) with \(H_p\), extend each truncated fitted inner model to the fixed disk. By the tail integral in the proof of Lemma~\ref{hyb:capture-matching} this changes the estimate by at most $C(r/\tau)^{2t}\tau^{-(2-2/\widehat p)}$ in $L^{\widehat p}$ for $1<\widehat p\le p$. Choose $\widehat p$ close enough to one using
\[
 \theta<s/t,
\]
so that the error is bounded by a positive power of this outer radius. The mean-zero Green inverse bounds the corresponding $C^0$ error. The fixed cut-offs in the exterior models give a bound of the same form.

The regular part of $G$ and its first spatial derivative can be replaced by their values at the center using the fractional moment estimate in the same lemma. These comparisons are done at each individual radius before summing.

On a mixed singular disk the radial Green calculation
\[
 \Delta\log(d^{2q}+R^{2q})=m_{q,R},
\]
with \(d=|z|\), now gives the factor multiplying \(H_p\), up to a multiplicative \(1+O(D^c|\log D|^N)\), as
\[
 \frac{1+(R/d)^{2s}}{(1+(r/d)^{2t})^2}.
\]
Here the signed Green combination is \(-Gm_{s,R}+2Gm_{t,r}\). The limiting coefficient \(H_p\) is a bounded factor times \(d^{2A-2}\).

On \(d\ge\tau\) subtract one from the displayed ratio and use that the absolute difference is at most \((R/d)^{2s}+2(r/d)^{2t}\). After multiplication by \(H_p\) the outer term has the \(L^p\) majorant \(C R^{2s}d^{-2t-2}\). Its norm is bounded by
\[
 C R^{2s}\tau^{2/p-2t-2},
\]
a positive power of \(R\): at \(p=1\) the exponent is
\[
 2s-2t\theta>0
\]
since
\[
 \tau=R^\theta,\qquad \theta<s/t.
\]
Choose \(p\) close enough to one for this strict inequality at every selected mixed singular point. The remaining term is bounded by \(C(r/\tau)^{2t}\). In the removed disk the limiting coefficient has norm $O(\tau^{2A-2+2/p})$.

In the transition, write the coefficient difference as $(1-b)$ times the full coefficient difference minus $bH_p$. The former is estimated from $d\ge\tau$ and the latter by the same disk bound with $2\tau$ in place of $\tau$.

At a pure singular point the factor is \(1+(R/d)^{2A}\), and multiplication by the limiting
\[
 H_p=O(d^{2(A+B-1)})
\]
makes its difference \(O(R^{2A})\) because \(B>1\). At an ordinary point the corresponding quadratic-zero cancellation is
\[
 d^2\{1+(R/d)^2\}-d^2=R^2.
\]
At a mixed regular point the coefficient comparison in \eqref{reg:actual-diffuse} bounds the error on $d\ge\tau$ by
\[
 C R^{3/p-2}|\log R|^N,\qquad 3/p-2>0.
\]
All other weights and remote factors are bounded on their matching domains. Together these give
\[
 \|A_D-H_p\|_p\le CD^c|\log D|^N,
\]
after decreasing \(c>0\).

Since \(w_u\) is already uniformly bounded and \(A_D\) is uniformly positive on a fixed open set away from the centers and the prescribed singularities, integration of the exact decomposition first bounds \(C_u\) above: the remaining term \(bf_2-\nu_2\) already has small \(L^1\) norm. It follows that
\[
 f_2-\nu_2=C_uH_pe^{w_u}+e_2,\qquad
 \|e_2\|_p\le CD^c|\log D|^N .
\]
Set $I_u=\int_M H_pe^{w_u}$, $Q_u=H_pe^{w_u}/I_u$, and
\[
 \delta_2=\int_M\nu_2-(\rho_2-\gamma).
\]
Integration gives $C_u I_u=\gamma-\delta_2-\int_M e_2$ and hence a positive lower bound for $C_u$, since $I_u$ is bounded above and below. The defect in the mean-zero scalar equation for $w_u$ is exactly
\[
 \Pi\bigl(2e_2-(f_1-\nu_1)-2(\delta_2+\textstyle\int_M e_2)Q_u\bigr).
\]
Its $L^p$ norm satisfies the same estimate because $H_p$ is bounded and $w_u$ is uniformly bounded. The scalar linearization from mean-zero $W^{2,p}$ to mean-zero $L^p$ is a compact perturbation of the mean-zero Laplacian.

Its kernel space is zero by the assumed nondegeneracy in $H^1$, and its inverse is uniform by continuity on the compact selected branch.

The normalized exponential has a quadratic Taylor remainder from $W^{2,p}$ to $L^p$ near this branch. The established $W^{2,p}$ convergence permits absorption of this remainder, yielding
\[
 \|w_u-w(p)\|_{W^{2,p}}\le CD^c|\log D|^N.
\]

\noindent\textit{Step 2. The projected family.}\quad The bounds for the individual concentrating profiles and their remaining tails, and the scalar comparison just established, are the hypotheses of Lemma~\ref{hyb:capture-matching}. It gives the actual common normalization, $|\mu|\leq Ce^{-c_0T}$ after decreasing
\[
 c_0>0,\qquad \|u-B_T\|_{C^0\cap H^1}\leq Ce^{-c_0T}T^N.
\]
Apply Lemma~\ref{journal:modulation-uniqueness} in
\[
 X=(\dot C^0\cap\dot H^1)^2
\]
with \(\theta\) the regular centers scaled by their own translation widths, the shape coordinates, and \(T\), and with
\[
 L_\theta=(L_\alpha).
\]
By \eqref{hyb:modulation-bounds} and the explicit cut-off profiles, the scaled parameter derivatives are uniformly bounded in $X$. At each active scale their values are bounded and their rescaled gradients are square-integrable.

The other disks contribute at most positive powers of their radii. The translation and shape diagonal blocks are $(4\pi/3)I$ and $(32\pi/3)I$, by the ordinary and quadratic profile moments.

At a mixed regular point the inner dilation from differentiating $r(B)$ is even and pairs to zero against the inner translations. Since
\[
 \partial_T\log R_a=-1/(2s_a),
\]
the first component of the $T$ tangent on the singular reference scale is $Z_{s_0,R_0}+1$ up to a constant and an error vanishing in the local weighted pairing; its centered radial pairing tends to $4\pi s_0/3$.

On the mixed regular outer scale the $T$ tangent is $Z_{0,B}+1-\tfrac12 B\cdot\partial_B\log m_{R,B}$ up to a constant and vanishing errors. In fact
\[
 \partial_T\log R=-1/4,
\]
and differentiation of $-2\log(R^4+|z^2-R^2B|^2)$ at fixed $z$ and $B$ gives this expression. The centered shape pairing eliminates the constant and the dilation, so the shape--$T$ column tends to $-(16\pi/3)B$.

All other off-diagonal limits vanish by parity, separation or the inner-to-outer ratio.

The resulting triangular matrix has a uniformly bounded inverse on the compact parameter set. In a fixed small ball of the initial scaled parameter coordinates, the rescaled rational profiles and the source pairings vary uniformly continuously.

On fixed rescaled compact sets differentiate using the bounded ratios of the varied and initial widths. Profile decay bounds the weighted tails uniformly. After subtraction of constants in the centered pairings the remote errors tend uniformly to zero.

Together with \eqref{hyb:modulation-bounds} this gives the required uniform control of the variation of $L_\theta D_\theta B_\theta$ with respect to $\theta$. The scaled parameter correction and the resulting remainder are $O(e^{-c_0T}(1+T)^N)$. The change in $T$ is $o(1)$, so the same bound holds for the parameters satisfying
\[
 L_\theta(u-B_\theta)=0.
\]
Choose the construction exponent $\delta_0$ smaller than the finitely many positive exponents in the preceding estimates. The parameters and remainder of the solution lie in the domain and uniqueness ball of Proposition~\ref{hyb:range}. Indeed, with
\[
 0<\delta_0<\min\{c_0,\bar\delta\},
\]
the ratio of the preliminary error bound to the ball radius is bounded by
\[
 Ce^{-(c_0-\delta_0)T}T^{N-C_0}\to0;
\]
the parameter ratio tends to zero as well. This is the elementary comparison in Lemma~\ref{journal:strict-capture}, and does not require the sharp radial mass law. The correction solves the projected equation in the uniqueness ball and hence agrees with the one constructed there.

The inner product estimates in the proof of invertibility give linear independence of the sources $S_\alpha$. An original solution has zero residual, hence all its multipliers vanish.

Conversely, when these multipliers vanish, the projected equation has zero full residual, with the same integral normalizations. Subsequence extraction and the finite selected set of shadow solutions give the assertion uniformly over the selected shadow solutions.
\end{proof}
\subsection{The scale equation and the local degree}
Let $\mathcal J_\mu$ be the normalized Toda functional and $j(T,z,B,\mu)=\mathcal J_\mu(B_T+\eta)$. At fixed $\mu$, centers and shapes put $Y=\partial_T B_T$. For $|\mu|\le e^{-\delta_0T}$ one has
\begin{equation}\label{hyb:energy-comparison}
 j_T-\partial_T\mathcal J_\mu(B_T)
       =O((\epsilon(T)^2+\mu^2)T^C).
\end{equation}
Indeed, differentiation of $L_\alpha\eta=0$ gives
\[
 j_T=\langle\mathcal F_\mu(B_T+\eta),Y\rangle_{\Hh}
        -\sum_\alpha c_\alpha(\partial_TL_\alpha)\eta .
\]
By self-adjointness the linear part of the difference of the first pairing from $\partial_T\mathcal J_\mu(B_T)$ is $\langle\eta,\partial_T\mathcal F_\mu(B_T)\rangle_{\Hh}$. The radial width derivatives of the patches are uniformly bounded. For the mixed regular patch, at fixed $B,\beta$, let
\[
 L=Q_0/(4R^4).
\]
Its two tangents are $\partial_T\log L-2\partial_T\log F$ and $\partial_T\log F-2\partial_T\log L$, each bounded by the positive comparability of the denominators.

The fixed cut-offs and surface mean subtraction preserve the bound on $Y$ in $C^0$. The probabilities $N_i/M_i$ have mass one and exponential $L^{p_0}$ bounds for a fixed
\[
 p_0>1.
\]
Proposition~\ref{hyb:range} gives
\[
 \sqrt T\|\eta\|_{\Hh}=o(1),
\]
so Lemma~\ref{journal:normalized-nonlinearity}, with $Y$ as the bounded test, bounds the Taylor error by $CT\|\eta\|_{\Hh}^2$. The differentiated residual bounds the linear term by
\[
 C(\epsilon(T)+|\mu|)T^C\|\eta\|_{\Hh}.
\]
The last sum is at most $C|c|\|\eta\|_{\Hh}$ by \eqref{hyb:modulation-bounds}. The correction and multiplier bounds in Proposition~\ref{hyb:range} now yield \eqref{hyb:energy-comparison}.
\begin{lemma}[Expansion of the energy derivative]
\label{hyb:trial-energy}
With $J_q=4\pi^2/\sin(\pi/q)$, the corrected approximate solution satisfies
\begin{equation}\label{hyb:trial-radial}
 \partial_T\mathcal J_\mu(B_T)
 =-\mu+\sum_{a\in P\cup Q}\frac{\kappa_aJ_{s_a}}{4s_a}R_a^2
    +\sum_{b\ {\rm mixed\ regular}}R_b^2
       \{\kappa_b\mathcal M(B_b)+4\pi t_b^{\rm sh}\cdot B_b\}
    +o(\mathcal L)+O(|\mu|\mathcal L T^C).
\end{equation}
Only terms of order $\mathcal L$ are required.
\end{lemma}
\begin{proof}
\noindent\textit{Step 1. Differentiation of the normalized functional.}\quad We differentiate the normalized functional before expanding any local term.

The derivative thus includes its spatial constants and both denominators exactly. First make the matching coefficients explicit. At every outer bubble,
\[
 d_a^{s_a}=\frac{e^{\psi_{1,a}(0)}}{4s_a^2}.
\]
At a mixed singular bubble define the regular function
\[
 \psi_{2,a}=\log h_2+\varphi+u^0_2
                   -(2(s_a-t_a)-2)\log|z|.
\]
Then $b_a=e^{c_2+\psi_{2,a}(0)}/(4t_a^2)$. At a mixed regular bubble the same definition has no logarithm, and
\[
 \beta_b=e^{c_2+\psi_{2,b}(0)},\qquad
 R_b^4=\frac{e^{\psi_{1,b}(0)}}{16}e^{-T},\qquad
 r_b^2=\frac{\beta_b(1+|B_b|^2)}4R_b^4.
\]
In particular these matching coefficients are independent of $T$ when the center and shape parameters are fixed. All the remote masses and the scalar shadow are already included in the functions $\psi_{i,a}$. Write
\[
 B_T=\widehat B_T-\int_M\widehat B_T,\qquad N_i=e^{c_i}h_i e^{\widehat B_{T,i}},
\]
as in the proof of Proposition~\ref{hyb:trial-comparison}. Put $M_i=\int_MN_i$, $\rho^0=(M,\rho_2)$, $N=(N_1,N_2)^t$, and
\[
 \mathscr R_T=-\Delta_g\widehat B_T-K(N-\rho^0).
\]
The functional on this mean-zero approximate solution is exactly
\[
 \mathcal J_\mu(B_T)
 =\tfrac12\|\widehat B_T\|_{\Hh}^2
       -\sum_i\rho_i\log M_i+\sum_i\rho_i c_i
       +\sum_i\rho_i\int_M\widehat B_{T,i}.
\]
Put $\widehat Y=\partial_T\widehat B_T$. The derivative of the Dirichlet term is $\int_M(N-\rho^0)\cdot\widehat Y+ \int_M K^{-1}\mathscr R_T\cdot\widehat Y$. Using
\[
 \int_MN_i\widehat Y_i=M_i'-c_i'M_i,
\]
the terms with $c_i'$ add up to $\sum_i(\rho_i-M_i)c_i'$ and the terms with $M_i'$ to $\sum_i(1-\rho_i/M_i)M_i'$. Since $c_1'=-1$ and $c_2'=0$, the first sum is $M_1-\rho_1=-\mu+(M_1-M)$. The mean terms leave $\sum_i(\rho_i-\rho_i^0)\int_M\widehat Y_i$. This gives the exact identity
\begin{align}
 \partial_T\mathcal J_\mu(B_T)
 ={}&-\mu+(M_1-M)
       +\int_M K^{-1}\mathscr R_T\cdot\partial_T\widehat B_T
       \notag\\
 &+\sum_i\frac{M_i-\rho_i}{M_i}M_i'
       +\mu\int_M\partial_T\widehat B_{T,1}.
 \label{hyb:exact-trial-variation}
\end{align}
Here no local mass has been normalized separately. We estimate the terms in this identity. For a fixed small coordinate disk, uniformly in the compact parameters,
\begin{equation}\label{hyb:trial-mass-order}
 |M_i-\rho_i^0|+|M_i'|\le C\mathcal L T^C,\qquad
 \int_M|\partial_T\widehat B_T|\le C\mathcal L T^C.
\end{equation}
These improve the residual norm bound only for the integrated masses.

At a singular bubble the linear terms are odd, the quadratic outer moments are $O(R^2)$, and inner moments and overlaps are respectively $O(r^2T^C)$ and $O(R^{2A}T^C)$ at a mixed singular bubble. At an ordinary bubble the quadratic moment is
\[
 O(R^2|\log R|).
\]
For the profile \eqref{hyb:gram}, $n_1=\Delta\log F$ and
\[
 n_2=\Delta\log Q_0.
\]
Both polynomials grow as positive multiples of $|z|^4$, so the divergence theorem gives mass $8\pi$ for each density. The profile
\[
 n_{\rm bg}=\beta(1+\beta|z|^2/4)^{-2}
\]
has mass $4\pi$ by radial integration. Its actual smooth factor in the background is $e^{J_2}$, which is also the factor multiplying $n_2$ in the trial density.

Therefore the concentrating mass of the second component is computed from $(n_2-n_{\rm bg})e^{J_2}$. The second moment of the first density is $O(R^2)$.

The second-density linear term vanishes by evenness, and the absolute second moment of $n_2-n_{\rm bg}$ on a fixed disk is $O(R^2+r^2(1+|\log r|))$ by the bounds at the two scales in Proposition~\ref{hyb:trial-comparison}.

The same division of the rational functions shows that their omitted tails outside a fixed disk are $O_d(R^4)$ in the first component and $O_d(R^2)$ for $n_2-n_{\rm bg}$.

The fixed annuli contribute $O(e^{-T}T^C)$ in the first component and at most $O(R^2)$ in a regular mixed second component.

The same bounds hold after one $T$ derivative: the full-plane masses are constant, parity is unchanged, and differentiating the density estimates in the fixed $z,B,\beta$ variables preserves the indicated majorants.

Finally the tails of the tangent of the approximate solution are $O(R^{2s}|z|^{-2s})$ for a radial bubble and $O(R^2|z|^{-2})$ for a mixed regular bubble. Integration proves the last estimate in \eqref{hyb:trial-mass-order}.

Since all outer profile exponents are at most $s_*$, these observations prove the displayed common bound. Consequently the last line of \eqref{hyb:exact-trial-variation} is $o(\mathcal L)+O(|\mu|\mathcal L T^C)$.

\noindent\textit{Step 2. The singular points.}\quad Consider a singular outer bubble with profile exponent $q$. Put
\[
 f_q(y)=1-\frac{2q}{(q-1)(1+|y|^{2q})},\qquad
 H_2(x)=\tfrac12x^tD^2\psi_1(0)x.
\]
Thus $g\cdot x+P=(g\cdot x)f_q(x/R)$. After applying $K^{-1}$ to the residual, its first component contains $-m_{q,R}(e^{J_1+P}-1-P-g\cdot x)$, apart from the inner and terms from the residual density estimated below. Here $J_1=\psi_1-\psi_1(0)$. Its quadratic term is $-E_1$, where
\[
 E_1=m_{q,R}\{H_2(x)+\tfrac12(g\cdot x)^2f_q(x/R)^2\}.
\]
In the first mass defect $M_1-M$, the same quadratic term occurs with sign $+E_1$: the terms $m_{q,R}P$ and $m_{q,R}g\cdot x$ are odd on the fixed disk, and the model mass omitted outside it is $O(R^{2q})$.

At fixed centers and shapes the leading tangent is $(Z_{q,R}+1,-(Z_{q,R}+1)/2)$. Consequently the first mass defect and its pairing with the residual multiplied by $K^{-1}$ combine as
\[
 \int E_1-\int (Z_{q,R}+1)E_1=-\int Z_{q,R}E_1.
\]
The contribution to the first line of \eqref{hyb:exact-trial-variation} is therefore
\begin{equation}\label{hyb:pinned-radial-moment}
 -\int m_{q,R}Z_{q,R}
       \{H_2(x)+\tfrac12(g\cdot x)^2f_q(x/R)^2\}+o(R^2).
\end{equation}
To justify the error, use the uniform $C^{2,\alpha}$ background bounds. The Taylor remainders and the cubic first-order terms are bounded by $Cm_{q,R}|z|^{2+\delta}$ with
\[
 0<\delta<\min(\alpha,1,2q-2)
\]
on the fixed disk. A first-order correction tangent is $O(|z|)$ and obeys the same bound against the quadratic term. Extending the quadratic radial integral to the plane discards only an $O(R^{2q})$ tail.

At a mixed singular bubble the two linear coupling terms in \eqref{hyb:linear-couplings}--\eqref{hyb:linear-coupling-bounds} are odd against every radial part of the tangent. Against a first-order tangent they have size $O(R^{2A+2})$ by \eqref{hyb:linear-coupling-moments}.

The terms involving $Q_1$ in the first Taylor square are bounded by $C|z|S|Q_1|$ and have the same integral bound. The additional inner radial part of the first tangent has absolute value at most $C\min(1,h)$ with
\[
 h=(r/|z|)^{2t}.
\]
Against the outer quadratic term the resulting error is at most
\[
 C\int |z|^2Sh=O(R^{2A+2})
\]
by the pointwise majorant for $Sh$ in the trial residual. The remaining errors from the two concentrating densities are bounded by
\[
 C\{\int I_0|z|^2+\int Sh\}=O(r^2+R^{2A}).
\]
For the term from the residual density the integral of $g_{\rm loc}\min(1,h)$ is at most $C(r^{2A}+r^{2t})(1+|\log r|)$ by splitting at $r$. For the angular terms from the residual density, splitting at $R$ gives
\[
 \int g_{\rm loc}|P|=O(R^{2A+1})
\]
: outside $R$ the integrand is bounded by $CR^{2s}d^{-2t-1}$ and $t>1$. Splitting at $r$ gives
\[
 \int g_{\rm loc}|Q_1| \le C(r^{2A+1}+r^{2t})(1+|\log r|),
\]
including the equality of the two radial exponents. Together these contributions are $O((r^2+R^{2A}+e^{-T})T^C)$. At a pure prescribed source the residual density is $O(|z|^{2(A+B)-2})$. Multiplication by the outer radial patch difference or its first-order correction has integral
\[
 O(R^{2A})=O(e^{-T}).
\]
All of these bounds persist after one logarithmic-scale derivative at the fixed parameters.

The remainder estimates can be made with a single positive exponent. Choose
\[
 0<\delta<\min_{a\in P\cup Q}
       \{\alpha,1,2s_a-2\}.
\]
Lemma~\ref{hyb:moment-estimates}, with $p=1$, bounds a Taylor remainder $O(|z|^{2+\alpha})$ and a cubic correction by
\[
 O(R_a^{2+\delta}).
\]
Multiplication by a bounded radial tangent, or by any of its first two logarithmic derivatives, preserves this bound. At a mixed singular point, the additional moments are bounded by
\[
 C\{r_a^2+R_a^{2A_a}+R_a^{2A_a+2}\}(1+T)^C=o(R_a^2).
\]
The first two terms follow respectively from the inner second moment and the overlap estimate; the last is the first moment of the linear couplings paired with a first-order tangent. Their strict exponent inequalities use
\[
 s_a/t_a>1,\qquad A_a>1.
\]
Differentiation in $T$ preserves their majorants by the same lemma. Hence the error and its first logarithmic derivative are uniformly $o(R_a^2)$ on the compact parameter set, before the different outer radii are compared.

For $h=1/q\in(0,1)$, make the substitution $w=|y|^{2q}$. The weighted moment with $k$ extra powers of $(1+w)^{-1}$ is
\[
 4\pi q\int_0^\infty w^h(1+w)^{-k-2}\,dw
 =4\pi q\,\frac{\Gamma(1+h)\Gamma(k+1-h)}{\Gamma(k+2)}.
\]
The value at $k=0$ is $J_q$ by the reflection formula. Dividing the expression by $J_q$ gives $a_k=\prod_{j=1}^k(j-h)/(k+1)!$. In particular $a_0=1$ and $2a_1-a_0=-h$. With $c=2/(1-h)$, the normalized moment of $(2/(1+w)-1)(1-c/(1+w))^2$ is
\[
 (2a_1-a_0)-2c(2a_2-a_1)+c^2(2a_3-a_2)
 =-h-\tfrac23(1-2h)+\tfrac13(2-h)=0.
\]
This proves
\[
 \int m_{q,1}|y|^2=J_q,\qquad
 \int m_{q,1}Z_{q,1}|y|^2=-J_q/q,\qquad
 \int m_{q,1}Z_{q,1}|y|^2f_q(y)^2=0 .
\]
Angular averaging in \eqref{hyb:pinned-radial-moment} now gives $\kappa J_qR^2/(4q)$. In particular the entire $|g|^2$ term cancels, even when the first-order coefficient at a singular point is nonzero.

\noindent\textit{Step 3. The regular points.}\quad At a regular mixed bubble, set $n_{\rm bg}=\beta(1+\beta|z|^2/4)^{-2}$. The exact residual is \eqref{hyb:gram-residual}, while
\[
 \partial_T\widehat B_1=1+\partial_T\log n_1,\qquad
 \partial_T\widehat B_2=\partial_T\log n_2
\]
on the bubble disk. The mass identity for $n_1$ cancels its constant profile contribution in $M_1-M$. Hence its contribution to \eqref{hyb:exact-trial-variation} is
\begin{align*}
 &-\int (\partial_T n_1)(e^{J_1}-1)\\
 &\hspace{8mm}-\int(n_2-n_{\rm bg})(\partial_T\log n_2)(e^{J_2}-1)
       +O_d(R^4T^C).
\end{align*}
Here $J_1=\kappa|z|^2/4+\Re(T_bz^2)+o(|z|^2)$, and $J_2=\psi_2-\psi_2(0)+2\log(1+\beta|z|^2/4)$. Here are the estimates on the differentiated densities. Put $U=R^2$, $u=|z|^2$, $x=u+r^2$, and $A_0=U^2+|z^2-UB|^2$. Then
\[
 A_0\asymp U^2+u^2,\qquad F=A_0+\beta U^2u.
\]
In $n_1=4Q_0/F^2$ the numerator other than $16U^2u$ is $O(U^2 A_0)$, and the error from replacing $F$ by $A_0$ in the leading term is $O(U^4u^2/A_0^3)$. Both give
\[
 |n_1-m_{R,B}|\le CU^2/A_0,
\]
with the same bound after two logarithmic derivatives at fixed $B,\beta,z$. On $u\le U$ use $Q_0\asymp U^2x$ and $F\asymp U^2$, so $n_2+n_{\rm bg}\le Cr^2/x^2$ and $|\partial_T^j\log n_2|\le C$ for $j=1,2$. On $u\ge U$ use
\[
 F/u^2=1+O(U/u),\qquad
 Q_0/(4U^2u(1+\beta u/4))=1+O(U/u).
\]
The errors and their first two $T$ derivatives satisfy the same bounds: each is a rational expression in $U,u,z^2,B,\beta$ with the denominator bounded below relative to its displayed scale. Thus
\[
 n_2-n_{\rm bg}=O(n_{\rm bg}U/u),\qquad \partial_T\log n_2=O(U/u)
\]
on this region, with the same bounds for one further derivative. We obtain
\[
 |(n_2-n_{\rm bg})\partial_T\log n_2|
 \le C\left\{\frac{r^2}{(r^2+|z|^2)^2}
             +\frac{R^4}{(R^2+|z|^2)^2}\right\}
\]
and the same majorant for its first $T$ derivative. The second moments of this bound and of $CU^2/A_0$ are $O_d(R^4(1+|\log R|))$ by polar integration in $u$.

All factors in the second integrand other than $e^{J_2}-1$ are even; its linear part therefore integrates to zero. This proves that the second integral is $O_d(R^4(1+|\log R|))$. The first integral may be replaced by $-\int (\partial_Tm_{R,B})J_1$ with an $o(R^2)$ error.

At fixed $B$, $R_T/R=-1/4$, and
\[
 V_B:=\partial_T\log m_{R,B}
   =Z_{0,B}(z/R)-\tfrac12 B\cdot Z_{{\rm sh},B}(z/R).
\]
The second moments are absolutely integrable. Differentiating them at fixed $B$ gives
\[
 \int m_{R,B}V_B|z|^2=-4R^2\mathcal M(B),\qquad
 \int m_{R,B}V_Bz^2=-4\pi R^2B.
\]
Thus the mixed contribution is $R^2\{\kappa\mathcal M(B)+4\pi\Re(T_bB)\}+o(R^2)$, as claimed. At an ordinary bubble the linear Taylor term again integrates to zero against its radial tangent, and the whole contribution is $O(Te^{-T})$.

On every fixed cut-off annulus the residual and the tangent of the approximate solution are both $O(R^2)$ for a mixed regular bubble, and both $O(R^{2q})$ for a radial bubble; their pairing is therefore $o(\mathcal L)$. The omitted tails of the first density have mass $O(e^{-T})$.

This accounts for the approximate solution, including all cut-off annuli. The estimates are uniform on the compact parameter set, and their differentiated versions follow from the same majorants. Summing \eqref{hyb:pinned-radial-moment} and the mixed contributions in \eqref{hyb:exact-trial-variation} proves \eqref{hyb:trial-radial}.
\end{proof}
Order the scaled parameters as centers, shapes and $T$. Scale each center by its own translation width. The pairing matrix between the sources and the parameter derivatives of the approximate solution tends to
\begin{equation}\label{hyb:chart}
 J_0=\begin{pmatrix}
 (4\pi/3)I&0&0\\
 0&(32\pi/3)I&-(16\pi/3)\mathbf B\\
 0&0&4\pi s_0/3
 \end{pmatrix},\qquad \det J_0>0.
\end{equation}
For the center diagonal blocks, the scaled test against a translation tangent is
\[
 \int m_{1,1}\,4y_i y_j/(1+|y|^2)^2=(4\pi/3)\delta_{ij}.
\]
On the mixed outer scale the shape sources pair with the first component in the variable $w=y^2-B$. The resulting inner products are
\[
 128\int w_iw_j(1+|w|^2)^{-4}\,dA_w=(32\pi/3)\delta_{ij}.
\]
The $T$ tangent in that first component is $1+Z_{0,B}-B\cdot Z_{\rm sh}/2$. Its constant and radial terms have zero limiting shape pairing by centering and parity in $w$. This gives the displayed shape--$T$ column. The reference diagonal is
\[
 \int m_{s_0,1}Z_{s_0,1}(Z_{s_0,1}+1) =4\pi s_0/3.
\]
A possible nested inner tangent tends to zero in the reference outer weighted integral because $r/R\to0$ and the outer measure has no point mass at zero.

The other off-diagonal limits are zero by separated supports and their centering, or by parity and scale separation on a mixed disk. The corresponding contribution of the outer shape derivative to the inner mean is $2 B_j/(1+|B|^2)$. The matched inner dilation contributes
\[
 -2\,\partial_{B_j}\log r=-2 B_j/(1+|B|^2)
\]
since $\partial_{\log r}D_{1,r}$ averages to $-2$ in the normalized inner measure. Thus those means cancel as well. The parameterization of the corrected solutions has the same orientation since
\[
 L_\alpha\partial_p\eta=-(\partial_pL_\alpha)\eta=o(1)
\]
in these scaled coordinates. The uniform source bounds are \eqref{hyb:modulation-bounds}, and the bounded scaled tangents follow by differentiating the explicit profiles as in the choice of parameters in Theorem~\ref{hyb:capture}.

For the shape equations we use the uniform $C^0$ convergence proved above. The derivative estimates apply to the parameterization.

Let $c_0$ be the radial multiplier at the singular reference point and $e_0=L_0u_T\to4\pi s_0/3>0$.
\begin{proposition}[The reduced scale equation]
\label{hyb:radial-row}
If $c^{\rm ctr}=0$, then for all shape parameters under consideration,
\begin{align}
 e_0c_0={}&-\mu+
 \sum_{a\in P\cup Q}\frac{\kappa_aJ_{s_a}}{4s_a}R_a^2
 +\sum_{b\ {\rm mixed\ regular}}\kappa_bR_b^2\mathcal H(B_b)
 \notag\\[-2pt]
 &\hspace{12mm}+o(\mathcal L)+O(\mu^2T^C+|\mu|\mathcal L T^C),
 \qquad \mathcal H(B)=\mathcal M(B)-B\cdot\nabla\mathcal M(B).
 \label{hyb:actual-radial}
\end{align}
The same formula holds at $\mathbf B=0$ for arbitrary center parameters.
\end{proposition}
\begin{proof}
At $c^{\rm ctr}=0$, the chain rule along the family of projected solutions gives
\[
 j_T=e_0c_0+\sum_b c_b^{\rm sh}\cdot L_b^{\rm sh}u_T .
\]
On this entire set, \eqref{hyb:shape-row} gives $c_b^{\rm sh}=-(3/(16\pi))R_b^2\Phi_b(B_b)+o(R_b^2)$ and \eqref{hyb:chart} gives
\[
 L_b^{\rm sh}u_T=-(16\pi/3)B_b+o(1).
\]
Their product is $R_b^2B_b\cdot\Phi_b(B_b)+o(R_b^2)$. Subtract it from \eqref{hyb:trial-radial} and use \eqref{hyb:energy-comparison}. The terms $4\pi R_b^2t_b^{\rm sh}\cdot B_b$ cancel, and the $\epsilon(T)^2T^C$ term is $o(\mathcal L)$. The two terms involving $\mu$ retain the orders in \eqref{hyb:actual-radial}.

At $\mathbf B=0$ the $T$ derivative of the approximate solution is radial on every regular source support. Its pairings with odd translations and shape derivatives of angular frequency two vanish exactly; the singular corrections have disjoint supports and constants pair to zero. Using
\[
 L_\alpha\partial_T\eta=-(\partial_TL_\alpha)\eta
\]
for $\alpha\ne0$ now gives
\[
 j_T=e_0c_0+O((\epsilon(T)^2+\mu^2)T^C),
\]
without using the equations for the centers. Since $\mathcal H(0)=\mathcal M(0)$, the asserted formula follows in this case as well.
\end{proof}
At occupied prescribed points the effective weight parameter $\beta_p$ is positive, so the residual density vanishes there. In \eqref{hyb:canonical-trace} and in the prescribed patches, the normalization of the coordinate has
\[
 \varphi(0)=0.
\]
For a prescribed point, the smooth part of the weight contributes $\Delta_g\log h_1^*(p_a)-4\pi\sum_p\alpha_{p,1}$, the conformal factor contributes $-2K_g(p_a)$, and the first Green field contributes $2M-N$ after the pole is removed.

The scalar part $-w/2$ contributes $-\gamma$ at the prescribed point. Since $N+\gamma=\rho_2$, we have
\[
 \kappa_a=\Delta_g\log h_1^*(p_a)-2K_g(p_a)
 -4\pi\sum_{p\in S}\alpha_{p,1}+2M-\rho_2.
\]
The selected weights make these traces positive and give $\kappa_b>4|T_b|$ at mixed regular points. To compute $\mathcal H(B)$, write $f_B(v)=|v+B|$. Differentiation under the integral is justified by the bounded gradient of $f_B$. Then
\[
 \mathcal H(B)
 =\int_{\C}\frac{v\cdot\nabla f_B(v)}{(1+|v|^2)^2}\,dA_v
 =\frac12\int_\C\frac{dA_v}{(1+|v|^2)|v+B|}>0 .
\]
In the integration by parts use
\[
 \Delta f_B=1/|v+B|
\]
off $-B$. The deleted circle about $-B$ contributes
\[
 O(\varepsilon),
\]
and the circle of radius $L$ at infinity contributes $O(L^{-1})$. The integrals are finite and continuous on compact $B$-sets. Shrink the center neighborhoods so that the strict trace inequalities persist throughout their closures. Since
\[
 R_a^2=d_a e^{-T/s_a},
\]
the coefficient of $\mathcal L$ in \eqref{hyb:actual-radial} is
\[
 \Gamma=\sum_{\substack{a\in P\cup Q\\s_a=s_*}}
             \frac{d_a\kappa_a J_{s_a}}{4s_a}
        +\mathbf1_{\{s_*=2\}}
          \sum_{b\ {\rm mixed\ regular}}d_b\kappa_b\mathcal H(B_b).
\]
A term with $s_a<s_*$ is exponentially smaller. At least one term with $s_a=s_*$ is present, and every such term is strictly positive. Continuity and compactness hence give
\[
 0<\gamma_0\le\Gamma\le K_0
\]
uniformly. At a solution of the full equation in the projected family, all multipliers vanish. Hence $j_T=0$. Using only \eqref{hyb:energy-comparison} and \eqref{hyb:trial-radial} first gives
\[
 |\mu|\le C\mathcal L+
          C(\mu^2+|\mu|\mathcal L)T^C+o(\mathcal L).
\]
In the enlarged domain $|\mu|\le e^{-\delta_0T}$ the coefficient $C(|\mu|+\mathcal L)T^C$ tends to zero, so absorption gives $|\mu|=O(\mathcal L)$. The two parameter errors in \eqref{hyb:actual-radial} are then $o(\mathcal L)$. Since $c_0=0$,
\begin{equation}\label{hyb:mass-law}
 \mu/\mathcal L=\Gamma+o(1).
\end{equation}
Theorem~\ref{hyb:capture} gives $|\mu|\le e^{-\delta_0T}$. The scale equation then gives $0<c\le\mu/\mathcal L\le C$ for every sufficiently concentrated solution. Its widths are comparable to the matched radii by the same construction.
\begin{lemma}[Degree of the reduced map]
\label{hyb:homotopy}
Let $\Omega\subset\R^{2k}$ be a bounded domain for the centers. Take a product of large two-dimensional balls for the shape parameters and $(T_-,T_+)$ for the scale. Assume $(C,Q,S)$ is continuous on the closure of this product domain, and suppose
\begin{itemize}
\item[\textup{(i)}]
$F_0\in C(\overline\Omega;\R^{2k})$, $0\notin F_0(\partial\Omega)$ and $\sup|C-F_0|<\min_{\partial\Omega}|F_0|$.
\item[\textup{(ii)}]
Whenever $C=0$, one has $Q_a\cdot B_a>0$ on each shape boundary sphere, $S>0$ at $T_-$, and $S<0$ at $T_+$.
\end{itemize}
Write $\deg_{\mathrm B}F_0:=\deg_{\mathrm B}(F_0,\Omega,0)$. The Brouwer degree of $(C,Q,S)$ on the product domain at zero is $-\deg_{\mathrm B}F_0$. In the absence of center variables, omit the center boundary condition and set $\deg_{\mathrm B}F_0=1$.
\end{lemma}
\begin{proof}
First homotope $S$ to $(T_-+T_+-2T)/(T_+-T_-)$ while keeping $C$ and $Q$ fixed. A zero on a scale endpoint would require $C=0$, where the two scale functions have the same strict sign. On a shape boundary the condition
\[
 Q_a\cdot B_a>0
\]
at $C=0$ excludes zeros, and on the center boundary the uniform distance to $F_0$ already excludes zeros. Keeping this scale function fixed, deform each $Q_a$ linearly to $B_a$. The shape condition at $C=0$ again excludes zeros on the shape boundary. Finally homotope $C$ to $F_0$.

At this step the shape and scale equations exclude zeros on their respective boundaries for every value of the center component. The three factors have degrees $\deg_{\mathrm B}F_0$, $1$, and $-1$.
\end{proof}
\begin{theorem}[The local degree contribution]
\label{hyb:transfer}
Fix a shadow solution $Z$ whose configuration contains at least one singular bubble. Assume the prescribed traces are positive and the strict trace inequalities hold at every mixed regular point.
\begin{itemize}
\item[\textup{(i)}]
As $\rho_1$ increases through $M$, this solution contributes
\begin{equation}\label{hyb:jump}
  \Delta d_Z=(-1)^{U+|P|}d_S(Z).
\end{equation}
\item[\textup{(ii)}]
For all sufficiently small $\mu\ne0$, there is a bounded isolating neighborhood containing all sufficiently concentrated solutions associated with this shadow solution.
\end{itemize}
\end{theorem}
\begin{proof}
\noindent\textit{Step 1. The finite-dimensional degree.}\quad On the positive side use $\mathcal L=\mu\tau$ with
\[
 0<\tau_-<1/K_0,\qquad \tau_+>1/\gamma_0.
\]
The increasing $T$ interval runs from $T_-=-s_*\log(\mu\tau_+)$ to $T_+=-s_*\log(\mu\tau_-)$. Choose $\delta_0$ as in Proposition~\ref{hyb:range}, also with $\delta_0<1/s_*$ and $\epsilon(T)=o(e^{-\delta_0T})$. On the closed scale interval,
\[
 \mu e^{\delta_0T}=\mu(\mu\tau)^{-s_*\delta_0}
 \leq \tau_-^{-s_*\delta_0}\mu^{1-s_*\delta_0}\to0 .
\]
Thus the projected equation is defined up to both endpoints. The bound $C(\epsilon(T)+|\mu|)T^C$ is strictly smaller than the contraction radius $e^{-\delta_0T}T^{C_0}$ for large $C_0$ and large $T$.

This excludes zeros on the complementary boundary over the whole interval. By \eqref{hyb:mass-law}, no sufficiently concentrated solution with this shadow limit has
\[
 \mu\le0.
\]
Take the local domain to be empty for $\mu<0$. For $\mu>0$ every sufficiently concentrated solution has $T_-<T<T_+$. For each shape parameter take the equation to be
\[
 -(16\pi/3)c_b^{\rm sh}/R_b^2=0.
\]
The multiplier is a negative scalar on a real two-dimensional space and thus preserves orientation. To choose the shape balls, differentiate $\mathcal M$ under the integral to obtain
\[
 \nabla\mathcal M(B)=\int_{\mathbb C}
   \frac{v+B}{|v+B|}\frac{dA_v}{(1+|v|^2)^2}.
\]
The value of the integrand at $v=-B$ is irrelevant. The integral on a fixed disk differs from the integral of $B/|B|$ against the same weight by $O(1/|B|)$, uniformly in the direction of $B$.

The remaining integral is bounded by twice the weight of the complement. The total weight is $\pi$. Since $\kappa_b-4|T_b|$ has a positive minimum on the closed center neighborhood, this gives
\[
 \Phi_b(B_b)\cdot B_b/|B_b|>c>0
\]
on sufficiently large fixed spheres. Where $C=0$, the normalized shape equation has left-hand side $\Phi_b(B_b)+o(1)$ by \eqref{hyb:shape-row}. The limiting shapes of the selected shadow solution are strictly inside the chosen balls.

For the center equations use $-c^{\rm ctr}/(3w_{\rm ctr})$, where $w_{\rm ctr}$ is its ordinary width or its mixed inner width. Proposition~\ref{hyb:stress-rows} gives
\[
 C=F_0+o(1)
\]
uniformly on the closed product, where $F_0$ is the center field evaluated on the selected scalar branch. If the center domain has a boundary, then $|F_0|$ is bounded away from zero there. Each negative normalization in a real two-dimensional block has positive determinant.

Set $S=e_0c_0$. Proposition~\ref{hyb:radial-row} applies at each point of the closed product satisfying $C=0$, with the shapes free to vary in their balls. Since
\[
 \mathcal L=\mu\tau,\qquad \mu T^C\to0
\]
uniformly, dividing \eqref{hyb:actual-radial} by $\mu$ gives $S/\mu=-1+\tau\Gamma+o(1)$ on that entire set. At $T_-$ this is at least
\[
 -1+\tau_+\gamma_0+o(1)>0.
\]
At $T_+$ it is at most $-1+\tau_-K_0+o(1)<0$. These are the scale conditions of Lemma~\ref{hyb:homotopy}, with the center and shape conditions just verified.

\noindent\textit{Step 2. The equation on the orthogonal complement.}\quad To justify the product with the complementary degree, take the center domain to be a sufficiently small coordinate ball about the center of the selected shadow solution, and write
\[
 p=(z,B,T).
\]
In the moving spaces $E_p$, write \(u=B_p+\eta_p+v\) with \(v\in E_p\) and \(\|v\|\leq \tfrac12 e^{-\delta_0T}T^{C_0}\). The graph of the projected solutions lies strictly inside the original uniqueness ball for small \(\mu>0\). Put
\[
 A_p=P_{E_p}D\mathcal F_\mu(B_p+\eta_p)|_{E_p}.
\]
The inverse estimate gives
\[
 \|A_p^{-1}\|\leq CT.
\]
On the stated fiber ball the projected equation has the form \(A_pv+N_p(v)\). By Lemma~\ref{journal:normalized-nonlinearity}, the third derivatives of the normalized energy are \(O(T^{3/2})\) in the energy norm along the segment from \(\eta_p\) to \(\eta_p+v\) relative to \(B_p\).

Indeed the probabilities at \(B_p\) have the fixed exponential \(L^{p_0}\) bound, and every point of this segment has \(\sqrt T\)-times its energy norm tending to zero. Taylor's formula based at \(\eta_p\) therefore gives
\[
 \|N_p(v)\|\leq CT^{3/2}\|v\|^2.
\]
For \(0\leq\sigma\leq1\), a zero of \(A_pv+\sigma N_p(v)\) there satisfies
\[
 \|v\|\leq CT^{5/2}\|v\|^2.
\]
Since
\[
 e^{-\delta_0T}T^{C_0+5/2}\to0,
\]
its only zero is \(v=0\). We may therefore set
\[
 \sigma=0
\]
and at the same time homotope the tangential components to their values at \(v=0\). All zeros during this homotopy are on the projected graph, where the boundary tests just proved apply.

For fixed $\mu$, the explicit formulas for the sources and their inner products show that the orthogonal projections depend smoothly on $p$ in operator norm.

Fix an interior parameter $p_*$ and let $U_p$ transport $E_{p_*}$ to $E_p$ along the straight segment. With $P(t)$ the orthogonal projection along the segment, solve
\[
 \dot U=[\dot P,P]U,\qquad U(0)=I.
\]
The generator is skew-adjoint. The identities
\[
 P\dot P P=0,\qquad [[\dot P,P],P]=\dot P
\]
show that $U(t)^{-1}P(t)U(t)$ is constant. The generator has finite rank because $\dot P$ has only off-diagonal blocks relative to $P$ and $\ker P$ is finite dimensional. The integral equation implies that $U_p-I$ is compact and norm-continuous in $p$.

In the fixed splitting spanned by the sources at $p_*$ and by $E_{p_*}$, write $\Psi(p,v)=B_p+\eta_p+U_pv$. Identify the finite coordinates with the ordered source basis at $p_*$. To express the equation in source coordinates use
\[
 D_p(c,v)=\sum_\alpha c_\alpha S_\alpha(p)+U_pv.
\]
Both $\Psi-I$ on bounded products and $D_p-I$ are compact, by compactness of the closed parameter product and norm continuity of $U_p-I$. The family $D_p$ is invertible and is joined to the identity through $D_{p_*+t(p-p_*)}$.

Multiplying the equation by $D_p^{-1}$ therefore preserves its degree in the fixed coordinates. For the source component of the derivative of $\Psi$, differentiating
\[
 L_p(\eta_p+U_pv)=0
\]
gives the finite block
\[
 G_p^{-1}\bigl(L_p\,d_pB_p-(d_pL_p)(\eta_p+U_pv)\bigr),
 \qquad (G_p)_{\alpha\beta}
       =\langle S_\alpha(p),S_\beta(p)\rangle_{\mathcal H}.
\]
Here $d_p$ denotes parameter differentiation. The other diagonal block of $D_p^{-1}d\Psi$ is the identity on $E_{p_*}$, and its upper right block vanishes. The matrix $G_p$ is positive definite.

In the scaled coordinates the second term in the finite block is $o(1)$ uniformly on the fiber ball, while the first has the positive determinant in \eqref{hyb:chart}.

Center widths are positive, so the sign agrees with that in the unscaled parameter coordinates. Thus the determinant of the local change of variables has positive sign.

\noindent\textit{Step 3. Disjointness of the neighborhoods.}\quad Choose the closed center neighborhoods so that the compact images consisting of the unordered centers of each type and the corresponding mean-zero scalar field are disjoint for distinct selected shadow solutions.

At a zero, the centers of each type and the scalar field determine the shadow solution. The choice follows by continuity of the scalar branches and finiteness. Use one labeling within each neighborhood.

Suppose two representations by $\Psi$ give the same function along a sequence
\[
 \mu\to0.
\]
Write $p$ and $\widetilde p$ for the two parameter choices, and use
tildes for quantities associated with $\widetilde p$. The corrections
are $o(1)$ in $H^1$. Since the trial fields have mean zero,
\[
 \|B_p-B_{\widetilde p}\|_{L^2(M)}
 \le C\|\nabla(B_p-B_{\widetilde p})\|_{L^2(M)}=o(1).
\]
Pass to a subsequence along which the centers and shapes converge,
and let $\Sigma$ be the union of the limiting blow-up sets.
Off $\Sigma$ the cut-off corrections tend locally to zero by
\eqref{hyb:matching}. For the limiting exterior fields $u^0,\widetilde u^0$ we obtain
\[
 B_p\to u^0,\qquad B_{\widetilde p}\to\widetilde u^0
 \quad\text{in }L^2_{\rm loc}(M\setminus\Sigma).
\]
Thus $u^0=\widetilde u^0$ off $\Sigma$. The Green logarithms are
locally integrable, so the equality holds in distributions on $M$.
Writing $p_a$ for the limiting centers, we have
\[
 -\Delta_g K^{-1}u^0
 =4\pi\sum_a\binom{s_a}{t_a}(\delta_{p_a}-1)
  +\binom{0}{\gamma(He^w/I-1)}.
\]
The scalar density has no atoms. Comparing the Dirac masses identifies
the unordered centers of each type. Subtracting their Green potentials
and using the mean-zero normalization gives
\[
 w=\widetilde w.
\]
For distinct shadow solutions this common limit would belong to two
disjoint compact images. Thus their neighborhoods are disjoint for small $|\mu|$.

\noindent\textit{Step 4. Uniqueness of the parameters.}\quad
On the scale interval of Step~1,
\[
 |\widetilde T-T|\le T_+-T_-=s_*\log(\tau_+/\tau_-).
\]
The coefficients in \eqref{hyb:matching} are bounded above and away
from zero. Hence, uniformly as $\mu\to0^+$,
\[
 C^{-1}\le\frac{\widetilde R_a}{R_a}\le C,\qquad
 C^{-1}\le\frac{\widetilde r_a}{r_a}\le C\quad(t_a>0).
\]
At the prescribed reference point, write $s_0$ for the outer exponent
and $R_0$ for the outer width of $B_p$. The inner correction cancels
in the first component of $K^{-1}B_p$, whose radial term is
\[
 -\log(R_0^{2s_0}+|z|^{2s_0}).
\]
Take a subsequence with $\widetilde R_0/R_0\to\lambda>0$ and set
$z=R_0y$. On fixed rescaled disks the cut-off is one, the smooth
exterior gradients are $O(R_0)$, and the first-order corrections are
$o(1)$ in Dirichlet norm. It follows that
\[
 \nabla_y\log(1+|y|^{2s_0})
 =\nabla_y\log(\lambda^{2s_0}+|y|^{2s_0}).
\]
At $|y|=1$ the radial derivatives give $s_0=2s_0/(1+\lambda^{2s_0})$,
so every limiting ratio is one. By Step~1 of the proof of
Lemma~\ref{hyb:trial-energy}, $d_0$ depends only on the centers and the
scalar branch. Hence
\[
 \frac{\widetilde d_0}{d_0}\to1,\qquad
 \widetilde T-T
 =s_0\log\frac{\widetilde d_0}{d_0}
  -2s_0\log\frac{\widetilde R_0}{R_0}=o(1).
\]
At a free ordinary point, $(K^{-1}B_p)_1$ has logarithmic part
$-\log(R^2+|z|^2)$. At a free mixed point the two logarithmic parts
are $-\log F$ and $-\log Q_0$ by \eqref{hyb:gram}. In both cases
the smooth remainders have bounded gradients after the tails are restored.
On the mixed inner scale,
\[
 z=ry,\qquad r^2=\frac{\beta R^4(1+|B|^2)}4.
\]
Since $B$ is bounded and $r/R=O(R)$, \eqref{hyb:gram} gives
\[
 \frac{Q_0(ry)}{4R^4r^2}
 =|y|^2+\frac{1+|B+(r/R)^2y^2|^2}{1+|B|^2}
 \longrightarrow 1+|y|^2
 \quad\text{in }C^1_{\rm loc}(\mathbb C).
\]
The corresponding statement on the outer scale is
\[
 F(Ry)/R^4\to 1+|y^2-B|^2
\]
in $C^1_{\mathrm{loc}}$. These limits hold along a subsequence in the compact shape set. Rescale by the first ordinary or inner width.

If the relative displacement is unbounded, the second rescaled gradient vanishes on compact sets. Indeed, outside its width the logarithmic gradient is $O(1/|z|)$, and the rescaled smooth error tends to zero.

This would contradict the nonzero gradient of the limiting ordinary profile, since the Dirichlet norm of the difference tends to zero.

For bounded relative displacement and a bounded positive width ratio, equality of the two limiting ordinary gradients fixes the displacement by their unique zero.

At this zero the second derivative is a negative scalar multiple of the identity. The absolute value of the scalar determines the dilation. Thus the free center displacements are $o(R)$ at ordinary points and $o(r)$ at mixed points.

On the mixed outer scale the center displacement is consequently $o(R)$, and the outer widths have ratio tending to one by the already identified $T$ and the outer matching coefficient, which depends only on the centers of each type and the scalar branch.

Equality of the outer inverse-Cartan gradients makes $1+|y^2-B|^2$ and the second polynomial proportional. Their quartic coefficients are both one, and their quadratic terms are $-2\operatorname{Re}(y^2\overline B)$, so the limiting shapes are equal.

Thus the scaled parameters differ by $o(1)$. Use the same reference widths along the segment joining them. The uniform derivative bounds for the sources and the approximate solutions give
\[
 d_p[L_p(u-B_p)]=-J_0+o(1)
\]
uniformly on the segment in the scaled coordinates. Put
$p_t=(1-t)p+t\widetilde p$. Since both parameter choices represent $u$,
\[
 0=L_{\widetilde p}(u-B_{\widetilde p})-L_p(u-B_p)
 =\left(\int_0^1
   \left.d_q[L_q(u-B_q)]\right|_{q=p_t}\,dt\right)(\widetilde p-p).
\]
The averaged matrix is invertible by \eqref{hyb:chart}, so
$\widetilde p=p$ and $\widetilde v=v$. At each representation, the same
scaled estimates and the bound on $\eta_p+U_pv$ give
\[
 \det d_p[L_p(u-B_p)]\ne0.
\]
The implicit function theorem gives $p=p(u)$ locally, and injectivity
makes the inverse single-valued. Its fiber component is
\[
 v=U_p^{-1}(u-B_p-\eta_p).
\]
Write $\Omega_\mu$ for the open parameter domain of Step~1 and
$r(T)=\tfrac12 e^{-\delta_0T}T^{C_0}$ for the fiber radius. The local
inverse is continuous, so the image
\[
 \mathcal U_{Z,\mu}
 =\{\Psi(p,v):p=(z,B,T)\in\Omega_\mu,\ v\in E_{p_*},\
                \|v\|_{\Hh}<r(T)\}
\]
is open. For fixed $\mu>0$ it is bounded by continuity on the compact
base $\overline\Omega_\mu$ and boundedness of $r$. If
$u_j=\Psi(p_j,v_j)\to u$ in $\Hh$, pass to a subsequence with
$p_j\to p\in\overline\Omega_\mu$. Then
\[
 v_j=U_{p_j}^{-1}(u_j-B_{p_j}-\eta_{p_j})
 \longrightarrow U_p^{-1}(u-B_p-\eta_p)=v
 \quad\text{in }E_{p_*}.
\]
Thus a boundary point has $p\in\partial\Omega_\mu$ or $\|v\|_{\Hh}=r(T)$.
At a zero the correction lies in the uniqueness ball of Step~2.
That step and the boundary signs in Step~1 give
\[
 \mathcal F_\mu(\Psi(p,v))=0
 \ \Longrightarrow\ (p,v)\in\Omega_\mu\times\{0\},
 \qquad 0\notin\mathcal F_\mu(\partial\mathcal U_{Z,\mu}).
\]
A transposition of regular centers of the same type has sign $(-1)^{2\cdot2}=1$.

\smallskip\noindent\textit{Step 5. The local degree.}\quad
We compute the degree in the coordinates of Step~2. At $u=B_p+\eta_p$,
\[
 Df_i(u)[\phi]
 =f_i(u)\left(\phi_i-\frac1{\rho_i}
                 \int_M f_i(u)\phi_i\,dV_g\right).
\]
For fixed $p$, the exponential estimate gives $f_i(u)\in L^{p_0}$
for some $p_0>1$. By H\"older's inequality and the compact Sobolev
embeddings, for $1<q<p_0$ we have
\[
 Df_i(u):\Hh\longrightarrow L^q(M)\quad\text{compact},
 \qquad \|G*h\|_{H^1(M)}\le C_q\|h\|_{L^q(M)}.
\]
Thus $\widehat A_p:=U_p^{-1}A_pU_p$ is an invertible compact perturbation
of the identity on $E_{p_*}$. Contracting $p$ to $p_*$ within the base gives
\[
 \ker\widehat A_{(1-t)p+tp_*}=\{0\}\quad(0\le t\le1),
 \qquad \widehat A_{p_*}=A_{p_*}.
\]
For fixed $\mu>0$, excise to the ball $B_a\subset E_{p_*}$ of radius
$0<a<\min_{[T_-,T_+]}r(T)$ centered at zero. By
Proposition~\ref{hyb:range}\textup{(ii)}, the self-adjoint operator
$A_{p_*}$ has index $\nu$ given by \eqref{hyb:index}. Hence
\[
 \deg_{\mathrm{LS}}(A_{p_*},B_a,0)=(-1)^\nu,\qquad
 d_{Z,\mu}
 =\deg_{\mathrm B}\bigl((C,Q,S),\Omega_\mu,0\bigr)(-1)^\nu
 =(-1)^{\nu+1}\deg_{\mathrm B}F_0,
\]
by the product formula and Lemma~\ref{hyb:homotopy}.
The angular terms in \eqref{hyb:index} are even, and the relative
dilation block is positive by \eqref{hyb:positive-schur}. Consequently,
\[
 n-1+(V+|Q|)=U+|P|-1+2(V+|Q|),\qquad
 \nu\equiv m_S+U+|P|-1\pmod 2.
\]
Step~2 of the proof of Proposition~\ref{fnd:shadow-degree} gives
\[
 d_S(Z)=(-1)^{m_S}\deg_{\mathrm B}F_0.
\]
On the positive side the product is therefore $(-1)^{U+|P|}d_S(Z)$. The negative-side domain is empty, giving \eqref{hyb:jump}. Finally Theorem~\ref{hyb:capture} and \eqref{hyb:mass-law} place every unbounded sequence with this limiting shadow solution in the isolating neighborhood on the positive side.

Its limiting centers and shape parameters lie in the interior, and its correction belongs to the projected family. The passage from these neighborhoods to the total jump in degree, including the contributions from regular blow-up points alone, uses Proposition~\ref{prop:global-localization}.
\end{proof}
\begin{remark}
The outer radii at the singular and regular points may decay at different rates.
\begin{itemize}
\item[\textup{(i)}]
The local identities give the shape equation at each radius. The energy expansion gives the common scale equation to order $\mathcal L$.
\item[\textup{(ii)}]
Uniform $C^0$ convergence and the boundary signs suffice to compute the degree, even when the shape equations have several zeros.
\end{itemize}
\end{remark}
\section{Consequences and examples}
\label{sec:algebra}
We first prove the polynomial identity in \eqref{eq:degree-candidate} and compute the coefficients for integer source strengths, including their eventual value on a torus.

Examples with noninteger source strengths illustrate the range of continuation to negative source strengths. We then compute the degree in the first subcritical strips and derive the local estimates for negative source strengths.
\begin{corollary}
\label{cor:coefficient}
Suppose that all source strengths are nonnegative integers. Let
\[
 \prod_{p\in S}Q_{\alpha_{p,1},\alpha_{p,2}}
       =\sum_{r,s\geq0}c_{r,s}x^ry^s,\qquad
 b_n=(-1)^n\binom{\chi-1}{n}\quad(n\geq0),\qquad b_n=0\quad(n<0).
\]
Then
\begin{equation}\label{eq:finite-convolution}
 d_{k,l}
 =\sum_{\substack{0\leq r\leq k\\0\leq s\leq l}}
 c_{r,s}\sum_{j=0}^{\min(k-r,l-s)}
       (-1)^j\binom{\chi}{j}b_{k-r-j}b_{l-s-j}.
\end{equation}
The binomial coefficients have their generalized meaning for negative upper entries.
\end{corollary}
\begin{proof}
Expand \((1-x)^{\chi-1}\), \((1-y)^{\chi-1}\), and \((1-xy)^\chi\) in \eqref{eq:degree-candidate} and compare the coefficients in \eqref{eq:analytic-identification}. All sums for fixed \(k,l\) are finite.
\end{proof}
\begin{corollary}
\label{cor:positive-genus}
Suppose all source strengths are nonnegative integers.
\begin{itemize}
\item[\textup{(i)}]
If $\chi\leq0$, then $d_{k,l}>0$ for every $k,l\geq0$. Hence \eqref{eq:toda} has a solution for every positive noncritical parameter pair.
\item[\textup{(ii)}]
If $\chi=0$, then
\[
 d_{k,l}=\prod_{p\in S}\frac{A_pB_p(A_p+B_p)}2
\]
whenever $k,l$ are at least the respective degrees in $x,y$ of $\prod_p Q_{\alpha_{p,1},\alpha_{p,2}}$.
\end{itemize}
\end{corollary}
\begin{proof}
\smallskip\noindent\textit{(i)}\quad For \(\chi\leq0\), the first two binomial factors in the second expression of \eqref{eq:degree-candidate} have positive coefficients. The remaining factor and all the \(Q\)'s have nonnegative coefficients and constant term one.

Their constant terms already give a strictly positive contribution to every coefficient. The nonzero-degree property and Corollary~\ref{thm:integer-degree} give the existence assertion.

\smallskip\noindent\textit{(ii)}\quad If \(\chi=0\), multiplication by \((1-x)^{-1}(1-y)^{-1}\) takes rectangular partial sums of the finite polynomial. Past its support the sum is its value at \((1,1)\). Apply Proposition~\ref{prop:quotient}.
\end{proof}
\subsection{The regular problem and the first two strips}
With \(S=\varnothing\), \eqref{eq:finite-convolution} reduces to
\[
 d_{k,l}
 =\sum_{j=0}^{\min(k,l)}(-1)^j\binom{\chi}{j}b_{k-j}b_{l-j}.
\]
For a torus this is one in every parameter rectangle. On the sphere the series is the polynomial
\[
 (1-x)(1-y)(1-xy)^2.
\]
For general integer singularities, the degrees in the first strip are given by
\[
 \sum_{l\geq0}d_{0,l}y^l=(1-y)^{\chi-1}
       \prod_{p\in S_2}(1+y+\cdots+y^{\alpha_{p,2}}),
\]
which agrees with the scalar degree normalization \(2\rho_2\). For the next strip, put \(T_a(y)=1+\cdots+y^a\) and \(P_2(y)=\prod_{p\in S_2}T_{\alpha_{p,2}}(y)\). Then
\begin{align}
 \sum_{l\geq0}d_{1,l}y^l
  =(1-y)^{\chi-1}\bigg\{&
     [1-(\chi-\ell)(1+y)]P_2(y)\notag\\
 &-\sum_{p\in S\setminus S_1}T_{\alpha_{p,2}+1}(y)
       \prod_{q\in S_2\setminus\{p\}}T_{\alpha_{q,2}}(y)\bigg\}.
       \label{eq:first-strip}
\end{align}
To see this, write
\[
 F(x,y)=(1-y)\{1-(1+y)x+yx^2\},\qquad
 P_{A,B}(x,y)=1-y^B+c_{A,B}(y)x+O(x^2),
\]
where
\[
 c_{A,B}(y)=
 \begin{cases}
 -1+y^{B+1},&A=1,\\
 0,&A>1.
 \end{cases}
\]
Substitution in \eqref{eq:degree-candidate} and comparison of the terms linear in \(x\) give \eqref{eq:first-strip}.
\subsection{An example with integer source strengths}
For a single source with strengths \((\alpha_{p,1},\alpha_{p,2})=(1,1)\),
\[
 Q_{1,1}=1+x+y+2xy+x^2y+xy^2+x^2y^2,
 \qquad d_{1,1}=\chi^2-5\chi+5.
\]
Thus the degree in \((4\pi,8\pi)^2\) is respectively \(-1,5,19\) for the sphere, torus, and a genus-two surface. For the torus the first few degrees are
\[
\begin{array}{c|rrrr}
\toprule
 k\backslash l&0&1&2&3\\
\midrule
0&1&2&2&2\\
1&2&5&6&6\\
2&2&6&8&8\\
3&2&6&8&8\\
\bottomrule
\end{array}
\]
The eventual value is \(8\). Adding a marked point with source strengths \((0,0)\) does not change the degree, since
\[
 Q_{0,0}=1,
\]
or equivalently \(P_{1,1}=F\). Interchanging the two components interchanges \(x,y\) and the corresponding source strengths, in agreement with the symmetry of \(K\).
\subsection{A torus with noninteger source strengths}
The polynomial quotient in Proposition~\ref{prop:quotient} and the positivity conclusion of Corollary~\ref{cor:positive-genus} concern integer source strengths. For noninteger source strengths, the degree can be negative or zero even on a torus. Consider one prescribed point on a torus, with
\[
 A=\sqrt2,\qquad B=\sqrt3,
 \qquad (\alpha_1,\alpha_2)=(\sqrt2-1,\sqrt3-1).
\]
Both source strengths are positive, and $1,\sqrt2,\sqrt3$ are linearly independent over $\mathbb Q$. Write $d(r,s)$ for the degree at $(\rho_1,\rho_2)=(4\pi r,4\pi s)$. Since $\chi=0$ and $\ell=1$, the series whose rectangular partial sums give this degree is
\[
 J_{A,B}(x,y)=\frac{P_{A,B}(x,y)}{F(x,y)}.
\]
For real $t,u$, define
\[
 C(t,u)=\sum_{k\geq0}\lceil t-k\rceil_+\lceil u-k\rceil_+,
 \qquad a_+=\max\{a,0\}.
\]
Indeed, expanding $F^{-1}=\sum_{i,j,k\geq0}x^{i+k}y^{j+k}$ shows that $C(t,u)$ is the sum of its coefficients at exponents strictly below $(t,u)$. In particular $C(t,u)=0$ if $t\leq0$ or $u\leq0$. Theorem~\ref{thm:main} and \eqref{eq:P} give
\begin{align*}
 d(r,s)={}&C(r,s)-C(r-A,s)-C(r,s-B)\\
 &+C(r-A,s-A-B)+C(r-A-B,s-B)\\
 &-C(r-A-B,s-A-B).
\end{align*}
Both parameter pairs considered below are noncritical. In units of $4\pi$, the sharp exclusion sets consist of the positive elements of
\[
 \{n+e:n\in\mathbb Z_{\geq0},\ e\in\{0,A,A+B\}\},\qquad
 \{n+e:n\in\mathbb Z_{\geq0},\ e\in\{0,B,A+B\}\}.
\]
The numbers $17/10$ and $29/10$ are nonintegral rationals, whereas $A$ and $B$ are irrational; also $A+B>3$. Thus neither parameter pair meets these sets. For $r=s=29/10$, the last three terms vanish and direct summation gives
\[
 d(29/10,29/10)=(9+4+1)-(6+2)-(6+2)=-2.
\]
For $r=17/10$ and $s=29/10$, the same three terms vanish and
\[
 d(17/10,29/10)=(6+2)-3-(4+1)=0.
\]
The first value implies existence of a solution. The second gives no conclusion on existence. In particular, the strict positivity in Corollary~\ref{cor:positive-genus} does not extend to arbitrary positive source strengths.
\subsection{Further results for signed source strengths}
\label{sec:all-signed-strips}
The scalar degree formula gives a direct extension to arbitrary source strengths greater than $-1$ when one parameter is sufficiently small. The other parameter need not be subcritical.
\begin{theorem}[Degree in the first subcritical strips]
\label{thm:all-signed-first-strip}
Assume that $\alpha_{p,i}>-1$ for every $p\in S$ and $i=1,2$, and put
\[
 a_1=\min\bigl(\{1\}\cup\{A_p:p\in S\}\bigr),\qquad
 a_2=\min\bigl(\{1\}\cup\{B_p:p\in S\}\bigr).
\]
If
\[
 0<\rho_1<4\pi a_1,\qquad \rho_2>0,
 \qquad \rho_2\notin\mathcal C_2(\alpha),
\]
then the solution set of \eqref{eq:toda} is compact and
\begin{equation}\label{eq:all-signed-first-strip}
 d_\alpha(\rho_1,\rho_2)
 =\mathscr D_\alpha\left(\frac{\rho_1}{4\pi},
                         \frac{\rho_2}{4\pi}\right).
\end{equation}
The same conclusion holds for
\[
 0<\rho_2<4\pi a_2,\qquad \rho_1>0,\qquad
 \rho_1\notin\mathcal C_1(\alpha).
\]
No rational independence assumption \eqref{fnd:independence} is required.
\end{theorem}
\begin{proof}
\smallskip\noindent\textit{Step 1. Compactness along the homotopy.} We consider the first set of hypotheses and prove a uniform bound for all zeros along the homotopy
\[
                   (t\rho_1,\rho_2),\qquad 0\leq t\leq1.
\]
In \eqref{eq:toda} the parameters are positive. The estimate here includes the endpoint $t=0$. We first consider a sequence of solutions with $0<t_k\leq1$ and $t_k\to t_\infty\in[0,1]$. Set
\[
 I_i^k=\int_Mh_i e^{u_i^k}\,dV_g,\qquad
 c_1^k=\log\frac{t_k\rho_1}{I_1^k},\qquad
 c_2^k=\log\frac{\rho_2}{I_2^k},\qquad
 f_i^k=e^{c_i^k}h_i e^{u_i^k}.
\]
Since the regular parts have zero mean, Jensen's inequality gives
\[
 I_i^k\geq\exp\left(\int_M\log h_i\,dV_g\right)>0.
\]
In particular, $c_i^k\leq C$; if $t_\infty=0$, then $c_1^k\to-\infty$. Passing to a subsequence, we obtain
\[
 f_i^k\,dV_g\rightharpoonup\nu_i,\qquad
 c_i^k\longrightarrow c_i\in[-\infty,C].
\]
The total masses are
\[
                  \nu_1(M)=t_\infty\rho_1,
                  \qquad \nu_2(M)=\rho_2.
\]
The local estimates in the proof of Proposition~\ref{fnd:signed-compactness} apply to this sequence. They also apply at $t_\infty=0$, since they require only the upper bounds for $c_i^k$, the bounded total masses, and an $L^p$ bound, with $p>1$, for the fixed weights. The Green representation
\[
                 u_i^k=\sum_jK_{ij}G*f_j^k
\]
and the lower bound for $G$ give $u_i^k\leq2G*(\mathbf1_D f_i^k)+C$ on a smaller disk. Choose fixed $p,q>1$ with $h_i\in L^{pq}$, and set $q'=q/(q-1)$. If the mass
\[
 m_i^k=\int_D f_i^k
\]
is below a fixed threshold strictly less than $2\pi/(pq')$, Jensen's inequality for the probability $f_i^k\mathbf1_D/m_i^k$ bounds the integral of $e^{2pq'G*(\mathbf1_D f_i^k)}$ on the smaller disk by a constant times $\sup_y\int_{D'}|z-y|^{-pq'm_i^k/\pi}\,dV_g(z)$. The case of zero local mass is immediate.

Together with the weight and the upper bound for $c_i^k$, H\"older's inequality then gives a local $L^p$ bound for $f_i^k$. Let $\mathcal B$ be the finite set of points where a limiting local mass reaches this threshold.

At every other point both local densities are bounded in $L^p$ and the near logarithmic kernels give bounds for both signs of $u_i^k$.

The separated Green kernels contribute bounded terms. Elliptic compactness gives local convergence away from the concentration points and the marked points.

For global $L^1$ convergence, truncate the Green kernel at distance $\varepsilon$ from its pole. The integral of the omitted logarithmic part in the observation variable is $o(1)$ as
\[
 \varepsilon\downarrow0,
\]
uniformly in the pole. The truncated kernels pass to the weak measure limit, giving
\[
 u_i^k\longrightarrow u_i=\sum_jK_{ij}G*\nu_j
 \quad\text{in }L^1(M).
\]
Thus a finite cover gives a global $C^0$ bound when
\[
 \mathcal B=\varnothing.
\]
At a point of $\mathcal B$, the reduction to a local Toda system also remains valid. In a conformal coordinate write $h_i=|z|^{2\alpha_i}b_i(z)$ and choose smooth functions $\xi_i^k$ with zero boundary values such that
\[
 \Delta\xi_i^k=
       \bigl(K(t_k\rho_1,\rho_2)^{\mathsf T}\bigr)_i e^\varphi.
\]
Then
\[
 v_i^k=u_i^k+c_i^k+2\alpha_i\log|z|-\xi_i^k
\]
satisfies the local Toda equation with coefficients
\[
 \widehat h_i^k=e^\varphi b_i e^{\xi_i^k},
\]
uniformly positive and bounded in $C^3$. Its component masses and boundary oscillations are uniformly bounded. The constants $c_i^k$ do not affect the latter bounds.

A positive local mass in this disk implies that at least one of the regular parts $v_i^k-2\alpha_i\log|z|$ has unbounded supremum; otherwise both densities on the disk have a common $L^p$ bound with $p>1$. The other points of the disk have uniform local bounds. The fields
\[
 U_i^k=v_i^k+\log\widehat h_i^k
\]
have bounded smooth forcing and are locally regular after subtracting the fixed source logarithms. Lemmas~\ref{fnd:local-unitary} and \ref{fnd:coarse-local} apply: the proofs allow a normalizing constant tending to $-\infty$, since they use mass bounds and boundary oscillations.

At an unmarked point both source strengths are zero. In particular, every nonzero local mass pair, written in units of $4\pi$ as $(x,y)$, satisfies
\begin{equation}\label{eq:first-strip-pohozaev}
                    x^2-xy+y^2=Ax+By,
\end{equation}
where $(A,B)=(A_p,B_p)$ at $p\in S$, and $(A,B)=(1,1)$ elsewhere. The local masses belong to the lists defining $\mathcal C_i(\alpha)$.

\smallskip\noindent\textit{Step 2. Exclusion of concentration.} We claim that
\[
 c_2> -\infty.
\]
If $c_2=-\infty$, the local bounds outside $\mathcal B$ imply that $\nu_2$ has no residual measure. All its mass is therefore concentrated at finitely many points. The local classification then gives
\[
 \rho_2=\nu_2(M)\in\mathcal C_2(\alpha),
\]
contrary to the hypothesis. Thus $\nu_2$ has the positive residual density $e^{c_2}h_2e^{u_2}$ away from $\mathcal B\cup S$. Its integrability at a point with local masses $(x,y)$ implies
\begin{equation}\label{eq:first-strip-diffuse-inequality}
                              2y-x\leq B.
\end{equation}
Indeed, the point masses contribute $-2(2y-x)\log r$ to the circular mean of $u_2$. For the remaining measures the logarithmic kernel on the disk has circular mean $\log\max\{r,|z|\}$ in the integration variable $z$. The measures assign zero mass to the center.

For $r<d$, divide by $|\log r|$, bound the part from $|z|<d$ by a constant times its mass, and use $O_d(1)$ on the complement. Sending first $r$ and then $d$ to zero gives $o(|\log r|)$ for each remainder. If
\[
 2y-x>B=1+\alpha_2,
\]
the circular mean of the logarithm of the second density is $(-2-\epsilon)\log r+o(|\log r|)$ for some $\epsilon>0$. Jensen's inequality on circles would bound the circular mean of the density below by $cr^{-2-\epsilon/2}$, which is not integrable.

On the other hand, nonnegativity of the limiting measures gives
\[
 0\leq x\leq\frac{t_\infty\rho_1}{4\pi}
       <a_1\leq A.
\]
Equations \eqref{eq:first-strip-pohozaev} and \eqref{eq:first-strip-diffuse-inequality} are incompatible with a nonzero pair $(x,y)$. In fact,
\[
 \begin{split}
 0&=x(x-A)+y(y-x-B),\\
 y-x-B&\leq-\frac{x+B}{2}<0.
 \end{split}
\]
The first summand is negative when $x>0$, and the second is negative when $y>0$. Therefore every local mass is zero.

The local estimates already established now give a uniform $C^0$ bound for $u^k$. This argument also covers $t_\infty=0$: in that case $\nu_1(M)=0$ and hence $x=0$ at every point.

\smallskip\noindent\textit{Step 3. Reduction to the scalar degree.} At $t=0$, the two equations and the mean-zero normalization give $u_1=-u_2/2$, while $u_2$ satisfies
\begin{equation}\label{eq:first-strip-scalar-endpoint}
 -\Delta_g u_2=2\rho_2
 \left(\frac{h_2e^{u_2}}{\int_Mh_2e^{u_2}\,dV_g}-1\right),
 \qquad \int_Mu_2\,dV_g=0.
\end{equation}
The scalar critical parameters, expressed in terms of $\rho_2$, form the set
\[
 (0,\infty)\cap4\pi
 \left\{n+\sum_{p\in E}B_p:
       n\in\mathbb Z_{\geq0},\ E\subset S\right\},
\]
which is contained in $\mathcal C_2(\alpha)$. Lemma~\ref{fnd:scalar-parameter} therefore bounds all the endpoint solutions, with shifted parameters $B_p=1+\alpha_{p,2}$ and parameter $\rho_2$. Choose $p>1$ such that $h_i\in L^p(M)$. On a bounded subset of
\[
 (C^0(M)\cap\{\int_Mu=0\})^2,
\]
the normalized densities depend continuously in $L^p$ on $(t,u)$, including at $t=0$, and have uniformly positive denominators. Since the Green operator maps $L^p$ into $W^{2,p}$, compactly embedded in $C^0$, the homotopy is continuous and compact on bounded sets.

The uniform solution bound just proved makes it admissible on one sufficiently large ball. Consequently its degree is the degree at $t=0$.

For the endpoint map, use the bounded linear change of variables $(u_1,u_2)\mapsto(a,b)=(u_1+u_2/2,u_2)$. Conjugation gives
\[
 (a,b)\longmapsto
 \left(a,\ b-2\rho_2G*
       \frac{h_2e^b}{\int_Mh_2e^b\,dV_g}\right).
\]
The reduction property and invariance under linear conjugation identify its degree with the scalar degree in \eqref{eq:first-strip-scalar-endpoint}. By Proposition~\ref{fnd:scalar-degree-internal}, that degree is the sum of the coefficients with exponent less than $\rho_2/(4\pi)$ in
\[
             (1-y)^{\chi-\ell}
                       \prod_{p\in S}(1-y^{B_p}).
\]
Finally, every positive $x$-exponent in $\mathscr J_\alpha(x,y)$ is at least $a_1$. Thus $\rho_1/(4\pi)<a_1$ selects only the terms with $x$-exponent zero, whose series is exactly
\[
 [x^0]\mathscr J_\alpha(x,y)
          =(1-y)^{\chi-\ell}\prod_{p\in S}(1-y^{B_p}).
\]
This proves \eqref{eq:all-signed-first-strip}. Interchanging the two components proves the other assertion.
\end{proof}
\begin{corollary}[Continuation from the first subcritical strips]
\label{cor:signed-strip-continuation}
Fix the marked points and the positive smooth factors. Let $t\mapsto(\alpha(t),\rho(t))$ be a continuous path satisfying
\[
 \alpha_{p,i}(t)>-1,\qquad \rho_i(t)>0,\qquad
 \rho_i(t)\notin\mathcal C_i(\alpha(t))\quad(i=1,2,\ 0\le t\le1).
\]
If the initial point satisfies either set of hypotheses of Theorem~\ref{thm:all-signed-first-strip}, then
\[
 d_{\alpha(t)}(\rho_1(t),\rho_2(t))
 =\mathscr D_{\alpha(t)}\left(\frac{\rho_1(t)}{4\pi},
                             \frac{\rho_2(t)}{4\pi}\right)
 \qquad(0\leq t\leq1).
\]
\end{corollary}
\begin{proof}
Proposition~\ref{fnd:signed-compactness} makes the path an admissible compact homotopy. Formula~\eqref{eq:signed-finite-sum} shows that the formal cumulative sum is locally constant along it: each possible cut-off equality belongs to one of the excluded sets.

Both quantities are therefore constant. They agree at $t=0$ by Theorem~\ref{thm:all-signed-first-strip}.
\end{proof}
\begin{proposition}[Subcritical parameters]
\label{prop:signed-subcritical}
Fix arbitrary source strengths $\alpha_{p,i}>-1$ and set $a_i=\min(\{1\}\cup\{1+\alpha_{p,i}:p\in S\})$. If $0<\rho_i<4\pi a_i$ for $i=1,2$, then
\[
 d_\alpha(\rho_1,\rho_2)
 =\mathscr D_\alpha(\rho_1/(4\pi),\rho_2/(4\pi))=1.
\]
\end{proposition}
\begin{proof}
Consider the homotopy $t\rho$, $0\leq t\leq1$. At a blow-up point let $(x,y)$ be the local masses in units of $4\pi$, and write $(A,B)$ for its two shifted parameters; at a regular point take
\[
 A=B=1.
\]
The total-mass bounds imply $x<A$ and $y<B$. For every nonzero pair with $x,y\geq0$, these inequalities give
\[
 x^2-xy+y^2<Ax+By,
\]
contrary to the local Pohozaev identity. Thus the solutions are uniformly bounded when $t$ lies in a compact subinterval of $(0,1]$.

To treat $t\to0$, write $I_i=\int_Mh_ie^{u_i}\,dV_g$ and $f_i=t\rho_i h_ie^{u_i}/I_i$. Jensen's inequality gives
\[
 I_i\geq\exp(\int_M\log h_i\,dV_g)>0.
\]
For each fixed finite $q>1$, the Green representation and the logarithmic-kernel estimate in the proof of Proposition~\ref{fnd:signed-compactness} give $\int_Me^{qu_i}\,dV_g\leq C_q$ for all sufficiently small $t$, since
\[
 \int_Mf_i\,dV_g=t\rho_i\to0.
\]
The fixed weights belong to $L^p(M)$ for some $p>1$. H\"older's inequality therefore gives $\|f_i\|_{L^{1+\eta}(M)}\leq Ct$ for some fixed $\eta>0$. The mean-zero normalization and elliptic estimates imply
\[
 \|u_1\|_{C^0(M)}+\|u_2\|_{C^0(M)}\leq Ct.
\]
Thus the homotopy is admissible on a common ball, including at $t=0$; its continuity and compactness follow from the same $L^p$ argument.

At $t=0$ the degree is that of the identity, namely one. Finally, every nonconstant monomial in \eqref{eq:real-jump-series} has first exponent at least $a_1$ or second exponent at least $a_2$. Only the constant term contributes to the rectangular sum.
\end{proof}
There are also supercritical examples covered by continuation. On a torus with all source strengths zero, $P_{1,1}=F$ and $\chi=0$, hence
\[
 \mathscr J_0=1.
\]
Corollary~\ref{cor:small-signed-orders} gives $d_\alpha(6\pi,10\pi)=1$ whenever all source strengths lie in $(-\varepsilon,0]$, for some $\varepsilon>0$. The same conclusion holds uniformly on each fixed compact subset of $((0,\infty)\setminus4\pi\mathbb N)^2$, with a possibly smaller $\varepsilon$. No explicit lower bound for $\varepsilon$ is asserted.
\begin{remark}[A component inaccessible by continuation from the first strips]
\label{rem:signed-chamber-not-reached}
Theorem~\ref{thm:all-signed-first-strip} and continuation from nonnegative source strengths do not cover every admissible configuration with signed source strengths. To see this, take one marked point and set
\[
 A=\frac{\sqrt2}{10},\qquad B=\frac{\sqrt3}{10},\qquad
 r=\frac{\rho_1}{4\pi}=\frac15,
 \qquad s=\frac{\rho_2}{4\pi}=\frac15.
\]
Here $1,A,B$ are linearly independent over $\mathbb Q$, and
\[
                    A<r<1,\qquad B<s<1,
                    \qquad r,s<A+B.
\]
\begin{itemize}
\item[\textup{(i)}]
With masses in units of $4\pi$, the quadratic in the Pohozaev identity has its minimum in $y$ at $y=(x+B)/2$. Setting the minimum equal to zero gives
\[
 3x^2-(4A+2B)x-B^2=0.
\]
The larger root occurs with $y>0$. Interchanging the coordinates gives the second maximum. Thus the two maxima are
\[
 \frac{2A+B+2\sqrt{A^2+AB+B^2}}3,\qquad
 \frac{A+2B+2\sqrt{A^2+AB+B^2}}3.
\]
\item[\textup{(ii)}]
Since $A,B<1/5$, the two coordinate maxima satisfy
\[
 \max x,\ \max y < (3+2\sqrt3)/15 < 1 .
\]
Thus the recursion admits no positive integer increment. The two reflections with zero increment preserve the orbit of the six initial pairs, so these pairs exhaust $\Sigma_p$ and
\[
 \frac{\mathcal C_1}{4\pi}
  =(0,\infty)\cap
       \{n,n+A,n+A+B:n\in\mathbb Z_{\geq0}\},
 \qquad
 \frac{\mathcal C_2}{4\pi}
  =(0,\infty)\cap
       \{n,n+B,n+A+B:n\in\mathbb Z_{\geq0}\}.
\]
The chosen parameters avoid these sets.
\item[\textup{(iii)}]
Along a continuous path with the same marked set and $\rho_i\notin\mathcal C_i$, the initial local masses and the masses at regular points exclude
\[
 r=A,\qquad s=B,\qquad r=1,\qquad s=1.
\]
The four differences therefore keep their signs. Such a path reaches neither $A,B\ge1$ nor either first subcritical strip. Here
\[
                         \mathscr D_\alpha(r,s)=1-1-1=-1.
\]
The direct calculation in Section~\ref{sec:wall-completion} proves that this is the degree of the Toda system.
\end{itemize}
\end{remark}
\subsection{Local estimates for negative source strengths}
\label{sec:negative-local-estimates}
For negative source strengths we need the following local estimates. They extend the model inverse to positive shifted parameters and identify the fractional and nonlocal terms that enter when a source strength is negative.

In Lemma~\ref{hyb:covariance}, each model probability with a fixed profile exponent $q>0$ has the required exponential $L^{p_0}$ bound. Indeed, the rescaling
\[
 m_{q,R}(Ry)=R^{-2}m_{q,1}(y)
\]
gives
\[
 \|m_{q,R}\|_{p_0}=C_{q,p_0}R^{2/p_0-2}.
\]
Near zero the density is integrable to the required power when $p_0(1-q)<1$ for $q<1$; at infinity the density decays like $|y|^{-2q-2}$. Fixed cut-offs and normalization of the cut-off masses preserve the bound. The matching law gives
\[
 \log(1/R)=O(T),
\]
so the bound is exponential in $T$. The lemma also requires an $L^1$ estimate for the difference between the normalized density and the model density by a positive power of the radius.
\begin{proposition}[Positive profile exponents in the model operator]
\label{signed:model-inverse}
Fix $A_a,B_a>0$ in Theorem~\ref{hyb:model-inverse}. Assume the following.
\begin{itemize}
\item[\textup{(i)}]
Each active scalar profile exponent at a prescribed blow-up point is noninteger. On the given compact family the linearized scalar shadow equation has a uniformly bounded inverse.
\item[\textup{(ii)}]
The normalized residual density is bounded in $L^{p_0}(M)$ for some $p_0>1$, uniformly on the family.
\item[\textup{(iii)}]
The centers are uniformly separated, the matching coefficients are bounded above and away from zero, and $\inf(\rho_2-N)>0$ for $\gamma=\rho_2-N$.
\end{itemize}
Then the theorem holds, including the inverse estimate, the index formula and the positivity of the small relative radial eigenvalues.
\end{proposition}
\begin{proof}
\smallskip\noindent\textit{Step 1. The estimates on each scale.}\quad We check the estimates in the proof of Theorem~\ref{hyb:model-inverse} that involve a singular profile exponent.

The local variance estimate and its tail bound use only $q>0$. Indeed, the radial probability measure has density comparable to $r^{2q-1}$ at zero and to $r^{-2q-1}$ at infinity in scaled coordinates.

Circular Poincar\'e and Cauchy--Schwarz for the circular mean therefore give the bound $CL^{-2q}(1+\log L)\|\nabla f\|_2^2$. Thus centered weighted compactness is valid for every profile exponent $q>0$.

The kernel functions of the model operator are bounded. A radial field, after subtraction of its limit at infinity, is bounded by $C\min\{1,(R/|x|)^{2q}\}$. Hence for every $1\leq k<\infty$ and every
\[
 0<\eta<\min\{2q,2/k\}
\]
its $L^k$ norm is $O(R^\eta)$. For the residual density $f_g$, H\"older's inequality gives
\[
 \int_M f_g|\Psi|^2
 \le \|f_g\|_{p_0}\|\Psi\|_{2p_0'}^2=o(1),
\]
with a positive power of $R$. The corresponding centered bilinear form has the same decay. This replaces boundedness or smoothness of the residual density in the estimates for these kernel functions. For the means used in the contrast terms, the Green kernels satisfy
\[
 \sup_z\|G(\cdot,z)\|_{2p_0'}<\infty
\]
by uniform integrability of the logarithmic singularity on a compact surface. Integrating in $z$ against any bubble probability preserves this bound by the integral triangle inequality. H\"older consequently gives a uniform second moment of its Green average against $f_g$.

The Green potential of $f_g$ itself is uniformly bounded by the $W^{2,p_0}$ estimate. Moreover, $H^1(M)\hookrightarrow L^{2p_0'}(M)$ is compact. Thus the term involving $f_g$ is uniformly compact in the limit scalar equation.

\smallskip\noindent\textit{Step 2. The interaction between the two scales.}\quad At a nested singular pair, $s=A+B>t=B>0$, and the matching equations give
\[
 \frac rR=C\exp\left(-\frac{A}{2Bs}T\right).
\]
The factor $C$ is bounded above and below. The outer radial field differs from its value at zero by $O(|y|^{2s})$ there; the inner radial field differs from its value at infinity by $O(|y|^{-2t})$ there.

Splitting the nested integrals between the two scales therefore gives a positive power of $r/R$, with a logarithm when two exponents coincide.

No estimate requires $s>1$ or $t>1$. The regular profiles and their translation and shape fields are unchanged. Consequently the error in the equation for the kernel functions, with possibly different positive constants, remains exponentially small in $T$.

\smallskip\noindent\textit{Step 3. The index.}\quad For a noninteger profile exponent $q>0$ the derivative with respect to dilation spans the finite-energy scalar kernel space.

The corresponding centered scalar form has negative index $2(\lceil q\rceil-1)$; in particular it is zero for $0<q<1$. The centered inverse and index proofs therefore apply.

To estimate the inner products of the mean Green potentials, we use the $L^{p_0}$ bound for each rescaled probability density. It bounds
\[
 \int(\log^-|y-z|)^2\,d\omega(y)
\]
uniformly in the translated pole $z$: use H\"older and the $L^{p_0'}$ norm of the squared logarithm on a unit disk. Its positive part is controlled by
\[
 C+\log(2+|y|)+\log(2+|z|),
\]
and the scaled probabilities have uniform finite second logarithmic moments. At a nested pair rescale by the outer radius. If
\[
 \delta=r/R\le1,
\]
these estimates bound the joint second moment of $\log|y-\delta z|$ uniformly, since
\[
 \log(2+|\delta z|)\le\log(2+|z|).
\]
Jensen also gives the same $L^2$ bound for the averaged kernels against either probability. Subtracting the logarithm of the outer radius therefore gives the $O(1)$ interaction remainder and the bounded centered variances needed for the resolvent.

For two factors at the same point use the radius of that point; the separated-disk interactions use bounded Green kernels. Finally, for differences of the local means the leading matrix of inner products is positive definite. Indeed, at each pair of nested scales,
\[
 \frac{c^2-cd}{s}+\frac{d^2}{t}
 =\frac{(c-d/2)^2}{s}
  +\left(\frac1t-\frac1{4s}\right)d^2,
 \qquad s>t>0.
\]
The normalized radial mean contrasts have the same triangular limiting matrix. The finite-rank inversion and the positive radial matrix after elimination in the original proof therefore give the stated conclusions.
\end{proof}
The $L^{p_0}$ hypothesis holds for the diffuse densities arising from fixed shifted parameters $A,B>0$ and a compact bounded scalar-shadow family. At a mixed singular point the density is
\[
 O(|x|^{2A-2}),
\]
whereas at a pure singular point it is $O(|x|^{2(A+B)-2})$. At a prescribed point with no bubble the exponent is
\[
 2\alpha_{p,2}>-2.
\]
Thus a common $p_0>1$ may be chosen for the finitely many prescribed points. This observation concerns the model operator. It supplies no estimate for the residual of an approximate nonlinear solution.
\begin{lemma}[A finite-part expansion for a negative singularity]
\label{signed:finite-part}
Let $0<q<1$, $\delta>0$, and $\phi\in C^2(\overline{D_\delta})$. Then
\begin{align*}
 \int_{D_\delta}m_{q,R}\phi
 ={}&4\pi q\phi(0)\\
 &+4q^2R^{2q}\left\{
   \int_{D_\delta}
    \frac{\phi(x)-\phi(0)-\nabla\phi(0)\cdot x}{|x|^{2q+2}}\,dA
       -\frac{\pi\phi(0)}{q\delta^{2q}}\right\}
   +o(R^{2q}).
\end{align*}
In particular, the coefficient of $R^{2q}$ is not determined by the derivatives of $\phi$ at the origin.
\end{lemma}
\begin{proof}
Radial integration gives
\[
 \int_{D_\delta}m_{q,R}
 =\frac{4\pi q\delta^{2q}}{R^{2q}+\delta^{2q}}
 =4\pi q-4\pi q\delta^{-2q}R^{2q}+o(R^{2q}),
 \qquad
 \int_{D_\delta}m_{q,R}\nabla\phi(0)\cdot x=0.
\]
Write $\phi_2(x)=\phi(x)-\phi(0)-\nabla\phi(0)\cdot x$. Since $|\phi_2(x)|\le C|x|^2$ and
\[
 R^{-2q}m_{q,R}(x)\le4q^2|x|^{-2q-2},
\]
the function on the right, multiplied by $|\phi_2|$, is integrable near zero precisely in the range used here. Dominated convergence therefore proves the expansion.

If $b$ is a nonzero nonnegative smooth function supported in an annulus contained in $D_\delta$, replacing $\phi$ by $\phi+\epsilon b$ preserves $\phi$ on a neighborhood of zero and changes the displayed coefficient by
\[
 4q^2\epsilon\int b|x|^{-2q-2}\,dA\ne0.
\]
\end{proof}
The same proof applies to a continuous $\phi$ for which
\[
 \phi(x)-\phi(0)-L(x)=O(|x|^\eta),
\]
where $L$ is linear and $\eta>2q$. In this case $L$ replaces $\nabla\phi(0)\cdot x$ in the finite-part integral. A fractional term of order at most $2q$ must instead be treated separately. For $\phi\in C_c^2(\mathbb R^2)$ one may equivalently write
\[
 \int_{\mathbb R^2}m_{q,R}\phi
 =4\pi q\phi(0)+4q^2R^{2q}\mathcal D_q[\phi]+o(R^{2q}),
\]
where, for any radial $\chi\in C_c^\infty(\mathbb R^2)$ equal to one near zero,
\[
 \mathcal D_q[\phi]
 =\int_{\mathbb R^2}
   \frac{\phi(x)-\phi(0)-\chi(x)\nabla\phi(0)\cdot x}
        {|x|^{2q+2}}\,dA.
\]
The integral is absolutely convergent and is independent of the radial cut-off. In particular, taking $\phi(x)=|x|^2$ on $D_\delta$ in the lemma gives
\[
 \int_{D_\delta}|x|^2m_{q,R}
 =\frac{4\pi q^2}{1-q}\delta^{2-2q}R^{2q}+o(R^{2q}).
\]
\begin{lemma}[The term from the residual density at a mixed singular point]
\label{signed:diffuse-cusp}
Let $0<A<1$ and suppose that, in a conformal coordinate centered at the origin, $w$ is continuous and satisfies
\[
 -\Delta w=2c|x|^{2A-2}(1+o(1))+O(1),\qquad c>0.
\]
Then its circular mean satisfies
\[
 \overline w(r)-w(0)
   =-\frac{c}{2A^2}r^{2A}+o(r^{2A}).
\]
Consequently, if $\psi_1=\psi_{\rm sm}-w/2$ with $\psi_{\rm sm}$ smooth, then
\[
 \overline{\psi_1}(r)-\psi_1(0)
   =\frac{c}{4A^2}r^{2A}+o(r^{2A}).
\]
\end{lemma}
\begin{proof}
Integrating the equation on $D_r$ gives
\[
 r\overline w'(r)
 =-\frac{c}{A}r^{2A}+o(r^{2A})+O(r^2).
\]
Integration from zero to $r$ proves the first assertion. The circular mean of a smooth function differs from its value at zero by $O(r^2)$, which is $o(r^{2A})$ because $A<1$.
\end{proof}
The fractional term in Lemma~\ref{signed:diffuse-cusp} must be retained in the radial calculation. Put $s=A+B$, $B>0$, and suppose $r^{2B}=bR^{2s}$, where $b=c/(4B^2)$. Direct integration on the plane gives
\begin{align*}
 \int_{\mathbb R^2}m_{s,R}Z_{s,R}|x|^{2A}
   &=-\frac{4\pi^2 A^2}{s\sin(\pi A/s)}R^{2A},\\
 \int_{\mathbb R^2}m_{s,R}Z_{s,R}(r/|x|)^{2B}
   &=\frac{4\pi^2B^2b}{s\sin(\pi A/s)}R^{2A}.
\end{align*}
Both integrals converge since $0<A<s$. To verify the identities, use $v=(|x|/R)^{2s}$ and the beta integral; the resulting integrals are respectively
\[
 \int_0^\infty\frac{v^{a}(1-v)}{(1+v)^3}\,dv
       =-\frac{\pi a^2}{\sin(\pi a)},\qquad
 \int_0^\infty\frac{v^{a-1}(1-v)}{(1+v)^3}\,dv
       =\frac{\pi(1-a)^2}{\sin(\pi a)},
 \qquad a=A/s.
\]
The two terms therefore satisfy
\[
 -\frac{c}{4A^2}\int m_{s,R}Z_{s,R}|x|^{2A}
 -\int m_{s,R}Z_{s,R}(r/|x|)^{2B}=0.
\]
This identity concerns these two terms only. It does not estimate the remaining terms in the reduced equation. In particular, a bound of order $R^{4A}$ for a quadratic correction is insufficient to obtain an $o(R^2)$ remainder when $A\leq1/2$.
\begin{lemma}[An exact radial family]
\label{signed:exact-radial}
Let $A,B>0$, $s=A+B$, and $c_0,c_1,c_2>0$. For $r=|x|>0$, set
\begin{align*}
 F(r)&=c_0+c_1r^{2A}+c_2r^{2s},\\
 Q(r)&=A^2c_0c_1r^{2A-2}+s^2c_0c_2r^{2s-2}
                   +B^2c_1c_2r^{2(A+s)-2},\\
 n_1(r)&=\frac{4Q(r)}{F(r)^2},\qquad
 n_2(r)=\frac{4A^2B^2s^2c_0c_1c_2r^{2(A+s-3)}F(r)}{Q(r)^2}.
\end{align*}
Then $u_i=\log n_i$ satisfy, in distributions on $\mathbb R^2$,
\[
 \Delta u_1+2e^{u_1}-e^{u_2}=4\pi(A-1)\delta_0,
 \qquad
 \Delta u_2+2e^{u_2}-e^{u_1}=4\pi(B-1)\delta_0.
\]
Their masses are
\[
 \int_{\mathbb R^2}n_1\,dA
 =\int_{\mathbb R^2}n_2\,dA=4\pi s.
\]
The equations include the Dirac masses at the origin.
\end{lemma}
\begin{proof}
\smallskip\noindent\textit{Step 1. The radial identities.}\quad The radial identity
\[
 \Delta\log\left(\sum_j d_jr^{2a_j}\right)
 =\frac{4\sum_{i<j}(a_i-a_j)^2d_id_j
                 r^{2(a_i+a_j)-2}}
        {(\sum_jd_jr^{2a_j})^2}
\]
gives $\Delta\log F=n_1$ for $r>0$. Write $Q=r^{2A-2}\widehat Q$, where
\[
 \widehat Q=A^2c_0c_1+s^2c_0c_2r^{2B}
                         +B^2c_1c_2r^{2s}.
\]
The same identity gives
\[
 \Delta\log Q=\Delta\log\widehat Q
 =\frac{4A^2B^2s^2c_0c_1c_2r^{2B-2}F}{\widehat Q^2}=n_2
 \quad(r>0).
\]
Taking the Laplacian of $\log n_1$ and $\log n_2$ proves the equations away from zero. The asymptotics at zero are
\[
 n_1(r)\sim\frac{4A^2c_1}{c_0}r^{2A-2},\qquad
 n_2(r)\sim\frac{4B^2s^2c_2}{A^2c_1}r^{2B-2}.
\]

\smallskip\noindent\textit{Step 2. The source terms and the masses.}\quad To check the distributional terms, the explicit power sums give $r(\log F)'=o(1)$ and
\[
 r(\log\widehat Q)'=o(1)
\]
as $r\downarrow0$. Their logarithms are bounded there and the radial integration by parts on punctured disks has zero limiting inner boundary contribution. Therefore on the whole plane
\[
 \Delta\log F=n_1,\qquad
 \Delta\log Q=n_2+4\pi(A-1)\delta_0 .
\]
In $\Delta\log n_2$ the explicit power $r^{2(A+s-3)}$ and $Q^{-2}$ contribute
\[
 4\pi(A+s-3)-8\pi(A-1)=4\pi(B-1)
\]
to the coefficient of $\delta_0$. This proves both claimed distributional equations. Both densities are integrable at zero by the asymptotics above. At infinity their highest powers give respectively $O(r^{-2B-2})$ and $O(r^{-2A-2})$. Finally, integration on annuli yields
\begin{align*}
 \int_{\mathbb R^2}n_1\,dA
   &=2\pi\left[ r(\log F)'\right]_{0}^{\infty}
     =2\pi(2s-0),\\
 \int_{\mathbb R^2}n_2\,dA
   &=2\pi\left[ r(\log Q)'\right]_{0}^{\infty}
     =2\pi\{2(A+s)-2-(2A-2)\}.
\end{align*}
\end{proof}
\begin{corollary}[A mixed radial degeneration]
\label{signed:exact-degeneration}
In Lemma~\ref{signed:exact-radial}, choose $c_0=R^{2s}$, $c_1=kR^{2s}$, and $c_2=1$, where $k>0$ is fixed. Put
\[
 \beta=\frac{B^2k}{s^2},\qquad
 \ell^{2B}=\frac{A^2k}{s^2}R^{2s},\qquad
 f_\beta(x)=\frac{4A^2\beta|x|^{2A-2}}
                          {(1+\beta|x|^{2A})^2}.
\]
As $R\to0$, $\ell/R\to0$, and
\[
 n_1\,dA\rightharpoonup4\pi s\delta_0,
 \qquad
 n_2\,dA\rightharpoonup4\pi B\delta_0+f_\beta\,dA.
\]
The rescaled densities satisfy
\[
 R^2n_1(Ry)\longrightarrow m_{s,1}(y),\qquad
 \ell^2n_2(\ell y)\longrightarrow m_{B,1}(y)
\]
locally uniformly on $\mathbb R^2\setminus\{0\}$.
\end{corollary}
\begin{proof}
The scale ratio is a positive constant times $R^{A/B}$. Write $\tau=R^{2s}$ and $Y=|y|$. Then
\[
 F(d)=\tau(1+kd^{2A})+d^{2s},\qquad
 \widehat Q(d)=s^2\tau(\ell^{2B}+d^{2B}+\beta d^{2s}).
\]
Consequently,
\begin{align*}
 R^2n_1(Ry)
 &=\frac{4s^2Y^{2A-2}\{(\ell/R)^{2B}+Y^{2B}+\beta R^{2A}Y^{2s}\}}
          {(1+kR^{2A}Y^{2A}+Y^{2s})^2},\\
 \ell^2n_2(\ell y)
 &=\frac{4B^2 (F(\ell Y)/\tau)Y^{2B-2}}
          {(1+Y^{2B}+\beta\ell^{2A}Y^{2s})^2}.
\end{align*}
Here $n_i$ are evaluated radially and $F(\ell Y)/\tau\to1$ locally uniformly because $\ell/R\to0$. These formulas yield the two rescaled limits on compact subsets away from zero. On such subsets in the original coordinate, $n_1\to0$ and $n_2\to f_\beta$ since
\[
 4B^2\ell^{2B}/\tau=4A^2\beta.
\]
At a fixed $d>0$ the derivatives $d\,\partial_d\log F,d\,\partial_d\log Q$ tend to $2s$ and $2s-2+2A\beta d^{2A}/(1+\beta d^{2A})$, respectively. Subtracting their respective limiting derivatives $0,2A-2$ at zero proves, by radial integration,
\[
 \lim_{R\to0}\int_{D_d}n_1\,dA=4\pi s,
 \qquad
 \lim_{R\to0}\int_{D_d}n_2\,dA
 =4\pi B+\frac{4\pi A\beta d^{2A}}{1+\beta d^{2A}}.
\]
The last summand is $\int_{D_d}f_\beta\,dA$. Since $\int f_\beta=4\pi A$ and both total masses equal $4\pi s$, these identities prove the weak convergence, with no loss of mass at infinity.
\end{proof}
For $A<1$ the local profile contains fractional powers of the radius. Its holomorphic representation permits matching to a prescribed scalar solution.
\begin{proposition}[Degree of a reduced equation]
\label{prop:signed-conditional-transfer}
Let $Z$ be a nondegenerate shadow solution for $(U,V,\mathcal P,\mathcal Q)$, and $m_S$ the Morse index of its scalar part. Assume the following for the finite-dimensional reduction near $Z$.
\begin{itemize}
\item[\textup{(i)}]
On the chosen closed product and near the graph of the correction, the $C^1$ parametrization is injective and preserves orientation. Along the graph the complementary linearization is invertible, with negative index
\[
 \nu\equiv m_S+U+|\mathcal P|-1\pmod2.
\]
\item[\textup{(ii)}]
The local degree of the center map is $\deg_{\mathrm B}F_0$ and the shape map has degree one, with
\[
 d_S(Z)=(-1)^{m_S}\deg_{\mathrm B}F_0.
\]
For all sufficiently small $\mu\Gamma>0$, a homotopy joins the reduced map to the product of the center, shape and scale maps. It has no boundary zero on the center and shape domains times
\[
 a|\mu|\le e^{-\lambda T}\le b|\mu|,
\]
where $0<a<\min(c,1/|\Gamma|)$ and $b>\max(C,1/|\Gamma|)$.
\item[\textup{(iii)}]
As $T$ increases towards concentration, the scale equation is
\[
 S(T,\mu)=-\mu+\Gamma e^{-\lambda T}
                   +o(|\mu|+e^{-\lambda T}),
 \qquad \lambda>0,\quad \Gamma\ne0.
\]
The error is uniform on the parameter neighborhoods whenever
\[
 0<a\le e^{-\lambda T}/|\mu|\le b<\infty\qquad(a,b\text{ fixed}).
\]
For $\mu\Gamma\le0$ the a priori estimates exclude sufficiently concentrated zeros. Every remaining zero satisfies
\[
 c|\mu|\le e^{-\lambda T}\le C|\mu|,\qquad
 0<c<C<\infty,
\]
with $c,C$ independent of $\mu$.
\item[\textup{(iv)}]
Every sequence blowing up at this critical value with limiting shadow solution $Z$ eventually lies in the neighborhoods of the reduction.
\end{itemize}
The local degree at $\mu>0$ minus that at $\mu<0$ is
\[
 \Delta d_Z=(-1)^{U+|\mathcal P|}d_S(Z),
\]
for either sign of $\Gamma$.
\end{proposition}
\begin{proof}
Write $\gamma=|\Gamma|$ and choose constants $0<a<\min\{c,1/\gamma\}$ and $b>\max\{C,1/\gamma\}$. If $\mu\Gamma>0$, take the scale interval whose endpoints satisfy
\begin{equation}\label{signed:conditional-scale-window}
 e^{-\lambda T_-}=b|\mu|,\qquad
 e^{-\lambda T_+}=a|\mu|.
\end{equation}
For small $|\mu|$, the scalar equation has sign $\sgn\Gamma$ at $T_-$ and sign $-\sgn\Gamma$ at $T_+$. Its one-dimensional degree is therefore $-\sgn\Gamma$.

This conclusion uses only the boundary signs and does not require uniqueness of a scale zero. At each fixed nonzero $\mu$ let $p$ denote the finite parameters in the closed product.

In the oriented coordinates of assumption \textup{(i)}, write the complementary equation at a fiber vector $v$ as $A_p v+E_p(v)$, where $A_p$ is the compressed linearization of the Toda gradient along the graph.

Its inverse is uniformly bounded on this compact product, and continuity of the derivative gives
\[
 \|E_p(v)\|=o(\|v\|)
\]
uniformly there. Choose a common positive fiber radius so that
\[
 \|A_p^{-1}E_p(v)\|\le\|v\|/2.
\]
The homotopy to $A_pv$ therefore has no complementary zero at
\[
 v\ne0.
\]
The finite-dimensional equations may be evaluated on the graph during this deformation because every zero of the full system then has $v=0$.

On the faces of the finite product, assumption \textup{(ii)} excludes those graph zeros.

The complementary degree is $(-1)^\nu$, since the linearization is self-adjoint in the Cartan inner product and is a compact perturbation of the identity. In these oriented product coordinates the degree is
\[
 d_{Z,\mu}=(-1)^\nu(-\sgn\Gamma)\deg_{\mathrm B}F_0
          =(-1)^{U+|\mathcal P|}\sgn\Gamma\,d_S(Z).
\]
On the opposite side the local contribution is zero by the stated exclusion. If $\Gamma>0$, the displayed contribution lies on the positive side.

If $\Gamma<0$, it lies on the negative side and is subtracted in the difference of the degrees on the two sides. Both cases give the asserted jump. By assumption~\textup{(iv)}, this gives the full local contribution at the critical parameter.
\end{proof}
For the degree calculation we need local families containing all blowing-up solutions at the chosen source strengths and parameters. The ordinary leading scale has a logarithmic factor. When every outer profile exponent is below one, the proof uses two-sided bounds for the scale term.
\section{The degree formula for all source strengths greater than \texorpdfstring{$-1$}{-1}}
\label{sec:all-orders}
In this section $A_p,B_p>0$ satisfy \eqref{fnd:independence}. We first extend the results for the shadow system to these shifted parameters. We then construct the local profiles and estimate the difference between each blowing-up solution and its profile. For the scale equation we distinguish
\[
 s_*>1,\qquad s_*=1,
\]
and $s_*<1$, where $s_*$ is the largest outer profile exponent. If $s_*<1$, we vary the weights away from the sources. At the end of Section~\ref{sec:wall-completion} compactness removes the rational independence assumption.

The constants in the local estimates may depend on the fixed source strengths. No uniform bound for the inverse as a scalar profile exponent approaches an integer is needed.
\subsection{Shadow systems and regular bubbles}
\begin{proof}[Proof of Proposition~\ref{univ:shadow}]
\smallskip\noindent\textit{Step 1. The effective weight parameters.} Let $b_p=1+\beta_p$. Substitution of the three possible pairs of local masses in \eqref{fnd:effective} gives
\[
 b_p=B_p,A_p+B_p,A_p,
\]
respectively. In each case $b_p>0$. Thus there is a fixed $p_0>1$ for the finite family such that every factor $|z|^{2\beta_p}$ belongs to $L^{p_0}$: when $\beta_p<0$, take
\[
 p_0(-\beta_p)<1.
\]
If $\beta_p=0$, the local factor is smooth and positive. If $\beta_p>0$, it vanishes. The moving ordinary zeros are smooth quadratic factors.

Only the second weight parameter $\beta_p$ enters the scalar integrability condition. The first modified weights enter the balance equations at free centers away from the marked points.

First, Proposition~\ref{fnd:scalar-cosets} holds for $\beta>-1$. For each fixed member of the sequence, a source obtained by an exact coincidence has total strength $\alpha=\beta+m>-1$ or a nonnegative integer strength. Its regular part is continuous and $F=O(|z-a|^{2\alpha})$. Consequently
\[
 \int_{\partial D_r(a)}\sqrt F\,ds=O(r^{\alpha+1})\longrightarrow0.
\]
The local Green representation also gives
\[
 r\sup_{|z-a|=r}|\nabla v(z)|\to0.
\]
Indeed, its singular integral is bounded by $C(1+r^{2\alpha+1})$ if $\alpha\ne-1/2$, with a logarithm in the equality case; multiplication by $r$ tends to zero. Thus the limiting source holonomy is the same matrix as before.

The curvature estimate, the integer central factors, and the outer flux calculation are unchanged. This proves the restriction on the local mass for $\beta>-1$, including coincidences and bounded smooth forcing.

\smallskip\noindent\textit{Step 2. Compactness.} Normalization and the mean-zero Green representation give the same upper bound for the normalizing constants as before.

The weights in the scalar equation are uniformly bounded in some $L^p$, $p>1$, because each fixed weight parameter satisfies $\beta_p>-1$ and the moving ordinary factors have quadratic zeros.

The estimate for small local mass therefore gives a finite concentration set and convergence away from that set.

If the normalization tended to zero, all scalar mass would be concentrated. The extended restriction on the local mass, applied to each whole cluster, would put $\gamma$ in the set excluded by \eqref{fnd:scalar-exclusion}. Hence the diffuse remainder is positive.

At a cluster its limiting density has the form $H_0|z|^{\lambda}e^\psi$, with $\lambda=2(\beta+u)-4n$ and $\psi$ bounded below. Integrability implies the strict inequality
\[
                         2n<1+\beta+u.
\]
On a sufficiently small disk the exponential of the Green potential of the residual measure belongs to every fixed finite $L^q$ space. Since
\[
 \lambda>-2,
\]
H\"older's inequality gives an $L^p$ density for some $p>1$. Elliptic regularity then bounds $\psi$ and yields
\[
 |z|\,|\nabla\psi(z)|\to0.
\]
To obtain the integral identities for each solution in the sequence, first remove a disk of radius $\epsilon$ about the prescribed point. With that member fixed, its density is $O(\epsilon^{2\beta})$ on the circle and the local Green estimate gives
\[
 \epsilon\sup_{|z|=\epsilon}|\nabla w|\to0.
\]
For either test field $x$ or $Jx$, the inner boundary term satisfies
\[
 \left|\text{inner boundary term}\right|
 \le C\epsilon^2\sup_{|z|=\epsilon}(|\nabla w|^2+F)\longrightarrow0
 \qquad(\beta>-1).
\]
The volume terms converge as well: $x\cdot\nabla\log V$ and $Jx\cdot\nabla\log V$ are bounded near the prescribed point for this fixed member. Let $\epsilon$ tend to zero before taking the sequence limit.

We thus obtain the two moments and the integral identity of Lemma~\ref{fnd:source-exclusion}. Their algebra uses only $b=1+\beta>0$ and rational independence, and still excludes every prescribed-source collision.

For separated centers, Lemma~\ref{fnd:separated-scalar-compactness} applies with each fixed weight parameter
\[
 \beta_p>-1,
\]
a quadratic zero at each ordinary center and the same residual mass $\gamma$. It gives uniformity on compact separated center sets for the indicated convergent positive smooth factors; no boundedness of the original singular weights is needed.

It remains to consider a pair of ordinary regular centers. Their disks stay away from the prescribed sources. Contributions to the potentials from the other disks have bounded derivatives by separation and the total mass bounds.

For convergence at a pair of colliding centers it suffices that the local coefficients converge in $C^1$. For the local integral identities take $C^3$ neighborhoods of the weights or use shrinking perturbations supported away from those disks.

Thus the arguments for a pair of colliding centers use the same local hypotheses. The induction in Theorem~\ref{fnd:common-weights} gives the compactness neighborhoods required for genericity. Only the local coefficients on those disks enter the regular pair expansions.

\smallskip\noindent\textit{Step 3. Nondegeneracy.} Let $X=C^0(M)\cap\{\int_Mw=0\}$ and use the fixed-point scalar map
\[
             \mathcal S(w,P)=w-2Gf(w,P).
\]
For the ordinary equations use the extension
\[
 \widehat F_i^A(w,P,Q)
       =F_{i,0}^A(P,Q)
           -\int_M\nabla_{p_i}G(p_i,x)f(w,P)(x)\,dV_g.
\]
At scalar zeros this is the original equation. On compact separated center neighborhoods these maps are $C^1$ on $X$ times the finite center space. To check the apparent singularity of the force kernel, write
\[
 E(x,p)=e^{-4\pi G(x,p)}.
\]
In a common source-free coordinate disk, $E(x,p)=|x-p|^2e^{-4\pi R(x,p)}$ is jointly smooth, and
\[
 \nabla_pG(p,x)E(x,p)=-\frac1{4\pi}\nabla_pE(x,p).
\]
The combined kernel and all its first center derivatives are therefore regular through $x=p$. At the prescribed sources they are bounded smooth multipliers of the fixed $L^p$ weight.

The scalar derivative is identity minus a compact map, and the center part is finite dimensional, so the joint derivative is Fredholm of index zero.

To evaluate the gradients in the center equations and differentiate their local Taylor coefficients, fix open sets $U_0\Subset U_1\Subset M\setminus S$ containing all moving centers in the compact separated parameter neighborhood.

Every scalar solution, and the parameter derivatives along a $C^1$ solution branch, belong to
\[
 X_{p,\eta}(U_0)
 =\left\{w\in W^{2,p}(M):\int_Mw=0\right\}
       \cap C^{2,\eta}(\overline{U_0}),
\]
with the sum of the two norms and uniform bounds for a fixed
\[
 0<\eta<1.
\]
Indeed the uniform $L^p$ bound on the weights gives a global $W^{2,p}$ estimate for the scalar equation.

On $U_1$ every coefficient, including each ordinary factor $E(x,p_i)$, is jointly smooth; its zero at $p_i$ causes no loss of regularity. Interior $W^{2,q}$ estimates for $q>2$, followed by Schauder estimates, give the claimed bound on $U_0$.

Differentiating the scalar equation gives the same estimates for parameter derivatives on each compact $C^1$ branch.

If the scalar derivative with fixed centers is invertible, the implicit function theorem in $C^0$ gives such a branch. Repeating the interior estimates supplies higher local derivatives for the coefficients in the integral identities and the traces.

Consequently $\nabla w(p_i)$ and the second derivatives at mixed centers depend continuously differentiably on the parameters along the scalar branch in this topology. These evaluations require the local estimates when $p$ is close to one.

Away from the scalar zero set we use the integral formula for $\widehat F_i^A$ on $C^0$; on the zero set it agrees with the gradient equation. The nested sets $U_0,U_1$ also keep these regularity estimates strictly interior.

Consider a functional annihilating the derivative in this $C^0$ space, including the variations of the weights.

Independent variations of the first derivatives of the first smooth factor force its center components to vanish. Represent the remaining scalar functional by a finite measure $\Lambda$ of zero total mass. Put
\[
 M_f\xi=f\left(\xi-\gamma^{-1}\int_Mf\xi\right).
\]
Annihilation of the scalar derivative gives
\[
 \Lambda=2M_fG\Lambda.
\]
The logarithmic Green potential of a finite measure belongs to every finite $L^r$ space. Since $f\in L^p$ for some $p>1$, the right-hand side belongs to $L^q$ for some $q>1$. Thus $G\Lambda\in W^{2,q}\subset C^0$. Annihilation of arbitrary variations of the second smooth factor now gives
\[
 M_fG\Lambda=0,
\]
and hence $\Lambda=0$. This proves universal surjectivity and permits the finite-dimensional Sard argument used in Lemma~\ref{fnd:genericity}.

We next perturb the weights to make the scalar kernel space trivial. The form with potential $f\in L^p$ is a compact perturbation of the Dirichlet form on mean-zero $H^1$; its kernel space is finite dimensional and, by elliptic bootstrap, consists of continuous functions smooth on $M\setminus S$.

In the proof using $\int_Bv\phi_i\phi_j$, Lemma~\ref{fnd:scalar-unique} is applied on the connected punctured surface $M\setminus S$, where the potential is locally bounded and smooth. Continuity extends the resulting vanishing across $S$. A globally bounded potential is not required.

The perturbation $v/f$ remains smooth because its support is a fixed disk disjoint from $S$ and all ordinary centers. Thus the same finite spectral perturbation makes every fixed-center scalar block invertible.

At these weights, eliminating the scalar variable gives the original reduced center equations, so the local index and its reduced orientation agree with those used in the calculation of the degree of the shadow system. Write
\[
 C=2G M_f
\]
on the mean-zero spaces. For $\lambda\ne0$, expanding $(C-\lambda I)^j\xi=0$ expresses $\xi$ as $C$ applied to a polynomial in $C$ acting on $\xi$. On $C^0$ its range lies in $W^{2,p}\subset H^1$. On $H^1$, use $H^1\subset L^r$ for finite $r$ sufficiently large that
\[
 1/q=1/p+1/r<1.
\]
Then $M_f\xi\in L^q$, so the range lies in $W^{2,q}\subset C^0$. Both Green mappings are compact into the corresponding ambient spaces and agree on their intersection. The generalized eigenspaces at every nonzero eigenvalue thus agree. In the Dirichlet inner product one has
\[
 \langle Cv,z\rangle_{\dot H^1}=2\int_M f(v-\langle v\rangle_f)(z-\langle z\rangle_f).
\]
The eigenvalues of $C$ greater than one count the negative directions of $I-C$, and give the same Leray--Schauder sign in $C^0$.

The scalar inverse also exists from mean-zero $L^p$ to $W^{2,p}$, since $-\Delta_g-2M_f$ is a compact perturbation of the mean-zero Laplacian with this same trivial kernel space.

The finite matrix after elimination is the derivative of the same reduced center equations after solving the scalar equation.

Lemma~\ref{fnd:cross-separation} also applies: a local variation of the first derivatives of the second factor at a mixed center changes the scalar and ordinary center equations by $O(r|\xi|)$ in the combined-kernel formulation and gives the prescribed direct mixed force.

The shrinking quadratic bumps in Theorem~\ref{fnd:common-weights} then have supports avoiding $S$ and the ordinary coordinates in the relevant fixed neighborhoods. The mixed potential is unchanged.

In the preceding fixed-point formulation, the perturbed joint maps converge in $C^1$ at rate $O(r^2)$ on fixed separated neighborhoods. In the scalar and combined-force numerators, the perturbation is multiplication by
\[
 e^{-b_r}-1=O(r^2)
\]
in supremum norm. The coefficients and their first center derivatives have fixed $L^p$ bounds, and all normalizing integrals have uniform positive lower bounds on these neighborhoods.

The direct ordinary force is unchanged there. The implicit function theorem at each of the finitely many invertible joint blocks gives
\[
 \|w_r-w\|_{C^0}+|P_r-P|+|Q_r-Q|=O(r^2)=o(r).
\]
The scalar equations also give $W^{2,p}$ convergence at this rate.

On fixed disks about the mixed coordinates, disjoint from $S$ and the ordinary centers, the density difference is $O(r^2)$ in supremum norm by the center and $C^0$ bounds above. Interior $W^{2,q}$ estimates for large fixed $q$ first give
\[
 w_r-w=O(r^2)
\]
in $C^{0,\eta}$ on smaller disks. The smooth center factors change by $O(r^2)$ there, while $\|e^{-b_r}-1\|_{C^{0,\eta}}=O(r^{2-\eta})$ by the supremum and gradient bounds for the bump. Schauder estimates for the difference equation on still smaller disks give
\[
       \|w_r-w\|_{C^{2,\eta}}=O(r^{2-\eta}),\qquad0<\eta<1.
\]
Here the disks and $\eta$ are fixed independently of $r$. Thus the mixed centers remain in the quadratic part of their bumps. In normalized conformal coordinates, $\Delta_g(b_r/2)=2L_l+o(1)$ at the continued center, while
\[
 (b_r)_{zz}=0,\qquad (b_r)_z=o(1)
\]
there. The trace increases by $2L_l+o(1)$, and the trace-free coefficient changes by $o(1)$ by the local $C^{2,\eta}$ estimate. The compactness induction uses $C^1$ convergence with bounded $C^2$ at a collision of two ordinary centers.

For the integral identities at pairs of centers, it uses coefficients unchanged near the previous mixed centers.

\smallskip\noindent\textit{Step 4. The shadow degree.} The homotopy to the decoupled system uses only the same two cluster moments and $b>0$. The homotopy for the degree of the center equations uses rational nonresonance and is independent of scalar regularity.

The proof of Proposition~\ref{fnd:shadow-degree} applies in $C^0$, since the fixed $L^p$ bound gives relative compactness of the Green potentials on bounded sets.

Separated scalar compactness makes its center homotopies admissible on a common scalar ball. The product rule then gives the formula in terms of the degree of the scalar equation.
\end{proof}
\begin{lemma}[Regular bubbles with integrable prescribed weights]
\label{univ:regular-only}
Suppose that the preceding hypotheses on the shadow system hold and that
\[
 \alpha_{p,i}>-1\quad(p\in S),\qquad
 \{\text{blow-up points}\}\cap S=\varnothing .
\]
If all bubbles are ordinary, assume also $L_0\ne0$. Then the construction and the local degree formulas of Section~\ref{regzero:section} apply.
\end{lemma}
\begin{proof}
\smallskip\noindent\textit{Step 1. Local estimates.} The weights are smooth and positive on the disks containing the bubbles, so the local profiles, integral identities, and equations for translations and shape parameters follow as before.

On their fixed complement the first density is bounded by $C\tau h_1$, where $\tau$ is the common first normalization; it is $O(\tau)$ in a fixed $L^p$ with $p>1$ close to one. The second residual density has a uniform $L^{p_0}$ bound for some
\[
 p_0>1.
\]
The associated bilinear form is compact on $H^1$, and the separated model inverse uses only this compactness and bounded logarithmic moments. Both follow from H\"older's inequality. For its local centered tests, take $1<q<p_0$ and
\[
 1/q=1/p_0+1/r.
\]
Put $W=G(n_*^g v)$. By H\"older's inequality, $H^1\hookrightarrow L^r$
and the elliptic estimate, for some $0<\theta<1$ we have
\[
 \|W\|_{C^{0,\theta}}\le C\|W\|_{W^{2,q}}
 \le C\|n_*^g\|_{p_0}\|v\|_{H^1}.
\]
At each free center $a$, subtract $W(a)$ in the centered pairing.
The profile moments in the proof of Lemma~\ref{regzero:inverse} give
a positive power of the local width, since
\[
 |W(x)-W(a)|\le C\|v\|_{H^1}d_g(x,a)^\theta.
\]
Higher derivatives are used only on the smooth free disks.
Lemma~\ref{hyb:covariance} compares the normalized measures.
Decrease $p_0>1$ so that $h_1,h_2\in L^{p_0}(M)$.
For a trial probability density $\nu$ and correction $\eta$ as in
\eqref{hyb:measure-stability}, put
\[
 p_0'=\frac{p_0}{p_0-1},\qquad E=\|\eta\|_{H^1},\qquad
 Z=\int_M e^\eta\nu\,dV_g .
\]
Since $\|\nu\|_{p_0}\le e^{CT}$ and $\sqrt T E=o(1)$, H\"older's inequality
and the exponential Sobolev estimate at exponent $p_0'T$ give
\[
 |Z-1|\le\|e^\eta-1\|_{L^T(\nu)}
 \le\|\nu\|_{p_0}^{1/T}\|e^\eta-1\|_{p_0'T}
 \le C\sqrt T E e^{CTE^2}=o(1).
\]
Thus $\|e^\eta/Z-1\|_{L^T(\nu)}=o(1)$, giving uniform convergence in the
local integral identities. For a blowing-up sequence choose $1<s<p_0$ and
$k$ with $1/s=1/p_0+1/k$. At an unused source write
\[
 c_i=\log\frac{\rho_i}{\int_M h_i e^{u_i}\,dV_g},\qquad
 \mathfrak n_i=h_i e^{u_i+c_i}.
\]
Choose a source disk $B$ with sufficiently small local mass. The boundary
values of $u_i+c_i$ give a uniform upper bound for its harmonic part.
On a smaller disk $B'$, the logarithmic Green estimate in Step~1 of the
proof of Proposition~\ref{fnd:signed-compactness} gives
\[
 \|e^{u_i+c_i}\|_{L^k(B')}\le C_k,\qquad
 \|\mathfrak n_i\|_{L^s(B')}\le \|h_i\|_{L^{p_0}(B')}C_k.
\]
Choose $1<p<s$ close enough to one for the regular-disk estimates.
With the notation of \eqref{reg:actual-diffuse}, these bounds give uniform
integrability of $|d_D|^p$ and $|e_1|^p$. Qualitative convergence off the
centers now gives
\[
 \|d_D-n_*^g\|_p+\|e_1\|_p=o(1).
\]
Subtract \eqref{reg:shadow-scalar} from \eqref{reg:actual-diffuse}.
The mean-zero elliptic estimate places $w_u$ in the neighborhood of
the nondegenerate branch required in Lemma~\ref{journal:scalar-neighborhood}\textup{(i)},
\[
 \|w_u-w(P)\|_{W^{2,p}}\le C\bigl(\|d_D-n_*^g\|_p+\|e_1\|_p\bigr)=o(1).
\]

\smallskip\noindent\textit{Step 2. Estimates outside the regular disks.}\quad
Let $X$ be the exterior of the fixed bubble disks and their cut-offs.
Truncate each fitted density on a smaller concentric disk about $a$.
For its measure $\xi_a$ and limiting mass $m_a$, separation gives,
for each fixed $j\ge0$,
\[
 \|G\xi_a-m_aG(\cdot,a)\|_{C^j(X)}
 \le C_j\left(|\xi_a(M)-m_a|
          +\int_M d_g(y,a)\,d\xi_a(y)\right).
\]
Scaling and decay bound the first moment by a positive power of the fitted
radius $r_a$. Lemma~\ref{reg:local-power} controls the mass and total
variation errors. On the omitted annuli the corresponding $L^p$ bound
follows by decay and interpolation. Put
\[
 D=\max_a r_a,\qquad Z_2^0=-4\pi\sum_{i=1}^U G(\cdot,p_i),\qquad
 H_P=h_2e^{Z_2^0},
\]
where $p_i$ are the fitted ordinary centers. The two masses at each mixed
center cancel in $Z_2^0$. Summing over the individual radii gives, for some
$c>0$ and $N\ge0$,
\[
 \|Z_2-Z_2^0\|_{C^0(X)}\le C D^c(1+|\log D|)^N .
\]
On $X$ we have $A_D=h_2e^{Z_2}$. The exterior exponentials are uniformly
bounded, and the first normalization is $\tau$. Hence
\[
 \|A_D-H_P\|_{L^p(X)}\le C\|h_2\|_{L^p(X)}D^c(1+|\log D|)^N,
 \qquad \|e_1\|_{L^p(X)}\le C\tau\|h_1\|_{L^p(X)} .
\]
On the regular disks the weights are smooth. Step~1 of the proof of
Theorem~\ref{regzero:transfer} gives the remaining coefficient and mass
estimates for Lemma~\ref{journal:scalar-neighborhood}\textup{(ii)}.

\smallskip\noindent\textit{Step 3. The projected equation.}\quad
By Step~4 of the same proof, at $\mu=0$ we have on $X$
\[
 -\Delta_g\widehat u^{\rm app}=K\binom{-4\pi m}{n_*^g-\rho_2},
 \qquad N_2=n_*^g,\qquad \mathscr R=-K\binom{N_1}{0},
\]
where $\mathscr R$ is the residual before mass normalization.
Write $\lambda=\log t$ and let $\mathcal D_a$ denote a center derivative
scaled by its translation width. The prescribed weights are fixed, and only
the separated Green fields and the scalar branch depend on the centers. Thus
\[
 \|D_*\mathscr R\|_{L^1(X)}+\|D_*\mathscr R\|_{L^p(X)}
 \le C\tau,\qquad D_*\in\{I,\partial_\lambda,\mathcal D_a\}.
\]
Use $L^p\hookrightarrow H^{-1}$ on $X$ and the original $L^{6/5}$ estimates
where required on the regular disks. The normalized residual in
\eqref{reg:normalized-residual} has zero mean in each component, so its
pairing with any tangent $Y$ satisfies
\[
 \int_M K^{-1}\mathcal R^{\rm app}\cdot\Pi Y\,dV_g
 =\int_M K^{-1}\mathcal R^{\rm app}\cdot Y\,dV_g .
\]
The mass and mean terms are retained in \eqref{reg:ordinary-trial-variation}.
For an ordinary center $a$ set $Y_a=\mathcal D_a\widehat u^{\rm app}$.
Since $\tau=t^2$,
\[
 \partial_\lambda\widehat u^{\rm app}=0\quad\text{on }X,\qquad
 \|Y_a\|_{L^\infty(X)}\le Ct,\qquad
 \left|\int_X K^{-1}\mathscr R\cdot Y_a\,dV_g\right|
 \le C\tau t=Ct^3 .
\]
Changes of order $O(\tau)$ in the unnormalized masses give at most $O(\tau\sqrt T)$ in the dual norm.

Indeed the normalized probabilities have mass one and a fixed exponential $L^{p_0}$ bound, hence bounded $L^{1+1/T}$ norm, and the dual Sobolev bound is $O(\sqrt T)$.

Their scale derivatives obey the same bound by the bounded logarithmic derivatives of the trial densities. These errors have the required orders.

The leading terms of the scale equations have orders $t$ and $t^2\log(1/t)$, respectively. The normalization and complementary homotopy may therefore use the same correction balls and the same source and degree orientations.
\end{proof}
\subsection{The local scale moments}
\begin{lemma}[Masses and radial variation of the matched local family]
\label{univ:local-scale-moments}
In the notation of Lemma~\ref{univ:profile-estimates}, fix $d_0>0$ and the coefficients of the first-order terms. Set
\[
 \tau=R^{2s}=d_0^s e^{-T},\qquad J_i\in C^2(\overline{D_d}),
 \qquad J_i(0)=0,\quad\nabla J_i(0)=0\quad(i=1,2).
\]
On any sufficiently small fixed disk $D_d$ the following hold.
\begin{itemize}
\item[\textup{(i)}]
The two masses satisfy
\[
 \int_{D_d}n_{1,R}=4\pi s+O(\tau),\qquad
 \int_{D_d}n_{2,R}=4\pi B+\int_{D_d}f+O(\tau).
\]
Their error terms have $T$ derivatives of order $O(\tau)$.
\item[\textup{(ii)}]
Write $\kappa_1=\Delta J_1(0)$. For $s>1$, put $J_s=\int_{\mathbb R^2}|y|^2m_{s,1}(y)\,dA_y$. Then
\begin{align*}
 -\int_{D_d}(\partial_Tn_{1,R})(e^{J_1}-1)
     &=\frac{\kappa_1J_s}{4s}R^2+o(R^2),\\
 \int_{D_d}(n_{2,R}-f)(e^{J_2}-1)
                         \partial_T\log n_{2,R}&=o(R^2).
\end{align*}
If $0<s<1$, both integrals are $O(\tau)$.
\end{itemize}
All estimates are uniform on compact sets of data if the second derivatives of the $J_i$ have a common modulus of continuity.
\end{lemma}
\begin{proof}
\smallskip\noindent\textit{Step 1. The two masses in (i).}\quad On a punctured disk use $\Delta\log F_R=n_{1,R}$ and $\Delta\log Q_R=n_{2,R}$. On the fixed boundary,
\[
 \log F_R=2s\log r+O(\tau),\qquad
 \log Q_R=\log\tau+(2s-2)\log r+\log D+O(\tau).
\]
The explicit formulas with positive denominators show that the errors and their spatial derivatives of fixed order are $O(\tau)$ there, also after one $T$ derivative. The inner circular fluxes tend to $0$ and $2\pi(2A-2)$, respectively.

Their limiting values are independent of $T$, since $F_R$ is bounded at zero and $Q_R/r^{2A-2}$ has a positive finite limit for fixed $R$. Thus the mass integrals are the outer fluxes minus these inner fluxes. Since
\[
 \Delta\log D=f
\]
and $\log D$ has zero inner flux, this gives $4\pi s+O(\tau)$ and $4\pi B+\int f+O(\tau)$ with differentiated errors of the same order.

\smallskip\noindent\textit{Step 2. Differentiation in scale.}\quad Differentiation of the explicit formulas gives
\[
 \partial_T\log n_{1,R}
 =-1-\frac{H_R}{D+H_R}
               +\frac{2\tau F_0}{\tau F_0+r^{2s}},\qquad
 \partial_T\log n_{2,R}
 =-\frac{\tau F_0}{\tau F_0+r^{2s}}
               +\frac{2H_R}{D+H_R}.
\]
Both are uniformly bounded. If $s>1$, the pointwise comparisons used in Lemma~\ref{univ:profile-estimates}, including their differentiated versions, imply
\[
 \int_{D_d}r^2|\partial_Tn_{1,R}-\partial_Tm_{s,R}|=o(R^2).
\]
Indeed the three error majorants are $r^3m_{s,R}$, $r^{2+2A}m_{s,R}$, and $r^2m_{s,R}(\ell/r)^{2B}$. For any $b>-2s$, their radial integrals use
\[
 \int_{D_d}r^b m_{s,R}\,dA
 \le C R^{\min(b,2s)}
             (1+|\log R|),
\]
by splitting the integral at $R$. The logarithm is needed only for $b=2s$. In the last term use $b=2-2B>-2s$ and multiply by
\[
 \ell^{2B}\asymp R^{2s}.
\]
This yields the exponents $\min(3,2s)$, $\min(2+2A,2s)$, and $\min(2+2A,4s)$, respectively, all strictly greater than two. Taylor expansion gives
\[
 J_1=\kappa_1r^2/4+\Re(az^2)+o(r^2).
\]
The harmonic quadratic term integrates to zero against $\partial_Tm_{s,R}$. Moreover
\[
 \int_{D_d}r^2\partial_Tm_{s,R}
 =-\frac{J_s}{s}R^2+O(\tau).
\]
The uniform tail bound for the rescaled scalar second moment controls the $o(r^2)$ remainder. Replacing $J_1$ by $e^{J_1}-1$ produces an error $o(R^2)$ by the same argument, since $J_1^2=O(r^4)$. This proves the first scale identity.

\smallskip\noindent\textit{Step 3. The second component.}\quad For the second identity we use the bound for the logarithmic derivative and $e^{J_2}-1=O(r^2)$. The pointwise estimate for $n_{2,R}-f$ gives
\[
 \int_{D_d}r^2|n_{2,R}-f|
 \le C\{\ell^{\min\{2,2B\}}
       +\ell^{\min\{2A+2,2B\}}\}(1+|\log\ell|).
\]
This is $o(R^2)$ when $s>1$: $\ell/R\to0$, $\ell^{2B}\asymp R^{2s}$, and $2A+2>2$. If $s<1$, then $B<s<1$. For the first integral use $|\partial_Tn_{1,R}|\le Cn_{1,R}$ and the same radial majorants. Its leading second moment is $O(R^{2s})$, and $r^2m_{s,R}(\ell/r)^{2B}$ has integral
\[
 O(R^{\min(2+2A,4s)}(1+|\log R|))=o(R^{2s}).
\]
For the second integral, $r^2 I$ and the diffuse difference in the pointwise bound of Lemma~\ref{univ:profile-estimates} each have integral $O(\ell^{2B})$: in their outer integrals the powers are $1-2B$ and $1+2A-2B$, both greater than $-1$. As $\ell^{2B}\asymp R^{2s}$, both scale integrals are $O(\tau)$.
\end{proof}
The first identity enters the normalized energy derivative after the first mass discrepancies have been summed over all disks. Use the unnormalized densities before subtracting the surface means; normalizing the densities makes the final equation independent of those means.

Since the normalized equation residual has zero integral, its energy pairing can equally be evaluated on the tangent before mean subtraction. Let $M_T$ be the total unnormalized first mass and $M$ its prescribed value at
\[
 \mu=0.
\]
At $\mu=0$, its normalization adds $(1-M/M_T)\mathfrak n_1$ to the first density residual. On a mixed singular disk the first field tangent is $1+\partial_T\log n_{1,R}$. Globally this average is exactly $1+M_T'/M_T$, since the first unnormalized normalization is $e^{-T}$ and a prime denotes $\partial_T$. It is therefore
\[
 1+O(\mathcal L)
\]
when $s_*>1$, where
\[
 \mathcal L=e^{-T/s_*}
\]
is the largest scale and the flux and moment bounds give
\[
 M_T-M,M_T'=O(\mathcal L).
\]
The first normalization thus contributes
\[
                    (M_T-M)+O(\mathcal L^2).
\]
The constant part cancels the sum of $-\int n_{1,R}(e^{J_1}-1)$ and the corresponding terms on the other disks, leaving the local expression
\[
 -\int(\partial_Tn_{1,R})(e^{J_1}-1)
\]
together with the $O(e^{-T})$ exterior and profile-mass terms. This is a global cancellation; with different outer profile exponents, the normalization remainder is $o(\mathcal L)$, not necessarily $o(R_a^2)$ for every individual bubble. The second field tangent on a mixed singular disk is $\partial_T\log n_{2,R}$.

If $N_T$ is the total unnormalized second mass, its global normalized tangent average is $N_T'/N_T$, since the multiplicative normalization of the second unnormalized density is independent of $T$. At a mixed regular disk the boundary formula
\[
 \log Q_0-\log(4R^4)=2\log d+\log(1+\beta d^2/4)+O(R^2)
\]
holds on the fixed circle, also after one radial or $T$ derivative. Put $n_{\rm bg}(z)=16\beta/(4+\beta|z|^2)^2$. Integration of the radial flux gives
\[
 \int_{D_d}n_2=4\pi+\int_{D_d}n_{\rm bg}+O_d(R^2).
\]
The background density on this disk is $n_{\rm bg}e^{J_2}$. By evenness and the second-moment bound in Lemma~\ref{reg:residual},
\[
 \int_{D_d}(n_2-n_{\rm bg})e^{J_2}=4\pi+O_d(R^2).
\]
At fixed matching data the error has $T$ derivative $O_d(R^2)$ by the differentiated bounds in the proof of Lemma~\ref{hyb:trial-energy}. Together with the other disk and exterior bounds this gives
\[
 N_T-\rho_2,N_T'=O(\mathcal L).
\]
Consequently the second normalization term is $o(\mathcal L)$. These are the two normalization contributions used with the local formulas above. When all outer profile exponents are below one, the same estimates give
\[
 M_T-M=O(e^{-T})
\]
and a total normalization contribution of order $e^{-T}$, which must be retained in the scale estimate.

On a fixed cut-off annulus for a radial or prescribed patch, the patch correction, its $T$ derivative and its residual are $O(\tau)$, so the direct energy pairing is $O(\tau^2)$. At a regular mixed point the formulas for $F_R$ and $Q_R$ bound the same factors by $O(R^2)$, giving
\[
 O(R^4)=O(\tau)
\]
in the energy pairing. The latter patch has outer profile exponent two, so this error is $o(\mathcal L)$. Exterior first-density and integrated profile-mass terms can still contribute $O(\tau)$.

They are negligible when the largest outer profile exponent exceeds one, and enter the leading-order estimate when all outer profile exponents are below one.
\subsection{The projected equation and the ordinary scale}
\begin{lemma}[Pure singular disks]
\label{univ:pure-residual}
At a pure singular point with $A,B>0$ and $A\ne 1$, set
\[
 q_d=A+B,\qquad R^{2A}=d_0^A\tau,\qquad \tau=e^{-T}.
\]
Keep $d_0>0$ fixed when differentiating in $T$. For some fixed $p>1$ sufficiently close to one, there is a local approximation satisfying
\[
 \|\mathscr R\|_{H^{-1}}
 +\|\partial_T\mathscr R\|_{H^{-1}}
 \le
 \begin{cases}
  C\tau,&0<A<1,\\
  C(R^{2/p}+\tau),&A>1,
 \end{cases}
\]
The fixed cut-off annuli contribute $O(\tau)$. The constants are uniform on compact sets of background data.
\end{lemma}
\begin{proof}
\smallskip\noindent\textit{Step 1. The residual density.}\quad The residual density is $g=O(r^{2q_d-2})$, and the scalar outer profile is $S=m_{A,R}$. Before adding the angular correction, the densities are $Se^{J_1}$ and $g(1+h)$, where
\[
 h=(R/r)^{2A}.
\]
The two excess densities in the unnormalized equation are
\[
 E_1=S(e^{J_1}-1),\qquad E_2=gh.
\]
With the convention
\[
 \mathscr R=-\Delta_g\widehat B_T-K(N-\rho^0)
\]
they contribute $-K(E_1,E_2)^{\mathsf T}$ to the differential-equation residual and $-(E_1,E_2)^{\mathsf T}$ to the residual before the Cartan factor. Here $J_1(0)=0$. The right-hand side of the Poisson equation for the background potential is $O(r^{2q_d-2})+O(1)$.

\smallskip\noindent\textit{Step 2. The choice of the exponent.}\quad If $A<1/2$, choose $2A<\nu<\min\{1,2q_d\}$. The local logarithmic potential estimate gives $J_1=O(r^\nu)$. If $1/2\le A<1$, then $q_d>1/2$ and
\[
 J_1=g_1\cdot x+O(r^\nu),\qquad
 2A<\nu<\min\{2,2q_d\}.
\]
For $1/2\le A<1$, add $(P,-P/2)$, where
\[
 P=-\frac{2A}{A-1}
              \frac{g_1\cdot x}{1+(r/R)^{2A}},\qquad
 (-\Delta-2S)P=2Sg_1\cdot x.
\]
Since $|P|\le Cr$, the first residual is bounded by $CSr^\nu$. Choose $p>1$ so close to one that $2-2/p<\nu-2A$ and $p(1-B)<1$. Splitting at $r=R$ gives
\[
 \|Sr^\nu\|_p\le C\tau,\qquad \|gh\|_p\le C\tau.
\]
In the range $A<1/2$ put $P=0$. In both cases
\[
 0<A<1,
\]
the second residual is $g\{(1+h)e^{-P/2}-1\}$. In addition to $gh$, it is bounded by $C(g|P|+gh|P|)$. On $r\le R$, the first of these is $O(r^{2q_d-1})$, while on $r\ge R$ it is $O(\tau r^{2B-1})$. For $1<p<2$ both have $L^p$ norm $O(\tau)$, because $2B-1+2/p>0$. The same bound holds for $gh|P|\le C\tau r^{2B-1}$.

\smallskip\noindent\textit{Step 3. The case $A>1$.}\quad If $A>1$, then
\[
 q_d>1,
\]
and the background is $C^2$ at the point. Use the same first-order correction. The first residual is bounded by $CSr^2$, whose $L^p$ norm is $O(R^{2/p}+\tau)$. The preceding second-component estimates remain valid for every $B>0$ after choosing $p$ sufficiently close to one.

This proves the stated bounds, since $L^p$ embeds into $H^{-1}$ on the compact surface. Logarithmic-scale differentiation preserves the rational majorants. On a fixed cut-off annulus $P$, the scalar profile correction, and their fixed spatial derivatives are $O(\tau)$.
\end{proof}
\begin{lemma}[An ordinary leading scale]
\label{univ:ordinary-dominated}
Fix a separated nondegenerate shadow solution with invertible scalar linearization and at least one prescribed bubble. Assume
\[
 U>0,\qquad V=0,\qquad s_a<1\quad\text{at each prescribed bubble}.
\]
Let $X$ be the ordered ordinary centers near their limiting positions $X^0$, and $F(X)$ their balance map on the scalar branch. Write
\[
 t=e^{-T/2},\qquad L=\log(1/t),\qquad s_i=t d_i(X),\qquad d_i(X)>0,
\]
and set
\[
 \kappa_i(X)=\Delta_g\log h_1(p_i)+2M-\rho_2-2K_g(p_i),
 \qquad
 \Lambda(X)=2\pi\sum_i d_i(X)^2\kappa_i(X),\qquad
 L_0=\Lambda(X^0).
\]
Take one radial reference source at a prescribed bubble in the source space for the mixed configuration.
\begin{itemize}
\item[\textup{(i)}]
For fixed sufficiently small $\delta>0$, the projected family exists on a fixed neighborhood of $X^0$ for $|\mu|\le e^{-\delta T}$. It is unique in the ball $\|\eta\|_{H^1}\le e^{-\delta T}T^{C_0}$, for $C_0$ sufficiently large. Moreover,
\[
 |X-X^0|\le C_0\sqrt t,\quad |\mu|\le C_1t^2L
 \quad\Longrightarrow\quad \|\eta\|_{\Hh}=O(t^{7/6}L).
\]
\item[\textup{(ii)}]
Suppose $L_0\ne0$. On the side $\mu L_0>0$, the local Leray--Schauder degree of the full Toda map on a bounded isolating neighborhood is
\[
 (-1)^{U+|P|}\operatorname{sgn}(L_0)d_S(Z).
\]
\item[\textup{(iii)}]
Under the same assumption $L_0\ne0$, consider a sequence of solutions blowing up at this configuration. If $|\mu|+\|\eta\|_{H^1}=O(t^\delta)$ after modulation for some $\delta>0$, then
\[
 |F(X)|+|X-X^0|=O(t^{5/3}L^3),\qquad
 \|\eta\|_{H^1}=O(t^{5/3}L^2),\qquad
 \mu=L_0t^2L+O(t^2).
\]
The sequence lies strictly inside the isolating domain on its side of the critical value.
\end{itemize}
\end{lemma}
\begin{proposition}[The global projected equation]
\label{univ:projected-construction}
Fix a compact family of separated shadow configurations with at least one prescribed singular bubble. Assume the following.
\begin{itemize}
\item[\textup{(a)}]
The residual second mass is $\gamma=\rho_2-4\pi\sum_a t_a>0$. The scalar branch is compact and its linearization with the centers fixed is uniformly invertible.
\item[\textup{(b)}]
The effective weight parameters give a uniform $L^{p_0}$ bound for the residual density. The positive coefficients in the matched scales are bounded above and away from zero.
\item[\textup{(c)}]
At each center use the Taylor coefficients of the corresponding shadow background. All active scalar profile exponents at the prescribed singular points are positive and noninteger.
\end{itemize}
Use Lemma~\ref{univ:profile-estimates} for the exact mixed singular profiles, Lemma~\ref{univ:pure-residual} for the pure singular approximations, and Section~\ref{reg:section} at regular points. Set $s_*=\max_a s_a$.
\begin{itemize}
\item[\textup{(i)}]
There are a mean-zero approximation $B_T$ and $0<\varepsilon_T\le Ce^{-aT}T^C$, $a>0$, such that, for $j=0,1$,
\[
 \|\partial_T^j\mathcal F_\mu(B_T)\|_{\Hh}
       \le C(\varepsilon_T+|\mu|)T^C,
 \qquad
 \|D\mathcal F_\mu(B_T)-\mathcal T_T\|
       \le C(e^{-aT}+|\mu|)T^C.
\]
If $s_*\ne1$, one may choose $\varepsilon_T$ so that
\[
 \varepsilon_T^2T^N=o(\mathcal L)
 \quad\hbox{for every fixed }N,
 \qquad
 \mathcal L=
 \begin{cases}
 e^{-T/s_*},&s_*>1,\\
 e^{-T},&s_*<1.
 \end{cases}
\]
\item[\textup{(ii)}]
Fix $\delta>0$ sufficiently small and $C_0$ sufficiently large. For $T$ sufficiently large and $|\mu|\le e^{-\delta T}$, the projected equation has a unique solution in the ball of radius $e^{-\delta T}T^{C_0}$:
\[
 \eta\in E_T,\qquad P_{E_T}\mathcal F_\mu(B_T+\eta)=0,
 \qquad
 \mathcal F_\mu(B_T+\eta)=\sum_\alpha c_\alpha S_\alpha.
\]
For some $a'>0$,
\[
 \|\eta\|_{H^1}+\|\partial_T\eta\|_{H^1}+|c|+|\partial_Tc|
       \le C(\varepsilon_T+|\mu|)T^C,
 \qquad
 \|\eta\|_{C^0}+\|\partial_T\eta\|_{C^0}\le Ce^{-a'T}.
\]
\item[\textup{(iii)}]
The projected second derivative has the index of the model in Proposition~\ref{signed:model-inverse}. The limiting orientation of the coordinates in \eqref{hyb:chart} is positive.
\item[\textup{(iv)}]
Suppose there are no free centers. Let $\varepsilon$ vary the first weight away from the prescribed points and their disks. At each finite sufficiently large $T$, the projected family is smooth in $(T,\mu,\varepsilon)$.
\end{itemize}
In \textup{(iv)} the mixed derivatives of the reduced energy exist. Their existence is asserted for fixed $T$.
\end{proposition}
\begin{proof}
\smallskip\noindent\textit{Step 1. Estimates for the approximation.} If $\mathscr R_T=-\Delta_gB_T-K(f_T-\rho)$ is the mean-zero residual in the two differential equations, then
\[
 \mathcal F_\mu(B_T)=G\mathscr R_T,
\]
with $G$ acting componentwise. If the density residual is written before its Cartan factor as $\mathfrak r_T=K^{-1}\mathscr R_T$, this identity is
\[
 \mathcal F_\mu(B_T)=KG\mathfrak r_T.
\]
The bounded Green map from $H^{-1}$ to $\Hh$ identifies the density-residual estimates below with the $\Hh$ estimates required by the inverse of the projected linearization on $E_T$.

The same assertion holds after a $T$ derivative, since the Green map is fixed. We first verify these approximation estimates.

At a mixed singular point, Lemma~\ref{univ:profile-estimates} and the exact residual identity give $O(R^{2/p}+\ell^{2/p}+\tau)$, allowing fixed logarithmic factors. At a pure point apply Lemma~\ref{univ:pure-residual} with the fixed noninteger $A>0$, $B>0$ and the exact matching law
\[
 R^{2A}=d_0^A e^{-T}.
\]
The background data range over a compact family. We choose the same $p$ at all points to satisfy the strict inequalities on the exponents in that lemma. At a regular mixed point the existing bound is
\[
 O(R^{4/p-2});
\]
at an ordinary point it is $O(R^{2/p-1})$ on a fixed center domain. The common normalization is
\[
 \tau\asymp e^{-T}.
\]
Only finitely many profile exponents occur, so one $p>1$ sufficiently close to one satisfies all local integrability conditions. If $s_*>1$, also choose
\[
 2/p-1>1/s_*.
\]
These inequalities imply that the square of the sum of the residual bounds is $o(e^{-T/s_*})$, even after multiplication by any fixed power of $T$. If all outer profile exponents are less than one, there are no regular bubbles. At each mixed point choose $ps_a<1$ and
\[
 pB_a<1.
\]
The estimates for
\[
 ps_a<1,\qquad pB_a<1
\]
in Lemma~\ref{univ:profile-estimates} and the exact residual identity then give $O(\tau)$; the pure negative residuals have this bound already. Their squares are $o(e^{-T})$. If $s_*=1$, the ordinary residuals on a fixed center neighborhood are only asserted to be exponentially small in $T$.

The prescribed remainders are $O(e^{-T})$ by the same strict integrability choices for each profile exponent below one. These are sufficient for the preliminary projected construction. The comparison at the ordinary leading scale is made on the refined domain in Lemma~\ref{univ:ordinary-dominated}.

On a fixed cut-off annulus a radial or prescribed patch gives an $O(\tau)$ residual. At a regular mixed point the formulas for $F_R$ and $Q_R$ give an $O(R^2)$ residual on this annulus, which is already included in its disk estimate $O(R^{4/p-2})$.

The exterior first density and its $T$ derivative are $O(\tau)$ in a fixed $L^p$. The unnormalized mass discrepancies and their $T$ derivatives are bounded by the same residual sum, by integration of the equations.

For a normalized probability density $v_T$ appearing here, there is $p_0>1$ such that
\[
 \|v_T\|_{p_0}\le Ce^{bT}.
\]
Indeed the source strengths at each prescribed point are $A-1$ and $B-1$, both greater than $-1$. The normalizing masses are bounded above and below. Interpolation with the $L^1$ norm gives a uniform $L^{1+1/T}$ bound and hence
\[
 \|v_T-1\|_{H^{-1}}\le C\sqrt T.
\]
The normalization errors therefore add only powers of $T$. All these statements remain valid after one $T$ derivative. To compute this derivative, retain the unnormalized densities
\[
 N_i=e^{c_i}h_i e^{\widehat B_{T,i}},
\]
their masses $M_i=\int_MN_i$, and the probabilities $v_i=N_i/M_i$. Then
\[
 \partial_Tv_i=v_i\left(\partial_T\log N_i
                   -\int_Mv_i\partial_T\log N_i\right),
 \qquad c_1'=-1,\quad c_2'=0.
\]
The explicit local formulas and fixed cut-off patches give
\[
 |\partial_T\log N_i|\le C;
\]
in particular $|\partial_Tv_i|\le Cv_i$. The constant $c_1'$ cancels in this probability derivative but is retained in the unnormalized mass and energy identities. Thus differentiating the denominators introduces no extra unestimated term.

The local density comparisons are a positive power of the local radii; the pure and regular comparisons have the same property. Thus the normalized densities differ from the model probabilities by $O(e^{-aT})$ in $L^1$, and all have exponentially bounded $L^{p_0}$ norms.

Lemma~\ref{hyb:covariance} gives the operator comparison. The scalar residual density has a uniform $L^{p_0}$ bound, so Proposition~\ref{signed:model-inverse} applies. Thus
\[
 P_{E_T}D\mathcal F_\mu(B_T)|_{E_T}:E_T\to E_T
\]
has an inverse of norm at most $CT$.

\smallskip\noindent\textit{Step 2. The projected equation.} For the nonlinear term we obtain a bound polynomial in $T$. Put
\[
 p_T=1+1/T
\]
and let $v_T$ be any of the normalized approximate densities. H\"older's inequality and the two-dimensional Sobolev and exponential inequalities give, when $\sqrt T\|\eta\|_{H^1}$ is bounded,
\[
 \int_Mv_Te^{|\eta|}|z_1z_2z_3|
 \le CT^{3/2}e^{CT\|\eta\|_{H^1}^2}
                         \prod_{j=1}^3\|z_j\|_{H^1}.
\]
To see the power, apply H\"older first with $p_T$ and then with four equal exponents to the exponential and the three functions. Their integrability exponents are $4(T+1)$.

Each function contributes $C\sqrt T$; the exponential contributes $e^{CT\|\eta\|_{H^1}^2}$. Jensen's inequality bounds the normalization $\int v_Te^\eta$ away from zero, since
\[
 |\int v_T\eta|\le C\sqrt T\|\eta\|_{H^1}.
\]
The third derivative of the logarithm of this normalization therefore has trilinear norm $CT^{3/2}$ on the small ball. In particular,
\[
 \|\mathcal N(\eta)-\mathcal N(\xi)\|_{\Hh}
 \le CT^{3/2}(\|\eta\|_{H^1}+\|\xi\|_{H^1})
                                      \|\eta-\xi\|_{H^1}.
\]
Choose $\delta<a$ and $C_0$ larger than the polynomial powers in the inverse residual bound. The inverse equation maps the ball of radius $e^{-\delta T}T^{C_0}$ to itself and is a contraction there, because its Lipschitz constant is bounded by $CT^{5/2}e^{-\delta T}T^{C_0}$.

Its fixed point obeys the sharper bound
\[
 C(\varepsilon_T+|\mu|)T^C
\]
by absorbing its quadratic remainder. This proves both existence and uniqueness on the enlarged ball. The matrix $(\langle S_\alpha,S_\beta\rangle_{\Hh})_{\alpha,\beta}$ and its inverse are uniformly bounded. Hence the multipliers satisfy the same estimate as the correction.

The local formulas and their source functionals have bounded first logarithmic-scale derivatives. The same holds for the projection $P_{E_T}$, because this matrix of inner products is uniformly invertible. Every source is centered on its own disk, so
\[
 \int_Mh_\alpha=\int_M\partial_Th_\alpha=0.
\]
The differentiated centered sources are products of the local measures with bounded factors from the rescaled kernel functions, up to the cut-off errors. The weighted variance estimate therefore bounds both their $H^{-1}$ norms and the $H^1$ norms of their Green representatives.

This proves differentiability in $T$ of the matrix of inner products and of $P_{E_T}$, with uniform bounds. Differentiate the projected equation at fixed $\mu$. The derivative of the trial residual has the same bound.

The term $(\partial_TP_{E_T})\sum_\alpha c_\alpha S_\alpha$ is $O(|c|)$, and the difference of second derivatives is bounded by a fixed power of $T$ times $\|\eta\|_{H^1}$. The derivative $\partial_T\eta$ need not lie in $E_T$; differentiating
\[
 L_\alpha\eta=0
\]
gives
\[
 L_\alpha\partial_T\eta=-(\partial_TL_\alpha)\eta.
\]
The inverse matrix of the inner products of the source vectors therefore bounds the component of $\partial_T\eta$ orthogonal to $E_T$ by
\[
 C\|\eta\|_{H^1}.
\]
Move this controlled component to the right-hand side and invert the projected linearization on $E_T$ to estimate the component in $E_T$. This proves the estimates for $\partial_T\eta$ and $\partial_Tc$.

\smallskip\noindent\textit{Step 3. Uniform convergence of the correction.} We use the integrability of the density to estimate the correction in $C^0$. For every sufficiently small $b_1>0$ there is a fixed $p_1>1$, chosen sufficiently close to one, such that
\[
 \|v_T\|_{p_1}\le Ce^{b_1T}.
\]
This follows by interpolating the $L^1$ and $L^{p_0}$ bounds. Choose $1<p<p_1$ and apply H\"older to $v_T(e^\eta-1)$; the remaining exponential factor is estimated in a fixed finite exponent.

Its norm is $O(\|\eta\|_{H^1})$ on the exponentially small ball. The same $L^{p_1}$ bound holds for the source densities. The equation for $\eta$ and elliptic $W^{2,p}$ estimates give
\[
 \|\eta\|_{W^{2,p}}
 \le Ce^{b_1T}
   \{\|\eta\|_{H^1}+|c|+|\mu|\}
       +\|\mathscr R\|_p.
\]
All local residuals are exponentially small in this fixed $L^p$, after choosing $p$ sufficiently close to one. Take $b_1$ smaller than their decay exponents and that of the correction. The Sobolev embedding $W^{2,p}\hookrightarrow C^0$ proves the first $C^0$ bound.

To differentiate the normalized correction, write a trial probability as $v$ and $Z=\int_M ve^\eta$. Since $\int_M\partial_Tv=0$,
\[
 Z_T=\int_M(\partial_Tv)(e^\eta-1)
                   +\int_Mve^\eta\eta_T .
\]
The derivative of $v(e^\eta/Z-1)$ has $L^p$ norm at most $Ce^{b_1T}(\|\eta\|_{H^1}+\|\eta_T\|_{H^1})$ by the same fixed-exponent H\"older estimate, the bound
\[
 |\partial_Tv|\le Cv,
\]
and the positive lower bound on $Z$. The differentiated source is a sum of $\partial_Tc_\alpha\,h_\alpha+c_\alpha\,\partial_Th_\alpha$; both source factors have the fixed $L^{p_1}$ bound with the chosen exponent loss.

The residual derivative is exponentially small in fixed $L^p$. Applying $W^{2,p}$ to this differentiated equation gives the asserted $C^0$ bound for $\partial_T\eta$.

\smallskip\noindent\textit{Step 4. Dependence on the parameters.} We prove smooth dependence for finite $T$ when there are no free centers. The scalar shadow is independent of the exterior first-weight perturbation.

Its values and the matched Taylor coefficients at the prescribed points are unchanged, since that perturbation is supported outside the patch disks. Thus $B_T$ and the centered source space $E_T$ are independent of $\mu$ and $\varepsilon$, and depend smoothly on $T$ for every finite $T$.

Their explicit rational profile formulas give this smoothness in $H^1$ and in the source $H^{-1}$ norms. Locally in $T$, use the smooth projections $P_{E_T}$ to identify these moving closed subspaces with one fixed subspace.

In these fixed coordinates the normalized nonlinear map is smooth: the fixed singular weights belong to some $L^{p_0}$, and the exponential Sobolev estimates justify its derivatives of every fixed finite order. Its derivative in the correction variable is invertible.

The parameter implicit function theorem therefore makes $\eta$ and its multipliers smooth in $(T,\mu,\varepsilon)$ at every finite sufficiently large $T$. This supplies, in particular, $C^2$ regularity for the reduced energy and the interchange of its mixed derivatives.

Only the first $T$ derivative is estimated uniformly in the preceding argument; finite-parameter smoothness does not claim further uniform derivative bounds as
\[
 T\to\infty.
\]

The change in the second derivative is exponentially small, hence $o(T^{-1})$, so the index is unchanged.

At every local scale the curved singular profile converges to the same radial scalar profile; the added first derivative parameters contribute terms tending to zero on that scale. Consequently the leading source--parameter pairings are those of \eqref{hyb:chart}. For the correction use
\[
 L_\alpha\partial_p\eta=-(\partial_pL_\alpha)\eta.
\]
The source derivatives are bounded in the scaled coordinates, so this term tends to zero. This proves the asserted orientation.
\end{proof}
\begin{proof}[Proof of Lemma~\ref{univ:ordinary-dominated}]
When differentiating in the centers, we work on the fixed separated scalar branch and keep the prescribed points fixed. Write
\[
 f=|F(X)|,\qquad m=|\mu|,
\]
$E=\|\eta\|_{\Hh}$, and $\lambda=\log t$. If $\widehat B$ is the approximation before subtracting the surface means, put
\[
 N_i=e^{c_i}h_i e^{\widehat B_i},\quad M_i=\int_M N_i,\quad
 \rho^0=(M,\rho_2),\quad
 \mathscr R=-\Delta_g\widehat B-K(N-\rho^0).
\]
Here and below integration is with respect to the specified surface area measure unless a local Euclidean density is used. The first constant is
\[
 c_1=2\lambda,
\]
and $c_2$ is independent of $\lambda$. At an ordinary center $\psi_i^A$ and $J_i$ are the background and its difference from its central value, as in the ordinary calculation.

Since the second residual density vanishes at that center, the smooth part of $u_1^0$ after removal of its self logarithm $-4\log|z|$ has Laplacian $2M-\rho_2$ there. The conformal area term contributes $-2K_g(p_i)$ to the Laplacian of the canonical background. Consequently
\[
 \Delta J_i(0)=\kappa_i(X)
\]
in the normalized conformal coordinates.

\smallskip\noindent\textit{Step 1. The contribution from the prescribed disks.} Choose a fixed $p>1$ close to one, with $ps_a<1$ at every mixed prescribed point and with the strict inequalities for pure singular bubbles in Lemma~\ref{univ:pure-residual}. On each prescribed disk, the unnormalized residual obeys
\[
 \|\mathscr R\|_p+\|\partial_\lambda\mathscr R\|_p+
 \|\partial_X\mathscr R\|_p\le Ct^2.                 \tag{*}
\]
For a mixed prescribed disk use the two products in the exact local residual after Lemma~\ref{univ:profile-estimates}. Both $J_j$ and $\partial_XJ_j$ are $O(r^2)$ by the matching of the constant and linear terms.

In the rational formulas the parameter derivatives of $\tau$, $F_0$, $D$ and $H_R$ are bounded by a constant times these quantities. This uses
\[
 \tau=d_0^se^{-T}
\]
with bounded logarithmic derivatives of $d_0$ and the fixed nonzero values of $P_{01},P_{02},W_0$ at the origin. Thus $\partial_Xn_1$ has the same majorant as $n_1$.

Write $n_2-f$ as the inner summand minus $f(2DH_R+H_R^2)/(D+H_R)^2$. Differentiating this difference preserves the majorant $I+r^{2A-2}\min\{1,(\ell/r)^{2B}\}$.

After multiplication by $r^2$, split at $R$ or $\ell$, respectively. The outer integrals are bounded by the $L^p$ norms of $t^2r^{-2s}$, $t^2r^{-2B}$ and $t^2r^{2A-2B}$. Their exponents are strictly integrable. The inner part of $r^2S(\ell/r)^{2B}$ is at most $Cr^{2A}$ on $r\le R$, and its norm there is
\[
 O(R^{2A+2/p})=O(t^2).
\]
The other inner terms are bounded by
\[
 C(R^{2/p}+\ell^{2/p})=O(t^2).
\]
The same rational estimates give the $\lambda$ derivative. At a pure disk the Poisson source of the differentiated background is still $O(r^{2(A+B)-2})+O(1)$, because the fixed-source scalar branch and its parameter derivative are bounded in $C^0$.

The logarithmic potential argument of Lemma~\ref{univ:pure-residual} therefore gives the same orders for the differentiated Taylor remainders, including the derivative of the linear term when $A\ge1/2$.

The differentiated angular correction has the same bounds as that correction. The radial integrals in that lemma prove $(*)$.

On fixed cut-off annuli all patch differences and these derivatives are $O(t^2)$ in smooth norms. Off the bubble disks the first unnormalized density and these derivatives are $O(t^2)$ in a fixed $L^p$, using the integrable prescribed weights there.

The second unnormalized density off the disks is the exact scalar residual density.

For the two masses one obtains, including one $\lambda$ or $X$ derivative of the discrepancies,
\[
 M_1-M=O(t^2L),\qquad M_2-\rho_2=O(t^2).
\]
On an ordinary disk use its smooth background, integrate the odd linear term of $m_{s_i}(e^{J_i}-1)$ in the translated disk, and use the quadratic moment $O(t^2L)$. Its second density has excess $O(t^2)$.

On mixed prescribed disks the unweighted exact profile masses and their derivatives follow by the fixed boundary fluxes of $\log F_R$ and $\log Q_R$ in Lemma~\ref{univ:local-scale-moments}; the flux errors are equally $O(t^2)$ under $X$ differentiation by the outer expansions.

The additional weighted density errors are controlled by $(*)$.

At a pure disk the first linear terms of the angular correction and the background are odd when they are present. The remaining error in the first density is bounded by $Cm_{A,R}r^\nu$ with
\[
 \nu>2A,
\]
integrable with the $O(t^2)$ bound. The excess in its second density is $O(t^2)$ in $L^1$ by the calculation for a pure singular bubble.

These statements hold after differentiation in coordinates translated with the ordinary centers and in fixed coordinates at the prescribed points: the discarded ordinary and pure linear moments vanish identically at each parameter.

Finally, the integral of the full residual density is $\gamma$, independent of $X$.

\smallskip\noindent\textit{Step 2. The correction on the ordinary disks.} On the ordinary disks the $L^{6/5}$ residual estimate is $C(t^{2/3}f+t^{5/3})$, and it also holds after one $\lambda$ derivative.

The normalization multiplies $N_i/M_i$ by $O(t^2L+m)$, with the same scale bound for the differentiated coefficient. The probability and its scale derivative have $H^{-1}$ norm $O(\sqrt L)$ by the interpolation in Proposition~\ref{univ:projected-construction}. Thus
\[
 \|\mathcal F_\mu(B)\|_{\Hh}
 +\|\partial_\lambda\mathcal F_\mu(B)\|_{\Hh}
 \le C(t^{2/3}f+t^{5/3}L+m\sqrt L).                 \tag{**}
\]
The $L^1$ comparisons to the model probabilities have a positive power on each disk, and all these probabilities have a fixed $L^{p_0}$ norm bounded by a power of $t^{-1}$ for some
\[
 p_0>1.
\]
Thus the comparison of bilinear forms and the signed model inverse bound the norm of the inverse of the projected linearization on $E_T$ by $CL$ throughout the fixed center neighborhood.

The contraction estimate of Proposition~\ref{univ:projected-construction} uses exponential smallness of the residual on this fixed neighborhood. This gives the unique projected family and the estimate
\[
 E\le C(t^{2/3}Lf+t^{5/3}L^2+mL^{3/2}).             \tag{***}
\]
The same estimate holds for a modulated true solution with the stated preliminary estimates, by applying this inverse on $E_T$ and absorbing $CL^{5/2}E^2$.

The preliminary estimates also put this solution in the enlarged uniqueness ball by taking its exponent smaller than both of the given decay exponents. On the refined center and mass domain $(***)$ gives $E=O(t^{7/6}L)$.

\smallskip\noindent\textit{Step 3. The reduced equations.} We test the residual against the parameter derivatives on $M$. An ordinary center tangent $Y=s_i\partial_{p_i}B$ is bounded. On its moving disk its components are $(4Z_j,-2Z_j)+O(t)$, where
\[
 Z_j(y)=y_j/(1+|y|^2).
\]
On the other disks it is $O(t)$. In particular this last estimate holds on prescribed disks: the source logarithms have fixed coefficients and centers, and the differentiated log ratios and exterior regular remainders are bounded there.

The tangent before subtraction of its mean has $L^1$ norm $O(t)$ as in the regular calculation. The first local pairing is $-4\pi s_iF_{i,j}$ because
\[
 \int_{\mathbb R^2}m_{s_i}Z_j z_k=\pi s_i\delta_{jk};
\]
the quadratic term of $e^{J_i}-1$ is even. Taylor remainders and the tail of this integral give $O(t^2 f+t^3 L)$. The contribution of the prescribed disks is $O(t^3)$ by $(*)$ and $Y=O(t)$ there. For each normalization term
\[
 \int N_kY_k=O(t):
\]
subtract the centered radial profile in the ordinary active pairing, use the bounded residual density in the ordinary disk, and use $Y=O(t)$ times the bounded masses in the other disks. The mean subtraction and the constant $\mu$ term have the same or smaller orders. It follows that
\[
 \langle\mathcal F_\mu(B),Y\rangle_{\Hh}
     =-4\pi s_iF_{i,j}+O(t^2f+t^3L+tm).
                                                               \tag{****}
\]
The derivative residual in either this scaled center direction or the $\lambda$ direction also satisfies
\[
 \|D\mathcal F_\mu(B)Y\|_{\Hh}
     \le C\sqrt L(tf+t^2L+m),                                    \tag{*****}
\]
with $Y=\partial_\lambda B$ for the scale direction. On the ordinary disks the derivative of $m_{s_i}$ in either scaled direction is bounded by $Cm_{s_i}$. At a fixed local coordinate $z$,
\[
 \partial_XJ_i(z,X)=\partial_X[\psi_i^A(z,X)-\psi_i^A(0,X)]=O(|z|),
\]
by the bounded mixed derivatives on the separated scalar branch. On the moving disk, differentiating at a fixed surface point adds a bounded coordinate velocity times
\[
 \nabla_zJ_i=F_i(X)+O(|z|).
\]
Multiplication by $s_i=O(t)$ gives $O(t(f+|z|))$; the other ordinary disks have only the $O(t|z|)$ term. For the scale derivative $J_i$ is fixed. Integrating the differentiated residual and using
\[
 \int m_{s_i}|z|=O(t),\qquad \int m_{s_i}|z|^2=O(t^2L)
\]
bounds its absolute mass by $C(tf+t^2L)$. The second excess and cut-off derivatives are $O(t^2)$ in absolute mass. Dividing the upper bound for the density by that bound leaves a polynomial pointwise bound in $t^{-1}$; interpolation to $L^{1+1/L}$ gives the factor $\sqrt L$ in the dual norm.

The prescribed disks and the complement are estimated separately by $(*)$, which also bounds their dual norms. The mass derivatives already established and
\[
 |s_i\partial_{p_i}\log N_k|+|\partial_\lambda\log N_k|\le C
\]
give the same $O(\sqrt L(t^2L+m))$ bound for the differentiated global normalization terms. This proves $(*****)$. For the scale, $\partial_\lambda\widehat B$ is bounded globally and has $L^1$ norm $O(t^2L)$.

On a mixed prescribed disk its two components, by the two logarithmic derivatives in Lemma~\ref{univ:local-scale-moments}, are $O(\min\{1,t^2r^{-2s}\}+\min\{1,t^2r^{-2B}\})$. The integral is $O(t^2)$ since
\[
 s,B<1.
\]
At a pure prescribed disk the same estimate with $A<1$ applies to the radial correction, and the angular correction and its derivative have $L^1$ norm $O(t^2)$.

Use the exact trial variation \eqref{reg:ordinary-trial-variation} with $M$ in place of $4\pi m$ and the two masses above. In its first pairing the ordinary scale tangent is $(-2(Z_0+1),Z_0+1)$, with
\[
 Z_0=(s_i^2-r^2)/(s_i^2+r^2).
\]
Its constant part exactly cancels the ordinary correction to $-2(M_1-M)$. The remaining prescribed and exterior mass corrections and residual pairings are $O(t^2)$ by the mass and residual estimates above. On the ordinary disks the remaining term is
\[
 2\int m_{s_i}Z_0(e^{J_i}-1).
\]
Using the quadratic moment in \eqref{reg:ordinary-trial-expansion} gives
\[
 \partial_\lambda J_\rho(B)
   =2\mu-4\pi\sum_i\kappa_i(X)s_i^2\log(1/s_i)
                 +O(t^2+t^2Lf^2+mt^2L+t^4L^2).
                                                               \tag{******}
\]
The two normalization terms in the exact variation have sizes $O(mt^2L+t^4L^2)$, so the trace uses the full $M$. For any bounded tangent vector $Y$ used as a test function, self-adjointness and the weighted estimate for the third derivative of the normalized logarithm give
\[
 \left|\langle\mathcal F_\mu(B+\eta)
                 -\mathcal F_\mu(B),Y\rangle_{\Hh}\right|
 \le \|D\mathcal F_\mu(B)Y\|_{\Hh}E+CLE^2.
\]
Indeed, with one bounded argument, the other two arguments in the weighted Hölder estimate each cost $C\sqrt L$. For a true solution, $(*****)$ and $(******)$ then give $m\le C(t^2L+LE^2+m^2)$ by Young's inequality, and $(****)$ gives
\[
 tf\le C(t^2f+t^3L+tm+t\sqrt L E+LE^2+m^2).
\]
Substitution of $(***)$ in the mass estimate, using the preliminary positive power to absorb $Cm^2L^4$, yields
\[
 m\le C(t^2L+t^{4/3}L^3f^2)
\]
and $E\le C(t^{2/3}Lf+t^{5/3}L^2)$. The center estimate divided by $t$ then has constant part $O(t^{5/3}L^3)$ and its remaining terms are bounded by
\[
 C(t+t^{2/3}L^{3/2})f
 +C(t^{1/3}+t^{4/3})L^3f^2
 +Ct^{5/3}L^6f^4 .
\]
All are absorbed using $f=o(1)$. Nondegeneracy gives
\[
 |X-X^0|\le Cf.
\]
For the sharp mass law use $(**)$ in the scale comparison: with the refined bounds its error is $o(t^2)$. In $(******)$ replace $s_i$ by $td_i(X)$ and use $L|X-X^0|=o(1)$ to obtain
\[
 \mu=L_0t^2L+O(t^2).
\]
This proves strict inclusion in the refined center and mass ranges for all sufficiently concentrated true solutions having the stated preliminary estimates.

\smallskip\noindent\textit{Step 4. The local degree.} Test the projected solution against the ordinary center tangents and
\[
 \partial_T B =-\tfrac12\partial_\lambda B.
\]
Their pairing matrix with the retained sources is uniformly invertible and has positive limiting determinant: the ordinary blocks are $(4\pi/3)I_2$ and the prescribed radial block is $4\pi s_0/3$, as in \eqref{hyb:chart} with no shape variables.

For the last block, $\partial_TB_1$ tends to $1+Z_{s_0,1}$ on the reference outer scale and the reference source is centered. The center tangent is $O(t)$ there, and the other cross pairings vanish in the limit by separation and the odd ordinary sources.

It follows that the test equations are equivalent to vanishing of every multiplier. Divide each two-dimensional center block by $-4\pi s_i$; this preserves orientation.

On the refined domain the estimates above give the center equations $F(X)+o(\sqrt t)$ and the scale equation
\[
               -\mu+L_0t^2L+o(t^2L).
\]
For the scale error use $(**)$ and $E=O(t^{7/6}L)$ so that both the linear correction and $LE^2$ are $o(t^2L)$. No derivative of $\eta$ is needed for these test equations. On $\mu L_0>0$, choose the small positive root $t_\mu$ with
\[
 |L_0|t_\mu^2\log(1/t_\mu)=|\mu|.
\]
Use a fixed small center ball in $(X-X^0)/\sqrt{t_\mu}$ and the interval $t_\mu/2<t<2t_\mu$. On its closure, with
\[
 t=t_\mu\tau,\qquad |\mu|/(t^2L)\to |L_0|/\tau^2
\]
uniformly. Use $C_1>4|L_0|$ and center radius less than $C_0/\sqrt2$ to remain strictly inside the refined domain. This product also lies in the enlarged domain since
\[
 |\mu|=O(t^2L)=o(t^{2\delta})
\]
for $\delta<1$. The center boundary is zero-free by the invertible $DF(X^0)$, and on the two scale ends in increasing $T$ the signs are $\operatorname{sgn}(L_0)$ and $-\operatorname{sgn}(L_0)$.

These facts give the product degree $\operatorname{sgn}\det DF(X^0)\,[-\operatorname{sgn}(L_0)]$ in the finite parameters, by straight homotopy to the linear center map and an affine scalar map with those endpoint signs.

Take a complementary ball of radius $C t^{7/6}L$, with $C$ larger than the correction bound in \textup{(***)}. Translate each fiber by the projected correction and, if necessary, shrink its radius around the graph by excision.

The equation on the complement has left-hand side $A v+O(L^{3/2}\|v\|^2)$, where the second derivative $A$ projected onto and restricted to the complement at the graph satisfies
\[
 \|A^{-1}\|\le C L.
\]
Along the homotopy that multiplies the remainder by a factor in $[0,1]$, every zero satisfies
\[
 \|v\|\le C L^{5/2}\|v\|^2,
\]
so all such zeros lie on the projected graph. The center and scale equations thus remain nonzero on the corresponding boundary faces when their tangential components are deformed to their values on the graph.

Trivialize the complements over the product of the center ball and scale interval by projection transport as in Theorem~\ref{hyb:transfer}.

For small $|\mu|$ these coordinates are unique. Indeed, if two representations agree, their approximations differ by $o(1)$ in $H^1$.

Their ordinary widths are comparable on the chosen scale interval. Rescale the first-component gradients on one ordinary disk by the first representation's width.

If the other center displacement divided by this width were unbounded, the other gradient would tend to zero on bounded rescaled sets, contrary to the nonconstant ordinary profile. Otherwise pass to a subsequential limit of the relative displacement and width ratio.

The two limiting gradients are identical. The unique critical point of $-2\log(1+|y|^2)$ fixes the displacement and its second derivative fixes the ratio, giving zero displacement and ratio one. Apply this at every ordinary disk. Their matching relations give a difference $o(1)$ in $T$.

In the common scaled coordinates the derivative of the orthogonality equations is then the negative of the invertible limiting matrix plus $o(1)$, uniformly on the segment between the two parameters.

Indeed, the derivatives of the source functionals are bounded and act on an $o(1)$ remainder.

Integration gives equality of the parameters. One fixed labeling near $X^0$ suffices; a permutation of equal types has positive determinant in the two-dimensional blocks.

For fixed nonzero $\mu$ the closed parameter product is compact and the fiber radii are positive and continuous. The derivative in the base variables of $L_\theta(u-B_\theta)$ is invertible on this neighborhood by the same matrix estimate along short segments.

The implicit function theorem thus recovers the parameters locally. The image of the open product is thus open in the original space.

For a convergent sequence in its closure, pass to a subsequence of the closed base parameters and subtract the continuous graph field.

The translated fibers then converge after the projection transport. A boundary zero would consequently lie over a base face or on a fiber sphere, both excluded above.

The matrix of inner products of the source vectors is positive definite. The matrix pairing these vectors with the parameter derivatives therefore preserves orientation in the same identity-plus-compact coordinates as in Theorem~\ref{hyb:transfer}.

The degree on this neighborhood in the original space equals the finite-dimensional degree times $(-1)^\nu$. All active outer and inner profile exponents are at most one, with equality only at ordinary centers, so all angular terms in \eqref{hyb:index} vanish. With
\[
 n=U+|P|+|Q|
\]
and $v=|Q|$, the exact index is
\[
 \nu=m_S+U+|P|+2|Q|-1.
\]
The allowed side therefore has degree $(-1)^{U+|P|}\operatorname{sgn}(L_0)d_S(Z)$. Under the preliminary estimate for the blowing-up solutions, their mass expansion excludes both the critical parameter and the opposite side.

Thus for $L_0>0$ this value is on the positive side, while for $L_0<0$ it is on the negative side and is subtracted. In both cases the difference of the degrees on the two sides is $(-1)^{U+|P|}d_S(Z)$.
\end{proof}
\subsection{Quantitative description of blowing-up solutions}
\begin{lemma}[A weighted scalar estimate]
\label{univ:weighted-scalar}
Fix $q>0$ with $q\notin\mathbb N$, $0<\kappa\leq\min\{1,q\}$, and $0<\zeta<2q$. Set
\[
 \omega_{\kappa,\zeta}(y)
   =|y|^{2\kappa-2}(1+|y|)^{-2\kappa-\zeta},
 \qquad \ell(y)=\log(2+|y|).
\]
\begin{itemize}
\item[\textup{(i)}]
For sufficiently large $L$, let $v$ be a weak solution on $D_L$, bounded at zero with $v(0)=0$ for its continuous representative. Then
\[
 \sup_{D_L}\frac{|v|}{\ell}
 \leq C\left\{
 \sup_{D_L\setminus\{0\}}
 \frac{|(-\Delta-2m_{q,1})v|}{\omega_{\kappa,\zeta}}
 +\operatorname{osc}_{\partial D_L}v\right\}.
\]
The equation holds in distributions on $D_L$, with no additional point mass at zero.
\item[\textup{(ii)}]
If the right-hand side of the estimate is finite, the powers of $m_{q,1}$ and $\omega_{\kappa,\zeta}$ at zero give a local $L^p$ right-hand side for the equation for some $p>1$. Thus $v$ is continuous there.
\item[\textup{(iii)}]
The constant is uniform for $q$ in compact subsets of
\[
 (0,\infty)\setminus\mathbb N
\]
and for the other exponents in compact subsets of the stated strict inequalities.
\end{itemize}
\end{lemma}
\begin{proof}
\smallskip\noindent\textit{Step 1. Normalization at the origin.}\quad The radial kernel function obtained by dilation is
\[
 Z_q=(1-r^{2q})/(1+r^{2q});
\]
direct differentiation gives its equation and $Z_q(0)=1$. Near zero its independent radial solution is logarithmically singular, by the identity
\[
 r(v_1v_2^\prime-v_1^\prime v_2)=\mathrm{const}
\]
for two radial solutions. For angular frequency $j\geq1$, write the radial coefficient as $r^jF(t)$ with $t=r^{2q}$. Its equation is
\[
 qtF''+(q+j)F'+\frac{2q}{(1+t)^2}F=0.
\]
The function $F(t)=(j-q)+2q/(1+t)$ solves this equation. The independent radial coefficient grows like $r^{-j}$ at zero by the same determinant identity, so the bounded-at-zero solution is proportional to
\[
 r^j\frac{(j+q)+(j-q)r^{2q}}{1+r^{2q}}.
\]
Its leading coefficient at infinity is $j-q\ne0$. Projecting any homogeneous field of logarithmic growth onto each angular frequency therefore eliminates all $j\ge1$. The condition $v(0)=0$ eliminates the remaining dilation.

\smallskip\noindent\textit{Step 2. The Green representation.}\quad The disk Green-kernel difference used in Lemma~\ref{hyb:scalar-inverse} is
\[
 \mathcal K_L(y,z)=\frac1{2\pi}
 \left\{\log\frac{|z|}{|y-z|}
          +\log\left|1-\frac{y\overline z}{L^2}\right|\right\}.
\]
Its integral against $\omega_{\kappa,\zeta}$ is bounded by $C\ell(y)$. To handle a fixed disk about zero, take $p_0>1$ with
\[
 p_0(\kappa-1)>-1
\]
and apply H\"older to the translated logarithm. Its $L^{p_0'}$ norm on the fixed disk is bounded uniformly when $y$ stays in a fixed larger disk. The unshifted logarithm is handled in the same way, and for $|y|$ large the direct logarithm there is bounded by $C\ell(y)+|\log|z||$.

At infinity the weight is comparable to $r^{-2-\zeta}$. The corresponding tail integral against $m_{q,1}\ell$ tends to zero after division by $\ell(y)$, uniformly in $L$ and $y$, since the profile decays as $r^{-2-2q}$. More explicitly, with
\[
 \mathcal K_L(y,z)=G_L(y,z)-G_L(0,z),
\]
for $M\geq2$ one has
\[
 \sup_{L\geq2M}\sup_{y\in D_L}
 \frac{1}{\ell(y)}\int_{M<|z|<L}
       |\mathcal K_L(y,z)|m_{q,1}(z)\ell(z)\,dz
       \leq CM^{-2q}(1+\log M)^2.
\]
For the direct logarithm, split into $|z|<|y|/2$,
\[
 |y|/2\leq|z|\leq2|y|,\qquad |z|>2|y|.
\]
On the first set divide the bound $C+\log(|y|/|z|)$ by $\ell(y)$; on the last use $C|y|/|z|$. On the middle annulus rescale by $|y|$; the logarithmic singularity has a uniformly bounded integral and
\[
 m_{q,1}(z)\ell(z)\leq C|y|^{-2-2q}\log(2+|y|).
\]
For the reflected logarithm, if $|y|\leq L/2$ or $|z|\leq L/2$, its argument stays a fixed distance from zero. On the remaining set
\[
 |y|,|z|>L/2,
\]
rescaling by $L$ gives a uniformly integrable logarithmic singularity, and division by $\ell(y)$ cancels the factor $\log L$ in the density bound. The contribution is $O(L^{-2q})$. These estimates also prove the bound of the full kernel difference against $\omega_{\kappa,\zeta}$ stated above.

If the asserted estimate for all sufficiently large $L$ fails, take $L_n\to\infty$ and normalize counterexamples to logarithmic norm one, with weighted residuals and boundary oscillations tending to zero. On each fixed disk choose $p>1$ with
\[
 p(\kappa-1)>-1,\qquad p(q-1)>-1
\]
whenever needed. The bound $|v_n|\le\ell$ and the forcing bound then give local $W^{2,p}$ compactness into $C^0$, including at zero. The description of the kernel space and the condition $v_n(0)=0$ imply that the limit vanishes.

In the representation by the difference of Green kernels, the difference of the harmonic terms at $y$ and $0$ is bounded by the boundary oscillation. The residual term divided by $\ell(y)$ tends to zero by the full kernel bound.

On a fixed disk, the $2m_{q,1}v_n$ term tends to zero by local uniform convergence and the same logarithmic integrability. Outside that disk the uniform kernel tail estimate makes the remaining quotient small. This contradicts
\[
 \sup_{D_{L_n}}|v_n|/\ell=1.
\]
On compact parameter families the same $p$ can be chosen and a subsequence of the profile exponents converges to a positive noninteger limit; the preceding bounds and the condition $v(0)=0$ still hold.
\end{proof}
\begin{lemma}[Fractional moments and contour fields]
\label{univ:fractional-contours}
Let $A,t>0$ and set
\[
 s=A+t,\qquad r=R^{s/t+o(1)},\qquad
 1<\theta<s/t,\qquad \tau=R^\theta,\qquad \sigma=R^{1/2}.
\]
Assume the density bounds at a mixed point in the proof of Lemma~\ref{hyb:singular-rate}, with $t$ in place of $B$ and any sufficiently small fixed $\varepsilon>0$. After enlarging the annuli by a fixed factor if necessary, assume also
\[
 |\nabla f_i|\leq Cd^{-1}
       \sup_{d/4<|z|<4d}f_i
       \quad(d/2<|z|<2d).
\]
For any
\[
 0<\eta_t<\min\{1,2t\},\qquad
 0<\eta_s<\min\{1,2s\},\qquad
 0<\eta_A<\min\{1,2A\},
\]
choose $\varepsilon>0$ sufficiently small. At $|z|=d$ use the term
\[
 \frac{I_i(d)}{\pi z},\qquad I_i(d)=\int_{D_{d/2}}f_i\,dA.
\]
Keep $d$ and $I_i(d)$ fixed in spatial derivatives. For the remaining field $\mathcal C$, the bounds on $d|\mathcal C|+d^2|\nabla\mathcal C|$ are as follows.
\begin{itemize}
\item[\textup{(i)}]
At $d=\tau$, subtract the term for $f_2$. The bound is
\[
 C\left\{(r/\tau)^{\eta_t}
              +(\tau/R)^{\eta_s}+\tau^{\eta_A}\right\}.
\]
\item[\textup{(ii)}]
At $d=\sigma$, subtract the terms for both components. The bound is
\[
 C\left\{(R/\sigma)^{\eta_s}
              +(r/\sigma)^{\eta_t}+\sigma^{\eta_A}\right\}.
\]
\end{itemize}
Each term is a positive power of $R$, possibly with a smaller exponent. Substitution of a limiting mass for $I_i(d)$ requires an estimate of their difference.
\end{lemma}
\begin{proof}
\smallskip\noindent\textit{Step 1. The mass on each annulus.}\quad The radial scalar density satisfies
\[
 \int_{\mathbb R^2}|y|^\eta m_{q,1}(y)\,dy<\infty
                  \quad(0<\eta<2q).
\]
For $|y|<d/2$ and $|z|=d$, the differences between the Cauchy kernel and its leading term, and between their first derivatives, are bounded by $Cd^{-1}(|y|/d)^\eta$ and $Cd^{-2}(|y|/d)^\eta$, respectively, for $0<\eta\leq1$.

Integration gives the fractional-moment terms. The same estimates follow directly from the stated pointwise bounds, with $2q$ replaced by $2q-\varepsilon$.

For the first density between $2\tau$ and $R$, direct integration gives the relative field bound
\[
 C\tau R^{-2s+\varepsilon}
          \int_{2\tau}^R d^{2s-2-\varepsilon}\,dd.
\]
It is bounded by $C(\tau/R)^{\eta_s}$ for every
\[
 \eta_s<\min\{1,2s-\varepsilon\};
\]
the equality case in the integral produces only a logarithm. The differentiated integral has one further inverse power of $d$, but is multiplied by $\tau^2$, and obeys the same weaker bound.

Inside $\tau/2$, the first enclosed mass is $O((\tau/R)^{2s-\varepsilon})$: integrate the intermediate bound, and on $|y|<r$ the inner bound contributes $O((r/R)^{2s-\varepsilon})$. On
\[
 \tau/2<|y|<2\tau,
\]
the density bound is $O(\tau^{-2}(\tau/R)^{2s-\varepsilon})$, so the annular kernel estimates below give the same relative bound for the field and derivative. For $|y|>R$, the relative field bound is
\[
 C\tau R^{2s-\varepsilon}\int_R^{d_0}u^{-2s-2+\varepsilon}\,du=O(\tau/R),
\]
and the derivative has the analogous bound
\[
 O((\tau/R)^2).
\]
Changing the dividing radii by fixed factors only alter the constants. The second inner tail between $\tau/2$ and $R/C_\varepsilon$ costs $O((r/\tau)^{2t-\varepsilon})$ in both relative norms by its density bound and the adjoining-annulus estimate below.

For the residual density from $R/C_\varepsilon$ outward, including the transition annulus, the field and derivative are bounded by
\[
 C\tau\int_{R/C_\varepsilon}^{d_0}d^{2A-2-\varepsilon}\,dd,
 \qquad
 C\tau^2\int_{R/C_\varepsilon}^{d_0}d^{2A-3-\varepsilon}\,dd.
\]
If the first integral diverges as $R\downarrow0$, then
\[
 \tau R^{2A-1-\varepsilon}\leq\tau^{2A-\varepsilon}
\]
because $\tau<R$; the analogous estimate for the second integral is $\tau^2R^{2A-2-\varepsilon}\leq\tau^{2A-\varepsilon}$. Convergent integrals are bounded by $C\tau$ or $C\tau^2$. At an equality of exponents, the logarithm is absorbed by the strict choice $\eta_A<\min\{1,2A-\varepsilon\}$.

\smallskip\noindent\textit{Step 2. The local Cauchy integral.}\quad On the observation disk $D_{d/8}(z)$ the Cauchy kernel itself has integral $O(d)$. For its holomorphic derivative, subtract $f_i(z)$ there; the constant has zero principal value by angular integration, and the remainder has integral bounded by
\[
 C\sup_{D_{d/8}(z)}|\nabla f_i|
        \int_0^{d/8}\frac{\varrho}{\varrho^2}\varrho\,d\varrho
 \le C\sup_{d/4<|y|<4d} f_i .
\]
The rest of the adjoining annulus has kernel and derivative bounded by $C/d$ and $C/d^2$, so integration gives the same relative bounds. Only this derivative estimate needs the annular gradient hypothesis. The remaining complex derivative is
\[
 \partial_{\bar z}\mathcal C[f_i]=f_i,
\]
with the stated normalization of the Cauchy kernel, and obeys the same estimate. At $\sigma$, the part of the first density inside $\sigma/2$ has fractional moment bounded by $CR^{\eta_s}$ via the inner, intermediate and outer majorants. On
\[
 |y|\ge\sigma/2
\]
the first outer tail gives $O((R/\sigma)^{2s-\varepsilon})$ in the relative norms, using the annular estimate at the observation point. For the second density on
\[
 |y|<R/C_\varepsilon,
\]
its inner-scale majorants give the moment $Cr^{\eta_t}$, and the difference kernel is used since this region is inside $\sigma/2$. On $R/C_\varepsilon\leq|y|<\sigma/2$ the diffuse bound from the transition estimate gives at most $C\sigma^{-\eta_A}\int_0^{\sigma/2}u^{2A+\eta_A-1-\varepsilon}\,du$ in the relative norms. On
\[
 |y|\ge\sigma/2,
\]
the second diffuse bound gives $O(\sigma^{\eta_A})$ by the annular estimate and the outer integrals in the preceding paragraph with lower radius comparable to $\sigma$. The subtracted mass at $\sigma$ includes precisely all these interior pieces. Finally
\[
 r/\tau=R^{s/t-\theta+o(1)},\qquad \tau/R=R^{\theta-1}
\]
and $R/\sigma=R^{1/2}$, so all the powers are positive.
\end{proof}
\begin{lemma}[Interpolation for singular error densities]
\label{univ:singular-interpolation}
Suppose the signed densities $e_R$ on $M$ satisfy, for some $p_0>1$, $c>0$, $K\ge0$ and $N<\infty$,
\[
 \|e_R\|_1\leq CR^c(1+|\log R|)^N,
 \qquad \|e_R\|_{p_0}\leq CR^{-K}.
\]
Then for some $1<p<p_0$ and $c'>0$,
\[
 \|e_R\|_p+\|Ge_R\|_{C^0\cap H^1}\leq CR^{c'}.
\]
For finitely many independent radii $R_a$, one may sum these estimates with $R=\max_a R_a$, taking common exponents for the finite family.
\end{lemma}
\begin{proof}
Choose $\vartheta\in(0,1)$ such that $(1-\vartheta)c-\vartheta K>0$. With $1/p=1-\vartheta+\vartheta/p_0$, interpolation gives
\[
 \|e_R\|_p\leq
 \|e_R\|_1^{1-\vartheta}\|e_R\|_{p_0}^{\vartheta}
 \leq CR^{(1-\vartheta)c-\vartheta K}
                  (1+|\log R|)^{N(1-\vartheta)}.
\]
Decrease the exponent to absorb the logarithm and use $G:L^p\to W^{2,p}\hookrightarrow C^0\cap H^1$. For a finite list, first choose one $p$ satisfying all the strict inequalities and then sum the individual estimates.
\end{proof}
These estimates apply to the singular-rate argument in Lemma~\ref{hyb:singular-rate}, with $t$ denoting its inner profile exponent $B$. Only the undifferentiated contour field is needed for its circular oscillation estimates.

After subtraction of the actual enclosed charges the fields of the subtracted point masses are radial. This argument therefore uses the enclosed masses directly. The derivative estimate of Lemma~\ref{univ:fractional-contours} uses the additional annular gradient hypothesis.

In the comparison on the inner disk the radius $L$ is $\tau/r$ in the coordinates $z=ry$. Choose $\kappa<\min(1,A,t)$ and $\zeta>0$ small. After division by $\omega_{\kappa,\zeta}$, the scaled first density and the smooth term are bounded by
\[
 C\{(r/R)^{2s-\varepsilon}L^{2s-\varepsilon+\zeta}
                  +r^2L^{2+\zeta}\}.
\]
At zero the first quotient contains the nonnegative power $|y|^{2(A-\kappa)}$. The limiting powers of $R$ in the display are $(2s-\varepsilon)(\theta-1)-\zeta(s/t-\theta)$ and $2\theta-\zeta(s/t-\theta)$, both positive for sufficiently small $\zeta$.

On the outer disk subtract the Dirichlet correction of Lemma~\ref{hyb:singular-rate} before fixing the fixed value at zero. Its $L^\infty$ norm has a positive power by the inner fit, the two tail bounds and fixed-exponent interpolation. With
\[
 a_0=r/R,
\]
the additional outer forcing has the factor $a_0^{2t}|y|^{2A-2}e^v$. Dividing by the weight is harmless at zero for $\kappa\le A$, and at infinity $e^v\le C|y|^{-4s+\varepsilon}$ gives a factor bounded by $C|y|^{2A-4s+\varepsilon+\zeta}$. Thus its weighted norm is $O(a_0^{2t})$ when $\zeta<2s-\varepsilon$.

For the pure outer fit, the source terms in the disk of scaled radius $R^{-1/2}$ have weighted bounds $CR^{A+B-(\varepsilon+\zeta)/2}$ and $CR^{1-\zeta/2}$. Indeed at $|y|\le C_\varepsilon$ the second density is bounded by $CR^{2(A+B)-\varepsilon}|y|^{2B-2}$ after multiplication by $R^2$, so take
\[
 \kappa<B.
\]
Outside this fixed disk use the exterior power $|y|^{2(A+B)-2-\varepsilon}$. Choose $\zeta$ and $\varepsilon$ small enough that both powers are positive. In all three applications the preliminary logarithmic norm is $o(1)$ by local convergence, mass exhaustion and the circle bounds.

If that norm is $\delta$, then the nonlinear remainder divided by the weight is $O(\delta^2)$: use $\kappa\le q$ at zero and $m_{q,1}\ell^2e^{\delta\ell}/\omega_{\kappa,\zeta}$ bounded at infinity for
\[
 \delta<2q-\zeta.
\]
Lemma~\ref{univ:weighted-scalar} therefore gives the quantitative bounds with a positive power.

For interpolation choose a fixed $p_0>1$ sufficiently close to one that
\[
 p_0(A-1)>-1,\qquad p_0(t-1)>-1,
\]
and $p_0(s-1)>-1$ whenever the respective power occurs. Take the exponent loss small enough that the resulting bounds still satisfy the strict inequalities on the exponents. The relation
\[
 r=R^{s/t+o(1)}
\]
makes all width factors polynomial in the inverse outer radius. This gives the polynomial $L^{p_0}$ density bounds used in Lemma~\ref{univ:singular-interpolation}.

For the approximations used here, each ordinary center is scaled by its own width and each free mixed center by its inner width; the other coordinates are bounded shapes and the common logarithmic scale. We estimate the centered source functionals with their cut-offs.

In the component where it is supported, and up to multiplication by the fixed Cartan matrix, a local source has the form
\[
 q_{a,\theta}=\chi_a m_{a,\theta}
       (Z_{a,\theta}-\overline Z_{a,\theta}),\qquad
 \overline Z_{a,\theta}
 =\frac{\int_M\chi_am_{a,\theta}Z_{a,\theta}}
        {\int_M\chi_am_{a,\theta}}.
\]
Thus $\int_M q_{a,\theta}=0$. The denominator is bounded above and away from zero, and $Z_{a,\theta}$ is uniformly bounded on the separated, bounded-shape parameter set. Thus $\|q_{a,\theta}\|_{\rm TV}\leq C$. With
\[
 \omega_{a,\theta}=\chi_am_{a,\theta}/\int\chi_am_{a,\theta},
\]
the weighted Poincar\'e inequalities in Sections~\ref{reg:section} and~\ref{hyb:section} give
\[
 \left|\int_M q_{a,\theta}v\right|
 \leq\left(\int_M\chi_am_{a,\theta}\right)
       \sqrt{\operatorname{Var}_{\omega_{a,\theta}}Z_{a,\theta}}
       \sqrt{\operatorname{Var}_{\omega_{a,\theta}}v}
 \leq C\|\nabla v\|_2.
\]
The source is centered separately on each disk. Equivalently, its Green energy is bounded because the constant logarithmic dilation term pairs with two zero integrals.

For the normalized sources of the approximation and the constants subtracted after cut-off, we use \eqref{reg:modulation-bounds} and \eqref{hyb:modulation-bounds}.

In the comparison with the approximation, a positive-power \(L^p\) source error has a Green potential bounded by a positive power of the concentration radius.

Its centering coefficient is of the same order, and the Green representative of the normalized local mean costs only \(O(\sqrt T)\). This correction is still uniformly bounded.

For a first scaled derivative, write \(a_\theta=\int_M\chi_a m_{a,\theta}\). The quotient rule gives
\[
 D\overline Z_{a,\theta}=a_\theta^{-1}\left[
   \int_M D(\chi_a m_{a,\theta})(Z_{a,\theta}-\overline Z_{a,\theta})
       +\int_M\chi_a m_{a,\theta}DZ_{a,\theta}\right],
 \qquad \int_M Dq_{a,\theta}=0 .
\]
At fixed prescribed centers, differentiating the scale preserves the local upper bound for the density, also for a noninteger positive shift \(q<1\): for the radial model
\[
 |\partial_{\log R}m_{q,R}|\leq C_qm_{q,R}
\]
and
\[
 |\partial_{\log R}Z_{q,R}|\leq C_q.
\]
At a free mixed point the outer derivative in the inner translation coordinate carries the inner-to-outer width ratio. The differentiated quadratic-profile density is dominated on the outer scale by \(C R^{-2}(1+|z/R|^2)^{-2}\), including at its density zero.

On each fixed annulus the cut-off terms decay as a positive power of the radius of concentration. Thus one can estimate the derivative against its own upper bound after subtracting the corresponding local mean of \(v\).

The same weighted variance argument and the \(L^1\) bound give both dual norms. Consequently, for the approximation sources,
\begin{equation}\label{journal:source-functional-bounds}
 \begin{aligned}
 |L_{\theta,\alpha}v|+|(D_\theta L_{\theta,\alpha})v|
     &\leq C\|v\|_{C^0},\\
 |L_{\theta,\alpha}v|+|(D_\theta L_{\theta,\alpha})v|
     &\leq C\|v\|_{H^1}.
 \end{aligned}
\end{equation}
The common radial source is a bounded linear combination of the finitely many local sources. 
Set $J_\theta=L_\theta D_\theta B_\theta$. The bounds on the scaled
tangents and the inverse pairing matrices, together with
\eqref{journal:source-functional-bounds}, give
\[
 \|L_\theta\|+\|D_\theta L_\theta\|
 +\|D_\theta B_\theta\|+\|J_\theta^{-1}\|\le C.
\]
In the rescaled variable $y$ of each source, write $\widetilde q_{a,\theta}$
for its signed measure and $A_R=\{R^{-1}\le |y|\le R\}$ for $R>1$.
The profile bounds give
\[
 \lim_{R\to\infty}\sup_{a,\theta}
       \int_{A_R^c}d|\widetilde q_{a,\theta}|=0.
\]
The scaled tangents are uniformly bounded. For fixed $R$, the rescaled
profiles and tangents converge on $A_R$ uniformly on compact scaled
parameter balls. A tangent from a nested scale has a constant limit $c$,
possibly zero. Its constant part cancels in the full pairing since
\[
 \int_M c\,q_{a,\theta}\,dV_g=0.
\]
In coordinates based on the reference widths, the entries of $J_\theta$
are thus uniformly approximated by a compact family of continuous
pairings. On bounded scale intervals continuity on compact sets gives
the same conclusion. Hence, uniformly in the base scale,
\[
 \lim_{h\downarrow0}\ \sup_{|\theta-\theta'|\le h}
            \|J_\theta-J_{\theta'}\|=0
\]
in each scaled chart. By Lemma~\ref{journal:modulation-uniqueness}\textup{(i)},
the preliminary bound at $\theta_0$ gives a unique local representation
\[
 u=B_\theta+\eta,\qquad L_\theta\eta=0,\qquad
 |\theta-\theta_0|+\|\eta\|_X\le C\|u-B_{\theta_0}\|_X .
\]
Use tildes for a second representation of the same field. The closed scale
interval bounds the ratios of corresponding widths above and away from zero, and
\[
 \|B_\theta-B_{\widetilde\theta}\|_{H^1}
 \le\|u-B_\theta\|_{H^1}+\|u-B_{\widetilde\theta}\|_{H^1}=o(1).
\]
As in Subsection~\ref{reg:isolation}, let $\sigma_a$ be the ordinary
width or the mixed inner width at a free center $a$. The unique zero of
the limiting radial gradient and its derivative there give
\[
 \frac{d_g(a,\widetilde a)}{\sigma_a}\to0,\qquad
 \frac{\widetilde\sigma_a}{\sigma_a}\to1 .
\]
At a mixed center, pass to a subsequence with limiting shapes
$B,\widetilde B$ and outer width ratio $\widetilde R/R\to\lambda>0$.
Equality of the outer gradients gives proportional quartics, both with
leading coefficient one. Comparing the quadratic and constant terms gives
\[
 1+|y^2-B|^2=\lambda^4+|y^2-\lambda^2\widetilde B|^2,\qquad
 B=\lambda^2\widetilde B,\qquad \lambda^4=1 .
\]
At a prescribed reference point, if present, equality of the limiting
radial gradients on a fixed annulus gives $\widetilde R_0/R_0\to1$,
as in Step~4 of the proof of Theorem~\ref{hyb:transfer}. The limiting
matching coefficients agree, so \eqref{reg:scales} or \eqref{hyb:matching} gives
\[
 \widetilde T-T=o(1),\qquad |\widetilde\theta-\theta|=o(1)
 \quad\text{in the scaled coordinates}.
\]
Both parameters therefore lie in one uniqueness ball of
Lemma~\ref{journal:modulation-uniqueness}\textup{(i)}, so they coincide.
Only the profile tangents and the local source estimates enter this comparison.
\begin{proposition}[Quantitative description for positive shifted parameters]
\label{univ:complete-capture}
Fix $A_p,B_p>0$ satisfying \eqref{fnd:independence}. Use the weights and disjoint neighborhoods of shadow solutions chosen above. Assume the following.
\begin{itemize}
\item[\textup{(a)}]
If the leading scale is ordinary, assume $L_0\ne0$ as in Lemma~\ref{univ:regular-only} when all bubbles are regular, and as in Lemma~\ref{univ:ordinary-dominated} when a prescribed point also blows up.
\item[\textup{(b)}]
Fix $\rho_2\notin\Gamma_2(\alpha)$ satisfying the scalar exclusions. Let an unbounded sequence of Toda solutions approach a critical value of the first parameter.
\end{itemize}
Pass to a subsequence whose norms tend to infinity. Then the following hold.
\begin{itemize}
\item[\textup{(i)}]
The errors in the local profiles, matched masses and Green potentials are bounded by positive powers of their outer radii times fixed logarithmic factors. The common matching equations have the same error bounds.
\item[\textup{(ii)}]
If a prescribed point blows up, the solutions eventually belong to the family of Proposition~\ref{univ:projected-construction} or Lemma~\ref{univ:ordinary-dominated}, according to the configuration.
\item[\textup{(iii)}]
If all blow-up points are regular, the solutions eventually belong to the corresponding family of Lemma~\ref{univ:regular-only}.
\end{itemize}
\end{proposition}
\begin{proof}
\smallskip\noindent\textit{Step 1. The local masses.} We apply the arguments of Proposition~\ref{hyb:local-capture} and Theorem~\ref{hyb:capture} for $A_p,B_p>0$. If the second component concentrated completely, its generic local mass list would give
\[
 \rho_2\in\Gamma_2(\alpha).
\]
Let $\mathcal B$ be the blow-up set and let $4\pi(x_a,y_a)$ be the
local mass pair at $a\in\mathcal B$. Since $\rho_2\notin\Gamma_2(\alpha)$,
the second component has a residual density $g$ of mass
\[
 \gamma=\int_M g\,dV_g=\rho_2-4\pi\sum_{a\in\mathcal B}y_a>0.
\]
By Step~3 of the proof of Proposition~\ref{fnd:signed-compactness},
the circular averages and the local Pohozaev identity force the first
component to concentrate completely. Thus the normalized densities $f_i^k$ satisfy
\[
 f_1^k\,dV_g\rightharpoonup4\pi\sum_{a\in\mathcal B}x_a\delta_a,
 \qquad
 f_2^k\,dV_g\rightharpoonup4\pi\sum_{a\in\mathcal B}y_a\delta_a+g\,dV_g .
\]
Fix $a\in\mathcal B$ and write $(x,y)=(x_a,y_a)$. Let $(A,B)$ be the
local shifts, with $A=B=1$ at a regular point. In a coordinate centered
at $a$, the effective second weight $H$ and the residual potential $w$ satisfy
\[
 H(z)\asymp |z|^{2(B+x-2y)-2},\qquad w=2Gg\ge-C .
\]
Since the second normalization has a finite positive limit and $g\in L^1(M)$,
\[
 g(z)\ge c|z|^{2(B+x-2y)-2},\qquad
 \int_0^\varepsilon r^{2(B+x-2y)-1}\,dr<\infty
 \quad\Longrightarrow\quad 2y-x<B .
\]
Lemma~\ref{fnd:local-unitary}\textup{(ii),(iii)} therefore leaves
\[
 (x,y)\in
 \begin{cases}
 \{(A,0),(A+B,B)\},&a\in S,\\
 \{(1,0),(2,1)\},&a\notin S .
 \end{cases}
\]
Lemma~\ref{reg:shadow-limit} identifies the limiting scalar field and
centers. For the weights of Proposition~\ref{univ:shadow}\textup{(i)},
they lie in one of the disjoint neighborhoods of nondegenerate shadow
solutions, and at each regular mixed point
\[
 \kappa_j>4|T_j|.
\]
If $\mathcal B\cap S=\varnothing$, Lemma~\ref{univ:regular-only} applies.
Otherwise fix $a\in\mathcal B\cap S$. On an isolating disk $D$, the
canonical fields $U_i$ and smooth forcings $\kappa_i$ of
Lemma~\ref{hyb:centered-singular-profiles} satisfy
\[
 f_i=e^{U_i},\qquad
 \|\kappa_i\|_{C^1(D)}+\int_D f_i\,dA
       +\operatorname{osc}_{\partial D}U_i\le C,\qquad i=1,2,
\]
by the local reduction in the proof of Proposition~\ref{hyb:local-capture}.
By \eqref{fnd:independence}, Lemma~\ref{fnd:local-unitary}\textup{(ii)}
gives $|z|^2f_i(z)\le C$ on a smaller disk.
Lemma~\ref{hyb:centered-singular-profiles} then gives the profiles
at outer radius $R$ and mixed inner radius $r$,
\[
 \begin{aligned}
 (x,y)=(A,0)   &:\quad m_{A,R} \quad\text{in the first component},\\
 (x,y)=(A+B,B)&:\quad m_{A+B,R},\quad m_{B,r},\qquad r/R\to0 .
 \end{aligned}
\]
On a fixed small circle about $a$, the finite positive normalization
and convergence away from the blow-up centers give $\sup_{|z|=d}|U_2|\le C_d$.
Set $s=A$ in the pure case and $s=A+B$ in the mixed case.
Lemma~\ref{hyb:singular-rate} gives, for some $c>0$ and $N<\infty$,
\[
 \begin{gathered}
 \|f_1-m_{s,R}\|_{L^1(D_{R^{1/2}})}
       +\|f_1\|_{L^1(D\setminus D_{R^{1/2}})}
   \le CR^c(1+|\log R|)^N,\\
 \|f_2-m_{B,r}\|_{L^1(D_{R^\theta})}
   \le CR^c(1+|\log R|)^N
       \quad\text{in the mixed case},\qquad 1<\theta<s/B .
 \end{gathered}
\]
Here $D$ has been reduced if necessary. Each comparison is made at its
own fitted radius before summing Green potentials. On source-free disks
Lemma~\ref{reg:local-power} gives the regular estimates. For the shape
$B_j$ at a mixed regular point, Lemma~\ref{reg:stress}\textup{(ii)} gives
\[
 \kappa_j>4|T_j|\quad\Longrightarrow\quad |B_j|\le C .
\]
Step~2 estimates the second density on the transition annuli in
\eqref{univ:capture-actual-tail}--\eqref{univ:capture-fitted-tail}.
The exponents depend on the fixed shifts. Choose $p_0>1$ close enough to
one that every local weight power $|z|^{2(q-1)}$ satisfies
\[
 p_0(q-1)>-1 .
\]
The pointwise bounds then give, for each actual or fitted density $F$
on its local comparison disk of outer scale $R$,
\[
 \|F\|_{L^{p_0}}\le CR^{-K_0}\qquad\text{for some fixed }K_0>0 .
\]
Lemma~\ref{univ:singular-interpolation} then controls the corresponding $L^p$ and $C^0\cap H^1$ Green errors for any $p>1$ sufficiently close to one. These estimates are made separately on each of the finitely many initially independent radii before summing.

\smallskip\noindent\textit{Step 2. The residual density.} Let $\nu_i$ be the truncated densities of the fitted bubbles and set
\[
 Z=K(G\nu_1,G\nu_2)^t,\qquad w_u=u_2-Z_2,\qquad
 D=\max_aR_a.
\]
The exterior first tail is estimated before matching as in the proof of Theorem~\ref{hyb:capture}. On a fixed circle surrounding a bubble, the estimate for $e^{c_1}$ is a positive power of the radius of that bubble.

Boundedness of the mean-zero fields away from the disks about the blow-up points uses the estimate for small local mass with the fixed $L^p$ bound for the weight at each marked point without blow-up.

Multiplying the small normalization by $h_1e^{u_1}$ on this complement gives the positive-power $L^p$ bound for $e_1=f_1-\nu_1$. At a mixed prescribed point put
\[
 d=|z|,\qquad s=A+B
\]
and $d_0=R^\theta$ with $1<\theta<s/B$. Set
\[
 k=s/B-\theta>0,\qquad \lambda_p=2-2/p.
\]
Choose a cut-off
\[
 0\leq b_a\leq1
\]
equal to one on $D_{d_0}$ and supported in $D_{2d_0}$. Take the fitted second concentrated measure on this disk to be $b_a m_{B,r}\,dA$, using the same cut-off as for the solution density. The cut-offs on the different mixed disks are disjoint and their sum is denoted by $b$.

In the local estimates we use Euclidean densities. The smooth conformal area factor and its inverse are uniformly bounded, so the local $L^p$ estimates also hold for the corresponding surface densities.

The limiting coefficient satisfies $H_p\asymp d^{2A-2}$, by the prescribed weight and the two limiting Green charges. Polar integration gives
\[
 \|H_p\|_{L^p(D_{d_0})}\leq C d_0^{2A-\lambda_p}
 \quad\text{if }\lambda_p<2A.
\]
The same estimate on $D_{2d_0}$ controls $\|b_a H_p\|_p$. Since
\[
 r=R^{s/B+o(1)},
\]
the annulus $d_0\leq d\leq2d_0$ lies in $C_\varepsilon r\leq d\leq R/C_\varepsilon$ for small $R$. Throughout this larger annulus the second density of the solution is bounded by $C r^{2B-\varepsilon}d^{-2B-2+\varepsilon}$ by Lemma~\ref{hyb:singular-rate}. For
\[
 0<\varepsilon<2B
\]
its $p$th-power integral from $d_0$ to $R/C_\varepsilon$ is bounded by $C r^{p(2B-\varepsilon)}d_0^{2-p(2B+2-\varepsilon)}$. Independently, $m_{B,r}(d)\leq C r^{2B}d^{-2B-2}$ for $d\geq d_0$, so its $p$th-power integral from $d_0$ to infinity is bounded by $C r^{2Bp}d_0^{2-p(2B+2)}$. In particular,
\begin{align}
 \|f_2\|_{L^p(d_0<d<R/C_\varepsilon)}
     &\leq C r^{2B-\varepsilon}d_0^{-(2B-\varepsilon)-\lambda_p}
      =O(R^{(2B-\varepsilon)k-\theta\lambda_p+o(1)}),
       \label{univ:capture-actual-tail}\\
 \|\mathbf 1_{\{d\geq d_0\}}m_{B,r}\|_p
     &\leq C r^{2B}d_0^{-2B-\lambda_p}
      =O(R^{2Bk-\theta\lambda_p+o(1)}).
       \label{univ:capture-fitted-tail}
\end{align}
The second bound applies to the omitted fitted density $(1-b_a)m_{B,r}$ as well as to the fitted density on the transition annulus. In fact the norm of the omitted density on $2d_0<d<3d_0$ is comparable to $r^{2B}d_0^{-2B-\lambda_p}$. Dividing by $d_0^{2A-\lambda_p}$ gives
\[
 R^{2s(1-\theta)+o(1)}\to\infty.
\]
Thus the fitted tail and the deleted limiting coefficient require separate estimates.

Choose a fixed loss $\varepsilon>0$ smaller than $\min(2A,2B)$ at every prescribed point and smaller than $2$ at free points. Choose $p>1$ sufficiently close to one that
\begin{equation}\label{univ:capture-common-p}
 \lambda_p<2A-\varepsilon,\qquad
 \theta\lambda_p<(2B-\varepsilon)k
\end{equation}
at all mixed prescribed points. At pure prescribed points require $\lambda_p<\min\{2B,2(A+B)-\varepsilon\}$. At free mixed points the estimates between the inner and outer scales use $(s,B,\theta)=(2,1,3/2)$, so impose
\[
 (3/2)\lambda_p<(2-\varepsilon)/2,
\]
together with
\[
 p<3/2,\qquad p\varepsilon<2.
\]
Require $p$ also to lie in the interpolation range of the fitted errors and strictly below the local integrability exponents of both physical weights at unused sources.

At $p=1$ every displayed upper bound for $\lambda_p$ is strictly positive and the list is finite. Thus one fixed $p>1$ satisfies them simultaneously.

The three mixed prescribed powers, for the deleted limiting coefficient, the second density on the intermediate annulus, and the tail of the fitted profile, respectively, are $\theta(2A-\lambda_p)$, $(2B-\varepsilon)k-\theta\lambda_p$, and $2Bk-\theta\lambda_p$.

Choose a common exponent below all strictly positive limiting powers in the estimates below as well. This leaves room to absorb fixed logarithms and the $o(1)$ in the radius relation.

Inside $D_{d_0}$ the difference $b_a f_2-\nu_{2,a}$ has a positive power in $L^p$ by the fitted inner-profile estimate and interpolation. Here $\nu_{2,a}$ denotes the local part of $\nu_2$. On $D_{2d_0}\setminus D_{d_0}$ apply the triangle inequality using \eqref{univ:capture-actual-tail} and \eqref{univ:capture-fitted-tail}. The second density on
\[
 d_0\le d\le R/C_\varepsilon
\]
is bounded in $L^p$ by \eqref{univ:capture-actual-tail}. Beyond $R/C_\varepsilon$, the transition and exterior bounds for the actual second density are
\[
 Cd^{2A-2-\varepsilon}
\]
at a mixed prescribed point and $Cd^{-\varepsilon}$ at a free mixed point. Their $L^p$ norms on disks of shrinking fixed radius are uniformly small by the strict integrability inequalities above.

At a pure source the bounds are $Cd^{2(A+B)-2-\varepsilon}$ for $d\geq R$ and $CR^{2A-\varepsilon}d^{2B-2}$ inside $D_R$. The latter has norm $O(R^{2(A+B)-\varepsilon-\lambda_p})$.

At a pure or ordinary point the second local mass is zero, so no cut-off is needed for that component. The ordinary second-density bound $C(R+d)^{2-\varepsilon}$ and the $L^p$ bound at marked points without blow-up give the remaining local integrability.

For each small fixed disk radius $d_*$, the bounds just proved make the $L^p$ norms of these diffuse envelopes on the disks at most $Cd_*^\nu$ for some
\[
 \nu>0,
\]
up to errors from the fitted profiles and the intermediate annuli that tend to zero. The limiting residual density obeys the same local integrability bounds.

Outside these disks the fitted cut-offs on the second component vanish, and convergence away from the centers, with the fixed $L^p$ bound for the weight at each marked point without blow-up, gives convergence in $L^p$ there.

First take the concentration limit and then let $d_*\downarrow0$. This proves strong $L^p$ convergence of $f_2-\nu_2$ to the diffuse shadow density.

The exact equation \eqref{journal:scalar-remainder-equation} and the bounded Green map place $w_u$ in the $W^{2,p}$ inverse-function neighborhood of the selected scalar shadow. This step uses no quantitative coefficient comparison.

\smallskip\noindent\textit{Step 3. Comparison of the coefficients.} Put $H_T=h_2 e^{Z_2}$ and $A_D=(1-b)H_T$. The exact decomposition is
\[
 f_2-\nu_2=C_u A_D e^{w_u}+(bf_2-\nu_2),
 \qquad C_u=\frac{\rho_2}{\int_M h_2e^{u_2}}.
\]
On a mixed prescribed disk, extend the fitted inner profile to the fixed disk in the Green representation. Its omitted density is controlled by \eqref{univ:capture-fitted-tail}. The bounded map $G:L^p\to C^0$ thus gives a error bounded by a positive power of the scale in the logarithm of the coefficient.

For the smooth part of the Green kernel, subtract its value at the fitted center and use a profile moment of order
\[
 0<\eta<\min\{1,2B,2s\}.
\]
The planar logarithmic potential on a smaller disk receives a contribution from beyond the fixed disk bounded by a constant times the tail moment with weight $1+\log(1+d)$, which has a positive power by radial decay.

The planar potentials $-\log(d^{2q}+a^{2q})$ for the profile with exponent $q$ and width $a$, with the signs in $Z_2=-G\nu_1+2G\nu_2$, therefore give
\[
 H_T=(1+O(D^c))H_p Q_T,\qquad
 Q_T=\frac{1+(R/d)^{2s}}{(1+(r/d)^{2B})^2}
 \quad(d\geq d_0),
\]
after decreasing the positive exponent $c$ to absorb logarithms. On $d\geq d_0$, $|H_p(Q_T-1)|\leq C(R^{2s}d^{-2B-2}+r^{2B}d^{2A-2B-2})$. The first term has $L^p$ norm $O(R^{2Bk-\theta\lambda_p})$. For the second the norm is bounded by
\[
 C r^{2B}d_0^{\min(0,2A-2B-\lambda_p)}
                  (1+|\log d_0|)^{1/p},
\]
whose limiting power is $\min\{2s,\,2Bk-\theta\lambda_p+2A\theta\}>0$. These estimates give a uniform $L^p$ bound for $H_pQ_T$ on $d\geq d_0$. Consequently its multiplicative $O(D^c)$ error is small in $L^p$ there. On the full disk use
\[
 A_D-H_p=(1-b_a)(H_T-H_p)-b_aH_p.
\]
The first term uses only the estimates on $d\geq d_0$ and the second uses the $H_p$ bound on $D_{2d_0}$. The separate $b_a f_2-\nu_{2,a}$ error has already been estimated on the core and the transition annulus.

At a pure source the extra coefficient is $O(R^{2A}d^{2B-2})$, which belongs to $L^p$ when $p(1-B)<1$. At a free ordinary center the limiting coefficient is $d^2 k(z)$, with bounded $k(z)$, so the finite-width difference before the fitted Green error is $R^2 k(z)$. At a free mixed center take
\[
 d_0=R^{3/2},
\]
use the same cut-off arrangement and apply the comparison following \eqref{reg:actual-diffuse}. For $d\geq d_0$ the ratio differs from one by $O(R^4d^{-4}+(R^2+a^2)d^{-2})$ before the fitted Green error. Since
\[
 a=R^{2+o(1)},
\]
these terms have norms $O(R^{3/p-2})$ and $O(R^{3/p-1})$, respectively. The limiting coefficient is bounded on the removed disk, which lies away from all marked points. The core and transition errors in the second fitted concentrated density are estimated as above, using its profile exponent $q=1$.

On the complement of the bubble disks, the comparison of the separated Green kernels uses the same profile moments and mass errors. Multiplication by the prescribed physical weight there preserves the bound in the fixed $L^p$.

The boundedness of $w_u$, already obtained from preliminary strong $L^p$ convergence, permits multiplication by $e^{w_u}$.

On a fixed open set away from all centers the limiting coefficient is positive, which gives the upper bound for $C_u$ by integrating the nonnegative unnormalized density there.

For the lower bound, integrate the exact truncated decomposition. The coefficient $(1-b)h_2e^{Z_2}$ has a uniform $L^1$ bound by the comparisons just proved, $w_u$ is bounded, and
\[
 \int(bf_2-\nu_2)=o(1).
\]
Since
\[
 \int(f_2-\nu_2)\to\gamma>0,
\]
this gives $C_u\geq c>0$. In the exact decomposition set
\[
 e_2=C_u(A_D-H_p)e^{w_u}+(bf_2-\nu_2).
\]
For the fitted mass error, a radial inner model with exponent $q$ has mass $4\pi q$ and its mass outside $D_{d_0}$ is $4\pi q/(1+(d_0/r)^{2q})$.

This is a positive power of the corresponding outer radius, both for prescribed mixed points and for free mixed points of inner profile exponent one. The sum of the uncut inner masses is $\rho_2-\gamma$. Combining the estimates at their own radii and then using $D$ gives
\begin{equation}\label{univ:capture-scalar-errors}
 \begin{split}
 f_2-\nu_2&=C_uH_pe^{w_u}+e_2,\\
 \|e_1\|_p+\|e_2\|_p+
     \left|\int_M\nu_2-(\rho_2-\gamma)\right|
     &\leq CD^c(1+|\log D|)^N .
 \end{split}
\end{equation}
The fixed effective weight parameters, separated fitted centers and compact shadow branch give continuity and boundedness of $H_p$ in $L^p$, with a uniform lower bound on a fixed open set away from the centers.

The scalar linearization on the branch has a uniform inverse from mean-zero $L^p$ to mean-zero $W^{2,p}$.

The kernel space is the same as in the $C^0$ realization by the Green mapping, and the inverse is uniform by compactness of the nondegenerate branch. With the preliminary $W^{2,p}$ convergence already established and
\[
 \gamma>0,
\]
Lemma~\ref{journal:scalar-neighborhood} then gives the bound in \eqref{univ:capture-scalar-errors} also for $\|w_u-w(p)\|_{W^{2,p}}$.

\smallskip\noindent\textit{Step 4. Matching and the projected equation.} We apply Lemma~\ref{hyb:capture-matching}, using on the prescribed disks the exact matched profiles of Lemma~\ref{univ:profile-estimates}, and the pure-source corrections where appropriate.

The local rates, the exterior first tail and \eqref{univ:capture-scalar-errors} give the comparison estimates for the concentrated densities and the scalar remainder.

At each mixed point, $r/R^\theta$ is bounded by a positive power of $R$ and the inner cut-off radius is $R^\theta$ itself. Thus the fixed rescaled inner annuli lie where the cut-off is one, and fixed rescaled outer annuli lie in the cut-off of the first concentrated density. The same statements hold with
\[
 \theta=3/2
\]
at free mixed points. The annular matching and the summed first mass error consequently give
\[
 D\asymp e^{-T/(2s_*)},
\]
where $T=-c_1$, and
\[
 |\mu|\le Ce^{-c_0T},\qquad
 \|u-B_T\|_{C^0\cap H^1}\le Ce^{-c_0T}T^N
\]
for some $c_0>0$. For clarity, the change from the matched radial pair to the exact singular correction preserves this norm estimate. In the inverse Cartan components the two corrections have the form $-\log(1+(\varrho/d)^{2q}F_j(z))$ by the explicit formulas, where
\[
 (\varrho,q)=(R,s)
\]
or $(\ell,B)$. For the first component $F_j=F_0$, and for the second $F_j=s^2|P_{01}|^2/(A^2D)$, where $D=|P_{02}|^2+k|z|^{2A}|P_{12}|^2$ in the formula for $Q_R$. They satisfy
\[
 |F_j-1|+d|\nabla F_j|\le Cd^\beta
\]
for
\[
 \beta=\min(1,2A)>0.
\]
Subtracting the radial logarithm and differentiating the difference gives respective bounds $Cd^\beta\min(1,(\varrho/d)^{2q})$ and $Cd^{\beta-1}\min(1,(\varrho/d)^{2q})$. Splitting polar integration at $\varrho$ shows that the $C^0\cap H^1$ cost is
\[
 O(\varrho^{\min(\beta,2q)}(1+|\log\varrho|)^{1/2}).
\]
A fixed cut-off and mean subtraction preserve the bound.

At the matched seed use the scaled free-center and shape coordinates and the common logarithmic scale.

The local formulas give bounded tangents in $C^0\cap H^1$ in these scaled coordinates and the positive limiting source--tangent determinant in \eqref{hyb:chart}, hence a uniformly bounded matrix inverse. The centered sources satisfy \eqref{journal:source-functional-bounds}. Rescale each pairing by the width of its source.

On fixed annuli the matched profiles and their scaled derivatives converge uniformly. The other scales give constants or zero. The centered sources have integrable tails at both ends.

These bounds give a common modulus of continuity for the matrix. Lemma~\ref{journal:modulation-uniqueness} in $C^0\cap H^1$ therefore puts the error in $E_T$ while preserving a positive exponential rate.

If the largest outer profile exponent is one, the preceding estimates give exactly the preliminary hypothesis of Lemma~\ref{univ:ordinary-dominated}.

The three inequalities in its proof give the refined center, remainder, and parameter estimates, and uniqueness in that refined construction identifies the solution with its family.

Thus this case is settled without assuming its sharp scale law. In the other cases, choose the preliminary parameter exponent smaller than the exponent in
\[
 |\mu|\leq Ce^{-c_0T},
\]
so that the projected construction is available. The preliminary bound supplies the smallness needed for the following absorption, which estimates the remainder by the actual residual and places it in the chosen uniqueness ball. Write
\[
 E=\|\eta\|_{H^1}
\]
after imposing the orthogonality conditions. The solution satisfies
\[
 P_{E_T}\{\mathcal F_\mu(B_T)
       +D\mathcal F_\mu(B_T)\eta+\mathcal N(\eta)\}=0.
\]
The inverse bound and the weighted nonlinear estimate in Proposition~\ref{univ:projected-construction} give
\[
 E\leq CT\|\mathcal F_\mu(B_T)\|_{\Hh}
                      +CT^{5/2}E^2.
\]
The preliminary bound gives $T^{5/2}E\to0$. Thus the explicit smallness condition in Lemma~\ref{journal:projected-contraction} holds, and its absorption estimate gives
\[
 E\le C(\varepsilon_T+|\mu|)T^C.
\]
Choose the construction exponent $\delta<\min(a,c_0)$ from the residual and true parameter estimates. Then
\[
 \frac{E}{e^{-\delta T}T^{C_0}}
 \le C\bigl(e^{-(a-\delta)T}+e^{-(c_0-\delta)T}\bigr)T^{C'}
 \longrightarrow0 .
\]
The true remainder therefore lies strictly inside the contraction ball. The norm is that of the variational gradient obtained by applying $G$ to the residual of the differential equation.

This absorption does not weaken the decay rate chosen for $\varepsilon_T$ or its squared-error estimate relative to the leading scale.

Uniqueness identifies the remainder with the constructed one. The projected solution is a true solution exactly when all the source coefficients vanish, since the centered source vectors are independent by the invertible source--tangent matrix.

This argument places the solution in the projected family before the scale is determined.
\end{proof}
\subsection{Variation of the scale equation away from the sources}
\begin{lemma}[Variation of the reduced energy away from the sources]
\label{universal:exterior-envelope}
Fix a nondegenerate solution of the scalar shadow equation, with no regular blow-up points and $0<s_a<1$ at each prescribed bubble. Set
\[
 M=4\pi\sum_a s_a,\qquad \tau=e^{-T}.
\]
Let $K$ be compact and disjoint from the fixed source disks and all marked points. Take $0\le b\in C_c^\infty(K)$, $b\not\equiv0$, and put
\[
        h_{1,\varepsilon}^*=h_{1,0}^*e^{\varepsilon b},
        \qquad h_{2,\varepsilon}^*=h_2^*.
\]
\begin{itemize}
\item[\textup{(i)}]
The shadow solution, matched radii, approximation $B_T$, and model space $E_T$ are independent of $\varepsilon$.
\item[\textup{(ii)}]
Let $J_{\mu,\varepsilon}$ be the Toda functional with parameters $(M+\mu,\rho_2)$. Suppose the projected solutions
\[
 u_{T,\mu,\varepsilon}=B_T+\eta_{T,\mu,\varepsilon},
 \qquad \eta_{T,\mu,\varepsilon}\in E_T,
 \qquad P_{E_T}\nabla J_{\mu,\varepsilon}(u_{T,\mu,\varepsilon})=0
\]
are $C^2$ in their parameters. Write
\[
 \mathcal E(T,\mu,\varepsilon)
       =J_{\mu,\varepsilon}(u_{T,\mu,\varepsilon}),
 \qquad S(T,\mu,\varepsilon)=\partial_T\mathcal E(T,\mu,\varepsilon).
\]
Then
\begin{align}
 \partial_\varepsilon\mathcal E
   &=-\int_M b f_{1,T,\mu,\varepsilon}\,dV_g,
       \label{universal:energy-envelope}\\
 \partial_\varepsilon S
   &=-\partial_T\int_M b f_{1,T,\mu,\varepsilon}\,dV_g,
       \label{universal:scale-envelope}\\
 \partial_\mu S
   &=-\partial_T\log\int_Mh_{1,\varepsilon}
                       e^{u_{1,T,\mu,\varepsilon}}\,dV_g.
       \label{universal:mass-envelope}
\end{align}
where
\[
 f_{1,T,\mu,\varepsilon}
  =\frac{(M+\mu)h_{1,\varepsilon}e^{u_{1,T,\mu,\varepsilon}}}
          {\int_Mh_{1,\varepsilon}e^{u_{1,T,\mu,\varepsilon}}\,dV_g}.
\]
\end{itemize}
\end{lemma}
\begin{proof}
There are no free centers. The scalar shadow equation depends on $h_2^*$ and the fixed source configuration, but not on $h_1^*$.

On each source disk the perturbation is zero. Hence the matching data, the local corrections in $B_T$, and the source functionals defining $E_T$ are independent of $\varepsilon$.

At fixed $T$, both $\partial_\varepsilon\eta$ and $\partial_\mu\eta$ belong to $E_T$. Their pairings with the gradient of $J_{\mu,\varepsilon}$ therefore vanish. Differentiation of the functional gives
\[
 \partial_\varepsilon\mathcal E
   =-(M+\mu)
      \frac{\int_M b h_{1,\varepsilon}e^{u_1}\,dV_g}
           {\int_M h_{1,\varepsilon}e^{u_1}\,dV_g},
 \qquad
 \partial_\mu\mathcal E
   =-\log\int_Mh_{1,\varepsilon}e^{u_1}\,dV_g.
\]
Differentiating with respect to $T$ proves the remaining identities. In particular, these formulas include the dependence of the complementary correction on the weights.
\end{proof}
\begin{proposition}[Two-sided bounds for the scale equation]
\label{universal:positive-scale-without-limit}
Let $j$ range over finitely many of the shadow solutions with fixed sources in Lemma~\ref{universal:exterior-envelope}. For each fixed $E<\infty$, assume the projected construction is uniform on $[0,E]$ and that
\begin{equation}\label{universal:exterior-c1}
 \tau^{-1}\partial_T f_{1,j,T,0,\varepsilon}
       \longrightarrow -W_{j,\varepsilon}
       \quad\hbox{uniformly for $(x,\varepsilon)\in K\times[0,E]$},
 \qquad
 W_{j,\varepsilon}=W_{j,0}e^{\varepsilon b}>0.
\end{equation}
Suppose also that
\begin{equation}\label{universal:base-scale-bound}
                  |S_j(T,0,0)|\leq C_j\tau
                  \qquad(T\geq T_j).
\end{equation}
\begin{itemize}
\item[\textup{(i)}]
There is a finite $E>0$ such that, for all $j$,
\begin{equation}\label{universal:positive-scale-bounds}
                  c\tau\leq S_j(T,0,E)\leq C\tau
                  \qquad(T\geq T_*),
\end{equation}
Here $c,C>0$ are independent of $j,T$. No limit of $S_j(T,0,E)/\tau$ is assumed.
\item[\textup{(ii)}]
At this $E$, assume in addition that the reduction exists for large $T$ and $|\mu|\le e^{-\delta T}$ with fixed $0<\delta<1$, and that
\begin{equation}\label{universal:mu-derivative}
                \partial_\mu S_j(T,\mu,E)=-1+o(1)
\end{equation}
uniformly in this range. Each sufficiently concentrated zero in this range satisfies
\[
                    \mu>0,\qquad c_1\mu\leq\tau\leq C_1\mu.
\]
For fixed small $\mu>0$, the scale equation has degree $-1$ on $a\mu<e^{-T}<b_1\mu$ for suitable fixed $0<a<b_1$.
\end{itemize}
\end{proposition}
\begin{proof}
By \eqref{universal:scale-envelope} and \eqref{universal:exterior-c1}, uniformly for $0\leq\varepsilon\leq E$,
\[
 \tau^{-1}\partial_\varepsilon S_j(T,0,\varepsilon)
       \longrightarrow\int_M bW_{j,0}e^{\varepsilon b}\,dV_g.
\]
Integrating this derivative from $0$ to $E$ gives
\begin{equation}\label{universal:scale-change}
 S_j(T,0,E)-S_j(T,0,0)
   =\tau\int_MW_{j,0}(e^{Eb}-1)\,dV_g+o(\tau).
\end{equation}
The integral tends to $+\infty$ as
\[
 E\to+\infty,
\]
since $b$ is positive on a nonempty open set and $W_{j,0}>0$ there. Choose one finite $E$ such that this integral exceeds $C_j+2$ for every member of the finite collection.

The construction and convergence now apply on the fixed compact interval $[0,E]$. Equations \eqref{universal:base-scale-bound} and \eqref{universal:scale-change} prove \eqref{universal:positive-scale-bounds}.

For large $T$, \eqref{universal:mu-derivative} gives $-3/2\leq\partial_\mu S_j\leq-1/2$ throughout the stated range. Integrating at fixed $T$ from $0$ to $\mu$ gives $S_j\geq c\tau+|\mu|/2>0$ if $\mu\leq0$. If $\mu>0$, it gives
\[
 c\tau-\tfrac32\mu\leq S_j(T,\mu,E)\leq C\tau-\tfrac12\mu.
\]
Consequently each zero satisfies $(2c/3)\tau\leq\mu\leq2C\tau$. Choose $0<a<1/(2C)$ and $b_1>3/(2c)$. On the closed interval $a\mu\leq\tau\leq b_1\mu$, one has
\[
 \mu/\tau^\delta\leq a^{-\delta}\mu^{1-\delta}\to0.
\]
Thus the whole parameter segment from $0$ to $\mu$ belongs to the reduction domain for each $T$ in this interval. At the endpoints
\[
 \tau=b_1\mu,\qquad \tau=a\mu
\]
the preceding inequalities give strict positive and negative signs, respectively. As $T$ increases, $\tau$ decreases. These endpoint signs give degree $-1$, irrespective of the number of zeros.
\end{proof}
\begin{lemma}[Choice of the weights at each critical value]
\label{universal:wallwise-weights}
Assume that $1,(A_p)_{p\in S},(B_p)_{p\in S}$ are linearly independent over $\mathbb Q$.
\begin{itemize}
\item[\textup{(i)}]
A critical value of the form
\[
 \frac{\rho_1^{\rm wall}}{4\pi}
  = n+\sum_{p\in\mathcal P}A_p
       +\sum_{p\in\mathcal Q}(A_p+B_p),
 \qquad n=U+2V,
\]
determines $n,\mathcal P,\mathcal Q$ uniquely. If every outer profile exponent in one configuration is less than one, then $n=0$. All configurations at that value have $U=V=0$ and the same source sets.
\item[\textup{(ii)}]
Let $\rho_2>0$ and $\rho_1^\pm>0$ satisfy
\[
 \rho_2\notin\Gamma_2(\alpha),\qquad
 \rho_1^-<\rho_1^{\rm wall}<\rho_1^+,\qquad \rho_1^\pm\notin\Gamma_1(\alpha).
\]
Suppose there is no other critical value of the first parameter between $\rho_1^-$ and $\rho_1^+$. Then
\[
             d_\alpha(\rho_1^+,\rho_2)
                       -d_\alpha(\rho_1^-,\rho_2)
\]
is independent of the choice of positive smooth weights.
\end{itemize}
In Proposition~\ref{universal:positive-scale-without-limit}, one may therefore perturb the weight separately at each critical value with all outer profile exponents less than one. The shadow systems at other critical values may change.
\end{lemma}
\begin{proof}
\smallskip\noindent\textit{(i)}\quad In the expression for the critical value, the coefficient of $B_p$ is the indicator of $\mathcal Q$, and the coefficient of $A_p$ is the indicator of $\mathcal P\cup\mathcal Q$.

On subtracting two such representations, rational independence forces equality of each coefficient. A regular outer bubble has profile exponent one or two. Thus a configuration with all outer profile exponents less than one has
\[
 U=V=0,
\]
and uniqueness gives the same conclusion for every configuration at this critical value.

\smallskip\noindent\textit{(ii)}\quad Keep the source strengths and marked points fixed. Join any two pairs of positive smooth factors by linear interpolation of their logarithms. The resulting factors have uniform positive lower bounds and form a compact family in $C^3$.

At both noncritical pairs, Lemma~\ref{fnd:sharp-compactness} bounds all solutions along the path. The normalized Green maps define a compact homotopy on a common ball. Both endpoint degrees, and hence their difference, are independent of the weights.

If $n=0$ at a critical value, the corresponding shadow systems have only fixed source centers, and their scalar equations depend only on the second weight. The exterior perturbation of the first weight leaves all these shadow systems unchanged.
\end{proof}
\begin{proposition}[The scale equation when all outer profile exponents are below one]
\label{universal:all-small-wall}
Fix generic source strengths and a critical value of the first parameter with unique source configuration $(\mathcal P,\mathcal Q)$ and zero regular mass. Assume
\[
 0<A_p<1\quad(p\in\mathcal P),\qquad
 0<A_p+B_p<1\quad(p\in\mathcal Q).
\]
Choose $h_2^*$ so that the scalar shadow equation has a finite compact set of nondegenerate solutions. Write
\[
 M=4\pi\left(\sum_{p\in\mathcal P}A_p
                    +\sum_{p\in\mathcal Q}(A_p+B_p)\right),
 \qquad \rho_1=M+\mu,\qquad \tau=e^{-T}.
\]
Use the matched exact profiles and the spaces $E_T$. Make the following assumptions.
\begin{itemize}
\item[\textup{(a)}]
For fixed $0<\delta<1$, the projected construction on $|\mu|\le\tau^\delta$ satisfies
\begin{equation}\label{universal:small-projection-estimates}
 \|\eta\|_{H^1}+\|\partial_T\eta\|_{H^1}
       \leq C(\tau+|\mu|)T^N,
 \qquad
 \|\eta\|_\infty+\|\partial_T\eta\|_{C^0(K)}=o(1).
\end{equation}
Here $K$ is any fixed compact set disjoint from the source disks and all marked points.
\item[\textup{(b)}]
The normalized mass and its $T$ derivative satisfy the estimates obtained by integrating the matched profiles. The construction depends smoothly on the weights.
\item[\textup{(c)}]
These estimates hold uniformly on each fixed compact interval of smooth perturbations of $h_1^*$ supported outside the source disks and away from the marked points.
\end{itemize}
Then $h_1^*$ can be chosen so that the following hold for every scalar shadow solution.
\begin{itemize}
\item[\textup{(i)}]
The reduced scale equation $S=\partial_T\mathcal E$ satisfies
\[
             c\tau\leq S(T,0)\leq C\tau,
             \qquad \partial_\mu S(T,\mu)=-1+o(1).
\]
There are no zeros for $\mu\le0$. Every zero in the projected range satisfies $\tau\asymp\mu>0$, and the scale degree is $-1$.
\item[\textup{(ii)}]
If the projected families contain every sufficiently concentrated sequence of solutions whose parameters tend to this critical value, the local degree jump at a scalar shadow solution $Z$ is
\[
                    (-1)^{|\mathcal P|}d_S(Z).
\]
\end{itemize}
\end{proposition}
\begin{proof}
\smallskip\noindent\textit{Step 1. The mass estimates.} Choose $p>1$ sufficiently close to one that
\[
 ps_a<1
\]
for every outer profile exponent $s_a$, and so that all singular densities belong locally to $L^p$. At a mixed source the inner profile exponent $B_p$ is smaller than the outer profile exponent, so $pB_p<1$ as well. The estimates for
\[
 ps_a<1,\qquad pB_p<1
\]
in Lemma~\ref{univ:profile-estimates}, applied to the exact residual with its quadratic Taylor remainders, give $O(\tau)$ for each mixed disk and for its logarithmic-scale derivative. The fixed-annulus errors and their $T$ derivatives are $O(\tau)$.

For a pure source, apply the case $0<A_p<1$ of Lemma~\ref{univ:pure-residual}, including its angular correction if $A_p\geq1/2$. In addition to
\[
 ps_a<1,
\]
choose $p$ sufficiently close to one to satisfy the strict bounds in that lemma for the linear term and its integrability at each of the finitely many pure sources. The bound is $O(\tau)$ both for the residual and for its $T$ derivative.

Let $\mathfrak n_i$ denote the unnormalized approximate densities. Boundary integration of the exact local identities
\[
 n_1=\Delta\log F,\qquad n_2=\Delta\log Q
\]
gives the limiting bubble masses with errors $O(\tau)$ at a fixed disk boundary. The Taylor remainders after matching are weighted by $n_{1,R}$ and by $n_{2,R}-f$, since the exterior residual density is $fe^{J_2}$.

The estimates for $s_a<1$ bound their second absolute moments and the corresponding moments after one $T$ derivative by $O(\tau)$.

At a pure source use the expansion in Lemma~\ref{univ:pure-residual}. For $A<1/2$ its first remainder is $O(r^\nu)$ for $\nu>2A$. For $1/2\leq A<1$ the linear term and its first-order correction are odd, so that their integrals against the radial outer density vanish.

The remaining first error is bounded by $C m_{A,R}r^\nu$ for $\nu>2A$, and has integral $O(\tau)$. The second error is bounded in $L^1$ by $C\tau$ by the same lemma, using in particular
\[
 \tau r^{2B-2}\in L^1.
\]
These bounds also hold after one $T$ derivative, since the differentiated linear terms remain odd and the differentiated radial factors have the same bounds. Consequently
\[
 \int_M\mathfrak n_1=M+O(\tau),\qquad
 \int_M\mathfrak n_2=\rho_2+O(\tau),
\]
and their derivatives in $T$ are $O(\tau)$. Global normalization adds an $O(\tau)$ multiple of each density. Its $L^1$ size is $O(\tau)$ and its $H^{-1}$ size is $O(\tau T^N)$.

The field $\partial_T B_T$ is uniformly bounded in $C^0$. For the exact mixed profiles this follows directly by logarithmic differentiation: the resulting expressions are sums of bounded ratios of the positive terms in $F$ and $Q$. For a scalar profile it follows from
\[
 \partial_T D_{q,R}=Z_{q,R}+1.
\]
Pairing the normalized approximate equation with $\partial_TB_T$ therefore gives
\begin{equation}\label{universal:trial-scale-bound}
               |\partial_TJ_{0}(B_T)|\leq C\tau.
\end{equation}
This estimate uses the residual itself and does not require the existence of a limiting finite-part coefficient.

\smallskip\noindent\textit{Step 2. The scale equation at the critical parameter.} At $\mu=0$, the correction in \eqref{universal:small-projection-estimates}, the residual bound and its $T$ derivative give
\begin{equation}\label{universal:energy-correction-scale}
 \partial_T\{J_0(B_T+\eta)-J_0(B_T)\}
                         =O(\tau^2T^{N_1})=o(\tau).
\end{equation}
Indeed, the chain rule and symmetry of the second derivative give the difference $\langle\partial_T\mathcal F_0(B_T),\eta\rangle_{\Hh}+\langle\mathcal N_T(\eta),\partial_T B_T\rangle_{\Hh}+\langle\mathcal F_0(B_T+\eta),\partial_T\eta\rangle_{\Hh}$. The first term is
\[
 O(\tau^2T^{N_1})
\]
by the differentiated residual estimate, and the second is $O(T\|\eta\|_{\Hh}^2)$ by Lemma~\ref{journal:normalized-nonlinearity} with the bounded argument $\partial_T B_T$. For the last term we use
\[
 \|\mathcal F_0(B_T+\eta)\|_{\Hh}\le C\tau T^C
\]
by the trial residual and the polynomial operator bound for the second derivative along the segment to the correction. The bound for $\partial_T\eta$ in \eqref{universal:small-projection-estimates} then gives the same order. Equations \eqref{universal:trial-scale-bound} and \eqref{universal:energy-correction-scale} prove $|S(T,0,0)|\leq C\tau$ for each unperturbed shadow solution.

To compute the mean subtraction, put
\[
 \tau_a=R_a^{2s_a}\asymp \tau.
\]
On a mixed disk the two corrections are linear combinations of $\log(1+H_R/D)$ and $\log(1+\tau_a F_0/r^{2s_a})$. Each and its $T$ derivative is bounded by $C\tau(r^{-2B}+r^{-2s_a})$.

On a pure disk use the corresponding bounds $C\tau r^{-2A}$ for the radial logarithm and $C\tau r^{1-2A}$ for the angular term. The same bounds hold after one $T$ derivative.

Since the outer and inner profile exponents are below one and the exterior fields have zero mean, the subtracted means and their derivatives are $O(\tau)$. Writing
\[
 m_1(T)=\int_M\widehat B_{T,1},
\]
one has $\tau I_1=e^{-m_1(T)}\int_M\mathfrak n_1 e^{\eta_1}$. The integrated profile estimates and the correction bounds give
\[
 e^{-T}I_1(T,\mu,\varepsilon)
       =M+O((\tau+|\mu|)T^{N_2}),
 \qquad
 \partial_T(e^{-T}I_1)
       =O((\tau+|\mu|)T^{N_2}).
\]
For the derivative in this estimate, use the exact formula
\[
 \partial_T(\tau I_1)=e^{-m_1}\left(
  \int_M \partial_T\mathfrak n_1\,e^{\eta_1}
  +\int_M \mathfrak n_1 e^{\eta_1}\partial_T\eta_1
  -m_1'\int_M \mathfrak n_1 e^{\eta_1}\right).
\]
Here $m_1'=O(\tau)$, and $|\partial_T\mathfrak n_1|\leq C\mathfrak n_1$ by the explicit scale derivatives in the projected construction. Remove $\int_M\partial_T\mathfrak n_1=O(\tau)$ from the first term. For both remaining correction integrals use
\[
 \int_M(\mathfrak n_1/\textstyle\int_M\mathfrak n_1)
       (|\eta_1|+|\partial_T\eta_1|)
 \leq C\sqrt T(\|\eta\|_{H^1}+\|\partial_T\eta\|_{H^1}),
\]
which follows by H\"older from the uniform $L^{1+1/T}$ probability bound and the two-dimensional Sobolev bound in $L^{T+1}$. Since
\[
 \|\eta\|_\infty=o(1),
\]
the exponential factors preserve these bounds. This proves the differentiated estimate above. Hence
\[
 \partial_T\log I_1=1+\frac{\partial_T(\tau I_1)}{\tau I_1}
                  =1+o(1)
\]
uniformly on $|\mu|\leq\tau^\delta$. On $K$ the uncorrected field and its $T$ derivative converge respectively to $u_1^0$ and zero, including the mean subtraction. The $C^0$ estimates on the complementary correction and its derivative give the same limits for $u_1$ there. At
\[
 \mu=0,
\]
$\partial_T f_1=f_1(\partial_Tu_1-\partial_T\log I_1)$, and so
\[
 \tau^{-1}\partial_T f_{1,T,0,\varepsilon}
 \longrightarrow-h_{1,0}e^{\varepsilon b}e^{u_1^0}
\]
uniformly for $(x,\varepsilon)\in K\times[0,E]$.

\smallskip\noindent\textit{Step 3. The sign of the scale equation and its degree.} Apply Proposition~\ref{universal:positive-scale-without-limit} to the finite family of scalar shadow solutions. It provides one finite exterior perturbation
\[
 h_1^*\mapsto h_1^*e^{Eb}
\]
and the positive upper and lower bounds for $S(T,0,E)$. The exact identity \eqref{universal:mass-envelope} gives
\[
 \partial_\mu S=-1+o(1).
\]
The same proposition excludes zeros for $\mu\leq0$, gives $\tau\asymp\mu$ at each zero, and gives degree $-1$. The whole closed interval used in that proposition lies strictly inside the wider construction range. In particular $\mu/\tau^\delta$ tends to zero uniformly on it.

Only the radial mode at the reference point remains. Its pairing with $\partial_T(B_T+\eta)$ tends to
\[
 4\pi s_{\rm ref}/3>0.
\]
For the correction this follows by differentiating its orthogonality and using the uniform bound for the differentiated source functional. Thus $S=0$ is equivalent to the remaining projected equation.

For the degree product, use the graph homotopy from the proof of Theorem~\ref{hyb:transfer} on the scale interval just obtained. There are no free center or shape parameters.

On the closed scale interval take the closed fiber ball about the graph of radius $\tfrac12e^{-\delta T}T^{C_0}$, increasing $C_0$ as in Proposition~\ref{univ:projected-construction}. On this interval
\[
 \|\eta\|_{H^1}=O(\tau T^N)
\]
and is $o(e^{-\delta T}T^{C_0})$. The graph and this fiber ball are inside the original uniqueness ball. The derivative of the projected equation at $B_T+\eta$, restricted to $E_T$, has an inverse of norm $O(T)$.

The fixed exponential $L^{p_0}$ profile bound makes Lemma~\ref{journal:normalized-nonlinearity} applicable along the segment from the graph to each fiber point, because
\[
 \sqrt T\|\eta+av\|_{H^1}=o(1)
\]
for $0\leq a\leq1$. Its remainder is at most $CT^{3/2}\|v\|^2$. Thus the homotopy from the complementary equation to its linear part has every zero at $v=0$, by
\[
 CT^{5/2}\|v\|=o(1).
\]
In the same homotopy we can fix the finite-dimensional component at its value on the graph. Projection transport and the positive reference source pairing give the same orientation as in the proof of Theorem~\ref{hyb:transfer}.

The scale equation has strict signs at both endpoints, and the complementary equation has no zero on the boundary sphere.

The product rule therefore gives the complementary degree times the degree of the scale equation. Since all outer and inner profile exponents are below one, their local scalar indices vanish. The complementary index is
\[
 \nu=m_S+|\mathcal P|+2|\mathcal Q|-1.
\]
Multiplication by the scale degree $-1$ gives $(-1)^{|\mathcal P|}d_S(Z)$. The hypothesis on all blowing-up sequences ensures that this gives the full contribution of $Z$ to the degree jump at
\[
 \rho_1=M.
\]
On this scale interval,
\[
 \tau\asymp\mu
\]
bounds the ratio of the two reference widths. Suppose two approximations with nearby parameters differ by $o(1)$ in $H^1$.

At the reference scale, the gradients of their first inverse-Cartan components have the radial limit obtained by comparing the approximations in Proposition~\ref{univ:wall-jump}.

That comparison uses only $s>0$. Thus the two $T$ values differ by $o(1)$, and Lemma~\ref{journal:modulation-uniqueness} proves injectivity. On this interval we can again choose the radius of the complementary ball strictly greater than the bound for the correction.
\end{proof}
\begin{remark}
Condition \eqref{universal:exterior-c1} follows, for example, from the exterior estimates
\[
 u_{1,T,0,\varepsilon}\to u_1^0,\qquad
 \partial_Tu_{1,T,0,\varepsilon}\to0,
 \qquad
 \frac{M}{\int_Mh_{1,\varepsilon}e^{u_1}\,dV_g}
          =\tau(1+o(1)),
\]
and $\partial_T\log\int_Mh_{1,\varepsilon}e^{u_1}\,dV_g=1+o(1)$. Then
\[
 W_{j,0}=h_{1,0}e^{u_{1,j}^0}\quad\text{on }K .
\]
For the projected solutions these estimates follow in the proof of Proposition~\ref{universal:all-small-wall}.
\end{remark}
\begin{remark}
By Lemma~\ref{universal:wallwise-weights}, rational independence gives $U=V=0$ for every configuration at this critical value. We may perturb the weight there independently of the shadow systems at other critical values.
\end{remark}

\section*{Acknowledgement}
The research of Y. Lee is supported by a National Research Foundation of Korea (NRF) grant funded by the Korea government (MSIT), No. RS-2026-25489601. 
The research of J. Wei is partially supported by the Hong Kong General Research Fund of the Research Grants Council of Hong Kong under the project ``New frontiers in singularity formations in nonlinear partial differential equations.'' The research of Wen Yang was supported by the National Key R\&D Program of China 2022YFA1006800, NSFC No. 12531010, FDCT No. 0070/2024/RIA1, and University of Macau research grants UMDF-TISF/2025/006/FST, SRG2023-00067-FST, MYRG-GRG2024-00082FSC-UMDF, MYRG-GRG2025-00051-FST, MYRG-GRG2026-00034-FSC.
The research of L. Zhang is partially supported by Simons Foundation Grant SFI-MPS-TSM-00013752.

We declare that AI assistance (OpenAI ChatGPT 6.2 Astra and Codex) was used in the preparation of this manuscript.

\appendix
\section{Auxiliary estimates for the blow-up analysis}
\label{journal:capture-topology}
The following estimates concern interpolation at singular points and the choice of parameters by orthogonality. Throughout, $G$ is the mean-zero Green operator on the surface $M$, and $\Pi f=f-\int_M f$; the area of $M$ is one.
\begin{lemma}[An estimate for the truncated inner profile]
\label{journal:inner-tail}
Let $r=R^{k+o(1)}$, $k>1$, $t>0$ and $1<\theta<k$. Suppose
\[
 |f(z)|\le Cr^{2t-\varepsilon}d^{-2t-2+\varepsilon},\qquad R^\theta<d=|z|<R.
\]
For fixed $p>1$ sufficiently close to one and fixed $\varepsilon>0$ sufficiently small, there is a $c>0$ such that $\|f\|_{L^p(R^\theta<d<R)}=O(R^c)$.
\end{lemma}
\begin{proof}
Choose \(0<\varepsilon<2t\), so that \(p(2t+2-\varepsilon)>2\). The radial integral is bounded by its lower-end contribution:
\[
 \int_{R^\theta<d<R}|f|^p
   \leq C r^{p(2t-\varepsilon)}
                   R^{\theta(2-p(2t+2-\varepsilon))}.
\]
Writing \(\lambda_p=2-2/p\), the exponent for the norm is
\[
 (2t-\varepsilon)(k-\theta)-\theta\lambda_p+o(1).
\]
It is strictly positive for sufficiently small fixed \(\varepsilon\) and fixed \(p>1\) sufficiently close to one. The corresponding fitted radial profile satisfies the same calculation with
\[
 \varepsilon=0.
\]
For finitely many bubbles choose one pair \(p,\varepsilon\) by using the minimum strict margin. The bound outside \(R\) must be estimated separately using the local bound for the residual density.
\end{proof}
\begin{proof}[Proof of Lemma~\ref{journal:scalar-neighborhood}]
\smallskip\noindent\textit{Step 1. Convergence in (i).}\quad The assumed strong $L^p$ convergence, followed by the bounded Green map, gives
\[
 w_u\to w_*
\]
strongly in $W^{2,p}$, where $w_*$ is the selected mean-zero shadow solution. The branch is continuous in $C^0$ and $H_p$ is continuous in $L^p$, so its normalized right-hand sides converge in $L^p$ by the uniform lower bound on the integrals. The mean-zero Green estimate then gives
\[
 w(p)\to w_*
\]
in $W^{2,p}$ as well, proving $w_u-w(p)\to0$ in this space.

\smallskip\noindent\textit{Step 2. The estimate in (ii).}\quad In particular $w_u$ is bounded in $C^0$, and $I_u=\int_MH_pe^{w_u}$ is bounded above and away from zero. Integrating the second displayed comparison determines its coefficient:
\[
 C_uI_u=\rho_2-\int_M\nu_2-\int_Me_2
       =\gamma+O(\delta_D).
\]
Put $\mathcal D_u=H_pe^{w_u}/I_u$. The normalized scalar residual is exactly
\[
 -\Delta_gw_u-2\gamma\Pi\mathcal D_u
 =2(C_uI_u-\gamma)\Pi\mathcal D_u+2\Pi e_2-\Pi e_1.
\]
Since \(W^{2,p}(M)\hookrightarrow C^0(M)\), the \(L^p\) bound for \(H_p\) and the positive lower bound of the denominator give a uniform \(L^p\) bound for \(\mathcal D_u\). Thus the displayed mean-zero defect has \(L^p\) norm \(O(\delta_D)\). Write
\[
 N_p(w)=H_pe^w/\int_MH_pe^w.
\]
Its derivative in direction \(v\) is \(N_p(w)(v-\int_MN_p(w)v)\). The \(C^0\) bound of the branch and the lower bound of the denominator also bound the Taylor remainder in \(L^p\) by \(C\|v\|_{C^0}^2\), uniformly on a fixed small neighborhood of the branch.

Indeed the exponential has this remainder uniformly in \(C^0\), and the numerator is then multiplied by the uniformly \(L^p\)-bounded \(H_p\). Set
\[
 v=w_u-w(p)
\]
and apply the uniformly bounded inverse of the linearization at \(w(p)\) to the scalar defect. This gives
\[
 \|v\|_{W^{2,p}}\leq C\delta_D+C\|v\|_{W^{2,p}}^2.
\]
The strong convergence already established permits absorption of the last term, hence the asserted estimate, also for a weight unbounded at a marked point.

When fitted-model errors are used to verify these hypotheses, apply Lemma~\ref{journal:independent-radii} when each individual error has the stated positive-power $L^1$ bound and a polynomial $L^\infty$ bound in its own radius.

At a prescribed point where the latter bound is replaced by a polynomial $L^{p_0}$ bound for one fixed
\[
 p_0>1,
\]
apply Lemma~\ref{univ:singular-interpolation}. Choose a single admissible $p>1$ for the finite list of errors and for the diffuse coefficient, before summing over the individual radii.

For cut-offs of the solution densities, estimate the remaining densities directly. The exponent in this interpolation may differ from that used in the projected construction.

For the local tails, separate any signed fitted-model error, the vanishing inner tail of Lemma~\ref{journal:inner-tail}, and the exterior pointwise bound. The latter is $Cd^{2A-2-\varepsilon}$ at a mixed singular point and $Cd^{-\varepsilon}$ at a mixed free point.

With the strict choices of $p,\varepsilon$, the $L^p$ integrals of these bounds on disks of radius $d_0$ about the fitted centers are bounded by $C d_0^\upsilon$ for some
\[
 \upsilon>0,
\]
uniformly even when the centers move. The tails and signed fitted errors already tend to zero in $L^p$. For fixed $d_0$ use convergence outside the corresponding disks about the limiting centers, then let $d_0$ decrease to zero using the uniform disk bounds. This gives strong $L^p$ convergence.

These estimates are uniform with respect to the ratios of the different outer radii and precede the sharp mass estimate.
\end{proof}
To prove exhaustion of the isolating neighborhoods, pass to a subsequence of any escaping solution sequence with one of the selected shadow limits.

On branches using the preceding lemma, the local and residual estimates first put the solution in the projected family after matching the scales and imposing the orthogonality conditions.

The appropriate scale law then puts it strictly inside the isolating scale interval on its allowed side.

If all bubbles are ordinary, testing the equation on the initial parameter interval first gives the mass estimate and hence places the parameters in the smaller interval used for the construction. This argument applies before identifying the solution in that smaller interval.

In either argument the fitted center and shape limits are interior, and the correction on the graph is the unique solution of the projected equation. Thus each such sequence has a subsequence in one of the neighborhoods on the corresponding side of
\[
 \mu=0.
\]
\subsection{Normalized nonlinearities and the projected equation}
We retain the normalizing integrals. All Sobolev spaces here have zero mean; we use the Dirichlet norm and, for pairs, the fixed Cartan norm.

The constants are uniform in the concentration parameter and may depend on the surface and the compact set of profile exponents. For the applications above, if
\[
 \mathscr R=-\Delta_gU-K(\Pi f_i)_i
\]
is the mean-zero differential residual, then its variational gradient is
\[
 \mathcal F_T(U)=G\mathscr R
\]
by \eqref{eq:toda-gradient-identity}. For a mean-zero unmultiplied density error $(r_i)_i$, the corresponding gradient error is $GK(r_i)_i$.

These fixed operators have the same bounds from $H^{-1}$ to $\Hh$ up to a fixed constant. For a finite collection of fixed local weights $r^{2(q-1)}$ with $q>0$, choose
\[
 p_0>1\quad\Longrightarrow\quad p_0(q-1)>-1
\]
for every weight. In polar coordinates these powers belong to $L^{p_0}$. Together with the explicit profile scales and the positive mass normalization, this verifies the exponential $L^{p_0}$ bound in Lemma~\ref{journal:normalized-nonlinearity} when the profile exponents are fixed.

The integrability margin can be chosen on a compact subset of $q>0$. The inverse of the projected operator also depends on the kernel space of the linearized equation, which changes when $q$ is an integer. Continuation in the source strengths uses the compactness theorems.
\subsection{Weighted compactness, moments and energy differentiation}
\begin{proof}[Proof of Lemma~\ref{hyb:moment-estimates}]
\smallskip\noindent\textit{Step 1. Moments of the radial density.}\quad Polar integration and $z=Ry$ give
\[
 \big\||z|^b m_{q,R}\big\|_p^p
 =2\pi(4q^2)^p R^{p(b-2)+2}
   \int_0^{d_0/R}
   \frac{y^{p(b+2q-2)+1}}{(1+y^{2q})^{2p}}\,dy.
\]
The hypothesis gives integrability at zero. At infinity the exponent is $p(b-2q-2)+1$, giving the three alternatives. For fixed $z$ put $x=(|z|/R)^{2q}$. Then $\partial_Tx=x$, and
\[
 \partial_T\log m_{q,R}=\frac{1-x}{1+x},\qquad
 \partial_T^2\log m_{q,R}=-\frac{2x}{(1+x)^2}.
\]
Both rational functions are bounded. Consequently $|\partial_T^j m_{q,R}|\leq C m_{q,R}$ for $j=1,2$, which gives the differentiated estimates. Uniformity of a noncritical branch over varying parameters uses a compact positive margin from its integrability and moment thresholds.

\smallskip\noindent\textit{Step 2. Terms containing the inner scale.}\quad The matching relation gives the pointwise inequality
\[
 Sh\leq C\frac{|z|^{2A-2}}{(1+(|z|/R)^{2s})^2}.
\]
Scaling by $R$ proves both overlap bounds; the integral at infinity converges since $s=A+t$. The same inequality holds after two logarithmic derivatives.

\smallskip\noindent\textit{Step 3. The first-order terms.}\quad For the coupling moments put $d=|z|$ and split at $r$ and $R$. On $r<d<R$, both products are bounded by $Cd^{2A-1}$, so that their first moments are bounded by
\[
 C\int_r^R d^{2A+1}\,dd\leq CR^{2A+2}.
\]
On $d<r$, both are bounded by $CR^{-2s}d^{2s-1}$, and their first moments are at most $CR^{-2s}r^{2s+2}$. Using $r^{2t}\asymp R^{2s}$ shows that this is bounded by $CR^{2A+2}$, since the difference of the corresponding exponents is $2A(A+1)/t>0$. On $d>R$, the two bounds are
\[
 C R^{2s}r^{2t}d^{-2s-2t-1},\qquad
 C r^{2t}d^{-2t-1}.
\]
After multiplication by $d$ and the polar area factor, their integrals are bounded respectively by $CR^{2A+2}$ and $CR^{2A+2}$. The second bound uses $t>1$, since then
\[
 \int_R^{d_0}d^{1-2t}\,dd=O(R^{2-2t}).
\]
Each first-order correction is a fixed linear form in $z$ times a rational radial factor. Since its coefficients are independent of $T$, two logarithmic derivatives preserve all these majorants.
\end{proof}
\begin{lemma}[Differentiation of the reduced energy]
\label{hyb:energy-differentiation}
Take the $T$ derivative with the mass defect $\mu$, the center parameters and the shape parameters fixed. The projected solution $u(T)=B_T+\eta(T)$ satisfies
\begin{equation}
 \frac{d}{dT}\mathcal J_\mu(u(T))
       -\frac{d}{dT}\mathcal J_\mu(B_T)
   =O\bigl((\epsilon(T)^2+\mu^2)T^C\bigr).
\end{equation}
The estimate is uniform on the domain of Proposition~\ref{hyb:range}. For finite $T$, the implicit function theorem in the varying complementary space gives differentiability. Differentiating the orthogonality conditions controls the term containing $\dot\eta$.
\end{lemma}
\begin{proof}
\smallskip\noindent\textit{Step 1. Differentiation of the projected equation.}\quad A dot denotes differentiation in $T$ and $\eta$ is the correction to $B_T$. Put $E=\|\eta\|_{\Hh}$ and
\[
 r_T=(\epsilon(T)+|\mu|)T^{C_1},
\]
increasing the fixed exponent $C_1$ when necessary. Proposition~\ref{hyb:trial-comparison} and the projected equation give
\[
 \|\mathcal F_\mu(B_T)\|_{\Hh}
 +\|\partial_T\mathcal F_\mu(B_T)\|_{\Hh}\leq C r_T,
 \qquad E+|c|\leq C T^{C_2}r_T.
\]
The explicit profiles give
\[
 \|\dot B_T\|_\infty\leq C.
\]
Indeed a logarithmic-scale derivative of every radial correction is bounded, and the logarithmic derivatives of $F$ and $Q_0$ in \eqref{hyb:gram} are bounded because these positive expressions are comparable to their sums of nonnegative radial leading terms on the compact shape set.

The differentiated first-order corrections are bounded as well. Subtracting the surface mean preserves this bound.

Lemma~\ref{journal:normalized-nonlinearity}, with its third argument equal to the bounded tangent $\dot B_T$, gives
\[
 |\langle\mathcal N_T(\eta),\dot B_T\rangle_{\Hh}|
       \leq CTE^2.
\]
The probabilities at $B_T$ have a fixed $L^{p_0}$ margin and at most exponential growth in that norm, and $\sqrt T E\to0$. Thus all hypotheses of that lemma hold. Both normalization denominators are included in this estimate.

\smallskip\noindent\textit{Step 2. The energy derivative.}\quad The chain rule and self-adjointness now give the exact decomposition
\begin{align*}
 &\frac{d}{dT}\mathcal J_\mu(B_T+\eta)
             -\frac{d}{dT}\mathcal J_\mu(B_T)\\
 &\quad=\langle\partial_T\mathcal F_\mu(B_T),\eta\rangle_{\Hh}
 +\langle\mathcal F_\mu(B_T+\eta)-\mathcal F_\mu(B_T)
                -D\mathcal F_\mu(B_T)\eta,\dot B_T\rangle_{\Hh}
 +\langle\mathcal F_\mu(B_T+\eta),\dot\eta\rangle_{\Hh}.
\end{align*}
For the last term, differentiate the orthogonality constraints and use the exact projected equation. This gives
\[
 \langle\mathcal F_\mu(B_T+\eta),\dot\eta\rangle_{\Hh}
 =-\sum_\alpha c_\alpha(\partial_TL_\alpha)\eta.
\]
By \eqref{hyb:modulation-bounds} its absolute value is at most $C|c|E$. The other two terms are bounded by $Cr_TE$ and $CTE^2$. Using $E+|c|\leq CT^{C_2}r_T$ in all three terms bounds their sum by $CT^C r_T^2=O((\epsilon(T)^2+\mu^2)T^{C'})$, proving the lemma.
\end{proof}
\subsection{Quantitative scalar profile comparison}
\begin{lemma}[From mass exhaustion to a quantitative profile estimate]
\label{univ:profile-rate}
Fix a noninteger $q>0$, $0<\kappa\leq\min(1,q)$, and $0<\zeta<2q$. Write
\[
 W_q(y)=\log\frac{4q^2}{(1+|y|^{2q})^2},\qquad
 m_q=|y|^{2q-2}e^{W_q},\qquad \ell(y)=\log(2+|y|).
\]
Let $L_n\to\infty$, and let $w_n$ be bounded near the origin, with $w_n(0)$ taken in its continuous representative, and solve distributionally on $D_{L_n}$
\[
 -\Delta w_n=2|y|^{2q-2}e^{w_n}+F_n,
 \qquad w_n(0)=W_q(0).
\]
Let $a_n\to0$. Suppose the following estimates hold.
\begin{itemize}
\item[\textup{(i)}]
We have $w_n\to W_q$ locally uniformly. The circular oscillations of $w_n-W_q$ are uniformly bounded and
\[
 \lim_{M\to\infty}\limsup_{n\to\infty}
 \sup_{M\leq d\leq L_n}
 \left|\int_{D_d}(|y|^{2q-2}e^{w_n}-m_q)\,dy\right|=0.
\]
\item[\textup{(ii)}]
The error satisfies
\[
 \sup_{D_{L_n}\setminus\{0\}}
       \frac{|F_n|}{\omega_{\kappa,\zeta}}
       +\operatorname{osc}_{\partial D_{L_n}}(w_n-W_q)
       \leq a_n,
\]
\end{itemize}
Then, for all sufficiently large $n$,
\[
 \sup_{D_{L_n}}\frac{|w_n-W_q|}{\ell}\leq Ca_n.
\]
The assertion remains valid after fitting the scalar radius by a factor tending to one so that the condition at zero holds.
\end{lemma}
\begin{proof}
\smallskip\noindent\textit{Step 1. The circular means.}\quad Put $e_n=w_n-W_q$ and denote circular means by a bar. Testing the distributional equation by radial functions that approximate the indicator of $D_d$ gives, for almost every $d$,
\[
 -2\pi d\bar e_n'(d)
 =2\int_{D_d}(|y|^{2q-2}e^{w_n}-m_q)\,dy
       +\int_{D_d}F_n\,dy.
\]
Both integrands are locally integrable at zero. There is no additive flux there, since it would be a point mass in the given distributional equation. Moreover, $\int |F_n|\leq C a_n$ by the integrability of $\omega_{\kappa,\zeta}$. For fixed large $M$, the mass exhaustion hypothesis then gives
\[
 \mathop{\rm ess\,sup}_{M\leq d\leq L_n}|d\bar e_n'(d)|
       \leq\delta_M+o_n(1),\qquad \delta_M\longrightarrow0.
\]
Integrating from $M$ to $d$ and using local convergence on $D_M$ and the uniform circular oscillation bound shows that
\[
 \limsup_{n\to\infty}\|e_n/\ell\|_\infty
       \leq C\delta_M+\frac{C}{\log(2+M)}.
\]
Consequently $E_n=\|e_n/\ell\|_\infty=o(1)$; the argument requires only local convergence and the preceding bounds.

\smallskip\noindent\textit{Step 2. The scalar equation.}\quad The equation for $e_n$ is
\[
 (-\Delta-2m_q)e_n=F_n+2m_q(e^{e_n}-1-e_n).
\]
Taylor's formula gives
\[
 |e^{e_n}-1-e_n|
 \leq\tfrac12E_n^2\ell^2(2+|y|)^{E_n}.
\]
The quotient $m_q/\omega_{\kappa,\zeta}$ is bounded at zero by $C|y|^{2(q-\kappa)}$, and at infinity by $C|y|^{-2q+\zeta}$. For large $n$ one has
\[
 E_n<(2q-\zeta)/2,
\]
so the tail of $\ell^2(2+|y|)^{E_n}m_q/\omega_{\kappa,\zeta}$ is bounded by a constant times $\ell^2|y|^{-(2q-\zeta)/2}$. Near zero, $q\ge\kappa$ and $\ell=O(1)$ give the bound, including at $q=\kappa$. Lemma~\ref{univ:weighted-scalar} therefore yields
\[
 E_n\leq Ca_n+CE_n^2.
\]
For sufficiently large $n$, $CE_n\leq1/2$, which proves the claim. For radius fitting take $\lambda_n=\exp((W_q(0)-w_n(0))/(2q))\to1$ and use
\[
 \widetilde w_n(y)=w_n(\lambda_n y)+2q\log\lambda_n
\]
on $D_{L_n/\lambda_n}$. Its forcing is $\lambda_n^2F_n(\lambda_n y)$, of the same weighted size because the weight is comparable at radii $r$ and $\lambda_n r$. The radial models are constant on the boundary circles, so the boundary oscillation is unchanged. Its actual density integral on $D_d$ is the old one on $D_{\lambda_n d}$.

The model integrals at these two radii differ by at most $C|\log\lambda_n|$, uniformly in $d$, by the bounded logarithmic derivative of the model mass. Thus local convergence, bounded circular oscillations and mass exhaustion are preserved, and
\[
 \widetilde w_n(0)=W_q(0).
\]
\end{proof}
\subsection{Exponential estimates and comparison of local degrees}

\begin{lemma}[Exponential estimates]
\label{fnd:exponential}
On a fixed surface, let $u$ have zero mean and let $E=\|\nabla u\|_2$. For $q\geq2$,
\[
 \|u\|_q\leq C\sqrt q\,E,\quad
 \|e^{|u|}\|_q\leq C e^{CqE^2},\quad
 \|e^u-1\|_q\leq C\sqrt q\,E e^{CqE^2}.
\]
The bounds also hold for $1\leq q<2$ after a fixed change of constants.
\end{lemma}
\begin{proof}
Integration by parts in the Green representation gives
\[
 u(x)=\int_M\nabla_yG(x,y)\cdot\nabla u(y)\,dV_g(y).
\]
On finitely many coordinate disks the kernel is bounded by $C|x-y|^{-1}$ near the diagonal. It is smooth and bounded on separated supports. In Young's inequality take
\[
 t=\frac{2q}{q+2},\qquad 1+\frac1q=\frac12+\frac1t.
\]
The $L^t$ norm of the local kernel is at most $C(2-t)^{-1/t}\le C\sqrt q$, since
\[
 2-t=4/(q+2),\qquad q^{1/t-1/2}=q^{1/q}
\]
is bounded. Approximation gives the first estimate for $u\in H^1$. Suppose $E>0$. Apply the first estimate with $q=2j$ to the expansion of $e^{cu^2/E^2}$. After integration its $j$-th term is bounded by $(2cC^2)^jj^j/j!$. Since
\[
 j!\ge(j/e)^j,
\]
the integral is uniformly bounded for fixed $c>0$ sufficiently small. Now $q|u|\le cu^2/E^2+Cq^2E^2$ gives the second bound after taking the $q$-th root. For the last bound apply H\"older's inequality with exponent $2q$ to
\[
 |e^u-1|\le |u|e^{|u|}.
\]
If $E=0$, then $u=0$ by the mean-zero condition.
\end{proof}

\begin{lemma}[Comparison of local degrees]
\label{lem:local-degree-spaces}
Let \(X_0=\{u\in C^0(M)^2:\int_M u_i\,dV_g=0,\ i=1,2\}\).
Suppose that a compact set of zeros of \(I-\mathcal T_\rho\) is isolated in both \(X_0\) and \(\Hh\). Its local Leray--Schauder degrees in these spaces agree.
\end{lemma}
\begin{proof}
To compare the two degrees, choose
\[
 h_i\in L^p,\qquad 1<q<p,\qquad b<\infty,\qquad
 \frac1q=\frac1p+\frac1b .
\]
On each bounded mean-zero $H^1$ ball Lemma~\ref{fnd:exponential} bounds $e^{u_i}$ in $L^{2b}$. Jensen's inequality and $\log h_i\in L^1$ bound the denominators away from zero.

Thus $f_i$ is bounded in $L^q$. The same estimate holds on bounded $C^0$ balls. From either space $\mathcal T_\rho$ is continuous into $W^{2,q}$. For $H^1$ this follows from
\[
 |e^u-e^v|\le |u-v|(e^u+e^v),
\]
H\"older's inequality and $H^1\hookrightarrow L^{2b}$.

Give
\[
 X=C^0\cap H^1
\]
the sum norm. Since $W^{2,q}\hookrightarrow X$ compactly in dimension two, $\mathcal T_\rho$ maps bounded sets of either space into relatively compact subsets of $X$. In each space factor $\mathcal T_\rho$ through $X$ by inclusion. Let a compact zero set be isolated in both spaces.

On its isolating neighborhoods smooth the two compact images and truncate in a smooth Laplace eigenbasis. This approximates both images uniformly in $X$ by maps with the same finite-dimensional range.

For sufficiently close approximations, compactness in $X$ keeps the zeros of the straight-line homotopies near the original zero set.

By excision and restriction to the common finite-dimensional range, both local degrees equal the degree of the same map with the same orientation.
\end{proof}

Differentiation of \(J_\rho\) with respect to the inner product on \(\Hh\) gives
\begin{equation}\label{eq:toda-gradient-identity}
 \nabla_{\Hh}J_\rho=I-\mathcal T_\rho.
\end{equation}

\section{Proofs of the local degree formulas}

\begin{proof}[Proof of Proposition~\ref{def:complete-transfer}]
\noindent\textit{Step 1. Local degree and blow-up at the critical parameter.} If
\[
 P=Q=\varnothing,
\]
the disks about the regular points stay away from the fixed zeros and $\gamma>0$. Theorem~\ref{regzero:transfer} applies by the strict trace inequalities at mixed points and, if $V=0$, by
\[
 L_0\ne0.
\]
For $V>0$, the local degree is $(-1)^Ud_S(Z)$ for $\mu>0$, and the local domain is empty for $\mu<0$. For $V=0$ write
\[
 t=t_\mu\tau,\qquad |\mu|=|L_0|t_\mu^2\log(1/t_\mu),
 \qquad \tfrac12\le\tau\le2.
\]
On the side where solutions exist, the normalized scale equation converges uniformly to $\operatorname{sgn}(L_0)(1-\tau^2)$. The endpoint signs give degree $-\operatorname{sgn}(L_0)$.

The index on the complement is $m_S+U-1$, and the degree of the center equations is $\operatorname{sgn}\det DF$. Thus the local degree is $(-1)^U\operatorname{sgn}(L_0)d_S(Z)$ on this side and zero on the other. Taking the degree for $\mu>0$ minus that for $\mu<0$ gives \eqref{eq:local-transfer} for either sign of $L_0$.

If
\[
 P\cup Q\ne\varnothing,
\]
the prescribed source strengths are positive and noninteger. The strict trace inequalities and the separation of the scalar branch give the other assumptions of Theorem~\ref{hyb:transfer}. Its local domain gives \eqref{eq:local-transfer} for $\mu>0$ and is empty for $\mu<0$.

The classification and Lemma~\ref{reg:shadow-limit} give the limiting configuration and shadow solution for each unbounded sequence. The description of these solutions and their scale equation also hold at
\[
 \mu=0.
\]
Since the leading coefficient is nonzero, an unbounded sequence at $\mu=0$ is impossible.

\noindent\textit{Step 2. Separation of the local domains.} Shrink the closed neighborhoods of the centers for the finitely many shadow solutions.

The unordered centers with their types, the choices at the singular points, and the continuous mean-zero scalar branches then have pairwise disjoint compact images.

At the center of each neighborhood these data determine its shadow solution. They also distinguish different configurations, since the pair of local masses at each blow-up point determines its type. If two local neighborhoods intersected along a sequence
\[
 \mu\to0,
\]
their approximate solutions would differ by \(o(1)\) in energy. After taking convergent center and scalar branch parameters, the approximate solutions converge away from all the limiting blow-up points to their exterior Green fields plus their scalar remainders.

The exterior limits then agree there and hence as distributions.

Applying \(-\Delta_g K^{-1}\) recovers the same two local masses at each blow-up point, since the scalar remainder gives zero mass to points. Subtracting the Green potentials of the point masses gives the same mean-zero scalar remainder.

This contradicts the disjointness of the compact images, even if the limiting parameters are on their boundaries.

The constants can be chosen uniformly over the finite family. In a nonempty local domain, the largest outer scale in the closed parameter interval tends to zero with $\mu$.

\noindent\textit{Step 3. Divergence of the norms.} Take a fixed small circle about one blow-up point and a concentric circle of radius $R_a^{1/2}$. The circular means of the approximate first component differ by $(2s_a-t_a)\log(1/R_a)+O(1)$.

Indeed, the outer logarithmic term has coefficient of $\log r$ equal to $-2(2s_a-t_a)$ and the corrections are bounded on both circles. The same formula applies at an ordinary point with
\[
 (s_a,t_a)=(1,0).
\]
Cauchy--Schwarz in the logarithmic radial variable bounds the Dirichlet energy below by $c\log(1/R_a)$. At a zero the correction on the complement is $o(1)$ in $H^1$.

Thus the $H^1$ norms tend to infinity. So do the $C^0$ norms, since normalization and the Green estimate give an $H^1$ bound for any $C^0$-bounded family of solutions.
\end{proof}
\begin{proof}[Proof of Proposition~\ref{univ:wall-jump}]
\noindent\textit{Step 1. Configurations and choice of weights.} The rational noninteger choice of $\rho_2/(4\pi)$ excludes complete concentration of the second component by \eqref{fnd:independence} and the generic mass list.

In an unbounded sequence the first component therefore concentrates completely. The integrability of the residual density leaves the local pairs, in $4\pi$ units,
\[
       (1,0),\quad(2,1),\quad(A_p,0),\quad(A_p+B_p,B_p).
\]
For a fixed critical value put $n_0=U+2V$. Subtracting two representations of this value gives a rational relation in \eqref{fnd:independence}, so the data $n_0,P,Q$ are determined by this critical value. The configurations to consider are precisely
\[
 0\le V\le\lfloor n_0/2\rfloor,\qquad U=n_0-2V,\qquad
 V+\sum_{p\in Q}B_p<\rho_2/(4\pi).
\]
The last inequality gives positive residual scalar mass. Lemma~\ref{reg:shadow-limit} identifies every limiting configuration with a solution of the corresponding shadow system, including configurations with no prescribed bubbles.

If a pair replacement has positive residual mass, it increments $V$ by one within this list and preserves the first mass.

Every scalar exclusion \eqref{fnd:scalar-exclusion} follows from rational independence, since moving the sum of the second local masses to its right side would express $\rho_2/(4\pi)$ as an integer plus an integral linear combination of the shifted parameters.

Proposition~\ref{univ:shadow} therefore applies to the whole finite family. It supplies the finite nondegenerate solution sets, invertible scalar blocks, and strict trace inequalities at regular mixed points. At each prescribed point, $A_p,B_p,A_p+B_p\notin\mathbb Z$.

For configurations where ordinary points determine the leading scale, we choose $L_0\ne0$ before constructing the solutions. Write
\[
 L_0=2\pi\sum_i d_i^2\kappa_i
\]
with the traces computed at the full critical mass $M$. Let $x_j$ run over the distinct free centers of the finitely many shadow solutions, including those obtained after a pair replacement.

Perturb $\log h_1^*$ in disjoint disks about the $x_j$ away from the prescribed points. Use smooth cut-offs of quadratic functions vanishing with their first derivatives at $x_j$.

At each shadow solution these perturbations preserve the scalar equation, the balance equations at the centers and the coefficients $d_i$, also for configurations with prescribed blow-up points. If $c_j$ is the Laplacian of the perturbation at $x_j$, its coefficient in $L_0$ is
\[
 2\pi\sum_{p_i=x_j}d_i^2>0
\]
whenever the solution uses $x_j$. The conditions $L_0=0$ therefore give finitely many proper affine hyperplanes. Choose arbitrarily small $c_j$ outside their union.

Nondegeneracy and the strict trace inequalities at regular mixed points persist by openness, and the compactness neighborhoods of Proposition~\ref{univ:shadow} exclude additional solutions. This choice applies also when all bubbles are ordinary of Lemma~\ref{univ:regular-only}.

If $P\cup Q=\varnothing$, Lemma~\ref{univ:regular-only} applies. We have the strict inequalities at regular mixed points and $L_0\ne0$ when all points are ordinary. In both cases the local degree for $\mu>0$ minus that for $\mu<0$ is $(-1)^Ud_S(Z)$. This includes
\[
 S=\varnothing.
\]
For $P\cup Q\ne\varnothing$ choose a prescribed point as the reference point for the radial parameter.

We distinguish three cases according to the largest outer profile exponent $s_*$.

\noindent\textit{Step 2. The case $s_*>1$.} If the largest outer profile exponent $s_*$ exceeds one, put
\[
 \mathcal L=e^{-T/s_*}.
\]
The square of the residual, including fixed powers of $T$, is $o(\mathcal L)$. Differentiating the energy after the complementary correction consequently changes its scale derivative by
\[
 o(\mathcal L)+O(\mu^2T^C),
\]
as in \eqref{hyb:energy-comparison}. Lemma~\ref{univ:local-scale-moments} replaces the mixed singular calculation in that identity. At a pure point with
\[
 q=A>1,\qquad B>0,
\]
the residual density is $g=O(r^{2(q+B)-2})$. The bounded local right-hand side first gives a $C^1$ scalar remainder. Since
\[
 2(q+B)-2>0,
\]
the density is then H\"older and the first corrected background admits a $C^{2,\eta}$ Taylor expansion for some
\[
 \eta>0.
\]
In \eqref{hyb:pinned-radial-moment} the first-order term is canceled by \eqref{hyb:firstjet}. The terms linear in that term have zero angular integral. The Taylor remainders have moment $O(R^{2+d})$ for
\[
 0<d<\min(\eta,1,2q-2).
\]
The weighted quadratic identities in Lemma~\ref{hyb:trial-energy} cancel the angular-square term and give the remaining coefficient $\kappa J_q R^2/(4q)$. For the second density, the radial overlap is $O(R^{2q}r^{2B-2})$, which is $O(R^{2q})$ in $L^1$ since $B>0$. For the first-order correction use
\[
 |P|\le Cr/(1+(r/R)^{2q}).
\]
Its product with $g$ integrates on $r\le R$ to $O(R^{2(q+B)+1})$, and on $R<r<d_0$ to at most $CR^{2q}\int_R^{d_0}r^{2B}\,dr$. Its overlap with the second radial correction costs at most the already bounded radial overlap.

Multiplication by the bounded second scale tangent preserves this $O(e^{-T})$ estimate. Every prescribed outer profile exponent below one contributes $O(e^{-T})$ by the residual and moment estimates for profile exponents below one. An ordinary bubble contributes $O(Te^{-T})$ by its truncated quadratic moment.

At a matched mixed singular point the trace is well defined even when the original residual density is unbounded. Indeed, if $g$ is that density and $f$ its matched exact local limit, then
\[
 \Delta J_1=\kappa_1+g-f,
 \qquad g-f=f(e^{J_2}-1)=O(|z|^{2A}).
\]
Thus $\Delta J_1(0)=\kappa_1(0)$. At a pure prescribed point of outer profile exponent greater than one the residual density vanishes, so that the same trace follows directly from the exterior equation. Indeed, away from the blow-up points, $\Delta_g u_1^0=2M-\rho_2+g$ and
\[
 \Delta_g\log h_1=\Delta_g\log h_1^*-4\pi\sum_{q\in S}\alpha_{q,1}.
\]
Removal of the prescribed logarithm contributes zero off the source, and the canonical area factor in an isothermal coordinate with $\varphi(0)=0$ contributes $-2K_g(p)$. Hence the trace in both cases is
\[
 \kappa_1=\Delta_g\log h_1^*(p)-2K_g(p)
         -4\pi\sum_{q\in S}\alpha_{q,1}+2M-\rho_2.
\]
At this critical value choose the additional trace perturbations of Corollary~\ref{fnd:extra-traces} simultaneously at all prescribed points where they are needed. Their proof uses the extended source exclusion and the compactness in Proposition~\ref{univ:shadow}.

The perturbations are supported away from the finitely many free centers, so the ordinary leading coefficients already chosen and the strict trace inequalities at the regular mixed points persist.

This makes a common choice of weights for all configurations at this critical value, including when some configurations have an ordinary leading scale and others have largest outer profile exponent greater than one.

Build the local models and use Proposition~\ref{univ:complete-capture} at these same weights to put concentrating true solutions in the family of solutions of the projected equation, before using the scale equation to locate its zeros.

For this finite family of shadow solutions, choose the exponent in the construction smaller than the positive exponents in both the approximation of the solutions and the estimate for the residual.

For the regular center and shape equations, apply Section~\ref{hyb:section}. The contours avoid the prescribed points, so contributions from the other points have uniformly smooth derivatives there. For \eqref{hyb:measure-stability} the normalized approximate densities satisfy
\[
 \|\nu\|_{p_0}\le e^{CT},
\]
by the pointwise profiles and the fixed exponent $p_0>1$. Let $p_0'=p_0/(p_0-1)$. By H\"older's inequality,
\[
 \|e^\eta-1\|_{L^q(\nu)}
 \le \|\nu\|_{p_0}^{1/q}\|e^\eta-1\|_{L^{qp_0'}} .
\]
For $q=T$ the first factor is bounded. The exponential estimate in \eqref{hyb:measure-stability} bounds the right-hand side by
\[
 C\sqrt T\|\eta\|_{H^1}e^{CT\|\eta\|_{H^1}^2}=o(1).
\]
The normalization factor differs from one by at most this error. The matched Taylor coefficients at the prescribed points are $C^1$ in the shadow parameters on the compact family.

The equations for the shape parameters therefore converge uniformly, as required for the degree. For $\rho_1=M+\mu$ the additional terms in \eqref{hyb:exact-trial-variation} are
\[
 -\mu-\mu M_T'/M_T+\mu\int_M\partial_T\widehat B_{T,1},\qquad
 M_T=\int_M\mathfrak n_1 .
\]
The local mass estimates give $M_T'=O(\mathcal L)$. Differentiating the local logarithmic terms and cut-offs gives
\[
 \int_M|\partial_T\widehat B_{T,1}|=O(\mathcal L T^C).
\]
At a zero of the center equations subtract the shape contribution as in \eqref{hyb:actual-radial}. The scale equation has left-hand side
\[
 -\mu+\mathcal L\Gamma+o(\mathcal L)
       +O(\mu^2T^C+|\mu|\mathcal L T^C).
\]
For each point of maximal outer profile exponent the coefficient is positive, since $\kappa_aJ_q>0$ at a prescribed point and $\kappa_b\mathcal H(B_b)>0$ at a regular mixed point. Thus $0<c\le\Gamma\le C$. At a zero,
\[
 |\mu|\le C\mathcal L+C(|\mu|+\mathcal L)T^C|\mu|.
\]
In the preliminary exponential range $C(|\mu|+\mathcal L)T^C=o(1)$. Absorption gives
\[
 |\mu|=O(\mathcal L).
\]
Substituting in the same equation makes both parameter errors $o(\mathcal L)$, so $\mu=\mathcal L(\Gamma+o(1))>0$. For $\mu>0$ fix small $a>0$ and large $b>0$ independently of $\mu$. Choose them so that
\[
 \mathcal L/\mu=a
\]
makes the scale component strictly negative, and $\mathcal L/\mu=b$ makes it strictly positive where the balance equations vanish throughout the shape balls. The bounds for its coefficient are uniform on that set. Choose the exponent in Proposition~\ref{univ:projected-construction} also with
\[
 \delta<1/s_*.
\]
On the closed scale interval,
\[
 \mu e^{\delta T}\leq a^{-\delta s_*}\mu^{1-\delta s_*}\to0,
\]
so the reduction and its error bounds hold on the whole interval. On the closed product, normalize each two-dimensional center component by $-1/(3w_{\rm ctr})$ and each shape component by $-16\pi/(3R_b^2)$. The factors preserve orientation in real dimension two.

By Proposition~\ref{hyb:stress-rows} the center component is uniformly close to the isolated center field on its closed ball. At a zero of the center equations, the shape component is $\Phi_b(B_b)+o(1)$. The gap $\kappa_b>4|T_b|$ and
\[
 \nabla\mathcal M(B)\to\pi B/|B|
\]
uniformly in the direction of large $B$ make this component point outward on a suitably large fixed shape sphere. The normalized scale component has the signs above at every zero of the center equations. Lemma~\ref{hyb:homotopy} therefore gives scale degree $-1$ and shape degree one.

In the local coordinates use a fiber radius $\tfrac12e^{-\delta T}T^{C_0}$ about the graph, strictly above the graph's correction bound for the chosen $\delta$ and sufficiently large $T$.

The graph and the fibers remain inside the original ball of uniqueness for the projected equation. The $O(T)$ inverse on the complement and the fixed $L^{p_0}$ exponential probability bound give
\[
 \|A^{-1}N(v)\|\leq CT^{5/2}\|v\|^2
\]
along the fiber homotopy of Theorem~\ref{hyb:transfer}. Since
\[
 CT^{5/2}\|v\|=o(1)
\]
on the entire fiber ball, each zero during the homotopy lies on the graph, where the preceding boundary estimates apply. The widths within each closed scale interval are comparable. We now prove that the coordinates are unique.

For the prescribed reference point $p$, with fixed profile exponent $s>0$, the first inverse-Cartan component has gradient equal to that of $-\log(R^{2s}+|z|^{2s})$ plus a uniformly bounded $C^1$ field and an $o(1)$ term in $L^2$. On a mixed disk,
\[
 (2\log n_{1,R}+\log n_{2,R})/3=-\log F_R+(2/3)\log|W_R|+\log4.
\]
The logarithms of the prescribed weights cancel the contribution of the power of $z$ in $W_R$. The remaining term $\log|W_0|$ and the matched background are smooth or $C^2$. Set $b=\min(1,2A)>0$. To replace $F_R=\tau F_0+|z|^{2s}$ by $\tau+|z|^{2s}$, use
\[
 F_0-1=O(r^b),\qquad \nabla F_0=O(r^{b-1}).
\]
The gradient of the difference of the logarithms is $O(r^{b-1}\min\{1,\tau r^{-2s}\})$. Its squared $L^2$ norm satisfies
\[
 \left\|\nabla\log\frac{F_R}{\tau+|z|^{2s}}\right\|_2^2
 \le CR^{2b}+C\tau^2\int_R^d r^{2b-4s-1}\,dr=O(R^{2b}),
 \qquad b<2s .
\]
At a pure point the additional correction is proportional to $P$ in Lemma~\ref{univ:pure-residual}. Its gradient tends to zero away from the origin and is bounded by
\[
 C/(1+(r/R)^{2s}).
\]
It therefore tends to zero in $L^2$. These estimates also apply on the fixed cut-off annulus. Subtracting the means preserves the gradients.

If two approximations differ by $o(1)$ in $H^1$, take a subsequential limit $\lambda\in(0,\infty)$ of their reference width ratio and rescale by the first width. Equality of the limiting gradients on a fixed rescaled annulus gives
\[
 \frac{|y|^{2s-1}}{1+|y|^{2s}}
   =\frac{|y|^{2s-1}}{\lambda^{2s}+|y|^{2s}}.
\]
Hence the prescribed width ratio tends to one. For each free point rescale the active inverse-Cartan gradient by the width of the ordinary bubble or the inner width of the mixed bubble.

All fixed source terms on its disk have smooth remote gradients, and the fitted local errors and projected corrections are $o(1)$ in $H^1$.

In the mixed case the second inverse-Cartan potential is $-\log Q_R$ up to such terms, and $Q_R/R^4$ is a constant multiple of $a^2+|z|^2$ to leading order on the inner scale by \eqref{reg:gram}. In either comparison the scaled gradient is that of a translated radial logarithmic profile.

An unbounded center shift in scaled units would give zero gradient on every fixed rescaled compact set, by the bound $C/(a+|z|)$ for $\nabla\log Q_R$ in the mixed case (with the analogous ordinary width bound).

A bounded shift is fixed by the unique zero of the limiting gradient, and its derivative there fixes the positive width ratio. The now aligned mixed center shift is
\[
 o(a)=o(R).
\]
On its outer scale the first inverse-Cartan potential is $-\log F_R$. If the outer width ratio tends to
\[
 \lambda_0>0,
\]
equality of the gradients makes $1+|y^2-B|^2$ and $\lambda_0^4+|y^2-\lambda_0^2B'|^2$ proportional. Their quartic coefficients give proportionality one, their quadratic coefficients give
\[
 B=\lambda_0^2B',
\]
and their constants give $\lambda_0=1$. With the centers and shapes now matching in the limit, the continuous positive factor $d_0$ in the reference law
\[
 R^{2s}=d_0^s e^{-T}
\]
matches as well. Thus the two $T$ values differ by $o(1)$. The scaled coordinate differences are all $o(1)$; the uniformly bounded inverse of the matrix of pairings with the parameter derivatives and the small remainders allow Lemma~\ref{journal:modulation-uniqueness} on the segment between them.

\noindent\textit{Step 3. The case $s_*=1$.} If $s_*=1$, genericity implies that some ordinary regular bubbles are present, no mixed regular bubbles are present, and all prescribed outer profile exponents are below one.

The projected and reduced equations are those of Lemma~\ref{univ:ordinary-dominated}. To identify the true solutions use the preliminary estimates from Proposition~\ref{univ:complete-capture}. With
\[
 t^2=e^{-T},\qquad L=\log(1/t),
\]
the scale equation has
\[
                 \mu=L_0t^2L+O(t^2).
\]
The prescribed contributions are $O(t^2)$. The refined estimates in Lemma~\ref{univ:ordinary-dominated} place the centers strictly in the shrinking domains and put the correction strictly within the neighborhood of the solution of the projected equation. To check both signs of $L_0$, fix $\mu L_0>0$ and write
\[
 t=t_\mu\omega,
\]
where $t_\mu$ is the small root of $|L_0|t_\mu^2\log(1/t_\mu)=|\mu|$. On
\[
 1/2\leq\omega\leq2
\]
the scale component divided by $|\mu|$ is $\operatorname{sgn}(L_0)(\omega^2-1)+o(1)$ uniformly, including the whole center ball in this refined scale. In increasing $T$ the two ends have signs $\operatorname{sgn}(L_0)$ and $-\operatorname{sgn}(L_0)$, with limiting magnitudes $3$ and $3/4$.

Its one-dimensional degree is $-\operatorname{sgn}(L_0)$. The center boundary stays zero-free and the complementary space and the parametrization are those of Lemma~\ref{univ:ordinary-dominated}. The mass law excludes blow-up for the opposite sign of $\mu$ and for $\mu=0$.

\noindent\textit{Step 4. The case $s_*<1$.} If every outer profile exponent is below one, there are no regular centers. Proposition~\ref{univ:projected-construction}, with the strict estimates for outer profile exponents below one, applies on each fixed compact interval of exterior first-weight perturbations.

The shadow systems, matching data and local patches are unchanged by these perturbations. Proposition~\ref{universal:all-small-wall} then selects the weight at which
\[
 ce^{-T}\leq S(T,0)\leq Ce^{-T},
\]
with no zeros for $\mu\leq0$ and scale degree $-1$. Apply Proposition~\ref{univ:complete-capture} at this final fixed weight. Every blowing-up sequence then belongs to the neighborhood on which this local degree was computed. The two noncritical strip degrees can be compared at this choice by Lemma~\ref{universal:wallwise-weights}.

\noindent\textit{Step 5. Local degrees and separation of the neighborhoods.} For a configuration with prescribed bubbles, let
\[
 n=U+V+|P|+|Q|,\qquad v=V+|Q|.
\]
The actual complementary Morse index supplied by Proposition~\ref{univ:projected-construction} and Proposition~\ref{signed:model-inverse} is
\[
 \nu=m_S+U+2V+|P|+2|Q|-1
      +2\sum_a(\lceil s_a\rceil-1)
      +2\sum_{t_b>0}(\lceil t_b\rceil-1).
\]
The parametrization by the chosen sources and parameters preserves orientation. Each shape block, when present, has degree one. Write the shadow local degree as
\[
 d_S(Z)=(-1)^{m_S}\deg_{\mathrm B}F_0,
\]
with center degree one if there are no free centers. For $s_*>1$ or $s_*<1$ the scale degree is $-1$. Thus the local degree for $\mu>0$ is $(-1)^{U+|P|}d_S(Z)$ and the local contribution for $\mu<0$ is zero.

For $s_*=1$ the local degree for $\mu L_0>0$ is $(-1)^{U+|P|}\operatorname{sgn}(L_0)d_S(Z)$, and the other side is empty. Taking the positive side minus the negative side gives the displayed jump in all three cases, including when $L_0<0$. The regular-only case was computed above.

For every selected shadow solution, take its closed center neighborhood sufficiently small that the compact image of its unordered centers with their types, prescribed states, and continuous mean-zero scalar branch is disjoint from the corresponding data for the other solutions.

The data distinguish the finitely many roots, so this is possible. The ordinary shrinking product is taken inside such a neighborhood. The scale interval is closed in each degree calculation. With
\[
 u=B_p+\eta_p+v_p,
\]
the equation $L_p(u-B_p)=0$ solves locally for $p$ by invertibility of the derivative with respect to $p$. The resulting local inverse, together with injectivity of the parametrization, gives a bounded open image in the energy space of the open parameter product and open fiber ball.

If a sequence in its image converges, a subsequence of its parameter values converges in the closed product.

Continuity of $B_p,\eta_p$ and of the subspace projections then gives a strong limit in the closed fiber ball. Thus a boundary zero must lie over a base face or on a fiber sphere, all of which were excluded.

If neighborhoods of distinct shadow solutions intersected along parameters tending to the critical value, their approximations would differ by $o(1)$ in $H^1$ and hence in $L^2$, by the mean-zero condition.

Pass to a subsequence of the center and scalar branch parameters in their closed domains.

Away from all the limiting blow-up points the corrections in the local construction vanish locally, so each approximation converges to its exterior Green field there.

The limits agree in local $L^2$. The logarithmic Green fields are locally integrable and there are only finitely many blow-up points, so the two limits agree as distributions.

Applying $-\Delta_gK^{-1}$ identifies the sums of the point masses in both components, since the extra term $(0,\gamma He^w/I-\gamma)^t$ has no point masses. Each type is recovered from its position and pair of charges, with $A_p,B_p$ fixed at each mark.

Subtracting the Green potentials of the point masses and using the mean normalization then identifies the same scalar remainder. This contradicts the disjoint closed shadow data even if the limiting parameters lie on their boundaries.

Equal-type label changes simultaneously permute two-dimensional blocks of parameters and equations. They preserve orientation, and one labeling is used for each unordered shadow datum.

For configurations with only regular blow-up points, the exterior fields distinguish the neighborhoods in the same way. For these weights every unbounded sequence therefore lies in the indicated neighborhoods for the corresponding sign of $\mu$.
\end{proof}
\begin{proof}[Proof of Proposition~\ref{fnd:scalar-degree-internal}]
\noindent\textit{Step 1. Reduction to one concentrating component.} Lemma~\ref{fnd:scalar-parameter} gives scalar compactness whenever the displayed parameter exclusions hold.

In particular for any scalar parameter below the first positive mass, the homotopy to zero has degree one, since the zero-parameter map is the identity.

First suppose that $1,(a_p)_{p\in S}$ are rationally independent. Use $A_p=a_p$ in an auxiliary Toda system. Choose
\[
 B_p>\gamma/(4\pi)+1
\]
outside the countably many proper affine hyperplanes that would violate \eqref{fnd:independence}. Choose $\rho_2/(4\pi)\in(0,1)\cap\mathbb Q$ and keep it fixed along the auxiliary first-parameter path.

Any positive second local mass at a regular point is at least $4\pi$, and the nonzero prescribed second masses are $4\pi B_p$ and $4\pi(A_p+B_p)$.

All exceed the total mass of the second component. Thus the limiting second measure has no point masses, and the only blow-up configurations have
\[
 V=0,
\]
\[
 Q=\varnothing.
\]
At a positive critical value of $\rho_1$, $4\pi(U+\sum_{p\in P}a_p)$ determines $(U,P)$ by rational independence.

The residual second parameter at each critical value of the auxiliary system is $\rho_2$. For any fixed separated configuration, its scalar weight has shifts $B_p$ or $A_p+B_p$ at the marks and shift two at each ordinary center. The marks and centers remain distinct.

All these shifts exceed one. The scalar homotopy of the coefficient from $2\rho_2$ to zero therefore meets no positive critical parameter of Lemma~\ref{fnd:scalar-parameter}. The fiber degree is one by the identity at zero. This calculation uses scalar compactness only.

Apply Proposition~\ref{univ:wall-jump} at this choice of $\rho_2$. A pair replacement would require a second local mass
\[
 4\pi>\rho_2,
\]
so the retained family is closed under all admissible replacements. If $\rho_2/(4\pi)$ belonged to the right-hand side of \eqref{fnd:scalar-exclusion} in normalized units, its rationality and \eqref{fnd:independence} would force all shift coefficients to vanish, leaving an integer. Thus the stronger shadow exclusions hold.

\noindent\textit{Step 2. The degree jump and its sum.} At each shadow solution the local jump multiplies its scalar Morse sign and reduced center sign by $(-1)^{U+|P|}$. The formulation of Proposition~\ref{fnd:shadow-degree} in Proposition~\ref{univ:shadow} multiplies the unordered center degree by the fixed-center scalar degree.

Here that last factor is the degree one just proved. Since $V=0$, the center degree is $\binom{\chi(M)-|S|}{U}$, with the equal-type labels already quotiented out. With the excision below, the whole jump is
\[
          (-1)^{U+|P|}\binom{\chi(M)-|S|}{U}.
\]
On either side of each critical value take endpoints avoiding $\Gamma_1$, which below $\gamma$ consists exactly of the critical values of the scalar equation since
\[
 B_p>\gamma/(4\pi)+1.
\]
The fixed second parameter avoids $\Gamma_2$. By Lemma~\ref{fnd:sharp-compactness}, log-linear homotopies of the positive smooth factors at each endpoint keep its degree unchanged.

Thus the choices made at each critical value in Lemma~\ref{universal:wallwise-weights} can all be compared with one fixed original choice, including the first weight $h$.

When $P=\varnothing$, use the description of the regular bubbling solutions and the nonzero coefficient of the ordinary scale equation. When
\[
 P\ne\varnothing,
\]
use the description of all bubbling solutions and the three scale equations in Proposition~\ref{univ:wall-jump}. In each case the nonzero scale law rules out an escaping sequence at the critical value.

With the selected weights fixed, take an energy ball containing all solutions at the critical parameter strictly in its interior.

A sequence of boundary zeros with parameters tending to the critical value would converge to a boundary zero there by compactness of the normalized Green map. The disjoint closed neighborhoods supplied in the calculation of the degree jump contain all escaping nearby zeros.

At an occupied point the approximate inverse-Cartan potential has the radial term $-\log(R^{2s}+|z|^{2s})$ up to a bounded energy term. On $LR<|z|<d$ its gradient has magnitude at least $c_s/|z|$. The energy thus diverges, also after the $o(1)$ graph correction.

Consequently all these local zeros are outside the fixed ball as the parameter tends to its critical value.

A missed zero outside the ball would have either a bounded subsequence converging to a zero at the critical value outside it, or a subsequence whose norms tend to infinity and whose members eventually lie in one of these neighborhoods.

Both alternatives are impossible. Excising the closed ball from the neighborhoods and applying additivity gives the sum of their local degrees. The two degrees on the ball cancel by homotopy.

Start with positive first parameter below both $\gamma$ and $4\pi\min(\{1\}\cup\{a_p:p\in S\})$. Since $\rho_2<4\pi$ and all
\[
 B_p>1,
\]
Proposition~\ref{prop:signed-subcritical} applies to the simultaneous descent of both parameters to zero. The starting degree is one. The binomial expansion of $(1-y)^{\chi(M)-|S|}$ assigns weight $(-1)^U\binom{\chi(M)-|S|}{U}$ to $y^U$; choosing the factors $-y^{a_p}$ exactly for $p\in P$ gives the remaining sign and exponent.

Adding the jumps at positive critical values strictly below $\gamma$ to the starting degree sums exactly the asserted scalar coefficients. The only exponent-zero contribution is the starting degree one.

\noindent\textit{Step 3. The scalar equation at $\rho_2=0$.} At the target first parameter $\gamma$, decrease $\rho_2$ to zero within its original range.

For positive sequence members the Green small-mass estimates of Proposition~\ref{fnd:signed-compactness} give local bounds away from the finite set of blow-up points, using only upper bounds for the masses and normalization factors.

In a disk containing only one blow-up point subtract the smooth constant forcing, add the normalizing constants and source logarithms as in that proof.

Lemma~\ref{fnd:local-unitary} applies with fixed generic source strengths at a marked point, or $(1,1)$ as the shifted parameters at an unmarked point. Its nonzero second masses all exceed the total second mass bound, so every second local mass is zero. Let $4\pi x>0$ be the first local mass, and set
\[
 a=A_p
\]
at a mark and $a=1$ at a free point. Since the second local mass vanishes, the Pohozaev identity gives
\[
 x^2=ax,\qquad x=a,\qquad 4\pi x=4\pi a.
\]
If the first normalizing factor has a positive subsequential limit, integrability of the residual mass requires $2x\le a$. Indeed, on circles of radius $d$,
\[
 \overline u_1(d)=-4x\log d+o(|\log d|).
\]
For the residual measure use the circular mean of the logarithmic kernel, $\log\max\{d,|z|\}$: after division by $|\log d|$ the contribution inside a fixed small disk is bounded by its mass, and the contribution outside tends to zero. Jensen's inequality on circles now gives
\[
 2x\le a,
\]
contradicting $x=a>0$. The first normalization factor would therefore tend to zero, so that the convergence away from the blow-up points leaves no residual mass in the first component. The entire first parameter would then be a scalar critical value, a contradiction.

The argument also applies along sequences with $\rho_2>0$ and $\rho_2\to0$, since the local classification only uses uniform mass, forcing and boundary-oscillation bounds after the source logarithms are removed. At
\[
 \rho_2=0,
\]
Lemma~\ref{fnd:scalar-parameter} supplies the bound. On bounded sets the normalized Green maps are jointly compact and continuous up to this endpoint by the uniform $L^p$ bound on the weights for a fixed $p>1$. At zero set
\[
 g=2\gamma h e^{u_1}/\int_M h e^{u_1}.
\]
The residual is $(u_1-Gg,u_2+Gg/2)$. Conjugating domain and range by
\[
 (u_1,u_2)\mapsto(u_1,u_2+u_1/2)
\]
gives the scalar residual in the first variable and the identity in the second. This identifies the computed auxiliary degree with the scalar degree.

\noindent\textit{Step 4. Continuation in the shifted parameters.} For arbitrary positive shifts at a noncritical parameter, take a small convex neighborhood of the target shifts. The labels $n+\sum_{p\in P}a_p$ below any bounded cut-off have a uniformly bounded nonnegative integer $n$.

At the target none of their positive values equals the cut-off. By shrinking the neighborhood keep all the finitely many strict tests fixed, and join the target to rationally independent shifts in it.

Lemma~\ref{fnd:scalar-parameter} gives a common bound along the homotopy of the shifted parameters, also for shifts of two. Hence the scalar degree equals the labelled coefficient sum at the target. In the shadow system, the added ordinary zeros stay distinct and away from the marks.
\end{proof}


\begin{thebibliography}{99}

\bibitem{BJLY2019}
D. Bartolucci, A. Jevnikar, Y. Lee and W. Yang,
Uniqueness of bubbling solutions of mean field equations,
\emph{J. Math. Pures Appl.} (9) \textbf{123} (2019), 78--126.
\href{https://doi.org/10.1016/j.matpur.2018.12.002}
{doi:10.1016/j.matpur.2018.12.002}.

\bibitem{BT2017}
D. Bartolucci and G. Tarantello,
Asymptotic blow-up analysis for singular Liouville type equations
with applications,
\emph{J. Differential Equations} \textbf{262} (2017), 3887--3931.
\href{https://doi.org/10.1016/j.jde.2016.12.003}{doi:10.1016/j.jde.2016.12.003}.

\bibitem{BYZ2026U}
D. Bartolucci, W. Yang and L. Zhang,
Asymptotic analysis and uniqueness of blowup solutions of
non-quantized singular mean field equations,
\emph{Math. Ann.} \textbf{395} (2026), Paper No.~103.
\href{https://doi.org/10.1007/s00208-026-03536-z}
{doi:10.1007/s00208-026-03536-z}.

\bibitem{BYZ2026N}
D. Bartolucci, W. Yang and L. Zhang,
Non degeneracy of blow-up solutions of non-quantized singular
Liouville-type equations and the convexity of the mean field entropy
of the Onsager vortex model with singular sources,
\emph{Ann. Sc. Norm. Super. Pisa Cl. Sci.} (5),
published online (2026).
\href{https://doi.org/10.2422/2036-2145.202501_030}
{doi:10.2422/2036-2145.202501\_030}.


\bibitem{CL2002}
C.-C. Chen and C.-S. Lin,
Sharp estimates for solutions of multi-bubbles in compact Riemann
surfaces,
\emph{Comm. Pure Appl. Math.} \textbf{55} (2002), 728--771.
\href{https://doi.org/10.1002/cpa.3014}{doi:10.1002/cpa.3014}.


\bibitem{CL2003}
C.-C. Chen and C.-S. Lin,
Topological degree for a mean field equation on Riemann surfaces,
\emph{Comm. Pure Appl. Math.} \textbf{56} (2003), 1667--1727.

\bibitem{CL2010}
C.-C. Chen and C.-S. Lin,
Mean field equations of Liouville type with singular data:
sharper estimates,
\emph{Discrete Contin. Dyn. Syst.} \textbf{28} (2010), 1237--1272.
\href{https://doi.org/10.3934/dcds.2010.28.1237}
{doi:10.3934/dcds.2010.28.1237}.


\bibitem{CL2015}
C.-C. Chen and C.-S. Lin,
Mean field equation of Liouville type with singular data: topological degree,
\emph{Comm. Pure Appl. Math.} \textbf{68} (2015), 887--947.
\href{https://doi.org/10.1002/cpa.21532}{doi:10.1002/cpa.21532}.

\bibitem{DEM2012}
M. del Pino, P. Esposito and M. Musso,
Nondegeneracy of entire solutions of a singular Liouville equation,
\emph{Proc. Amer. Math. Soc.} \textbf{140} (2012), 581--588.
\href{https://doi.org/10.1090/S0002-9939-2011-11134-1}{doi:10.1090/S0002-9939-2011-11134-1}.

\bibitem{LLWYZ2026}
Y. Lee, H. Li, J. Wei, W. Yang and L. Zhang,
The topological degree of the regular $SU(3)$ Toda system,
preprint (2026).
\href{https://arxiv.org/abs/2611.20013}{arXiv:2611.20013}.

\bibitem{LLY2022}
Y. Lee, C.-S. Lin and S. Yan,
Degree counting formula for non-symmetric Toda systems of rank two,
\emph{J. Funct. Anal.} \textbf{283} (2022), Paper No. 109642.
\href{https://doi.org/10.1016/j.jfa.2022.109642}{doi:10.1016/j.jfa.2022.109642}.

\bibitem{LLWY2018}
Y. Lee, C.-S. Lin, J. Wei and W. Yang,
Degree counting and shadow system for Toda system of rank two: one bubbling,
\emph{J. Differential Equations} \textbf{264} (2018), 4343--4401.
\href{https://doi.org/10.1016/j.jde.2017.12.018}{doi:10.1016/j.jde.2017.12.018}.

\bibitem{LLYZ2020}
Y. Lee, C.-S. Lin, W. Yang and L. Zhang,
Degree counting for Toda system with simple singularity: one point blow up,
\emph{J. Differential Equations} \textbf{268} (2020), 2163--2209.
\href{https://doi.org/10.1016/j.jde.2019.09.016}{doi:10.1016/j.jde.2019.09.016}.


\bibitem{Li1999}
Y.~Y. Li,
Harnack type inequality: the method of moving planes,
\emph{Comm. Math. Phys.} \textbf{200} (1999), 421--444.
\href{https://doi.org/10.1007/s002200050536}
{doi:10.1007/s002200050536}.

\bibitem{LWWZ2018}
C.-S. Lin, J. Wei, W. Yang and L. Zhang,
On rank-2 Toda systems with arbitrary singularities: local mass and new estimates,
\emph{Anal. PDE} \textbf{11} (2018), 873--898.
\href{https://doi.org/10.2140/apde.2018.11.873}{doi:10.2140/apde.2018.11.873}.

\bibitem{LWY2012}
C.-S. Lin, J. Wei and D. Ye,
Classification and nondegeneracy of $SU(n+1)$ Toda system with singular sources,
\emph{Invent. Math.} \textbf{190} (2012), 169--207.
\href{https://doi.org/10.1007/s00222-012-0378-3}{doi:10.1007/s00222-012-0378-3}.

\bibitem{LWZ2015}
C.-S. Lin, J. Wei and L. Zhang,
Classification of blowup limits for \(SU(3)\) singular Toda systems,
\emph{Anal. PDE} \textbf{8} (2015), 807--837.
\href{https://doi.org/10.2140/apde.2015.8.807}{doi:10.2140/apde.2015.8.807}.


\bibitem{LWZhao2012}
C.-S. Lin, J. Wei and C. Zhao,
Sharp estimates for fully bubbling solutions of a $SU(3)$ Toda system,
\emph{Geom. Funct. Anal.} \textbf{22} (2012), 1591--1635.
\href{https://doi.org/10.1007/s00039-012-0193-4}
{doi:10.1007/s00039-012-0193-4}.


\bibitem{Tarantello2023}
G. Tarantello,
On the blow-up analysis at collapsing poles for solutions of singular
Liouville-type equations,
\emph{Comm. Partial Differential Equations} \textbf{48} (2023), 150--181.
\href{https://doi.org/10.1080/03605302.2022.2139725}{doi:10.1080/03605302.2022.2139725};
\href{https://arxiv.org/abs/2206.15309}{arXiv:2206.15309}.

\bibitem{WWZ2023}
J. Wei, L. Wu and L. Zhang,
Estimates of bubbling sequences of \(SU(3)\) Toda systems at critical
parameters: Part 2,
\emph{J. London Math. Soc.} \textbf{107} (2023), 2150--2196.
\href{https://doi.org/10.1112/jlms.12736}{doi:10.1112/jlms.12736}.

\bibitem{WZ2020}
L. Wu and L. Zhang,
Estimates of bubbling solutions of \(SU(3)\) Toda systems at critical
parameters. Part 1,
\emph{SIAM J. Math. Anal.} \textbf{52} (2020), 3738--3776.
\href{https://doi.org/10.1137/20M131031X}{doi:10.1137/20M131031X}.

\bibitem{Zhang2009}
L. Zhang,
Asymptotic behavior of blowup solutions for elliptic equations with
exponential nonlinearity and singular data,
\emph{Commun. Contemp. Math.} \textbf{11} (2009), 395--411.
\href{https://doi.org/10.1142/S0219199709003417}
{doi:10.1142/S0219199709003417}.

\end{thebibliography}
\end{document}